\documentclass[twoside]{report} 
\usepackage{latexsym}
\usepackage{amsfonts}
\usepackage{amssymb}
\usepackage{graphicx}

\begin{document}

\pagenumbering{roman}

\thispagestyle{empty}

\vspace*{2cm}

\begin{center}
{\huge Rank Equalities Related to Generalized  \\ 
Inverses of Matrices and Their Applications}
\end{center}

\vspace*{2cm}

\begin{center}
{\LARGE Yongge Tian}
\end{center}

\newpage

\pagenumbering{roman}
\vspace*{0.3cm}
\begin{center} {\Large{\bf ACKNOWLEDGEMENTS}}\end{center} 
\addcontentsline{toc}{part}{ACKNOWLEDGEMENTS}
\vspace*{0.5cm}

Many thanks go to my thesis advisor
George P. H. Styan for his guidance, invaluable discussions and thorough
understanding during the preparation of this thesis. 

\medskip 

Sincere thanks also go to Professor  Ronald J. Stern 
for his advice and encouragement.

\newpage

\vspace*{0.3cm}
\begin{center} {\Large {\bf ABSTRACT }} \end{center} 
\addcontentsline{toc}{part}{ABSTRACT}

\vspace*{1cm}

This paper is divided into two parts. In the first part, we develop a general
 method for expressing ranks of matrix expressions that involve Moore-Penrose 
inverses, group inverses, Drazin inverses, as well as  weighted  Moore-Penrose 
inverses of matrices. Through this method we establish a variety of valuable 
rank equalities related to generalized inverses of matrices mentioned above. 
Using them, we characterize many matrix equalities in the theory of 
generalized inverses of matrices and their applications. 

In the second part, we consider maximal and minimal possible ranks of matrix expressions that 
involve variant matrices, the fundamental work is concerning extreme ranks of  the two linear matrix expressions $A - BXC$ and  $A - B_1X_1C_1 - B_2X_2C_2$. 
As applications, we present a wide range of their consequences and 
applications in matrix theory. \\


{\em Key  words}: Rank, inner inverse, Moore-Penrose inverse, 
 group inverse, Drazin inverse, reverse order law, EP matrix,   block matrix,  Schur complement, matrix expression,  matrix equation,  solution.

\tableofcontents

\newpage

\pagenumbering{arabic}

\setcounter{page}{1}

\pagestyle{myheadings}

\markboth{YONGGE  TIAN }
{1. INTRODUCTION AND PRELIMINARIES}

\chapter{\bf Introduction and preliminaries}

\noindent This is a comprehensive work on ranks of matrix expressions
involving Moore-Penrose inverses,  group inverses, Drazin inverses, 
as well as weighted Moore-Penrose inverses. In the theory of generalized inverses 
of matrices and their applications, there are numerous matrix expressions and
  equalities that involve these three kinds of
 generalized inverses of matrices. Now we propose such a problem: Let
  $p(\, A_1^{\dagger}, \ \cdots, \ A_k^{\dagger} \,)$ and  $q(\, B_1^{\dagger},
  \ \cdots, \ B_l^{\dagger} \,)$ be two matrix expressions involving 
  Moore-Penrose inverses of matrices. Then determine necessary and sufficient conditions
   such that
 $p(\, A_1^{\dagger}, \ \cdots, \ A_k^{\dagger} \,) = q(\, B_1^{\dagger}, \
 \cdots, \ B_l^{\dagger} \,)$ holds. A seemingly trivial condition for this equality to hold is apparently 
 $$ 
 {\rm rank\,}[ \, p(\, A_1^{\dagger}, \ \cdots, \ A_k^{\dagger} \,) -
 q(\, B_1^{\dagger}, \ \cdots, \ B_l^{\dagger} \,) \,] = 0. \eqno (1.1)
 $$
 However, if we can reasonably find a formula for expressing the rank
 of the left-hand side of (1.1), then we can derive immediately from (1.1)
 nontrivial conditions for 
$$
p(\, A_1^{\dagger}, \ \cdots, \ A_k^{\dagger} \,) =
  q(\, B_1^{\dagger}, \ \cdots, \ B_l^{\dagger} \,)
$$ 
to hold. This work has a far-reaching influence to many problems in the theory of 
generalized inverses of matrices and their applications. 
This consideration motivates us to make a thorough investigation to this work.  In fact, the author has successfully 
used this idea to establish necessary and sufficient conditions 
such that
 $(ABC)^{\dagger} = C^{\dagger}B^{\dagger}A^{\dagger}$,  $(ABC)^{\dagger} =
 (BC)^{\dagger}B(AB)^{\dagger},$ and $ (A_1 A_2 \cdots \ A_k )^{\dagger}
 = A_k^{\dagger} \cdots A_2^{\dagger}A_1^{\dagger}$ (cf. \cite{Ti2}, \cite{Ti4}). But the methods used in 
those papers are somewhat restricted and not applicable to various kinds of matrix expressions. In this 
paper, we shall develop a general and complete method for establishing rank equalities for matrix 
expressions involving Moore-Penrose inverses, group inverses,  Drazin inverses, as well as weighted 
Moore-Penrose inverses of matrices. Using these rank formulas, we shall characterize various equalities for generalized inverses of matrices, and then present their applications in the theory of generalized inverses of matrices.

The matrices considered in the this paper are mainly over the complex number field ${\cal C}.$
Let $ A \in {\cal C }^{ m \times n }$. We use $ A^*$,  $r( A ) $ and $ R( A )$ to stand for the conjugate transpose, the rank and the range (column space) of $ A $, respectively.

It is well known that the Moore-Penrose inverse of matrix $ A $ is defined to be the unique solution $ X $ of 
the following four Penrose equations
$$
(1)  \ \ \ AXA = A,  \ \ \ \  (2)  \ \ \ XAX = X,  \ \  \ \ (3) \ \ \  (AX)^* = AX,  
 \ \ \ \  (4) \ \ \ (XA) ^* = XA, 
$$  
and is often denoted by $ X = A^{\dagger}$. In addition, a matrix $ X $ that
 satisfies  the first equation above is called an inner inverse of $ A $, and often  denoted by 
$A^-$.  A matrix $ X $ that satisfies  the second equation above is called an outer inverse of $ A $, and often  denoted by $A^{(2)}$. For simplicity,  we use $ E_A $ and  $ F_A $ to stand for 
the two projectors 
$$
 E_A = I - AA^{\dagger} \qquad  {\rm and} \qquad  F_A = I - A^{\dagger}A
$$ 
induced by $ A $. As to various basic properties concerning Moore-Penrose inverses of matrices, see, e.g., Ben-Israel and Greville \cite{BeG} (1980), Campbell and Meyer \cite{CM2} (1991), 
Rao and Mitra \cite{RM} (1971).   

Let $ A \in {\cal C }^{ m \times m}$ be given with ${\rm Ind\,}A = k$, the smallest positive 
integer such that $r(A^{k+1}) = r(A^k)$. The Drazin inverse of matrix $ A $ is defined to be 
the unique solution $ X $ of  the following three  equations
$$
(1)  \ \ \ A^kXA = A^k,   \ \ \ \ \  (2)  \ \ \ XAX = X,   \ \ \ \ \ (3) \ \ \  AX = XA,   
$$  
and is often denoted by $ X = A^{D}$. In particular, when ${\rm Ind\,}A = 1$, the Drazin inverse 
of matrix $ A $ is call  the group inverse of $ A $, and is often  denoted by 
$A^{\#}$.

Let  $ A \in {\cal C}^{m \times n}$. The weighted Moore-Penrose inverse of 
$ A \in {\cal C}^{m \times n} $ with respect to the two positive definite
 matrices $ M \in {\cal C}^{m \times m} $ and $N \in {\cal C}^{n \times n}$
 is defined to be the unique solution of the following four matrix equations
$$
(1)  \ \ \ AXA = A,  \ \ \ \ \ (2)  \ \ \ XAX = X,  \ \ \ \ \ (3) \ \ \  (MAX)^* = MAX,,  
 \ \ \ \ \ (4) \ \ \  (NXA)^* = NXA,
$$  
and this $ X $ is often denoted by $ X = A^{\dagger}_{M,N}$. In particular, when
 $ M = I_m$ and $ N = I_n$, $A^{\dagger}_{M,N}$ is the conventional Moore-Penrose inverse $ A^{\dagger}$ of $ A$. Various basic properties concerning Drazin inverses,   group inverses  and 
Weighted Moore-Penrose inverses of matrices can be found  in Ben-Israel and 
Greville \cite{BeG}, Campbell and Meyer \cite{CM2},  Rao and Mitra \cite{RM}.

It is well known that generalized inverses of matrices are a powerful tool for
establishing various rank equalities  for matrices. Two seminal references
are  the paper \cite{MaS1} by Marsaglia and Styan (1974)  and the paper \cite{MJ4} by Meyer (1973). 
In those two papers, some fundamental rank equalities and inequalities related to generalized inverses of block matrices were established 
and a variety of consequences and applications of these rank equalities and inequalities
 were considered. Since then, the main results in those two papers have widely been applied to 
dealing  with various problems in the theory of generalized inverses of matrices and its 
applications. To some extent, this paper could be regarded as a summary and extension of all 
work related to those two remarkable papers.

We next list some key results in those  papers, which will be intensively applied in this monograph.

\medskip

\noindent {\bf Lemma 1.1} (Marsaglia and Styan \cite{MaS1}, 
Meyer \cite{MJ4}).\, {\em Let $ A \in {\cal C}^{ m \times n },
 \, B \in {\cal C}^{ m \times k}, \,  C \in {\cal C}^{ l \times n }$ and $
  D \in {\cal C}^{ l \times k } $. Then
$$ 
\displaylines{ 
\hspace*{1cm} r[ \, A, \ B \, ] = r( A ) + r(\, B - AA^{\dagger}B \, ) = r( B ) 
+ r(\, A - BB^{\dagger}A \, ),
   \hfill (1.2)
\cr
\hspace*{1cm} r \left[ \begin{array}{c}  A  \\ C  \end{array} \right]  = r( A ) + 
r(\, C - CA^{\dagger}A \, ) = r( C ) + r(\, A - AC^{\dagger}C \, ),  \hfill (1.3)
 \cr 
\hspace*{1cm} 
  r\left[ \begin{array}{cc}  A  & B  \\ C &  0 \end{array} \right]  =
   r( B ) + r( C ) + r( E_BAF_C ) = r( B ) + r( C ) + r[\,(\, I_m - BB^{\dagger}\, )A(\,I_n  - C^{\dagger}C  \,)\, ],
 \hfill (1.4)
\cr
\hspace*{1cm}  r \left[ \begin{array}{cc}  A  & B  \\ C &  D \end{array} \right]
  = r( A ) + r \left[ \begin{array}{cc}  0  &  E_A B \\  CF_A  &  S_A
  \end{array} \right],  \hfill (1.5)
\cr
\hspace*{1cm}  
  r \left[ \begin{array}{cc}  A  & B  \\ C &  D \end{array} \right]
  = \left[ \begin{array}{cc} A \\ C
  \end{array} \right] + r[ \, A, \ B \,]  -r(A) + r[J(D)], \hfill (1.6)
\cr}
$$ 
where $ S_A = D - CA^{\dagger}B $ is the Schur complement of $A $ in
 $ M =  \left[ \begin{array}{cc}  A  & B  \\ C &  D \end{array} \right],$
and 
$$ \displaylines{ 
\hspace*{1cm} 
J(D) =[ \, I - (CF_A)(CF_A)^{\dagger}\,]
S_A[ \, I - (E_AB)^{\dagger}(E_AB) \,],\hfill 
\cr }
$$
called the rank complement of $ D $  in  $ M $.  In particular$,$ if $ R( B ) \subseteq R ( A )$ and 
$ R( C^* )  \subseteq R ( A^* ),$ then 
 $$\displaylines{ 
\hspace*{1cm}  
  r \left[ \begin{array}{cc}  A  & B  \\ C &  D \end{array} \right]
  = r( A ) + r( \, D - CA^{\dagger}B \, ). \hfill (1.7)
\cr}
$$ 
 The six rank equalities in {\rm (1.2)}---{\rm (1.7)} are also
 true when replacing $A^{\dagger}$ by any inner inverse $ A^-$ of $ A$.
} 

\medskip

\noindent {\bf Lemma 1.2}\cite{MaS1}.\, {\em Let $ A \in {\cal C}^{ m \times n },
 \, B \in {\cal C}^{ m \times k}, \,  C \in {\cal C}^{ l \times n }$ and $
  D \in {\cal C}^{ l \times k } $. Then

{\rm (a)} \ $ r[ \, A, \ B \, ] = r( A ) + r( B ) \Leftrightarrow 
R(A) \cap R(B) = \{ 0 \} \Leftrightarrow R[(E_AB)^*] = R(B^*)
\Leftrightarrow R[(E_BA)^*] = R(A^*).$

{\rm (b)} \ $ r \left[ \begin{array}{c}  A  \\ C  \end{array} \right]  \
= r( A ) +  r(C)  \Leftrightarrow R(A^*) \cap R(C^*) = \{ 0 \}
\Leftrightarrow R(CF_A) = R(C) \Leftrightarrow R(AF_C) = R(A). $

{\rm (c)} \ $ r[ \, A, \ B \, ] = r(A)   \Leftrightarrow R(B) \subseteq R(A) \Leftrightarrow E_AB =0.$

{\rm (d)} \ $ r \left[ \begin{array}{c}  A  \\ C  \end{array} \right]  \
= r( A )  \Leftrightarrow R(C^*) \subseteq  R(A^*)
\Leftrightarrow CF_A = 0.$

{\rm (e)} \ $ r \left[ \begin{array}{cc}  A & B   \\ C  & 0  \end{array}
\right] = r( A ) + r(B) + r(C)$ $ \Leftrightarrow $ $  R(A) \cap R(B) = \{ 0 \}$
and $R(A^*) \cap R(C^*) = \{ 0 \}. $

{\rm (f)} \ $ r \left[ \begin{array}{cc}  A & B   \\ C  & D  \end{array}
\right] = r( A )$ $ \Leftrightarrow $ $  R(B) \subseteq R(A)$ and
$R(C^*) \subseteq R(A^*)$ and $ D = CA^{\dagger}B. $ } 

\medskip

\noindent {\bf Lemma 1.3}\cite{MaS1}({\em Rank cancellation rules}).\, {\em Let  $ A \in {\cal C}^{m \times n}, \  B \in
 {\cal C}^{ m \times k}$ and $ C \in {\cal C}^{l \times n }$ be given$,$ and 
suppose  that 
$$
\displaylines{ 
\hspace*{2cm}
R(AQ) = R(A) \ \ \ \   and  \ \ \ \ R[(PA)^*] = R(A^*).\hfill 
\cr 
\hspace*{0cm}
Then \hfill
\cr
\hspace*{2cm}
 r[\, AQ, \ B \, ] = r[\, A, \ B \, ],  \qquad  r \left[ \begin{array}{c}
 PA  \\ C  \end{array} \right] =
 r \left[ \begin{array}{c}  A  \\ C  \end{array} \right].  \hfill (1.8)
\cr
\hspace*{0cm}
In \ particular, \hfill
\cr 
\hspace*{2cm}
 r[\, AA^*, \ B \, ] = r[\, A, \ B\, ],  \qquad r \left[ \begin{array}{c}
  A^*A  \\ C  \end{array} \right] =
 r \left[ \begin{array}{c}  A  \\ C  \end{array} \right].  \hfill (1.9)
 \cr}
$$}
{\bf Lemma 1.4}\cite{MaS1}.\,  {\em Let $ A,  \, B  \in
{\cal C}^{m \times n }.$ Then 

{\rm (a)} \ $ r(\, A \pm B \, ) \geq  r \left[ \begin{array}{c}  A  \\ B  
\end{array} \right]  + r[\, A, \ B \, ] - r(A) - r(B).$ 

{\rm (b)} \ If $ R(A) \cap R(B) = \{ 0 \},$  then $ r(\, A + B\, )
= r \left[ \begin{array}{c}  A  \\ B  \end{array} \right].$

{\rm (c)} \, If $ R(A^*) \cap R(B^*) = \{ 0 \},$  then $  r( \, A + B \,) = r[\, A, \ B \, ].$

{\rm (d)} \  $r(\, A + B \, ) = r(A) + r(B ) \Leftrightarrow  r(\, A - B \, ) = r(A) + r(B ) 
\Leftrightarrow  R(A) \cap R(B) = \{ 0 \},$ and $ R(A^*) \cap R(B^*) = \{ 0 \}.$ }

\medskip

In addition, we shall also use in the sequel the following several basic rank formulas, which 
are either well known or  easy to prove.

\medskip

\noindent {\bf Lemma 1.5.}\, {\em Let  $ A \in {\cal C}^{m \times n}, \,  B \in
 {\cal C}^{ n \times m}$ and $ N \in {\cal C}^{m \times m}.$ Then  
$$
\displaylines{ 
\hspace*{2cm} r(\, A - ABA \, ) = r(A)  + r( \, I_n - BA\,) - n  = r(A)
  + r( \, I_m -  AB \,) - m,  \hfill (1.10)
\cr
\hspace*{2cm}
 r(\, N \pm N^{k+1} \, ) = r(N)  + r( \, I_n \pm N^k \,) - m,  \ \ \ 
for \ all \  k \geq 1,  \hfill (1.11)
\cr
\hspace*{2cm} 
r( \, I_m - N^2 \,) =  r( \, I_m + N \,)  +  r( \, I_m -  N \,) - m.
 \hfill (1.12)
\cr}
$$}
{\bf Lemma 1.6.}\, {\em Let $ A,  \ B  \in
{\cal C}^{m \times n }$. Then
$$\displaylines{ 
\hspace*{2cm}
 r \left[ \begin{array}{cc}  A & B  \\ B & A  \end{array} \right]
= r( \, A + B \, ) + r( \, A - B \, ). \hfill (1.13)
\cr}
$$}
{\bf Proof.}\, Follows from the following decomposition 
$$
\frac{1}{2} \left[ \begin{array}{cc}  I_m & I_m \\ I_m & -I_m  \end{array}
\right] \left[ \begin{array}{cc}  A & B  \\ B & A  \end{array} \right]
\left[ \begin{array}{cc}  I_n & I_n  \\ I_n & -I_n  \end{array} \right]
 = \left[ \begin{array}{cc}  A+B & 0  \\ 0 & A - B \end{array} \right].
$$ 
{\bf Lemma 1.7.}\, {\em Let  $ A \in {\cal C}^{m \times m}.$ Then  
$$
r( \, A - A^{2k+1} \,) =  r( \, A + A^k \,)  +  r( \, A -  A^k \,) - r(A)
= r(A) +  r( \, I_m + A^{k-1} \,)  +  r( \, I_m -  A^{k-1} \,) - 2m.
 \eqno (1.14)
$$
In particular (Anderson and Styan \cite{AnS}) 
$$\displaylines{
\hspace*{1cm}
r( \, A - A^3 \,) =  r( \, A + A^2 \,)  +  r( \, A -  A^2 \,) - r(A)
= r(A) +  r( \, I_m + A \,)  +  r( \, I_m -  A \,) - 2m.
 \hfill (1.15)
\cr}
$$}
{\bf Proof.}\, Replace $ B$ in (1.13) by $ A^k$ and simplify to yield (1.14). \qquad  $\square$  

\medskip

\noindent {\bf Lemma 1.8.}\, {\em Let  $ A \in {\cal C}^{m \times m}$ and  $  \lambda_1, \, \lambda_2,\, 
 \cdots, \, \lambda_k \in {\cal C}$ with $ \lambda_i \neq  \lambda_j$  for $ i \neq j $. Then for any positive integer $t_1, t_2, \cdots, t_k,$ the following rank equality
$$ 
\displaylines{
\hspace*{1cm}
r[\, ( \lambda_1I - A)^{t_1}  ( \lambda_2I - A)^{t_2} \, \cdots \, ( \lambda_kI - A)^{t_k}\,] \hfill
\cr
\hspace*{1cm}
 =  r[\, ( \, \lambda_1I - A \,)^{t_1}] +  r[( \, \lambda_1I - A \,)^{t_2}] + \cdots +  r[( \, \lambda_kI - A \,)^{t_k}\, ]
 - (k-1)m.
\hfill (1.16)
\cr}
$$
always holds.  This result can also alternatively be stated that for any two polynomials
 $p(x)$ and $ q(x)$ without common roots$,$ there is 
$$
r[\, p(A)q(A) \,] = r[\, p(A) \,] +  r[\, q(A) \,] - m.  \eqno (1.17) 
$$}

This paper is divided into two parts with 31 chapters. They organize as follows.

In Chapter 2, we establish  several universal rank formulas for matrix expressions that 
involve Moore-Penrose inverses of matrices. These rank formulas will serve as a basic tool for developing the content in the subsequent chapters.

In Chapter 3, we present  a set of  rank formulas related to sums, differences 
and  products of idempotent matrices. Based on them, we shall reveal a series of new and nontrivial 
properties for idempotent matrices.

In Chapter 4, we extend the results in Chapter 3 to some matrix expressions that 
involve both idempotent matrices and general matrices. In addition, we shall also establish 
 a group of new rank formulas related to involutory matrices and then consider their consequences.

In Chapter 5, we establish a set of rank formulas related to outer inverses of a  
matrix. Some of them will be applied in the subsequent chapters.  

In Chapter 6, we examine various relationships between a matrix and its Moore-Penrose inverse
using the rank equalities obtained in the preceding chapters. We also consider in the chapter how 
characterize some special types of matrices, such as, EP matrix, conjugate EP matrix, bi-EP matrix, 
star-dagger matrix, power-EP matrix, and so on.

In Chapter 7, we discuss various rank equalities for matrix expressions that involve two or 
more Moore-Penrose inverses, and then use them to characterize various matrix equalities that 
involve Moore-Penrose inverses.

In Chapter 8, we investigate various kind of reverse order laws for Moore-Penrose inverses 
of products of two or three matrices using the rank equalities established in the preceding  chapters.

In Chapter 9, we investigate Moore-Penrose inverses of $ 2 \times 2$ block matrices, as well as 
 $ n \times n$ block matrices using the rank equalities established in the preceding  chapters.

In Chapter 10, we investigate  Moore-Penrose inverses 
of sums of matrices using the rank equalities established in the preceding  chapters.

In Chapter 11, we study the relationships between Moore-Penrose inverses of block 
circulant matrices and  sums of matrices. Based on them and the results in Chapter 9, we shall 
present a group of expressions for Moore-Penrose inverses of sums of matrices.

In Chapter 12, we present a group of formulas for expressing ranks of submatrices 
in the  Moore-Penrose inverse of a matrix.

In Chapters 13---17, our work is concerned with rank equalities for Drazin inverses, 
  group inverses,  and weighted Moore-Penrose inverses of matrices and their applications.
Various kinds of problems examined in Chapters 6---12 for Moore-Penrose inverses of matrices are almost 
considered  in these five chapters for Drazin inverses, group inverses, and weighted Moore-Penrose 
inverses of matrices.   

In Chapter 18, we present maximal and minimal ranks of the matrix expression 
$A - BXC$ with respect to the variant matrix $ X$, and then consider rank and range  
invariance of $A - BXC$ with respect to the variant matrix $ X$. In addition, 
we also consider shorted matrices of a matrix relative to a given matrix set.

In Chapter 19,  we determine maximal and minimal ranks of the matrix expression 
$A - B_1X_1C_1 - B_2X_2C_2 $ with respect to the variant matrices $ X_1$ and $X_2$  under
 the conditions $R (B_1) \subseteq R (B_2)$ and $R (C_2^T) \subseteq R(C_1^T)$.

In Chapter 20,  we determine maximal and minimal ranks of the matrix expression 
$A_1 - B_1XC_1$ subject to a consistent matrix equation $ B_2XC_2 = A_2$, and consider some related topics.
 
In Chapters 21---23,  we present maximal and minimal ranks of the Schur complement 
$D - CA^-B$  with respect to $ A^-$ and then consider various related topics, including 
some problems on generalized inverses of sums and products of matrices. 

In Chapters 24---26, we determine extreme ranks of submatrices in solutions of the matrix equation $ BXC = A$, and extreme  ranks of two real 
matrices $ X_0$ and $ X_1$ in solutions to the complex matrix equation $ B( \, X_0 + X_1\, )C = A$, as well as  extreme ranks of solutions of the matrix
equation $ B_1XC_1 + B_2YC_2 = A$.       

In Chapters 27---31, we consider extreme ranks of some general matrix expressions. The basic work is concerning extreme 
ranks of $A - B_1X_1C_1 - B_2X_2C_2$ with respect to the two independent variant matrices $ X_1$ and $X_2$. The work 
is partially extended to some linear matrices expressions with more than two independent variant matrices. Quite a lot of 
consequences are derived from these results, including extreme ranks of $ A_1 - B_1XC_1$ subject to a pair of 
consistent matrix equations $B_2XC_2= A_2$ and $B_3XC_3= A_3$;   extreme ranks of $ A - BX - XC $ subject to a 
consistent matrix equations $BXC= D$; extreme ranks of a quadratic matrix expression 
$ A - (\, A_1 - B_1X_1C_1\,)D(\, A_2 - B_2X_2C_2\,)$ with respect to $ X_1 $ and $ X_2$. In addition, we also present many rank formulas for 
matrix expressions involving generalized inverses of matrices in these chapters.

\pagestyle{myheadings}
\markboth{YONGGE  TIAN }
{2. BASIC RANK FORMULAS}

\chapter{Basic rank formulas for Moore-Penrose inverses}

The first and most fundamental rank formula used in the sequel is given below. 

\medskip

\noindent {\bf Theorem 2.1.}\, {\em Let $ A \in {\cal C}^{ m \times n },
\, B \in  {\cal C}^{ m \times k}, \,  C \in {\cal C}^{ l \times n }$ and $  D
 \in {\cal C}^{ l \times k } $ be given. Then the rank of the Schur complement
 $ S_A = D - CA^{\dagger}B $ satisfies the equality 
 $$\displaylines{ 
\hspace*{1.5cm} 
  r( \,  D - CA^{\dagger}B \, ) = r \left[ \begin{array}{cc}  A^*AA^*
  & A^*B  \\ CA^* &  D \end{array} \right] - r( A ).  \hfill (2.1)
 \cr}
$$ }
{\bf Proof.}\, It is obvious that 
$$\displaylines{ 
\hspace*{1.5cm}
R( A^*B ) \subseteq R( A^* ) =  R( A^*AA^* ), \ \  {\rm and} \ \  R( AC^*) \subseteq R( A )
 = R (AA^*A ). \hfill
\cr}
$$ 
Then it follows by (1.7) and  a well-known basic property 
$A^*( A^*AA^* )^{\dagger}A^* = A^{\dagger}$(see \cite{RM} pp. 69) that
$$ \displaylines{ 
\hspace*{1.5cm} 
 r \left[ \begin{array}{cc}  A^*AA^*  & A^*B  \\ CA^* &  D \end{array} \right] 
= r( A^*AA^* ) + r[ \, A -
  C A^*( A^*AA^* )^{\dagger}A^*B \, ] =  r ( A ) + r( \,  D - CA^{\dagger}B \, ), \hfill
\cr}
$$ 
establishing (2.1). \qquad  $ \Box$ 

\medskip

 The significance of (2.1) is in that the rank of the Schur complement $
 S_A = D - CA^{\dagger}B $ can be evaluated by a block matrix formed by $ A, \ B, \ C$ and $ D $ in it, where
 no restrictions are imposed  on $ S_A$ and no Moore-Penrose inverses appear
 in the right-hand side of (2.1). Thus (2.1) in fact provides us a
 powerful tool to express ranks of matrix expressions that involve
 Moore-Penrose inverses of matrices.

Eq.\,(2.1) can be extended to various general formulas. We next present some of 
them, which will widely be used in the sequel. 

\medskip

\noindent  {\bf Theorem 2.2.}\, {\em Let $ A_1, \,  A_2, \, B_1, \, B_2, \, C_1,
 \ C_2 $ and $D $ are matrices such that expression
 $ D - C_1A_1^{\dagger}B_1 - C_2A_2^{\dagger}B_2 $ is defined. Then
 $$ \displaylines{ 
\hspace*{1.5cm}
 r( \, D - C_1A_1^{\dagger}B_1 - C_2A_2^{\dagger}B_2 \, ) =
  r \left[ \begin{array}{ccc}  A_1^*A_1A_1^*
  & 0 & A_1^*B_1   \\  0 & A_2^*A_2 A_2^*  & A_2^*B_2   \\ C_1 A_1^*  &
  C_2A_2^* &  D \end{array} \right] - r( A_1 ) -  r( A_2 ).  \hfill (2.2)
 \cr
\hspace*{0cm}
 In particular, \ if \hfill
\cr
\hspace*{1.5cm}
 R(B_1) \subseteq R(A_1),  \ \ \  R(C_1^*) \subseteq R(A_1^*),
 \ \  R(B_2) \subseteq R(A_2) \ \ \  and   \ \ \ R(C_2^*) \subseteq R(A_2^*), \hfill
\cr
\hspace*{0cm}
then \hfill
\cr
\hspace*{1.5cm} 
 r( \, D - C_1A_1^{\dagger}B_1 - C_2A_2^{\dagger}B_2 \, ) =
 r \left[ \begin{array}{ccc}  A_1 & 0 & B_1   \\  0 & A_2  & B_2   \\ C_1  &
  C_2 &  D \end{array} \right] -r( A_1 ) -  r( A_2 ).
  \hfill (2.3)
\cr}
$$}
\hspace*{0.3cm}Let 
$$ \displaylines{ 
\hspace*{1.5cm} 
 C = [\, C_1, \ C_2 \, ], \qquad  B =  \left[ \begin{array}{c}  B_1  \\ B_2  
\end{array} \right] \qquad {\rm and}  \qquad  A = \left[ \begin{array}{cc} 
 A_1 & 0 \\ 0 &  A_2 \end{array} \right].\hfill 
\cr}
$$
Then  (2.1) can be written as  (2.2), and  (2.3) follows from  (1.6). 
 \qquad  $ \Box$ 

\medskip

If the matrices in  (2.2) satisfy certain conditions, the block matrix in
 (2.2) can easily be reduced to some simpler forms. Below are some of
 them. 

\medskip

\noindent  {\bf Corollary 2.3.}\, {\em Let $ A \in {\cal C}^{ m \times n },
\, B \in {\cal C}^{ m \times k}, \,  C \in {\cal C}^{ l \times n }$ and $
  N \in {\cal C}^{ k \times l} $ be given. Then
  $$ \displaylines{ 
\hspace*{1.5cm}
 r( \, N^{\dagger} - CA^{\dagger}B \, ) = r \left[ \begin{array}{ccc}
 AA^*A - A(BNC)^*A & AC^* & 0  \\
   B^*A &  0  & N   \\ 0  &  N  &  0  \end{array} \right] - r( A ) -  r(N).
   \hfill (2.4)
\cr} 
$$     
In particular$,$ if $ R( B^*A ) \subseteq R(N )$ and $ R(CA^*) \subseteq R( N^* ),$
then
$$
\displaylines{ 
\hspace*{2cm}
 r( \, N^{\dagger} - CA^{\dagger}B \, ) =  r[ \, AA^*A - A(BNC)^*A \, ]
 +  r(N) - r( A ). \hfill (2.5)
\cr}
$$ 
If $ R( B^*A ) \subseteq R(N ), \ R(CA^*) \subseteq R( N^*), \
  R( BN ) \subseteq R(A)$  and \ $ R[ (NC)^*] \subseteq R(A^*),$ then 
$$
\displaylines{ 
\hspace*{1.5cm} 
r( \, N^{\dagger} - CA^{\dagger}B \, ) =  r( \, A - BNC \, ) +  r(N) -
 r( A ). \hfill (2.6)
\cr}
$$}
\hspace*{0cm}{\bf Theorem 2.4.}\, {\em Let $ A_t, \, B_t, \, B_t, \, C_t( t = 1, \, 2 , \,
 \cdots , \, k )$ and $D $ are matrices such that expression
 $ D - C_1A_1^{\dagger}B_1 - \cdots - C_kA_k^{\dagger}B_k $ is defined.
 Then
 $$ \displaylines{ 
\hspace*{1.5cm}
 r( \, D - C_1A_1^{\dagger}B_1 -  \cdots - C_kA_k^{\dagger}B_k \, )
 = r \left[ \begin{array}{cc} A^*AA^*  & A^*B   \\ C A^*  &  D \end{array}
 \right] - r( A ),  \hfill (2.7)
\cr}
 $$
where $ A = {\rm diag}( \, A_1, \,  A_2, \, \cdots, \, A_k \, ), \
 B^* = [ \, B_1^*, \,  B_2^*, \, \cdots, \, B_k^* \, ] $ and $
 C = [ \, C_1, \,  C_2, \, \cdots, \, C_k \, ].$  } 

\medskip

\noindent  {\bf Theorem 2.5.}\, {\em Let $A, \, B, \, C, \, D, \, P  $ and $Q $ are
 matrices such that expression $ D - CP^{\dagger}A Q^{\dagger}B $ is defined.
 Then
 $$ 
\displaylines{ 
\hspace*{1.5cm}
 r( \, D - CP^{\dagger}A Q^{\dagger}B \, ) = r \left[ \begin{array}{ccc}
 P^*AQ^*  & P^*PP^*  & 0   \\
     Q^*QQ^* & 0 & Q^*B   \\ 0 &  CP^* &  -D \end{array} \right] -
  r(P) -  r(Q).  \hfill (2.8)
 \cr
In \ particular, \ if \hfill
\cr
\hspace*{1.5cm}
R(A) \subseteq R(P),   \ \ \ R(A^*) \subseteq R(Q^*),  \ \ \  R(B)
 \subseteq R(Q) \ \ \ and  \ \ \ R(C^*) \subseteq R(P^*),\hfill 
\cr
\hspace*{0cm}
then \hfill
\cr
\hspace*{1.5cm}
 r( \, D - CP^{\dagger}A Q^{\dagger}B \, ) = r \left[ \begin{array}{ccc} A
 & P  & 0 \\ Q & 0 & B \\ 0 &  C &  -D \end{array} \right] - r(P) - r(Q).
 \hfill (2.9)
 \cr}
$$} 
\hspace*{0cm}{\bf Proof.}\,  Note that 
 \begin{eqnarray*}
 r( \, D - CP^{\dagger}A Q^{\dagger}B \, )
 &  = & r \left[ \begin{array}{cc} A  & A Q^{\dagger}B \\ CP^{\dagger}A  &
 D \end{array} \right] - r( A ) \\
 & = & r \left( \, \left[ \begin{array}{cc} A & 0 \\ 0  &  D \end{array}
 \right] + \left[ \begin{array}{cc} A & 0 \\ 0  & C \end{array} \right]
 \left[ \begin{array}{cc} 0 & P \\ Q  &  0 \end{array} \right]^{\dagger}
  \left[ \begin{array}{cc} A & 0 \\ 0  &  B \end{array} \right] \right) -
  r( A ).
 \end{eqnarray*}
Applying (2.1) to it  and then simplifying yields  (2.8). Eq.(2.9) is 
derived from (2.8) by the rank cancellation law (1.8). \qquad$ \Box$ 

\medskip

\noindent {\bf Theorem 2.6.}\, {\em  Suppose that the matrix expression
 $ S = D - C_1P_1^{\dagger}A_1 Q_1^{\dagger}B_1 - C_2P_2^{\dagger}A_2
 Q_2^{\dagger}B_2 $ is defined.  Then
 $$ \displaylines{ 
\hspace*{1.5cm}
 r( S ) = r \left[ \begin{array}{ccccc}  P_1^*A_1Q_1^*  & 0  & P_1^*P_1P_1^*
  & 0 & 0  \\
   0  &   P_2^*A_2Q_2^* &  0 &  P_2^*P_2P_2^* & 0 \\
  Q_1^*Q_1Q_1^* & 0 & 0 & 0 & Q_1^*B_1   \\
  0 & Q_2^*Q_2Q_2^*   & 0 & 0 &  Q_2^*B_2  \\
 0 & 0 & C_1P_1^* & C_2P_2^* & -D \end{array} \right] - d,  \hfill (2.10)
 \cr}
$$
 where $ d = r(P_1) + r(P_2) + r(Q_1) + r(Q_2).$  In particular$,$ if
$$ 
\displaylines{ 
\hspace*{1.5cm}
R(A_i) \subseteq R(P_i), \ \ \ R(A^*_i) \subseteq R(Q^*_i), \ \ \ R(B_i)
\subseteq R(Q_i) \ \ and  \ \ R(C^*_i) \subseteq R(P^*_i), i = 1, \ 2, \hfill
\cr
\hspace*{0cm}
then \hfill
\cr
\hspace*{1.5cm}
 r( S ) = r \left[ \begin{array}{ccccc} A_1  & 0  & P_1  & 0 & 0  \\
  0  & A_2 &  0 &  P_2 & 0 \\
  Q_1 & 0 & 0 & 0 & B_1   \\
  0 & Q_2  & 0 & 0 & B_2  \\
 0 & 0 & C_1 & C_2 & -D \end{array} \right] - r(P_1) -  r(Q_1) - r(P_2) -
 r(Q_2).  \hfill (2.11)
\cr
\hspace*{0cm}
Moreover, \hfill
\cr
\hspace*{0.8cm}
 r( \, D^{\dagger} - CP^{\dagger}A Q^{\dagger}B \, ) =
  r \left[ \begin{array}{cccc}  D^*DD^*  & 0 & 0   & D^* \\ 0 &  P^*A Q^*
  & P^*PP^* & 0  \\ 0 & Q^*QQ^* & 0 & Q^*B  \\ D^* & 0 & CP^* & 0
   \end{array} \right] -
  r(P) -  r(Q) - r(D).  \hfill (2.12)
\cr} 
$$}
\hspace*{0cm}{\bf Proof.}\,  Writing $ S $ as  
 $$\displaylines{ 
\hspace*{2cm}
 S = D - [\, C_1, \ C_2 \,] \left[ \begin{array}{cc} P_1 & 0 \\ 0  &  P_2
 \end{array} \right]^{\dagger}
 \left[ \begin{array}{cc} A_1 & 0 \\ 0  & A_2 \end{array} \right] 
 \left[ \begin{array}{cc} Q_1 & 0 \\ 0  &  Q_2 \end{array} \right]^{\dagger}
  \left[ \begin{array}{c} B_1 \\ B_2 \end{array} \right], \hfill 
\cr}
$$
and then applying  (2.8) to it produce  (2.10). Eq.\,(2.11) is derived
  from  (2.10) by the rank cancellation law (1.8). Eq.\,(2.12) is a special 
case of (2.10). \qquad  $ \Box$ 

\medskip

It is easy to see that a general rank formula for 
 $$\displaylines{ 
\hspace*{2cm}
 D - C_1P_1^{\dagger}A_1Q_1^{\dagger}B_1 -
  C_2P_2^{\dagger}A_2 Q_2^{\dagger}B_2 - \cdots -
  C_kP_2^{\dagger}A_kQ_k^{\dagger}B_k   \hfill 
\cr}
 $$ 
 can also be established by the similar method for deriving  (2.10).
 As to some other general matrix expressions, such  as
 $$ \displaylines{ 
\hspace*{2cm}
 S_k = A_0P_1^{\dagger}A_1P_2^{\dagger}A_2 \cdots P_k^{\dagger}A_k \hfill 
\cr}
 $$ 
 and their linear combinations, the formulas for expressing their ranks can
 also be established. However they are quite tedious in form, we do not
 intend to give them here.

\markboth{YONGGE  TIAN }
{3.  RANK EQUALITIES FOR IDEMPOTENT MATRICES}

\chapter{Rank equalities for idempotent matrices}

\noindent A square matrix $ A $ is said to be idempotent if $A^2$ = $ A $.  If we consider 
it as  a matrix equation $A^2$ = $ A $, then its general solution can be written as 
$A = V( V^2)^{\dagger}V,$ where $V $ is an arbitrary square complex matrix. This assertion can easily be 
verified. In fact, $A = V( V^2)^{\dagger}V,$  apparently satisfies $A^2$ = $ A $. Now for any 
matrix $ A $ with $A^2$ = $ A $, we let $ V = A $. Then $V( V^2)^{\dagger}V =  
A( A^2)^{\dagger}A= A A^{\dagger}A= A$. Thus $A = V( V^2)^{\dagger}V$ is indeed the general 
solution the  idempotent equation $A^2$ = $A$. This fact clearly implies that any matrix 
expression that involves idempotent matrices could be regarded as a conventional  matrix 
expression that involves  Moore-Penrose inverses of matrices. Thus the formulas in  Chapter 2 
are all applicable to determine ranks of  matrix expressions that involve 
idempotent matrices. However because of speciality of idempotent matrices, the rank equalities 
related to idempotent matrices can also be deduced by various elementary methods. The results 
in the chapter are originally derived by the rank formulas in Chapter 2, we later also find 
some elementary methods to establish them. So we only show these results in these elementary 
methods inn this chapter. 

\medskip

\noindent {\bf Theorem 3.1.}\, {\em Let $ P, \, Q \in {\cal C}^{m \times m}$
 be two idempotent matrices. Then the difference $ P-Q $ satisfies the rank equalities  
 $$
\displaylines{ 
\hspace*{2cm}
 r( \, P - Q \, ) = r \left[ \begin{array}{c}  P \\ Q \end{array} \right] + 
r[ \, P,  \ Q \, ] - r(P) -  r(Q), \hfill (3.1) 
\cr
\hspace*{2cm}
 r( \, P - Q \, ) =  r( \, P - PQ \, ) + r(\, PQ - Q \, ), \hfill (3.2)
 \cr
\hspace*{2cm}
 r( \, P - Q \, ) = r( \, P - QP \, ) + r(\, QP - Q \, ).  \hfill (3.3)
\cr }
$$ }
{\bf Proof.}\, Let $ M = \left[ \begin{array}{ccc} -P  & 0  & P  
\\ 0 & Q  & Q  \\ P & Q & 0  \end{array} \right] $. Then it is easy to see by
 block elementary operations  of matrices that
 $$ \displaylines{ 
\hspace*{2cm}
 r(M ) = r \left[ \begin{array}{ccc} -P & 0 & 0  \\ 0  & Q  & 0  \\ 0 & 0 & P-Q 
\end{array} \right] = r(P) + r(Q) + r( \, P - Q \,).  \hfill 
 \cr}
$$
 On the other hand, note that $ P^2 = P$ and $ Q^2= Q$. It is also easy to find by
 block elementary operations of matrices that 
$$\displaylines{ 
\hspace*{2cm}
 r(M)  = r \left[ \begin{array}{ccc} -P  & 0  & P  \\
 -QP & 0  & Q  \\ P & Q & 0  \end{array} \right] =  r \left[ \begin{array}{ccc} 
0  & 0 & P \\
    0 & 0 & Q  \\ P & Q  & 0  \end{array} \right] = r \left[ \begin{array}{c}  
P \\ Q \end{array} \right] +
  r[ \, P,  \ Q \,].  \hfill 
\cr}
$$
Combining the above two equalities yields 
(3.1). Consequently applying 
(1.2) and (1.3)  to
 $[ \, P,  \ Q \,] $ and $  \left[ \begin{array}{c}  P \\ Q \end{array}
 \right] $ in  (3.1)
 respectively yields 
 $$ 
  r[ \, P, \ Q \, ] = r( P ) + r(\, Q - PQ \, ),  \eqno (3.4) 
 $$
 $$
 r[ \, P, \ Q \, ] = r(Q) + r(\, P - QP \, ), \eqno (3.5) 
 $$
 $$
  r \left[ \begin{array}{c}  P  \\ Q  \end{array} \right]  = r( P ) + r(\, Q - QP \,), \eqno (3.6)
 $$ 
 $$
  r \left[ \begin{array}{c}  P  \\ Q  \end{array} \right] = r( Q ) + r( \, P - PQ \,). \eqno (3.7)
 $$ 
 Putting  (3.4) and (3.7) in  (3.1) produces  (3.2), putting (3.5) and (3.6) 
in (3.1) produces (3.3).  \qquad  $ \Box$ 

\medskip

\noindent {\bf Corollary 3.2.}\, {\em Let $ P, \, Q \in {\cal C}^{m \times m}$
 be two idempotent matrices.  Then   

{\rm (a)} \ $R(\, P - PQ \, ) \cap R( \, PQ - Q \,) = \{ 0 \}$ and
 $ R[ (\, P - PQ\, )^* ] \cap R[ ( \, PQ - Q \,)^*] = \{ 0 \}.$

{\rm (b)} \  $ R(\, P - QP\, ) \cap R( \, QP - Q \,) = \{ 0 \}$ and
 $ R[ (\, P - QP \, )^*] \cap R[ ( \, QP - Q \,)^*] = \{ 0 \}.$

{\rm (c)} \  If $ PQ = 0 $ or $ QP = 0,$ then $ r( \, P - Q \, ) =
  r(P) + r(Q),$ i.e.$,$ $ R(P) \cap R(Q) = \{ 0 \}$ and   $ R(P^*) \cap R(Q^*) = \{ 0 \}.$ 

{\rm (d)} \  If both $ P $ and $ Q $ are Hermitian idempotent$,$ then 
$r( \, P - Q \, ) = 2r[ \, P,  \ Q \, ] - r(P) -  r(Q).$ 
 }

\medskip

\noindent {\bf Proof.}\, Parts (a) and (b) follows from applying Lemma 1.4(d) to 
(3.2) and (3.3). Part (c) is a direct consequence of  (3.2) and (3.3).
Part (d) follows from  (3.1). \qquad  $ \Box$ 

\medskip

On the basis of  (3.1), we can easily deduce the following known result due to Hartwig and Styan \cite{HaSt2} on 
the rank subtractivity  two idempotent matrices.    

\medskip

\noindent {\bf Corollary 3.3.}\, {\em Let $ P, \, Q \in {\cal C}^{m \times m}$
 be two idempotent matrices.  Then the following statements are equivalent$:$  

 {\rm (a)} \  $r( \, P - Q \, ) = r(P) - r(Q),$ i.e.$,$  $ Q \leq_{rs} P$.  

 {\rm (b)} \ $ r \left[ \begin{array}{c}  P \\ Q \end{array} \right] = r[ \, P,  \ Q \, ]
 = r(P). $

 {\rm (c)} \ $R(Q) \subseteq R(P)$ and $ R(Q^*) \subseteq R(P^*)$.

 {\rm (d)} \ $ PQ = QP = Q.$  

 {\rm (e)} \ $ PQP = Q.$  }

\medskip

\noindent {\bf Proof.}\, The equivalence of Parts (a) and (b) follows immediately from applying (3.1). 
The equivalence  of Parts (b), (c), (d) and (e) can trivially be verified by  (1.2) and (1.3). \qquad  $ \Box$ 

\medskip

\noindent {\bf Corollary 3.4.}\, {\em Let $ P, \, Q \in {\cal C}^{m \times m}$ be two idempotent matrices. 
   Then the following statements are equivalent$:$  

{\rm (a)} \ The difference $ P - Q $ is nonsingular. 

 {\rm (b)} \ $ r \left[ \begin{array}{c}  P \\ Q \end{array} \right] = r[ \, P,  \ Q \, ] = r(P) + r(Q) = m.$ 

{\rm (c)} \ $ R(P) \oplus R(Q) = R(P^*) \oplus R(Q^*)  = {\cal C}^m.$}

\medskip 

\noindent {\bf Proof.}\,  Follows directly from  (3.1). \qquad  $ \Box$ 

\medskip

Notice that if a matrix  $P$ is  idempotent, the $ I_m - P$ is also idempotent. Thus replacing 
$ P $ in (3.1) by $ I_m - P$, we get the following.     

\medskip

\noindent {\bf Theorem 3.5.}\, {\em Let $ P, \, Q \in {\cal C}^{m \times m}$ be two idempotent matrices. Then the rank of $ I_m - P - Q $ satisfies the equalities
 $$
\displaylines{ 
\hspace*{2cm}
 r( \, I_m  - P - Q \, ) = r(PQ) + r(QP ) - r(P) - r(Q) + m,  \hfill (3.8)
\cr
\hspace*{2cm}
r( \, I_m  - P - Q \, ) = r( \, I_m  - P - Q  + PQ \, ) + r(PQ),  \hfill (3.9)
\cr
\hspace*{2cm}
r( \, I_m  - P - Q \, ) = r( \, I_m  - P - Q  + QP \, ) + r(QP).  \hfill (3.10)
\cr}
$$}
{\bf Proof.}\, Replacing $ P $ in  (3.1) by
$I_m - P$ yields
$$ \displaylines{ 
\hspace*{2cm}
 r( \, I_m  - P - Q \, ) = r \left[ \begin{array}{c} I_m - P \\ Q
 \end{array} \right] + r[ \,I_m  - P,  \ Q \, ] - r( \,I_m  -P \,) -
  r(Q). \hfill (3.11) 
 \cr}
$$ 
It follows  by  (1.2) and (1.3) that 
 $$ \displaylines{ 
\hspace*{2cm}
 r[ \,I_m  - P,  \ Q \, ] = r( \,I_m  - P \, ) +
 r[ \, Q -  (\, I_m  - P \,)Q \, ] = m - r(P) + r(PQ), \hfill 
\cr
\hspace*{0cm}
and \hfill
\cr 
\hspace*{2cm}
 r \left[ \begin{array}{c} I_m - P \\ Q \end{array} \right]
 = r( \,I_m  - P \, ) +
  r[ \, Q -  Q(\, I_m  - P \,)\, ] = m - r(P) + r(QP).\hfill 
\cr} 
$$
Putting them in (3.11) produces (3.8). On the other hand, replacing 
$ P $ in  (3.2) and (3.3) by
  $I_m  - P$  produces
 \begin{eqnarray*}
 r[ \, (I_m  - P) - Q \, ] & = & r[ \, (I_m - P) - (I_m  - P)Q \, ]
 + r[ \, (I_m  - P)Q - Q \, ] \\
 & = & r( \, I_m - P - Q + PQ \, ) + r( PQ ), 
 \end{eqnarray*}
and
 \begin{eqnarray*}
 r[ \, (I_m  - P) - Q \, ] & = & r[ \, (I_m - P) - Q(I_m  - P) \, ]
 + r[ \, Q(I_m  - P) - Q \, ] \\
 & = & r( \, I_m - P - Q + QP \, ) + r(QP), 
 \end{eqnarray*}
both of which are exactly (3.9) and (3.10). \qquad  $ \Box$ 

\medskip

\noindent {\bf Corollary 3.6.}\, {\em Let $ P, \, Q \in {\cal C}^{m \times m}$
 be two idempotent matrices.  Then

{\rm (a)} \ $ R(\, I_m - P - Q + PQ \, ) \cap R(PQ) = \{ 0 \}$ and
 $ R[ (\,I_m -P - Q + PQ \, )^* ] \cap R[(PQ )^*] = \{ 0 \}.$

{\rm (b)} \ $ R(\, I_m - P - Q + QP \, ) \cap R(QP) = \{ 0 \}$ and
 $ R[ (\,I_m -P - Q + QP \, )^*] \cap R[(QP)^*] = \{ 0 \}.$

{\rm (c)} \ $ P + Q = I_m  \ \Leftrightarrow \ PQ = QP = 0 \ and
  \  R(P) \oplus R(Q) =  R(P^*) \oplus R(Q^*) = {\cal C}^m.$

{\rm (d)} \ If $ PQ = QP = 0,$ then $ r( \, I_m - P - Q \, ) =  m - r(P) -
r(Q)$.

{\rm (e)} \ $ I_m - P - Q $ is nonsingular if and only if
$ r(PQ) = r(QP ) = r(P) =r(Q)$.

{\rm (f)} \ If both $P$ and $ Q $ are Hermitian idempotent$,$ then 
$r( \, I_m - P - Q \, ) = 2r(PQ) -r(P) - r(Q) + m $.} 

\medskip

\noindent{\bf Proof.}\, Parts (a) and (b) follow from applying Lemma 1.4(d) to 
 (3.9) and (3.10).  Note from  (3.8)---(3.10) that $ P + Q = I_m $ is 
equivalent to $ PQ = QP = 0$ and  $r(P ) + r(Q) = m$. This assertion is also 
equivalent to $ PQ = QP = 0$ and  $R(P) \oplus R(Q) =  R(P^*) \oplus R(Q^*) 
= {\cal C}^m,$  which is Part (c). Parts (d), (e) and (f) follow from (3.8). 
\qquad  $ \Box$ 

\medskip

As for the rank of sum of two idempotent matrices, we have the following several results.

\medskip
  
\noindent {\bf Theorem  3.7.}\, {\em Let $ P, \, Q \in {\cal C}^{m \times m}
$ be  two idempotent matrices. Then the sum $ P + Q $ satisfies the rank equalities
 $$
\displaylines{ 
\hspace*{2cm}
 r( \, P + Q \,) = r\left[ \begin{array}{cc} P & Q \\ Q & 0 \end{array}
 \right] - r(Q) = r\left[ \begin{array}{cc} Q & P \\ P & 0 \end{array}
  \right] - r(P),   \hfill (3.12) 
\cr
\hspace*{2cm} 
 r( \, P + Q \,) = r( \, P - PQ - QP + QPQ \, ) + r( Q ), \hfill (3.13)
\cr
\hspace*{2cm}
 r( \, P + Q \,) =  r( \, Q - PQ - QP + PQP \, ) + r(P). \hfill (3.14)
\cr }
$$}
{\bf Proof.}\, Let $ M = \left[ \begin{array}{ccc} P  & 0  & P  
\\ 0 & Q  & Q  \\ P & Q & 0  \end{array} \right] $. Then it is easy to see by block elementary operations of matrices that
 $$ \displaylines{ 
\hspace*{2cm}
 r(M ) = r \left[ \begin{array}{ccc} P & 0 & 0  \\ 0  & Q  & 0  \\ 0 & 0 & -P-Q \end{array} \right] 
= r(P) + r(Q) + r( \, P + Q \,). \hfill 
 \cr}
$$
 On the other hand, note that $ P^2 = P$ and $ Q^2= Q$. It is also easy to find by block elementary operations  of matrices that 
 $$\displaylines{ 
\hspace*{0.5cm}
 r(M) =  r \left[ \begin{array}{ccc} P  & 0  & P  \\
 -QP & 0  & Q  \\ P & Q & 0  \end{array} \right]  = r \left[ \begin{array}{ccc} 2P  & 0 & P \\
    0 & 0 & Q  \\ P & Q  & 0  \end{array} \right] = r \left[ \begin{array}{ccc} 2P 
 & 0 & 0 \\ 0 & 0 & Q  \\ 0 & Q  & \frac{1}{2}P  \end{array} \right] = r \left[ \begin{array}{cc}
  P  & Q \\ Q & 0 \end{array} \right] + r(P), \hfill
\cr
\hspace*{0cm}
and \hfill
\cr
\hspace*{0.5cm}
 r(M)  = r \left[ \begin{array}{ccc} 0  & -PQ  & P  \\ 0  & Q  & Q  \\ P & Q & 0  \end{array} \right] 
= r \left[ \begin{array}{ccc} 0  & 0 & P \\
    0 & 2Q & Q  \\ P & Q  & 0  \end{array} \right] = r \left[ \begin{array}{ccc} 0  & 0 & P \\
    0 & 2Q & 0  \\ P & 0  &  \frac{1}{2}Q  \end{array} \right]
    = r \left[ \begin{array}{cc}  Q  & P \\  P & 0  \end{array} \right]
     + r(Q). \hfill
\cr}
$$
The combination of the above three rank equalities yields the two equalities in  (3.12). 
Consequently applying (1.4) to the two block matrices
 in (3.12) yields  (3.13) and (3.14), respectively. \qquad  $ \Box$ 

\medskip

\noindent {\bf Corollary 3.8.}\, {\em Let $ P, \, Q \in {\cal C}^{m \times m}$
 be two idempotent matrices. 

{\rm (a)}\, If $ PQ = QP,$  then
$$\displaylines{ 
\hspace*{2cm}
r( \, P + Q \,) = r[\, P ,\ Q \,] =
r \left[ \begin{array}{c} P \\ Q \end{array} \right], \hfill (3.15)
\cr
\hspace*{0cm}
or equivalently, \hfill
\cr
\hspace*{2cm}
R(Q) \subseteq R( P + Q )  \ \ and  \ \  R(Q^*) \subseteq R( P^* + Q^* ).
\hfill (3.16)
\cr}
$$

{\rm (b)}\, If $ R( Q) \subseteq R(P)$ or $ R(Q^*) \subseteq R(P^*),$  then
 $r( \, P + Q \,) = r(P)$.  }  

\medskip

\noindent {\bf Proof.}\, If $PQ= QP$, then (3.13) and (3.14) reduce to
$$\displaylines{ 
\hspace*{2cm}
r( \, P + Q \, ) = r(\, P - PQ \, ) + r(Q) = r( \, Q - PQ \, ) + r(P). \hfill
\cr}
$$
Combining them with  (3.4) and (3.7) yields  (3.15). The equivalence of
 (3.15) and (3.16) follows from a simple fact that
$$\displaylines{ 
\hspace*{2cm}
r \left[ \begin{array}{c} P \\ Q \end{array} \right]
= r \left[ \begin{array}{c}  P + Q  \\ Q \end{array} \right] \ \ {\rm and} \ \
r[ \, P, \  Q \, ] = r[\, P + Q, \ Q \, ], \hfill 
\cr}
$$
as well as  Lemma 1.2(c) and (d). The result in Part (b) follows immediately from
 (3.12). \qquad  $\Box$ 

\medskip

\noindent {\bf Corollary 3.9.}\, {\em Let $ P, \, Q \in {\cal C}^{m \times m}$
 be two idempotent matrices. Then the following five statements are
 equivalent$:$

 {\rm (a)} \  The sum $ P+Q $ is nonsingular.

 {\rm (b)} \ $ r \left[ \begin{array}{c}  P \\ Q \end{array} \right] = m \ \
 and  \ \ R \left[ \begin{array}{c}  P \\ Q \end{array} \right] \cap R
 \left[ \begin{array}{c}  Q \\ 0  \end{array} \right] = \{ 0 \} $.

 {\rm (c)} \ $ r[\, P, \ Q \,] = m \ \ and  \ \ R \left[ \begin{array}{c}
  P^* \\ Q^* \end{array} \right] \cap R \left[ \begin{array}{c}  Q^* \\ 0
   \end{array} \right] = \{ 0 \}  $.

 {\rm (d)} \ $ r \left[ \begin{array}{c}  Q \\ P \end{array} \right] = m \ \
  and  \ \ R \left[ \begin{array}{c}  Q \\ P \end{array} \right] \cap
  R \left[ \begin{array}{c} P \\ 0  \end{array} \right] = \{ 0 \}$.

 {\rm (e)} \ $ r[\, Q, \ P \,] = m \ \ and  \ \ R \left[ \begin{array}{c}
 Q^* \\ P^* \end{array} \right] \cap R \left[ \begin{array}{c}  P^* \\ 0
  \end{array} \right] = \{ 0 \} $. } 

\medskip

 \noindent {\bf Proof.}\,  In light of  (3.12), the sum $ P + Q $ is
 nonsingular if and only if
 $$\displaylines{ 
\hspace*{2cm}
  r\left[ \begin{array}{cc} P & Q \\ Q & 0 \end{array} \right] = r(Q) + m,
  \hfill (3.17)
\cr
\hspace*{0cm}
or \ equivalently \hfill
\cr 
\hspace*{2cm}
 r\left[ \begin{array}{cc} Q & P \\ P & 0 \end{array} \right] = r(P) + m.
  \hfill (3.18)
 \cr
\hspace*{0cm}
Observe \  that \hfill
\cr 
\hspace*{2cm}
 r\left[ \begin{array}{cc} P & Q \\ Q & 0 \end{array} \right] \leq
 r \left[ \begin{array}{c}  P \\ Q \end{array} \right] +
 r \left[ \begin{array}{c}  Q \\ 0 \end{array} \right] \leq m + r(Q),  \hfill
 \cr
\hspace*{2cm} 
r\left[ \begin{array}{cc} P & Q \\ Q & 0 \end{array} \right]
   \leq r[ \,  P, \ Q \, ] +
 r [\, Q, \ 0 \,] \leq m + r(Q), \hfill
  \cr
\hspace*{2cm} 
 r\left[ \begin{array}{cc} Q & P \\ P & 0 \end{array} \right] \leq
  r \left[ \begin{array}{c}  Q \\ P \end{array} \right] +
  r \left[ \begin{array}{c}  P \\ 0 \end{array} \right] \leq m + r(P), \hfill
 \cr
\hspace*{2cm} 
  r\left[ \begin{array}{cc} Q & P \\ P & 0 \end{array} \right]
  \leq r[ \,  Q, \ P \, ] +
 r [\, P, \ 0 \,] \leq m + r(P). \hfill
 \cr}
$$
 Combining them with  (3.17) and (3.18) yields the equivalence of
  Parts (a)---(e). \qquad $ \Box$ 

\medskip

 \noindent {\bf Theorem 3.10.}\, {\em Let $ P, \, Q \in
 {\cal C}^{m \times m} $ be two idempotent matrices. Then 

{\rm (a)}\, The rank of $ I_m  + P - Q $ satisfies the  equality
 $$\displaylines{ 
\hspace*{2cm} 
 r( \, I_m  + P - Q \,) = r(QPQ) - r(Q) + m. \hfill (3.19) 
 \cr}
$$

{\rm (b)}\, The rank of $ 2I_m  -  P - Q $ satisfies the
 two equalities
 $$
\displaylines{ 
\hspace*{2cm} 
 r( \, 2I_m - P - Q \,)=  r( \, Q - QPQ \, ) - r(Q) + m, \hfill (3.20)
 \cr
\hspace*{2cm} 
 r( \, 2I_m - P - Q \,) = r( \, P - PQP \, ) - r(P) + m.  \hfill (3.21) 
\cr}
$$}
{\bf Proof.}\, Replacing $ Q $ in (3.12) by the idempotent
 matrix  $ I_m  - Q$  and  applying (1.4) to it yields
 \begin{eqnarray*}
 r( \, I_m  + P - Q \,)  &= & r\left[ \begin{array}{cc} P &  I_m  - Q \\ I_m  - Q & 0 \end{array} \right] 
 - r( \,I_m  - Q \,) \\
 & = & r( \,I_m  - Q \,)  + r[ (\, I_m -( \,I_m  - Q \,) \, )P(\, I_m - ( \,I_m  - Q \,) \,) ] \\
 & = & m  - r(Q)  + r(QPQ), 
 \end{eqnarray*} 
 establishing (3.19). Further, replacing $ P $ and $ Q $ in (3.12) by
 $ I_m  - P$ and  $ I_m  - Q$, we also by (1.4) find that
 \begin{eqnarray*}
 r( \, 2I_m  - P - Q \,)  &= & r\left[ \begin{array}{cc} I_m  - P &  I_m  - Q \\ I_m  - 
Q & 0 \end{array} \right] - r( \,I_m  - Q \,) \\
 & = & r( \,I_m  - Q \,) +
 r[ \, (\, I_m -( \,I_m  - Q \,) \, )( \,I_m  - P \,)(\, I_m - ( \,I_m  - Q \,) \,)\, ] \\
 & = & m  - r(Q)  + r(\, Q - QPQ \,), 
 \end{eqnarray*}
 establishing (3.20). Similarly, we can show (3.21). \qquad  $ \Box$ 

\medskip

 \noindent {\bf Corollary 3.11.}\, {\em Let $ P, \, Q \in
 {\cal C}^{m \times m} $ be two idempotent matrices. 

 {\rm (a)} \ If $ R(P) \subseteq R(Q) $ and $R(P^*) \subseteq R(Q^*),$ then $ P $ and $Q $ satisfy the two rank equalities
$$\displaylines{ 
\hspace*{2cm}
 r(\, I_m + P - Q \, ) = m  + r(P) - r(Q), \hfill (3.22)
\cr
\hspace*{2cm}
 r(\, 2I_m - P - Q \, ) = m + r( \,Q - P \,) - r(Q). \hfill (3.23)
\cr}
$$ 

{\rm (b)} \ $ I_m  + P - Q $ is nonsingular $ \Leftrightarrow$  $r(QPQ) = r(Q).$

{\rm (c)} \ $ 2I_m  - P - Q $ is nonsingular $\Leftrightarrow$
$r(\, P - PQP \,) = r(P)$  $ \Leftrightarrow$ $ r(\, Q - QPQ \,) = r(Q).$

{\rm (d)} \ $ Q - P = I_m $ $ \Leftrightarrow$ $r(QPQ) + r(Q ) = m.$ } 

\medskip

\noindent {\bf Proof.}\, The two conditions $ R(P) \subseteq R(Q) $ and $R(P^*) \subseteq R(Q^*)$ 
are equivalent to $ QP = P = PQ $. In that case,  (3.19) reduces to 
 (3.22), (3.20) and (3.21) reduce to  (3.23). The results in Parts (a)---(c) are direct
consequences of (3.19).  \qquad $ \Box$

\medskip

we next consider the rank of $ PQ - QP$ for two idempotent matrices $P $ and $ Q $.

\medskip

\noindent {\bf Theorem  3.12.}\, {\em Let $ P, \, Q \in {\cal C}^{m \times m}
$ be two idempotent matrices. Then the  difference $ PQ-QP $ satisfies the five rank equalities
$$   
\displaylines{ 
\hspace*{1cm}
 r( \, PQ - QP \, )  =  r( \, P - Q \, ) +  r( \, I_m -  P - Q \, ) - m,
 \hfill (3.24)
 \cr
 \hspace*{1cm} 
 r( \, PQ - QP \, )  = r( \, P - Q \, ) + r(PQ) + r(QP) - r(P) - r(Q),
 \hfill (3.25)
 \cr
\hspace*{1cm}
 r( \, PQ - QP \, ) =  r \left[ \begin{array}{c}  P \\ Q \end{array}
  \right] + r[  \, P, \ Q \, ] + r(PQ) + r(QP ) - 2r(P) -  2r(Q),
  \hfill (3.26)
 \cr
\hspace*{1cm}
 r( \, PQ - QP \, )  =  r( \, P - PQ \, ) + r(\,  PQ - Q \, )
 + r(PQ) + r(QP ) - r(P) - r(Q), 
\hfill  (3.27)
 \cr
\hspace*{1cm}
r( \, PQ - QP \, )  =  r( \, P - QP \, ) + r(\, QP - Q \, ) + r(PQ)
 + r(QP ) - r(P) - r(Q). \hfill (3.28)
\cr }
$$
In particular$,$ if both $P$ and $Q$ are Hermitian idempotent$,$ then 
 $$ \displaylines{ 
\hspace*{1cm}
 r( \, PQ - QP \, ) = 2r[\, P, \ Q \,] + 2r(PQ) - 2r(P) - 2r(Q).
 \hfill (3.29)
\cr}
$$}
{\bf Proof.}\, It is easy to verify that that $ PQ - QP = ( \, P - Q \, )
 ( \,P + Q - I_m \, )$. Thus  the rank of  $ PQ-QP $ can be expressed as
 $$ \displaylines{ 
\hspace*{1cm}
 r(\, PQ- QP \, ) = r [ \,( \, P - Q \, )( \, P + Q - I_m  \, ) \, ]
 = r\left[ \begin{array}{cc} I_m &   P + Q - I_m \\ P - Q  & 0
 \end{array} \right] - m.     \hfill (3.30)
\cr} 
$$ 
 On the other hand, it is easy to verify the factorization 
 $$\displaylines{ 
\hspace*{1cm}
 \left[ \begin{array}{cc} I_m &   P + Q - I_m \\ P - Q  & 0  \end{array}
 \right] =
 \left[ \begin{array}{cc} I_m &   2P- I_m \\ 0  & I_m  \end{array} \right]
 \left[ \begin{array}{cc} 0 &   P + Q - I_m \\ P - Q  & 0  \end{array} \right]
 \left[ \begin{array}{cc} I_m &  0 \\  2Q- I_m  & I_m  \end{array} \right].  \hfill 
 \cr
\hspace*{0cm}
Hence \hfill
\cr
\hspace*{1cm}
 r\left[ \begin{array}{cc} I_m &   P + Q - I_m \\ P - Q  & 0  \end{array}
 \right] =r( \, P - Q \, ) +  r( \, I_m -  P - Q \, ).  \hfill 
\cr}
$$ 
 Putting it in (3.30) yields  (3.24). Consequently putting (3.8) in
 (3.24) yields (3.25); putting  (3.1) in (3.25) yields  (3.26);
  putting  (3.2)  and (3.3) respectively in (3.25) yields  (3.27)
  and (3.28). \qquad  $ \Box$ 

\medskip

 \noindent {\bf Corollary 3.13.}\, {\em Let $ P, \, Q
 \in {\cal C}^{m \times m}$ be two idempotent matrices. Then the
 following five statements are equivalent$:$

 {\rm (a)} \ $ PQ = QP $.

 {\rm (b)} \  $ r( \, P - Q \, ) + r( \, I_m -  P - Q \, ) = m $.

 {\rm (c)} \ $ r(\, P - Q \, ) =  r(P) + r(Q) - r(PQ) -r(QP)$. 

 {\rm (d)} \ $ r( \, P - PQ \,) = r(P) -r(PQ)  \ \ and  \ \
 r( \, Q - PQ \,) = r(Q) -r(PQ), $ i.e.$,$  $ PQ \leq_{rs} P $ and $ PQ \leq_{rs} Q $.  

 {\rm (e)} \ $ r( \, P - QP \,) = r(P) -r(QP)  \ \ and  \ \
 r( \, Q - QP \,) = r(Q) -r(QP),$ i.e.$,$  $ QP \leq_{rs} P $ and $ QP \leq_{rs} Q $.

 {\rm (f)} \ $ r \left[ \begin{array}{c}
 P \\ Q \end{array} \right] = r(P) + r(Q) - r(PQ)  \ \ and \ \ 
 r[  \, P, \ Q \, ] = r(P) + r(Q) - r(QP). $

{\rm (g)} \ $ r \left[ \begin{array}{c}
 P \\ Q \end{array} \right] = r(P) + r(Q) - r(QP)  \ \ and \ \ 
 r[  \, P, \ Q \, ] = r(P) + r(Q) - r(PQ). $ } 

\medskip

\noindent {\bf Proof.}\, Follows immediately from  (3.24)---(3.28). \qquad
 $ \Box$

\medskip

 \noindent {\bf Corollary 3.14.}\, {\em Let $ P, \, Q
 \in {\cal C}^{m \times m}$ be two idempotent matrices. Then the following
 three statements are equivalent:

 {\rm (a)} \ $ r( \, PQ - QP \, ) = r( \, P - Q \, )$.

 {\rm (b)} \ $ I_m - P - Q $ is nonsingular.

{\rm (c)} \  $  r(PQ) =r(QP) =r(P) = r(Q)$. } 

\medskip

 \noindent {\bf Proof.}\, The equivalence of Parts (a) and (b) follows
 from (3.24). The  equivalence of Parts (b) and (c) follows from
 Corollary 3.6(e). \qquad  $ \Box$ 

\medskip

 \noindent {\bf Corollary 3.15.}\, {\em Let $ P, \, Q \in 
 {\cal C}^{m \times m}$ be two idempotent matrices. Then the following
 three statements are equivalent:

 {\rm (a)} \ $ PQ - QP$  is nonsingular.  

 {\rm (b)} \  $ P - Q $ and  $ I_m - P - Q $ are nonsingular. 

 {\rm (c)} \ $ R(P) \oplus R(Q) = R(P^*) \oplus R(Q^*) = {\cal C}^m $ and
 $ r(PQ) =r(QP) =r(P) = r(Q)$ hold.  } 

\medskip

 \noindent {\bf Proof.}\, The equivalence of Parts (a) and (b) follows
 from  (3.24). The  equivalence of Parts (b) and (c) follows from
  Corollaries 3.4(e) and 3.6(e). \qquad  $ \Box$ 

\medskip

A group of analogous rank equalities can also be derived for $PQ + QP$, where 
$P $ and $ Q $ are  two idempotent matrices $P $ and $ Q $.  

\medskip

\noindent {\bf Theorem  3.16.}\, {\em Let $ P, \, Q \in {\cal C}^{m \times m}
 $ be two idempotent matrices. Then $ PQ + QP $ satisfies the rank equalities
 $$
\displaylines{ 
\hspace*{1cm}
 r( \, PQ + QP \, ) =  r( \, P + Q \, ) +  r( \, I_m - P - Q \, ) - m ,
 \hfill (3.31)
 \cr
\hspace*{1cm}
 r( \, PQ + QP \, ) =  r( \, P + Q \, ) + r(PQ) + r(QP ) - r(P) - r(Q),
  \hfill (3.32)
 \cr
 \hspace*{1cm}
 r( \, PQ + QP \, ) = r( \, P - PQ - QP + QPQ \, ) + r(PQ) + r(QP ) - r(P),
  \hfill (3.33)
 \cr
\hspace*{1cm}
 r( \, PQ + QP \, ) = r( \, Q - PQ - QP  + PQP \, ) + r(PQ) + r(QP ) - r(Q).
  \hfill (3.34)
 \cr }
$$ } 
{\bf Proof.}\, Note that $ PQ +QP = (\,P  + Q \, )^2 -
 (\, P + Q \, )$. Then applying  (1.11) 
to it, we directly obtain (3.31). Consequently, putting (3.8) in
 (3.21) yields  (3.32), putting  (3.13) and (3.14) respectively
 in (3.32) yields (3.33) and (3.34). \qquad $ \Box$ 

 \medskip 









 \noindent {\bf Corollary 3.17.}\, {\em Let $ P, \, Q \in
 {\cal C}^{m \times m}$ be two idempotent matrices. Then the following four
 statements are equivalent$:$

 {\rm (a)} \ $ r(\,  PQ + QP \, ) =  r( \,P + Q  \, ) $.  

 {\rm (b)} \ $I_m - P - Q$ is nonsingular. 

 {\rm (c)} \ $ r(PQ) = r( QP) = r(P) = r(Q)$. 

 {\rm (d)} \ $ r(\,  PQ - QP \, ) =  r( \,P - Q \,) $.} 

\medskip

 \noindent {\bf Proof.}\, The equivalence of Parts (a) and (b) follows
 from ( 3.31),  and the  equivalence of  Parts (b)---(d) comes from
 Corollary 3.14.  \qquad  $\Box$ 

\medskip

 \noindent {\bf Corollary 3.18.}\, {\em Let $ P, \ Q \in
 {\cal C}^{m \times m}$ be two idempotent matrices. Then the following two
 statements are equivalent$:$

 {\rm (a)} \ $PQ + QP $ is nonsingular.  

 {\rm (b)} \  $ P + Q $ and  $ I_m - P - Q $ are nonsingular.  } 

\medskip

 \noindent {\bf Proof.}\,  Follows directly from (3.31). \qquad
 $ \Box$ 

\medskip

 Combining the two rank equalities in  (3.24) and (3.31), we obtain the
 following. 

\medskip

 \noindent {\bf Corollary 3.19.}\, {\em Let $ P, \, Q \in {\cal C}^{m \times m}
  $ be two idempotent matrices. Then both of them satisfy the following rank
  identity
 $$ \displaylines{ 
\hspace*{1cm}
 r(\, P + Q \,) + r( \, PQ - QP \, ) =  r(\, P - Q \,) + r( \, PQ + QP \, ).
  \hfill (3.35)
\cr} 
$$}
{\bf Thoerem 3.20.}\, {\em Let $ P, \ Q
 \in {\cal C}^{m \times m}$
 be two idempotent matrices. Then
$$ \displaylines{ 
\hspace*{1cm}
 r[\, (\, P - Q \,)^2 - ( \,  P - Q \, ) \,]
 = r(\,I_m - P + Q \,) + r(\, P - Q \, ) - m .  \hfill (3.36)
\cr 
\hspace*{1cm}
r[\, (\, P - Q \,)^2 - ( \,  P- Q  \, ) \,]
 =  r(PQP) - r(P) + r( \, P - Q \, ).  \hfill (3.37)
\cr }
$$ } 
{\bf Proof.}\, Eq.\,(3.36) is derived from (1.11).  
According to   (3.19), we have $r(\,I_m - P + Q \,) =  r(PQP) - r(P) + m$. 
Putting it in  (3.36) yields  (3.37). \qquad  $ \Box$

\medskip

 \noindent {\bf Corollary 3.21} (Hartwig and Styan \cite{HaSt2}).\, {\em Let $ P, \, Q \in
 {\cal C}^{m \times m}$ be two idempotent matrices. Then the following five
 statements are equivalent$:$

 {\rm (a)} \ $ P-Q$ is  idempotent.  

 {\rm (b)} \ $ r(\,I_m - P + Q \,)  = m - r(\, P - Q \, ).$

 {\rm (c)} \ $ r( \, P - Q \, ) = r(P) - r(Q),$ i.e.$,$  $Q \leq_{rs} P$. 

 {\rm (d)} \  $ R(Q) \subseteq R(P)$ and  $ R(Q^*) \subseteq R(P^*)$.

 {\rm (e)} \  $PQP = Q.$  } 

\medskip

 \noindent {\bf Proof.}\, The equivalence of Parts (a) and (b) follows
 immediately from  (3.36), and the  equivalence of Parts (c), (d) and (e)
  is from  Corollary 3.3(d). The equivalence of Parts (a) and (e) follows
  from a direct matrix computation.  \qquad  $ \Box$ 

\medskip

 In Chapter 4, we shall also establish a rank formula for
 $(\, P - Q \,)^3 - ( \, P - Q \, )$ and consider tripotency of
  $ P - Q$, where $ P, \ Q $ are  two idempotent  matrices. 

\medskip

 \noindent {\bf Theorem  3.22.}\, {\em Let $ P, \, Q \in
  {\cal C}^{m \times m} $ be two idempotent matrices. Then
  $ I_m - PQ$ satisfies the  rank equalities
 $$ \displaylines{ 
\hspace*{2cm}
  r( \, I_m - PQ \, ) = r( \, 2I_m - P - Q \, )= 
 r[ \, (\, I_m - P \,) +  (\, I_m - Q \, \, ) \, ].   \hfill (3.38)
\cr}
 $$}
{\bf Proof.}\, According to (1.10) we have
 $$\displaylines{ 
\hspace*{2cm}
  r( \, I_m - PQ \, ) = r( \, Q - QPQ \, ) -r(Q) + m.  \hfill 
 \cr}
$$ 
 Consequently putting  (3.20) in it yields  (3.38). \qquad   $ \Box$ 

\medskip

 \noindent {\bf Corollary 3.23.}\, {\em Let $ P, \, Q \in
 {\cal C}^{m \times m} $ be two idempotent matrices. Then the sum $ P + Q $
 satisfies the  rank identities
 $$ \displaylines{ 
\hspace*{2cm}
 r(\, P + Q \,) = r( \,  P + Q  - PQ \, ) = r( \,  P + Q  - QP \, ).
 \hfill (3.39)
 \cr}
$$ 
 In particular$,$ if $ PQ = QP,$ then 
 $$ \displaylines{ 
\hspace*{2cm}
 r(\, P + Q \,) = r( P ) + r(Q) -r( PQ ).  \hfill (3.40) 
 \cr}
$$ } 
{\bf Proof.}\, Replacing $ P $ and $ Q$ in (3.38) by  two
 idempotent matrices $ I_m - P$  and $ I_m - Q$ immediately yields (3.39).
 If $ PQ = QP,$ then we know by  (3.13) and (3.14) that
 $$ \displaylines{ 
\hspace*{2cm}
 r(\, P + Q \,) = r( \,  P - PQ \, ) + r( Q ) = r( \,  Q - QP \, ) + r(P),
 \hfill (3.41)
 \cr}
$$  
 and by Corollary 3.13 we also know that  
 $$ \displaylines{ 
\hspace*{2cm}
 r( \,  P - PQ \, ) = r(P) - r(PQ) \ \ 
{\rm and} \ \  r( \,  Q - QP \, )
 = r(Q) - r(QP). \hfill (3.42)
 \cr}
$$ 
Putting (3.42) in (3.41) yields  (3.40).  \qquad   $ \Box$ 

\medskip

 \noindent {\bf Corollary  3.24.}\, {\em Let $ P, \, Q \in
 {\cal C}^{m \times m} $ be two idempotent matrices. Then
 $$ \displaylines{ 
\hspace*{2cm}
 r[ \, PQ - (PQ)^2 \,] = r( \, I_m - PQ \, ) + r(PQ) - m
 = r( \, 2I_m - P - Q \, ) + r(PQ) - m .\hfill (3.43)
 \cr} 
$$
 In particular$,$ the following statements are equivalent$:$

 {\rm (a)} \ $PQ$ is idempotent.

 {\rm (b)} \ $r( \, I_m - PQ \, ) = m -  r(PQ)$. 

 {\rm (c)} \ $ r( \, 2I_m - P - Q \, )  = m - r(PQ) $. } 

\medskip

 \noindent {\bf Proof.} \ Applying  (1.11) to
  $PQ - (PQ)^2$ gives the first equality in (3.43). The second one
  follows from (3.38).  \qquad  $ \Box$ 

\medskip

\noindent {\bf Corollary 3.25.}\, {\em Let $ P, \, Q \in
 {\cal C}^{m \times m} $ be two idempotent matrices. Then
 $$ \displaylines{ 
\hspace*{2cm} 
 r( \, I_m - P - Q + PQ \,) = m -  r(P) -r(Q) + r(QP). \hfill 
\cr}
 $$ }
{\bf Proof.}\, This follows from replacing $ A $ in (1.4) by
 $ I_m $. \qquad  $ \Box$

\medskip

Notice that if a matrix $ A $ is idempotent, then $A^*$ is also idempotent. 
Thus we can easily find the following.

\medskip

\noindent {\bf Corollary 3.26.}\, {\em Let $ P \in
 {\cal C}^{m \times m}$ be an idempotent matrix. Then

 {\rm (a)} \ $ r( \, P - P^* \,) = 2r[\, P, \ P^* \,] - 2r(P)$.  

 {\rm (b)} \ $ r(\,I_m - P - P^* \,)  = r(\,I_m + P - P^* \,) = m$.   

 {\rm (c)} \ $ r( \, P + P^* \, ) = r( \, PP^* + P^*P \, ) = r[\, P, \ P^* \,],$ i.e.$,$ $R(P) \subseteq R( \, P + P^* \,)$ and $R(P^*) \subseteq R( \, P + P^*\,)$.

 {\rm (d)} \  $ r( \, PP^* - P^*P \,) = r(\, P - P^* \, ).$ }

\medskip

 \noindent {\bf Proof.}\, Part (a) follows from (3.1).  Part (b) follows from  (3.8) and (3.22). 
 Part (c) follows from  (3.31). Part (d) follows
  from  (3.24) and Part (b).  \qquad  $ \Box$ 

\medskip

The results in the preceding theorems and corollaries can easily be
 extended to matrices with properties
 $P^2 = \lambda P $ and $ Q^2 = \mu Q$, where $ \lambda \neq 0 $ and $ \mu 
 \neq 0 $. In fact, observe that
 $$ 
 \left( \, \frac{1}{\lambda}P \, \right)^2 = \frac{1}{\lambda^2}P^2 = \frac{1}{\lambda}P, \qquad
  \left( \, \frac{1}{\mu}Q \, \right)^2 = \frac{1}{\mu^2}Q^2 = \frac{1}{\mu}Q. 
 $$        
 Thus both $ P/\lambda$  and $Q/\mu$ are idempotent. In that case, applying the results in 
 the previous theorems and corollaries, one may establish a variety of  rank equalities and their consequences
  related to such kind of matrices. For example,
$$
\displaylines{ 
\hspace*{2cm} 
 r( \, \mu P - \lambda Q \, ) = r \left[ \begin{array}{c}  P \\ Q \end{array} \right] 
+ r[ \, P,  \ Q \, ] - r(P) -  r(Q), \hfill 
\cr
\hspace*{2cm}  
 r( \, \mu P + \lambda Q \, ) = r \left[ \begin{array}{cc}  P  & Q \\ Q & 0 \end{array} \right] - r(Q)
  = r \left[ \begin{array}{cc} Q & P \\ P & 0 \end{array} \right] - r(P), \hfill 
\cr
\hspace*{2cm}
 r(\, \lambda \mu I_m - \mu P - \lambda Q \, ) = r(PQ) + r(QP ) - r(P) - r(Q) + m, \hfill 
\cr
\hspace*{2cm} 
 r( \, PQ - QP \, )  =  r( \, \mu P - \lambda Q \, ) +  r( \, \lambda \mu I_m - \mu P - \lambda Q \, ) - m, \hfill 
\cr
\hspace*{2cm} 
 r( \, PQ + QP \, )  =  r( \, \mu P + \lambda Q \, ) +  r( \, \lambda \mu I_m - \mu P - \lambda Q \, ) - m, \hfill 
\cr
\hspace*{2cm} 
  r( \, \lambda \mu I_m - PQ \, ) = r(\, 2\lambda \mu I_m - \mu P - \lambda Q \, ),  \hfill 
\cr }
$$
and so on. We do not intend to present them in details.

\markboth{YONGGE  TIAN }
{4.  MORE RANK EQUALITIES FOR IDEMPOTENT MATRICES}

\chapter{More on rank equalities for idempotent matrices}

\noindent The rank equalities in Chapter 3 can partially be extended to matrix expressions that 
involve idempotent matrices and  general  matrices. In addition, they can also be applied to 
establish rank equalities related to involutory matrices. The corresponding results are 
presented in this chapter.

\medskip

\noindent {\bf Theorem  4.1.}\, {\em Let $ A \in  {\cal C}^{m \times n}$
 be given$,$ $ P \in {\cal C}^{m \times m}$ and $ Q \in {\cal C}^{n \times n} $
   be two idempotent matrices. Then the  difference $ PA-AQ $ satisfies
   the two rank equalities
 $$
\displaylines{ 
\hspace*{1.5cm} 
 r( \, PA - AQ \, ) = r \left[ \begin{array}{c}  PA \\ Q \end{array} \right]
  + r[ \, AQ, \ P \,] - r(P) -  r(Q), \hfill (4.1)
\cr
\hspace*{1.5cm}
 r( \, PA - AQ \, ) = r( \, PA - PAQ \, ) + r(\,PAQ  - AQ \, ).  \hfill (4.2)  
\cr }
$$  }
{\bf Proof.}\, Let $ M = \left[ \begin{array}{ccc} -P  & 0  & PA  
\\ 0 & Q  & Q  \\ P & AQ & 0  \end{array} \right] $. Then it is easy to see by the block 
elementary operations of matrices that
 $$\displaylines{ 
\hspace*{1.5cm}  
 r(M ) = r \left[ \begin{array}{ccc} -P & 0 & 0  \\ 0  & Q  & 0 \\ 0 & 0 & PA-AQ \end{array} \right] 
= r(P) + r(Q) + r( \, PA - AQ \,). \hfill (4.3)
\cr} 
$$
 On the other hand, note that $ P^2 = P$ and $ Q^2= Q$. It is also easy to
 find by  block elementary operations of matrices that
 $$\displaylines{ 
\hspace*{1.5cm} 
 r(M)  = r \left[ \begin{array}{ccc} 0  & PAQ  & PA  \\
 0  & Q  & Q  \\ P & AQ & 0  \end{array} \right] =  r \left[ \begin{array}{ccc} 0  & 0 & PA \\
    0 & 0 & Q  \\ P & AQ  & 0 \end{array} \right] = r \left[ \begin{array}{c} 
 PA \\ Q \end{array} \right] + r[ \, AQ,  \ P \,]. \hfill (4.4)
 \cr}
$$
 Combining  (4.3) and (4.4) yields  (4.1). Consequently applying
  (1.2) and (1.3)  to $[ \, AQ,  \ P \,] $ and
  $  \left[ \begin{array}{c}  PA \\ Q \end{array} \right] $ in (4.1)
  respectively yields  (4.2). \qquad  $ \Box$

\medskip

\noindent {\bf Corollary 4.2.}\, {\em Let $ A \in  {\cal C}^{m \times n}$
 be given$,$ $ P \in {\cal C}^{m \times m}$ and $ Q \in {\cal C}^{n \times n} $
   be two idempotent matrices. Then

{\rm (a)} \ $R( \, PA - PAQ\, ) \cap R( \, PAQ - AQ\,) = \{ 0\}$ and
 $R[( \, PA - PAQ\, )^*] \cap R[( \, PAQ - AQ\,)^*] = \{ 0\}.$

{\rm (b)} \ If $ PAQ = 0,$ then $ r( \, PA - AQ \, ) = r(PA) + r(AQ),$
 or$,$ equivalently $ R(PA)\cap R(AQ) = \{ 0 \}$ and
  $ R[(PA)^*] \cap R[(AQ)^*] = \{ 0 \}.$

{\rm (c)} \ $ PA = AQ \Leftrightarrow  PA(\, I  - Q \,) = 0   \ and \ ( \, I - P \,)AQ =0 \Leftrightarrow  R( AQ) \subseteq R(P)  \ and \  R[(PA)^*] \subseteq R(Q^*).$  
}

\medskip

\noindent {\bf Proof.}\, Part (a) follows from applying  Lemma 1.4(d) to 
 (4.2). Parts (b) and (c) are direct consequences  of (4.2). \qquad
 $ \Box $  

\medskip

\noindent {\bf Corollary 4.3.}\, {\em Let $ A, \, P, \, Q \in
 {\cal C}^{m \times m}$ be given with  $ P, \ Q $ being two idempotent matrices.
 Then the following three statements are equivalent$:$

 {\rm (a)} \ $ PA-AQ $ is nonsingular.

 {\rm (b)} \ $ r \left[ \begin{array}{c}  PA \\ Q \end{array} \right]
 = r[ \, AQ,  \ P \, ] = r(P) + r(Q) = m $.

 {\rm (c)} \ $ r(PA) = r(P), \ r(AQ) = r(Q)$ and
 $ R(AQ) \oplus R(P) = R[(PA)^*] \oplus R(Q^*)  = {\cal C}^m.$ } 

\medskip

\noindent {\bf Proof.}\, Follows from (4.1). \qquad  $ \Box $  

\medskip

Based on Corollary 4.2(c), we find an interesting result on the general solution of  a matrix equation. 
 
\medskip

\noindent {\bf Corollary 4.4.}\, {\em Let $ P \in {\cal C}^{m \times m}$ and $ Q \in {\cal C}^{n \times n} $ be two idempotent matrices. Then the general solution of the matrix equation 
$ PX =XQ $ can be written  in the two forms 
$$ \displaylines{ 
\hspace*{2cm}
X = PUQ + (\, I_m - P\,)V(\, I_m - Q\,), \hfill (4.5)
\cr
 \hspace*{2cm}
X = PW + WQ - 2PWQ, \hfill (4.6)
\cr}
$$  
where $ U, \ V, \ W  \in {\cal C}^{m \times n}$ are arbitrary.}

\medskip

\noindent {\bf Proof.}\, According to Corollary 4.2(c), 
the matrix equation $ PX = XQ $ is equivalent to the pair of matrix equations
$$\displaylines{ 
\hspace*{2cm}  
 PX(\, I  - Q \,) = 0   \  \ \ {\rm and}  \ \ \ ( \, I - P \,)XQ =0. \hfill (4.7)
\cr}
$$
Solving the pair of equations, we can find that both  (4.5) and  (4.6) are the general 
solutions of $ PX = XQ$. The process is somewhat tedious. Instead, we give here a  
direct verification.  Putting (4.5) in $ PX $ and $ XQ$,  we  get 
$$ \displaylines{ 
\hspace*{2cm} 
PX = PUQ \ \ \  \ {\rm and}  \ \ \ \ XQ = PUQ.  \hfill
\cr}
$$ 
Thus  (4.5) is solution of $ PX = XQ$. On the other hand, suppose that $ X_0$ is a solution of 
$ PX = XQ$ and  let $ U = V = X_0$ in (4.5). Then (4.5) becomes 
$$\displaylines{ 
\hspace*{2cm}  
X = PX_0Q + (\, I_m - P\,)X_0(\, I_m - Q\,) = PX_0Q + X_0 - PX_0 - X_0Q + PX_0Q = X_0,  \hfill 
\cr}
$$   
which implies that any solution of $ PX = XQ$ can be expressed by  (4.5). Hence  (4.5) is
 indeed the general solution of the equation $ PX = XQ $.  Similarly 
we can verify that  (4.6) is also a general solution to $ PX =XQ$. \qquad $ \Box $

\medskip

As one of the basic linear  matrix equation, $ AX = XB $ was examined 
(see, e.g., Hartwig \cite{Ha0},  Horn and Johnson \cite{HJ}, Parker \cite{Par},  Slavova 
et al \cite{SBZ}). In  general cases, solutions of $ AX = XB$ can only 
be determined through canonical forms of $ A$ and $B$. The result in Corollary 4.4 
manifests that for idempotent matrices $ A$ and $ B $, the 
general solution of  $ AX = XB$ can directly be written in 
$A$ and $B$. Obviously, the result in Corollary 4.4 is also 
valid for an operator equation of the form $ AX = XB$ when 
both $ A $ and $ B $ are idempotent operators.

\medskip

\noindent {\bf Theorem  4.5.}\, {\em Let $ A \in  {\cal C}^{m \times n}$ be
 given$,$  $ P \in {\cal C}^{m \times m}$ and $ Q \in {\cal C}^{n \times n} $
 be two idempotent matrices.  Then the sum $ PA + AQ $ satisfies the rank equalities
 $$
\displaylines{ 
\hspace*{2cm} 
 r( \, PA + AQ \, ) =
 r \left[ \begin{array}{cc}  PA & AQ \\ Q & 0 \end{array} \right] - r(Q)
 = r \left[ \begin{array}{cc} AQ & P \\ PA & 0 \end{array} \right]
 - r( P ),  \hfill (4.8)
 \cr
\hspace*{2cm}
 r( \, PA + AQ \, ) = r \left[ \begin{array}{c}  AQ - PAQ \\
 PA \end{array} \right]  = r[ \, PA - PAQ , \  AQ \, ].     \hfill (4.9)
\cr }
$$}
{\bf Proof.}\, Let $ M = \left[ \begin{array}{ccc} P  & 0  & PA
 \\ 0 & Q  & Q  \\ P & AQ & 0 \end{array} \right] $. Then it is easy to see
 by block elementary operations of matrices that 
 $$\displaylines{ 
\hspace*{2cm}  
 r(M ) = r \left[ \begin{array}{ccc} P & 0 & 0  \\ 0  & Q  & 0 \\ 0 & 0 &
 PA+AQ \end{array} \right] = r(P) + r(Q) + r( \, PA + AQ \,).\hfill
\cr}
$$
 On the other hand, note that $ P^2 = P$ and $ Q^2= Q$. We also obtain by block elementary operations  of matrices that 
\begin{eqnarray*}
 r(M) & = & r \left[ \begin{array}{ccc} P  & -PAQ  & PA  \\
 0  & 0  & Q  \\ P & AQ & 0  \end{array} \right] \\
 & = & r \left[ \begin{array}{ccc} 2P  & 0 & PA \\
    0 & 0 & Q  \\ P & AQ  & 0 \end{array} \right]
 = r \left[ \begin{array}{ccc} 2P  & 0 & 0 \\
    0 & 0 & Q  \\ 0 & AQ  & -\frac{1}{2}PA \end{array} \right]
 =  r(P) + r \left[ \begin{array}{cc}  PA & AQ \\ Q & 0 \end{array} \right],
 \end{eqnarray*}
and
 \begin{eqnarray*}
 r(M) & = & r \left[ \begin{array}{ccc} 0  & -PAQ  & PA  \\
 0  & Q  & Q  \\ P & AQ & 0  \end{array} \right] \\
 & = & r \left[ \begin{array}{ccc} 0  & 0 & PA \\
    0 & 2Q & Q  \\ P & AQ  & 0 \end{array} \right]
    = r \left[ \begin{array}{ccc} 0  & 0 & PA \\
    0 & 2Q & 0  \\ P & 0  & -\frac{1}{2}AQ \end{array} \right]
    =  r(Q) +  r \left[ \begin{array}{cc} AQ & P \\ PA & 0 \end{array}   \right].
 \end{eqnarray*}
 Combining the above three rank equalities for $ M $  yields  (4.8). Consequently applying  (1.2) and (1.3)  to the two block matrices in (4.8) yields  (4.9). \qquad  $\Box$ 

\medskip

\noindent {\bf Corollary 4.6.}\, {\em Let $ A \in  {\cal C}^{m \times n}$ be
 given$,$ $ P \in {\cal C}^{m \times m}$ and $ Q \in {\cal C}^{n \times n} $
 be two idempotent matrices.

{\rm (a)}\, If $ PAQ = 0,$  then  $ r( \, PA + AQ \, ) = r( PA ) + r(AQ),$
 or equivalently $ R(PA) \cap R(AQ) = \{ 0 \}$ and $ R[(PA)^*] \cap
  R[(AQ)^*] = \{ 0 \}.$

{\rm (b)}\,  $ PA + AQ  = 0  \Leftrightarrow 
 PA = 0 \ and  \ AQ = 0.$ 

{\rm (c)}\,  The general solution of the matrix equation $ PX + XQ  = 0$ is 
$ X = ( \, I - P \,)U( \, I - Q \,), $   where $ U \in {\cal C}^{m \times n}$ is 
arbitrary.} 

\medskip

\noindent {\bf Proof.}\, If $ PAQ = 0,$   then $ r( \, PA - AQ \, )
 = r( PA ) + r(AQ)$ by Theorem 4.1(b). Consequently $ r( \, PA + AQ \, )
 = r( PA ) + r(AQ)$ by Lemma 1.4(d). The result in Part (b) follows from (4.9).
According to (b), the equation  $PX + XQ  = 0$  is equivalent to the pair of matrix  equations  $PX = 0 $ and $  XQ  = 0$. 
According  to Rao and Mitra \cite{RM}, and Mitra \cite{Mi4},  the common general solution of the pair of equation 
is exactly $ X = ( \, I - P \,)U( \, I - Q \,),$ where $ U \in {\cal C}^{m \times n}$ is arbitrary. \qquad  $ \Box$

\medskip

\noindent {\bf Corollary 4.7.}\, {\em Let $ A, \ P, \ Q \in
 {\cal C}^{m \times m}$ be given with  $ P, \ Q $ being two idempotent matrices.
 Then the following five statements are equivalent$:$

 {\rm (a)}\, The sum $ PA +AQ $ is nonsingular.

 {\rm (b)} \  $ r\left[ \begin{array}{cc} PA & AQ \\ Q & 0 \end{array} \right]
  = m + r(Q). $

 {\rm (c)} \ $ r \left[ \begin{array}{cc} AQ & P \\ PA & 0 \end{array}
 \right] = m + r(P). $

 {\rm (d)} \ $ r[ \, PA, \  AQ \, ] = m \ and \ R \left[ \begin{array}{c}
 (PA)^* \\ (AQ)^* \end{array} \right] \cap R \left[ \begin{array}{c} Q^* \\
 0 \end{array} \right] = \{ 0 \}.$

 {\rm (e)} \ $ r\left[ \begin{array}{c} AQ \\ PA \end{array} \right]
 = m $ \ and  \ $ R \left[ \begin{array}{c} AQ \\ PA \end{array} \right]
 \cap R \left[ \begin{array}{c} P \\ 0 \end{array} \right]
 = \{ 0 \}. $  }  

\medskip

\noindent {\bf Proof.}\, Follows from (4.8).  \qquad  $ \Box$ 

\medskip

\noindent {\bf Theorem  4.8.}\, {\em Let $ A \in {\cal C}^{m \times n}$ be
 given$,$ $ P \in {\cal C}^{m \times m}$ and $ Q \in {\cal C}^{n \times n} $
 be two idempotent matrices.  Then the  rank of  $  A - PA - AQ $ satisfies the  equalities
 $$
 r( \,  A - PA - AQ \, ) = r \left[ \begin{array}{cc}  A & P \\ Q & 0
 \end{array} \right] + r(PAQ) - r(P) -  r(Q) = 
r( \, A - PA - AQ + PAQ \, ) + r(PAQ). \eqno (4.10) 
$$
In particular$,$ 

{\rm (a)} \ $PA + AQ = A  \Leftrightarrow  (\, I - P\,)A(\, I - Q \,) = 0$  and $  PAQ =0$.

{\rm (b)}\, The general solution of the matrix equation $ PX + XQ  = X$ is $ X = ( \, I - P \,)U Q + V( \, I - Q \,), $   where $ U, \ V \in {\cal C}^{m \times n}$ are arbitrary.}

\medskip

 \noindent {\bf Proof.}\, According to (4.1), we first find that 
 $$
r( \,  A - PA - AQ \, )= r[ \, (\, I - P\,) A - AQ \, ] = 
r \left[ \begin{array}{c} (\, I - P\,)A  \\ Q \end{array} \right] + r[\, AQ, \ I - P \,] - r(\, I - P\,) - r(Q).
$$
According to (1.2) and (1.3), we also get 
$$
r \left[ \begin{array}{c} (\, I - P\,)A  \\ Q \end{array} \right] = 
r \left[ \begin{array}{cc}  A & P \\ Q & 0 \end{array} \right] - r(P),  \ \ {\rm and}  \ \ 
r[\, AQ, \ I - P \,] = r(PAQ) + r(\, I - P\,).   
$$
Combining the above three yields  the first equality in  (4.10).   Consequently applying (1.4) to the block matrix in it yields the second  equality in  (4.10). 
Part (a) is a direct consequence of  (4.10), Part (a) follows from Corollary 4.4.  
\qquad  $ \Box$ 

\medskip

If replacing $ P$ and $ Q $ in Theorem 4.5 by $I_m - P$ and $ I_m - Q $,
 we can also obtain two rank equalities for  $ 2A - PA - AQ$. For simplicity
 we omit them here. 

\medskip

\noindent {\bf Theorem  4.9.}\, {\em Let $ A \in {\cal C}^{m \times n}$ be
 given$,$ $ P \in {\cal C}^{m \times m}$ and $Q \in {\cal C}^{n \times n} $
 be two idempotent matrices. Then the rank of  $A - PAQ $ satisfies the
 equality
 $$
 r( \,  A - PAQ \, ) =  r \left[ \begin{array}{ccc}  A & AQ  & P \\
 PA & 0 & 0 \\ Q & 0 & 0 \end{array} \right] - r(P) -  r(Q)= r \left[ \begin{array}{cc} (\, I - P \,) A(\, I - Q \,)  
& (\, I - P \,)AQ \\
 PA(\, I - Q \,) & 0  \end{array} \right].   \eqno (4.11)
 $$
In particular$,$ 

{\rm (a)} \ $ PAQ = A   \Leftrightarrow  (\, I - P \,) A(\, I - Q \,) = 0, \ \ 
 (\, I - P \,)AQ  = 0  \ \ and \ \  PA(\, I - Q \,)=0 \Leftrightarrow 
PA = A  \ and  \ AQ = A.$ 

{\rm (b)}\,  The general solution of the matrix equation $ PXQ = X$ is $ X = PUQ,$   where $ U \in {\cal C}^{m \times n}$ is arbitrary.
 }

\medskip

\noindent {\bf Proof.}\, Note that $ P^2 = P$ and $ Q^2 = Q$. It is easy to
 find that
 \begin{eqnarray*}
  r \left[ \begin{array}{ccc} A  & AQ  & P  \\ PA  & 0  & 0 \\ Q & 0 & 0
  \end{array} \right]
 & = & r \left[ \begin{array}{ccr} A & 0  & P \\  0 & -PAQ & -P  \\ Q & -Q  &
 0 \end{array} \right] \\
 & = & r \left[ \begin{array}{crr} A & 0  & P  \\
    -PAQ & 0 & -P  \\ 0 & -Q  & 0 \end{array} \right] \\
  & =  & r \left[ \begin{array}{crr} A -PAQ &  0 & 0 \\ 0  & 0  & -P \\
  0 & -Q  & 0 \end{array} \right] = r( \, A - PAQ \, ) + r(P) + r(Q),
 \end{eqnarray*}
as required for the first equality in  (4.11). Consequently applying  (1.4) to 
its left side yields the second one in  (4.11). Part (a) is a direct consequence of  (4.11), Part (b) can trivially be verified. \qquad  $ \Box$ 

\medskip

Applying  (4.1) to powers of difference of two idempotent matrices, we also
 find following several results. 

\medskip

 \noindent {\bf Theorem 4.10.}\, {\em Let $ P, \, Q \in
 {\cal C}^{m \times m}$ be two idempotent
 matrices. Then

{\rm (a)}\,  $ ( \,P - Q \, )^3$ satisfies the  two rank equalities
$$
\displaylines{ 
\hspace*{2cm} 
 r[ \, ( \, P - Q \, )^3 \,] = r \left[ \begin{array}{c}  P - PQP \\ Q
 \end{array} \right] + r[ \, Q - QPQ, \ P \,] - r(P) -  r(Q), \hfill (4.12)
\cr
\hspace*{2cm} 
 r[ \, ( \, P - Q \, )^3 \,] =  r[\, P - PQP - PQ + (PQ)^2 \,] +
 r[\, Q - QPQ - PQ + (PQ)^2 \,]. \hfill (4.13)
\cr}
$$ 
In particular$,$

{\rm (b)}\,  If $(PQ)^2 = PQ,$ then
$$ \displaylines{ 
\hspace*{2cm} 
r[ \, ( \, P - Q \, )^3 \,] = r( \, P - PQP \,) +  r( \,  QPQ - Q \,).
\hfill (4.14)
\cr}
$$

{\rm (c)} $ r[\, ( \, P - Q \, )^3 \, ] = r( \, P - Q \, ),$ i.e.$,$
 ${\rm Ind}( \, P - Q \, ) \leq 1,$  if and only if
 $$ 
\displaylines{ 
\hspace*{2cm}  
 r\left[ \begin{array}{c} P - PQP \\ Q \end{array} \right]
 = r \left[ \begin{array}{c} P \\ Q \end{array} \right], \ \ \
 and \ \ \ r[\, Q- QPQ, \ P \, ] =  r[\, Q, \ P \, ],  \hfill (4.15)
 \cr}
$$
or$,$ equivalently$,$
$$  \displaylines{ 
\hspace*{2cm} 
 R\left( \left[ \begin{array}{c} P - PQP \\ Q \end{array} \right]^* \right)
 = R\left( \left[ \begin{array}{c} P \\ Q \end{array} \right]^* \right), \ \ \
 and \ \ \ R[\, Q- QPQ, \ P \, ] =  R[\, Q, \ P \, ].  \hfill (4.16)
\cr} 
$$

{\rm (d)} \  $ ( \, P - Q\, )^3 = 0  \Leftrightarrow 
 r \left[ \begin{array}{c}  P - PQP \\ Q
 \end{array} \right] = r(Q)$ and $r[ \, Q - QPQ, \ P \,] = r(P) \Leftrightarrow 
 R(\, Q - QPQ \, ) \subseteq R(P)$ \  and  $R[( \, P - PQP \,)^*] \subseteq
 R(Q^*). $} 

\medskip

\noindent {\bf Proof.}\, Since $ P^2 = P$ and $ Q^2 = Q$, it is easy to
 verify that
 $$ \displaylines{ 
\hspace*{2cm} 
 ( \, P - Q \, )^3  = P(\, I_m - QP \,) - (\, I_m - QP \,)Q. \hfill (4.17)
 \cr}
$$ 
Letting $ A = I_m - QP$ and applying (4.1) and (4.2) to (4.17) immediately
 yields (4.12) and (4.13). The results in Parts (b)---(d) are  natural consequences
 of (4.13).  \qquad $\Box$  

\medskip

 \noindent {\bf Corollary 4.11.}\, {\em Let $ P, \, Q \in
 {\cal C}^{m \times m}$ be two idempotent matrices. Then
 $$\displaylines{ 
\hspace*{2cm} 
 r[ \, ( \, P - Q \, )^3 - ( \, P - Q \, ) \,]
 = r \left[ \begin{array}{c} PQP \\ Q \end{array} \right] +
 r[ \, QPQ, \ P \,] - r(P) -  r(Q). \hfill (4.18)
\cr} 
$$ 
In particular$,$

 {\rm (a)} \ $ P - Q$ is  tripotent $\Leftrightarrow$ $ R(QPQ) \subseteq R(P) $
 and $R[(PQP)^*] \subseteq R(Q^*).$

 {\rm (b)}\, If $ PQ = QP,$ then $ P - Q $ is tripotent.}

\medskip

 \noindent {\bf Proof.}\, Observe from  (4.17) that 
 $$\displaylines{ 
\hspace*{2cm}  
 ( \, P - Q \, )^3  - ( \, P - Q \, ) = -PQP + QPQ. \hfill
\cr} 
 $$ 
 Applying  (4.1) to it immediately yields  (4.18). The results in
  Parts (b) and (c) are natural consequences of  (4.18).  \qquad  $\Box$

\medskip

 \noindent {\bf Corollary 4.12.}\, {\em  A matrix  $A \in
 {\cal C}^{m \times m}$ is tripotent if and only if
  it can factor as  $ A = P - Q,$ where $ P $ and $ Q $ are two
  idempotent matrices with $PQ = QP$. }

\medskip

 \noindent {\bf Proof.}\, The ``if'' part comes from  Corollary 4.11(b).
 The `` only if'' part follows from a decomposition  of $ A $
 $$\displaylines{ 
\hspace*{2cm} 
 A = \frac{1}{2} ( \, A^2 + A  \, ) - \frac{1}{2} ( \, A^2 - A  \, ), \hfill
\cr}
 $$ 
where $ P = \frac{1}{2} ( \, A^2 + A  \, )$ and
  $ Q = \frac{1}{2}( \, A^2 - A  \,)$ are two idempotent matrices with
  $ PQ= QP$. \qquad $\Box$ 

\medskip

 The rank equality (4.12) can be extended to
the  matrix $ ( \, P - Q\, )^5,$  where both $ P $ and $ Q $ are idempotent. In fact,
 it is easy to verify 
 $$ \displaylines{ 
\hspace*{2cm} 
 ( \, P - Q \, )^5  = P(\, I_m - QP \,)^2 - (\, I_m - QP \,)^2Q. \hfill
 \cr}
$$ 
 Hence by (4.1) it follows that 
 $$\displaylines{ 
\hspace*{2cm} 
 r[ \, ( \, P - Q \, )^5 \,] = r \left[ \begin{array}{c}
 P(\, I_m - QP \,)^2 \\ Q \end{array} \right] +
 r[ \, (\, I_m - QP \,)^2Q, \ P \,] - r(P) -  r(Q).\hfill
\cr} 
$$ 
 Moreover, the above work can also be extended to
 $( \, P - Q \, )^{2k+1}( \, k = 3, \ 4, \ \cdots ),$ where both $ P $ and
 $ Q $ are idempotent. 

\medskip

 Applying (4.1) to $ PQ - QP$, where both $ P $ and $ Q $ are idempotent,
  we also obtain the following.

\medskip

 \noindent {\bf Corollary 4.13.}\, {\em  Let $ P, \, Q
 \in {\cal C}^{m \times m}$ be  two idempotent matrices. Then
 $$
\displaylines{ 
\hspace*{2cm} 
  r( \, PQ - QP   \,) = r \left[ \begin{array}{c}  PQ \\ P \end{array} \right]
   + r[ \,QP, \ P \,] - 2r(P), \hfill (4.19)
 \cr
\hspace*{2cm} 
 r( \, PQ - QP \,) = r \left[ \begin{array}{c}  QP \\ Q \end{array} \right]
 + r[ \,PQ, \ Q \,] - 2r(Q), \hfill (4.20)
 \cr
\hspace*{2cm}
  r( \, PQ - QP \,) = r( \,PQ - PQP \,) +  r( \, PQP -QP \,) , \hfill (4.21)
  \cr
\hspace*{2cm}
  r( \, PQ - QP \,) = r( \,PQ - QPQ \,) + r( \, QPQ -QP \,). \hfill (4.22)
 \cr}
$$
In particular$,$ if both $ P $ and $ Q $ are Hermitian idempotent$,$ then  
$$
\displaylines{ 
\hspace*{2cm}
  r( \, PQ - QP \,) = 2r( \,PQ - PQP \,) =  2r( \, PQ -QPQ \,). \hfill (4.23)
\cr} 
$$ } 
\hspace*{0.4cm} The rank equality (4.23) was proved by B\'{e}rub\'{e}, Hartwig and Styan \cite{BHS}. 

\medskip

\noindent {\bf Corollary 4.14.}\, {\em  Let $P, \, Q \in
 {\cal C}^{m \times m}$ be two idempotent matrices. Then
 $$ 
\displaylines{ 
\hspace*{2cm}
  r[ \, ( \,P - PQ \,)  + \lambda ( \, PQ -Q \,) ] = r(\, P - Q \, )
   \hfill (4.24)
\cr} 
$$ 
 holds for all $ \lambda \in {\cal C}$ with $ \lambda \neq 0$. In particular$,$
 $$\displaylines{ 
\hspace*{2cm}
  r( \,P + Q - 2PQ \,) = r( \,P + Q - 2QP \,) = r(\, P - Q \,). \hfill (4.25)
 \cr}
$$ } 
{\bf Proof.}\, Observe that 
 $$
\displaylines{ 
\hspace*{2cm} 
 ( \,P - PQ \,)  + \lambda (\, PQ - Q \,) = P (\, P + \lambda Q \, ) -
 (\, P + \lambda Q \, )Q. \hfill
 \cr}
$$
Thus it follows by (4.1) that 
 \begin{eqnarray*}
  r[ \, ( \,P - PQ \,)  + \lambda ( \, PQ -Q \,) ]
  &= & r \left[ \begin{array}{c}  P (\, P + \lambda Q \, ) \\
  Q \end{array} \right] + r[ \,(\,  P + \lambda Q \, )Q, \ P \,]
  - r(P) -r(Q) \\
 &= & r \left[ \begin{array}{c}  P \\ Q \end{array} \right]
 + r[\, \lambda Q, \ P \,] - r(P) -r(Q) \\
 &= & r \left[ \begin{array}{c}  P \\ Q \end{array} \right]
 + r[\,P , \ Q \,] - r(P) -r(Q).
 \end{eqnarray*}
 Contrasting it with  (3.1) yields  (4.24). Setting $ \lambda = -1 $ we have (4.25). 
\qquad  $\Box$

\medskip

 Replacing $ P $ by $ I_m - P$ in  (4.24), we also obtain the following. 

\medskip

 \noindent {\bf Corollary 4.15.}\, {\em  Let $ P, \, Q
 \in {\cal C}^{m \times m}$ be two idempotent matrices. Then
 $$ \displaylines{ 
\hspace*{2cm}  
 r( \, I_m - P - Q + \lambda PQ \, ) = r( \, I_m - P - Q \, ) \hfill
 \cr}
$$
holds for all $ \lambda \in {\cal C}$ with $ \lambda \neq 1$. } 

\medskip

In the remainder of this chapter, we apply the results in Chapter 3 to establish various rank equalities related to involutory matrices. A matrix $ A $ is said to be involutory if its square is identity, i.e., $ A^2 = I $. As two special types of matrices, involutory matrices and idempotent matrices are closely linked. As a matter of fact, for any involutory matrix $ A $, the two corresponding matrices $(\,  I + A \,)/ 2$ and $(\,  I - A \,)/ 2$ are idempotent. Conversely, for any idempotent matrix $ A$, the two corresponding matrices $ \pm ( \,I - 2A \,)$ are involutory. Based on the basic fact, all the results in Chapter 3 and this  chapter on idempotent matrices can dually be extended to involutory matrices. We next list  some of them.   
   
\medskip

\noindent {\bf Theorem 4.16.}\,  {\em Let $ A, \, B \in {\cal C}^{m \times m}$
be two involutory matrices. Then the ranks of  $ A + B $  and  $  A - B  $ satisfy the equalities  
 $$
\displaylines{ 
\hspace*{2cm}
 r( \, A + B  \, ) = r \left[ \begin{array}{c}   I + A \\  I - B
  \end{array} \right] +
r[ \,  I + A,  \  I - B  \, ] - r(\, I + A \, ) -  r(\, I - B \,),
\hfill (4.26)
\cr
\hspace*{2cm}
r( \,  A + B  \, ) =  r[\, ( \, I + A \,)( \, I + B \,)\,] +
r[\, ( \, I - A \,)( \, I - B \,)\,], \hfill (4.27)
 \cr
\hspace*{2cm}
 r( \, A - B  \, ) = r \left[ \begin{array}{c}   I + A \\  I + B
 \end{array} \right] +
r[ \,  I + A,  \  I + B  \, ] - r(\, I + A \, ) -  r(\, I + B \,),
\hfill (4.28)
\cr
\hspace*{2cm}
 r( \,  A - B  \, ) =  r[\, ( \, I + A \,)( \, I - B \,)\,] +
  r[\, ( \, I - A \,)( \, I + B \,)\,]. \hfill (4.29)
 \cr}
$$}
{\bf Proof.}\, Notice that  both $ P = ( \, I + A \,)/2 $
and $Q = ( \, I - B \,)/2$ are idempotent when $A $ and $ B $ are involutory.  In that 
case, 
$$\displaylines{ 
\hspace*{2cm} 
 r( \,  P  - Q  \, ) =  r\left[  \frac{1}{2}( \, I + A \,) -
 \frac{1}{2}( \, I - B \,) \right] =  r( \,  A + B  \, ),\hfill
\cr
\hspace*{0cm}
and \hfill
\cr
\hspace*{2cm} 
r \left[ \begin{array}{c} P \\  Q  \end{array} \right] +
r[ \, P,  \  Q  \,] - r(P) - r(Q) =  r \left[ \begin{array}{c}   I + A \\  I - B  \end{array} \right] + 
r[ \,  I + A,  \  I - B  \, ] - r(\, I + A \, ) -  r(\, I - B \,).\hfill
\cr}
$$
Putting them in (3.1) produces  (4.26). Furthermore we have 
 $$\displaylines{ 
\hspace*{2cm} 
r( \,  P  - PQ  \, ) = r\left[\, ( \, I + A \,) \left( \, I -  \frac{1}{2}( \,
I - B \,)\right) \, \right]= r[\, ( \, I + A \,)( \, I + B \,)\,], \hfill
\cr
\hspace*{2cm} 
r( \, PQ - Q \, ) = r \left[  \left( \frac{1}{2}( \, I + A \,) - I  \right)
( \, I - B \,) \right] = r[\, ( \, I - A \,)( \, I - B \,) \,]. \hfill
\cr}
$$ 
Putting them in (3.2) yields  (4.27). moreover, if $ B $ is involutory,
then $-B $ is also involutory. Thus replacing
 $ B $ by $-B$ in  (4.26) and (4.27) yields  (4.28) and (4.29).  \qquad $ \Box$

\medskip

\noindent {\bf Corollary 4.17.}\,  {\em Let $ A, \, B \in {\cal C}^{m
\times m}$ be two involutory matrices.

{\rm (a)}\,  If $ (\, I + A \, )(\, I - B \,) = 0$ or  $ (\, I - B \, )(\,
I  + A \,) = 0, $ then
 $$ \displaylines{ 
\hspace*{2cm} 
r( \, A + B  \, ) =  r(\, I + A \, ) + r(\, I - B \,). \hfill (4.30) 
\cr}
$$
 
{\rm (b)}\,  If $ (\, I + A \, )(\, I + B \,) = 0$ or  $ (\, I + B \, )(\,
 I  + A \,) = 0, $ then
 $$\displaylines{ 
\hspace*{2cm} 
r( \, A - B  \, ) =  r(\, I + A \, ) + r(\, I + B \,). \hfill (4.31) 
\cr}
$$}
{\bf Proof.}\,  The condition $ (\, I + A \, )(\, I - B \,) = 0$
is equivalent to $ I + A = B + BA$ and $ I - B = AB - A.$ In that case,   
$(\, I + A \, )(\, I + B \,)  = I  + A + B + AB = 2( \, I + A \,).$ 
and $(\, I - A \, )(\, I - B \,)  = I -  B  - A + AB = 2( \, I  - B  \,).$
Thus (4.27) reduces to (4.30). Similarly we show that under $ (\, I - B \, )
(\,  I  + A \,) = 0$, the rank equality (4.30)
also holds. The result in Part (b) is obtained by replacing $ B $ in Part (a)
 by $-B$. \qquad $ \Box$

\medskip 

\noindent {\bf Corollary 4.18.}\,  {\em Let $ A, \, B \in {\cal C}^{m \times
 m}$ be two involutory matrices. Then

{\rm (a)}\, The sum $ A + B $ is nonsingular if and only if   
 $$
R( \, I + A \, ) \cap R(\, I -  B \, ) = \{ 0 \}, \ \  
R( \, I + A^* \, ) \cap R(\, I -  B^* \, ) = \{ 0 \},  
\ \ and  \ \   r(\, I + A \, ) + r(\, I - B \,) = m.
$$

{\rm (b)}\, The difference $ A - B $ is nonsingular if and only if   
 $$
R( \, I + A \, ) \cap R(\, I + B \, ) = \{ 0 \}, \ \ R( \, I + A^* \, ) \cap 
R(\, I +  B^* \, ) = \{ 0 \}, \ \ and  \ \   r(\, I + A \, ) + r(\, I + B \,) = m.
$$}
{\bf Proof.}\, Follows immediately from  (4.26) and (4.27).
\qquad $ \Box $

\medskip

\noindent {\bf Theorem 4.19.}\,  {\em Let $ A, \, B \in {\cal C}^{m \times
m}$ be two involutory matrices. Then
 $ A + B $  and  $  A - B  $ satisfy the rank equalities  
 $$
\displaylines{ 
\hspace*{0.5cm}
r( \,  A + B  \, ) =  r[\, ( \, I + A \,)( \, I + B \,) \,] + r[\, ( \, I +
B \,)( \, I + A \,)\,]
- r(\, I + A \, ) -  r(\, I + B \,) + m, \hfill (4.32)
 \cr
\hspace*{0.5cm}
r( \,  A - B  \, ) =  r[\, ( \, I + A \,)( \, I - B \,) \,] + r[\, ( \, I -
B \,)( \, I + A \,)\,]
- r(\, I + A \, ) -  r(\, I - B \,) + m. \hfill (4.33)
 \cr }
$$}
{\bf Proof.}\, Putting  $ P = ( \, I + A \,)/2 $ and $Q = ( \,
I + B \,)/2$ in (3.8) and simplifying
 yields (4.32). Replacing  $ B $ by $-B$ in (4.32) yields  (4.33).  \qquad $
 \Box$

\medskip

The combination of  (4.27) with (4.32) produces the following rank equality 
$$  \displaylines{ 
\hspace*{2cm} 
r[\, ( \, I + B \,)( \, I + A \,) \,] = r(\, I + B \, ) +  r(\, I + A \,) -
 m +
r[\, ( \, I - A \,)( \, I - B \,) \,]. \hfill (4.34)
\cr}
$$
{\bf Theorem 4.20.}\,  {\em Let $ A, \ B \in {\cal C}^{m \times m}$
be  two involutory matrices. Then
$$ \displaylines{ 
\hspace*{2cm} 
r( \,  AB - BA  \, ) =  r( \, A + B \,) + r( \, A -  B \,) - m.  \hfill (4.35)
\cr
\hspace*{0cm}
In particular, \hfill
\cr
\hspace*{2cm} 
AB = BA \Leftrightarrow r( \, A + B \,) + r( \, A -  B \,) = m.  \hfill (4.36)
\cr}
$$ }
{\bf Proof.}\, Putting  $ P = ( \, I + A \,)/2 $ and $Q = ( \, I
- B \,)/2$ in  (3.24) and simplifying
 yields (4.35). \qquad $ \Box$

\medskip

 Putting  the formulas (4.26)---(4.29), (4.32) and (4.33) in  (4.35)
 may yield some other rank equalities for $ AB - BA$. We leave them to the reader. 

\medskip

\noindent {\bf Theorem 4.21.}\,  {\em Let $ A, \, B \in {\cal C}^{m \times m}$
 be  two involutory matrices. Then

 {\rm (a)} \ $ r\left[ \left( \frac{A + B}{2} \right)^2 - \frac{A + B}{2} \right] =  r( \,
 I -  A - B \,) + r( \, A + B \,) - m.$

{\rm (b)} \  $ r\left[ \left( \frac{A - B}{2} \right)^2 - \frac{A - B}{2} \right] =  r( \,
 I -  A + B \,) + r( \, A - B \,) - m. $ \\
In particular$,$ 

{\rm (c)} \  $ \frac{1}{2}(\, A + B \,) \  is \ idempotent \Leftrightarrow r( \, I -  A - B \,)
 + r( \, A + B \,) = m \Leftrightarrow r( \, A + B \,) =  r( \, I + A \,) -  r( \, I - B \,).$

{\rm (d)} \ $ \frac{1}{2}(\, A - B \,) \  is \ idempotent \Leftrightarrow r( \, I -  A + B \,)
 + r( \, A - B \,) = m \Leftrightarrow r( \, A - B \,) =  r( \, I + A \,) -  r( \, I + B \,).
$ }
 
\medskip

\noindent {\bf Proof.}\, Putting  $ P = ( \, I + A \,)/2 $ and $Q = ( \, I -
B \,)/2$ in  (3.32) and simplifying  yields Part (a).  Replacing $B$ by  $-B$ we get 
 Part (b). Part (c) and (d) follow from Parts (a) and (b), and Corollary 3.21(c). 
 \qquad $ \Box$

\medskip

\noindent {\bf Theorem 4.22.}\,  {\em Let $ A, \, B \in {\cal C}^{m \times m}$
be two involutory matrices. Then
$$\displaylines{ 
\hspace*{2cm} 
 r( \, 3I -  A - B  - AB \,) = r( \, 2I -  A - B \,). 
\hfill (4.39) 
\cr}
$$ }
{\bf Proof.}\, Putting  $ P = ( \, I + A \,)/2 $ and $Q = ( \, I +
B \,)/2$ in  (3.34) and simplifying
 yields  (4.39). \qquad $ \Box$

\medskip

\noindent {\bf Theorem 4.23.}\,  {\em Let $ A \in {\cal C}^{m \times m}$ be
an involutory matrix. Then

{\rm (a)} \  $ r( \, A - A^* \,) =  2r[ \, I +  A, \  I + A^*  \,] - 2r( \,
I + A \,)
 =  r[ \, I -  A, \  I - A^*  \,]- 2r( \, I - A \,).$
 
{\rm (b)} \  $ r( \, A + A^* \,) = m.$
 
{\rm (c)} \ $ r( \, AA^* - A^*A \,) =  r( \, A - A^* \,).$ } 

\medskip

\noindent {\bf Proof.}\, Putting  $ P = ( \, I \pm A \,)/2 $ and $Q = ( \, I
\pm A^* \,)/2$ in Corollary 3.26 and simplifying
 yields the desired results. \qquad $ \Box$

\medskip

\noindent {\bf Theorem 4.24.}\,  {\em Let $ A \in {\cal C}^{m \times m}$
and $ B \in {\cal C}^{n \times n}$ be two involutory matrices$,$
and $ X \in {\cal C}^{m \times m}$.   Then  $ AX - XB $ satisfies the rank
equalities
 $$
\displaylines{ 
\hspace*{0.5cm}
 r( \, AX - XB  \, ) = r \left[ \begin{array}{c}   (\, I_m + A \,)X \\  I_n
 + B \end{array} \right] +
r[ \,  X(\, I_n + B \,),  \  I_m + A  \, ] - r(\, I_m + A \, ) -
 r(\, I_n + B \,),
\hfill (4.40)
\cr
\hspace*{0.5cm}
 r( \,  AX - XB  \, ) =  r[\, ( \, I_m + A \,)X( \, I_n - B \,)\,] +
  r[\, ( \, I_m - A \,)X( \, I_n + B \,)\,]. \hfill (4.41)
 \cr}
 $$
In particular$,$
$$\displaylines{ 
\hspace*{0.5cm} 
AX = XB \Leftrightarrow ( \, I_m + A \,)X( \, I_n - B \,) = 0  \ \ and
 \ \ ( \, I_m - A \,)X( \, I_n + B \,)  = 0. \hfill (4.42)
\cr}
$$
}
{\bf Proof.}\, Putting  $ P = ( \, I_m + A \,)/2 $
and $Q = ( \, I_n + B \,)/2$ in (4.1) and (4.2) yields  (4.40) and (4.41).
The equivalence  in (4.42) follows from (4.41). \qquad $\Box$

\medskip

\noindent {\bf Theorem 4.25.}\,  {\em Let $ A \in {\cal C}^{m \times m}$
and $ B \in {\cal C}^{n \times n}$ be two involutory matrices.
 Then the general solution of the matrix equation $ AX = XB $  is
$$
X = V + AVB, \eqno (4.43)
$$
where $ V \in {\cal C}^{m \times n}$ is arbitrary. }

\medskip

\noindent {\bf Proof.}\, We only give the verification. Obviously
the matrix $X $ in (4.43) satisfies $ AX = AV + VB$ and $ XB = VB + AV$.
Thus $ X $ is a solution of $AX = XB$. On the other hand, for any solution
$ X_0$ of $ AX = XB$, let $ V = X_0/2$ in (4.43). Then we get $ V = AX_0B
= X_0$, that is, $ X_0$ can be represented by  (4.43). Thus  (4.43) is the
general  solution of the matrix equation $ AX = XB.$  \qquad $ \Box $

\medskip

\noindent {\bf Theorem 4.26.}\,  {\em Let $ A \in {\cal C}^{m \times m}$
be an involutory matrix$,$  and $ X \in {\cal C}^{m \times m}$. Then

{\rm (a)} \ $ AX - XA $ satisfies the rank equalities
$$
\displaylines{ 
\hspace*{2cm}
 r( \, AX - XA  \, ) = r \left[ \begin{array}{c}   (\, I + A \,)X \\
 I + A \end{array} \right] +
r[ \,  X(\, I + A \,),  \  I + A  \, ] - 2r(\, I + A \, ),
\hfill
\cr
\hspace*{2cm}
 r( \,  AX - XA  \, ) =  r[\, ( \, I + A \,)X( \, I - A \,)\,] +
  r[\, ( \, I - A \,)X( \, I + A \,)\,], \hfill 
 \cr
 \hspace*{0cm}
In \ particular, \hfill
\cr
\hspace*{2cm}
AX = XA \Leftrightarrow ( \, I + A \,)X( \, I - A \,) = 0  \ \ and
 \ \ ( \, I - A \,)X( \, I + A \,)  = 0.\hfill  
\cr}
$$

{\rm (b)}\, The general solution of the matrix equation $ AX = XA $ is
$$\displaylines{ 
\hspace*{2cm} 
X = V + AVA,  \hfill
\cr}
$$
where $ V \in {\cal C}^{ m \times m}$ is arbitrary.} 

\markboth{YONGGE  TIAN }
{5. RANK EQUALITIES FOR OUTER INVERSES OF MATRICES}

\chapter{Rank equalities for outer
inverses of matrices }

\noindent An outer inverse of a matrix $ A $ is the solution to the matrix 
equation $ XAX = X $, and is often denoted by $ X = A^{(2)}$. The 
collection of all outer inverses of $ A $ is often denoted by $ A\{2\}$. 
Obviously, the Moore-Penrose inverse, the Drazin inverse, the group inverse, and the weighted Moore-Penrose inverse of a matrix are 
naturally outer inverses of the matrix. If outer inverse of a 
matrix is also an inner inverse the matrix,  it is called a 
reflexive inner inverse of the matrix, and is often denoted by 
$ A^-_r$. The collection of all reflexive inner inverses of 
a matrix $ A $  is denoted by $ A\{1,\, 2\}$.  As one of important kinds of generalized inverses of matrices, outer 
inverses of matrices  and their applications have well been examined in the literature (see, e.g., 
\cite{BeG, CM2, GH, HLG, NC, WY,WW}). In this chapter, we shall 
establish several basic rank equalities related to differences and sums of outer inverses of a matrix, and then consider their various consequences.  The results obtained in 
this chapter will also be applied in the subsequent chapters. 

\medskip

\noindent {\bf Theorem 5.1.}\, {\em Let $  A \in {\cal C}^{m \times n}$ be
given$,$ and $ X_1, \, X_2 \in A\{ 2 \}$.  Then the difference $ X_1 - X_2$
satisfies the following three rank equalities
$$
\displaylines{ 
\hspace*{2cm} 
r( \, X_1 - X_2 \, ) = r \left[ \begin{array}{c} X_1 \\ X_2 \end{array} 
\right] 
+ r[ \, X_1,  \ X_2 \, ] - r(X_1) -  r(X_2), \hfill ( 5.1) 
\cr
\hspace*{2cm} 
r( \, X_1 - X_2 \, ) =  r( \, X_1 - X_1AX_2 \, ) + r(\, X_1AX_2  - X_2 \, ),
 \hfill (5.2)
\cr
\hspace*{2cm} 
r( \, X_1 - X_2 \, ) =  r( \, X_1 - X_2AX_1 \, ) + r(\, X_2AX_1  - X_2 \, ).
   \hfill (5.3)
\cr}
$$ }
{\bf Proof.}\, Let $ M = \left[ \begin{array}{ccc} -X_1  & 0  & X_1  \\ 0 & X_2  & X_2 
 \\ X_1 & X_2 & 0  \end{array} \right] $. Then it is easy to see by  block elementary operations of matrices
that 
$$ \displaylines{
\hspace*{2cm} 
r(M ) = r \left[ \begin{array}{ccc} -X_1 & 0 & 0  \\ 0  & X_2  & 0  \\ 0 & 0 &  X_1 - X_2 \end{array}
 \right] = r(X_1) + r(X_2) + r( \, X_1 - X_2 \,).\hfill (5.4) 
\cr}
$$
On the other hand, note that $ X_1AX_1 = X_1$ and $ X_2AX_2 = X_2$. Thus 
$$
\left[ \begin{array}{ccc} I_n & 0  & X_1A  \\ 0  & I_n  & 0  \\ 0 & 0 & I_n  \end{array} \right]
 \left[ \begin{array}{ccc} -X_1  & 0  & X_1  \\ 0 & X_2  & X_2  \\ X_1 & X_2 & 0  \end{array}\right]
\left[ \begin{array}{ccc} I_m & 0  & 0  \\ 0  & I_m  & 0  \\ 0 & -AX_2 & I_m  \end{array} \right] 
= \left[ \begin{array}{ccc} 0  & 0  & X_1  \\ 0 & 0  & X_2  \\ X_1 & X_2 & 0
 \end{array} \right],
$$
which implies that 
$$\displaylines{ 
\hspace*{2cm} 
 r(M) = r \left[ \begin{array}{c} X_1 \\ X_2 \end{array} \right] + r[ \, X_1,  \ X_2 \, ].\hfill (5.5) 
\cr}
$$  
Combining  (5.4) and (5.5) yields  (5.1). Consequently applying  (1.2) and (1.3) 
 to the two block matrices  in (5.1) respectively and noticing that $ A \in X_1\{ 2\}$ and $ A \in X_2\{ 2\}$, we can
write (5.1) as   (5.2) and (5.3).  \qquad $ \Box$ 

\medskip
 
It is obvious that if $ A = I_m$ in Theorem 5.1, then  $ X_1, \, X_2 \in
I_m\{ 2\}$ are actually two idempotent matrices. In that case,  (5.1)---(5.3)
reduce to the results in Theorem 3.1.  

\medskip

\noindent {\bf Corollary  5.2.}\, {\em Let $ A \in {\cal C}^{m \times n}$
 be  given$,$ and $ X_1, \, X_2 \in A\{ 2 \}$.  Then 
  
{\rm (a)} \ $R( \, X_1 - X_1AX_2 \, ) \cap R( \, X_1AX_2 - X_2 \,) = \{ 0\}$
 and $R[( \, X_1 - X_1AX_2 \, )^*] \cap R[( \, X_1AX_2 - X_2 \,)^*]
 = \{ 0\}.$

{\rm (b)} \ $R( \, X_1 - X_2AX_1 \, ) \cap R( \, X_2AX_1 - X_2 \,) = \{ 0\}$
 and $R[( \, X_1 - X_2AX_1 \, )^*] \cap R[( \, X_2AX_1 - X_2 \,)^*]
 = \{ 0\}.$

{\rm (c)}\, If $ X_1AX_2 = 0 $ or  $ X_2AX_1 = 0,$  then $ r( \, X_1 - X_2 \, ) =  r(X_1) + r(X_2).$ }

\medskip

\noindent {\bf Proof.}\, The results in Parts (a) and (b) follow
 immediately from applying Lemma 1.4(d)  to  (5.2) and (5.3).
Parts (c) is a direct consequence of  (5.2) and (5.3).  \qquad $ \Box$ 

\medskip

\noindent {\bf Corollary 5.3.}\, {\em Let $  A \in {\cal C}^{m \times n}$ be
given$,$ and $ X_1, \, X_2 \in A\{ 2 \}$.  Then  the following five statements are 
equivalent$:$ 

{\rm (a)} \ $ r( \, X_1 - X_2 \, ) =  r(X_1) - r(X_2),$  i.e.$,$  $ X_2 \leq_{rs} X_1.$   

{\rm (b)} \ $ r\left[ \begin{array}{c} X_1 \\ X_2 \end{array} \right] 
 = r[ \, X_1,  \ X_2 \, ] = r(X_1).$

{\rm (c)} \ $R(X_2) \subseteq R(X_1)$ and $ R(X_2^*) \subseteq R(X_1^*).$ 

{\rm (d)} \ $ X_1AX_2 = X_2 $ and $ X_2AX_1 = X_2.$

{\rm (e)} \ $ X_1AX_2AX_1 = X_2.$ }

\medskip

\noindent {\bf Proof.}\,  The equivalence of Parts (a) and (b) follows directly  
from (5.1). The equivalence of Parts (b), (c) and (d) follows directly  
from Lemma 1.2(c) and (d). Combining  the two equalities in  Part (d) yields the 
equality in Part (e). Conversely, suppose that $ X_1AX_2AX_1 = X_2$ holds. Pre- and 
post-multiplying $ X_1A $ and $ AX_1$  to it yields  $ X_1AX_2AX_1 = X_1AX_2 
= X_2AX_1$. Combining it with $ X_1AX_2AX_1 = X_2$ yields the two rank equalities in 
Part (d).  \qquad $ \Box$        

\medskip

\noindent {\bf Corollary 5.4.}\, {\em Let $  A \in {\cal C}^{m \times m}$ be
given$,$ and $ X_1, \, X_2 \in A\{ 2 \}$.  Then  the following three statements are 
equivalent$:$ 

{\rm (a)}\, The  difference $ X_1 - X_2 $ is nonsingular.

{\rm (b)} \ $ r \left[ \begin{array}{c}  X_1 \\ X_2 \end{array} \right]
 = r[ \, X_1,  \ X_2 \, ] = r(X_1) + r(X_2) = m$. 

{\rm (c)}  \ $R(X_1) \oplus R(X_2) = R(X_1^*) \oplus R(X_2^*)  = {\cal C}^m.$  } 

\medskip

\noindent {\bf Proof.}\, A trivial consequence of  (5.1).  \qquad $ \Box$ 

\medskip

\noindent {\bf Corollary  5.5.}\, {\em Let $ A \in {\cal C}^{m \times n}$ be  
given$,$ and $ X \in A\{ 2 \}$.  Then 
$$\displaylines{ 
\hspace*{2cm} 
r( \, A  - AXA \, ) = r(A) - r(AXA), \ \ i.e., \ \ AXA \leq_{rs} A.  \hfill (5.6)
\cr
\hspace*{0cm}
In \ particular, \hfill
\cr
\hspace*{2cm} 
 AXA = A, \ i.e., \   X \in A \{1, \, 2 \} \Leftrightarrow  r( A) = r(X). 
\hfill (5.7)
\cr}
$$ }
{\bf Proof.}\, It is easy to verify that both $ A $ and 
$ AXA$ are outer inverses of $A^{\dagger}$. Thus by  (5.1) we obtain
$$ \displaylines{ 
\hspace*{2cm} 
r( \, A  - AXA \, ) = r \left[ \begin{array}{c} A \\ AXA \end{array} \right]
 + r[ \, A,  \ AXA \, ] - r(A) -  r(AXA) = r(A) - r(AXA), 
\cr}
$$
the desired in  (5.6). \qquad $\Box $

\medskip

\noindent {\bf Corollary  5.6.}\, {\em Let $ A \in {\cal C}^{m \times m}$ be  
given$,$ and $ X \in A\{ 2 \}$.  Then 
$$ \displaylines{ 
\hspace*{1cm} 
r( \, AX  - XA \, ) =  r \left[ \begin{array}{c} X \\ XA \end{array} \right]
 + r[ \, X, \  AX \,] - 2r(X) =  r( \, XA - XA^2X \,) +  r( \, XA^2X - AX \,) . \hfill (5.8)
\cr
\hspace*{0cm}
In particular, \hfill
\cr 
\hspace*{2cm}  
AX = XA \Longleftrightarrow R(AX ) = R(X)  \ \ and \ \ R[(AX )^*] = R(X^*). \hfill (5.9)
\cr}
$$} 
{\bf Proof.}\, It is easy to verify that both $ AX $ and $ XA $ are idempotent when $ X \in A\{ 2 \}$. Thus we find by  (3.1), (1.2) and (1.3) that
\begin{eqnarray*}
r( \, AX  - XA \, ) & = & r \left[ \begin{array}{c} AX \\ XA \end{array} \right]
 + r[ \, AX,  \ XA \, ] - r(AX) -  r(XA) \\
& = & r \left[ \begin{array}{c} X \\ XA \end{array} \right]
 + r[ \, X,  \ AX \, ] - r(AX) -  r(XA) \\
& = & r( \, XA - XA^2X \,) +  r( \, XA^2X - AX \,),
\end{eqnarray*}
as required for  (5.8). Eq.\,(5.9) is a direct consequence of (5.8). 
\qquad $\Box $  

\medskip

\noindent {\bf Corollary 5.7.}\, {\em Let $ A \in {\cal C}^{m \times n}$ be  
given$,$ and $ X_1, \, X_2 \in A\{ 2 \}$.  Then 
$$\displaylines{ 
\hspace{2cm}
r( \, AX_1A  - AX_2A \, ) =  r \left[ \begin{array}{c} X_1A \\ X_2A
\end{array} \right] + r[ \, AX_1, \ AX_2 \, ] - r(X_1) -  r(X_2) . \hfill (5.10)
\cr
\hspace*{0cm}
In \ particular, \hfill
\cr 
\hspace{2cm}
AX_1A  = AX_2A  \Leftrightarrow  X_1AX_2AX_1 = X_1 \  and \   
X_2AX_1AX_2 = X_2. \hfill (5.11)
\cr}
$$}
{\bf Proof.}\, Notice that Both $ AX_1A$ and $AX_2A$ are outer inverses of  $A^{\dagger}$ when  $ X_1, \ X_2 \in A\{ 2 \}$. Moreover, observe that 
$r(AX_1A) = r(AX_1)= r(X_1A) = r(X_1)$, and  $r(AX_2A) = r(AX_2) = r(X_2A) = r(X_2).$ Thus it follows from (5.1) that
\begin{eqnarray*}
r( \, AX_1A  - AX_2A \, ) & = & r \left[ \begin{array}{c} AX_1A \\ AX_2A
\end{array} \right] + r[ \, AX_1A,  \ AX_2A \, ] - r(AX_1A) -  r(AX_2A) \\
& = & r \left[ \begin{array}{c} X_1A \\ X_2A \end{array} \right]
 + r[ \, AX_1,  \ AX_2 \, ] - r(X_1) -  r(X_2), 
 \end{eqnarray*}
as required for (5.10). The verification of (5.11) is trivial, hence is omitted. \qquad $\Box $

\medskip

\noindent {\bf Corollary 5.8.}\, {\em Let $ A \in {\cal C}^{m \times n}$ be  
given$,$ and $ X_1, \, X_2 \in A\{ 2 \}$.  Then  the following five statements are 
equivalent$:$ 

{\rm (a)} \ $ r( \, AX_1A - AX_2A \, ) =  r(AX_1A) - r(AX_2A),$ i.e.$,$  
$ AX_2A \leq_{rs} AX_1A.$ 

{\rm (b)} \ $ \left[ \begin{array}{c} X_1A \\ X_2A \end{array} \right] 
 = r[ \, AX_1,  \ AX_2 \, ] = r(X_1).$

{\rm (c)} \ $R(AX_2) \subseteq R(AX_1)$ and $ R[(X_2A)^*] \subseteq R[(X_1A)^*].$ 

{\rm (d)} \ $ AX_1AX_2A = AX_2A$ and $ AX_2AX_1A = AX_2A.$

{\rm (e)} \ $ AX_1AX_2AX_1A = AX_2A.$ }

\medskip

\noindent {\bf Proof.}\,  Follows form Corollary 5.3 by noticing that Both $ AX_1A$ 
and $AX_2A$ are outer inverses of  $A^{\dagger}$ when  $ X_1, \ X_2 \in A\{ 2 \}$.
\qquad $\Box $ 

\medskip

\noindent {\bf Theorem  5.9.}\, {\em Let $ A \in {\cal C}^{m \times n}$
be given$,$ and $ X_1, \, X_2 \in A\{ 2 \}$.  Then the sum $ X_1 + X_2$
satisfies the  rank equalities
$$
\displaylines{ 
\hspace*{2cm} 
r( \, X_1  + X_2 \, ) = r \left[ \begin{array}{cc} X_1  & X_2 \\ X_2  & 0
 \end{array} \right] - r(X_2)
 = r \left[ \begin{array}{cc} X_2  & X_1 \\ X_1  & 0 \end{array} \right]
  - r(X_1),  \hfill (5.11)
\cr
\hspace*{2cm} 
r( \, X_1 + X_2 \, ) =  r[\, ( \, I_n -  X_2A \,)X_1( \, I_m - AX_2 \,)\ ] +
 r(X_2), \hfill (5.12)
\cr
\hspace*{2cm} 
r( \, X_1 + X_2 \, ) =  r[\, ( \, I_n -  X_1A \,)X_2( \, I_m - AX_1 \,)\ ]
+ r(X_1). \hfill (5.13)
\cr}
$$ }
{\bf Proof.}\, Let $ M = \left[ \begin{array}{ccc} X_1  & 0 & X_1
  \\ 0 & X_2  & X_2  \\ X_1 & X_2 & 0  \end{array} \right] $. Then it is easy
 to see by block elementary operations  that
$$\displaylines{
\hspace*{2cm}   
r(M ) = r \left[ \begin{array}{ccc} X_1 & 0 & 0  \\ 0  & X_2  & 0
\\ 0 & 0 & -( \, X_1 + X_2 \,) \end{array} \right] = r(X_1) + r(X_2) +
  r( \, X_1 + X_2 \,). \hfill (5.14)
\cr}
$$
On the other hand, note that $ X_1AX_1 = X_1$ and $ X_2AX_2 = X_2$. Thus  
$$
\left[ \begin{array}{ccc} I_n & 0  & X_1A  \\ 0  & I_n  & 0  \\ 0 & 0 & I_n
\end{array} \right]
 \left[ \begin{array}{ccc} X_1  & 0  & X_1  \\ 0 & X_2  & X_2  \\ X_1 & X_2 &
  0  \end{array}\right]
\left[ \begin{array}{ccc} I_m & 0  & 0  \\ 0  & I_m  & 0  \\ 0 & -AX_2 & I_m
  \end{array} \right]
= \left[ \begin{array}{ccc} 2X_1  & 0  & X_1  \\ 0 & 0  & X_2  \\ X_1 & X_2 &
0  \end{array} \right],
$$
which implies that 
$$
r(M) = r\left[ \begin{array}{ccc} 2X_1  & 0  & X_1  \\ 0 & 0  & X_2  \\ X_1 & X_2 &
0  \end{array} \right] = r\left[ \begin{array}{ccc} 2X_1  & 0  & 0  \\ 0 & 0
 & X_2  \\ 0 & X_2 & -\frac{1}{2}X_1  \end{array} \right] 
= r \left[ \begin{array}{cc} X_1  & X_2 \\ X_2 & 0  \end{array} \right] + r(X_1).
  \eqno (5.15)
$$  
Combining  (5.14) and (5.15) yields the first equality in  (5.11). By
symmetry, we have the second equality in  (5.15). Applying  (1.3) to the two
block matrices in (5.11), respectively, and noticing that $ A \in \{ X_1^-\}$ 
and $ A \in \{ X_2^- \}$, we then can write  (5.11) as  (5.12) and (5.13). 
\qquad $ \Box$ 

\medskip

\noindent {\bf Corollary 5.10.}\, {\em Let $ A \in {\cal C}^{m \times n}$ be  
given$,$ and $ X_1, \, X_2 \in A \{ 2 \}$.

{\rm (a)} \ If $ X_1AX_2 = X_2AX_1,$ then $ r( \, X_1 + X_2 \, ) = r \left[ \begin{array}{cc} X_1 \\ X_2 \end{array} \right] = r[\, X_1, \ X_2\, ].$ 

{\rm (b)} \ If $ X_1AX_2 = X_2AX_1 = 0,$  then $r( \, X_1 + X_2 \, ) = 
 r(X_1) + r(X_2).$ } 

\medskip

\noindent {\bf Proof.}\, Under $ X_1AX_2 = X_2AX_1$, we find from
 (5.12) and (5.13) that
$$\displaylines{
\hspace*{2cm}   
r( \, X_1 + X_2 \, ) =  r(\, X_1 - X_1AX_2 \, ) + r(X_2) 
= r(X_1) + r(\, X_1AX_2 - X_2 \,).  \hfill  
\cr}
$$ 
Note by (1.2) and (1.3) that
$$\displaylines{ 
\hspace*{2cm} 
r \left[ \begin{array}{cc} X_1 \\ X_2 \end{array} \right]
=  r(\, X_1 - X_1AX_2 \, ) + r(X_2),  \ \  {\rm and} \ \
r[\, X_1, \ X_2\, ] = r(X_1) + r(\, X_1AX_2 - X_2 \,). \hfill  
\cr}
$$
Thus we have the results in Part (a). Part (b) follows immediately from (5.12).
\qquad $ \Box$

\medskip
 
\noindent {\bf Corollary  5.11.}\, {\em Let $ A \in {\cal C}^{m \times m}$
 be  given, and $ X_1, \, X_2
\in A\{ 2 \}$. Then the following five statements are equivalent$:$

{\rm (a)}\,  The sum $ X_1+ X_2 $ is nonsingular.

{\rm (b)} \ $ r \left[ \begin{array}{c}  X_1 \\ X_2 \end{array} \right] = m \ \ and  
\ \ R \left[ \begin{array}{c} X_1 \\ X_2 \end{array} \right] \cap R \left[ \begin{array}{c} 
 X_2 \\ 0  \end{array} \right] = \{ 0 \} $.

{\rm (c)} \ $ r[\, X_1, \ X_2 \,] = m \ \ and  \ \ R \left[ \begin{array}{c}  X_1^* \\ X_2^* 
\end{array} \right] \cap R \left[ \begin{array}{c} X_2^* \\ 0  \end{array} \right] = \{ 0 \}  $.
 
{\rm (d)} \ $ r \left[ \begin{array}{c}  X_2 \\ X_1 \end{array} \right] = m \ \ and  \ \ 
R \left[ \begin{array}{c}  X_2 \\ X_1 \end{array} \right] 
\cap R \left[ \begin{array}{c} X_1 \\ 0  \end{array} \right] = \{ 0 \}$.

{\rm (e)} \  $ r[\, X_2, \ X_1 \,] = m \ \ and  \ \ R \left[ \begin{array}{c}
 X_2^* \\ X_1^* \end{array} \right] \cap R \left[ \begin{array}{c}
  X_1^* \\
 0  \end{array} \right] = \{ 0 \} $. }

\medskip
 
\noindent {\bf Proof.}\, Follows immediately from  (5.11).
\qquad $ \Box$ 

\medskip

\noindent {\bf Corollary 5.12.}\, {\em Let $ A \in {\cal C}^{m \times n}$ be
 given$,$ and $ X \in A\{ 2 \}$. Then
$$\displaylines{
\hspace*{2cm}  
r( \, A  + AXA \, ) = r(A). \hfill (5.16)
\cr}
$$
holds for all $X \in A\{ 2 \}$. }

\medskip

\noindent {\bf Proof.}\, Notice that Both $A$ and $AX_2A$ are outer inverses
 of  $A^{\dagger}$ when  $ X \in A\{ 2 \}$. Thus  (5.16) follows from  (5.11). \qquad $\Box $ 

\medskip

\noindent {\bf Theorem  5.13.}\, {\em Let $  A \in {\cal C}^{m \times n}$ be
given$,$ and $ X_1, \, X_2
\in A\{ 2 \}$.  Then the difference $ X_1 - X_2$ satisfies the rank equalities
$$ 
\displaylines{ 
\hspace*{1cm} 
r[\, ( \, X_1 - X_2 \, )A( \, X_1 - X_2 \, ) - ( \, X_1 - X_2 \, ) \,]
 = r(\, I_m - AX_1 + AX_2 \, ) + r( \, X_1 - X_2 \, ) -m , \hfill (5.17)
\cr
\hspace*{1cm} 
r[\, ( \, X_1 - X_2 \, )A( \, X_1 - X_2 \, ) - ( \, X_1 - X_2 \, ) \,] =
r( X_1AX_2AX_1 \, ) - r(X_1) + r( \, X_1 - X_2 \, ). \hfill (5.18)
\cr}
$$ }
{\bf Proof.}\, Letting $ X = X_1 - X_2 $ and applying (1.10)
yields
$$ 
r( \, XAX  - X \, )  =  r( \, I_m - AX \,) +  r(X) - m, 
$$
which is  (5.17). Note that $ AX_1$ and $ AX_2$ are idempotent. It turns
 out by  (3.19) that
$$                               
r(\, I_m - AX_1 + AX_2 \, ) = r(AX_1AX_2AX_1) -r(AX_1) + m = r(X_1AX_2AX_1)
-r(X_1) + m .
$$    
Putting it in  (5.17) yields  (5.18). \qquad $ \Box$

 \medskip                                                         

\noindent {\bf Corollary 5.14.}\, {\em Let $A \in {\cal C}^{m \times n}$ be given$,$ and $ X_1, \, X_2 
\in A\{ 2 \}$. Then the following five statements are equivalent$:$

{\rm (a)} \  $ X_1 - X_2 \in A\{ 2 \}$.  
 
{\rm (b)} \ $ r(\,I_m - AX_1 + AX_2 \,)  = m - r( \, X_1 - X_2 \, ).$

{\rm (c)} \ $ r( \, X_1 - X_2 \, ) = r(X_1) - r(X_2),$ i.e.$,$  $ X_2 \leq_{rs} X_1$.

{\rm (d)}  \ $ R(X_2) \subseteq R(X_1)$ and  $ R(X_2^*) \subseteq R(X_1^*)$.

{\rm (e)} \  $X_1AX_2AX_1 = X_2.$  } 

\medskip

\noindent {\bf Proof.}\, The equivalence of Parts (a) and (b) follows immediately from  (5.17). 
The  equivalence of Parts (c), (d) and (e) is from Corollary  5.3. We next 
show the equivalence of Parts (a) and (e). It is easy to verify that  
$$\displaylines{
\hspace*{2cm}
( \, X_1 - X_2 \, )A( \, X_1 - X_2 \, ) - ( \, X_1 - X_2 \, ) = - X_1AX_2 - X_2AX_1 + 2X_2, \hfill
\cr}
$$ 
Thus $ X_1 - X_2 \in A\{ 2 \}$ holds if and only if 
$$\displaylines{
\hspace*{2cm}
X_1AX_2  +  X_2AX_1  = 2X_2. \hfill (5.19)
\cr}
$$
Pre- and post-multiplying $ X_1A $ and $ AX_1$ to it, we get 
$$\displaylines{
\hspace*{2cm} 
X_1AX_2AX_1 = X_1AX_2 \ \ {\rm and} \ \   X_1AX_2AX_1 = X_2AX_1. \hfill (5.20)
\cr}
$$ 
Putting them in (5.19) yields Part (e). Conversely, if Part (e) holds,
then  (5.20) holds. Combining Part (e) with 
 (5.20) leads to   (5.19), which is equivalent to $ X_1 - X_2 \in
 A\{ 2 \}$. \qquad $ \Box$ 

\medskip

The problem considered in Corollary 5.14 could be regarded as an extension of 
the work in Corollary 3.21, which was examined by Getson and Hsuan cite{GH}.
In that monograph, they only gave a sufficient condition for $ X_1 - X_2 \in
A\{ 2 \}$ to hold when $ X_1, \ X_2 \in A\{ 2 \}$.  Our result in Corollary 5.14 
is a complete conclusion on this problem. 
 
\medskip

\noindent {\bf Theorem  5.15.}\, {\em Let $  A \in {\cal C}^{m \times n}$ be
given$,$ and $ X_1, \, X_2
\in A\{ 2 \}$.  Then the sum $ X_1 + X_2$ satisfies the two rank
equalities
$$\displaylines{
\hspace*{0cm}
r[\, ( \, X_1 + X_2 \, )A( \, X_1 + X_2 \, ) - ( \, X_1 + X_2 \, ) \,]
= r(\, I_m - AX_1 - AX_2 \, ) + r( \, X_1 + X_2 \, ) -m , \hfill (5.21)
\cr
\hspace*{0cm}
r[\, ( \, X_1 + X_2 \, )A( \, X_1 + X_2 \, ) - ( \, X_1  + X_2 \, ) \,]
= r( X_1AX_2 \, ) + r( X_2AX_1 \, ) + r( \, X_1 +  X_2 \, ) - r(X_1) -
r(X_2). \hfill (5.22)
\cr}
$$ } 
{\bf Proof.}\, Letting $ X = X_1+  X_2 $ and applying  (1.10) to
 $XAX  - X $ yields  (5.21). Note that $ AX_1$ and $ AX_2$ are idempotent.
  It turns out by (3.8) that
$$ \displaylines{
\hspace*{2cm}
r(\, I_m - AX_1 - AX_2 \, ) = r( X_1AX_2 \, ) + r( X_2AX_1 \, ) - r(X_1)
 - r(X_2) + m. \hfill
\cr}
$$    
Putting it in (5.21) yields (5.22).  \qquad $ \Box$ 

\medskip

\noindent {\bf Corollary 5.16.}\, {\em Let $A \in {\cal C}^{m \times n}$ be
given$,$ and $ X_1, \, X_2 \in A\{ 2 \}$. Then the following four statements
are equivalent$:$

{\rm (a)} \ $ X_1 + X_2 \in A\{ 2 \}$.  
 
{\rm (b)}  \ $  X_1AX_2  + X_2AX_1  = 0.$

{\rm (c)}  \ $ r(\,I_m - AX_1 - AX_2 \,)  = m - r( \, X_1 + X_2 \, ).$

{\rm (d)} \  $  X_1AX_2 = 0,$ and $  X_2AX_1  = 0.$ }

\medskip

\noindent {\bf Proof.}\, The equivalence of Parts (a) and (b) follows
immediately from expanding
$( \, X_1  + X_2 \, )A( \, X_1 + X_2 \, ) - ( \, X_1 + X_2 \, ).$
The  equivalence of (a) and (c) is from (5.21). We next show the equivalence of 
(b) and (d). Pre- and post-multiplying $ X_1A $ and $ AX_1$ to $ X_1AX_2  +  X_2AX_1  = 0,$ we get
$$
X_1AX_2  + X_1AX_2AX_1 = 0 , \ \ \ {\rm and}  \ \ \   X_1AX_2AX_1  + X_2AX_1 = 0, 
$$
which implies  that $ X_1AX_2 = X_2AX_1$. Putting them in Part (b) yields (d). 
Conversely, if Part (d) holds, then  Part (b) naturally holds. \qquad $ \Box$ 

\markboth{YONGGE  TIAN }
{6. RANK EQUALITIES FOR A MATRIX AND ITS MOORE-PENROSE INVERSE}

\chapter{Rank equalities for a matrix and its  Moore-Penrose inverse}

\noindent In this chapter, we shall establish a variety of rank equalities 
related to a matrix and its Moore-Penrose inverse, and then use them to characterize various specified matrices, such as, EP matrices, conjugate EP matrices, bi-EP matrices, star-dagger matrices, and so on. 

A matrix $ A $ is said to be EP (or Range-Hermitian) if $R( A ) = R( A^*)$. EP matrices have some nice properties, meanwhile they are quite inclusive. Hermitian matrices,  normal matrices, as well as nonsingular matrices are special cases of EP matrices. As a class of important matrices, EP matrices and their applications have well be examined in the 
literature.  One of the basic and nice properties related to an EP matrix $ A$ is $AA^{\dagger} = A^{\dagger}A$, see, e.g., 
Ben-Israel and Greville \cite{BeG}, Campbell and Meyer \cite{CM2}.  This equality motivates us to consider the rank of 
$AA^{\dagger} - A^{\dagger}A$, as well as its various extensions.          
  
\medskip

\noindent {\bf Theorem 6.1.}\, {\em Let $ A \in {\cal C}^{m \times m}$ be
given.
Then the rank of $ AA^{\dagger} - A^{\dagger}A$ satisfies the following rank
equalities
$$\displaylines{
\hspace*{2cm}
r(  \, AA^{\dagger} -  A^{\dagger}A \, ) = 2r[\, A, \ A^* \, ] - 2r( A )
=2r( \, A -  A^2A^{\dagger} \, ) = 2r( \, A -  A^{\dagger}A^2 \, ).
\hfill (6.1)
\cr}
$$
In particular$,$ 

{\rm (a)}  $ AA^{\dagger} =  A^{\dagger}A  \Leftrightarrow r[\, A, \ A^* \, ] = r( A )
\Leftrightarrow A = A^2A^{\dagger} \Leftrightarrow
 A = A^{\dagger}A^2  \Leftrightarrow R( A ) = R( A^* ), \
  i.e., \ A \ is \ EP.$

{\rm (b)}  $ AA^{\dagger} - A^{\dagger}A$ is nonsingular
 $ \Leftrightarrow r[\, A, \ A^*\,] = 2r( A ) = m
 \Leftrightarrow R( A ) \oplus R( A^* ) = {\cal C}^m. $ } 

\medskip

\noindent {\bf Proof.}\,  Note that $AA^{\dagger}$ and $A^{\dagger}A $ are
idempotent matrices. Then applying (3.1), we first obtain
$$\displaylines{
\hspace*{2cm}
r( \, AA^{\dagger} -  A^{\dagger}A \,) = r\left[ \begin{array}{c}
AA^{\dagger}
 \\ A^{\dagger}A \end{array} \right] + r[\, AA^{\dagger}, \ A^{\dagger}A \,]
   - r(AA^{\dagger}) - r(A^{\dagger}A). \hfill (6.2)
\cr}
$$
Observe that $ r(AA^{\dagger}) = r(A^{\dagger}A) =r(A),$ and 
$$\displaylines{
\hspace*{2cm}
r\left[ \begin{array}{c} AA^{\dagger} \\ A^{\dagger}A \end{array} \right] 
= r\left[ \begin{array}{l} A^{\dagger} \\ A \end{array} \right]
= r\left[ \begin{array}{l} A^* \\ A \end{array} \right], \qquad
r[\, AA^{\dagger}, \ A^{\dagger}A \,] = r[\, A, \ A^{\dagger} \,]
= r[\, A, \ A^*\,]. \hfill
\cr}
$$ 
Thus (6.2) reduces to the first rank equality in (6.1). Consequently
applying (1.2) to $[\, A, \ A^* \,]$  in  (6.1) yields the other two
rank equalities in  (6.1). The equivalence in Part (a) are well-known
results on a  EP matrix, which  now is a direct consequence of (6.1).
It remains to show Part (b). If $ r[\, AA^{\dagger} - A^{\dagger}A \, ]
= m,$ then $ r[\, A, \ A^*\, ] =  r[\, AA^{\dagger}, \ A^{\dagger}A \,]
= r[\, AA^{\dagger}- A^{\dagger}A, \ A^{\dagger}A \,] = m.$
 Putting it in (6.1), we obtain $2r( A ) = m$. Conversely,
 if $ r[\, A, \ A^*\, ] = 2r( A ) = m$, then we immediately have
 $ r(  \, AA^{\dagger} -  A^{\dagger}A \, ) = m$ by  (6.1).
 Hence the first equivalence in Part (b)  is true. The second equivalence
 is obvious.  \qquad $\Box$ 

\medskip

Another  group of rank equalities related to EP matrix is given below, which is motivated by 
 a work of Campbell and Meyer \cite{CM1}. 

\medskip

\noindent {\bf Theorem 6.2.}\, {\em Let $ A \in {\cal C}^{m \times m}$ and $ 0 \neq k \in {\cal C}$  be
given. Then 

{\rm (a)}  \ $r[\, AA^{\dagger}(\, A + kA^{\dagger}\,) - (\, A + kA^{\dagger}\,) AA^{\dagger} \, ] 
= 2r[\, A, \ A^* \, ] - 2r( A ).$

{\rm (b)} \  $r[  \, A^{\dagger}A(\, A + kA^{\dagger}\,) - (\, A + kA^{\dagger}\,)A^{\dagger}A \, ] 
= 2r[\, A, \ A^* \, ] - 2r( A ).$

{\rm (c)} \  $r[  \, AA^{\dagger}(\, A + kA^*\,) - (\, A + kA^*\,) AA^{\dagger} \, ] 
= 2r[\, A, \ A^* \, ] - 2r( A ).$

{\rm (d)} \  $r[  \, A^{\dagger}A(\, A + kA^* \,) - (\, A + kA^* \,)A^{\dagger}A \, ] 
= 2r[\, A, \ A^* \, ] - 2r( A ).$  
 
{\rm (e)}\, The following five statements are equivalent$:$ 

\hspace*{0.5cm} {\rm (1)} \ $A$ is EP.  

\hspace*{0.5cm} {\rm (2)} \  $AA^{\dagger}(\, A + kA^{\dagger}\,) = (\, A + kA^{\dagger}\,) AA^{\dagger}.$

\hspace*{0.5cm} {\rm (3)} \  $A^{\dagger}A(\, A + kA^{\dagger}\,) = (\, A + kA^{\dagger}\,)A^{\dagger}A.$ 

\hspace*{0.5cm} {\rm (4)} \  $AA^{\dagger}(\, A + kA^* \,) = (\, A + kA^*\,) AA^{ \dagger}.$

\hspace*{0.5cm} {\rm (5)} \  $A^{\dagger}A(\, A + kA^*\,) = (\, A + kA^* \,)A^{\dagger}A.$ }

\medskip 

\noindent {\bf Proof.}\, We only show Part (a). Notice that both $AA^{\dagger}$ and $A^{\dagger}A $ are idempotent matrices. Thus by (4.1) we get
$$
\displaylines{
\hspace*{2cm}
r[\, AA^{\dagger}(\, A + kA^{\dagger}\,) - (\, A + kA^{\dagger}\,) AA^{\dagger} \, ] \hfill
\cr
\hspace*{2cm} 
 =  r \left[ \begin{array}{c} AA^{\dagger}(\, A + kA^{\dagger}\,)  \\  AA^{\dagger} \end{array} \right]  + 
r[\, (\, A + kA^{\dagger}\,) AA^{\dagger}, \  AA^{\dagger} \, ] - r(AA^{\dagger}) - r(A^{\dagger}A) \hfill
\cr
\hspace*{2cm}
 =  r \left[ \begin{array}{c}  A + kAA^{\dagger}A^{\dagger}  \\  A^* \end{array} \right]  + 
r[\, A^2A^{\dagger} + kA^{\dagger} , \  A \, ] - 2r(A)  \hfill
\cr 
\hspace*{2cm} 
=  r \left[ \begin{array}{c}  A   \\  A^* \end{array} \right]  + 
r[\, A^* , \  A \, ] - 2r(A) = 2r[\, A, \ A^* \, ] - 2r( A ), \hfill
\cr}
$$
establishing Part (a). \qquad $ \Box$  

\medskip

When $ k = 1$, the corresponding result in Theorem 6.2(e) was established by Campbell and Mayer \cite{CM1}. 

\medskip 

\noindent {\bf Theorem 6.3.}\, {\em Let $ A \in {\cal C}^{m \times m}$ and $ 0 \neq k \in {\cal C}$  be
given. Then 

{\rm (a)} \ $r[\, AA^{\dagger}(\, AA^* + kA^*A\,) - (\, AA^* + kA^*A\,) 
A^{\dagger}A \, ] = 2r[\, A, \ A^* \, ] - 2r( A ).$

{\rm (b)} \ $AA^{\dagger}(\, AA^* + kA^*A\,) = (\, AA^* + kA^*A\,) 
A^{\dagger}A \Leftrightarrow$  $ A $ is EP. } 

\medskip 

\noindent {\bf Proof.}\, Follows from (4.1) by noting that both $AA^{\dagger}$ and $A^{\dagger}A $ 
are idempotent matrices. \qquad $\Box$

\medskip 

In an earlier paper by Meyer \cite{MJ2} and a recent paper by Hartwig and Katz \cite{HK}, they established 
 a necessary  and  sufficient condition for a block triangular matrix to be EP. 
Their work now can be extended to  the following general settings.    

\medskip

\noindent {\bf Corollary 6.4.}\, {\em Let $ A \in  {\cal C}^{m \times m},$ $ B \in  {\cal C}^{m \times k},$
 and $  D \in  {\cal C}^{k \times k}$ be given$,$ and let
 $ M = \left[ \begin{array}{cc} A & B \\ 0 & D  \end{array} \right]$.
 Then
$$\displaylines{
\hspace*{2cm}
r(\, MM^{\dagger} -  M^{\dagger}M \,) = 2r\left[ \begin{array}{cccc}
A & A^* & B & 0 \\ 0 & B^* & D & D^*  \end{array} \right] - 2r(M).
\hfill (6.3)
\cr}
$$
In particular$,$

{\rm (a)}\, If both $ A $ and $ D $ are EP$,$ then 
$$
\displaylines{
\hspace*{2cm}
r(\, MM^{\dagger} -  M^{\dagger}M \,) = 2r[\, A, \ B \, ] +
2r\left[ \begin{array}{c} B \\ D \end{array} \right] -
2r\left[ \begin{array}{cc} A & B \\ 0 & D  \end{array} \right].\hfill  
\cr}
$$ 

{\rm (b)}\, If $R(B) \subseteq R(A)$ and $R( B^*) \subseteq R(D^*),$  then 
$$
\displaylines{
\hspace*{2cm}
r(\, MM^{\dagger} -  M^{\dagger}M \,) = 2r[\, A, \ A^* \, ] +
2r[\, D, \ D^* \, ] - 2r(A) - 2r(D).\hfill  
\cr}
$$ 

{\rm (c) (Meyer \cite{MJ2}, Hartwig and Katz \cite{HK})} \ $ M $ is EP if and only if both $ A $ and $ D $ are EP$,$  and $R( B) \subseteq R(A)$ and $R( B^*) \subseteq R(D^*).$ In that case$,$ 
 $ MM^{\dagger} = \left[ \begin{array}{cc} AA^{\dagger} &0 \\ 0 & DD^{\dagger}  \end{array} \right]$.
}

\medskip

\noindent {\bf Proof.}\,  Follows immediately from Theorem 6.1 by putting 
$ M$ in it. \qquad $\Box$

\medskip

 \noindent {\bf Corollary 6.5.}\, {\em Let
$$
 M = \left[ \begin{array}{cccc} A_{11} & A_{12} & \cdots & A_{1k}  \\
  & A_{22} & \cdots & A_{2k} \\ & & \ddots & \vdots \\
  & & & A_{kk}   \end{array} \right] \in {\cal C}^{m \times n}, \qquad  A_{ij} \in 
{\cal C}^{m_i \times m_j} 
$$
be given. Then $ M $ is $ EP $ if and only if  $A_{11}, \ A_{22} \cdots, \ 
A_{nn}$ are $EP,$ and 
$$
R(A_{ij}) \subseteq R(A_{ii}), \ \ and  \ \  R(A_{ij}^*) \subseteq R(A_{jj}^*),
 \ \ \ \ i,  \ j = 1, \ \cdots, \ n.   
$$
In that case$,$ $MM^{\dagger} ={\rm diag}( \, A_{11}A_{11}^{\dagger}, \
 A_{22}A_{22}^{\dagger}, \  \cdots, \ A_{nn}A_{nn}^{\dagger} \,)$.}

\medskip

\noindent {\bf Proof.}\,  Follows from Theorem 6.1(a) by putting $ M$ in it. 
\qquad $\Box$

\medskip

We leave the verification of the following result to the reader. Let $ A, \ B  \in  
{\cal C}^{m \times m}$ be given, and let
 $ M = \left[ \begin{array}{cc} A & B \\ B & A  \end{array} \right]$.
Then
$$\displaylines{
\hspace*{1cm}
r(\, MM^{\dagger} -  M^{\dagger}M \,) \hfill
\cr
\hspace*{1cm}
 = 2r[ \, A+ B, \ (A + B)^* \,] + 
2r[\, A - B , \ (A - B)^* \,] - 2r( \,A + B \,) - 2r( \,A - B  \,) \hfill
\cr
\hspace*{1cm} 
= r[\,( \,A + B \,)( \,A+ B \,)^{\dagger} - ( \,A + B \,)^{\dagger}(\, A + B \,) \,] 
 + r[\,( \,A - B \,)( \,A  - B \,)^{\dagger} - ( \,A - B \,)^{\dagger}( \,A - B \,)  \,]. \hfill
\cr}
$$
In particular, $ M $ is EP  if and only if $ A \pm B $ are EP.

A parallel concept to EP matrices is so-called  conjugate EP matrices. A matrix $ A $ is said to be conjugate EP if $R( A ) = R( A^T)$. If matrices considered are real, then  EP matrices and conjugate EP matrices are identical.     
 Much similar to EP matrices, conjugate EP matrices also have some nice 
properties.  One of the  basic and nice properties related to a conjugate 
EP matrix  $ A$ is $AA^{\dagger} = \overline{A^{\dagger}A}$ (see the 
series work \cite{MI0}, \cite{MI1}, \cite{MI2}, \cite{MI3}, and \cite{MI4} 
by Meenakshi and Indira).  This equality motivates us to find the following results.        
  
\medskip

\noindent{\bf Theorem 6.6.}\, {\em Let $ A \in {\cal C}^{m \times m} $ be
given. Then
$$
r(\, AA^{\dagger} - \overline{ A^{\dagger}A} \,) = 2r[\, A, \ A^T \,]
- 2r( A ). \eqno (6.4)
$$ 
In particular$,$

{\rm (a)} \  $ AA^{\dagger} = \overline{A^{\dagger}A} 
\Leftrightarrow
 r[\, A, \ A^T \,] = r( A ) \Leftrightarrow
R( A ) = R( A^T),$ \ i.e.$,$ $ A $  is conjugate EP.   

{\rm (b)} \  $ AA^{\dagger} -  \overline{ A^{\dagger}A}$ is
nonsingular $ \Leftrightarrow  r[\, A, \ A^T \, ] = 2r( A ) = m
\Leftrightarrow R( A ) \oplus R( A^T) = {\cal C}^m. $ } 

\medskip

\noindent {\bf Proof.}\,  Since $ A^{\dagger}A $ and $\overline{A^{\dagger}A}$ are idempotent, 
applying (3.1) to $ AA^{\dagger} - \overline{ A^{\dagger}A}$, we  obtain
\begin{eqnarray*}
r( \, AA^{\dagger} - \overline{ A^{\dagger}A} \,) 
& = & r\left[ \begin{array}{c} AA^{\dagger} \\ \overline{ A^{\dagger}A}
\end{array} \right]+ r[\, AA^{\dagger}, \ \overline{ A^{\dagger}A} \,]
 - r(AA^{\dagger}) - r(\overline{ A^{\dagger}A} ) \\
& = & r\left[ \begin{array}{c} A^{\dagger} \\ \overline{A} \end{array} \right] +
 r[\, A, \ \overline{ A^{\dagger}} \,]  - 2r(A) \\
& = & r\left[ \begin{array}{c} A^{*} \\ \overline{A} \end{array} \right] + 
r[\, \overline{A}, \ A^* \,]  - 2r(A) = 2r[\, A, \ A^T \,] - 2r( A ), 
\end{eqnarray*}
which is exactly (6.4). The results in Parts (a) and (b) follow
immediately from (6.4).
\qquad \ $ \Box$ 

\medskip

\noindent {\bf Corollary 6.7.}\, {\em Let $ A \in  {\cal C}^{m \times m}$,
$ B \in  {\cal C}^{m \times k},$
 and $  D \in  {\cal C}^{k \times k}$ be given$,$ and denote
 $ M = \left[ \begin{array}{cc} A & B \\ 0 & D  \end{array} \right]$.
 Then
$$\displaylines{
\hspace*{2cm}
r(\, MM^{\dagger} -  \overline{ M^{\dagger}M } \,) = 2r\left[ \begin{array}{cccc}
A & A^T & B & 0 \\ 0 & B^T & D & D^T  \end{array} \right] - 2r(M). \hfill
\cr}
$$
In particular$,$ $ M $ is conjugate EP if and only if $ A $ and $ D $ are 
con-EP$,$ and $R( B) \subseteq R(A)$ and $R( B^T) \subseteq R(D^T).$ }

\medskip

\noindent {\bf Proof.}\, Follows from Theorem 6.6 by putting $ M $ in it. 
\qquad \ $ \Box$ 

\medskip

The work in Theorem 6.1 can be extended to matrix expressions that involve powers of a matrix. 

\medskip

\noindent {\bf Theorem 6.8.}\, {\em Let $ A \in {\cal C}^{m \times m} $ be
given and $k$ be an integer with $k \geq 2.$ Then
$$\displaylines{
\hspace*{2cm}
r( \, A^kA^{\dagger} -  A^{\dagger}A^k \, ) = r \left[ \begin{array}{c}
 A^k  \\ A^* \end{array} \right]  + r[\, A^k, \ A^*\, ] - 2r(A).  \hfill (6.5)
\cr}
$$
In particular$,$

{\rm (a)} \ $r( \, A^kA^{\dagger} -  A^{\dagger}A^k \, ) =
  r[\, A, \ A^*\, ] - 2r(A),$ if $ r(A) = r(A^2).$    

{\rm (b)} \   $ A^kA^{\dagger} = A^{\dagger}A^k \Leftrightarrow  
r \left[ \begin{array}{c} A^k  \\ A^* \end{array} \right] =  r[\, A^k, \ A^*\, ] = r(A) \Leftrightarrow
 R(A^k) \subseteq  R(A^*) \ and \  R[(A^k)^*] \subseteq  R(A). $ 

{\rm (c)} \ $A^kA^{\dagger} -  A^{\dagger}A^k $ is nonsingular
  $ \Leftrightarrow  r \left[ \begin{array}{c}  A^k  \\ A^* \end{array}
  \right] =  r[\, A^k, \ A^*\, ] =  2r(A) = m  \Leftrightarrow
  r(A^k) = r( A )$  and $ R(A) \oplus R(A^*) = {\cal C}^m.$ 

{\rm (d)}  \ $ r(A) = r(A^2)$  and $A^kA^{\dagger} = A^{\dagger}A^k \Leftrightarrow$  $ A $ is EP. }

\medskip

\noindent {\bf Proof.}\, Writing $A^kA^{\dagger} -  A^{\dagger}A^k
= -[\, (A^{\dagger}A)A^{k-1} - A^{k-1}(AA^{\dagger}) \, ]$ and applying
Eq. (4.1) to  it, we obtain
\begin{eqnarray*}
r( \, A^kA^{\dagger} -  A^{\dagger}A^k \,)
& = & r\left[ \begin{array}{c} A^{\dagger}A^k \\ AA^{\dagger} \end{array}
\right] + r[\, A^kA^{\dagger}, \ A^{\dagger}A \,]  - r(AA^{\dagger}) -
 r(A^{\dagger}A) \\
& = & r\left[ \begin{array}{c} A^k \\ A^{\dagger} \end{array} \right]
+ r[\, A^k, \ A^{\dagger} \,] - 2r(A) = r\left[ \begin{array}{c} A^k
\\ A^* \end{array} \right] + r[\, A^k, \ A^* \,]  - 2r(A),
\end{eqnarray*}
as required for (6.5). The results in Parts (a)---(d) follow
immediately from  (6.5). \qquad $ \Box$

\medskip

>From the result in Theorem 6.8(a) we can extend the concept of EP matrix to power case: 
A square matrix $ A$ is said to be {\em k-power-EP} if $ R(A^k) \subseteq R(A^*)$ and 
 $R[(A^k)^*] \subseteq  R(A)$, where $ k \geq 2.$  It is believed that power-EP matrices, as 
a special type of matrices,  might also have some more interesting properties. But we do not intend further to discuss power-EP matrices and the related topics in this monograph. As an exercise, we leave the verification of the following result to the reader.

As an application, now we let $ M = \left[ \begin{array}{cc} A & B \\ B & A  \end{array} 
\right],$ where $ A, \ B  \in  {\cal C}^{m \times m}$.
Then
$$\displaylines{
\hspace*{1cm}
r(\, M^kM^{\dagger} -  M^{\dagger}M^k \,) \hfill
\cr
\hspace*{1cm}
 = 2r[\, (A+ B)^k , \ (A + B)^* \,] + 
2r[\, (A - B)^k , \ (A - B)^* \,] - 2r( \,A + B \,) - 2r( \,A - B  \,) \hfill
\cr
\hspace*{1cm} 
= r[\,( \,A + B \,)^k( \,A+ B \,)^{\dagger} - ( \,A + B \,)^{\dagger}(\, A + B \,) \,] + r[\,( \,A - B \,)( \,A  - B \,)^{\dagger} - ( \,A - B \,)^{\dagger}( \,A - B \,)  \,]. \hfill
\cr}
$$
In particular, $ M $ is $k$-power-EP if and only if $ A \pm B $ are $k$-power-EP.

\medskip

In general for any square matrix $ A $  and a  polynomial $ p(x)$, there is   
$$\displaylines{
\hspace*{2cm}
r[ \, p(A)A^{\dagger} -  A^{\dagger}p(A) \, ] = r \left[ \begin{array}{c}
 A^k  \\ A^* \end{array} \right]  + r[\, A^k, \ A^*\, ] - 2r(A),  \hfill (6.6)
\cr}
$$
and $ p(A)A^{\dagger} = A^{\dagger}p(A)$ holds if  and only if $ R[p(A)] \subseteq  R(A^*)$  
and  $ R[p(A)^*] \subseteq  R(A). $ 
 
\medskip

\noindent {\bf Theorem 6.9.}\, {\em Let $ A \in {\cal C}^{m \times m} $ be
given and  $k$ be an integer with $ k \geq 2$. Then

{\rm (a)}  \ $r[ \, A(A^k)^{\dagger} -  (A^k)^{\dagger}A \, ] = r \left[ \begin{array}{c}
 A^k  \\ A^kA^* \end{array} \right]  + r[\, A^k, \ A^*A^k \, ] - 2r(A^k).$

{\rm (b)} \ $r[ \, A(A^k)^{\dagger} -  (A^k)^{\dagger}A \, ] = 2 r[\, A, \ A^* \, ] - 2r(A),$ 
if $ r(A) = r(A^2).$ 

{\rm (c)} \ $ A(A^k)^{\dagger} = (A^k)^{\dagger}A  \Leftrightarrow  
r \left[ \begin{array}{c} A^k  \\ A^kA^* \end{array} \right] = r[\, A^k, \ A^*A^k \, ] = r(A^k) 
\Leftrightarrow R(A^k) = R(A^*A^k) \ and \ R[(A^k)^*] = R[A(A^k)^*].$ 

{\rm (d)}\, If  $ A(A^k)^{\dagger} = (A^k)^{\dagger}A,$  then  $A^kA^{\dagger} = A^{\dagger}A^k.$

{\rm (e)} \ $ r(A) = r(A^2)$ and $A(A^k)^{\dagger} = (A^k)^{\dagger}A \Leftrightarrow$ A is EP.  
 }

\medskip

\noindent {\bf Proof.}\, It follows by (2.2) and block elementary operations that 
\begin{eqnarray*}
r[ \, A(A^k)^{\dagger} -  (A^k)^{\dagger}A \, ]
& = & r\left[ \begin{array}{ccc} (A^k)^*A^k(A^k)^* & 0 &(A^k)^*  \\  0 & - (A^k)^*A^k(A^k)^* &
 (A^k)^*A  \\  A(A^k)^* &  (A^k)^* & 0 \end{array} \right] - 2r(A^k) \\
& = & r\left[ \begin{array}{ccc} (A^k)^*A^k(A^k)^* & (A^k)^*A^{k-1}(A^k)^*   &(A^k)^*  \\  0 &
 0 & (A^k)^*A  \\  A(A^k)^* &  (A^k)^* & 0 \end{array}
 \right] - 2r(A^k) \\
& = &  r\left[ \begin{array}{c} A^k \\ A^kA^* \end{array} \right]
+ r[\, A^k, \ A^*A^k \,] - 2r(A^k), 
\end{eqnarray*}
establishing Part (a). Parts (b), (c) and (e) follow immediately  from Part (a).
 Combining Part (c) and Theorem 6.8(a) yields the implication in Part (d). \qquad $\Box$ 

\medskip

\noindent {\bf Theorem 6.10.}\, {\em Let $ A \in {\cal C}^{m \times m}$  be
given. Then 

{\rm (a)} \ $r[\, A(\, AA^{\dagger} - A^{\dagger}A\,) - (\, AA^{\dagger} - A^{\dagger}A\,)A \,] 
= r \left[ \begin{array}{ccc} A & A^2  & A^*  \\
A^2 & 0  & 0 \\ A^* & 0 & 0 \end{array} \right] - 2r(A).$

{\rm (b)} \ $A$ commutes with $ AA^{\dagger} - A^{\dagger}A \Leftrightarrow 
 r \left[ \begin{array}{ccc} A & A^2  & A^*  \\ A^2 & 0  & 0 \\ A^* & 0 & 0 \end{array} 
\right] = 2r(A).$

{\rm (c)} \ $r[\, A(\, AA^{\dagger} - A^{\dagger}A\,) - (\, AA^{\dagger} - A^{\dagger}A\,)A \,] 
= 2r[\, A, \ A^* \, ] - 2r(A),$ if $ r(A) = r(A^2)$. 

{\rm (d)} \ $r(A) = r(A^2)$  and $A$ commutes with $AA^{\dagger} - A^{\dagger}A
\Leftrightarrow  A $ is EP. } 

\medskip 

\noindent {\bf Proof.}\, Notice that $ A(\, AA^{\dagger} - A^{\dagger}A\,) - (\, AA^{\dagger} - A^{\dagger}A\,)A = A^2A^{\dagger} + A^{\dagger}A^2 -2A.$ Thus according to (2.2), we get 
\begin{eqnarray*}
r( \, A^2A^{\dagger} + A^{\dagger}A^2 -2A \,)
 & = & r \left[ \begin{array}{ccc} A^*AA^*  & 0  & A^*  \\
0  & A^*AA^*  & A^*A^2 \\ A^2A^* & A^* & 2A \end{array} \right] - 2r(A) \\
 & = & r \left[ \begin{array}{ccc} A^*A  & 0  & A^*  \\
0  & AA^*  & A^2 \\ A^2 & A^* & 2A \end{array} \right] - 2r(A) \\
 & = & r \left[ \begin{array}{ccc} 0  & 0  & A^*  \\
-A^3 & AA^*  & A^2 \\ -A^2 & A^* & 2A \end{array} \right] - 2r(A) \\
 & = & r \left[ \begin{array}{ccc} 0  & 0  & A^*  \\
0 & 0  & -A^2 \\ -A^2 & A^* & 2A \end{array} \right] - 2r(A) \\
&= & r \left[ \begin{array}{ccc} A & A^2  & A^*  \\ A^2 & 0  & 0 \\ A^* & 0 & 0 \end{array} \right] - 2r(A),
\end{eqnarray*}
establishing Part (a).  Parts (b), (c) and (d) are are direct consequences of Part (a). \qquad $\Box$ 

\medskip

\noindent {\bf Theorem 6.11.}\, {\em Let $ A \in {\cal C}^{m \times m}$ be given and  $k$ be an integer with $ k \geq 2$.  Then 

{\rm (a)} \ $r[\, (AA^{\dagger})(A^*A) -  (A^*A)(AA^{\dagger}) \, ] = 2r[\, A, \ A^*A^2 \, ] - 2r( A ).$

{\rm (b)}  \ $r[\, (A^{\dagger}A)(AA^*) -  (AA^*)(A^{\dagger}A) \, ] 
= 2r\left[ \begin{array}{c} A \\ A^2A^* \end{array} \right] - 2r( A ).$

{\rm (c)} \ $r[\, (AA^{\dagger})(A^*A)^k -  (A^*A)^k(AA^{\dagger}) \, ] = 2r[\, A, \ (A^*A)^kA \, ] - 2r( A ).$

{\rm (d)} \ $r[\, (A^{\dagger}A)(AA^*)^k -  (AA^*)^k(A^{\dagger}A) \, ] 
= 2r\left[ \begin{array}{c} A \\ A(AA^*)^k \end{array} \right] - 2r( A ).$

{\rm (e)} \ $AA^{\dagger}$ commutes with $ A^*A \Leftrightarrow R(A^*A^2) \subseteq R(A)$. 

{\rm (f)} \ $A^{\dagger}A$ commutes with $ AA^* \Leftrightarrow R[A(A^2)^*] \subseteq R(A^*)$. 

{\rm (g)} \ $r[\, (AA^{\dagger})(A^*A) -  (A^*A)(AA^{\dagger}) \, ] = 
r[\, (A^{\dagger}A)(AA^*) -  (AA^*)(A^{\dagger}A) \, ] = 2r[\, A, \ A^* \, ] - 2r( A ),$ if 
$r(A) = r(A^2)$. 

{\rm (h)} \ $r(A) = r(A^2)$ and $AA^{\dagger}$ commutes with $ A^*A \Leftrightarrow r(A) = r(A^2)$ and 
$A^{\dagger}A$ commutes with $ AA^* \Leftrightarrow A $ is EP. }   

\medskip 

\noindent {\bf Proof.}\, Notice that both $AA^{\dagger}$ and $A^{\dagger}A $ are idempotent. 
The two rank equalities  in Parts (a)---(d) can trivially be derived from  (4.1).
The results in Parts (e)---(h) are direct consequences of  Parts (a) and (b).   \qquad $\Box$

\medskip

\noindent {\bf Theorem 6.12.}\, {\em Let $ A \in {\cal C}^{m \times m}$ given and  $k$ be an integer with $ k \geq 2$. Then 

{\rm (a)} \ $r(\, A^{\dagger}A^*A -  A^*AA^{\dagger} \, ) = r[\, A(A^*A) -  (A^*A)A \, ].$

{\rm (b)} \   $r(\, A^{\dagger}AA^* -  AA^*A^{\dagger} \, ) = r[\, A(AA^*) -  (AA^*)A \, ].$

{\rm (c)} \ $r[\, A^{\dagger}(A^*A)^k -  (A^*A)^kA^{\dagger} \, ] = r[\, A(A^*A)^k -  (A^*A)^kA \, ].$

{\rm (d)} \  $r[\, A^{\dagger}(AA^*)^k -  (AA^*)^kA^{\dagger} \, ] = r[\, A(AA^*)^k -  (AA^*)^kA \, ].$

{\rm (e)}\, The following statements are equivalent$:$

\hspace*{0.5cm} {\rm (1)}  $A^{\dagger}$ commutes with $ AA^*.$

\hspace*{0.5cm} {\rm (2)}  $A^{\dagger}$ commutes with $ A^*A.$

\hspace*{0.5cm} {\rm (3)}  $A^{\dagger}$ commutes with $ (AA^*)^k.$

\hspace*{0.5cm} {\rm (4)}  $A^{\dagger}$ commutes with $ (A^*A)^k.$

\hspace*{0.5cm} {\rm (5)\cite{EI}} $A$ commutes with $ AA^*.$

\hspace*{0.5cm} {\rm (6)\cite{EI}} $A$ commutes with $ A^*A.$

\hspace*{0.5cm} {\rm (7)\cite{EI}} $A$ commutes with $ (AA^*)^k.$

\hspace*{0.5cm} {\rm (8)\cite{EI}} $A$ commutes with $ (A^*A)^k.$

\hspace*{0.5cm} {\rm (9)\cite{EI}} $A$ is normal$,$ i.e.$,$ $AA^* = A^*A.$ }

\medskip 

\noindent {\bf Proof.}\, Notice that $A^{\dagger}A^*A -  A^*AA^{\dagger}=- (\,  A^* -  
A^{\dagger}A^*A \,)$. Thus by (2.1) we find that
\begin{eqnarray*}
r(\,  A^* -  A^{\dagger}A^*A \, ) & = &  r\left[ \begin{array}{cc} A^*AA^* & A^*A^*A  \\ 
A^* & A^* \end{array} \right] - r(A) \\
& = &r\left[ \begin{array}{cc} 0 &  A^*A^*A - A^*AA^*  \\ A^* & 0 \end{array} \right] - r(A)
 = r(\, AA^*A - A^*A^2 \,),
\end{eqnarray*}
establishing Part (a). Similarly we can establish Parts (b)---(c). The equivalence of (1) and (5),  
(2) and (6), (3) and (7),  (4) and (8) in Part (c) follow from the four  formulas in Parts (a)---(b). 
The equivalence of (5)---(9) in Part (c) were presented in \cite{EI}.  \qquad $\Box$   

\medskip

A square matrix $ A $ is said to be bi-EP, if $A $ and its Moore-Penrose inverse 
$A^{\dagger}$ satisfy $(AA^{\dagger})(A^{\dagger}A) = (A^{\dagger}A) (AA^{\dagger})$. This special type of matrices were examined by
 Campbell and Meyer \cite{CM1}, Hartwig and Spindelb\"{o}ck \cite{HaSp1}, \cite{HaSp2}.

Just as for EP matrices and conjugate EP matrices, bi-EP matrices can also be characterized by a 
rank equality. 

\medskip

\noindent {\bf Theorem 6.13.}\, {\em Let $ A \in {\cal C}^{m \times m} $ be
given. Then

{\rm (a)} \ $r[ \, (AA^{\dagger})(A^{\dagger}A) -  (A^{\dagger}A) (AA^{\dagger}) \, ] =  2r[\, A, \ A^* \,] + 2r( A^2 ) -4r( A ). $

{\rm (b)}\ $r[ \, A^2 - A^2(A^{\dagger})^2A^2 \, ]  = r[\, A, \ A^*\,] + r( A^2 ) -
2r( A ).$

{\rm (c)}\  $r[ \, (AA^{\dagger})(A^{\dagger}A) -  (A^{\dagger}A) (AA^{\dagger}) \, ] =
  2r[\, A, \ A^* \,] - 2r( A ),$ if $ r( A) = r( A^2 ).$ 

{\rm (d)} \ $r[ \, A^2 - A^2(A^{\dagger})^2A^2 \, ]  = r[\, A, \ A^*\,]  - r( A ),$ 
if $ r( A) = r( A^2 ).$ 

{\rm (e)}\, The following four statements are equivalent$:$

\hspace*{0.5cm} {\rm (1)} \ $ (AA^{\dagger})(A^{\dagger}A) = (A^{\dagger}A) (AA^{\dagger}),$  i.e.$,$ 
 $ A$  is bi-EP. 

\hspace*{0.5cm} {\rm (2)} \ $(A^{\dagger})^2 \in \{ (A^2)^- \}.$  

\hspace*{0.5cm} {\rm (3)} \ $ r[\, A, \ A^*\,] = 2r( A ) - r( A^2 ).$ 

\hspace*{0.5cm} {\rm (4)} \ $ {\rm dim}[R(A)\cap R(A^*)] = r(A^2).$

{\rm (f)\cite{HaSp2}} \ $ A $ is bi-EP and $r( A ) = r( A^2 )$ $\Leftrightarrow$ 
$ A^2(A^{\dagger})^2A^2 = A^2$ and $r( A ) = r( A^2 )$ $\Leftrightarrow$  $A$  is EP.}
 
\medskip

\noindent  {\bf Proof.}\,  Note that both $AA^{\dagger} $ and $ A^{\dagger}A $ are
 Hermitian idempotent and $ R(A^{\dagger}) = R(A^*)$. We find  by (3.29) that 
 \begin{eqnarray*}
 r[ \, (AA^{\dagger})(A^{\dagger}A) -  (A^{\dagger}A) (AA^{\dagger}) \, ] & = & 2r[\, AA^{\dagger}, \ A^{\dagger}A \, ] +
  2r[\,(AA^{\dagger})(A^{\dagger}A)\,]  -  2r(AA^{\dagger}) - 2r(A^{\dagger}A) \\
& = &  2r[ \, A, \ A^* \,]  +2r( A^{\dagger} A^{\dagger}) - 4r(A) \\
& = &  2r[\, A, \ A^* \,] + 2r( A^2 ) -4r( A ), 
\end{eqnarray*}
establishing Part (a). Applying (2.8) and then the rank cancellation laws 
in (1.8)  to $ A^2 - A^2(A^{\dagger})^2A^2$, we obtain
\begin{eqnarray*}  
r[ \, A^2 - A^2(A^{\dagger})^2A^2 \, ] & = & r \left[ \begin{array}{ccc} 
A^*A^* &  A^*AA^*  & 0  \\
A^*AA^* & 0  & A^*A^2 \\
0 & A^2A^* & -A^2 \end{array} \right] - 2r(A) \\ 
& = & r \left[ \begin{array}{ccc} 
A^*A^* &  A^*AA^*  & 0  \\
A^*AA^* & A^*A^2A^*  & 0 \\
0 &  0 & -A^2 \end{array} \right] - 2r(A) \\ 
& = & r \left[ \begin{array}{cc} 
A^*A^* &  A^*AA^* \\
A^*AA^* & A^*A^2A^* \end{array} \right]  + r(A^2) - 2r(A) \\ 
& = & r \left[ \begin{array}{cc} 
A^*A^* &  A^*A \\ AA^* & A^2 \end{array} \right]  + r(A^2) - 2r(A) \\ 
& = & r \left( \, [\, A, \  A^* \, ]^* [\, A, \ A^* \,]  \right) + r(A^2) - 2r(A) \\ 
& = & r[\, A,\ A^* \,] + r(A^2) - 2r(A), 
\end{eqnarray*}
as required for Part (b). The equivalence of  (1)---(3) in Part (e) follows from the two formulas in
 Parts (a) and (b)  The equivalence of (3) and (4)  in Part (e) follows from a well-known rank formula
 $ r[\, A, \  B \,] = r(A) + r(B) - {\rm dim}[\,R(A)\cap R(B)\,]$.  \qquad  $ \Box$ 

\medskip 

The above work can also be extended to the conjugate case. 

\medskip 

\noindent {\bf Theorem 6.14.}\, {\em Let $ A \in {\cal C}^{m \times m} $ be
given. Then
$$
r[ \, (AA^{\dagger})( \overline{A^{\dagger}A }) -  (\overline{A^{\dagger}A})(AA^{\dagger}) \, ] =
  2r[\, A, \ A^T \,] + 2r( A\overline{A} ) -4r( A ). 
$$
In particular$,$
$$
(AA^{\dagger})( \overline{A^{\dagger}A }) = (\overline{A^{\dagger}A})(AA^{\dagger}) \Leftrightarrow
r[\, A, \ A^T \,] = 2r(A) - r( A\overline{A} ).
$$ } 
{\bf Proof.} \ Follows from (3.29) by noticing that both 
$ AA^{\dagger}$ and  $ \overline{A^{\dagger}A}$ are idempotent.  \qquad $\Box$. 

\medskip

Based on the above results, a parallel concept to bi-EP matrix now can be introduced: A square matrix $ A$ is said to be {\em conjugate bi-EP} if $(AA^{\dagger})( \overline{A^{\dagger}A }) = (\overline{A^{\dagger}A})(AA^{\dagger})$. The properties and applications of this special type of matrices remain to further study. 

We next consider rank equalities related to  star-dagger matrices. A square matrix $ A$ is said to be star-dagger if $ A^*A^{\dagger}= A^{\dagger}A^*$. 
This special type of matrices were well investigated by Hartwig and Spindelb\"{o}ck \cite{HaSp2}, and later 
by Meenakshi and Rajian \cite{MR}.

\medskip

\noindent  {\bf Theorem 6.15.}\, {\em Let $ A \in {\cal C}^{ m \times m } $
be given. Then

 {\rm (a)} \ $ r( \, A^*A^{\dagger} -  A^{\dagger}A^* \, )
 = r( \, AA^*A^2 - A^2A^*A \, ).$

 {\rm (b)} \ $ r( \, AA^*A^{\dagger}A -  AA^{\dagger}A^*A \, )
 = r( \, AA^*A^2 - A^2A^*A\, ).$

 {\rm (c)} \ $ r( \, A^*A^{\dagger} -  A^{\dagger}A^* \, )
 = r( \, AA^* - A^*A\, ),$ if $ A $ is EP.

{\rm (d)}\, The following statements are equivalent
(Hartwig and Spindelb\"{o}ck \cite{HaSp2}):

\hspace*{0.5cm} {\rm (1)} \ $ A^*A^{\dagger} = A^{\dagger}A^*,$ i.e.$,$ $ A $ is
 star-dagger. 

\hspace*{0.5cm} {\rm (2)} \ $ AA^*A^{\dagger}A =  AA^{\dagger}A^*A$.

\hspace*{0.5cm} {\rm (3)} \ $ AA^*A^2 = A^2A^*A.$

{\rm (e)\cite{MJ4}} \ $ A $ is both EP and star-dagger $\Leftrightarrow A$ is normal.}

\medskip

\noindent  {\bf Proof.}\, We find by  (2.2) that 
 \begin{eqnarray*}
 r( \, A^*A^{\dagger} - A^{\dagger}A^* \, )
  & = & r \left[ \begin{array}{ccc} -A^*AA^*  & 0  & A^*   \\
  0 &  A^*AA^* & A^*A^*  \\ A^*A^* &  A^* & 0 \end{array} \right] - 2r(A)\\
 &= & r \left[ \begin{array}{ccc} 0 &  0 & A^* \\ (A^*)^2AA^* - A^* A(A^*)^2 & 0 \\ 0 &
 0 & A \\  A & 0  & 0  \end{array} \right] - 2r(A)
 = r( \, AA^*A^2 - A^2A^*A\, ),
 \end{eqnarray*}
 as required for Part (a). Similarly we can establish Part (b) by  (2.2). 
The formula in in Part (c)  is derived from Part (a), and Part (d) is 
direct consequences of Parts (a) and (b). Part (e) comes from Part (c).  \qquad  $ \Box$ 

\medskip

As pointed out by Hartwig and Spindelb\"{o}ck \cite{HaSp2}, the class of star-dagger matrices are 
quite inclusive.  Normal matrix, partial isometry (i.e., $A^{\dagger} = A^*$),  idempotent matrix, 
2-nilpotent matrix, power Hermitian  matrix (i.e., $A^* = A^k$),  and so on are all special 
cases of star-dagger matrices, this assertion can easily be seen from the statement (3) in Theorem 6.15(d).  

The results in Theorem 6.15 can be extended 
to general cases. Below are three of them. Their  proofs are much similar to 
that of Theorem 6.15 and are, therefore, omitted.    
 
\medskip

\noindent  {\bf Theorem 6.16.}\, {\em Let $ A \in {\cal C}^{ m \times m } $
be given. Then

 {\rm (a)} \ $ r( \, A^*AA^*A^{\dagger} -  A^{\dagger}A^*AA^* \, )
 = r[ \, (AA^*)^2A^2 - A^2(A^*A)^2 \, ].$

 {\rm (b)} \ $ r[\, (AA^*)^2A^{\dagger}A -  AA^{\dagger}(A^*A)^2 \, ]
 = r[ \, (AA^*)^2A^2 - A^2(A^*A)^2 \, ].$

 {\rm (c)} \ $ A^*AA^*A^{\dagger} = A^{\dagger}A^*AA^* \Leftrightarrow 
(AA^*)^2A^{\dagger}A = AA^{\dagger}(A^*A)^2 \Leftrightarrow  (AA^*)^2A^2 = A^2(A^*A)^2.$ 

{\rm (d)}\, If $ A^*A^{\dagger} = A^{\dagger}A^*,$ then $ A^*AA^*A^{\dagger} = A^{\dagger}A^*AA^*.$ 
}

\medskip

 \noindent {\bf Theorem 6.17.}\, {\em Let $ A \in {\cal C}^{ m \times m } $
 be given and $k $ be an integer with $k \geq 2$. Then

{\rm (a)} \ $ r[ \, (A^*)^kA^{\dagger} - A^{\dagger}(A^*)^k \,] = r( \, AA^*A^{k+1} -
 A^{k+1}A^*A \,). $

{\rm (b)} \ $ r[ \, A(A^*)^kA^{\dagger}A - AA^{\dagger}(A^*)^kA \,] = r( \, AA^*A^{k+1} -
 A^{k+1}A^*A \,). $

{\rm (c)} \ $ (A^*)^kA^{\dagger} = A^{\dagger}(A^*)^k \Leftrightarrow
 A(A^*)^kA^{\dagger}A = AA^{\dagger}(A^*)^kA  \Leftrightarrow
   AA^*A^{k+1} = A^{k+1}A^*A.$ 

{\rm (d)}\, If $A^{k+1} = A,$  or $A^{k+1} =0,$ or  $AA^* = A^*A,$ or $AA^*A = A,$ 
then  $(A^*)^kA^{\dagger} = A^{\dagger}(A^*)^k$ holds. } 

\medskip

In general for any square matrix $ A $  and a  polynomial $ p(x)$, there is   
$$\displaylines{
\hspace*{2cm}
r[ \, p(A^*)A^{\dagger} -  A^{\dagger}p(A^*) \, ] = r [\, AA^*p(A)A - Ap(A)A^*A \,]. \hfill (6.7)
\cr}
$$
In particular, $ p(A^*)A^{\dagger} = A^{\dagger}p(A^*)$ holds if  and only if $  AA^*p(A)A = Ap(A)A^*A$. 
 
\medskip

\noindent {\bf Theorem 6.18.}\, {\em Let $ A \in {\cal C}^{ m \times m } $
 be given and $k $ be an integer with $k \geq 2$. Then

{\rm (a)}  \ $ r[ \, A^*(A^k)^{\dagger} - (A^k)^{\dagger}A^* \,] =
 r[ \, A^k(A^k)^*A^{k + 1} - A^{k + 1}(A^k)^*A^k \,]. $

{\rm (b)} \  $ A^*(A^k)^{\dagger} = (A^k)^{\dagger}A^* \Leftrightarrow 
A^k(A^k)^*A^{k + 1} = A^{k + 1}(A^k)^*A^k$.}

\medskip


Next are several results on ranks of matrix expressions involving powers of
 the Moore-Penrose inverse of a matrix. 

\medskip

\noindent  {\bf Theorem 6.19.}\, {\em Let $ A \in {\cal C}^{ m \times m } $ be given.
 Then

 {\rm (a)} \ $ r[\,I_m \pm A^{\dagger} \,] = r( \, A^2 \pm AA^*A \, ) -
 r( A ) + m .$

 {\rm (b)} \  $ r[\,I_m -  (A^{\dagger})^2 \,] = r( \,A^2 +  AA^*A \, ) +
 r( \, A^2 -  AA^*A \, ) - 2r( A )  + m .$

 {\rm (c)}  \ $ r[\,I_m \pm A^{\dagger} \,] = r( \, A \pm AA^* \, ) - r( A )
 + m, $  if $ A $ is EP.

 {\rm (d)}  \ $ r[\,I_m -  (A^{\dagger})^2 \,] = r( \,A + AA^* \, ) +
 r( \, A - AA^* \, ) - 2r(A) + m,$  if $ A $ is EP.

 {\rm (e)} \  $ r[\,I_m \pm A^{\dagger} \,] = r( \, A \pm A^2 \, ) - r( A )
  + m,$ if $ A $ is Hermitian.

{\rm (f)}  \ $ r[\,I_m -  (A^{\dagger})^2 \,] = r( \,A + A^2 \, ) +
 r( \, A -  A^2 \, ) - 2r( A ) + m,$
 if $ A $ is Hermitian. } 

\medskip

 \noindent {\bf Proof.}\, By (2.1) we easily obtain
\begin{eqnarray*} 
  r(\,I_m -  A^{\dagger} \,)& = & 
r\left[ \begin{array}{cc}  A^*AA^* &  A^* \\  A^*  & I_m  \end{array} \right] -r( A ) \\
& = & r\left[ \begin{array}{cc}  A^*AA^* - A^*A^* &  0 
\\ 0  & I_m  \end{array} \right] -r( A ) =  r( \, AA^*A - A^2 \, )+ m  - r( A ),
\end{eqnarray*} 
and
\begin{eqnarray*} 
  r(\,I_m + A^{\dagger} \,)& = & 
r\left[ \begin{array}{cc}  -A^*AA^* &  A^* \\  A^*  & I_m  \end{array} \right] -r( A ) \\
& = & r\left[ \begin{array}{cc}  -A^*AA^* - A^*A^* &  0 
\\ 0  & I_m  \end{array} \right] -r( A ) =  r( \, AA^*A + A^2 \, )+ m - r( A ).
\end{eqnarray*} 
Both of the above  are exactly Part (a). Next applying  (1.12) to
 $ I_m -  (A^{\dagger})^2 $ we obtain
\begin{eqnarray*} 
r[\,I_m -  (A^{\dagger})^2 \,]  & = &  r(\,I_m +  A^{\dagger} \,) + r(\,I_m - A^{\dagger} \,) - m \\
  & = & r( \,A^2 +  AA^*A \, ) +  r( \, A^2 -  AA^*A \, ) - 2r( A )  + m,
 \end{eqnarray*} 
establishing Part (b). The results in
 Parts (c)---(f) follow directly from Parts (a) and (b).
  \qquad $ \Box$

\medskip

\noindent  {\bf Theorem 6.20.}\, {\em Let $ A \in {\cal C}^{ m \times m } $ be given.
 Then

{\rm (a)} \ $ r[\, A^{\dagger} \pm (A^{\dagger})^2 \,] =
  r( \, A^2 \pm AA^*A \, ) = r[ \, A \pm A(A^{\dagger})^2A \, ].$

 {\rm (b)} \ $ r[\, A^{\dagger} - (A^{\dagger})^2 \,]  =
 r( \, A - AA^* \, ) = r( \, A - A^*A \, ),$ if $ A $ is EP.

 {\rm (c)} \  $ r[\, A^{\dagger} \pm (A^{\dagger})^2 \,]  =
  r( \, A \pm A^2 \, ),$ if $ A $ is Hermitian.

{\rm (d)} \  $r[\, A^{\dagger} - (A^{\dagger})^2 \,] = r( A^{\dagger}) - r[\,(A^{\dagger})^2\,],$ i.e.$,$ 
$(A^{\dagger})^2 \leq_{rs} A^{\dagger} \Leftrightarrow  r( \, AA^*A - A^2 \, ) = r( A)  - r( A^2),$ i.e.$,$ 
$A^2 \leq_{rs} AA^*A.$

{\rm (e)} \  $(A^{\dagger})^2  = A^{\dagger} \Leftrightarrow
 AA^*A = A^2  \Leftrightarrow  A = (AA^{\dagger})(A^{\dagger}A)
 \Leftrightarrow (A^{\dagger})^2 \in \{ A^- \}.$

{\rm (f)} \ $ (A^{\dagger})^2  = A^{\dagger} \Leftrightarrow
 AA^* = A^*A = A,$ if $ A $ is EP.

{\rm (g)} \  $ (A^{\dagger})^2 =  A^{\dagger} \Leftrightarrow
 A^2 = A,$  \ if $ A $ is Hermitian. }

\medskip

\noindent {\bf Proof.}\, It follows first from  (1.11) that 
$$\displaylines{
\hspace*{2cm}
r [\, A^{\dagger} - (A^{\dagger})^2 \,] = r( \, I_m - A^{\dagger} \, )
 + r( A )- m, \hfill 
\cr
\hspace*{2cm}
r [\, A^{\dagger} + (A^{\dagger})^2 \,] = r( \, I_m + A^{\dagger} \, )
 + r( A )- m. \hfill 
\cr}
$$ 
Then we have the first two equalities in Part (a) by Theorem 6.15(a).
Note that
$$\displaylines{
\hspace*{2cm}
A[ \, A^{\dagger} \pm (A^{\dagger})^2 \,]A =  A \pm A(A^{\dagger})^2A
\ \ {\rm and}  \ \ A^{\dagger} [ \, A \pm A(A^{\dagger})^2A \,]A^{\dagger}
= A^{\dagger} \pm (A^{\dagger})^2. \hfill
\cr}
$$
It follows that
$$\displaylines{
\hspace*{2cm}
r[ \, A^{\dagger} \pm (A^{\dagger})^2 \,]
= r[\,  A \pm A(A^{\dagger})^2A \,]. \hfill
\cr}
$$
Thus we have the second equality in Part (a).  Parts (b)---(g) follow from
Part (a). \qquad  $ \Box$

\medskip

\noindent  {\bf Theorem 6.21.}\, {\em Let $ A \in {\cal C}^{ m \times m } $
be given. Then

{\rm (a)} \ $ r[\, A^{\dagger} - (A^{\dagger})^3 \, ]
 = r( \, A^2 + AA^*A \, ) + r( \,  A^2 - AA^*A \, ) - r( A ).$ 

{\rm (b)} \ $ r[\, A^{\dagger} - (A^{\dagger})^3 \, ]
 =  r( \, A + AA^* \, ) + r( \,  A - AA^* \, ) - r( A ),$  if $ A $ is EP.

{\rm (c)} \ $ r[\, A^{\dagger} - (A^{\dagger})^3 \, ]
 =  r( \, A + A^2 \, ) + r( \,  A - A^2 \, ) - r( A )
 = r( \,  A - A^3 \, ),$ if $ A $ is Hermitian.

{\rm (d)} \ $ (A^{\dagger})^3 = A^{\dagger} \Leftrightarrow
  r( \,A^2 + AA^*A  \, ) + r( \, A^2 - AA^*A \, ) = r( A )
  \Leftrightarrow R( \, AA^*A + A^2 \, ) \cap R( \, AA^*A - A^2 \, ) = \{0\}
  $ and $R[( \, AA^*A + A^2 \, )^*] \cap R[( \, AA^*A - A^2 \, )^*]
  = \{0\}.$ }

\medskip

 \noindent {\bf Proof.}\, Applying the rank equality (1.15) to 
$A^{\dagger} - (A^{\dagger})^3$, we obtain 
 $$\displaylines{
\hspace*{2cm}
 r[\, A^{\dagger} - (A^{\dagger})^3 \, ] = r[\, A^{\dagger} +
 (A^{\dagger})^2 \, ]
 + r[\, A^{\dagger} - (A^{\dagger})^2 \, ] - r( A ).\hfill
 \cr}
$$ 
 Then putting Theorem 6.20(a) in it yields Part (a). The  results in
 Parts (b)---(d) follow
 all from Part (a). \qquad  $ \Box$ 

\medskip

\noindent  {\bf Theorem 6.22.}\, {\em Let $ A \in {\cal C}^{ m \times n } $
be given. Then

{\rm (a)} \ $ r(\, A^{\dagger} - A^* \, ) = r( \, A - AA^*A \, ).$ 

{\rm (b)} \  $ r(\, A^{\dagger} - A^*AA^* \,) = r( \, A - AA^*AA^*A \, ).$ 

{\rm (c)} \  $ r[\, A^{\dagger} - (A^k)^* \,] =
 r( \, A - A^kA^*A \, ) = r( \, A - AA^*A^k \, ).$ \\
In particular$,$

{\rm (d)} \  $ A^{\dagger} = A^*  \Leftrightarrow  AA^*A = A, \ i.e.,   \ A $ is partial isometry.

{\rm (e)} \ $ A^{\dagger} = A^*AA^*  \Leftrightarrow  AA^*AA^*A = A.$ 

{\rm (f)} \  $ A^{\dagger} = (A^k)^*  \Leftrightarrow 
 A^kA^*A = AA^*A^k = A.$ 
}

\medskip

\noindent {\bf Proof.}\, Follows from (2.1). \qquad $\Box$

\medskip

\noindent  {\bf Theorem 6.23.}\, {\em Let $ A \in {\cal C}^{ m \times m } $
be idempotent. Then

{\rm (a)}  \ $ r(\, A - A^{\dagger} \, ) = 2r[ \, A, \  A^* \, ] - 2r(A).$ 

{\rm (b)} \ $ r(\, 2A - AA^{\dagger} -A^{\dagger}A  \, ) = 2r[ \, A, \  A^* \, ] - 2r(A).$ 

{\rm (c)} \ $r[\, (AA^{\dagger})(AA^*) -  (AA^*)(AA^{\dagger}) \, ] = 
r[\, (A^{\dagger}A)(A^*A) -  (A^*A)(A^{\dagger}A) \, ] = 2r[\, A, \ A^* \, ] - 2r( A ).$ 

{\rm (d)}  \ $ r(\, A^{\dagger} - AA^{\dagger} A^{\dagger} A  \,) 
 = 2r[ \, A, \  A^* \, ] - 2r(A).$ 

{\rm (e)} \ $ r(\, A - AA^{\dagger}A^{\dagger} A  \,) 
 = r[ \, A, \  A^* \, ] - r(A).$ 

{\rm (f)} \ $ A^{\dagger} $ commutes with $A^*$. 

{\rm (g)} \ $ A^{\dagger} $ commutes with $A^*AA^*$. 

{\rm (h)} \ $ (AA^*)^2A^{\dagger}A =AA^{\dagger}A (A^*A)^2.$ \\
In particular$,$

{\rm (i)} \  $ A = A^{\dagger} \Leftrightarrow AA^{\dagger} + A^{\dagger}A = 2A \Leftrightarrow
 (AA^{\dagger})(AA^*) = (AA^*)(AA^{\dagger}) \Leftrightarrow (A^{\dagger}A)(A^*A) = (A^*A)(A^{\dagger}A)
\Leftrightarrow  A^{\dagger} = AA^{\dagger} A^{\dagger} A \Leftrightarrow 
 A = AA^{\dagger} A^{\dagger} A \Leftrightarrow  A $ is Hermitian.}

\medskip

\noindent {\bf Proof.}\,  Note that $ A, \ A^{\dagger} \in A\{2\}$ when 
$ A $ is idempotent. Thus we have by (5.1) that
\begin{eqnarray*}
r( \, A -  A^{\dagger} \,) &=& r\left[ \begin{array}{c}
A \\ A^{\dagger} \end{array} \right] + 
r[\, A, \ A^{\dagger} \,]   - r(A) - r(A^{\dagger}) \\
&=& r\left[ \begin{array}{c} A \\ A^* \end{array} \right] + 
r[\, A, \ A^* \,]  - 2r(A) = 2r[\, A, \ A^* \,]  - 2r(A),
\end{eqnarray*}
establishing Part (a).  Part (b) follows from Theorem 6.10(c), Part (c) follows from Theorem 6.11(c), 
and Parts (d) and (e) follow from Theorem 6.13(a) and (b). Parts (f)---(h) follow from Theorems 
6.15 and 6.16. Part (i) is a direct consequence of Parts (a)---(e). \qquad $\Box$

\medskip

\noindent  {\bf Theorem 6.24.}\, {\em Let $ A \in {\cal C}^{ m \times m } $
be tripotent$,$ that is$,$ $A^3 = A$. Then

{\rm (a)} \ $ r(\, A - A^{\dagger} \, ) = 2r[ \, A, \  A^* \, ] - 2r(A).$ 

{\rm (b)} \ $ r(\, A^2A^{\dagger} - A^{\dagger}A^2  \,) 
 = 2r[ \, A, \  A^* \, ] - 2r(A).$ 

{\rm (c)}  \ $ r[\, A(A^2)^{\dagger} - (A^2)^{\dagger}A \,] 
 = 2r[ \, A, \  A^* \, ] - 2r(A).$ 

{\rm (d)}  \ $r[\, A(\, AA^{\dagger} - A^{\dagger}A\,) - (\, AA^{\dagger} - A^{\dagger}A\,)A \,] 
= 2r[\, A, \ A^* \, ] - 2r(A)$.

{\rm (e)} \  $r[\, (AA^{\dagger})(AA^*) -  (AA^*)(AA^{\dagger}) \, ] = 2r[\, A, \ A^* \, ] - 2r( A ).$

{\rm (f)} \  $ r[\, (A^{\dagger}A)(A^*A) -  (A^*A)(A^{\dagger}A) \, ] = 2r[\, A, \ A^* \, ] - 2r( A )$. 

{\rm (g)} \ $ r[\, (AA^{\dagger})(A^{\dagger}A) - (A^{\dagger}A)(AA^{\dagger}) \,] 
 = 2r[ \, A, \  A^* \, ] - 2r(A).$ 

{\rm (h)} \ $ r[\, A^2 - A^2(A^{\dagger})^2A^2  \,]  = r[ \, A, \  A^* \, ] - r(A).$ 

{\rm (i)} \ $ (A^*)^2A^{\dagger} = A^{\dagger}(A^*)^2.$

{\rm (j)}\, The following nine statements are equivalent$:$

\hspace*{0.5cm} {\rm (1)}   $ A = A^{\dagger}.$
 
\hspace*{0.5cm} {\rm (2)}  $A^2A^{\dagger} = A^{\dagger}A^2.$

\hspace*{0.5cm} {\rm (3)}  $ A(A^2)^{\dagger} = (A^2)^{\dagger}A.$

\hspace*{0.5cm} {\rm (4)}  $ A(\, AA^{\dagger} - A^{\dagger}A\,) 
= (\, AA^{\dagger} - A^{\dagger}A\,)A.$ 

\hspace*{0.5cm} {\rm (5)}  $ (AA^{\dagger})(AA^*) = (AA^*)(AA^{\dagger})$.

\hspace*{0.5cm} {\rm (6)}  $ (A^{\dagger}A)(A^*A) = (A^*A)(A^{\dagger}A)$.

\hspace*{0.5cm} {\rm (7)} $ (AA^{\dagger})(A^{\dagger}A) = (A^{\dagger}A)(AA^{\dagger}).$
 
\hspace*{0.5cm} {\rm (8)}   $A^2 = A^2(A^{\dagger})^2A^2.$
 
\hspace*{0.5cm} {\rm (9)}  $ R(A) = R(A^*), \  i.e.,  \  A  \  is  \ EP.$ }

\medskip

\noindent {\bf Proof.}\,  Note that $ A, \ A^{\dagger} \in A\{2\}$ when 
$ A $ is tripotent. Thus we have by (5.1) that
\begin{eqnarray*}
r( \, A -  A^{\dagger} \,) &=& r\left[ \begin{array}{c}
A \\ A^{\dagger} \end{array} \right] + 
r[\, A, \ A^{\dagger} \,]   - r(A) - r(A^{\dagger}) \\
&=& r\left[ \begin{array}{c} A \\ A^* \end{array} \right] + 
r[\, A, \ A^* \,]  - 2r(A) = 2r[\, A, \ A^* \,]  - 2r(A),
\end{eqnarray*}
establishing Part (a).  Parts (b)---(h) follow respectively from (6.5), Theorem 6.9(b), 
 Theorem 6.10(c), Theorem 6.11(e), and Theorem 6.13(a) and (b). Part (i) follows from 
Theorem 6.17(d). Part (j) is a direct consequence of Parts (a)---(h). \qquad $\Box$

\medskip
 
The following result is motivated by a problem of Rao and Mitra \cite{RM} on the 
nonsingularity of a matrix of the form $ I + A - A^{\dagger}A$.

\medskip

\noindent  {\bf Theorem 6.25.}\, {\em Let $ A \in {\cal C}^{ m \times m }$ and 
$ 1 \neq \lambda \in {\cal C}$ be given. 
Then

{\rm (a)} \ $ r(\, I_m + A - A^{\dagger}A \, ) = 
r(\, I_m + A - AA^{\dagger} \, ) = r(A^2)  - r(A) +m.$ 

{\rm (b)} \ $ r(\, I_m - A - A^{\dagger}A \, ) = 
r(\, I_m - A - AA^{\dagger} \, ) = r(A^2)  - r(A) +m.$ 

{\rm (c)} \ $ r(\, \lambda I_m + A - A^{\dagger}A \, ) = r(\, \lambda I_m + A - AA^{\dagger} \, ) = 
 r[\, (\, \lambda - 1 \,)I_m + A \,].$ 

{\rm (d)} \ $ r(\, \mu I_m - A \,) = r( \, \mu I_m  + I_m - A - A^{\dagger}A \, ) = 
r( \, \mu I_m + I_m  - A - AA^{\dagger} \, ),$ when  $ \mu \neq 0$.              

{\rm (e)\cite{RM}} \ $I_m + A - A^{\dagger}A$ is nonsingular 
 $ \Leftrightarrow  $ $I_m + A - AA^{\dagger}$  is nonsingular 
$ \Leftrightarrow r(A^2) = r(A)$.}

\medskip

\noindent {\bf Proof.}\, Applying (2.1) and then (1.8), we find that
\begin{eqnarray*}
r( \, I_m  + A -  A^{\dagger}A \,) 
&=& r\left[ \begin{array}{cc} A^*AA^* & A^*A \\ 
A^*  & I_m + A \end{array} \right] - r(A) \\
&=& r\left[ \begin{array}{cc} AA^* & A \\ 
A^*  & I_m + A   \end{array} \right] - r(A) \\
&=& r\left[ \begin{array}{cc} 0 & A \\  -AA^*  & I_m \end{array} \right] - r(A) \\
&=& r\left[ \begin{array}{cc} A^2A^* & 0 \\  0 & I_m \end{array} \right] - r(A) = m + r(A^2) - r(A). 
\end{eqnarray*}
By symmetry, we also get $ r(\, I_m + A - AA^{\dagger} \, ) =  m + r(A^2)  - r(A).$
Both of them are the result in Part (a). Replace $ A$ by $-A$ to yield Part (b). Part (c) is derived  by (2.1) and (1.11).
Replace $ \lambda$ by $ \mu + 1$ in Part (c) to yield Part (d).  Part (e)  is a direct consequence of Part (a). \qquad $\Box$ 

\medskip

The result in Theorem 6.26(d) reveals an interesting fact that a square matrix $A$  has the same nonzero eigenvalues as 
 the  matrix $A + A^{\dagger}A - I_m $ (or $ A + AA^{\dagger} - I_m$) has. Of course, this result 
is trivial when $ A$ is nonsingular. 

\medskip 

\noindent  {\bf Theorem 6.26.}\, {\em Let $ A \in {\cal C}^{ m \times m }$ be given. 
Then

{\rm (a)}  \ $ r(\, I_m + A^k - A^{\dagger}A \, ) = 
r(\, I_m + A^k - AA^{\dagger} \, ) = r(A^{k+1})  - r(A) +m.$ 

{\rm (b)}  \ $ r(\, I_m - A^k - A^{\dagger}A \, ) = 
r(\, I_m - A^k - AA^{\dagger} \, ) = r(A^{k+1})  - r(A) +m.$ 

{\rm (c)}  \ $I_m + A^k - A^{\dagger}A$ is nonsingular $ \Leftrightarrow $ $I_m + A^k - AA^{\dagger}$  is nonsingular 
  $ \Leftrightarrow $ $I_m - A^k - A^{\dagger}A$ is nonsingular $ \Leftrightarrow $ $I_m - A^k - AA^{\dagger}$  is nonsingular 
$ \Leftrightarrow r(A^{k+1}) = r(A),$ i.e.$,$ $r(A) = r(A^2).$   

{\rm (d)} \ $ r(\, I_m + A^k - A^{\dagger}A \, ) = m - r(A) \Leftrightarrow  A^{k+1} =0.$ }

\medskip

\noindent {\bf Proof.}\, We only show the first equality in Part (a). Applying 
(2.1) and then  (1.8), we find that
\begin{eqnarray*}
r( \, I_m  + A^k -  A^{\dagger}A \,) 
&=& r\left[ \begin{array}{cc} A^*AA^* & A^*A \\ 
A^*  & I_m + A^k \end{array} \right] - r(A) \\
&=& r\left[ \begin{array}{cc} AA^* & A \\ 
A^*  & I_m + A^k  \end{array} \right] - r(A) \\
&=& r\left[ \begin{array}{cc} 0 & A \\  -A^kA^*  & I_m \end{array} \right] - r(A) \\
&=& r\left[ \begin{array}{cc} A^{k+1}A^* & 0 \\  0 & I_m \end{array} \right] - r(A) = m + r(A^{k+1}) - r(A). \qquad \Box 
\end{eqnarray*}

The rank equalities in Theorem 6.26 are still valid when replacing the Moore-Penrose inverse of $ A$ by any inner 
inverse  $ A^-$ of $A$.  We shall prove this in Chapter 23.

\medskip 

\noindent  {\bf Theorem 6.27.}\, {\em Let $ A \in {\cal C}^{ m \times m }$ be given. Then

{\rm (a)}  \ $ r(\, I_m  - AA^{\dagger} - A^{\dagger}A \, ) = 2r(A^2)  - 2r(A) +m.$ 

{\rm (b)}  \ $ I_m - AA^{\dagger} - A^{\dagger}A$ is nonsingular $ \Leftrightarrow$ $ r(A^2) = r(A)$. 

{\rm (c)} \  $AA^{\dagger} + A^{\dagger}A = I_m \Leftrightarrow $ $m$ is even$,$  $ r(A) =m/2$ and $ A^2 = 0$. 
}

\medskip

\noindent {\bf Proof.}\, Apply (3.8) to $  I_m  - AA^{\dagger} - A^{\dagger}A$ to yield
\begin{eqnarray*}
r( \, I_m  - AA^{\dagger} - A^{\dagger}A  \,) 
&=&  r(AA^{\dagger}A^{\dagger}A) + r(A^{\dagger}AAA^{\dagger}) -  r(AA^{\dagger}) - r(A^{\dagger}A) + m \\
&=&  2r(A^2) -  2r(A) + m, 
\end{eqnarray*}
as required in Part (a). Part (b) is obvious from Part (a). According to the rank formula in Part (a), 
 the equality $AA^{\dagger} + A^{\dagger}A = I_m$ holds if and only if $  2r(A^2)  - 2r(A) +m = 0$. 
This rank equality implies that $m$  must be even and $ r(A^2) = r(A) - m/2$. Contracting this rank equality  with 
the Frobenius rank inequality $ r(A^2) \geq 2r(A) - m$, we get Part (c). \qquad $\Box$

\markboth{YONGGE  TIAN }
{7. RANK EQUALITIES FOR  MATRICES AND THEIR MOORE-PENROSE INVERSES}

\chapter{Rank equalities for matrices and their  Moore-Penrose inverses}

\noindent We consider in this chapter ranks of various matrix expressions that involve
two or more matrices and their  Moore-Penrose inverses, and present their
various consequences, which can  reveal a series of intrinsic properties related to 
Moore-Penrose inverses of matrices. Most of the results obtained in this chapter are 
new and are not considered before. 

\medskip

\noindent {\bf Theorem 7.1.}\, {\em Let $ A, \, B  \in {\cal C}^{ m \times m }
 $ be  given. Then

{\rm (a)} \ The rank of $AA^{\dagger}B -  BA^{\dagger}A$ satisfies
$$\displaylines{
\hspace*{1.5cm}
r(  \, AA^{\dagger}B -  BA^{\dagger}A \, )
 = r \left[ \begin{array}{c} A \\  A^*B  \end{array} \right]
 + r[ \, A, \ BA^* \, ] - 2r ( A ). \hfill (7.1)
\cr}
$$

{\rm (b)} \ $AA^{\dagger}B =  BA^{\dagger}A \Leftrightarrow r \left[ \begin{array}{c} A \\  A^*B  \end{array} \right] = r[ \, A, \ BA^* \, ] = r ( A ) \Leftrightarrow
  R( BA^* ) \subseteq R( A ) \ and \ R( B^*A ) \subseteq R( A^* ).$ 

{\rm (c)} \ $ AA^{\dagger}B -  BA^{\dagger}A$ is nonsingular $ \Leftrightarrow  r \left[ \begin{array}{c} A \\  A^*B  \end{array}
 \right] =  r[ \, A, \ BA^* \, ] =  2r ( A ) = m  \Leftrightarrow
  R( A ) \oplus  R( BA^* ) = {\cal C}^m  \ and \ R( AB^* ) = R(A)
 \Leftrightarrow R( A^* ) \oplus  R( B^*A ) = {\cal C}^m \ and \
   R( A^*B) = R( A^* ).$ }
\medskip

\noindent  {\bf Proof.}\, Note that $AA^{\dagger} $ and $ A^{\dagger}A $ are idempotent and 
$ R(A^{\dagger}) = R(A^*)$. We have by Eq. (4.1) that
\begin{eqnarray*}
r( \, AA^{\dagger}B -  BA^{\dagger}A \,) & = & r \left[ \begin{array}{c} AA^{\dagger}B \\ 
A^{\dagger}A  \end{array} \right] + r[ \,BA^{\dagger}A , \  AA^{\dagger} \, ] - r (AA^{\dagger} ) - r(A^{\dagger}A) \\
& = & r \left[ \begin{array}{c} A^{\dagger}B \\ A  \end{array} \right]
 + r[ \,BA^{\dagger} , \  A \, ] - 2r(A) \\
 & = & r \left[ \begin{array}{c} A^*B \\ A \end{array} \right]
 + r[ \,BA^*, \  A \, ] - 2r(A),
\end{eqnarray*}
establishing  (7.1). The results in Parts (b) and (c) follow from it. \qquad  $ \Box$ 

\medskip

Clearly the results in Theorems 6.1 and 6.7 are special cases of the above
 theorem.

\medskip

 \noindent {\bf Theorem 7.2.}\, {\em Let $ A \in {\cal C}^{ m \times n},  \,
 B \in {\cal C}^{ m \times k}$ and $C \in {\cal C}^{ l \times n} $ be given.
 Then

 {\rm (a)} \ $ r(  \, AA^{\dagger} -  BB^{\dagger} \, )
 = 2r[ \, A, \ B \, ] - r ( A ) - r( B ).$

{\rm (b)} \  $ r(  \, A^{\dagger}A -  C^{\dagger}C \, )
 = 2 r \left[ \begin{array}{c} A \\ C  \end{array} \right]
  - r ( A ) - r( C ).$

{\rm (c)} \ $ r(  \, AA^{\dagger} +  BB^{\dagger} \, )
 = r[ \, A, \ B \, ],$ that is$,$ $ R(  \, AA^{\dagger} +  BB^{\dagger} \, ) = R[ \, A, \ B \, ].$

{\rm (d)} \ $ r(  \, A^{\dagger}A  + C^{\dagger}C \, )
 =  r \left[ \begin{array}{c} A \\ C  \end{array} \right],$ that is$,$ $ R( \, A^{\dagger}A + 
 C^{\dagger}C \, ) = R[ \, A^*, \ C^* \, ].$ \\
In particular$,$

{\rm (e)} \ $ AA^{\dagger} =  BB^{\dagger} \Leftrightarrow
  R(A) = R( B ).$

{\rm (f)}  \ $ A^{\dagger}A =  C^{\dagger}C  \Leftrightarrow
  R(A^*) = R(C^*).$

{\rm (g)} \  $ r(  \, AA^{\dagger} -  BB^{\dagger} \, ) = r(AA^{\dagger}) -  r(BB^{\dagger}) 
\Leftrightarrow  R(B) \subseteq R(A).$

{\rm (h)} \ $ r( \, A^{\dagger}A -  C^{\dagger}C \, ) =  r(A^{\dagger}A) -  r(C^{\dagger}C) 
\Leftrightarrow  R(C^*) \subseteq R(A^*). $

{\rm (i)} \ $ r( \, AA^{\dagger} -  BB^{\dagger} \, ) = m
  \Leftrightarrow  r[ \, A, \ B \, ]
  = r( A ) + r( B )
  = m    \Leftrightarrow  R(A) \oplus R(B) = {\cal C}^m.$

{\rm (j)} \ $ r( \, A^{\dagger}A -  C^{\dagger}C \, ) = n
  \Leftrightarrow r \left[ \begin{array}{c} A \\ C  \end{array} \right]
  = r ( A ) + r(C) = n  \Leftrightarrow  R(A^*) \oplus R(C^*)
  = {\cal C}^n.$ }

\medskip

\noindent  {\bf Proof.}\, Note that $AA^{\dagger}$, $ A^{\dagger}A, $ $ BB^{\dagger}$, 
and $ C^{\dagger}C$ are all idempotent. Thus we can easily derive  by   (3.1) and (3.12) 
the four rank equalities in Parts (a)---(d).  The results in Parts (e)---(j) are direct consequences 
of Parts (a) and (b). \qquad  $ \Box$

\medskip

\noindent  {\bf Theorem 7.3.}\, {\em Let $ A \in {\cal C}^{ m \times n},
\, B \in {\cal C}^{ k \times m}$ be given. Then
$$
\displaylines{
\hspace*{1.5cm}
 r(  \, AA^{\dagger}B^{\dagger}B -
  B^{\dagger}BAA^{\dagger} \, )  = 2r[ \, A, \ B^* \, ] + 2r( BA ) - 2r ( A )
   - 2r( B ). \hfill (7.2)
\cr
In \ particular \hfill
\cr
\hspace*{0cm}
(AA^{\dagger})(B^{\dagger}B)
 = (B^{\dagger}B)(AA^{\dagger}) \Leftrightarrow r[ \, A, \ B^* \, ]
 = r ( A ) + r( B ) - r( BA ) \Leftrightarrow {\rm dim}[R(A)\cap R(B^*)]
 = r(BA). \hfill (7.3)
\cr}
$$}
{\bf Proof.}\, Note that $AA^{\dagger}$, $ A^{\dagger}A, $
$ BB^{\dagger}$, and $ B^{\dagger}B$ are Hermitian idempotent.  Thus we find by (3.29) that
\begin{eqnarray*}
r[ \, (AA^{\dagger})(B^{\dagger}B) -  (B^{\dagger}B) (AA^{\dagger}) \, ] & = &2 r[\, AA^{\dagger},\  B^{\dagger}B \, ] +
  2r[\, (AA^{\dagger})(B^{\dagger}B) \,]  -  2r(AA^{\dagger}) -
  2r(B^{\dagger}B) \\
& = &  2r[\, A, \ B^* \,] + 2r( BA ) -2r( A ) - 2r(B) , 
\end{eqnarray*}
as required for  (7.2). The result in (7.3) is a direct consequence
 of (7.2). \qquad  $ \Box$ 

\medskip

Replace $ B $ by $B^*$ in (7.2) to yield an alternative formula
$$ 
\displaylines{
\hspace*{1.5cm}
 r(\, AA^{\dagger}BB^{\dagger} -
  BB^{\dagger}AA^{\dagger} \, )  = 2r[ \, A, \ B \, ] + 2r( B^*A ) - 2r ( A )
   - 2r( B ). \hfill (7.4)
\cr}
$$  
Some interesting consequences can be derived from (7.2) and (7.4). For example,  let $B = I_m - A $ in (7.2). Then  we
get by (1.11)  
$$
\displaylines{
\hspace*{1.5cm}
 r[ \, AA^{\dagger}( I_m - A)^{\dagger}(I_m - A) -
  ( I_m - A)^{\dagger}( I_m - A)AA^{\dagger} \, ] \hfill
\cr
\hspace*{1.5cm}
= 2r[ \, A, \ I_m - A^* \, ] + 2r(A - A^2 ) - 2r( A ) - 2r(I_m - A ) \hfill
\cr
\hspace*{1.5cm}
 = 2r[ \, A, \ I_m - A^* \, ] - 2m \leq 0. \hfill
\cr}
$$
Because the rank of a matrix is nonnegative, the above inequality in fact implies 
 that $r[ \, A, \ I_m - A^* \, ]= m$ and $ AA^{\dagger}$ commutes with 
$ ( I_m - A)^{\dagger}(I_m - A)$ for any square matrix $ A $. Based on this result, one can easily see 
that  $ AA^{\dagger}$ also commutes with $ ( \lambda I_m - A)^{\dagger}( \lambda I_m - A)$ for any 
 $\lambda \neq 0$.  However, it is curious that  $ AA^{\dagger}$ does not commute with 
$ (I_m - A )(I_m - A)^{\dagger}$ in general. In fact,  we find by (7.4) that 
$$
\displaylines{
\hspace*{1.5cm}
 r[ \, AA^{\dagger}(I_m - A)( I_m - A)^{\dagger} -
  (I_m - A)( I_m - A)^{\dagger}AA^{\dagger} \, ] \hfill
\cr
\hspace*{1.5cm}
= 2r[ \, A, \ I_m - A \, ] + 2r(A - A^*A ) - 2r( A ) - 2r(I_m - A ) \hfill
\cr
\hspace*{1.5cm}
= 2m + 2r(A - A^*A ) - 2r( A ) - 2r(I_m - A ) \hfill
\cr
\hspace*{1.5cm}
= 2r(A - A^*A ) - 2r( A  - A^2 ). \hfill
\cr}
$$
Thus $ AA^{\dagger}$  commutes with $ (I_m - A )(I_m - A)^{\dagger}$ if and only if $ r(A - A^*A ) = r( A  - A^2 )$. 

Next replacing $ A $ and $ B $ by  $I_m + A $  and $I_m - A $ in (7.2), respectively,  we can get 
$$
\displaylines{
\hspace*{1.5cm}
 r[ \, ( I_m + A)( I_m + A)^{\dagger}( I_m - A)^{\dagger}(I_m - A) -
  ( I_m - A)^{\dagger}( I_m - A)( I_m + A)( I_m + A)^{\dagger} \, ] \hfill
\cr
\hspace*{1.5cm}
= 2r[ \, I_m +  A, \ I_m - A^* \, ] + 2r(I_m - A^2 ) - 2r( I_m +  A ) - 2r(I_m - A ) \hfill
\cr
\hspace*{1.5cm}
 = 2r[ \, I_m + A, \ I_m - A^* \, ] - 2m \leq 0. \hfill
\cr}
$$
This inequality implies that $r[ \,I_m +  A, \ I_m - A^* \, ]= m$ and $( I_m + A)( I_m + A)^{\dagger}$ commutes with 
$ ( I_m - A)^{\dagger}(I_m - A)$ for any square matrix $ A $. By (7.4) we also find that 
$$
\displaylines{
\hspace*{1.5cm}
 r[ \, ( I_m + A)( I_m + A)^{\dagger}( I_m - A)( I_m - A)^{\dagger} - 
 ( I_m - A)(I_m - A)^{\dagger}( I_m + A)( I_m + A)^{\dagger} \, ] \hfill
\cr
\hspace*{1.5cm}
= 2r[ \, I_m +  A, \ I_m - A \, ] + 2r[(I_m - A^*)( I_m + A) - 2r( I_m +  A ) - 2r(I_m - A ) \hfill
\cr
\hspace*{1.5cm}
 = 2m + 2r[(I_m - A^*)( I_m + A)] - 2r( I_m +  A ) - 2r(I_m - A ) \hfill
\cr
\hspace*{1.5cm}
= 2r[(I_m - A^*)( I_m + A)] - 2r( I_m - A^2). \hfill
\cr}
$$
Thus $(I_m + A)( I_m + A)^{\dagger}$ commutes with $ (I_m - A )(I_m - A)^{\dagger}$ if and only if 
$r[(I_m - A^*)( I_m + A)] = r( I_m - A^2) $.

A general result is that for any two polynomials $ p(\lambda)$ and  $ q(\lambda)$ without common roots, we have according to (7.2) and (1.17)  
$$
\displaylines{
\hspace*{1.5cm}
 r[ \, p(A)p^{\dagger}(A)q^{\dagger}(A)q(A) - q^{\dagger}(A)q(A)p(A)p^{\dagger}(A)\, ] \hfill
\cr
\hspace*{1.5cm}
= 2r[ \, p(A), \ q(A^*) \, ] + 2r[p(A)q(A)] - 2r[p(A)] - 2r[q(A)] \hfill
\cr
\hspace*{1.5cm}
 = 2r[ \,  p(A), \ q(A^*) \, ] - 2m \leq 0. \hfill
\cr}
$$ 
This implies that $ r[ \,  p(A), \ q(A^*) \, ] = m$ and 
$[p(A)p^{\dagger}(A)][q^{\dagger}(A)q(A)] = [q^{\dagger}(A)q(A)][p(A)p^{\dagger}(A)]$, that is, 
$p(A)p^{\dagger}(A)$ commutes with $q^{\dagger}(A)q(A)$. On the other hand,  the fact 
$r[ \,  p(A), \ q(A^*) \, ] = m$ can also alternatively be stated that 
the for any square matrix $ A $ and  any two polynomials $ p(\lambda)$ and  $ q(\lambda)$ without 
common roots the Hermitian matrix $ p(A)p(A^*) +  q(A^*)q(A)$ is always positive 
definite.  

\medskip

Observe that the Moore-Penrose inverses of $ [\, A, \ B \, ] $ and $ \left[ \begin{array}{c} A \\ C \end{array} \right]$ 
can be expressed as 
$$
\displaylines{
\hspace*{2cm}
[ \, A, \ B \, ]^{\dagger} = [ \, A, \ B \, ]^*\left( [ \, A, \ B \, ]
[ \, A, \ B \, ]^* \right)^{\dagger} = \left[ \begin{array}{c} A^*(\, AA^* + BB^* \,)^{\dagger} \\
  B^*(\, AA^* + BB^* \,)^{\dagger}  \end{array} \right], \hfill
\cr
\hspace*{2cm} \left[ \begin{array}{c} A \\ C \end{array} \right]^{\dagger} = 
\left( \left[ \begin{array}{c} A \\ C \end{array} \right]^* \left[ \begin{array}{c} A \\ 
 C \end{array} \right] \right)^{\dagger}\left[ \begin{array}{c} A \\ C \end{array} \right]^*  
= [\, (\, A^*A + C^*C \,)^{\dagger}A^*,  \  (\, A^*A + C^*C \,)^{\dagger}C^* \, ].\hfill
\cr}
$$
Based on the two expressions  we can find a series of rank equalities related to  $ [\, A, \ B \, ] $ and 
$ \left[ \begin{array}{c} A \\ C \end{array} \right]$  and their consequences.  
 
\medskip

\noindent  {\bf Theorem 7.4.}\, {\em Let $ A \in {\cal C}^{ m \times n},
\, B \in {\cal C}^{ k \times m}$ be given. Then

{\rm (a)} \ $r[ \, AA^*( \, AA^* +  BB^* \, )^{\dagger}A - A \, ]  = r ( A )
   + r( B) - r[ \, A, \ B \, ].$

{\rm (b)} \ $r[ \, AA^{\dagger}( \, AA^{\dagger} +  BB^{\dagger} \, )^{\dagger}A - A \, ]  = r ( A )
   + r( B) - r[ \, A, \ B \, ].$

{\rm (c)} \ $r[ \, A( \, A^*A +  C^*C \, )^{\dagger}A^*A - A \, ]  = r ( A )
   + r( C) - r \left[ \begin{array}{c} A \\  C  \end{array} \right].$

{\rm (d)}  \ $r[ \, A( \, A^{\dagger}A +  C^{\dagger}C \, )^{\dagger}A^{\dagger}A - A \, ]  = r ( A )
   + r( C) - r \left[ \begin{array}{c} A \\  C  \end{array} \right].$

{\rm (e)} \ $r[ \, A^*( \, AA^* +  BB^* \, )^{\dagger}B \, ]  = r ( A )
   + r( B) - r[ \, A, \ B \, ].$

{\rm (f)}  \ $r[ \, A^{\dagger}( \, AA^{\dagger} +  BB^{\dagger} \, )^{\dagger}B \, ]  = r ( A )
   + r( B) - r[ \, A, \ B \, ].$

{\rm (g)} \ $r[ \, A( \, A^*A +  C^*C \, )^{\dagger}C^* \, ]  = r ( A )
   + r( C) - r \left[ \begin{array}{c} A \\  C  \end{array} \right].$

{\rm (h)} \  $r[ \, A( \, A^{\dagger}A +  C^{\dagger}C \, )^{\dagger}C^{\dagger} \, ]  = r ( A )
   + r( C) - r \left[ \begin{array}{c} A \\  C  \end{array} \right].$

{\rm (i)} The following five statements are equivalent$:$

\hspace*{ 0.5cm} {\rm (1)}  $AA^*( \, AA^* +  BB^* \,)^{\dagger}A = A.$

\hspace*{ 0.5cm} {\rm (2)}  $AA^{\dagger}( \, AA^{\dagger} +  BB^{\dagger} \,)^{\dagger}A = A.$

\hspace*{ 0.5cm} {\rm (3)}  $A^*( \, AA^* +  BB^* \, )^{\dagger}B =0. $

\hspace*{ 0.5cm} {\rm (4)} $ A^{\dagger}( \, AA^{\dagger} +  BB^{\dagger} \, )^{\dagger}B = 0.$ 

\hspace*{ 0.5cm} {\rm (5)} $ r[ \, A, \ B \, ] = r ( A )+ r( B), \  i.e., \   R(A) \cap R(B) = \{0\}.$

{\rm (j)}\,The following five statements are equivalent$:$

\hspace*{0.5cm} {\rm (1)}  $A( \, A^*A +  C^*C \,)^{\dagger}A^*A = A.$

\hspace*{ 0.5cm} {\rm (2)}  $A( \, A^{\dagger}A +  C^{\dagger}C \,)^{\dagger}A^{\dagger}A = A.$

\hspace*{0.5cm} {\rm (3)}  $ A( \, A^*A +  C^*C \, )^{\dagger}C^* =0. $

\hspace*{ 0.5cm} {\rm (4)}  $ A( \, A^{\dagger}A +  C^{\dagger}C \, )^{\dagger}C^{\dagger} = 0.$ 

\hspace*{0.5cm} {\rm (5)} $ r \left[ \begin{array}{c} A \\  C  \end{array} \right] = r ( A )
   + r( C) \ i.e., \  R(A^*) \cap R(C^*) = \{0\}. $
}

\medskip

\noindent  {\bf Proof.}\, We only show Parts (a) and (b).  Note that $ R(A) \subseteq R( \, AA^* + BB^*\,).$ Thus we find by (1.7) that 
\begin{eqnarray*}
r[ \, AA^*( \, AA^* +  BB^* \, )^{\dagger}A - A \, ]
 &=&  r \left[ \begin{array}{cc} AA^* +  BB^*  & A \\ AA^* & A \end{array}
 \right] - r( \, AA^* +  BB^* \, ) \\
&=&  r \left[ \begin{array}{cc} BB^*  & 0 \\ 0 & A \end{array}
 \right] - r[ \, A, \  B \, ] = r ( A ) + r( B) - r[ \, A, \ B \, ],
\end{eqnarray*}
as required for Part (a).  Similarly
\begin{eqnarray*}
r[ \, AA^{\dagger}( \, AA^{\dagger} +  BB^{\dagger} \, )^{\dagger}A - A \, ]
 &=&  r \left[ \begin{array}{cc} AA^{\dagger} +  BB^{\dagger}  & A \\ AA^{\dagger} & A \end{array}
 \right] - r( \,  AA^{\dagger} +  BB^{\dagger}  \, ) \\
&=&  r \left[ \begin{array}{cc}  BB^{\dagger}   & 0 \\ 0 & A \end{array}
 \right] - r[ \, A, \  B \, ] = r ( A ) + r( B) - r[ \, A, \ B \, ],
\end{eqnarray*}
as required for Part (b). \qquad  $ \Box$ 
 
\medskip

\noindent  {\bf Theorem 7.5.}\, {\em Let $ A \in {\cal C}^{m \times n},
 \, B \in {\cal C}^{ m \times k}$ and $ C \in {\cal C}^{ l \times n} $ be
 given. Then

{\rm (a)}\, The rank of $ BB^{\dagger}A -  AC^{\dagger}C$ satisfies 
$$\displaylines{
\hspace*{1.5cm}
r(  \, BB^{\dagger}A -  AC^{\dagger}C \, ) = r \left[ \begin{array}{c} B^*A
\\  C  \end{array} \right] +  r[ \, AC^*, \ B \, ] - r( B ) - r( C ). \hfill
\cr}
$$

{\rm (b)} \ $BB^{\dagger}A=  AC^{\dagger}C \Leftrightarrow 
r \left[ \begin{array}{c} B^*A \\  C  \end{array} \right] =  r( C )  \  and  \  
r[ \, AC^*, \ B \, ] = r( B ) \Leftrightarrow  R( AC^* ) \subseteq R( B ) \ and
  \  R( A^*B ) \subseteq R( C^* ).$ 

{\rm (c)} \ $ BB^{\dagger}A -  AC^{\dagger}C $ is nonsingular
  $\Leftrightarrow  r \left[ \begin{array}{c} B^*A \\  C  \end{array}
 \right] = r[ \, AC^*, \ B \, ] = r( B ) + r( C ) = m.$ }

\medskip

 \noindent {\bf Proof.} \ Follows from (4.1) by noticing that both $ BB^{\dagger}$ and $ C^{\dagger}C$ are idempotent. \qquad $ \Box$

\medskip

It is well known that the matrix equation $ BXC = A$ is solvable if and only
if $BB^{\dagger}AC^{\dagger}C  = A$. This leads us to 
consider the rank of $ A - BB^{\dagger}AC^{\dagger}C$.  

\medskip

\noindent {\bf Theorem 7.6.}\, {\em Let $ A \in {\cal C}^{ m \times n}, 
 \, B \in {\cal C}^{ m \times k}$ and $C \in {\cal C}^{ l \times n} $ be
 given. Then
$$\displaylines{
\hspace*{1.5cm}
r(  \, A - BB^{\dagger}AC^{\dagger}C \, )
 = r \left[ \begin{array}{ccc} A & AC^*  & B \\
  B^*A & 0 & 0 \\ C & 0 & 0    \end{array} \right] - r( B ) - r( C ),
  \hfill (7.5)
\cr
\hspace*{0cm}
and \hfill
\cr
\hspace*{1.5cm}
r(  \, 2A - BB^{\dagger}A - AC^{\dagger}C \, )
 = r \left[ \begin{array}{ccc} A & AC^*  & B \\
  B^*A & 0 & 0 \\ C & 0 & 0    \end{array} \right] - r( B ) - r( C ).
 \hfill (7.6)
\cr
\hspace*{0cm}
In \ particular, \hfill
\cr
\hspace*{1.5cm}
 BB^{\dagger}AC^{\dagger}C = A \Leftrightarrow
   BB^{\dagger}A + AC^{\dagger}C = 2A  \Leftrightarrow
   R( A) \subseteq R( B ) \ \ and  \ \ R(A^*) \subseteq R(C^*).
   \hfill (7.7)
\cr}
$$ } 
{\bf Proof.}\, Applying (2.8) and the rank cancellation law (1.8) 
to $ A - BB^{\dagger}AC^{\dagger}C $ produces
 \begin{eqnarray*}
 \lefteqn {r( \, A - BB^{\dagger}AC^{\dagger}C  \, ) } \\
  & = & r \left[ \begin{array}{ccr} B^*AC^*  & B^*BB^*  & 0 \\
   C^*CC^*  &  0 & C^*C  \\ 0 &  BB^* & -A  \end{array} \right] - r(B) -
   r( C )  \\
 & = &  r \left[ \begin{array}{ccr} B^*AC^*  & B^*B  & 0 \\ CC^*  &  0 & C
  \\ 0 &  B & -A  \end{array} \right] - r(B) - r( C )
  = r \left[ \begin{array}{ccr} 0 & 0  & B^*A \\ 0 &  0 & C  \\ AC^* &  B &
   -A  \end{array} \right] - r(B) - r( C ),
 \end{eqnarray*}
as required for (7.5). In the same way we can show  (7.6). The result
  in  (7.7) is well known.  \qquad  $ \Box$ 

\medskip

\noindent  {\bf Theorem 7.7.}\, {\em Let $ A \in {\cal C}^{ m \times n},  \,
 B \in {\cal C}^{ m \times k}$ and $ C \in {\cal C}^{ l \times n} $ be given$,$ and 
let $M = \left[ \begin{array}{cc} A & B \\  C & 0 \end{array} \right].$
 Then
 $$
  r(  \, A - BB^{\dagger}A - AC^{\dagger}C \, ) = r(M) + r( CA^*B ) - r( B ) - 
r( C ), \eqno (7.8)
 $$
 that is$,$  the block matrix $M$ satisfies the rank equality   
 $$ 
 r(M) = r( B ) + r( C ) - r( CA^*B ) + r( \, A - BB^{\dagger}A - AC^{\dagger}C \,).
   \eqno (7.9)
 $$ }
{\bf Proof.}\, Applying (2.2) and (1.8) to
 $ A - BB^{\dagger}A - AC^{\dagger}C$ yields
 \begin{eqnarray*}
r( \, A - BB^{\dagger}A - AC^{\dagger}C  \, ) & = & r \left[ \begin{array}{ccc} B^*BB^*  & 0  & B^*A \\
   0 &  C^*CC^*  & C^*C  \\  BB^*  &  AC^* & A  \end{array} \right] -
   r(B) - r( C )  \\
 & = &  r \left[ \begin{array}{ccc}  B^*B  & 0  & B^*A \\  0 &  CC^*  & C  \\
   B  &  AC^* & A  \end{array} \right] - r(B) - r( C )  \\
 & = & r \left[ \begin{array}{ccc} 0  &  -B^*AC^*   & 0 \\ 0  &  0 & C  \\
   B &  0 & A  \end{array} \right] - r(B) - r( C ) \\
&= & r \left[ \begin{array}{cc} A & B \\  C & 0 \end{array}
 \right] + r( B^*AC^* ) - r( B ) - r( C ),
  \end{eqnarray*}
as required for  (7.8). \qquad  $ \Box$

\medskip

\noindent  {\bf Theorem 7.8.}\, {\em Let $ A \in {\cal C}^{ m \times n},  \,
 B \in {\cal C}^{ m \times k}$ and $ C \in {\cal C}^{ l \times n} $ be given$,$ and 
let $M = \left[ \begin{array}{cc} A & B \\  C & 0 \end{array} \right].$ Then

{\rm (a)}\,  The rank of $ A - A(E_BAF_C)^{\dagger}A$ satisfies  $$\displaylines{
\hspace*{2cm}
 r[  \, A - A(E_BAF_C)^{\dagger}A \, ] = r(A) + r(B) + r(C) - r(M), \hfill (7.10)
 \cr
\hspace*{0cm}
that \ is, \hfill
\cr
\hspace*{2cm}
r\left[ \begin{array}{cc} A & B \\  C & 0 \end{array} \right]
 =  r(A) + r(B) + r(C) - r[ \, A - A(E_BAF_C)^{\dagger}A \, ], \hfill (7.11)
\cr}
$$
where $ E_B = I - BB^{\dagger}$ and  $ F_C = I - C^{\dagger}C$.

{\rm (b)}\, In particular$,$
$$\displaylines{
\hspace*{2cm}
(E_BAF_C)^{\dagger} \in \{ A^-\} \hfill (7.12)
\cr
\hspace*{0cm}
 holds \ if \ and \ only \ if \hfill
\cr
\hspace*{1cm}
r(M) = r(A) + r(B) + r(C), \ \ i.e., \ \  R(A) \cap R(B) = \{ 0\}  \ and \  R(A^*) \cap R(C^*)
 = \{ 0\}. \hfill (7.13)
\cr}
$$

{\rm (c)}  \ $r[  \, A - A(E_BA)^{\dagger}A \, ] = r(A) + r(B) - r[\, A , \ B \, ].$ 

{\rm (d)} \  $r[  \, A - A(AF_C)^{\dagger}A \, ] = r(A) + r(C) - r\left[ \begin{array}{c} A \\ C \end{array} \right].$ 

{\rm (e)}  \ $A(E_BA)^{\dagger}A  = A \Leftrightarrow r[\, A , \ B \, ] = r(A) + r(B), \ i.e., \ R(A) \cap R(B) = \{ 0\}.$ 

{\rm (f)}  \ $ A(AF_C)^{\dagger}A = A  \Leftrightarrow  r\left[ \begin{array}{c} A \\ C \end{array}\right]
= r(A) + r(C), \ i.e., \ R(A^*) \cap R(C^*) = \{ 0\}.$ }

\medskip

\noindent {\bf Proof.}\, Let $ N = E_BAF_C $. Then it is easy to verify that 
$ N^*NN^* = N^*AN^*.$ In  that case, applying  (2.1), and then (1.2) and (1.3) to 
$ A -  A(E_BAF_C)^{\dagger}A$ yields        
\begin{eqnarray*}
 r[ \, A -  A(E_BAF_C)^{\dagger}A \, ] = r[ \, A -  AN^{\dagger}A \, ]
 & = & r \left[ \begin{array}{cc} N^*NN^*  & N^*A \\ AN^*  & A
 \end{array} \right] -  r(M) \\
& = & r \left[ \begin{array}{cc} N^*NN^* - N^*AN^*  & 0 \\ 0  & A
 \end{array} \right] -  r(N) \\
& = & r \left[ \begin{array}{cc} 0  & 0 \\ 0  & A
 \end{array} \right] -  r(N) \\
 & = & r(A) - r(N)  = r(A) + r(B) + r(C) - r(M),
\end{eqnarray*}
as required for (7.10).  The equivalence of   (7.12) and (7.13) follows immediately from (7.11). \qquad  $ \Box$ 

\medskip

It is known that for any $ B $ and $ C $, the matrix $(E_BAF_C)^{\dagger}$ is always an outer inverse 
of $ A $ (Greville \cite{Gr2}). Thus the rank formula (7.11) can also be derived from  (5.6). 
  
In the remainder of this chapter, we establish various  rank equalities related
to  ranks of  Moore-Penrose inverses of block  matrices, and  then present their consequences.

\medskip

\noindent {\bf Theorem 7.9.}\, {\em Let $ A \in {\cal C}^{ m \times n},
 \, B \in {\cal C}^{ m \times k}$ and $ C \in {\cal C}^{ l \times n} $ be
 given. Then

 {\rm (a)} \  $ r \left(  \, [\, A, \ B \, ]^{\dagger} -
  \left[ \begin{array}{c} A^{\dagger} \\ B^{\dagger}  \end{array} \right] \,
  \right) = r[ \, AA^*B, \ BB^*A \,].$

 {\rm (b)} \ $ r \left( \,  \left[ \begin{array}{c} A \\ C  \end{array}
 \right]^{\dagger} - [\, A^{\dagger}, \ C^{\dagger} \, ] \,  \right)
 = r\left[ \begin{array}{c} AC^*C \\ CA^*A  \end{array} \right]. $ 

 {\rm (c)} \ $ r \left(  \, [\, A, \ B \, ]^{\dagger}[\, A, \ B \, ] -
  \left[ \begin{array}{c} A^{\dagger} \\ B^{\dagger}  \end{array} \right]
[\, A, \ B \, ] \,  \right) = r[ \, AA^*B, \ BB^*A \,].$

 {\rm (d)} \ $ r \left( \, \left[ \begin{array}{c} A \\ C  \end{array}
 \right] \left[ \begin{array}{c} A \\ C  \end{array}
 \right]^{\dagger} - \left[ \begin{array}{c} A \\ C  \end{array}
 \right][\, A^{\dagger}, \ C^{\dagger} \, ] \,  \right)
 = r\left[ \begin{array}{c} AC^*C \\ CA^*A  \end{array} \right]. $ \\
In  particular$,$

{\rm (e)} \  $ [\, A, \ B \, ]^{\dagger}
 = \left[ \begin{array}{c} A^{\dagger} \\ B^{\dagger}  \end{array} \right]  
 \Leftrightarrow [\, A, \ B \, ]^{\dagger}[\, A, \ B \, ] = 
\left[ \begin{array}{c} A^{\dagger} \\ B^{\dagger}  \end{array} \right]
[\, A, \ B \, ] \Leftrightarrow  A^*B = 0.$

 {\rm (f)} \  $ \left[ \begin{array}{c} A \\ C  \end{array} \right]^{\dagger}
 = [\, A^{\dagger}, \ C^{\dagger} \, ] \Leftrightarrow 
\left[ \begin{array}{c} A \\ C  \end{array} \right] \left[ \begin{array}{c} A \\ C  
\end{array} \right]^{\dagger} = \left[ \begin{array}{c} A \\ C  \end{array}
 \right][\, A^{\dagger}, \ C^{\dagger} \, ] \Leftrightarrow CA^* = 0.$ } 

\medskip

\noindent   {\bf Proof.}\,  Let $ M = [ \, A, \ B \, ].$ Then  it follows 
by (2.7) that
 \begin{eqnarray*}
\lefteqn {r\left(  \,[\, A, \ B \, ]^{\dagger} -
 \left[ \begin{array}{c} A^{\dagger} \\ B^{\dagger}  \end{array} \right] \,
 \right)  } \\
 & = &  r\left(  \, [\, A, \ B \, ]^{\dagger} -
 \left[ \begin{array}{c} I \\ 0 \end{array} \right] A^{\dagger} -
 \left[ \begin{array}{c} 0 \\ I \end{array} \right] B^{\dagger} \, \right) \\
& = & r \left[ \begin{array}{cccc} 
 -M^*MM^*  & 0 & 0 & M^*   \\ 0 &  A^*AA^* &  0  &  A^*  \\ 0  & 0  &
  B^*BB^*  & B^* \\
  \left[ \begin{array}{c} A^* \\ B^* \end{array} \right] &
  \left[\begin{array}{c} A^* \\ 0 \end{array} \right] &
   \left[ \begin{array}{c} 0 \\ B^* \end{array} \right] &
   \left[ \begin{array}{c} 0 \\ 0 \end{array}
 \right] \end{array} \right] - r(M) - r( A ) - r( B )  \\ 
 & = & r \left[ \begin{array}{cccc} -M^*MM^*  & 0  & 0 & M^*  \\
 - A^*AA^* & 0  & 0 &  A^*  \\
  - B^*BB^*  & 0  & 0  & B^* \\ 0 & A^* & 0  & 0  \\ 0 & 0 & B^* & 0 
  \end{array} \right] - r(M) - r( A ) - r( B )  \\ 
& = & r \left[ \begin{array}{cc} M^*MM^* & M^*  \\ A^*AA^* &  A^*  \\
 B^*BB^*  & B^*  \end{array} \right] - r(M) \\ 
& = & r \left[ \begin{array}{ccc} MM^*M  & AA^*A  & BB^*B  \\ M & A & B
 \end{array} \right] - r(M) \\
 & = & r \left[ \begin{array}{ccc} 0  & AA^*A - MM^*A  & BB^*B -MM^*B  \\
  M & 0  & 0  \end{array} \right] - r(M) \\
 & = & r[\,  AA^*A - MM^*A, \ BB^*B -MM^*B \,]
 =  r[ \, AA^*B, \ BB^*A \, ],
 \end{eqnarray*}
as required in Part (a). Similarly, we can show Parts (b), (c) and (d).  The results in
 Parts (e) and (f) follow immediately from Parts (a)---(d). \qquad  $\Box$ 

\medskip

A general result is  given below, the proof is omitted. 

\medskip

\noindent {\bf Theorem 7.10.}\, {\em Let $ A =  [\, A_1, \, A_2, \, \cdots,
 \, A_k \, ] \in {\cal C}^{ m \times n}$  be given$,$ and denote
 $ M = \left[ \begin{array}{c} A_1^{\dagger} \\ \vdots \\A_k^{\dagger}
 \end{array} \right]$. Then 
$$\displaylines{
\hspace*{2cm}
 r( \,  A^{\dagger} - M \, ) = r( \,  A^{\dagger}A - MA \, ) = r[ \, N_1N_1^*A_1, \,
 N_2N_2^*A_2, \,  \cdots, \, N_kN_k^*A_k  \, ], \hfill (7.14)
\cr}
$$
where $ N_i = [\, A_1,  \, \cdots, \, A_{i-1}, \, A_{i+1}, \, \cdots, \, A_k \,], 
\ \ i = 1, \, 2, \, \cdots, \, k.$ 
In particular$,$
$$\displaylines{
\hspace*{2cm} 
A^{\dagger} = M \Leftrightarrow A^{\dagger}A = MA \Leftrightarrow A_iA^*_j= 0 
\ \ \  for \ all \  i \neq j. \hfill (7.15)
\cr}
$$ } 
\hspace*{0.4cm} Let $[\, A, \ B\,]^{\dagger} = \left[ \begin{array}{c} G_1 \\ G_2 \end{array} \right]$ and 
$ \left[ \begin{array}{c} A \\ C  \end{array} \right]^{\dagger} = [\, H_1, \ H_2\,].$ 
We next consider the relationships between $ A $ and $G_1$, $A$ and $H_1$.       

\medskip

\noindent {\bf Theorem 7.11.}\, {\em Let $ A \in {\cal C}^{ m \times n},  \, B
 \in {\cal C}^{ m \times k}$ and $ C \in {\cal C}^{ l \times n} $ be given. Then

{\rm (a)} \ $ r( \, A^{\dagger} - [ \, I_n, \ 0 \,][\, A, \ B \,]^{\dagger} \, )  
=  r( \, B^{\dagger} - [ \, 0, \  I_k \,][\, A, \ B \,]^{\dagger} \, )  = r (A^*B).$ 

{\rm (b)}  \ $ r \left(  A^{\dagger} - \left[ \begin{array}{c} A \\ C  \end{array}
 \right]^{\dagger}\left[ \begin{array}{c} I_m \\ 0 \end{array}
 \right] \right) 
= r \left(  C^{\dagger} - \left[ \begin{array}{c} A \\ C  \end{array} \right]^{\dagger} 
\left[ \begin{array}{c} 0  \\ I_l \end{array}
 \right] \right)
=  r( CA^*).$

{\rm (c)}  \ $ r( \, AA^{\dagger} - [ \, A, \ 0 \,][\, A, \ B \,]^{\dagger} \, )  
=  r( \, BB^{\dagger} - [ \, 0, \  B \,][\, A, \ B \,]^{\dagger} \, ) = r (A^*B).$ 

{\rm (d)} \ $ r \left(  A^{\dagger}A - \left[ \begin{array}{c} A \\ C  \end{array}
 \right]^{\dagger}\left[ \begin{array}{c} A \\ 0  \end{array}
 \right]   \right) 
= r \left(  C^{\dagger}C - \left[ \begin{array}{c} A \\ C  \end{array} 
 \right]^{\dagger}\left[ \begin{array}{c} 0  \\ C  \end{array}
 \right] \right) =  r( CA^*).$ \\
In  particular$,$

{\rm (e)} \ $ [\, I_m, \ 0 \, ][\, A, \ B \, ]^{\dagger} = A^{\dagger} \Leftrightarrow 
[ \, 0, \  I_k \,][\, A, \ B \,]^{\dagger} = B^{\dagger} \Leftrightarrow 
 [\, A, \ 0 \, ][\, A, \ B \, ]^{\dagger} = AA^{\dagger} \Leftrightarrow 
[ \, 0, \  B \,][\, A, \ B \,]^{\dagger} = B^{\dagger}B \Leftrightarrow A^*B = 0.$ 

 {\rm (f)} \  $ \left[ \begin{array}{c} A \\ C  \end{array} \right]^{\dagger}\left[ \begin{array}{c} I_m \\ 0 \end{array}
 \right]  = A^{\dagger} \Leftrightarrow 
\left[ \begin{array}{c} A \\ C  \end{array} \right]^{\dagger}\left[ \begin{array}{c} 0 \\ I_l  \end{array} \right]  
= C^{\dagger} \Leftrightarrow \left[ \begin{array}{c} A \\ C  \end{array}
 \right]^{\dagger}\left[ \begin{array}{c} A \\ 0  \end{array} \right]  = A^{\dagger}A \Leftrightarrow 
\left[ \begin{array}{c} A \\ C  \end{array} \right]^{\dagger}\left[ \begin{array}{c} 0 
\\ C  \end{array} \right]  = C^{\dagger}C 
\Leftrightarrow CA^* =0.$ }

\medskip

\noindent  {\bf Proof.}\,  We only prove Parts (a) and (c). Let $ M = [\, A, \ B \, ].$ Then it follows
 by  (2.7), (1.8) and block elementary operations of matrices  that
\begin{eqnarray*}
 r( \, A^{\dagger} - [\, I_m, \ 0 \, ]M^{\dagger} \,) 
& = & r \left[ \begin{array}{ccc} 
 -A^*AA^*  & 0 & A^*   \\ 0 &  M^*MM^* &  M^*  
\\ A^* & \left[\, I_m, \ 0 \, \right] M^* & 0 \end{array} \right]
   - r(A) - r( M )   \\
 & = & r \left[ \begin{array}{cccc} A^*AA^*  & -A^*MM^* & A^*  \\ 0 & 0 & M^*  \\ A^* & A^* & 0 \end{array} \right]
  - r( A ) - r(M)  \\
 & = & r \left[ \begin{array}{cccc} 0 & -A^*BB^* & A^*  \\ 0 & 0 & M^*  \\ A^* & 0 & 0 \end{array} \right]
  - r( A ) - r(M) = r(A^*BB^*) = r(A^*B),
 \end{eqnarray*}
establishing the first equality in Part (a). Similarly
\begin{eqnarray*}
r( \, AA^{\dagger} - [\, A, \ 0 \, ]M^{\dagger} \,)& = & r \left[ \begin{array}{ccc} 
 -A^*AA^*  & 0 & A^*   \\ 0 &  M^*MM^* &  M^*  
\\ AA^* & \left[\, A, \ 0 \, \right] M^* & 0 \end{array} \right]
   - r(A) - r(M) \\
 & = & r \left[ \begin{array}{ccc} A^*A  & -A^*MM^* & A^*  \\ 0 & 0 & M^*  \\ A & AA^* & 0 \end{array} \right]
  - r( A ) - r(M)  \\
 & = & r \left[ \begin{array}{ccc} 0 & -A^*AA^* + A^*MM^*  & A^*  \\ 0 & 0 & M^*  \\ A^* & 0 & 0 \end{array} \right]
  - r( A ) - r(M) \\
& = & r \left[ \begin{array}{ccc} 0 & A^*BB^*  & A^*  \\ 0 & 0 & A^* \\ 0 & 0 & B^*  \\ A^* & 0 & 0 \end{array} \right]
  - r( A ) - r(M) = r(A^*BB^*) = r(A^*B),
 \end{eqnarray*}
establishing the first equality in Part (c). \qquad   $ \Box$ 

\medskip

\noindent {\bf Theorem 7.12.}\, {\em Let $ A \in {\cal C}^{ m \times n},  \, B
 \in {\cal C}^{ m \times k}$ and $ C \in {\cal C}^{ l \times n} $ be given.
 Then

 {\rm (a)} \ $ r( \, [\, A, \ B \, ][\, A, \ B \, ]^{\dagger} -
 ( \, AA^{\dagger} + BB^{\dagger} \, ) \, ) =
  r[ \, A, \ B \, ] + 2r (A^*B) - r ( A ) - r( B ).$ 

 {\rm (b)} \ $ r \left( \left[ \begin{array}{c} A \\ C  \end{array}
 \right]^{\dagger}
 \left[ \begin{array}{c} A \\ C  \end{array} \right] -
 ( \, A^{\dagger}A + C^{\dagger}C \, )  \right)
  = r\left[ \begin{array}{c} A \\ C  \end{array} \right]  +
  2r( CA^*)- r ( A ) - r( C ).$ \\
In particular$,$

 {\rm (c)} \ $ [\, A, \ B \, ][\, A, \ B \, ]^{\dagger}
 = AA^{\dagger} + BB^{\dagger} \Leftrightarrow  A^*B = 0
 \Leftrightarrow  [\, A, \ B \, ]^{\dagger}
 = \left[ \begin{array}{c} A^{\dagger} \\ B^{\dagger}  \end{array} \right].$

 {\rm (d)}  \ $ \left[ \begin{array}{c} A \\ C  \end{array} \right]^{\dagger}
 \left[ \begin{array}{c} A \\ C  \end{array} \right]
 = A^{\dagger}A + C^{\dagger}C  \Leftrightarrow 
  CA^* = 0  \Leftrightarrow \left[ \begin{array}{c} A \\ C  \end{array} \right]^{\dagger}
   = [\, A^{\dagger}, \ C^{\dagger} \,]. $ } 

\medskip

\noindent  {\bf Proof.}\,  Let $ M = [\, A, \ B \, ].$ Then it follows
 by (2.7), (1.8) and block elementary operations of matrices  that
\begin{eqnarray*}
 \lefteqn{ r( \, MM^{\dagger} -  AA^{\dagger} -  BB^{\dagger} \,) } \\
  & = & r \left[ \begin{array}{cccc} 
 -M^*MM^*  & 0 & 0 & M^*   \\ 0 &  A^*AA^* &  0  &  A^*  \\ 0  & 0  &
  B^*BB^*  & B^* \\ MM^* & AA^* & BB^* & 0 \end{array} \right]
   - r(M) - r( A ) - r( B )  \\
 & = & r \left[ \begin{array}{cccc} -M^*M  & 0  & 0 & M^*  \\ 0 &  A^*A &
  0 &  A^*  \\ 0  & 0  & B^*B  & B^* \\  M & A & B & 0  \end{array} \right]
   - r(M) - r( A ) - r( B )  \\
& = & r \left[ \begin{array}{ccccc} 0 & 0  & 0 & 0  & A^*  \\
0 & 0 & 0 & 0  & B^*  \\ A^*A & A^*B  & A^*A & 0 &  A^*  \\ B^*A & B^*B & 0
 & B^*B  & B^* \\ A & B & A & B & 0  \end{array} \right]
   - r(M) - r( A ) - r( B )  \\
& = & r \left[ \begin{array}{ccccc} 0 & 0  & 0 & 0  & A^*  \\
0 & 0 & 0 & 0  & B^*  \\ 0 & 0 & 0 & -A^*B & 0  \\  0 & 0
 & -B^*A  &  0 &  0  \\ A & B & 0 & 0 & 0  \end{array} \right]
   - r(M) - r( A ) - r( B )  \\
 & = &  r(M) + 2r (A^*B) - r (A) - r(B),
 \end{eqnarray*}
as required in Part (a). In the same way, we can show Part (b). We know from
 Part (a) that
 $$
  MM^{\dagger} = AA^{\dagger} +BB^{\dagger} \Leftrightarrow 
 r[ \, A, \ B \, ] = r(A) + r(B) - 
 2r(A^*B). \eqno (7.16) 
 $$ 
 On the other hand, observe from  (1.2) that
\begin{eqnarray*}  
 r[ \, A, \ B \, ] &= &r(A) + r( \, B - AA^{\dagger}B \,) \\
 &\geq & r(A) + r(B ) - r( AA^{\dagger}B) \\
& = & r(A) +r(B) - r(A^*B) \\
 &\geq & r(A) + r(B ) - 2r(A^*B).
\end{eqnarray*}
Thus  (7.16) is also equivalent to $A^*B = 0$. In the similar manner,
 we can show Part (d).  \qquad   $ \Box$ 

\medskip

 A general result is given below, the proof is omitted.  

\medskip

\noindent  {\bf Corollary 7.13.}\, {\em Let $ A = [\, A_1, \, A_2, \, \cdots, \, A_k \, ]
 \in {\cal C}^{ m \times n}$ be given. Then
$$ 
 r[ AA^{\dagger} -  (\, A_1A_1^{\dagger} + \cdots +
 A_kA_k^{\dagger} \,) ]
 = r\left[ \begin{array}{cccc} 0 & A_1^*A_2 & \cdots & A_1^*A_k \\ A_2^*A_1
 & 0 & \cdots &  A_2^*A_k  \\ \vdots  & \vdots  &
 \ddots & \vdots \\  A_k^*A_1 & A_k^*A_2 & \cdots & 0 \end{array} \right]
 + r( A ) - r ( A_1 ) - \cdots - r(A_k). \eqno (7.17)
$$ 
In particular$,$
$$ 
 AA^{\dagger} = A_1A_1^{\dagger} + \cdots + A_kA_k^{\dagger}
 \Leftrightarrow
 A_i^*A_j  = 0, \ \ for \ all \ i \neq j. \eqno (7.18)
$$ } 

\noindent {\bf Theorem 7.14.}\, {\em Let $ A \in {\cal C}^{ m \times n},
  \ B \in {\cal C}^{ m \times k}$ and $ C \in {\cal C}^{ l \times n} $ be
  given. Then

 {\rm (a)}\  $ r \left(  \, [\, A, \ B \, ]^{\dagger} -
  \left[ \begin{array}{c} (E_BA)^{\dagger} \\
 (E_AB)^{\dagger}  \end{array} \right] \, \right)
 = r(A) + r(B) - r[ \, A, \ B \,].$

 {\rm (b)} \ $ r \left( \,  \left[ \begin{array}{c} A \\ C  \end{array}
 \right]^{\dagger} -[\, (AF_C)^{\dagger}, \ (CF_A)^{\dagger} \, ] \, \right)
 =  r(A) + r(C) - r\left[ \begin{array}{c} A \\ C  \end{array} \right]. $ \\
In  particular$,$

 {\rm (c)} \ $  [\, A, \ B \, ]^{\dagger} =
  \left[ \begin{array}{c} (E_BA)^{\dagger} \\ (E_AB)^{\dagger}  \end{array}
\right] \Leftrightarrow r[ \, A, \ B \,] = r(A) + r(B),  \ i.e., \  
R(A) \cap R(B) = \{0\}.$

 {\rm (d)} \  $  \left[ \begin{array}{c} A \\ C  \end{array} \right]^{\dagger}
  = [\, (AF_C)^{\dagger}, \ (CF_A)^{\dagger} \, ]
   \Leftrightarrow  r\left[ \begin{array}{c} A \\ C  \end{array} \right]
 =  r(A) + r(C), \ i.e., \   R(A^*) \cap R(C^*) = \{0\}.$ } 

\medskip

\noindent  {\bf Proof.}\,  Let $ M = [ \, A, \ B \, ].$ Then it follows by  (2.7) 
and (1.8) and block elementary operations of matrices  that
 \begin{eqnarray*}
 \lefteqn{ r\left(  \, [\, A, \ B \, ]^{\dagger} - \left[ \begin{array}{c}
 (E_BA)^{\dagger} \\(E_AB)^{\dagger}  \end{array} \right] \, \right) } \\
 & = &  r\left(  \, [\, A, \ B \, ]^{\dagger} -
 \left[ \begin{array}{c}  I  \\ 0  \end{array} \right] (E_BA)^{\dagger}
- \left[ \begin{array}{c} 0  \\ I  \end{array} \right] (E_AB)^{\dagger}
 \, \right) \\
 & = & r \left[ \begin{array}{cccc} 
 -M^*MM^*  & 0 & 0 & M^*   \\ 0 &  (E_BA)^*(E_BA)(E_BA)^* &  0  &  (E_BA)^*
  \\ 0  & 0 &
 (E_AB)^*(E_AB)(E_AB)^*  & (E_AB)^*  \\ \left[ \begin{array}{c} A^* \\
  B^* \end{array} \right] & \left[ \begin{array}{c} (E_BA)^* \\ 0 \end{array}
  \right] & \left[ \begin{array}{c} 0 \\ (E_AB)^* \end{array} \right] &
   \left[ \begin{array}{c} 0 \\ 0 \end{array} \right] \end{array} \right] \\
  & & \ - r(M) - r(E_BA) - r(E_AB) \\ 
 & = & r \left[ \begin{array}{cccc} 
 -M^*MM^*  & 0 & 0 & M^*   \\ 0 &  (E_BA)^*A(E_BA)^* &  0  & (E_BA)^*
  \\ 0  & 0 &
 (E_AB)^*B(E_AB)^*  & (E_AB)^*  \\ A^* & (E_BA)^* & 0 & 0 \\
 B^* & 0 &  (E_AB)^* &  0 \end{array} \right] \\
  & & \ - r(M) - r(E_BA) - r(E_AB) \\ 
& = & r \left[ \begin{array}{cccc}
 0 & M^*A (E_BA)^* &  M^*B (E_AB)^* & M^*  \\ 0 &  (E_BA)^*A(E_BA)^* &  0  &
 (E_BA)^*  \\ 0  & 0  &
 (E_AB)^*B(E_AB)^*  & (E_AB)^*  \\ A^* & 0 & 0 & 0 \\ B^* & 0 & 0 & 0
 \end{array} \right] - r(M) - r(E_BA) - r(E_AB)  \\
 & = & r \left[ \begin{array}{ccc} 
  M^*A (E_BA)^* &  M^*B (E_AB)^* & M^*  \\ (E_BA)^*A(E_BA)^* &  0 &
  (E_BA)^* \\ 0 &(E_AB)^*B(E_AB)^*  & (E_AB)^*  \end{array} \right]
  - r(E_BA) - r(E_AB)  \\
 & = & r \left[ \begin{array}{ccc} 
  0 & 0 & M^*  \\ 0  &  0 &  (E_BA)^* \\ 0  & 0  & (E_AB)^*  \end{array}
  \right] - r(E_BA) - r(E_AB)  \\
 & = & r[\, M, \ E_BA, \ E_AB \, ] - r(E_BA) - r(E_AB) \\
 & = & r[\, A, \ B \, ] - r(E_BA) - r(E_AB)
 = r(A) + r(B) - r[ \, A, \ B \,],
 \end{eqnarray*}
as required for Part (a). Similarly, we can show Part (b). The results in
  Parts (c) and (d) follow immediately from Parts (a) and (b). \qquad  $\Box$

\medskip

 A general  result is given below, its proof is much similar 
to that of Theorem 7.14 and is, therefore, omitted. 

\medskip

\noindent {\bf Theorem 7.15.}\, {\em Let $ A = [\, A_1, \, A_2, \, \cdots, \,
A_k \, ] \in {\cal C}^{ m \times n}$ be given. Then
$$
r \left(  \, [\, A_1, \ A_2, \ \cdots, \ A_k \, ]^{\dagger} -
  \left[ \begin{array}{c} (E_{N_1}A_1)^{\dagger} \\ (E_{N_2}A_2)^{\dagger} \\
  \vdots \\ (E_{N_k}A_k)^{\dagger} \end{array} \right] \, \right)
 = r(A_1) + r(A_2) + \cdots + r(A_k) - r(A), \eqno (7.19)
$$
where $ N_i = [\, A_1, \ \cdots, \ A_{i-1}, \ A_{i+1}, \cdots,\ A_k \,], \ i =
1, \ 2, \ \cdots, \ k.$ In particular$,$
$$
[\, A_1, \ A_2, \ \cdots, \ A_k \, ]^{\dagger} =
\left[ \begin{array}{c} (E_{N_1}A_1)^{\dagger} \\ (E_{N_2}A_2)^{\dagger} \\
  \vdots \\ (E_{N_k}A_k)^{\dagger} \end{array} \right]
\Leftrightarrow r(A) = r(A_1) + r(A_2) + \cdots + r(A_k).\eqno (7.20)
$$ } 
{\bf Theorem 7.16.}\, {\em Let $ A \in {\cal C}^{ m \times n},
  \ B \in {\cal C}^{ m \times k}$ and $ C \in {\cal C}^{ l \times n} $ be
  given. Then

 {\rm (a)} \  $ r \left(  \, [\, A, \ B \,][\, A, \ B \, ]^{\dagger} -
 A(E_BA)^{\dagger}  - B(E_AB)^{\dagger} \, \right) = r(A) + r(B) -
 r[ \, A, \ B \,].$

 {\rm (b)}  \ $ r \left( \,  \left[ \begin{array}{c} A \\ C  \end{array}
 \right]^{\dagger} \left[ \begin{array}{c} A \\ C  \end{array}
 \right] - (AF_C)^{\dagger}A - (CF_A)^{\dagger}C \, \right)
  =  r(A) + r(C) - r\left[ \begin{array}{c} A \\ C  \end{array} \right]. $ \\
In  particular$,$

 {\rm (c)}  \ $ [ \,A, \ B\,][\, A, \ B \, ]^{\dagger} = A(E_BA)^{\dagger} +
 A(E_AB)^{\dagger}  \Leftrightarrow   R(A) \cap R(B) = \{0\}.$

 {\rm (d)} \ $  \left[ \begin{array}{c} A \\ C  \end{array} \right]^{\dagger}
 \left[ \begin{array}{c} A \\ C  \end{array} \right] = (AF_C)^{\dagger}A +
 (CF_A)^{\dagger}C   \Leftrightarrow   R(A^*) \cap R(C^*) = \{0\}.$ } 

\medskip

The proof of Theorem 7.16 is much similar to that of  Theorem 7.14 and is, therefore, omitted. 

\medskip

\noindent {\bf Theorem 7.17.}\, {\em Let $ A \in {\cal C}^{m \times n}, \,
 B \in {\cal C}^{ m \times k}$ and $C \in {\cal C}^{l \times n} $ be given. Then

{\rm (a)} \ $ r\left( \,  [\, A, \ B \,]^{\dagger}[\, A, \ B \,]-
 \left[ \begin{array}{cc} A^{\dagger}A &
   0 \\ 0 & B^{\dagger}B \end{array} \right] \, \right)
    =  r ( A ) + r( B ) -  r[\, A, \ B \, ].$

{\rm (b)} \ $ r\left( \,  [\, A, \ B \,]^{\dagger}[\, A, \ B \,]-
 \left[ \begin{array}{cc}  (E_BA)^{\dagger}(E_BA) & 0 \\ 0 & (E_AB)^{\dagger}(E_AB) \end{array} \right] \, \right)
    =  r ( A ) + r( B ) -  r[\, A, \ B \, ].$

{\rm (c)} \ $ r \left( \left[ \begin{array}{c} A \\ C  \end{array} \right]
 \left[ \begin{array}{c} A \\ C  \end{array} \right]^{\dagger} -
 \left[ \begin{array}{cc} AA^{\dagger} &
   0 \\ 0 & CC^{\dagger} \end{array} \right]  \right)
   =  r ( A ) +  r( C ) - r\left[ \begin{array}{c} A \\ C  \end{array}
    \right].$

{\rm (d)} \ $ r \left( \left[ \begin{array}{c} A \\ C  \end{array} \right]
 \left[ \begin{array}{c} A \\ C  \end{array} \right]^{\dagger} -
 \left[ \begin{array}{cc} (AF_C)(AF_C)^{\dagger} &
   0 \\ 0 & (CF_A)(CF_A)^{\dagger} \end{array} \right]  \right)
   =  r ( A ) +  r( C ) - r\left[ \begin{array}{c} A \\ C  \end{array}
    \right].$ \\
In  particular$,$ 

{\rm (e)} \ $ [\, A, \ B \, ]^{\dagger}[\, A, \ B \, ]
 =  \left[ \begin{array}{cc} A^{\dagger}A &
   0 \\ 0 & B^{\dagger}B \end{array} \right] \Leftrightarrow  [\, A, \ B \,]^{\dagger}[\, A, \ B \,]=
 \left[ \begin{array}{cc}  (E_BA)^{\dagger}(E_BA) & 0 \\ 0 & (E_AB)^{\dagger}(E_AB) \end{array} \right]
 \Leftrightarrow R(A) \cap R(B) = \{0\}.$

{\rm (f)} $ \left[ \begin{array}{c} A \\ C \end{array} \right]
 \left[ \begin{array}{c} A \\ C  \end{array} \right]^{\dagger} =
 \left[ \begin{array}{cc} AA^{\dagger} &   0 \\ 0   & CC^{\dagger}
 \end{array} \right]\Leftrightarrow \left[ \begin{array}{c} A \\ C  \end{array} \right]
 \left[ \begin{array}{c} A \\ C  \end{array} \right]^{\dagger} =
 \left[ \begin{array}{cc} (AF_C)(AF_C)^{\dagger} &
   0 \\ 0 & (CF_A)(CF_A)^{\dagger} \end{array} \right] \Leftrightarrow  R(A^*) \cap R(C^*) = \{0\}. $} 

\medskip

\noindent  {\bf Proof.}\,  Let $ M = [\, A, \ B \,] $ and $N =
\left[ \begin{array}{cc} A & 0 \\ 0 & B \end{array} \right].$
Then we find by Theorem 7.2(b) that
\begin{eqnarray*}
 r\left( \,  [\, A, \ B \,]^{\dagger}[\, A, \ B \,]-
 \left[ \begin{array}{cc} A^{\dagger}A & 0 \\ 0 & B^{\dagger}B
 \end{array} \right] \, \right) & = & r( \, M^{\dagger}M - N^{\dagger}N \,)\\
 & = & 2 r\left[ \begin{array}{c} M \\ N  \end{array} \right] - r(M) - r(N)\\
& = & 2 r\left[ \begin{array}{cc} A & B \\ A & 0  \\ 0 & B  \end{array}
 \right] - r[\, A, \ B \, ] - r(A) - r(B)\\
& = & r( A ) + r(B) - r[\, A, \ B \, ],
\end{eqnarray*}
as required for Part (a). Similarly we can show Parts (b)---(d). Parts (e) and (f) are direct
consequences of  Parts (a) and (b). \qquad $ \Box$ 

\medskip

A general result is given below, and its proof is omitted.  

\medskip

\noindent  {\bf Corollary 7.18.}\, {\em Let $ A = [\, A_1, \, A_2, \, \cdots, \, A_k \,]
 \in {\cal C}^{ m \times n}$ be given. Then
$$\displaylines{
\hspace*{1cm}
r [ \, A^{\dagger}A - {\rm diag}( \,  A_1^{\dagger}A_1, \
  A_2^{\dagger}A_2, \ \cdots, \
  A_k^{\dagger}A_k \, ) \, ]
  =  r ( A_1 ) + r ( A_2 ) + \cdots +  r ( A_k ) -  r(A). \hfill (7.21)
\cr
In \ particular, \hfill
\cr
\hspace*{1cm}
A^{\dagger}A = {\rm diag}( \,  A_1^{\dagger}A_1, \
 A_2^{\dagger}A_2, \ \cdots, \  A_k^{\dagger}A_k \, )
  \Leftrightarrow  r(A) =  r ( A_1 ) + r ( A_2 ) + \cdots +
  r ( A_k ). \hfill (7.22)
\cr}
$$ } 
{\bf Theorem 7.19.}\, {\em Let $ A \in {\cal C}^{ m \times n}, \, B \in
 {\cal C}^{ m \times k}$ and $C \in {\cal C}^{l \times n} $ be given$,$ and 
let $M = \left[ \begin{array}{cc} A & B \\  C & 0 \end{array} \right].$  Then 
 $$\displaylines{
\hspace*{2cm}
 r\left(\, A - [\, A , \ 0 \,]M^{\dagger}\left[ \begin{array}{c} A \\ 0  \end{array}
   \right] \right)  = r(A) + r(B) + r(C) - r(M), \hfill (7.23)
\cr
\hspace*{0cm}
or \ alternatively \hfill
\cr
\hspace*{2cm}
 r(M) = r(A) + r(B) + r(C) - r\left(\, A - [\, A , \ 0 \,] M^{\dagger} \left[ \begin{array}{c} A \\ 0  \end{array}  \right] \, \right).  \hfill (7.24)
\cr
\hspace*{0cm}
In \ particular, \hfill
\cr
\hspace*{2cm}
 A[\, I_n, \ 0 \, ] \left[ \begin{array}{cc} A & B \\ C & 0 \end{array}
  \right]^{\dagger} \left[ \begin{array}{c} I_m \\  0 \end{array} \right] A = A
  \hfill (7.25) 
\cr
\hspace*{0cm}
holds \ if \ and \ only \ if \hfill
\cr
\hspace*{1cm}
 r(M) = r(A) + r(B) + r(C), \ \ i.e., \ \  R(A) \cap R(B) = \{ 0 \} \ \ and  \ \  R(A^*) \cap R(C^*) = \{ 0 \}.
 \hfill (7.26)
 \cr}
$$ }
{\bf Proof.}\, It follows by (2.1) that 
 \begin{eqnarray*}
  r\left( \, A - [\, A, \ 0 \, ]M^{\dagger}\left[ \begin{array}{c} A \\
   0 \end{array} \right] \, \right)
 & = & r \left[ \begin{array}{cc} M^*MM^* & M^*\left[ \begin{array}{c} A \\
 0 \end{array} \right]
  \\ \left[\, A, \  0 \, \right]M^*  & A \end{array} \right]  - r( M ) \\
 & = & r \left[ \begin{array}{cc} M^*MM^* -  M^* \left[ \begin{array}{c} I \\
  0 \end{array} \right]A\left[ \, I, \ 0 \,\right]M^* & 0 \\ 0 & A \end{array}
  \right] - r( M ) \\
 & = & r  \left( \, M^* \left[ \begin{array}{cc} 0 & B \\ C & 0 \end{array}
 \right]M^* \, \right) + r( A )
 - r( M ) \\
& = & r \left[ \begin{array}{cc} BB^*A +  AC^*C & BB^*B \\ CC^*C & 0 \end{array}
 \right] + r( A ) - r( M ) \\
&= & r( A ) + r( B ) + r( C ) - r( M ), 
 \end{eqnarray*}
 as required for  (7.23). \qquad $\Box$ 

\medskip

\noindent  {\bf Theorem 7.20.}\, {\em Let $ A \in {\cal C}^{ m \times n}, \, B \in
 {\cal C}^{ m \times k}$ and $C \in {\cal C}^{l \times n} $ be given$,$ and 
let $M = \left[ \begin{array}{cc} A & B \\  C & 0 \end{array} \right].$  Then 
 $$\displaylines{
\hspace*{0.5cm}
 r\left(\, (E_BAF_C)^{\dagger} - [\, I_n, \ 0 \,]M^{\dagger}\left[ \begin{array}{c} I_m \\ 0  \end{array}
   \right] \right)  = r\left[ \begin{array}{c} A \\ C  \end{array}
   \right] + r[\, A, \ B \,] + r(B) + r(C) - 2r(M), \hfill (7.27)
\cr
or \ alternatively \hfill
\cr
\hspace*{0.5cm}
 r(M) = \frac{1}{2}r\left[ \begin{array}{c} A \\ C  \end{array}
   \right] + \frac{1}{2}r[\, A, \ B \,] + \frac{1}{2}r(B) + \frac{1}{2}r(C) - 
\frac{1}{2}r\left(\, (E_BAF_C)^{\dagger} - 
[\, I_n, \ 0 \,]M^{\dagger}\left[ \begin{array}{c} I_m \\ 0  
\end{array} \right] \right). \hfill (7.28)
\cr
In \ particular, \hfill
\cr
\hspace*{1cm}
 [\, I_n, \ 0 \, ]M^{\dagger} \left[ \begin{array}{c} I_m \\  0 \end{array} \right] 
= (E_BAF_C)^{\dagger} \Leftrightarrow r(M) = r\left[ \begin{array}{c} A \\ C  \end{array}
   \right] + r(B) =  r[\, A, \ B \,] + r(C). \hfill (7.29)
\cr}
$$}
{\bf Proof.}\, Let  $ P = [\, I_n, \ 0 \,]$ and $Q = 
\left[ \begin{array}{c} I_m \\ 0 \end{array} \right].$ It follows by (2.2) that 
\begin{eqnarray*}
\lefteqn{r\left[\, (E_BAF_C)^{\dagger} - PM^{\dagger}Q \right]} \\
 & = & r \left[ \begin{array}{ccc}  (E_BAF_C)^*A(E_BAF_C)^*  & 0 & (E_BAF_C)^* \\ 
0 &  -M^*MM^* & M^*Q \\  (E_BAF_C)^*  & PM^* & 0 \end{array} \right] - r(E_BAF_C) - r(M) \\ 
& = & r \left[ \begin{array}{ccc} 0  & 0 & (E_BAF_C)^* \\ 
0 &  -M^*MM^* + M^*QAPM^*  & M^*Q \\  (E_BAF_C)^*  & PM^* & 0 \end{array} \right] + r(B) + r(C) - 2r(M) \\ 
& = & r \left[ \begin{array}{ccc} 0  & 0 & (E_BAF_C)^* \\ 
0 &  M^* \left[ \begin{array}{cc} 0 & B \\ C & 0
  \end{array} \right]M^*  & M^*Q \\  (E_BAF_C)^*  & PM^* & 0 \end{array} \right] + r(B) + r(C) - 2r(M) \\ 
& = & r \left[ \begin{array}{cccc} 0  & 0  & 0 & E_BAF_C \\ 
0 & 0 & 0 & A \\ 0 & 0 & 0 & C \\ E_BAF_C  & A  & B & 0 \end{array} \right] + r(B) + r(C) - 2r(M) \\ 
& = & r\left[ \begin{array}{c} A \\ C  \end{array} \right] + r[ \, A, \ B \,] + r(B) + r(C) - 2r(M),
\end{eqnarray*}
establishing  (7.27). \qquad $\Box$ 

\medskip

\noindent  {\bf Theorem 7.21.}\, {\em Let $ A \in {\cal C}^{ m \times n}, \,
B \in {\cal C}^{ m \times k}$ and $C \in {\cal C}^{l \times n} $ be given.
Then
$$\displaylines{
\hspace*{1.0cm}
 r\left(\, A + [\, 0 , \ B \,] \left[ \begin{array}{cc} A & B \\ C & 0
  \end{array} \right]^{\dagger} \left[ \begin{array}{c} 0 \\  C \end{array}
   \right] \, \right)  = r[\,  A, \ B \,] +
   r\left[ \begin{array}{c} A \\  C \end{array} \right]  -
   r\left[ \begin{array}{cc} A & B \\ C & 0 \end{array} \right], \hfill(7.30)
\cr
or \ alternatively \hfill
\cr
\hspace*{1.0cm}
 r\left[ \begin{array}{cc} A & B \\ C & 0 \end{array} \right] = 
r[\,  A, \ B \,] +   r\left[ \begin{array}{c} A \\  C \end{array} \right] -
r\left(\, A + [\, 0 , \ B \,] \left[ \begin{array}{cc} A & B \\ C & 0
  \end{array} \right]^{\dagger} \left[ \begin{array}{c} 0 \\  C \end{array}
   \right] \, \right). \hfill (7.31)
\cr}
$$ }
{\bf Proof.}\, Let $M = \left[ \begin{array}{cc} A & B \\ C & 0 \end{array}
 \right]$. Then it follows by  (2.1) and block elementary operation that 
\begin{eqnarray*}
 r\left( \, A + [\, 0, \ B \,]M^{\dagger}\left[ \begin{array}{c} 0 \\ C
  \end{array} \right] \, \right) & = & r \left[ \begin{array}{cc} M^*MM^* & M^*\left[ \begin{array}{c} 0 \\
  C \end{array} \right]
  \\ \left[\, 0, \ B \,\right]M^*  & -A \end{array} \right] - r( M ) \\
 & = & r \left[ \begin{array}{cc} M^*MM^* -  M^* \left[ \begin{array}{c} 0
 \\ I \end{array} \right]
 C\left[  \, I, \ 0 \, \right] M^*  & M^*\left[ \begin{array}{c} 0 \\  C
 \end{array} \right] \\
  \left[\, A, \ B \, \right]M^* & -A \end{array} \right] - r( M ) \\
 & = & r \left[ \begin{array}{cc} M^*MM^* -  M^* \left[ \begin{array}{cc} A
 & B  \\ C & 0  \end{array} \right]M^* & M^*\left[ \begin{array}{c} A \\  C
 \end{array} \right]  \\  \left[\, A, \ B \, \right] M^* & -A \end{array}
  \right] - r( M ) \\
& = & r \left[ \begin{array}{cc} 0 & M^*\left[ \begin{array}{c} A \\  C
 \end{array} \right]
  \\ \left[\, A, \ B \, \right]M^* & -A \end{array} \right] - r( M ) \\
 & = & r \left[ \begin{array}{cc} 0 & \left[ \begin{array}{c} A \\  C
 \end{array} \right]
  \\ \left[\, A, \ B \, \right] & -A \end{array} \right] - r( M )
  = r[\,  A, \ B \, ]+ r\left[ \begin{array}{c} A \\  C \end{array} \right]
  - r(M),
 \end{eqnarray*} 
as requird for  (7.30). \qquad $\Box$ 

\medskip

It is easy to derive from (1.6) that 
$$\displaylines{
\hspace*{2cm}
r\left[ \begin{array}{cc} A & B  \\  C & 0  \end{array} \right] 
\geq r\left[ \begin{array}{c} A  \\  C \end{array} \right]  + r[\, A , \ B\, ]
 - r(A).\hfill
\cr}
$$
Now replacing $A $ by $ A - BXC $ in the above inequality, where $ X $ is
arbitrary, we obtain
$$\displaylines{
\hspace*{2cm}
 r(\, A - BXC \, ) \geq  r[\, A , \ B\, ] + r\left[ \begin{array}{c} A  \\
  C \end{array} \right]  - r\left[ \begin{array}{cc} A  & B  \\  C  & 0
  \end{array} \right]. \hfill (7.32)
\cr}
$$
This rank inequality implies that the quantity in the right-hand side of
 (7.32) is a lower bound for the rank of $ A - BXC $ with respect to
the choice of $ X $. Combining  (7.30) and (7.32), we immediately obtain
$$\displaylines{
\hspace*{2cm}
\min_X r(\, A - BXC \, ) =  r[\, A , \ B\, ] + r\left[ \begin{array}{c} A  \\
  C \end{array} \right]  - r\left[ \begin{array}{cc} A  & B  \\  C  & 0
  \end{array} \right], \hfill (7.33)
\cr}
$$
and a matrix satisfying  (7.33) is given by
$$\displaylines{
\hspace*{2cm}
X = - [\, 0, \ I_l \,] \left[ \begin{array}{cc} A  & B  \\  C  & 0
  \end{array} \right]^{\dagger}
  \left[ \begin{array}{c} 0  \\  I_k \end{array} \right]. \hfill (7.34)
\cr}
$$ 
{\bf Theorem 7.22.}\, {\em Let $ A \in {\cal C}^{m \times n}, \,
 B \in {\cal C}^{ m \times k}$ and $D \in {\cal C}^{l \times k} $ be given.
 Then

 {\rm (a)} \ $ r\left( \, \left[ \begin{array}{cc} A & B \\
 0 & D \end{array} \right]\left[ \begin{array}{cc} A & B \\
 0 & D \end{array} \right]^{\dagger} -
  \left[ \begin{array}{cc} AA^{\dagger} & 0 \\
  0 &  DD^{\dagger} \end{array} \right] \, \right)
  =  r ( D ) - r(A) + 2r[\, A, \ B \, ] - r\left[ \begin{array}{cc} 
A & B \\  0 & D \end{array} \right].$

{\rm (b)} \  $ r\left( \, \left[ \begin{array}{cc} A & B \\
 0 & D \end{array} \right]^{\dagger}\left[ \begin{array}{cc} A & B \\
 0 & D \end{array} \right] -
  \left[ \begin{array}{cc} A^{\dagger}A & 0 \\
 0  & D^{\dagger}D \end{array} \right] \, \right)
 =   r ( A )  - r(D) +  2r\left[ \begin{array}{c} B \\  D \end{array}
  \right] - r\left[ \begin{array}{cc} A & B \\
 0 & D \end{array} \right]. $ \\
In  particular$,$

{\rm (c)} \ $\left[ \begin{array}{cc} A & B \\
 0 & D \end{array} \right]\left[ \begin{array}{cc} A & B \\
 0 & D \end{array} \right]^{\dagger} =
  \left[ \begin{array}{cc} AA^{\dagger} & 0 \\
 0 & DD^{\dagger} \end{array} \right]
  \Leftrightarrow R(B) \subseteq R(A). $

{\rm (d)} \ $ \left[ \begin{array}{cc} A & B \\
 0 & D \end{array} \right]^{\dagger}\left[ \begin{array}{cc} A & B \\
 0 & D \end{array} \right] =
  \left[ \begin{array}{cc} A^{\dagger}A & 0 \\
 0 & D^{\dagger}D \end{array} \right]  \Leftrightarrow
  R(B^*) \subseteq R(D^*).$   }

\medskip

\noindent  {\bf Proof.}\,  Let $ M =
\left[ \begin{array}{cc} A & B \\ 0 & D \end{array} \right]$ and $N =
\left[ \begin{array}{cc} A & 0 \\ 0 & D \end{array} \right].$
Then we find by Theorem 7.2(a) that
\begin{eqnarray*}
r( \, MM^{\dagger} - NN^{\dagger} \,)
 & = & 2 r[\,  M, \ N \, ] - r(M) - r(N) \\
& = & 2 r\left[ \begin{array}{cccc} A & B & A & 0 \\ 0 & D & 0 & D
 \end{array} \right] - r(M) - r(A) - r(D)\\
& = & 2r[\, A, \ B \, ] + r(D) - r(M) - r(A), 
\end{eqnarray*}
as required for Part (a). Similarly we can show Part (b). Observe that  
$$\displaylines{
\hspace*{1cm}
r ( D ) - r(A) + 2r[\, A, \ B \, ] - r\left[ \begin{array}{cc} 
A & B \\ 0 & D \end{array} \right] = \left(  r[\, A, \ B \, ] - r(A) \right) + 
\left(  r ( D )  + r[\, A, \ B \, ] - r\left[ \begin{array}{cc} 
A & B \\  0 & D \end{array} \right] \right), \hfill
 \cr
and \hfill
\cr
\hspace*{1cm}
r ( A )  - r(D) +  2r\left[ \begin{array}{c} B \\  D \end{array}
  \right] - r\left[ \begin{array}{cc} A & B \\
 0 & D \end{array} \right] = \left(  r\left[ \begin{array}{c} B \\  D \end{array}
  \right] - r(D) \right) + \left(  r(A) + r\left[ \begin{array}{c} B \\  D \end{array}
  \right]  -  r\left[ \begin{array}{cc} A & B \\
 0 & D \end{array} \right] \right).   \hfill
\cr}
$$
Thus Parts (c) and (d) follow. \qquad $ \Box$ 

\medskip

A general result is given below, and  the proof is omitted for simplicity.

\medskip

\noindent {\bf Theorem 7.23.}\, {\em Let
$$
 M = \left[ \begin{array}{cccc} A_{11} & A_{12} & \cdots & A_{1k}  \\
  & A_{22} & \cdots & A_{2k} \\ & & \ddots & \vdots \\
  & & & A_{kk}   \end{array} \right] \in {\cal C}^{m \times n}
$$
be given$,$ and let $ A = {\rm diag}( \, A_{11}, \, A_{22}, \, \cdots, \, A_{kk}
 \, )$. Then

{\rm (a)} \ $ r\left[ \, MM^{\dagger} -
{\rm diag}( \,  A_{11}A_{11}^{\dagger}, \ \cdots, \ A_{kk}A_{kk}^{\dagger}
\, ) \, \right] = 2r[\, M, \ A \, ] -
 r(M) - r(A).$

{\rm (b)} \ $ r\left[ \,M^{\dagger}M  -{\rm diag}( \, A_{11}^{\dagger}A_{11},
 \ \cdots, \ A_{kk}^{\dagger}A_{kk} \,) \, \right] =
2r\left[ \begin{array}{c} M \\  A \end{array} \right] - r(M) - r(A).$ 

{\rm (c)} \ $ MM^{\dagger} = {\rm diag}( \, A_{11}A_{11}^{\dagger}, \, \cdots, 
\ A_{kk}A_{kk}^{\dagger}
\, )  \Leftrightarrow  R(M) = R(A) \Leftrightarrow
R(A_{ij}) \subseteq R(A_{ii}), \, j = i +1, \, \cdots, \, k, \, i = 1, \, \cdots, \,
k-1.$

{\rm (d)} \ $ M^{\dagger}M = {\rm diag}( \, A_{11}^{\dagger}A_{11},
 \ \cdots, \ A_{kk}^{\dagger}A_{kk} \,) \Leftrightarrow  R(M^*) = R(A^*)
  \Leftrightarrow
R(A_{ij}) \subseteq R(A_{jj}), \ j = 2, \ \cdots, \ k, \  i = 1, \ \cdots, \
j-1.$  }
 
\medskip

By Theorem 7.2(a) and (b), we can also establish the following.

\medskip

\noindent {\bf Theorem 7.24.}\, {\em Let $ A \in {\cal C}^{m \times n}, \,
 B \in {\cal C}^{ m \times k}$ and $D \in {\cal C}^{l \times k} $ be given.
 Then

 {\rm (a)} \  $ r\left( \, \left[ \begin{array}{cc} A & B \\
 C & D \end{array} \right]\left[ \begin{array}{cc} A & B \\
 C & D \end{array} \right]^{\dagger} -
  \left[ \begin{array}{cc} AA^{\dagger} & 0 \\
  0 &  DD^{\dagger} \end{array} \right] \, \right)
  = 2r[\, A, \ B \, ]+ 2r[\, C, \ D \, ] - r(A) - r(D)- 
r\left[ \begin{array}{cc} A & B \\  C & D \end{array} \right].$

{\rm (b)} \ $ r\left( \, \left[ \begin{array}{cc} A & B \\
 C & D \end{array} \right]^{\dagger}\left[ \begin{array}{cc} A & B \\
 C & D \end{array} \right] - \left[ \begin{array}{cc} A^{\dagger}A & 0 \\
 0  & D^{\dagger}D \end{array} \right] \, \right)
 = 2r\left[ \begin{array}{c} A \\  C \end{array}
  \right] + 2r\left[ \begin{array}{c} B \\  D \end{array}
  \right] - r(A) - r(D)  - r\left[ \begin{array}{cc} A & B \\
 C & D \end{array} \right]. $ \\
In  particular$,$

{\rm (c)} \  $\left[ \begin{array}{cc} A & B \\
 C & D \end{array} \right]\left[ \begin{array}{cc} A & B \\
 C & D \end{array} \right]^{\dagger} =
  \left[ \begin{array}{cc} AA^{\dagger} & 0 \\
 0 & DD^{\dagger} \end{array} \right] \Leftrightarrow 
r\left[ \begin{array}{cc} A & B \\ C & D \end{array} \right] 
= r(A) + r(D), \ R(B) \subseteq R(A)$ and $ R(C) \subseteq R(D)$.

{\rm (d)} \ $ \left[ \begin{array}{cc} A & B \\
 C & D \end{array} \right]^{\dagger}\left[ \begin{array}{cc} A & B \\
 C & D \end{array} \right] =
  \left[ \begin{array}{cc} A^{\dagger}A & 0 \\
 0 & D^{\dagger}D \end{array} \right] \Leftrightarrow 
r\left[ \begin{array}{cc} A & B \\ C & D \end{array} \right] = r(A) + r(D), 
\ R(C^*) \subseteq R(A^*)$ and $R(B^*) \subseteq R(D^*)$.}

\markboth{YONGGE  TIAN }
{8. REVERSE ORDER LAWS FOR MOORE-PENROSE INVERSES}

\chapter{Reverse order laws for  Moore-Penrose inverses}

\noindent Reverse order laws for generalized inverses of  products of matrices have been an attractive topic in the theory of generalized inverses of matrices, for
 these laws can reveal essential relationships between generalized
 inverses of  products of matrices and generalized inverses of each matrix in
 the products. Various results on reverse order laws
related to inner inverses, reflexive inner inverses,
Moore-Penrose  inverses,  group inverses, Drazin inverses, and weighted
Moore-Penrose inverses of products  of matrices have widely been established by lot of 
authors (see, e.g.,  \cite{BG1, BG2, BG3, Cl1, DW, Er,  GP, Gr1, Ha4,  SS1, SS2, SW, Ti2,  Ti4, We1, We2, 
WHG}). In this chapter,
we shall present some rank equalities related to products of 
Moore-Penrose inverses of matrices, and then derive from them various types of reverse order 
laws for  Moore-Penrose inverses of  products of matrices. 

\medskip

\noindent {\bf Theorem 8.1.}\, {\em Let $ A \in {\cal C}^{ m \times n}$ and
 $ B \in {\cal C}^{ n \times p}$ be 
given. Then
$$
 r( \, AB - ABB^{\dagger}A^{\dagger}AB \,) = r( \, B^{\dagger}A^{\dagger} - 
B^{\dagger}A^{\dagger}ABB^{\dagger}A^{\dagger}\,) = r[\, A^*, \ B \, ] +
 r(AB)- r(A)- r(B). \eqno (8.1)
$$
In particular$,$ the following seven statements are equivalent$:$

{\rm (a)} \  $B^{\dagger}A^{\dagger} \in \{ ( AB )^-_r \},$  i.e.$,$   
 $ B^{\dagger}A^{\dagger}$  is a reflexive inner inverse of $AB$.  

{\rm (b)} \  $ r[\, A^*, \  B \, ] = r(A) + r(B) - r(AB).$

{\rm (c)} \ $ {\rm dim}[ R(A) \cap R(B^*)] = r(AB).$

{\rm (d)} \ $r(\, B - A^{\dagger}AB \, ) = r(B) - r( A^{\dagger}AB),$ i.e.$,$  
$A^{\dagger}AB \leq_{rs} B$.  

{\rm (e)} \  $ r(\, A - ABB^{\dagger} \, ) = r(A) - r(ABB^{\dagger}),$ i.e.$,$  
$ABB^{\dagger} \leq_{rs} A$.  

{\rm (f)\cite{HaSp2}} \ $B^{\dagger}A^{\dagger}ABB^{\dagger}A^{\dagger} = B^{\dagger}A^{\dagger}.
$

{\rm (g)\cite{HaSp2}} \ $ AA^{\dagger}B^{\dagger}B = B^{\dagger}BAA^{\dagger}.$}

\medskip

\noindent  {\bf Proof.}\, Applying  (2.8) and (1.7) to
 $ AB - ABB^{\dagger}A^{\dagger}AB,$ we obtain
 \begin{eqnarray*}
r( \, AB - ABB^{\dagger}A^{\dagger}AB \,) & = & r \left[ \begin{array}{ccc}
  B^*A^*  & B^*BB^* & 0  \\ A^*AA^* & 0 &  A^*AB  \\ 0 &  ABB^* & -AB
  \end{array} \right] - r( A )  - r( B ) \\
  & = &  r \left[ \begin{array}{ccc}
  B^*A^*  & B^*B & 0  \\ AA^* & 0 & AB  \\ 0 &  AB & -AB  \end{array} \right]
  - r( A )  - r( B ) \\
  & = &  r \left[ \begin{array}{cc}
  B^*A^*  & B^*B  \\ AA^* & AB  \end{array} \right]  + r( AB) - r( A )
   - r( B ) \\
  & = &  r \left(  \, \left[ \begin{array}{c} A \\ B^*  \end{array} \right]
   [\, A^*, \ B \, ] \, \right) + r( AB) - r( A )  - r( B ) \\
  & = & r(\, [\, A^*, \  B \, ]^*[\, A^*, \  B \, ] \,) + r(AB)- r(A)- r(B) \\
 & = & r[\, A^*, \  B \, ] + r(AB)- r(A)- r(B).
 \end{eqnarray*}
Thus we have the first part of (8.1). Replace $ A$ by $B^{\dagger}$  and $ B$ by $A^{\dagger}$ 
and simplify to yield the second part of   (8.1).  The equivalence of Parts (a), (b) and (f) 
follows  immediately from (8.1). The equivalence of Parts (b) and (c) follows from the well-known 
rank  formula
$$
r[\, A^*, \ B \, ] = r(A) + r(B ) - {\rm dim}[R(A^*)\cap R(B)].
$$
The equivalence of Parts (b), (d) and (e) follows from  (1.2) and (1.3). The equivalence of Parts (b) 
and (g) follows from  (7.2). \qquad $\Box$ 

\medskip

The rank formula (8.1) was established by Baksalary and Styan \cite{BaSt} in an alternative form 
$$
 r(AE_BF_AB) = r[\, A^*, \ B \, ] + r(AB)- r(A)- r(B). \eqno (8.1')
$$  
Observe that 
$$ 
A( \, I  - BB^{\dagger} \,)(\,  I -  A^{\dagger}A \,)B =  - AB + 
ABB^{\dagger}A^{\dagger}AB. 
$$  
Thus $(8.1')$ is exactly (8.1). Some extensions and applications of  
 $(8.1')$ in mathematical statistics were also considered by Baksalary and 
Styan \cite{BaSt}. But in this monograph we only consider the application of 
 (8.1) to the reverse order law $B^{\dagger}A^{\dagger} \in \{ ( AB )^- \}.$
 In addition, the results in Theorem 8.1 can also be extended to a product of 
$ n $ matrices. The corresponding results were presented by the author 
in \cite{Ti4}.  

\medskip

As an application of (8.1), we let $ B = I_m - A$ in (8.1). Then 
$$\displaylines{
\hspace*{2cm}
 r[ \, (A - A^2)  - (A - A^2)(I - A)^{\dagger}A^{\dagger}(A - A^2) \,] \hfill
\cr
\hspace*{2cm}
= r[\, A^*, \  I_m - A  \, ] +  r(A -A^2) - r(A) - r( I_m - A) \hfill
\cr
\hspace*{2cm}
= r[\, A^*, \  I_m - A  \, ] - m \leq  0. \hfill
\cr}
$$  
This inequality implies that $r[ \, A^*, \ I_m - A \,] = m $ and   $ (I - A)^{\dagger}A^{\dagger}$ is a reflexive
 inner inverse of the matrix $ A - A^2$. By symmetry,  $ (I - A)^{\dagger}A^{\dagger}$ is also  a reflexive
 inner inverse of the matrix $ A - A^2$. 

Replacing  $ A$ and $ B $ in (8.1) by $ I_m + A$ and $ I_m - A $, respectively, we then get 
$$\displaylines{
\hspace*{2cm}
 r[ \, (I_m - A^2)  - (I_m - A^2)(I - A)^{\dagger}(I_m + A)^{\dagger}(I_m - A^2) \,] \hfill
\cr
\hspace*{2cm}
= r[\, I_m +  A^*, \  I_m - A  \, ] +  r( I_m  -A^2 ) - r(I_m + A ) - r( I_m - A) \hfill
\cr
\hspace*{2cm}
= r[\, I_m +  A^*, \  I_m - A  \, ] - m \leq  0. \hfill
\cr}
$$
This inequality implies that $r[ \, I_m + A^*, \ I_m - A \,] = m $ and  $ (I - A)^{\dagger}(I_m + A)^{\dagger}$ is a 
reflexive inner inverse of the matrix $ I_m - A^2$. By symmetry,  $ (I + A)^{\dagger}(I_m - A)^{\dagger}$ is also 
 a reflexive inner inverse of the matrix $ I_m - A^2$.

In general, for any two polynomials $ p(\lambda)$ and  $ q(\lambda)$ without common roots, 
we find by (8.1) and (1.17) the following   
$$
\displaylines{
\hspace*{1.5cm}
 r[ \, p(A)q(A) - p(A)q(A)q^{\dagger}(A)p^{\dagger}(A)p(A)q(A) \, ] \hfill
\cr
\hspace*{1.5cm}
= r[ \, p(A^*), \ q(A) \, ] + r[p(A)q(A)] - r[p(A)] - r[q(A)] \hfill
\cr
\hspace*{1.5cm}
 = r[ \,  p(A^*), \ q(A) \, ] - m \leq 0. \hfill
\cr}
$$ 
This implies that $ r[ \,  p(A^*), \ q(A) \, ] = m$ and 
$ q^{\dagger}(A)p^{\dagger}(A)$ is a reflexive inner inverse of the matrix $ p(A)q(A)$.
By symmetry, $ p^{\dagger}(A)q^{\dagger}(A)$ is also a reflexive inner inverse of the 
matrix $ p(A)q(A)$.

\medskip

\noindent {\bf Theorem 8.2.}\, {\em Let $ A \in {\cal C}^{ m \times n}$ and 
$ B \in {\cal C}^{ n \times p}$ be given. Then

{\rm (a)} \ $ r[ \, (AB)(AB)^{\dagger} - (AB)(B^{\dagger}A^{\dagger}) \,] 
= r[\, B, \ A^*AB \, ] - r(B) = r( \, A^*AB - BB^{\dagger}A^*AB \,).$ 
 
{\rm (b)}  \ $ r[ \, (AB)^{\dagger}(AB) - (B^{\dagger}A^{\dagger})(AB) \,] 
= r\left[ \begin{array}{c} A \\ 
ABB^*  \end{array} \right] - r(A)= r( \, ABB^* - ABB^*A^{\dagger}A \,).$\\
In  particular$,$

{\rm (c)} \ $ (AB)(AB)^{\dagger} = (AB)(B^{\dagger}A^{\dagger})  \Leftrightarrow
 A^*AB = BB^{\dagger}A^*AB \Leftrightarrow R(A^*AB)
 \subseteq R(B)  \Leftrightarrow B^{\dagger}A^{\dagger} \subseteq \{ \,(AB)^{(1,2,3)} \}$.
 
{\rm (d)} \ $ (AB)^{\dagger}(AB) = (B^{\dagger}A^{\dagger})(AB) \Leftrightarrow
 ABB^* = ABB^*A^{\dagger}A \Leftrightarrow R(BB^*A^* ) 
\subseteq R(A^*)$ $ \Leftrightarrow  B^{\dagger}A^{\dagger} \subseteq \{ \,(AB)^{(1,2,4)} \}$. 
 
{\rm (e)}\, The following four statements are equivalent$:$

\hspace*{0.3cm} {\rm (1)}  $(AB)^{\dagger}= B^{\dagger}A^{\dagger}.$

\hspace*{0.3cm} {\rm (2)}  $(AB)(AB)^{\dagger} = (AB)(B^{\dagger}A^{\dagger}) \ and \ 
 (AB)^{\dagger}(AB) = (B^{\dagger}A^{\dagger})(AB)$.

\hspace*{0.3cm} {\rm (3)}  $ A^*AB = BB^{\dagger}A^*AB \ and \ ABB^* = ABB^*A^{\dagger}A$.
 
\hspace*{0.3cm} {\rm (4)} $R(A^*AB) \subseteq R(B) \  and  \ R(BB^*A^*) \subseteq R(A^*)$.}  

\medskip

\noindent {\bf Proof.}\, Let $ N = AB$. Then by (2.1), (2.7) and (1.8),
it follows that
\begin{eqnarray*}
r( \, NN^{\dagger} - NB^{\dagger}A^{\dagger} \,) & = & r \left[ \begin{array}{cc}
 N^*NN^*  & N^*  \\ NN^* &  NB^{\dagger}A^{\dagger}  \end{array} \right] - r( N )  \\
& = & r \left[ \begin{array}{cc}
 0  & N^* -  N^*NB^{\dagger}A^{\dagger}  \\ NN^* & 0  \end{array} \right] - r( N )  \\
& = & r( \, N^* - N^*NB^{\dagger}A^{\dagger} \,) \\
& = & r \left[ \begin{array}{ccc}
 B^*A^*  & B^*BB^* &  0  \\ A^*AA^* & 0 &  A^*  \\ 0 &  N^*NB^* & -N^*  \end{array} \right] 
- r( A )  - r( B ) \\
& = & r \left[ \begin{array}{ccc}
 B^*A^*  & B^*B &  0  \\ 0 & 0 &  A^*  \\ B^*A^*AA^* &  N^*N & 0  \end{array} \right] 
- r( A )  - r( B ) \\
& = &  r \left[ \begin{array}{cc}
 B^*A^*  & B^*B  \\ B^*A^*AA^* & N^*N  \end{array} \right] - r( B ) \\
& = &  r \left[ \begin{array}{cc}
 B^*B  & B^*A^*AB  \\ AB & AA^*AB  \end{array} \right] - r( B ) \\
 & = & r \left( \,  \left[ \begin{array}{c} B^* \\ A  \end{array} \right]
[\,  B,  \ A^*AB \, ] \, \right) - r( B )  = r[\, B, \ A^*AB \, ] - r(B),
\end{eqnarray*}
as required for the first equality in Part (a). Applying  (1.2) to it the block matrix in it yields the 
second  equality in Part (a).  Similarly, we can establish Part (b). The results in
Parts (c) and (d) are direct consequences of Parts (a) and (b). The result  in Part (e) follows directly
from Parts (c) and (d).  \qquad $\Box$

\medskip

The result in Theorem 8.2(e) is well known, see, e.g., Arghiriade \cite{Ar}, Rao and Mitra \cite{RM}, 
 Ben-Israel and Greville \cite{BeG}, Campbell and  Meyer \cite{CM2}. Now it can be regarded as a direct 
consequence of some rank equalities related to Moore-Penrose inverses of products of two  
 matrices. We next present another group rank equalities related to Moore-Penrose inverses of products 
of two matrices, which can also help to establish necessary and sufficient conditions for 
$(AB)^{\dagger}= B^{\dagger}A^{\dagger}$.

\medskip

\noindent {\bf Theorem 8.3.}\, {\em Let $ A \in {\cal C}^{m \times n}$ and 
$ B \in {\cal C}^{n \times p}$ be 
 given. Then 
 
{\rm (a)}  \ $ r[ \, ABB^{\dagger} - (AB)(AB)^{\dagger}A \,] 
= r[\, B, \ A^*AB \, ] - r(B).$ 

{\rm (b)}  \ $ r[ \, A^{\dagger}AB - B(AB)^{\dagger}(AB) \,] 
= r\left[ \begin{array}{c} A \\ ABB^*  \end{array} \right] - r(A).$

 {\rm (c)}  \ $ r[ \, A^*ABB^{\dagger} - BB^{\dagger}A^*A \,] 
= 2r[\, B, \ A^*AB \, ] - 2r(B).$ 
 
{\rm (d)}  \ $ r[ \, A^{\dagger}ABB^* - BB^*A^{\dagger}A \,] 
= 2 r\left[ \begin{array}{c} A \\ ABB^* 
  \end{array} \right] - 2r(A).$\\
In  particular$,$

{\rm (e)} \ $(AB)(AB)^{\dagger}A = ABB^{\dagger}$ 
$ \Leftrightarrow$ $ A^*ABB^{\dagger} = BB^{\dagger}A^*A$ 
 $ \Leftrightarrow$ $R(A^*AB) \subseteq R(B) \Leftrightarrow B^{\dagger}A^{\dagger} \subseteq \{ \,(AB)^{(1,2,3)} \}$.

{\rm (f)} \ $ A^{\dagger}AB = B(AB)^{\dagger}(AB)$ $ \Leftrightarrow$ 
$A^{\dagger}ABB^* = BB^*A^{\dagger}A$ $ \Leftrightarrow$ 
$ R(BB^*A^*) \subseteq R(A^*) \Leftrightarrow  B^{\dagger}A^{\dagger} \subseteq \{ \,(AB)^{(1,2,4)} \}$.

{\rm (g)}\, The following three statements are equivalent (Greville \cite{Gr1})$:$

\hspace*{0.5cm} {\rm (1)}  $(AB)^{\dagger}= B^{\dagger}A^{\dagger}.$

\hspace*{0.5cm} {\rm (2)} $(AB)(AB)^{\dagger}A = ABB^{\dagger} \ and \ A^{\dagger}AB = B(AB)^{\dagger}(AB)$.

\hspace*{0.5cm} {\rm (3)}  $ A^*ABB^{\dagger} = BB^{\dagger}A^*A \ and \ 
 A^{\dagger}ABB^* = BB^*A^{\dagger}A$. }
 
\medskip

\noindent {\bf Proof.}\, We only show Part (b). Note that both $AA^{\dagger}$ 
and $ (AB)^{\dagger}(AB)$ are idempotent. We have by (3.1) that
$$\displaylines{
\hspace*{1cm}
r[ \, A^{\dagger}AB -  B(AB)^{\dagger}(AB) \,] \hfill
\cr
\hspace*{1cm}
 =  r \left[ \begin{array}{c} A^{\dagger}AB \\ (AB)^{\dagger}(AB)  \end{array} \right]
 + r[ \, B(AB)^{\dagger}(AB), \  A^{\dagger}A \, ]
 - r ( A^{\dagger}A) - r[(AB)^{\dagger}(AB)]  \hfill
\cr
\hspace*{1cm}
=r(AB) + r[ \, B(AB)^* , \  A^* \, ] -r(A) -r(AB)\\
=  r[ \, BB^*A^* , \  A^* \, ] -r(A), \hfill
\cr}
$$
as required.  \qquad  $\Box$ 

\medskip

\noindent {\bf Theorem 8.4.}\, {\em Let $ A \in {\cal C}^{m \times n}$ and 
$ B \in {\cal C}^{n \times p}$ be 
 given. Then 
 
{\rm (a)} \ $ r[ \, A^{\dagger} - B(AB)^{\dagger} \,]  = r\left[ \begin{array}{c} A \\ ABB^*  \end{array} \right] - r(AB).$ 

{\rm (b)} \ $ r[ \, B^{\dagger} - (AB)^{\dagger}A \,] = r[\, B, \ A^*AB \, ] - r(AB).$ \\
In  particular$,$

{\rm (c)} \ $ A^{\dagger} = B(AB)^{\dagger}   \Leftrightarrow  R(A^*) = R( BB^*A^*).$ 

{\rm (d)} \ $ B^{\dagger} = (AB)^{\dagger}A   \Leftrightarrow R(A^*) = R( A^*AB).$}

\medskip
\noindent {\bf Proof.}\, We only show Part (a).  According to  (2.2), we  find that 
\begin{eqnarray*}
 r[ \, A^{\dagger} - B(AB)^{\dagger} \,] & = & 
r \left[ \begin{array}{ccc} A^*AA^* &  0 & A^* \\ 0 &  -(AB)^*AB(AB)^* & (AB)^* 
\\ A^* & B(AB)^* & 0  \end{array} \right] - r(A) - r(AB) \\
& = & r \left[ \begin{array}{ccc} A^*AA^* &  0 & A^* \\ (AB)^*AA^*  &  & (AB)^* 
\\ A^* & B(AB)^* & 0  \end{array} \right] - r(A) - r(AB) \\
& = & r \left[ \begin{array}{ccc} 0 &  0 & A^* \\ 0  &  & (AB)^* 
\\ A^* & B(AB)^* & 0  \end{array} \right] - r(A) - r(AB) \\
& = & r\left[ \begin{array}{c} A \\ ABB^*  \end{array} \right] - r(AB), 
\end{eqnarray*}
establishing Part (a). \qquad  $\Box$ 

\medskip

We next consider ranks of  matrix 
expressions involving Moore-Penrose inverses of products of three matrices, and then present their 
consequences related to reverse order laws.
  
\medskip  

\noindent {\bf Theorem 8.5.}\, {\em Let $ A \in {\cal C}^{ m \times n}, \,  
B \in {\cal C}^{ n \times p}$  and $  C \in {\cal C}^{p \times q}$ be given$,$
 and let $ M = ABC$. Then
$$
  r[ \, M - M(BC)^{\dagger}B(AB)^{\dagger}M \,]
 = r \left(  \, \left[ \begin{array}{c} (BC)^* \\ A  \end{array} \right] B
 [\, (AB)^*, \ C \,] \,  \right)  + r(M)- r(AB)- r(BC). \eqno (8.2)
$$
In particular$,$ 
$$
 (BC)^{\dagger}B(AB)^{\dagger} \in \{ (ABC)^- \}  \Leftrightarrow
  r \left(\, \left[ \begin{array}{c} (BC)^* \\ A  \end{array} \right]B[\,
  (AB)^*, \ C \,] \, \right) = r(AB)+ r(BC) - r(M).\eqno (8.3)
$$ }
{\bf Proof.}\, Applying (2.8) and the rank cancellation law (1.8) to
 $ M - M(BC)^{\dagger}B(AB)^{\dagger}M,$ we obtain
\begin{eqnarray*}
 \lefteqn{ r[ \,M - M(BC)^{\dagger}B(AB)^{\dagger}M \,] } \\ 
 & = & r \left[ \begin{array}{ccc} (BC)^*B(AB)^* & (BC)^*(BC)(BC)^* & 0
  \\ (AB)^*(AB)(AB)^* & 0 &  (AB)^*M  \\  0 &  M(BC)^* & -M  \end{array}
  \right] - r( AB )  - r( BC ) \\
 & = & r \left[ \begin{array}{ccc} (BC)^*B(AB)^*  & (BC)^*(BC) & 0  \\
 (AB)(AB)^* & 0 & M  \\  0 &  M & -M  \end{array} \right] - r( AB )
  - r( BC ) \\
 & = & r \left[ \begin{array}{cc} (BC)^*B(AB)^*  & (BC)^*(BC)  \\
 (AB)(AB)^* & M  \end{array} \right]+ r(M) - r( AB )  - r( BC ) \\
 & = & r \left(  \, \left[ \begin{array}{c} (BC)^* \\ A  \end{array} \right]
  B [\, (AB)^*, \ C \,]\,  \right)  + r(M)- r(AB)- r(BC).
\end{eqnarray*}
 Thus we have  (8.2) and (8.3).  \qquad $\Box$ 

\medskip

As an application of (8.2), we consider the matrix product $ M  = (I_m + A)A(I_m - A) = A -  A^3 $.  Then 
$$\displaylines{
\hspace*{2cm}
 r[ \, (A -  A^3) - (A - A^3)(A - A^2)^{\dagger}A( A +  A^2)^{\dagger}(A - A^3) \,] \hfill
\cr
\hspace*{2cm}
= r \left(  \left[ \begin{array}{c} (A - A^2)^* \\ I_m + A  \end{array} \right]A
 [\, (A +A^2)^*, \ I_m - A  \,] \right) + r(A - A^3)- r( A +  A^2) -  r( A - A^2) \hfill
\cr
\hspace*{2cm}
=  r \left( \left[ \begin{array}{c} (A - A^2)^* \\ I_m + A  \end{array} \right]A
 [\, (A +A^2)^*, \ I_m - A  \,]  \right) - r(A)  \leq  0. \hfill
\cr}
$$  
Notice that the rank of a matrix is nonnegative. The above inequality in fact 
implies that 
$$
 r \left( \left[ \begin{array}{c} (A - A^2)^* \\ I_m + A  \end{array} \right]A
 [\, (A +A^2)^*, \ I_m - A  \,] \right) = r(A) \ \  {\rm and}  \ \  
(A - A^3)(A - A^2)^{\dagger}A( A +  A^2)^{\dagger}(A - A^3) =(A -  A^3),
$$ 
that is, the matrix product $ (A - A^2)^{\dagger}A( A +  A^2)^{\dagger}$ is 
an inner inverse of the matrix $ A - A^3$. By symmetry,  $(A + A^2)^{\dagger}A( A - A^2)^{\dagger}$
 is also an inner inverse of the matrix $ A - A^3$.  

In general, for any three polynomials $ p_1(\lambda), \, p_2(\lambda)$ and $ p_3(\lambda)$  
without common roots and a square matrix $A$,  we let $ p(A) = p_1(A)p_2(A)p_3(A)$. Then 
we can find by (8.2) and (1.17) the following   
$$
\displaylines{
\hspace*{0.5cm}
 r[ \, p(A) - p(A)[p_2(A)p_3(A)]^{\dagger}p_2(A)[p_1(A)p_2(A)]^{\dagger}p(A) \, ] \hfill
\cr
\hspace*{0.5cm}
= r \left(  \left[ \begin{array}{c} p_2(A^*)p_3(A^*) \\ p_1(A)  \end{array} \right]p_2(A)
 [\, [p_1(A^*)p_2(A^*), \ p_3(A)  \,]  \right) + r[p(A)]- r[p_1(A)p_2(A)] - r[p_2(A)p_3(A)]   \hfill
\cr
\hspace*{0.5cm}
= r \left(  \left[ \begin{array}{c} p_2(A^*)p_3(A^*) \\ p_1(A)  \end{array} \right]p_2(A)
 [\, [p_1(A^*)p_2(A^*), \ p_3(A)  \,]  \right) - r[p_2(A)] \leq 0.   \hfill
\cr}
$$  
This  implies that
$$
r \left(  \left[ \begin{array}{c} p_2(A^*)p_3(A^*) \\ p_1(A)  \end{array} \right]p_2(A)
 [\, [p_1(A^*)p_2(A^*), \ p_3(A)  \,]  \right) = r[p_2(A)]
$$
and 
$$
p(A)[p_2(A)p_3(A)]^{\dagger}p_2(A)[p_1(A)p_2(A)]^{\dagger}p(A) = p(A).
$$
Thus $[p_2(A)p_3(A)]^{\dagger}p_2(A)[p_1(A)p_2(A)]^{\dagger}$ is an inner  inverse  of 
the matrix product $ p_1(A)p_2(A)p_3(A)$. By symmetry, 
$$[p_1(A)p_2(A)]^{\dagger}p_2(A)[p_2(A)p_3(A)]^{\dagger}, \ \ \ \
 [p_1(A)p_2(A)]^{\dagger}p_1(A) [p_1(A)p_3(A)]^{\dagger}
$$
$$ [p_1(A)p_3(A)]^{\dagger}p_1(A) [p_1(A)p_2(A)]^{\dagger}, \ \ 
 [p_1(A)p_3(A)]^{\dagger}p_3(A) [p_2(A)p_3(A)]^{\dagger}, \ \   
 [p_2(A)p_3(A)]^{\dagger}p_3(A)[p_1(A)p_3(A)]^{\dagger}
$$
are all inner inverses of the matrix product $p_1(A)p_2(A)p_3(A)$.

\medskip

\noindent {\bf Theorem 8.6.}\, {\em Let $ A \in {\cal C}^{ m \times n}, \,
  B \in {\cal C}^{ n \times p}$  and $ C \in {\cal C}^{p \times q}$ be given$,$
  and let $ M = ABC$. Then

{\rm (a)}\,  The rank of $ M^{\dagger} - (BC)^{\dagger}B(AB)^{\dagger}$ satisfies the equality
$$
\displaylines{
\hspace*{2cm}
r[ \, M^{\dagger} - (BC)^{\dagger}B(AB)^{\dagger} \,]
 = r \left(  \, \left[ \begin{array}{c} (BC)^* \\ M^*A  \end{array} \right]
 B[\, (AB)^*, \ CM^* \,] \,\right) - r(M).  \hfill (8.4)
\cr}
$$

{\rm (b)}\, The following three statements are equivalent$:$

\hspace*{0.5cm} {\rm (1)} $ (ABC)^{\dagger} = (BC)^{\dagger}B(AB)^{\dagger}.$

\hspace*{0.5cm} {\rm (2)}  $ r\left[ \begin{array}{cc} MM^*M & M(BC)^*(BC) \\
      (AB)(AB)^*M & ABB^*BC  \end{array} \right] = r(ABC).$

\hspace*{0.5cm} {\rm (3)}  $ ABB^*BC = AB( BCM^{\dagger}AB)^* BC.$

{\rm (c)}\, If $ r(ABC) = r(B),$ then  
$$ \displaylines{
\hspace*{2cm}
(ABC)^{\dagger} = (BC)^{\dagger}B(AB)^{\dagger} \ \  and \ \ 
(ABC)^{\dagger} = (B^{\dagger}BC)^{\dagger}B^{\dagger}(ABB^{\dagger})^{\dagger}. \hfill (8.5) 
\cr}
$$ } 
{\bf Proof.}\, Applying  (2.12) to $ M^{\dagger} - 
(BC)^{\dagger}B(AB)^{\dagger},$ we obtain 
 \begin{eqnarray*}
 \lefteqn{ r[ \,M^{\dagger} - (BC)^{\dagger}B(AB)^{\dagger} \,] } \\
 & = & r \left[ \begin{array}{cccc} M^*MM^*  & 0 & 0 & M^* \\  0 &
  (BC)^*B(AB)^* & (BC)^*(BC)(BC)^*  & 0
 \\ 0 &  (AB)^*(AB)(AB)^* & 0 &  (AB)^* \\  M^* & 0 & (BC)^* & 0  \end{array}
  \right] -r(M)- r( AB )  - r( BC ) \\
 & = & r \left[ \begin{array}{cccc} M^*MM^*  & -M^*(AB)(AB)^* & 0 & 0 \\
  -(BC)^*(BC)M^*  & (BC)^*B(AB)^*
 & 0 & 0  \\ 0 & 0 & 0 &  (AB)^* \\ 0  &  0 & (BC)^* & 0  \end{array} \right]
 - r(M)- r( AB )  - r( BC ) \\
 & = & r \left[ \begin{array}{cc} (BC)^*B(AB)^* & (BC)^*(BC)M^* \\
 M^*(AB)(AB)^* & M^*MM^*  \end{array} \right] - r(M) \\
 & = & r \left(\, \left[ \begin{array}{c} (BC)^* \\ M^*A  \end{array} \right]
 B[\, (AB)^*, \ CM^* \,] \, \right) - r(M),
 \end{eqnarray*}
as required for (8.4).  Then the equivalence of  Statements (1) and (2) in  Part (b) follows
immediately from (8.4), and the equivalence of  Statements (2) and (3) in  Part (b) follows
from Lemma 1.2(f).  If $ r(ABC ) = r(B)$, then
$$\displaylines{
\hspace*{1cm}
r \left[ \begin{array}{cc} MM^*M & M(BC)^*(BC)\\ (AB)(AB)^*M &
ABB^*BC  \end{array} \right] \geq  r(MM^*M ) = r(M). \hfill
\cr}
$$
On the other hand, 
$$\displaylines{
\hspace*{1cm}
r \left[ \begin{array}{cc} MM^*M & M(BC)^*(BC)\\ (AB)(AB)^*M &
ABB^*BC  \end{array} \right] = r\left( \,
\left[ \begin{array}{c} MC^* \\ AB  \end{array} \right] B^* [\, A^*M ,
\ BC \,] \, \right) \leq r(B) = r(M). \hfill
\cr}
$$
Thus we have
$$\displaylines{
\hspace*{1cm}
r \left[ \begin{array}{cc} MM^*M & M(BC)^*(BC)\\ (AB)(AB)^*M &
ABB^*BC  \end{array} \right] = r(M). \hfill
\cr}
$$
Thus according to the statements (1) and (2) in  Part (b), we know that the first
equality in (8.5) is true. The second  equality in  (8.5) follows from
writing $ ABC = ABB^{\dagger}BC$ and then applying the
 first equality to it. \qquad $\Box$ 

\medskip

\noindent {\bf Theorem 8.7.}\, {\em Let $ A \in {\cal C}^{ m \times n}, \,
  B \in {\cal C}^{ n \times p}$  and $C \in {\cal C}^{p \times q}$ be given$,$
  and let $ M = ABC$. Then
$$\displaylines{
\hspace*{1cm}
r[ \, B^{\dagger} - (AB)^{\dagger}M(BC)^{\dagger} \,]
 = r \left[ \begin{array}{cc} M & (AB)(AB)^* \\ (BC)^*(BC) &
 (BC)^*B(AB)^* \end{array} \right] + r(B) - r(AB)  - r(BC). \hfill (8.6)
\cr
\hspace*{0cm}
In \ particular, \hfill
\cr
\hspace*{1cm}
B^{\dagger} = (AB)^{\dagger}M(BC)^{\dagger}  \Leftrightarrow 
  r\left[ \begin{array}{cc} M & (AB)(AB)^* \\ (BC)^*(BC) & (BC)^*B(AB)^* \end{array} \right]= r(AB)  + r(BC) - r(B). 
\hfill
\cr}
$$}
{\bf Proof.}\, Applying (2.11) to $ B^{\dagger} -
(AB)^{\dagger}M(BC)^{\dagger},$ we obtain 
 \begin{eqnarray*}
 \lefteqn{ r[ \,B^{\dagger} - (AB)^{\dagger}M(BC)^{\dagger} \,] } \\
 & = & r \left[ \begin{array}{cccc} B^*BB^*  & 0 & 0 & B^* \\  0 &
  (AB)^*M(BC)^* & (AB)^*(AB)(AB)^*  & 0
 \\ 0 &  (BC)^*(BC)(BC)^* & 0 &  (BC)^* \\  B^* & 0 & (AB)^* & 0  \end{array}
  \right] -r(B) - r( AB )  - r( BC ) \\
 & = & r \left[ \begin{array}{cccc} 0  & 0 & 0 & B^* \\  0 &
  (AB)^*M(BC)^* & (AB)^*(AB)(AB)^*  & 0
 \\ 0 &  (BC)^*(BC)(BC)^* &(BC)^* B(AB)^* & 0 \\  B^* & 0 & 0 & 0  \end{array}
  \right] -r(B) - r( AB )  - r( BC ) \\
 & = & r \left[ \begin{array}{cc} (AB)^*M(BC)^* & (AB)^*(AB)(AB)^* \\ 
(BC)^*(BC)(BC)^* &(BC)^*B(AB)^* \end{array}
  \right] + r(B) - r( AB )  - r( BC ) \\
& = & r \left[ \begin{array}{cc} M & (AB)(AB)^* \\ 
(BC)^*(BC) &(BC)^*B(AB)^* \end{array}
  \right] + r(B) - r( AB )  - r( BC ), \\
\end{eqnarray*}
as required for (8.6).  \qquad $\Box$

\medskip

\noindent {\bf Theorem 8.8.}\, {\em Let $ A \in {\cal C}^{ m \times n}, \,
  B \in {\cal C}^{ n \times p}$  and $C \in {\cal C}^{p \times q}$ be given. Then

{\rm (a)}  \ $ (ABC)^{\dagger}
= (A^{\dagger}ABC)^{\dagger}B(ABCC^{\dagger})^{\dagger}$.

{\rm (b)} \ $ (ABC)^{\dagger} = [(AB)^{\dagger}ABC]^{\dagger}B^{\dagger}[ABC(BC)^{\dagger}]^{\dagger}$.

{\rm (c) (Cline \cite{Cl1})} \ $ (AB)^{\dagger} = (A^{\dagger}AB)^{\dagger}(ABB^{\dagger})^{\dagger}$.

{\rm (d)}\, If $ R(C ) \subseteq R[(AB)^*] $ and 
$ R(A^* ) \subseteq R(BC), $ then
$ (ABC)^{\dagger} = C^{\dagger}B^{\dagger}A^{\dagger}.$ } 

\medskip

\noindent {\bf Proof.}\, Write $ABC$ as  $ ABC = A(A^{\dagger}ABCC^{\dagger})C$. Then it is 
evident that
$$ 
r( A^{\dagger} ABCC^{\dagger}) = r(ABC), \ \ R[( ABCC^{\dagger})^{\dagger}] \subseteq R(C), \ \ 
 {\rm and}  \ \ R[(( A^{\dagger}ABC)^{\dagger})^*] \subseteq R(A^*).
$$
Thus by (8.5), we find that 
\begin{eqnarray*}
 (ABC)^{\dagger} = [  \, A( A^{\dagger} ABCC^{\dagger})C \, ]^{\dagger} =  (A^{\dagger}ABC)^{\dagger}A^{\dagger}ABCC^{\dagger}(ABCC^{\dagger})^{\dagger} = (A^{\dagger}ABC)^{\dagger}B(ABCC^{\dagger})^{\dagger},
\end{eqnarray*}
as required for Part (a). On the other hand, we can write $ ABC $ as
$ ABC = (AB)B^{\dagger}(BC)$. Applying the equality in Part(a) to it yields
$$ 
 (ABC)^{\dagger} = [ \, (AB)B^{\dagger}(BC) \, ]^{\dagger} = 
[(AB)^{\dagger}ABC]^{\dagger}B^{\dagger}[ABC(BC)^{\dagger}]^{\dagger},
$$
as required for Part (b). Let $B$ be identity matrix and replace $C$
by $B$ in the result in Part (a). Then
we have the result in Part (c). The two conditions in Part (d) are
equivalent to
$$ 
(AB)^{\dagger}ABC = C,  \ \ {\rm and } \ \ ABC(BC)^{\dagger} = A.
$$
In that case, the result in Part (b) reduces to the result in Part (d). \qquad $ \Box$

\medskip
  
In the remainder of this chapter we consider the relationship of
$ (ABC)^{\dagger}$ and the reverse order product $ C^{\dagger}B^{\dagger}A^{\dagger}$, and
present necessary and
sufficient conditions for  $(ABC)^{\dagger} =
C^{\dagger}B^{\dagger}A^{\dagger}$ to hold.
Some of the  results were  presented  by the 
first author in \cite{Ti2} and \cite{Ti4}.

\medskip
   
\noindent {\bf Lemma 8.9}\cite{Ti4}.\, {\em Suppose that $ A_1,\, A_2, \, A_3, \, B_1 $
and $ B_2$ satisfy the
the following range inclusions
$$  
R(B_i) \subseteq R(A_{i+1}),  \ \ \   and  \ \ \ R(B_i^*) \subseteq R(A_{i}^*),
  \qquad  i = 1, \ 2. \eqno (8.7) 
$$ 
Then 
$$ 
 \left[ \begin{array}{ccc} 0  &  0  & A_1  \\ 0 & A_2 & B_1  \\
  A_3 & B_2 & 0 \end{array} \right]^{\dagger} =
  \left[ \begin{array}{ccc} A_3^{\dagger}B_2A_2^{\dagger}B_1 A_1^{\dagger}
  &  -A_3^{\dagger}B_2A_2^{\dagger}  &  A_3^{\dagger} \\
   - A_2^{\dagger}B_1 A_1^{\dagger} & A_2^{\dagger} & 0 \\  A_1^{\dagger} &
   0 & 0 \end{array} \right]. \eqno (8.8)
$$ }
{\bf Proof.}\, The range inclusions in (8.7) are equivalent to 
$$
A_{i+1}A_{i+1}^{\dagger}B_i = B_i, \ \  \ {\rm and} \ \ \
 B_iA_{i}^{\dagger}A_i = B_i, \qquad  i = 1, \ 2.
$$ 
In that case, it is easy to verify that the block matrix in the right-hand
side of  (8.8) and the given block matrix in the left-hand side of  (8.8) satisfy
the four Penrose equations. Thus  (8.8) holds.  \qquad $ \Box$ 

\medskip

\noindent {\bf Lemma 8.10.}\, {\em Let $ A \in {\cal C}^{ m \times n}, \,
  B \in {\cal C}^{ n \times p}$  and $C \in {\cal C}^{p \times q}$ be given.
  Then the product $C^{\dagger}B^{\dagger}A^{\dagger}$ can be written as 
$$\displaylines{
\hspace*{1.5cm}
C^{\dagger}B^{\dagger}A^{\dagger} =[ \, I_q, \ 0, \ 0 \,]
\left[ \begin{array}{ccc} 0  &  0  & AA^*  \\ 0 & B^*BB^* & B^*A^*  \\
  C^*C & C^*B^* & 0 \end{array} \right]^{\dagger} \left[ \begin{array}{ccc} I_m \\ 0 \\ 0 
\end{array} \right] : = PJ^{\dagger}Q,  \hfill (8.9)
\cr}
$$
where  the block matrices $P, \,J$ and $Q$ satisfy 
$$ \displaylines{
\hspace*{1.5cm}
r( J ) = r(A) + r(B ) + r(C ),  \ \ \ R(QA) \subseteq R(J), \ \ \ and \ \ \
R[(AP)^*] \subseteq R(J^*). \hfill (8.10)
\cr}
$$  } 
{\bf Proof.}\, Observe that 
$$\displaylines{
\hspace*{1.5cm}  
R(B^*A^*) \subseteq R(B^*BB^*),  \ \  R(AB) \subseteq R(AA^*),  \ \ 
R(C^*B^*) \subseteq R(C^*C),  \ \  R(BC) \subseteq R(BB^*B), 
\hfill 
\cr} 
$$
as well as the three basic equalities on the Moore-Penrose inverse of a matrix
$$\displaylines{
\hspace*{1.5cm}
N^{\dagger} = N^*( N^*NN^*)^{\dagger}N^*, \qquad
N^{\dagger} = ( N^*N)^{\dagger}N^*, \qquad
N^{\dagger} = N^*(NN^*)^{\dagger}.\hfill
\cr}$$
Thus we find by (8.8) that
 \begin{eqnarray*}
\left[ \begin{array}{ccc} 0  & 0  & AA^*  \\ 0 & B^*BB^* & B^*A  \\
  C^*C & C^*B^* & 0 \end{array} \right]^{\dagger} 
&=& \left[ \begin{array}{cccc} &
 (C^*C)^{\dagger}C^*B^*(BB^*B)^{\dagger}B^*A^*(AA^*)^{\dagger} & *  & *
\\ & * & * & 0 \\ & * & 0 & 0 \end{array} \right] \\
&=& \left[ \begin{array}{cccc} & C^{\dagger}B^{\dagger}A^{\dagger} & *  & *  
\\ & * & * & 0 \\ & * & 0 & 0 \end{array} \right].
\end{eqnarray*}
Hence we have (8.9). The properties in  (8.10) are obvious.
\qquad $ \Box$ 

\medskip

\noindent {\bf Theorem 8.11.}\, {\em Let $ A \in {\cal C}^{ m \times n}, \,
  B \in {\cal C}^{ n \times p}$  and $C \in {\cal C}^{p \times q}$ be given
  and let $ M = ABC$.
Then 
$$\displaylines{
\hspace*{1.5cm}
r[ \, M -  M(C^{\dagger}B^{\dagger}A^{\dagger})M \, ] =
 r \left[ \begin{array}{ccc} -M^*  &  0  & C^*C  \\ 0 & BB^*B & BC  \\
  AA^* & AB & 0 \end{array} \right] - r(A) - r(B) - r(C) + r(M).  \hfill (8.11)
\cr
\hspace*{0cm}
In \ particular, \hfill
\cr
\hspace*{1.5cm}
C^{\dagger}B^{\dagger}A^{\dagger} \in \{ (ABC)^- \}  \Leftrightarrow r \left[ \begin{array}{ccc} -M^*  &  0  & C^*C  \\ 0 & BB^*B & BC  \\
  AA^* & AB & 0 \end{array} \right] = r(A) + r(B) + r(C) - r(M).  \hfill (8.12)
\cr}
$$ }
{\bf Proof.}\, It follows from  (1.7) and (8.9) that 
\begin{eqnarray*}
r[ \, M -  M(C^{\dagger}B^{\dagger}A^{\dagger})M \, ] & = & r( \, M -  MPJ^{\dagger}QM \, ) \\
& = & r\left[ \begin{array}{cc} J  &  QM  \\ MP & M  \end{array}
\right] - r(J)  \\
& = & r\left[ \begin{array}{cc} J -QMP  &  0  \\ 0 & M  \end{array} \right] -r(J)  \\
& = & r( \,  J -QMP \, )  + r(M) - r(J) \\
& = & r \left[ \begin{array}{ccc} -M  &  0  & AA^*  \\ 0 & B^*BB^* & B^*A^* \\
  C^*C & C^*B^* & 0 \end{array} \right] + r(M) - r(J) \\
& = & r \left[ \begin{array}{ccc} -M^*  &  0  & C^*C  \\ 0 & BB^*B & BC  \\
  AA^* & AB & 0 \end{array} \right] + r(M) - r(A) - r(B) - r(C),
\end{eqnarray*} 
as required for (8.11).   \qquad $ \Box$ 

\medskip

Applying (8.11) to the matrix product $ M  = (I_m + A)A(I_m - A) = A -  A^3$, we can find  that  
$$\displaylines{
\hspace*{2cm}
 r[ \, M - M(I_m - A)^{\dagger}A^{\dagger}( I_m  +  A)^{\dagger}M \,] = 0.
\cr}
$$
Thus $ (I_m - A)^{\dagger}A^{\dagger}( I_m  +  A)^{\dagger}$ is an inner inverse of the matrix $ A - A^3$. 
 By symmetry,  $ (I_m - A)^{\dagger}A^{\dagger}( I_m  +  A)^{\dagger}$ is also  an inner inverse of the matrix 
$ A - A^3$. We leave the verification of the rank equality to the reader. 

Applying (8.11) to the matrix product $ M  = (I_m + A)A(I_m - A) = A -  A^3$, we can find  that  
$$\displaylines{
\hspace*{2cm}
 r[ \, M - M(I_m - A)^{\dagger}A^{\dagger}( I_m  +  A)^{\dagger}M \,] = 0.
\cr}
$$
Thus $ (I_m - A)^{\dagger}A^{\dagger}( I_m  +  A)^{\dagger}$ is an inner inverse of the matrix $ A - A^3$. 
 By symmetry,  $ (I_m - A)^{\dagger}A^{\dagger}( I_m  +  A)^{\dagger}$ is also  an inner inverse of the matrix 
$ A - A^3$. We leave the verification of the rank equality to the reader. 

In general, for any three polynomials $ p_1(\lambda), \ p_2(\lambda)$ and $ p_3(\lambda)$  
without common roots and a square matrix $A$,  we let $ p(A) = p_1(A)p_2(A)p_3(A)$. Then 
  one can find by (8.11) and (1.17) the following   
$$
 r[ \, p(A) - p(A)p_3^{\dagger}(A)p_2^{\dagger}(A)p_1^{\dagger}(A)p(A) \, ] = 0. 
$$
Thus $p_3^{\dagger}(A)p_2^{\dagger}(A)p_1^{\dagger}(A)$ is an inner inverse 
 inverse of the matrix product $ p_1(A)p_2(A)p_3(A)$.

\medskip

\noindent {\bf Theorem 8.12.}\, {\em Let $ A \in {\cal C}^{ m \times n}, \,
  B \in {\cal C}^{ n \times p}$  and $C \in {\cal C}^{p \times q}$ be given
  and let $ M = ABC$.
Then 
$$\displaylines{
\hspace*{1.5cm}
r( \, M^{\dagger} -  C^{\dagger}B^{\dagger}A^{\dagger} \, ) =
 r \left[ \begin{array}{ccc} -MM^*M  &  0  & MC^*C  \\ 0 & BB^*B & BC  \\
  AA^*M & AB & 0 \end{array} \right] - r(B) - r(M).  \hfill (8.13)
\cr
\hspace*{0cm}
In \ particular, \hfill
\cr
\hspace*{1.5cm}
(ABC)^{\dagger} = C^{\dagger}B^{\dagger}A^{\dagger} \Leftrightarrow 
 r\left[ \begin{array}{ccc} -MM^*M  &  0  & MC^*C  \\ 0 & BB^*B & BC  \\
  AA^*M & AB & 0 \end{array} \right] = r(B) + r(ABC).  \hfill (8.14)
\cr}
$$ }
{\bf Proof.}\, Notice that
$$\displaylines{
\hspace*{1.5cm}
 C^{\dagger}CM^{\dagger}AA^{\dagger} = M^{\dagger} \ \ \ {\rm and}  \ \ 
 C^{\dagger}C( C^{\dagger}B^{\dagger}A^{\dagger}) AA^{\dagger}
 = C^{\dagger}B^{\dagger}A^{\dagger}. \hfill
\cr}
$$
We first get the following 
\begin{eqnarray*}
r( \, M^{\dagger} -  C^{\dagger}B^{\dagger}A^{\dagger} \, ) & = &
r( \, CM^{\dagger}A -  CC^{\dagger}B^{\dagger}A^{\dagger}A \, ) \\
& = & r( \, CM^{\dagger}A -  CPJ^{\dagger}QA \, ) \\
& = & r[ \, C M^*(M^*MM^*)^{\dagger}M^*A -  CPJ^{\dagger}QA \, ] \\
& = & r \left( [\, CM^*, \ CP\, ] \left[ \begin{array}{cc} -M^*MM^* & 0 \\ 0 & J \end{array} 
\right]^{\dagger} \left[ \begin{array}{c} M^*A \\ QA  \end{array} \right] \right).
\end{eqnarray*} 
Observe from (8.10) that 
$$
R \left[ \begin{array}{c} M^*A \\ QA  \end{array} \right] \subseteq 
R \left[ \begin{array}{cc} -M^*MM^* & 0 \\ 0 & J \end{array} 
\right] \ \ \ {\rm and} \ \ \ 
R (\, [\, CM^*, \ CP\, ]^* \,)  \subseteq 
R\left( \left[ \begin{array}{cc} -M^*MM^* & 0 \\ 0 & J \end{array} 
\right]^* \right).
$$
Thus we find  by  (1.7) that
\begin{eqnarray*}
r \lefteqn{ \left( [\, CM^*, \ CP\, ] \left[ \begin{array}{cc} -M^*MM^* & 0 \\ 0 & J \end{array} 
\right]^{\dagger} \left[ \begin{array}{c} M^*A \\ QA  \end{array} \right] \right) } \\
& = & r\left[ \begin{array}{ccc} -MM^*M  &  0  & M^*A  \\ 0 & J & QA  \\
  CM^* & CP & 0 \end{array} \right] - r\left[ \begin{array}{cc} -M^*MM^* & 0 \\ 0 & J  
\end{array} \right] \\
& = & r\left[ \begin{array}{ccccc} -MM^*M  &  0 & 0 &  0 & M^*A  \\ 0 & 0 & 0 & AA^* & A  \\
 0 & 0 & B^*BB^*  & B^*A^* & 0 \\ 0 & C^*C & C^*B^* & 0 & 0  \\ CM^* & C & 0 & 0 & 0 \end{array} \right] 
- r(M) - r(J) \\
& = & r\left[ \begin{array}{ccccc} -MM^*M  & 0 & 0 &  0 -M^*AA^* & 0  \\ 0 & 0 & 0 & 0 & A  \\
 0 & 0 & B^*BB^*  & B^*A^* & 0 \\ -C^*CM^* & 0  & C^*B^* & 0 & 0  \\ 0 & C & 0 & 0 & 0 
\end{array} \right] - r(M) - r(J) \\
& = & r\left[ \begin{array}{ccc} -MM^*M  & 0 & -M^*AA^*  \\ 0 & B^*BB^*  & B^*A^*
 \\ -C^*CM^* &  C^*B^* & 0 
\end{array} \right] - r(M) - r(B).
\end{eqnarray*} 
The results in (8.13) and (8.14) follow from it.  \qquad $ \Box$

\medskip

\noindent {\bf Corollary 8.13.}\, {\em Let $ A \in {\cal C}^{ m \times n}, \,
  B \in {\cal C}^{ n \times p}$  and $C \in {\cal C}^{p \times q}$ be given$,$ and let $ M = ABC$. If 
$$ \displaylines{
\hspace*{1.5cm}
R(B ) \subseteq R(A^*)  \ \ and \ \ R(B^* ) \subseteq R(C),  \hfill (8.15)
\cr
\hspace*{0cm}
then \hfill
\cr
\hspace*{1.5cm}
r( \, M^{\dagger} -  C^{\dagger}B^{\dagger}A^{\dagger} \, ) =
 r \left[ \begin{array}{ccc} B \\ BCC^* \end{array} \right] 
 + r[ \,  B, \ A^*AB \, ] - 2r(B).  \hfill (8.16)
\cr
In \ particular, \hfill
\cr
\hspace*{1.5cm}
(ABC)^{\dagger} = C^{\dagger}B^{\dagger}A^{\dagger} \Leftrightarrow 
 R(A^*AB ) \subseteq R(B)  \ \ and \ \ R[(BCC^*)^*] \subseteq R(B^*). \hfill (8.17)
\cr}
$$ }
{\bf Proof.}\, Eq.\,(8.15) is equivalent to $ A^{\dagger}AB = B$ and 
$ BCC^{\dagger} = B$. Thus we can reduce  (8.13) by block elementary operations  to    
\begin{eqnarray*}
r( \, M^{\dagger} -  C^{\dagger}B^{\dagger}A^{\dagger} \, ) & =& 
 r \left[ \begin{array}{ccc} -BCM^*AB  &  0  & BCC^* \\ 0 & BB^*B & B  \\
  A^*AB & B & 0 \end{array} \right] - 2r(B) \\
& =&  r \left[ \begin{array}{ccc} 0  & BCC^*B^* & BCC^* \\ 0 & BB^*B & B  \\
  A^*AB & B & 0 \end{array} \right] - 2r(B) \\
& =&  r \left[ \begin{array}{ccc} 0  & 0 & BCC^*B \\ 0 & 0 & B  \\
  A^*AB & B & 0 \end{array} \right] - 2r(B) \\
 &= & r \left[ \begin{array}{ccc} B \\ BCC^* \end{array} \right] 
+ r[ \,  B, \ A^*AB \, ] - 2r(B).  \qquad \Box
\end{eqnarray*} 

\noindent  {\bf Corollary 8.14.}\, {\em Let $ B \in {\cal C}^{ m \times n}$
be given$,$ $ A \in {\cal C}^{m \times m}$ and  $ C \in {\cal C}^{n\times n}$ be two
 invertible matrices. Let $ M = ABC.$  Then
$$\displaylines{
\hspace*{1.5cm}
  r[ \, (ABC)^{\dagger} - C^{-1}B^{\dagger}A^{-1} \, ] 
= r\left[ \begin{array}{c} B \\ BCC^* \end{array} \right]
  + r[ \,  B, \ A^*AB \, ] - 2r(B), \hfill (8.18)
\cr
\hspace*{0cm}
and \hfill
\cr
\hspace*{1.5cm}
r[ \, (ABC)^{\dagger} - C^{-1}B^{\dagger}A^{-1} \, ]
  = r\left[ \begin{array}{c} M \\ MC^*C \end{array} \right] +
 r[ \,  M, \ AA^*M \,] - 2r(M). \hfill (8.19)
\cr
\hspace*{0cm}
In \ particular, \hfill
\cr
\hspace*{1.5cm}
 (ABC)^{\dagger} = C^{-1}B^{\dagger}A^{-1} \Leftrightarrow R(AA^*B )
 = R(B) \ \ and \ \  R( CC^*B^* ) = R( B^* ), \hfill (8.20)
\cr
\hspace*{0cm}
and \hfill
\cr
\hspace*{1.5cm}
 (ABC)^{\dagger} = C^{-1}B^{\dagger}A^{-1}  \Leftrightarrow  R(AA^*M ) = R(M) \ \ and \ \ 
  R(C^*CM^* ) = R( M^*).  \hfill (8.21) 
\cr} 
$$ }
{\bf Proof.}\, Follows immediately from Corollary 8.13.   \qquad $ \Box$

\medskip

\noindent {\bf Theorem 8.15.}\, {\em Let $ B \in {\cal C}^{m \times n}$ be given$,$  $ A \in {\cal C}^{m \times m }$ and $ C \in {\cal C}^{n \times n}$ be two invertible matrices. Let $ M = ABC$. Then
 
{\rm (a)} \ $ r( \, MM^{\dagger} - ABB^{\dagger}A^{-1} \,)
= r[\, B, \ A^*AB \, ] - r(B).$

{\rm (b)} \  $ r( \, M^{\dagger}M - C^{-1}B^{\dagger}BC \,) 
= r\left[ \begin{array}{c} B \\ BCC^*  \end{array} \right] - r(B).$

{\rm (c)} \ $ r( \, M^{\dagger} - C^{-1}B^{\dagger}A^{-1} \, ) = 
r( \, MM^{\dagger} - ABB^{\dagger}A^{-1} \,) + r( \, M^{\dagger}M - C^{-1}B^{\dagger}BC \,).$ \\
In  particular$,$

{\rm (d)} \ $ MM^{\dagger} = ABB^{\dagger}A^{-1}
\Leftrightarrow R(A^*AB) = R(B).$

{\rm (e)} \ $ M^{\dagger}M = C^{-1}B^{\dagger}BC
\Leftrightarrow R(CC^*B^*) = R(B^*).$

{\rm (f)} \ $M^{\dagger} = C^{-1}B^{\dagger}A^{-1} \Leftrightarrow MM^{\dagger} = ABB^{\dagger}A^{-1}
 \Leftrightarrow M^{\dagger}M = C^{-1}B^{\dagger}BC.$
}  

\medskip

\noindent  {\bf Proof.}\, Observe that 
$$
r(\, MM^{\dagger} - ABB^{\dagger} A^{-1} \, ) 
= r(\, MM^{\dagger}A - ABB^{\dagger} \, ) \ \ {\rm and}  \ \ r(\, M^{\dagger}M - C^{-1}B^{\dagger}BC \, ) 
= r(\, CM^{\dagger}M - B^{\dagger}BC \, ).
$$
Applying (7.4) to both of them yields Parts (a) and (b). Contrasting (8.18) with Parts (a) and (b) yields Part 
(c).  \qquad  $\Box$ 

\medskip

Finally we present author interesting result on the Moore-Penrose inverse of a triple matrix product.

\medskip

\noindent {\bf Corollary 8.16.}\, {\em Let $ A \in {\cal C}^{ m \times n}, \
  B \in {\cal C}^{ n \times p}$  and $C \in {\cal C}^{p \times q}$ be given$,$  and suppose that 
$ R(B ) \subseteq R(A^*)$  and $R(B^* ) \subseteq R(C).$   Then
$$
r[ \, (ABC)^{\dagger} - (\, I_q  - (C^{\dagger}F_B)(C^{\dagger}F_B)^{\dagger} \,) 
C^{\dagger}B^{\dagger}A^{\dagger}(\, I_m - (E_BA^{\dagger})^{\dagger}(E_BA^{\dagger}) \,) \,] 
= m + q - r(A) - r(C), \eqno (8.22)
$$
where $E_B = I_n -BB^{\dagger}$ and $F_B = I_p -B^{\dagger}B$. In particular$,$ the equality 
$$
\displaylines{
\hspace*{1.5cm}
 (ABC)^{\dagger} = [\, I_q  - (C^{\dagger}F_B)(C^{\dagger}F_B)^{\dagger} \,] 
C^{\dagger}B^{\dagger}A^{\dagger}[\, I_m - (E_BA^{\dagger})^{\dagger}(E_BA^{\dagger}) \,] \hfill (8.23)
\cr}
$$
holds if and only if $ r(A) = m$ and $ r(C) = q $. Thus if both $ A $ and $ C$ are nonsingular matrices$,$ then 
the identity     
$$
\displaylines{
\hspace*{1.5cm}
 (ABC)^{\dagger} = [\, I_q  - (C^{-1}F_B)(C^{-1}F_B)^{\dagger} \,] 
C^{-1}B^{\dagger}A^{-1}[\, I_m - (E_BA^{-1})^{\dagger}(E_BA^{-1}) \,] \hfill (8.24)
\cr}
$$
holds.}

\medskip

\noindent {\bf Proof.}\, Let $ M = ABC$ and $ N = [\, I_q  - (C^{\dagger}F_B)(C^{\dagger}F_B)^{\dagger} \,] 
C^{\dagger}B^{\dagger}A^{\dagger}[\, I_m - (E_BA^{\dagger})^{\dagger}(E_BA^{\dagger}) \,].$ 
 Then it is easy to verify that under $ R(B ) \subseteq R(A^*)$  and $R(B^* ) \subseteq R(C)$, $N$ is an outer 
inverse of $M$. Hence by (5.1) we get 
$$
\displaylines{
\hspace*{1.5cm}
r( \, M^{\dagger} -  N \, ) =
 r \left[ \begin{array}{c}  M^{\dagger}  \\ N \end{array} \right] 
 + r[ \,   M^{\dagger}, \ N  \, ] - r(M) - r(N).  \hfill (8.25)
\cr}
$$  
Simplifying the ranks of the matrices  in (8.25) by (8.23) and (1.2)---(1.4), we can eventually get
$$
\displaylines{
\hspace*{0.5cm}
r \left[ \begin{array}{c}  M^{\dagger}  \\ N \end{array} \right] = m + r(B) - r(C),  \ \ \ r[ \,   M^{\dagger}, \ N  \, ] = q  + r(A) - r(C), \ \ \ r(M) = r(B),   \ \ \ r(N ) = r(B). \hfill   
\cr}
$$
The tedious processes are omitted here. Putting them in (8.25) we have (8.22), and then (8.23) and (8.24). \qquad 
$ \Box$ 

\medskip

It is expected that the identity (8.24) can help to establish various equalities for Moore-Penrose inverses of 
 block matrices.

\markboth{YONGGE  TIAN }
{9. MOORE-PENROSE INVERSES OF BLOCK MATRICES}

\chapter{Moore-Penrose inverses of  block matrices}

\noindent In this chapter we establish some rank equalities related to  factorizations
 of $ 2 \times 2$ block matrices and then deduce from them various expressions of Moore-Penrose inverses 
 for $ 2 \times 2$  block matrices, as well as for  $ m \times n$  block matrices. Some of the results 
in this chapter appear in the author's recent paper \cite{Ti5}. In fact, any  $ 2 \times 2$ block matrix can 
simply factor as various types of products of block matrices. In that case, applying the rank equalities in Chapter 8 to 
them one can establish many new rank equalities related to the block matrix, and consequently, derive from them 
various expressions for the Moore-Penrose inverse of the block matrix. We begin this work first  with a bordered 
matrix.             

\medskip
   
\noindent  {\bf Theorem 9.1.}\, {\em    Let  $ M = \left[ \begin{array}{cc} A & B \\ C  & 0  \end{array} \right]$ 
 be a given bordered matrix over the field of  complex numbers$,$  where $ A \in {\cal C}^{ m \times n}, \, B \in
 {\cal C}^{ m \times k}$ and $C \in {\cal C}^{ l \times n},$ and factor $ M$  as 
$$ 
 M = \left[ \begin{array}{cc} A & B \\ C  & 0 \end{array} \right] 
 = \left[ \begin{array}{cc} I_m  & E_BAC^{\dagger}  \\ 0  & I_l
 \end{array} \right] \left[ \begin{array}{cc} E_BAF_C  & B \\ C  & 0
 \end{array} \right] \left[ \begin{array}{cc} I_n  &  0 \\ B^{\dagger}A & I_k \end{array} \right] : = PNQ,  \eqno (9.1) 
 $$
where $E_B = I_m - BB^{\dagger}$ and $F_C = I_n - C^{\dagger}C$. Then

{\rm (a)}\, The rank of $ M^{\dagger} - Q^{-1}N^{\dagger}P^{-1}$ satisfies the equality   
$$\displaylines{
\hspace*{1.5cm}
 r( \, M^{\dagger} - Q^{-1}N^{\dagger}P^{-1} \, )
  = r\left[ \begin{array}{c} A \\ C \end{array} \right] + r[ \, A, \ B \, ]
  + r(B) + r(C) - 2r(M). \hfill (9.2)
\cr}
 $$

{\rm (b)}\, The following four statements are equivalent$:$

\hspace*{0.5cm}{\rm (1)}\ The Moore-Penrose inverse of $ M $ can  be expressed as 
$ M^{\dagger} = Q^{-1}N^{\dagger}P^{-1},$ that is$,$  
 $$ \displaylines{
\hspace*{2cm} 
  M^{\dagger} =  \left[ \begin{array}{cc} (E_BAF_C)^{\dagger} & C^{\dagger} -
  (E_BAF_C)^{\dagger}AC^{\dagger}
  \\ B^{\dagger} -  B^{\dagger}A(E_BAF_C)^{\dagger}  &
  -B^{\dagger}AC^{\dagger} +
 B^{\dagger}A(E_BAF_C)^{\dagger}AC^{\dagger}  \end{array} \right].
 \hfill (9.3)
\cr} 
$$ 

\hspace*{0.5cm}{\rm (2)} \ $[\, I_n, \ 0 \,]M^{\dagger}\left[ \begin{array}{c} I_m \\ 0  \end{array}
   \right]=(E_BAF_C)^{\dagger}.$

\hspace*{0.5cm}{\rm (3)} \ $A, \ B$ and $ C $ satisfy the  rank additivity
 condition
$$ \displaylines{
\hspace*{2cm}
 r( M ) =  r \left[ \begin{array}{c} A \\ C \end{array} \right] + r(B)
 =  r[ \, A, \ B \, ] + r(C). \hfill (9.4)
\cr}
$$

\hspace*{0.5cm}{\rm (4)}\, The two conditions hold
$$\displaylines{
\hspace*{2cm}
R \left[ \begin{array}{c} A \\ C \end{array} \right] \cap
 R\left[ \begin{array}{c} B \\ 0 \end{array} \right] = \{0\} \ \  and
  \ \ R ([ \, A, \ B \, ]^*) \cap R( [\, C, \ 0 \,]^*) = \{0\}.
  \hfill (9.5)
\cr}
$$ }
{\bf Proof.}\, It follows first by  (9.1) and (8.19) that
$$\displaylines{
\hspace*{1.5cm}
 r( \, M^{\dagger} - Q^{-1}N^{\dagger}P^{-1} \, )
  = r\left[ \begin{array}{c} M \\ MQ^*Q \end{array} \right] +
  r[ \, M, \ PP^*M \, ] - 2r(M). \hfill (9.6)
\cr}
$$
The ranks of the two block matrices in  (9.6) can simplify to 
\begin{eqnarray*}
r[\, M, \ PP^*M \, ]
& = & r \left[ \, \left[ \begin{array}{cc} A  & B \\ C & 0 \end{array} \right],
\ \left[ \begin{array}{cc} I_m  & E_BAC^{\dagger} \\ 0 & I_l \end{array}
\right] \left[ \begin{array}{cc} I_m & 0 \\ (E_BAC^{\dagger})^*  & I_l
 \end{array} \right]\left[ \begin{array}{cc} A  & B \\ C & 0 \end{array}
 \right] \, \right] \\
& = & r \left[ \begin{array}{cccc} A  & B & A + (E_BAC^{\dagger})
(E_BAC^{\dagger})^*A + E_BAC^{\dagger}C  & B \\
C & 0 & C + (E_BAC^{\dagger})^*A  & 0 \end{array}
\right] \\
& = & r \left[ \begin{array}{ccc} A & B & AC^{\dagger}
(E_BAC^{\dagger})^*A + AC^{\dagger}C  \\
C & 0 &(E_BAC^{\dagger})^*A \end{array}
\right] \\
& = & r \left[ \begin{array}{ccc} A & B & AC^{\dagger}C  \\
C & 0 & 0 \end{array}
\right] \\
& = & r \left[ \begin{array}{ccc} A & B & 0  \\
C & 0 & -C \end{array}
\right]  = r[\, A, \ B \,] + r(C),
\end{eqnarray*}
and
\begin{eqnarray*}
r\left[ \begin{array}{c} M \\ MQ^*Q \end{array} \right]
& = & r\left[ \begin{array}{c}
 \left[ \begin{array}{cc} A  & B \\ C & 0 \end{array} \right]
\\
 \left[ \begin{array}{cc} A  & B \\ C & 0 \end{array} \right]
\left[ \begin{array}{cc} I_n  & (B^{\dagger}A)^* \\ 0 & I_k \end{array} \right]
 \left[ \begin{array}{cc} I_n  & 0 \\ B^{\dagger}A & I_k \end{array} \right]
    \end{array} \right] \\
& = & r\left[ \begin{array}{cc} A  & B \\ C & 0 \\ 
A + A(B^{\dagger}A)^*(B^{\dagger}A) + BB^{\dagger}A  & B +
A(B^{\dagger}A)^* \\
 C + C(B^{\dagger}A)^*B^{\dagger}A & C(B^{\dagger}A)^* \end{array}
 \right] \\
& = & r \left[ \begin{array}{cc} A  & B \\ C & 0
\\  A + A(B^{\dagger}A)^*(B^{\dagger}A)  &  A(B^{\dagger}A)^* \\
 C(B^{\dagger}A)^*B^{\dagger}A & C(B^{\dagger}A)^* \end{array}
 \right] \\
& = & r \left[ \begin{array}{cc} A  & B \\ C & 0 \\  A  &  0 \\
 0  &  0 \end{array}
 \right] = r\left[ \begin{array}{c} A \\ C  \end{array} \right]
 + r(B). 
\end{eqnarray*}
Putting both of them in  (9.6) yields  (9.2). Notice that
$ (E_BAF_C)^{\dagger}B^{\dagger} = 0$ and $ C^{\dagger}(E_BAF_C)^{\dagger}
 = 0$ always hold. Then it is easy to verify that
$$
\displaylines{
\hspace*{1.5cm}
N^{\dagger} = \left[ \begin{array}{cc} E_BAF_C & B \\ C   & 0 \end{array}
 \right]^{\dagger} = \left[ \begin{array}{cc} (E_BAF_C)^{\dagger}
  & C^{\dagger}  \\ B^{\dagger} & 0 \end{array} \right]. \hfill
\cr}
$$
Putting it in (9.2), we get 
$$ 
\displaylines{
\hspace*{1.5cm} 
  Q^{-1}N^{\dagger}P^{-1} =  \left[ \begin{array}{cc} (E_BAF_C)^{\dagger} & C^{\dagger} -
  (E_BAF_C)^{\dagger}AC^{\dagger}
  \\ B^{\dagger} -  B^{\dagger}A(E_BAF_C)^{\dagger}  &
  -B^{\dagger}AC^{\dagger} +
 B^{\dagger}A(E_BAF_C)^{\dagger}AC^{\dagger}  \end{array} \right].
 \hfill 
\cr} 
$$ 
The equivalence of the statements (1) and (3) in Part (b) follows from  (9.2). The equivalence of the statements 
(2) and (3) in Part (b) comes from  (7.29).  The equivalence of the statements 
(3) and (4) in Part (b) is obvious.  \qquad  $ \Box$ 

\medskip

The expression (9.3) for $M^{\dagger}$ is well known when $M$ satisfies (9.4) (see, i.e. \cite{Mi3} and \cite{RM}). 
The rank equality (9.2) further reveals a fact that (9.3) is not only sufficient but also necessary. Various 
consequences can be derived from Theorem 9.1 when the matrix $M$ in it satisfies some more restrictions. 
Here we only present one that is well known. 

\medskip 

\noindent {\bf  Corollary 9.2.}\, {\em Let $ A \in {\cal C}^{ m \times n},
\, B \in {\cal C}^{ m \times k}$ and $C \in {\cal C}^{ l \times n} $ be
given. If $ R(A) \cap R(B) = \{ 0 \} $ and $R(A^*) \cap R(C^*) = \{ 0 \},$  then 
$$
\displaylines{
\hspace*{1.5cm}
\left[ \begin{array}{cc} A & B \\ C   & 0 \end{array}
 \right]^{\dagger}
 = \left[ \begin{array}{cc} (E_BAF_C)^{\dagger} & C^{\dagger} -
  (E_BAF_C)^{\dagger}AC^{\dagger}
  \\ B^{\dagger} -  B^{\dagger}A(E_BAF_C)^{\dagger}  & 0  \end{array}
\right]. \hfill (9.7)
\cr}
$$}
{\bf Proof.}\, Under $ R(A) \cap R(B) = \{ 0 \} $ and $R(A^*) \cap R(C^*) = \{ 0 \},$
  the rank equality (9.4) naturally holds.
In that case, we know by Theorem 7.8 that $A(E_BAF_C)^{\dagger}A = A$. Thus
 (9.3) reduces to  (9.7). \qquad $ \Box$ 

\medskip

Besides the factorization (9.1), we can generally factor $ M $ as 
$$
\displaylines{
\hspace*{1.5cm}
M = \left[ \begin{array}{cc} A & B \\ C & 0 \end{array} \right]
  = \left[ \begin{array}{cc} I_m  & Y \\ 0 & I_l  \end{array} \right]
 \left[ \begin{array}{cc} A - BX - YC  & B \\ C   & 0 \end{array} \right]
 \left[ \begin{array}{cc} I_n  & 0 \\ X & I_k \end{array}
 \right] : = PNQ,  \hfill (9.8)
\cr}
$$
where $ X $ and $Y$  are arbitrary and $ P $ and $ Q $ are nonsingular. Clearly (9.1) is a special case of 
(9.8).  In that case, applying the first equality in (8.5) to (9.8) we obtain the following. 

\medskip 

\noindent  {\bf Theorem 9.3.}\, {\em Let $ A \in {\cal C}^{ m \times n}, \,
B \in {\cal C}^{ m \times k}$ and $C \in {\cal C}^{ l \times n} $ be given.
Then
$$
\displaylines{
\hspace*{1.5cm}
\left[ \begin{array}{cc} A & B \\ C   & 0 \end{array}\right]^{\dagger}
 = \left[ \begin{array}{cc} A-YC & B \\ C & 0 \end{array}
\right]^{\dagger}
\left[ \begin{array}{cc} A -BX - YC & B \\ C   & 0 \end{array}
\right]\left[ \begin{array}{cc} A-BX & B \\ C & 0 \end{array}
\right]^{\dagger}, \hfill (9.9)
\cr}
$$
where $ X $ and $Y$  are arbitrary. In particular$,$
$$
\displaylines{
\hspace*{1.5cm}
\left[ \begin{array}{cc} A & B \\ C & 0 \end{array}\right]^{\dagger}
 = \left[ \,  [\, AF_C, \ B \,]^{\dagger}, \ \left[ \begin{array}{c}
 C^{\dagger} \\ 0 \end{array} \right] \right]
\left[ \begin{array}{cc} A - BB^{\dagger}A -
 AC^{\dagger}C & B \\ C
& 0 \end{array} \right] \left[ \begin{array}{c}
\left[ \begin{array}{c} E_BA \\ C \end{array} \right]^{\dagger} 
\\ \left[ \, B^{\dagger}, \ 0 \,\right] \end{array} \right].   \hfill (9.10)
\cr}
$$}
{\bf Proof.}\, Under (9.8) we have by (8.5) that $ M =  (PNQ)^{\dagger} =(NQ)^{\dagger}N(PN)^{\dagger}.$
Written in an explicit form, it is (9.9). Now let $ X = B^{\dagger}A$ and
$ Y = AC^{\dagger}$ in  (9.9). Then  (9.9) becomes
$$\displaylines{
\hspace*{1.5cm}
\left[ \begin{array}{cc} A & B \\ C & 0 \end{array}\right]^{\dagger}
 = \left[ \begin{array}{cc} AF_C  & B \\ C   & 0 \end{array}
 \right]^{\dagger} \left[ \begin{array}{cc} A - BB^{\dagger}A -
 AC^{\dagger}C  & B \\ C  & 0 \end{array} \right]
 \left[ \begin{array}{cc} E_BA  & B \\ C   & 0 \end{array} \right]^{\dagger}. \hfill
\cr}
$$
Note that 
$$
\displaylines{
\hspace*{1.5cm}
[\, AF_C ,\ B \, ][\, C, \ 0 \,]^* = 0, \ \ \ {\rm and }  \ \ \
\left[ \begin{array}{c} B \\ 0  \end{array}\right]^*
\left[ \begin{array}{c} E_BA \\ C \end{array}\right] = 0. \hfill
\cr}
$$
Then it follows by Theorem 7.9(e) and (f) that
$$
\displaylines{
\hspace*{1.5cm}
\left[ \begin{array}{cc} A & B \\ C & 0 \end{array}\right]^{\dagger}
  = [\, [\, AF_C ,\ B \, ]^{\dagger}, \ [\, C, \ 0 \,]^{\dagger} \,], \ \ {\rm and}  \ \  
\left[ \begin{array}{cc} E_BA  & B \\ C  & 0 \end{array} \right]^{\dagger}
  =  \left[ \begin{array}{c}
   \left[ \begin{array}{c} E_BA \\ C  \end{array} \right]^{\dagger}
   \\ \left[ \begin{array}{c} B  \\ 0 \end{array} \right]^{\dagger} \end{array} \right]. \hfill
\cr}
$$
Thus we have (9.10). \qquad $\Box$ 

\medskip

Eq.\,(9.10) manifests that the Moore-Penrose inverse of a bordered matrix can be intermediately determined by 
 the Moore-Penrose inverse of  $B$,  $C$,  $[\, AF_C ,\ B \, ]$ and $\left[ \begin{array}{c} E_BA \\ C  \end{array} 
\right].$  Observe that 
$$
\displaylines{
\hspace*{0.5cm}
[ \, AF_C, \ B \, ]^{\dagger} = [ \, AF_C, \ B \, ]^*\left( [ \, AF_C, \ B \, ]
[ \, AF_C, \ B \, ]^* \right)^{\dagger} = \left[ \begin{array}{c} (AF_C)^* \left[\, (AF_C)(AF_C)^* + BB^* \, 
\right]^{\dagger} \\
  B^* \left[\, (AF_C)(AF_C)^* + BB^* \, \right]^{\dagger}  \end{array} \right], \hfill
\cr}
$$
\begin{eqnarray*}
\left[ \begin{array}{c} E_BA \\ C \end{array} \right]^{\dagger} &=& \left( \left[ \begin{array}{c} E_BA \\ C 
\end{array} \right]^* \left[ \begin{array}{c} E_BA \\  C \end{array} \right] \right)^{\dagger} 
\left[ \begin{array}{c} E_BA \\ C \end{array} \right]^*  \\
& = &[\, ( \, (E_BA)^*(E_BA) + C^*C \,)^{\dagger} (E_BA)^* ,  \   (\,  (E_BA)^* (E_BA) +
 C^*C \, )^{\dagger}C^*  \, ]. 
\end{eqnarray*}
Inserting them in (9.10) we get 
$$
\displaylines{
\hspace*{1cm}
\left[ \begin{array}{cc} A & B \\ C & 0 \end{array}\right]^{\dagger}
 = \left[ \begin{array}{cc} (AF_C)^*\left[\, (AF_C)(AF_C)^* + BB^* \,\right]^{\dagger} &  C^{\dagger} \\
  B^* \left[\, (AF_C)(AF_C)^* + BB^* \,\right]^{\dagger} &  0  \end{array} \right] 
\left[ \begin{array}{cc} A - BB^{\dagger}A -  AC^{\dagger}C & B \\ C
& 0 \end{array} \right]  \hfill
\cr
\hspace*{3cm}
 \times  \left[ \begin{array}{cc} \left[\, (E_BA)^*(E_BA) + C^*C \,\right]^{\dagger} (E_BA)^* &  
\left[\,  (E_BA)^* (E_BA) +
 C^*C \,\right]^{\dagger}C^*  \\  B^{\dagger} & 0  \end{array} \right], \hfill 
\cr}
$$
which could be regarded as a general expression for the Moore-Penrose inverse a boredered matrix 
when no restriction is posed on it. Moreover, this expression reveals another interesting that the Moore-Penrose 
inverse a boredered matrix can factor as  a product of three boredered matrices although the 
Moore-Penrose inverse of the boredered matrix is not boredered in general. 

Another well-known factorization for a bordered matrix is  
$$
\displaylines{
\hspace*{2cm}
 M = \left[ \begin{array}{cc} A & B \\ C  & 0 \end{array} \right]  = \left[ \begin{array}{cc} I_m  & 0 \\ CA^{\dagger}  & I_l
 \end{array} \right] \left[ \begin{array}{cc} A & E_AB \\ CF_A  & -CA^{\dagger}B
 \end{array} \right] \left[ \begin{array}{cc} I_n  &
 A^{\dagger}B  \\ 0 & I_k \end{array} \right].  \hfill  (9.11) 
\cr}
$$
But it can be considered as a special case of the Schur factorization of a $ 2 \times 2 $ block matrix, 
$$
\displaylines{
\hspace*{2cm}
 M = \left[ \begin{array}{cc} A & B \\ C  & D \end{array} \right]  = \left[ \begin{array}{cc} I_m  & 0 \\ CA^{\dagger}  & I_l
 \end{array} \right] \left[ \begin{array}{cc} A & E_AB \\ CF_A  & S_A
 \end{array} \right] \left[ \begin{array}{cc} I_n  &
 A^{\dagger}B  \\ 0 & I_k \end{array} \right] : = PNQ,  \hfill (9.12) 
\cr}
 $$
where $S_A = D - CA^{\dagger}B$. We next present a rank equality related to (9.12) and derive its consequences. 

\medskip

\noindent {\bf Theorem 9.4.}\, {\em Let $ M$ be given by {\rm (9.12)}$,$ where $ A \in {\cal C}^{ m \times n},
\, B \in {\cal C}^{ m \times k}, \, C \in {\cal C}^{ l \times n}$  and
 $D \in {\cal C}^{ l \times k}$. Then the rank of  $ M^{\dagger} - Q^{-1}N^{\dagger}P^{-1}$ satisfies 
the equality  
 $$\displaylines{
\hspace*{2cm}
 r( \, M^{\dagger} - Q^{-1}N^{\dagger}P^{-1} \, )
  = r\left[ \begin{array}{cc} A & 0 \\ 0 & C  \\ B & D
   \end{array} \right] + r\left[ \begin{array}{ccc} A & 0 & B \\ 0 & C & D
   \end{array} \right] - 2r(M). \hfill (9.13) 
\cr}
$$ 
In particular$,$ the Moore-Penrose inverse of $ M $ in  {\rm (9.12)} can be
 expressed as $ M^{\dagger} = Q^{-1}N^{\dagger}P^{-1},$ that is$,$  
$$\displaylines{
\hspace*{2cm}
M^{\dagger}  =  \left[ \begin{array}{cc} I_n
   & -A^{\dagger}B  \\ 0  & I_k \end{array} \right]
   \left[ \begin{array}{cc} A & E_AB \\ CF_A  & S_A \end{array}
   \right]^{\dagger}
 \left[ \begin{array}{cc} I_m  & 0 \\ -CA^{\dagger}  & I_l \end{array}\right] 
\hfill (9.14)
\cr}
$$
holds if and only if  $A, \ B, \ C $ and $ D $ satisfy
$$\displaylines{
\hspace*{1.5cm}
  r\left[ \begin{array}{cc} A & 0 \\ 0 & C  \\ B & D  \end{array} \right]
  = r(M)  \ \ and \ \
  r\left[ \begin{array}{ccc} A & 0  & B \\ 0 & C & D \end{array} \right]
   = r(M), \hfill (9.15) 
 \cr
\hspace*{0cm}
or \ equivalently \hfill
\cr
\hspace*{2cm}
  R \left[ \begin{array}{c} A \\ 0 \end{array} \right] \subseteq R
  \left[ \begin{array}{cc} A & B \\ C  & D \end{array} \right]  \ \ \ and
  \ \ \
  R \left[ \begin{array}{c} A^* \\ 0 \end{array} \right] \subseteq
  R\left[ \begin{array}{cc} A^* & C^* \\ B^*  & D^* \end{array} \right] .
  \hfill (9.16)
 \cr}
$$
} 
{\bf Proof.}\, It follows by  (9.12) and (8.19) that
$$\displaylines{
\hspace*{2cm}
 r( \, M^{\dagger} - Q^{-1}N^{\dagger}P^{-1} \, )
  = r\left[ \begin{array}{c} M \\ MQ^*Q \end{array} \right] +
  r[ \, M, \ PP^*M \, ] - 2r(M). \hfill
 \cr}
 $$
The ranks of the two block matrices in it can reduce to 
\begin{eqnarray*}
\lefteqn{ r[\, M, \ PP^*M \, ]} \\
& = & r \left[ \, \left[ \begin{array}{cc} A  & B \\ C & D \end{array} \right],
\ \left[ \begin{array}{cc} I  & 0  \\ CA^{\dagger} & I \end{array}
\right] \left[ \begin{array}{cc} I & (CA^{\dagger})^* \\ 0  & I
 \end{array} \right]\left[ \begin{array}{cc} A  & B \\ C & D \end{array}
 \right] \, \right] \\
& = & r \left[ \begin{array}{cccc} A  & B &  A + (CA^{\dagger})^*C & B + 
(CA^{\dagger})^*D  \\ C & D & CA^{\dagger}A + (CA^{\dagger})(CA^{\dagger})^*C + C  & 
CA^{\dagger}B + (CA^{\dagger})(CA^{\dagger})^*D + D  \end{array} \right] \\
& = & r \left[ \begin{array}{cccc} A  & B & (CA^{\dagger})^*C &  
(CA^{\dagger})^*D  \\ C & D & (CA^{\dagger})(CA^{\dagger})^*C + C &
 CA^{\dagger}B + (CA^{\dagger})(CA^{\dagger})^*D   \end{array} \right] \\
& = & r \left[ \begin{array}{cccc} A  & B & 0  &  0  \\ C & D &  C &
 CA^{\dagger}B  \end{array} \right] \\
& = & r \left[ \begin{array}{ccc} A  & B & 0 \\ C & D &  C  \end{array} \right] 
= r \left[ \begin{array}{ccc} A  & 0 & B \\ 0 & C &  D  \end{array} \right].
\end{eqnarray*}
Similarly we can get
$$\displaylines{
\hspace*{1.5cm}
r\left[ \begin{array}{c} M \\ MQ^*Q \end{array} \right] = 
 r\left[ \begin{array}{cc} A  & 0 \\ 0 & B \\  C  &  D \end{array} \right]. \hfill
\cr}
$$
Thus we have (9.13). Eqs.\,(9.14)---(9.16) are direct consequences of 
(9.12).
\qquad $ \Box$ 

\medskip

When $ D = 0$ in Theorem 9.4, we the following.

\medskip

\noindent  {\bf Corollary 9.5.}\, {\em Let $ M $ be given by {\rm (9.11)} where $ A \in {\cal C}^{ m \times n},
\, B \in {\cal C}^{ m \times k}$ and $ C \in {\cal C}^{ l \times n}$. Then the rank of 
$ M^{\dagger} - Q^{-1}N^{\dagger}P^{-1}$ satisfies the equality
$$ 
\displaylines{
\hspace*{2cm}
 r[ \, M^{\dagger} - Q^{-1}N^{\dagger}P^{-1} \, ]
  = r\left[ \begin{array}{c} A \\ C \end{array} \right] +
  r[\, A , \ B \, ] + r(B ) + r(C) - 2r(M). \hfill (9.17) 
\cr}
$$ 
In particular$,$ the Moore-Penrose inverse of $ M $ can be expressed as
$$\displaylines{
\hspace*{2cm} 
M^{\dagger} = Q^{-1}N^{\dagger}P^{-1} =  \left[ \begin{array}{cc} I_n
   & -A^{\dagger}B  \\ 0  & I_k \end{array} \right]
   \left[ \begin{array}{cc} A & E_AB \\ CF_A  & - CA^{\dagger}B
   \end{array} \right]^{\dagger}
 \left[ \begin{array}{cc} I_m  & 0 \\ -CA^{\dagger}  & I_l \end{array}
 \right], \hfill (9.18)
\cr}
$$
if and only if $A, \ B$ and $ C $ satisfy the following rank additivity condition
$$\displaylines{
\hspace*{2cm} 
 r( M ) =  r \left[ \begin{array}{c} A \\ C \end{array} \right] + r(B)
 =  r[ \, A, \ B \, ] + r(C). \hfill (9.19)
\cr} 
$$}

Eq.\,(9.19) shows that we have another expression for the Moore-Penrose inverse of a bordered matrix $M$ when it 
satisfies the rank additivity condition (9.19) (the first one is in (9.3)).

Clearly  the matrix $ N$ in (9.12) can be written as 
$$\displaylines{
\hspace*{1.5cm}
  N = \left[ \begin{array}{cc} A & 0 \\ 0 & 0
   \end{array} \right]  + \left[ \begin{array}{cc} 0 & E_AB \\ CF_A  & S_A
   \end{array} \right] =  N_1 +  N_2. \hfill  (9.20)
\cr}
$$
Then it is easy to verify that
$$\displaylines{
\hspace*{1.5cm}
  N^{\dagger} = \left[ \begin{array}{cc} A^{\dagger} & 0 \\ 0 & 0
   \end{array} \right]  + \left[ \begin{array}{cc} 0 & E_AB \\ CF_A  & S_A
   \end{array} \right]^{\dagger} =  N_1^{\dagger}  + N_2^{\dagger}. \hfill  (9.21)
\cr}
$$ 
Thus if we can find $ N_2^{\dagger}$, then we can give
 the expression of  $N^{\dagger}$ in (9.21). This consideration motivates  us to find the following set of results on Moore-Penrose inverses of block  matrices.

\medskip

\noindent  {\bf Lemma 9.6.}\, {\em Let $ A \in {\cal C}^{ m \times n},
\, B \in {\cal C}^{ m \times k}, \, C \in {\cal C}^{ l \times n}$  and
  $ D \in {\cal C}^{ l \times k}$ be given. Then the
 rank additivity condition
$$\displaylines{
\hspace*{2cm}
r \left[ \begin{array}{cc} A & B \\ C  & D \end{array} \right]
= r \left[ \begin{array}{c} A  \\ C  \end{array} \right] +
r \left[ \begin{array}{c} B \\ D \end{array} \right]  = r[\, A , \ B \,]
+ r[\, C, \ D \,] \hfill (9.22) 
\cr}
$$
is equivalent to the two range inclusions  
$$\displaylines{
\hspace*{2cm}
  R \left[ \begin{array}{c} A \\ 0 \end{array} \right] \subseteq R
  \left[ \begin{array}{cc} A & B \\ C  & D \end{array} \right],  \qquad 
 R \left[ \begin{array}{c} A^* \\ 0 \end{array} \right] \subseteq
  R\left[ \begin{array}{cc} A^* & C^* \\ B^*  & D^* \end{array} \right],
  \hfill(9.23)
 \cr}
$$  
and the rank additivity condition 
$$\displaylines{
\hspace*{2cm}
 r \left[ \begin{array}{cc} 0 & E_AB \\ CF_A   & S_A \end{array} \right]
 = r \left[ \begin{array}{c} E_AB \\ S_A \end{array} \right] + r (CF_A)
 = r[ \, CF_A, \ S_A \, ] + r( E_AB),   \hfill (9.24)
\cr}
$$
where $ S_A = D - CA^{ \dagger } B.$ } 

\medskip

\noindent  {\bf Proof}.\, Let
$$\displaylines{
\hspace*{2cm}
 V_1 = \left[ \begin{array}{c} A \\ C \end{array} \right], \ \ \ \
 V_2 = \left[ \begin{array}{c} B \\ D \end{array} \right],
 \ \ \ \ W_1 = [ \, A, \ B \, ],
  \ \ \ \ W_2 = [ \, C , \ D \, ]. \hfill (9.25)
\cr}
$$ 
Then (9.22) is equivalent to 
$$ \displaylines{
\hspace*{2cm}
R( V_1 ) \cap R( V_2 ) = \{ 0 \} \qquad  {\rm and } \qquad R
( W_1^* ) \cap R( W_2^*) = \{ 0 \}.  \hfill (9.26)
 \cr}$$
In that case,  we easily find 
$$\displaylines{
\hspace*{2cm}
 r \left[ \begin{array}{cc} W_1 &  A \\ W_2 & 0 \end{array} \right] = r[ \, W_1, \ A \, ]
 + r[ \, W_2 ,\ 0 \, ] = r( W_1 ) + r( W_2)
 = r( M),  \hfill
\cr
\hspace*{0cm}
and \hfill
\cr
\hspace*{2cm}
r\left[ \begin{array}{cc}   V_1 & V_2 \\ A & 0 \end{array} \right]
= r\left[ \begin{array}{c}   V_1 \\ A \end{array} \right]
+ r \left[ \begin{array}{c}  V_2 \\ 0 \end{array} \right]
= r( V_1 ) + r( V_2) = r( M),  \hfill
\cr}
$$
both of which are equivalent to the two inclusions in  (9.23).
On the other hand, observe that
$$\displaylines{
\hspace*{2cm}
\left[ \begin{array}{c}  E_AB \\ S_A \end{array} \right]
= \left[ \begin{array}{c}  B - AA^{ \dagger }B \\
D - CA^{\dagger}B \end{array} \right] = \left[ \begin{array}{c}  B \\ D
\end{array} \right]-
\left[ \begin{array}{c} A \\ C  \end{array} \right]A^{ \dagger }B
= V_2 - V_1 A^{ \dagger }B,  \hfill
\cr\hspace*{0cm}
and \hfill
\cr
\hspace*{2cm}
[ \, CF_A, \ S_A \, ]= [ \, C - CAA^{ \dagger }, \ D - CA^{\dagger}B \, ]
= [\, C , \ D\, ] - CA^{\dagger }[ \, A, \ B \, ] =  W_2 - CA^{\dagger }W_1.   \hfill
\cr}
$$  
Thus according to  (9.26) and Lemma 1.4(b) and (c), we find that
$$\displaylines{
\hspace*{2cm}
r \left[ \begin{array}{c} E_AB \\ S_A \end{array} \right] =
r( \,V_2 - V_1A^{\dagger}B \, )
= r  \left[ \begin{array}{c}V_2 \\ V_1A^{ \dagger }B \end{array} \right]
= r( V_2),  \hfill (9.27)
\cr
\hspace*{0cm}
and \hfill
\cr
\hspace*{2cm}
 r[ \, CF_A, \ S_A \, ] = r( \, W_2 - CA^{\dagger}W_1\, )
= r[ \, W_2, \ CA^{\dagger }W_1 \, ] = r( W_2).  \hfill(9.28)
\cr}
$$
From both of them and the rank formulas in  (1.2), (1.3), (1.5),
(9.23), (9.27) and (9.28), we derive the following two equalities
$$\displaylines{
\hspace*{1cm}
  r \left[ \begin{array}{cc} 0 & E_AB \\ CF_A & S_A  \end{array} \right]
  =  r( M) - r( A ) = r( V_1 ) + r(V_2 ) - r( A ) = r \left[ \begin{array}{c}
    E_AB \\ S_A \end{array} \right] + r( CF_A),  \hfill
\cr
\hspace*{0cm}
and \hfill
\cr
\hspace*{1cm}
 r \left[ \begin{array}{cc} 0 & E_AB \\ CF_A & S_A \end{array} \right]
 = r( M) - r( A ) =r( W_1 ) + r( W_2 ) - r( A ) = r[ \, CF_A, \ S_A \, ] +
 r( E_A B ).  \hfill
\cr}$$
Both of them are exactly the rank additivity condition  (9.24). Conversely,
adding $ r( A ) $ to the three sides of (9.24) and then applying
(1.2), (1.3)  and (1.5) to the corresponding result we first obtain
$$\displaylines{
\hspace*{0.5cm}
r \left[ \begin{array}{cc}  A & E_AB \\  CF_A & S_A \end{array} \right] =
 r\left[ \begin{array}{c} A \\ CF_A \end{array} \right] +
 r \left[ \begin{array}{c} E_AB \\ S_A \end{array} \right] =
r[ \, A, \ E_AB \, ] + r[ \, CF_A, \ S_A \, ]. 
\hfill  (9.29)
\cr}
$$
On the other hand, the two inclusions in  (9.23) are also equivalent to 
$$\displaylines{
\hspace*{0.5cm}
r( M) = r \left[ \begin{array}{ccc} A & B & A \\ C & D & 0 \end{array}
\right] = r\left[ \begin{array}{ccc} A & B & 0 \\ C & D & C \end{array} \right],
 \ \ {\rm and} \ \  r( M) = r \left[ \begin{array}{cc} A & B  \\ C & D  \\ A & 0
   \end{array} \right] =  \left[ \begin{array}{cc} A & B  \\ C & D
    \\ 0 & B  \end{array} \right].\hfill
\cr}
$$ 
Applying (1.5) to the right-hand sides of the above two equalities and
then combining them with (9.29), we  find
$$\displaylines{
\hspace*{0.5cm}
r( M) = r\left[ \begin{array}{ccc} A & E_AB & 0 \\ CF_A & S_A & C
\end{array} \right] =  r[ \, A, \ E_AB \, ] + r[ \, CF_A, \ S_A, \ C \, ] =
r[\, A, \ B \, ] + r[ \, C, \ D \,],  \hfill
\cr
\hspace*{0cm}
and \hfill
\cr
\hspace*{0.5cm}
 r( M ) = r\left[ \begin{array}{cc} A & E_AB \\ CF_A & S_A \\ 0
 & B  \end{array} \right]
 = r\left[ \begin{array}{c} A \\ CF_A \end{array} \right] +
 r\left[ \begin{array}{c} E_AB \\ S_A \\ B \end{array} \right]
 = r\left[ \begin{array}{c} A \\ C \end{array} \right] +
  r\left[ \begin{array}{c} B \\ D \end{array} \right]. \hfill
\cr}
$$
Both of them are exactly  (9.22).  \qquad $\Box $ 

\medskip

Similarly we can establish the following. 

\medskip

\noindent {\bf Lemma 9.7.}\, {\em  The rank additivity condition
{\rm (9.22)} is equivalent to the following four  conditions
$$\displaylines{
\hspace*{1.5cm}
 R \left[ \begin{array}{c} A \\ 0  \end{array} \right] \subseteq
 R \left[ \begin{array}{cc} A & B \\ C & D \end{array} \right], \qquad
 R \left[ \begin{array}{c}  A^* \\ 0   \end{array} \right] \subseteq
 R \left[ \begin{array}{cc} A^* &  C^* \\ B^* & D^* \end{array} \right],
 \hfill(9.30) 
\cr
\hspace*{0cm}
and \hfill
\cr
\hspace*{1.5cm}
r \left[ \begin{array}{cc} S_D & BF_D \\ E_DC & 0 \end{array} \right] =
r \left[ \begin{array}{c}S_D \\ E_DC \end{array} \right] + r( BF_D ) =
r[ \, S_D, \ BF_D \,] + r( E_DC ), \hfill (9.31)
\cr}
$$ 
where $ S_D= A- BD^{\dagger}C. $ } 

\medskip

\noindent {\bf Theorem 9.8.}\,  {\em Suppose that the block matrix $M $
in {\rm (9.11)} satisfies the rank additivity condition {\rm (9.22)}$,$
then the Moore-Penrose inverse of $M $ can be expressed in the two
forms
$$\displaylines{
\hspace*{1.5cm}
M^{ \dagger } =  \left[ \begin{array}{cc} H_1 -H_2CA^{ \dagger } -
A^{ \dagger }BH_3 +
A^{ \dagger }BJ^{ \dagger }(D)CA^{ \dagger } &
H_2 -A^{ \dagger }BJ^{ \dagger }(D) \\
H_3 - J^{ \dagger }(D)CA^{ \dagger } & J^{ \dagger }(D)
 \end{array} \right], \hfill(9.32)
 \cr
\hspace*{0cm}
and \hfill
\cr
\hspace*{1.5cm}
 M^{ \dagger } =  \left[ \begin{array}{cc} J^{ \dagger }(A) &
 J^{ \dagger }(C) \\
 J^{ \dagger }(B) & J^{ \dagger }(D) \end{array} \right] = \left[ \begin{array}{cc}
  (E_{B_2}S_DF_{C_2})^{ \dagger } & (E_{D_2}S_BF_{A_1})^{ \dagger } \cr
 ( E_{A_2}S_CF_{D_1})^{ \dagger } & (E_{C_1}S_AF_{B_1})^{ \dagger }
 \end{array} \right] ,
  \hfill (9.33)
\cr
where \hfill
\cr
\hspace*{2cm}
S_A = D - CA^{\dagger}B, \ \ \ S_B = C - DB^{\dagger}A,  \ \ \
S_C = B- AC^{\dagger}D,  \ \ \ S_D= A- BD^{\dagger}C, \hfill
\cr
\hspace*{3cm}
 A_1 = E_BA, \qquad A_2 = AF_C,   \qquad  B_1 = E_AB, \quad B_2 = BF_D,  \hfill
\cr
\hspace*{3cm}
 C_1 = CF_A,  \qquad   C_2 = E_DC,  \qquad  D_1 = E_CD, \quad D_2 = DF_B,  \hfill
\cr
\hspace*{2cm}
 H_1 = A^{\dagger} + C_1^{\dagger}[\, S_AJ^{\dagger}(D)S_A - S_A \, ]B_1^{\dagger}, \ \ \ H_2 = C_1^{\dagger}
[ \, I - S_AJ^{\dagger}(D) \, ],  \ \ \  H_3 = [ \,  I- J^{\dagger}(D)S_A \, ]B_1^{\dagger}. \hfill
\cr}
$$} 
{\bf Proof.}\, Lemma 9.6 shows  that the rank additivity condition in
 (9.22) is equivalent to  (9.23) and (9.24). It follows from
Theorem 9.4 that under (9.23), the Moore-Penrose inverse of $ M $ can be
expressed as (9.13). On the other hand, It follows from Theorem 9.1 that
 under  (9.24) the Moore-Penrose inverse of $N_2$ in  (9.20) can be
 written as
$$
 N_2^{\dagger} = \left[ \begin{array}{cc}  C_1^{\dagger}[ \,
 S_AJ^{\dagger}(D)S_A - S_A \, ]B_1^{\dagger} & C_1^{\dagger} -
 C_1^{\dagger} S_AJ^{\dagger}(D)  \cr  B_1^{\dagger}-
 J^{\dagger}(D)S_AB_1^{\dagger} & J^{\dagger}(D) \end{array} \right],
   \eqno (9.34)
$$
where $B_1 = E_AB, \, C_1 = CF_A $ and $J(D) = E_{C_1}S_AF_{B_1}.$ Now
substituting (9.34) into (9.21) and then  (9.21) into (9.13),
we get
$$
 M^{\dagger} =Q^{-1} N^{\dagger} P^{-1} = 
Q^{-1} \left[ \begin{array}{cc}   A^{\dagger} +
 C_1^{\dagger}[ \ S_AJ^{\dagger}(D)S_A - S_A \ ]B_1^{\dagger} &
 C_1^{\dagger} - C_1^{\dagger}S_AJ^{\dagger}(D)  \cr  B_1^{\dagger}-
 J^{\dagger}(D)S_AB_1^{\dagger} & J^{\dagger}(D) \end{array} \right]P^{-1}.
 \eqno (9.35)
$$
Written in a $ 2\times 2 $ block matrix,  (9.35) is (9.32). In the same 
way,  we can also decompose
$M$
in (9.11) into the other three forms,
$$ \displaylines{
\hspace*{3cm}
 M =\left[ \begin{array}{cc} I_m  & 0 \\ CB^{\dagger}  & I_l
 \end{array} \right] \left[ \begin{array}{cc} E_BA & B \\ S_B & DF_B
 \end{array} \right] \left[ \begin{array}{cc} I_n  &
  0 \\ B^{\dagger}A & I_k \end{array} \right],  \hfill
\cr
\hspace*{3cm}
 M =\left[ \begin{array}{cc} I_m  & AC^{\dagger} \\ 0 &  I_l \end{array}
 \right] \left[ \begin{array}{cc} AF_C & S_C \\ C & E_CD \end{array}
 \right] \left[ \begin{array}{cc} I_n  & C^{\dagger}D  \\ 0 &
 I_k \end{array} \right], \hfill
\cr
and \hfill
\cr
\hspace*{3cm}
M =\left[ \begin{array}{cc} I_m  &  BD^{\dagger} \\ 0 & I_l \end{array}
 \right] \left[ \begin{array}{cc} S_D & BF_D \\ E_DC & D \end{array}
 \right] \left[ \begin{array}{cc} I_n  & 0 \\  D^{\dagger}C  &
 I_k \end{array} \right]. \hfill
\cr}
$$
Based on the above decompositions of $ M $ we can also find that under
 (9.22) the Moore-Penrose inverse of $ M $ can also be expressed as
$$
M^{\dagger} =  \left[ \begin{array}{cc} * & J^{\dagger}(C) \cr * & *
\end{array} \right]
=  \left[ \begin{array}{cc}  * & * \cr J^{\dagger}( B ) & *
\end{array} \right] =
 \left[ \begin{array}{cc}  J^{\dagger}( A ) & * \cr * & *
 \end{array} \right]. \eqno (9.36)
$$
Finally from the uniqueness of the Moore-Penrose inverse of  a
matrix and the  expressions in (9.32) and (9.36), we obtain  (9.33).
\qquad   $ \Box $ 

\medskip

Some fundamental properties on the Moore-Penrose inverse
 of $M$ in  (9.11) can be derive from  (9.32) and (9.33). 

\medskip

\noindent {\bf Corollary 9.9.}\, {\em Denote the Moore-Penrose inverse of
$M$ in Eq. {\rm (9.11)} by
$$ 
M^{\dagger} = \left[ \begin{array}{cc} G_1 & G_2 \cr G_3 & G_4
\end{array} \right],
 \eqno (9.37)
$$
where $ G_1, \, G_2 , \, G_3 $ and $ G_4 $ are $ n \times m , \, n \times l,
\, k \times m $ and $ k \times l $ matrices$,$ respectively. If $M$ in
{\rm (9.11)}  satisfies the rank additivity condition {\rm (9.22)},
then the submatrices in $ M $ and $M^{\dagger}$ satisfy the  rank
equalities
$$  \displaylines{
\hspace*{4cm}
r( G_1) = r( V_1) + r( W_1 ) - r( M) +r( D ), \hfill (9.38)
\cr
\hspace*{4cm}
 r( G_2) = r( V_1) + r( W_2 ) - r( M) +r( B ), \hfill (9.39)
\cr
\hspace*{4cm}
 r( G_3) = r( V_2) + r( W_1 ) - r( M) +r( C ), \hfill (9.40)
\cr
\hspace*{4cm}
 r( G_4) = r( V_2) + r( W_2 ) - r( M) +r( A ), \hfill (9.41)
\cr
\hspace*{2cm}
 r( G_1) + r( G_4) = r( A ) + r( D ),   \qquad
 r( G_2) + r( G_3) = r( B ) + r( C ), \hfill (9.42)
\cr}
$$
where $V_1, \, V_2, \,  W_1 $ and $ W_2$ are defined in {\rm (9.25)}. Moreover$,$
 the products of $ MM^{\dagger}$ and $M^{\dagger}M$ have the forms
$$
 MM^{\dagger} = \left[ \begin{array}{cc}  W_1W_1^{\dagger} & 0 \cr 0 &  
W_2W_2^{\dagger} \end{array} \right], \qquad  M^{\dagger}M =
 \left[ \begin{array}{cc} V_1^{\dagger}V_1 & 0 \cr 0 &  V_2^{\dagger}V_2
  \end{array} \right]. \eqno (9.43)
$$ } 

\noindent {\bf Proof.} \  The four rank equalities in (9.38)---(9.42)
can directly be derived from the expression in (9.33) for $ M^{\dagger} $
and the rank formula (1.6). The two equalities in (9.42) come from the
 sums of (9.38) and (9.41),  (9.39) and (9.40), respectively. The
 two results in (9.43) are derived from (9.22) and Theorem 7.16(c)
 and (d).  \qquad  $ \Box $ 

\medskip

The rank additivity condition (9.22) is a quite weak restriction
to a $ 2 \times 2 $ block matrix. As a  matter of fact, any matrix with its rank great then 1 satisfies a rank additivity condition 
as in (9.22) when its rows and columns are properly permuted. We next present a group of consequences of Theorem 9.8. 

\medskip

\noindent {\bf Corollary 9.10.}\, {\em If the block matrix $M$ in
 {\rm (9.11)} satisfies {\rm (9.23)} and the following two conditions
$$
 R(C_1) \cap R( S_A) = \{ 0 \}  \ \ \  and  \ \ \
 R( B_1^*) \cap R(S_A^*)= \{ 0 \}, \eqno (9.44)
 $$
then the Moore-Penrose inverse of $M$ can be expressed as 
\begin{eqnarray*}
\left[ \begin{array}{cc}   A & B \cr C & D \end{array} \right]^{\dagger}
& = & Q^{-1}
\left[ \begin{array}{cc} A^{\dagger} & C_1^{\dagger} -
C_1^{\dagger}S_AJ^{\dagger}(D)  \cr  B_1^{\dagger}-
J^{\dagger}(D)S_AB_1^{\dagger} & J^{\dagger}(D) \end{array} \right]P^{-1} \\
& = &  \left[ \begin{array}{cc}  A^{\dagger} -H_2CA^{ \dagger } -
A^{ \dagger }BH_3 +
A^{ \dagger }BJ^{ \dagger }(D)CA^{ \dagger } &
H_2 -A^{ \dagger }BJ^{ \dagger }(D) \cr
H_3 - J^{ \dagger }(D)CA^{ \dagger } & J^{ \dagger }(D) \end{array} \right],
\end{eqnarray*}
where $ C_1, \, B_1, \, H_2, \, H_3 $ and $ J( D )$ are as in  {\rm (9.32)},
$P $ and $ Q $ are as in {\rm (9.13)}. } 

\medskip

\noindent {\bf Proof.}\,  The conditions in (9.44) imply that the block matrix
$ N_2$ in  (9.20) satisfies the following rank additivity condition
$$
 r(N_2 ) = r( E_AB) + r( CF_A) +r( S_A),
$$
which is a special case of  (9.24).  On the other hand, under  (9.44)
if follow by Theorem 7.7 that $S_A J^{\dagger}(D)S_A = S_A$. Thus  (9.35)
 reduces to the desired result in the corollary.   \qquad  $ \Box $ 

\medskip

\noindent {\bf Corollary 9.11.}\, {\em If the block matrix $M$ in
{\rm (9.11)} satisfies {\rm (9.23)} and the two conditions
$$ 
R(BS_A^*) \subseteq R(A) \ \ \  and \ \ \ R(C^*S_A) \subseteq R(A^*),
\eqno (9.45)
 $$
then the Moore-Penrose inverse of $M$ can be expressed as
\begin{eqnarray*} 
M^{\dagger} & = &  \left[ \begin{array}{cc}  I_n & -A^{\dagger}B \cr 0 &
  I_k \end{array} \right]
   \left[ \begin{array}{cc}  A^{\dagger} & ( CF_A)^{\dagger} \cr
   (E_AB)^{\dagger} &
  S_A^{\dagger} \end{array} \right]  \left[ \begin{array}{cc}  I_m & 0 \cr
  -CA^{\dagger} & I_l \end{array} \right] \\
& = &  \left[ \begin{array}{cc}  A^{\dagger} -
A^{\dagger}B(E_AB)^{\dagger} -
(CF_A)^{\dagger}CA^{\dagger} + A^{\dagger}BS_A^{\dagger}CA^{\dagger} &
(CF_A)^{\dagger} - A^{\dagger}BS_A^{\dagger} \\ ( E_AB)^{\dagger} -
S_A^{\dagger}CA^{\dagger} & S_A^{\dagger} \end{array} \right],
\end{eqnarray*} 
where $ S_A = D- CA^{\dagger}B. $ } 

\medskip

\noindent {\bf Proof.}\, Clearly  (9.45) are equivalent to
$ (E_AB)S^*_A= 0 $ and $S^*_A(CF_A)= 0 $. In that case,
 (9.24) is satisfied, and  $ N^{\dagger}=
 \left[ \begin{array}{cc}  A^{\dagger} & (CF_A)^{\dagger} \cr
 (E_AB)^{\dagger} & S_A^{\dagger} \end{array} \right]$
in (9.35). \qquad  $ \Box $ 

\medskip

\noindent {\bf Corollary 9.12}(Chen and Zhou \cite{CZ}).\, {\em If the block matrix $M$ in
 {\rm (9.11)} satisfies the following four
conditions
$$ 
R(B) \subseteq R(A),  \qquad   R(C^*) \subseteq R( A^*), \qquad  R(C )
\subseteq R(S_A), \qquad R( B^* ) \subseteq R(S_A^*),   \eqno (9.46)
 $$
then the Moore-Penrose inverse of $M$ can be expressed as
\begin{eqnarray*} 
\left[ \begin{array}{cc}  A & B \cr C & D \end{array} \right]^{\dagger}
& = &  \left[ \begin{array}{cc} I_n & -A^{\dagger}B
\cr 0 & I_k \end{array} \right]  \left[ \begin{array}{cc}  A^{\dagger} & 0
\cr 0  & S_A^{\dagger}\end{array} \right]
\ \left[ \begin{array}{cc} I_m & 0 \cr -CA^{\dagger} & I_l \end{array}
\right] \\
& = &\left[ \begin{array}{cc} A^{\dagger} + A^{\dagger}BS_A^{\dagger}
CA^{\dagger} & - A^{\dagger}BS_A^{\dagger} \cr - S_A^{\dagger}CA^{\dagger} &
S_A^{\dagger} \end{array} \right],
\end{eqnarray*} 
where $ S_A = D- CA^{\dagger}B. $ } 

\medskip

\noindent {\bf Proof.}\,  It is easy to verify that under the conditions in (9.46),
 the rank of $M$ satisfies the rank additivity condition (9.22). In that
 case, $ N^{\dagger} = \left[ \begin{array}{cc}  A^{\dagger} & 0 \cr 0  &
 S_A^{\dagger} \end{array} \right] $ in (9.35).   \qquad $ \Box $ 

\medskip

\noindent {\bf Corollary 9.13.} \ {\em If the block matrix $M$ in
 {\rm (9.11)} satisfies the  four conditions
$$
 R(A) \cap R(B)= \{ 0 \},  \ \   R(A^*) \cap R(C^*) = \{ 0 \}, \ \ 
 R(D) \subseteq R(C),  \ \  R(D^*) \subseteq R(B^*),   \eqno (9.47)
 $$
then the Moore-Penrose inverse of $M$ can be expressed as
\begin{eqnarray*}
 \left[ \begin{array}{cc} A & B \cr C & D \end{array} \right]^{\dagger}
 & = & \left[ \begin{array}{cc}  I_n &
 -A^{\dagger}B \cr 0 & I_k \end{array} \right]
 \left[ \begin{array}{cc}  A^{\dagger} -
 C_1^{\dagger}S_AB_1^{\dagger}  & C_1^{\dagger}
  \cr B_1^{\dagger}  & 0 \end{array} \right]
  \left[ \begin{array}{cc}  I_m & 0 \cr
  -CA^{\dagger} &  I_l \end{array} \right] \\
& = & \left[ \begin{array}{cc}  A^{\dagger} - A^{\dagger}BB_1^{\dagger} -
C_1^{\dagger}CA^{\dagger} -C_1^{\dagger}S_AB_1^{\dagger}  &
C_1^{\dagger} \cr  B_1^{\dagger} & 0  \end{array} \right],
\end{eqnarray*} 
where $ S_A = D- CA^{\dagger}B, \ B_1 = E_AB $ and $C_1 = CF_A. $ } 

\medskip

\noindent {\bf Proof.}\, It is not difficult to verify by (1.5) that
 under (9.47) the rank of $ M$ satisfies  (9.22). In that case,
  $ J(D) = 0$ and  $ N^{\dagger } =  \left[ \begin{array}{cc} A^{\dagger} -
  C_1^{\dagger}S_AB_1^{\dagger} & C_1^{\dagger} \cr B_1^{\dagger} & 0
  \end{array} \right] $ in (9.35).  \qquad  $ \Box $  

\medskip

\noindent {\bf Corollary 9.14.}\, {\em If the block matrix $M$ in
 {\rm (9.11)} satisfies the four conditions
$$
R(A) \cap R(B) = \{ 0 \},    \qquad   R(A^*) \cap R( C^* ) = \{ 0 \},
\eqno (9.48)
$$
$$ 
R( S_A) \subseteq N(C^*),   \qquad  R( S^*_A ) \subseteq N(B ),
\eqno (9.49)
 $$
then the Moore-Penrose inverse of $M$ can be expressed as
\begin{eqnarray*}
 \left[ \begin{array}{cc}  A & B \cr C & D \end{array} \right]^{\dagger}
 & = &  \left[ \begin{array}{cc}  I_n &
 -A^{\dagger}B \cr 0 & I_k \end{array} \right]
 \left[ \begin{array}{cc}  A^{\dagger} &
 (CF_A)^{\dagger}
  \cr (E_AB)^{\dagger} & S_A^{\dagger} \end{array} \right]
  \left[ \begin{array}{cc}  I_m & 0 \cr -  CA^{\dagger} & I_l   \end{array}
  \right]  \\
& = & \left[ \begin{array}{cc}  A^{\dagger} - A^{\dagger}B(E_AB)^{\dagger} -
(CF_A)^{\dagger}CA^{\dagger} & (CF_A)^{\dagger} \cr ( E_AB)^{\dagger} &
S_A^{\dagger} \end{array} \right],
\end{eqnarray*} 
where $ S_A = D- CA^{\dagger}B. $ } 

\medskip

\noindent {\bf Proof.}\, Clearly (9.49) is equivalent to $
C^{\dagger}S_A = 0 $ and $ S_AB^{\dagger} = 0$, as well as $S_A^{\dagger }C
= 0$ and $ BS_A^{\dagger} = 0 $. From them and (9.48), we also find
$$ 
(CF_A)^{\dagger}S_A = 0  \ \ \ and \ \ \  S_A(E_AB)^{\dagger} = 0.
\eqno (9.50)
$$
Combining (9.48) and (9.50) shows that $M$ satisfies (9.23) and
(9.24). In that case,
$$ 
N^{\dagger} =  \left[ \begin{array}{cc}  A^{\dagger} &
 (CF_A)^{\dagger} \cr (E_AB)^{\dagger} & S_A^{\dagger} \end{array} \right]
$$ 
in (9.35).   \qquad  $\Box $  

\medskip

\noindent {\bf Corollary 9.15.}\, {\em If the block matrix $M$ in
 {\rm (9.11)} satisfies
the  rank additivity condition
$$ 
r( M_1 ) = r( A ) + r( B ) + r( C ) + r ( D ),  \eqno (9.51) 
$$
then the Moore-Penrose inverse of $M$ can be expressed as
$$
 \left[ \begin{array}{cc}  A & B \cr C & D \end{array} \right]^{\dagger} =
  \left[ \begin{array}{cc}  ( E_BAF_C)^{\dagger}
 & (E_DCF_A)^{\dagger} \cr (E_ABF_D)^{\dagger} & (E_CDF_B)^{\dagger}
  \end{array} \right]. \eqno (9.52)
$$ } 
{\bf Proof.}\, Obviously  (9.51) is a special case of  (9.22).
On the other hand,  (9.51) is also equivalent to the following four
conditions
$$ \displaylines{
\hspace*{2cm}
R(A) \cap  R(B) = \{0\}, \ \   R(C) \cap  R(D) = \{0\}, \ \ 
 R(A^*) \cap  R(C^*) = \{0\}, \ \ R(B^*) \cap  R(D) = \{0\}.   \hfill
\cr
In \ that \ case, \hfill
\cr
\hspace*{2cm}
 R(A_1^*)  = R(A^*), \qquad    R(A_2)  = R(A), \qquad
 R(A_1^*)  = R(A^*), \qquad   R(B_2)  = R(B), \hfill
 \cr
\hspace*{2cm}
 R(C_1)  = R(A), \qquad R(C_2^*)  = R(C^*), \qquad  
 R(D_1^*)  = R(D^*), \qquad   R(A_2)  = R(D), \hfill
\cr}
$$
by Lemma 1.2(a) and (b). Then it turns out by Theorem 7.2(c) and (d) that 
$$ \displaylines{
\hspace*{2cm}
A_1^{\dagger}A_1 =  A^{\dagger}A, \qquad
A_2A_2^{\dagger} =  AA^{\dagger},  \qquad
B_1^{\dagger}B_1 =  B^{\dagger}B,  \qquad  B_2B_2^{\dagger} =
 BB^{\dagger}, \hfill
\cr
\hspace*{2cm}
 C_1C_1^{\dagger} =  CC^{\dagger}, \qquad
 C_2^{\dagger}C_2 =  C^{\dagger}C, \qquad
 D_1^{\dagger}D_1 = D^{\dagger}D, \qquad
D_2D_2^{\dagger} =  DD^{\dagger}. \hfill
\cr}
$$
Thus (9.33) reduces to (9.52)
 \qquad  $ \Box $  

\medskip

\noindent {\bf Corollary 9.16.}\, {\em If the block matrix $M$ in
Eq.\, {\rm (9.11)} satisfies $ r( M) = r( A ) + r ( D ) $ and
$$ 
R(B) \subseteq R( A ),  \ \ \ R(C ) \subseteq R(D ),  \ \ \ R(C^*)
\subseteq R(A^*),  \ \ \ R( B^* ) \subseteq R(D^* ),
$$
then the Moore-Penrose inverse of $M$ can be expressed as
$$ 
 \left[ \begin{array}{cc}  A & B \cr C & D \end{array} \right]^{\dagger} =
  \left[ \begin{array}{cc}  ( \,A - BD^{\dagger}C \,)^{\dagger} &
  -A^{\dagger}B( \, D - CA^{\dagger}B \,  )^{\dagger} \cr -( \,D -
  CA^{\dagger}B \, )^{\dagger}CA^{\dagger}  & ( \, D- CA^{\dagger}B \,
  )^{\dagger}  \end{array} \right].
$$ } 
{\bf Corollary 9.17.}\, {\em If the block matrix $M$ in  {\rm (9.11)} 
satisfies $ r( M ) = r( A ) + r ( D ) $ and the following four
conditions
$$
 R( A )= R( B ), \ \ \ R( C )= R( D ), \ \ \  R(A^*)= R(C^*), \ \ \
 R( B^* )= R(D^*),
$$
then the Moore-Penrose inverse of $M$ can be expressed as
$$
\left[ \begin{array}{cc}   A & B \cr C & D \end{array} \right]^{\dagger}
=  \left[ \begin{array}{cc}  S_D^{\dagger} &
  S_B^{\dagger}  \cr  S_C^{\dagger}  &   S_A^{\dagger} \end{array} \right]  =
  \left[ \begin{array}{cc}  ( \, A - BD^{\dagger}C \, )^{\dagger} &
  ( \, C - DB^{\dagger}A
   \, )^{\dagger} \cr ( \, B - AC^{\dagger}D  \, )^{\dagger} &
   ( \, D- CA^{\dagger}B  \, )^{\dagger} \end{array} \right].
$$} 
\hspace*{0.3cm} The above two corollaries can directly be derived from  (9.32) and (9.33),
the proofs are omitted here. 

\medskip

Without much effort, we can extend the results in Theorem 9.8 to
$ m \times n$  block matrices when they satisfy rank additivity conditions. 

Let
$$
 M = \left[ \begin{array}{cccc} A_{11} & A_{12} & \cdots & A_{1n} \\
                 A_{21} & A_{22} & \cdots & A_{2n} \\
                 \cdots & \cdots & \cdots & \cdots \\
                 A_{m1} & A_{m2} & \cdots & A_{mn}  \end{array} \right] 
   \eqno (9.53) 
$$
be an $ m \times n $ block matrix 
, where $ A_{ij} $ is an $ s_i \times t_j $ matrix $(1 \leq i \leq m,\ 1
\leq j
\leq n ),$ and suppose that $ M $ satisfies the following rank additivity
condition
$$ 
r( M ) = r( W_1) + r( W_2) + \cdots + r( W_m) = r( V_1) + r( V_2 ) + \cdots
+ r( V_n ), \eqno (9.54)
 $$ 
where 
$$
 W_i = [ \, A_{i1}, \, A_{i2}, \, \cdots, \, A_{in} \, ], \qquad V_j =
 \left[ \begin{array}{c}  A_{1j} \cr A_{2j} \cr \vdots \cr A_{mj}
 \end{array} \right], \qquad 1 \leq i \leq m, \ 1 \leq j \leq n. \eqno (9.55)
$$
For convenience of representation, we adopt the following notation. Let
$ M = ( A_{ij} ) $ be given in (9.53), where $ A_{ij}\in {\cal C}^{
s_i \times t_j}, \
1 \leq i \leq m, \ 1 \leq j \leq n,$ and $ \sum_{i=1}^m s_i = s, \
 \sum_{i=1}^n t_i = t.$ For each $ A_{ij}$ in $M$ we associate three block
 matrices as follows
$$
 B_{ij} = [ \, A_{i1}, \, \cdots, \, A_{i,j-1}, \, A_{i,j+1}, \, \cdots, \, A_{in}
 \, ], \eqno(9.56)
$$
$$
C^*_{ij} = [ \,  A^*_{1j},\, \cdots, \ ,A^*_{i-1,j}, \, A^*_{i+1,j}, \ \cdots, \
A^*_{mj} \, ], \eqno (9.57)
$$ 
$$
D_{ij} = \left[ \begin{array}{cccccc}  A_{11} & \cdots & A_{1,j-1} & A_{1,j+1} &
\cdots & A_{1n} \cr \vdots & & \vdots & \vdots & & \vdots \cr
 A_{i-1,1} & \cdots & A_{i-1,j-1} & A_{i-1,j+1} & \cdots & A_{i-1,n} \cr 
A_{i+1,1} & \cdots & A_{i+1,j-1} & A_{i+1,j+1} & \cdots & A_{i+1,n} \cr 
\vdots & & \vdots &  \vdots & & \vdots \cr
A_{m1} & \cdots & A_{m,j-1} & A_{m,j+1} & \cdots & A_{mn} \end{array}
\right]. \eqno (9.58)
$$
The symbol $ J( A_{ij}) $ stands for 
$$ 
J( A_{ij} ) = E_{\alpha_{ij}}S_{D_{ij}}F_{\beta_{ij}}, \qquad  1 \leq i \leq
m, \ 1 \leq j \leq n, \eqno (9.59) 
$$
where $ \alpha_{ij} = B_{ij}F_{D_{ij}}, \ \beta_{ij} = E_{D_{ij}}C_{ij}$ and
$  S_{D_{ij}} = A_{ij}- B_{ij}D_{ij}^{\dagger}C_{ij}$ is  the Schur
complement of $ D_{ij}$ in $ M $.  We call the matrix $ J(A_{ij})$
the {\it rank complement} of $ A_{ij} $ in $ M $. Besides we partition the
Moore-Penrose inverse of $M $ in (9.53) into the form
$$
 M^{\dagger} =  \left[ \begin{array}{cccc}  G_{11} & G_{12} & \cdots & G_{1m}
 \cr  G_{21} & G_{22} & \cdots & G_{2m} \cr   
  \cdots & \cdots & \cdots & \cdots \cr
  G_{n1} & G_{n2} & \cdots & G_{nm} \end{array} \right], \eqno (9.60) 
$$
where $ G_{ij} $ is a $ t_i \times s_j $ matrix, $1 \leq i \leq n, \ 1
\leq j \leq m$. 

\medskip

Next we build two groups of block permutation matrices as follows
$$ \displaylines{
\hspace*{3cm}
P_1= I_s, \qquad  P_i =  \left[ \begin{array}{cccccccc} 0 & & & I_{s_i} & &
& \cr I_{s_1} & \ddots & & & & & \cr
                        & \ddots & \ddots & & & & \cr
                      & & I_{s_{i-1}}& 0 &  & & \cr
                      & & & & I_{s_{i+1}}  & & \cr
                      & & & & & \ddots & \cr
                      & & & & & & I_{s_m } \end{array} \right],  \hfill (9.61)
\cr
\hspace*{3cm}
 Q_1= I_s,  \qquad  Q_j =  \left[ \begin{array}{cccccccc}  0 & I_{t_1} & & &
 & & \cr
                        & \ddots & \ddots & & & & \cr
                       &  & \ddots & I_{t_{j-1}} & & &  \cr
                       I_{t_j} & & & 0 &  & & \cr
                      & & & & I_{t_{j+1}}  & & \cr
                      & & & & & \ddots & \cr
                      & & & & & & I_{t_n} \end{array} \right], \hfill (9.62)
\cr}
$$ 
where $ 2 \leq i \leq m, \ 2 \leq j \leq n.$  Applying (9.61) and (9.62)
 to $M$ in (9.53) and $M^{\dagger} $ in (9.60) we have the following
 two groups of expressions
$$
P_iM Q_j = \left[ \begin{array}{cc}  A_{ij} & B_{ij} \cr C_{ij} & D_{ij}
\end{array} \right], \qquad  1 \leq i \leq m, \ 1 \leq j \leq n,
\eqno (9.63)
$$
and
$$
Q^T_jM^{\dagger} P_i^T =  \left[ \begin{array}{cc}  G_{ji} & * \cr * & *
\end{array} \right], \qquad 1 \leq i \leq m, \ 1 \leq j \leq n. \eqno (9.64)
$$
These two equalities show that we can use two block permutation matrices to
permute $ A_{ij} $ in $ M $  and the corresponding block $ G_{ji} $ in
$ M^{\dagger}$ to the upper left corners of $ M $ and $ M^{\dagger }$,
respectively. Observe that $ P_i $ and $ Q_j $ in  (9.61) and (9.62)
are all orthogonal matrices. The Moore-Penrose inverse of  $ P_iMQ_j$  in
 (9.63)
can be expressed as $ ( P_iMQ_j)^{\dagger} = Q_j^TM^{\dagger}P_i^T. $
Combining  (9.63) with (9.64), we have the following simple result  
$$ 
 \left[ \begin{array}{cc}  A_{ij} & B_{ij} \cr C_{ij} & D_{ij}
 \end{array} \right]^{\dagger} =
  \left[ \begin{array}{cc}  G_{ji} & * \cr * & * \end{array} \right], \qquad 1 \leq i
  \leq m, \ 1 \leq j \leq n. \eqno (9.65)
$$ 
If the block matrix $ M $ in  (9.53) satisfies the rank additivity
condition (9.54), then the $ 2 \times 2 $
block matrix on the right-hand side of (9.63) naturally satisfies the
following rank additivity condition
$$
 r \left[ \begin{array}{cc}  A_{ij} & B_{ij} \cr C_{ij} & D_{ij}
 \end{array} \right] =
 r  \left[ \begin{array}{cc}  A_{ij}
\cr C_{ij} \end{array} \right] + r \left[ \begin{array}{c} B_{ij} \cr D_{ij}
 \end{array} \right] = r[ \, A_{ij}, \  B_{ij} \, ] +
 r[ \, C_{ij}, \ D_{ij} \, ],  \eqno (9.66)
$$ 
where $ 1 \leq i \leq m,  \ 1 \leq j \leq n.$ Hence combining  (9.65) and
(9.66) with Theorems 9.8 and 9.9, we obtain the following general result. 

\medskip

\noindent {\bf Theorem 9.18.}\, {\em  Suppose that the $ m \times n $ block
matrix $ M $ in {\rm (9.53)} satisfies the rank additivity condition {\rm (9.54)}. Then 

 {\rm (a)}\,  The Moore-Penrose inverse of $M $ can be expressed as
$$
 M^{\dagger} =  \left[ \begin{array}{cccc}
 J^{\dagger}( A_{11}) & J^{\dagger}( A_{21})& \cdots
 & J^{\dagger}( A_{m1}) \cr
J^{\dagger}( A_{12}) & J^{\dagger}( A_{22})& \cdots & J^{\dagger}( A_{m2})
\cr  \cdots & \cdots & \cdots & \cdots \cr
J^{\dagger}( A_{1n}) & J^{\dagger}( A_{2n})& \cdots & J^{\dagger}( A_{mn})
    \end{array} \right], \eqno (9.67)
$$
where $ J( A_{ij}) $ is defined in {\rm (9.59)}.

{\rm (b)}\, The rank of the block entry $ G_{ji} = J^{\dagger}( A_{ij})
 $ in $ M^{\dagger} $ is
$$ r( G_{ji} ) = r[J( A_{ij})] = r( W_i) + r( V_j) - r( M ) + r( D_{ij} ),
\eqno (9.68)
$$
where $ 1 \leq i \leq m, \ 1 \leq j \leq n, \  W_i $ and $V_j $ are defined
in {\rm (9.55)}.

{\rm (c)} \  $ MM^{\dagger} $ and $M^{\dagger}M $ are two block diagonal
 matrices
$$ \displaylines{
\hspace*{2cm} 
MM^{\dagger} = {\rm diag}( \, W_1W_1^{\dagger}, \, W_2W_2^{\dagger}, \, \cdots,
\, W_mW_m^{\dagger} \,  ), \hfill (9.69)
\cr
\hspace*{2cm}
 M^{\dagger}M = {\rm diag}( \, V_1^{\dagger}V_1, \, V_2^{\dagger}V_2, \, \cdots,
  \, V_n^{\dagger}V_n  \, ), \hfill (9.70)
\cr}
$$  
written in explicit forms$,$ {\rm (9.69)} and {\rm (9.70)} are equivalent to
$$ \displaylines{
\hspace*{2cm}
A_{i1}G_{1j} + A_{i2}G_{2j} + \cdots + A_{in}G_{nj} =
 \left\{ \begin{array}{ccc}  W_iW_i^{\dagger} &  & i = j
\\ 0 &  &  i \neq j \end{array} \right.  \quad  i, \, j = 1, \, 2, \, \cdots ,
\, m, \hfill
\cr
\hspace*{2cm}
G_{i1}A_{1j} + G_{i2}A_{2j} + \cdots + G_{im}A_{mj} =
\left\{ \begin{array}{ccc}  V_i^{\dagger}V_i & &   i = j
\\ 0  &  & i \neq j \end{array} \right.  \quad  i, \, j = 1, \, 2, \,  \cdots,
 \, n. \hfill
\cr}
$$}     

In addition to the expression given in (9.67) for $ M^{\dagger}$, we can
also derive some other
expressions for $ G_{ij} $ in $ M^{\dagger} $ from (9.32). But they are quite
complicated in form, so we  omit them here. 

Various consequences can be derived from  (9.67) when the
submatrices in $ M $ satisfies some additional conditions, or $ M $ has
some particular patterns, such as triangular forms, circulant forms and
tridiagonal forms. Here we only give one consequences. 

\medskip

\noindent {\bf Corollary 9.19.}\, {\em If the block matrix $ M $ in
 {\rm (9.53)} satisfies the following rank additivity condition
$$
 r( M ) = \sum_{i=1}^m \sum_{j=1}^n r( A_{ij} ), \eqno (9.71)
 $$
then the Moore-Penrose inverse of $M $ can be expressed as
$$ 
M^{\dagger} =  \left[ \begin{array}{ccc}
( E_{B_{11}}A_{11}F_{C_{11}})^{\dagger} & \cdots &
( E_{B_{m1}}A_{m1}F_{C_{m1}})^{\dagger} \cr  \vdots &  & \vdots \cr 
( E_{B_{1n}}A_{1n}F_{C_{1n}})^{\dagger} &
\cdots & ( E_{B_{mn}}A_{mn}F_{C_{mn}})^{\dagger}  \end{array} \right],
\eqno (9.72)
$$
where $B_{ij}$ and $ C_{ij} $ are defined in  {\rm (9.56)} and
{\rm (9.57)}. }

\medskip

\noindent {\bf Proof.}\, In fact, (9.71) is equivalent to 
$$  \displaylines{
\hspace*{2cm}
R( A_{ij} ) \cap R( B_{ij} ) = \{ 0 \}, \qquad R( A_{ij}^* )
\cap R( C_{ij}^* )= \{ 0 \}, \qquad 1 \leq i \leq m, \ 1 \leq j \leq n, \hfill
\cr
\hspace*{2cm}
R( C_{ij} ) \cap R( D_{ij} ) = \{ 0 \},  \qquad  R(B_{ij}^*)
\cap R( D_{ij}^* )= \{ 0 \},  \qquad 1 \leq i \leq m, \ 1 \leq j \leq n. \hfill
\cr}
$$
We can get from them $J(A_{ij}) = E_{B_{ij}}A_{ij}F_{C_{ij}}. $ Putting
them in  (9.67) yields  (9.72). \qquad  $ \Box $ 

\medskip

The results so far established in the chapter are mainly based on the factorizations (9.1) and (9.12).      
However, if a given block matrix has certain special pattern such that we can factor it in some particular methods, then we can establish some new rank equalities through the special factorizations of the block matrix.    
  Here  we present some examples on the Moore-Penrose inverse of some special block matrices. 

\medskip

\noindent {\bf Theorem 9.20.}\, {\em Let $ A, \, B \in {\cal C}^{ m \times n}, \, 0 \neq p 
\in {\cal C},$ and let  
$$ 
 M = \left[ \begin{array}{cc} A & p^2B \\ B   & A  \end{array} \right], \qquad   N = \left[ \begin{array}{cc} A  + pB & 0 \\  0  & A - pB  \end{array} \right].
$$ 

{\rm (a)}\, If $  |p| = 1, $  then 
$$ \displaylines{
\hspace*{3cm} 
\left[ \begin{array}{cc} A & p^2B \\ B & A  \end{array} \right]^{\dagger} = P_{2n} \left[ \begin{array}{cc} ( \, A  + pB \, )^{\dagger}  & 0 \\  0  & (\, A - pB \,)^{\dagger}  \end{array} \right] P_{2m},  
\hfill (9.73)
\cr
where \hfill
\cr 
\hspace*{3cm} 
P_{2t} = P_{2t}^{-1}= \frac{1}{\sqrt{2}} \left[ \begin{array}{cc} I_t & pI_t \\ p^{-1}I_t & -I_t  \end{array} \right], \qquad  t = m, \ n.  \hfill
\cr
In \ particular,\hfill
\cr
\hspace*{1cm} 
\left[ \begin{array}{rr} A & -B \\ B & A  \end{array} \right]^{\dagger} = \frac{1}{2}
 \left[\begin{array}{cc} I_n & iI_n \\ -iI_n & -I_n \end{array} \right]
\left[ \begin{array}{cc} ( \, A  + iB \, )^{\dagger}  & 0 \\  0  & (\, A - iB \,)^{\dagger}  \end{array} \right] 
\left[\begin{array}{cc} I_m & iI_m \\ -iI_m & -I_m  \end{array} \right].  
\hfill
}
$$  

{\rm (b)}\, If $  |p| \neq 1, $  then 
$$ 
r(\, M^{\dagger} - P_{2n}N^{\dagger}P_{2m} \, )= 2 r\left[ \begin{array}{c} A  \\ B  \end{array} \right] +
 2r[ \, A, \ B \, ] -  2 r(\, A + pB \,) - 2 r(\, A - pB \,).    \eqno (9.74)
$$ 

{\rm (c)}\, Under $ |p| \neq 1,$ {\rm (9.73)} holds if and only if $ R(A) \subseteq R( \, A + pB \, )$ 
and  $ R(A^*) \subseteq R( \, A^* - \overline{p}B^* \, )$. } 

\medskip  
  
\noindent {\bf Proof.}\, It is easy to verify that the block matrix $ M$ can factor as 
$M = P_{2m}NP_{2n}$, where $ P_{2m}$ and $ P_{2n}$ are nonsingular with $ P_{2m}^2 = I_{2m}$ and $P_{2n}^2 = I_{2n}$. 
In that case, we find by Theorem 8.14 that 
$$  \displaylines{
\hspace*{1cm}
r(\, M^{\dagger} - P_{2n}^{-1}N^{\dagger}P_{2m}^{-1} \, )= 
r(\, M^{\dagger} - P_{2n}N^{\dagger}P_{2m} \, )=  r\left[ \begin{array}{c} N  \\ NP_{2n}P^*_{2n} \end{array} \right] +
 r[ \, N, \ P_{2m}^*P_{2m} N \, ] - 2r(N),   \hfill (9.75)
\cr
where \hfill
\cr
\hspace*{1.5cm}   
 P_{2m}^*P_{2m} = \frac{1}{2} \left[ \begin{array}{cc} I_m & \overline{p}^{-1}I_m \\ \overline{p}I_m & -I_m  \end{array} \right] \left[ \begin{array}{cc} I_m & pI_m \\ p^{-1}I_m & -I_m  \end{array} \right]
 =\frac{1}{2} \left[ \begin{array}{cc} ( \,1 + |p|^{-2} \,)I_m & ( \, p - \overline{p}^{-1} \, )I_n \\
 ( \, \overline{p} - p^{-1} \, ) I_m  &  ( \,1 + |p|^2 \,)I_m \end{array} \right],  \hfill
 \cr
\hspace*{1.5cm}
 P_{2n}P_{2n}^* = \frac{1}{2} \left[ \begin{array}{cc} I_n & pI_n \\  p^{-1}I_n & -I_n \end{array} \right]
 \left[ \begin{array}{cc} I_n & \overline{p}^{-1}I_n \\ \overline{p}I_n  & -I_n  \end{array} \right]
 =\frac{1}{2} \left[ \begin{array}{cc} ( \,1 + |p|^2 \,)I_n & ( \, \overline{p}^{-1} - p \,)I_n \\
 ( \, p^{-1} - \overline{p} \, ) I_n  &  ( \,1 + |p|^{-2} \,)I_n \end{array} \right]. \hfill 
\cr}
$$
If  $ |p| = 1$, then $ P_{2m}^*P_{2m} = I_{2m}$ and $ P_{2n}P_{2n}^* = I_{2n}$. Thus the right-hand side of (9.75) becomes zero, which implies that $  M^{\dagger} = P_{2n}N^{\dagger}P_{2m}$, the desired result in (9.73). 
 If $ |p| \neq 1,$ then we find 
$$ 
r[\, N, \ P_{2m}^*P_{2m} N \,] = r \left[ \begin{array}{cccc} 0 &  A - pB  &  A + pB  & 0 
  \\  A + pB  & 0 &  0 & A - pB  \end{array} \right] =  2r[ \, A, \ B \, ]. 
$$  
Similarly $ r\left[ \begin{array}{c} N  \\ NP_{2n}P^*_{2n} \end{array} \right] = 2r\left[ \begin{array}{c} A  \\ B \end{array} \right]. $  Thus  (9.75) becomes (9.74).   \qquad $\Box$ 

\medskip

\noindent {\bf Theorem 9.21.}\, {\em Let $ M = \left[ \begin{array}{cc} A & A \\ A   & A + B  \end{array} \right]$. 
where $ A, \ B \in {\cal C}^{ m \times n}.$ Then $ M$ factor as
$$
\displaylines{
\hspace*{2cm}
 M = PNQ = \left[ \begin{array}{cc} I_m & 0 \\ I_m & I_m \end{array} \right] \left[ \begin{array}{cc} A & 0 \\  0  & B  \end{array} \right] \left[ \begin{array}{cc} I_n & I_n \\ 0 & I_n \
\end{array} \right]. \hfill
\cr
\hspace*{0cm}
In \ that \ case \hfill
\cr
\hspace*{2cm} 
r(\, M^{\dagger} - Q^{-1}N^{\dagger}P^{-1} \, ) = 2 r\left[ \begin{array}{c} A  \\ B  \end{array} \right] +
 2r[ \, A, \ B \, ] -  2 r(A) - 2 r(B).   \hfill
\cr}
$$ 
In particular$,$ the Moore-Penrose inverse inverse of $ M^{\dagger}$ can be expressed as 
 $$ \displaylines{
\hspace*{2cm}
 M^{\dagger} = Q^{-1}N^{\dagger}P^{-1} = \left[ \begin{array}{cr} A^{\dagger} + B^{\dagger}  & - B^{\dagger} \\  - B^{\dagger}  & B^{\dagger} \end{array} \right], \hfill
\cr}
$$ 
if and only if $R(A)  = R(B)$ and $R(A^*) = R(B^*)$. } 

\medskip

We leave the verification of Theorem 9.21 to the reader. For the block matrix $ M $ in Theorem 9.21, we can also factor 
it, according to (9.12),  as 
$$
\displaylines{
\hspace*{2cm}
 M = PNQ = \left[ \begin{array}{cc} I_m & 0 \\ AA^{\dagger} & I_m \end{array} \right] \left[ \begin{array}{cc} A & 0 
\\  0  & B  \end{array} \right] \left[ \begin{array}{cc} I_n & A^{\dagger}A \\ 0 & I_n \
\end{array} \right]. \hfill
\cr}
$$
In that case
$$
\displaylines{
\hspace*{2cm} 
r(\, M^{\dagger} - Q^{-1}N^{\dagger}P^{-1} \, ) =  r\left[ \begin{array}{c} A  \\ B  \end{array} \right] +
 r[ \, A, \ B \, ] -  2r(B).   \hfill
\cr}
$$ 
Thus, we see that the Moore-Penrose inverse inverse of $ M^{\dagger}$ can be expressed as 
$$ 
\displaylines{
\hspace*{2cm}
 M^{\dagger} = Q^{-1}N^{\dagger}P^{-1} = \left[ \begin{array}{cc} A^{\dagger} 
+ A^{\dagger}AB^{\dagger}AA^{\dagger}  & - A^{\dagger}AB^{\dagger} \\  - B^{\dagger}AA^{\dagger}  & B^{\dagger} 
\end{array} \right], \hfill
\cr}
$$ 
if and only if $R(A) \subseteq R(B)$ and $R(A^*) \subseteq R(B^*)$.  

Another interesting example is concerning  the Moore-Penrose inverse of the $ k\times k$  block matrix
$$
M =  \left[ \begin{array}{cccc} A & B & \cdots  & B \\  B & A & \cdots  & B \\  
\vdots & \vdots & \ddots & \vdots \\ B & B & \cdots & A  \end{array} \right]_{ k\times k}, \eqno (9.76)
$$
where both $A$ and $B$ are $ m \times n$ matrices. It is easy to verify that 
$$ 
\displaylines{
\hspace*{1cm}
M = P_mNQ_n =  \left[ \begin{array}{cccc} I_m & -I_m & \cdots  & -I_m \\  I_m  & I_m & \cdots  & 0 \\  
\vdots & \vdots & \ddots & \vdots \\ I_m & 0 & \cdots & I_m  \end{array} \right]
 \left[ \begin{array}{cccc} A + (k - 1)B  &  &   &  \\   & A - B  &   &  \\  
 &  & \ddots &  \\  &  &  & A - B  \end{array} \right] \hfill
\cr
\hspace*{4cm}
\times \left[ \begin{array}{cccc} I_n/k & I_n/k  & \cdots  & I_n/k \\  -I_n/k  & (k-1)I_n/k & \cdots  & -I_n/k \\  
\vdots & \vdots & \ddots & \vdots \\ -I_n/k &  -I_n/k & \cdots & (k-1)I_n/k \end{array} \right], \hfill (9.77)
\cr}
$$
where $ P $ and $ Q $ are nonsingular, and both of them satisfy 
$$ 
P^*P = \left[ \begin{array}{cccc} kI_m & 0 & \cdots  & 0 \\  0  & 2I_m & \cdots  & I_m \\  
\vdots & \vdots & \ddots & \vdots \\ 0 & I_m & \cdots & 2I_m  \end{array} \right],  \ \ \
QQ^* = \left[ \begin{array}{cccc} I_n/k & 0 & \cdots  & 0 \\  0  & (k-1)I_n/k & \cdots  & -I_n/k \\  
\vdots & \vdots & \ddots & \vdots \\ 0 & -I_n/k & \cdots & (k-1)I_n/k  \end{array} \right]. 
$$
Note that $ N $ is diagonal. Then it is easy to  verify that $ r[ \, N , \ P^*PN \, ] =  r(N) $ and 
$ r\left[ \begin{array}{c} N  \\ NQQ^*  \end{array} \right]  = r(N)$. Thus we have 
$ M^{\dagger} = Q^{-1}_nN^{\dagger}P^{-1}_m$ according to (8.20). Furthermore, one can verify that 
$P_mQ_m = I_{km}$ and $P_nQ_n = I_{kn}$ when $ m = n$. Hence we can write $ M^{\dagger}$ as  
$ M^{\dagger} = P_nN^{\dagger}Q_m$. Written in an explicit form 
$$
M^{\dagger}=  \left[ \begin{array}{cccc} S & T & \cdots  & T \\  T & S & \cdots  & T \\  
\vdots & \vdots & \ddots & \vdots \\ T & T & \cdots & S  \end{array} \right]_{ k\times k}, \eqno (9.78)
$$ 
where 
$$
S = \frac{1}{k}[\, A + (k-1)B \,]^{\dagger} +  \frac{k-1}{k}(\, A - B \,)^{\dagger} , \qquad  
T = \frac{1}{k}[\, A + (k-1)B \,]^{\dagger} -  \frac{1}{k}(\, A - B \,)^{\dagger}. \eqno (9.79)
$$   
The expression (9.78) illustrates that $M^{\dagger}$ has the same pattern as $M$. To find  $M^{\dagger}$, 
what we actually need to do is to find $ [\, A + (k-1)B \,]^{\dagger}$ and $ (\, A - B \,)^{\dagger},$ and then put them in (9.78).  Some interesting subsequent results can be derived from (9.77) and (9.78). For example   
$$\displaylines{
\hspace*{0cm}
 r(\, MM^{\dagger} - M^{\dagger}M \,) = 2r[\,  A + (k-1)B, \, ( \,A + (k-1)B\,)^* \,] + 
2(k-1)r[\, A - B,  \, ( \,A - B \,)^* \,]- 2 r[\, A + (k-1)B \,]  \hfill
\cr
\hspace*{2.5cm}
- 2(k-1)r( \,A - B \,).  \hfill (9.80)           
\cr}
$$
In particular, $ M $ in (9.76) is EP if and only if both $ A + (k-1)B$ and $  A - B$ are EP. We leave the verification of the result to the reader. One can also find $r[\, (MM^{\dagger})(M^{\dagger}M) - 
(M^{\dagger}M)(MM^{\dagger}) \,]$ and  $ r(\, M^*M^{\dagger} - 
M^{\dagger}M^* \,)$ and so on for $ M $ in (9.76).

A more general work than (9.78) is to consider Moore-Penrose inverses of block circulant matrices. This topic was examined by Smith in \cite{Sm1} and some nice properties on   Moore-Penrose inverses of block circulant matrices  were presented there. Much similar to what we have done for $M$ in (9.76), through block factorization, we can also simply find a general expression for Moore-Penrose inverses of block circulant matrices, and derive from them various consequences. We shall present the corresponding results in Chapter 11.   

  
\medskip

\noindent {\bf Remark.}\, It should be pointed out that many results similar to those in Theorems 9.20 and 9.21, as well as in (9.78), (9.83)---(9.86), can be trvially established. 
In fact, properly choosing block matrices $ P, \ N$ and  $Q$  and then applying Theorems 
8.12, 8.13 and 8.14 to them, one can find out various rank equalities related to the Moore-Penrose inverses of 
the block matrices. Based on those rank equalities, one can further derive necessary and sufficient conditions for  
$(PNQ)^{\dagger} = Q^{\dagger}N^{\dagger}P^{\dagger}$ or $(PNQ)^{\dagger} = Q^{-1}N^{\dagger}P^{-1}$ to hold for 
these block matrices. We hope the reader to try this method and find some more interesting or unexpected conclusions 
about  Moore-Penrose inverses of block matrices.              

In addition to the methods mentioned above for finding  Moore-Penrose inverses of block matrices, 
   another possible tool is the identity (8.25) for the Moore-Penrose inverse of  product of 
three matrices.

\markboth{YONGGE  TIAN }
{10. RANK EQUALITIES FOR MOORE-PENROSE INVERSES OF  SUMS OF MATRICES}

\chapter{Rank equalities for Moore-Penrose inverses of sums of matrices}

\noindent In this chapter, we establish rank equalities related to 
 Moore-Penrose inverses of sums of matrices and consider their various 
consequences.  

\medskip

\noindent  {\bf Theorem 10.1.}\, {\em Let $ A, \, B
\in {\cal C}^{ m \times n}$ be given and let $ N = A + B $. Then
$$\displaylines{
\hspace*{1cm}
 r[ \, N - N( \, A^{\dagger} + B^{\dagger} \, )N  \,]
 = r \left[ \begin{array}{cc} AB^*A& AA^*B + AB^*B \\
  BA^*A + BB^*A & BA^*B \end{array} \right] + r(N) - r(A) - r( B ). \hfill (10.1)
\cr
\hspace*{0cm}
In \ particular, \hfill
\cr
\hspace*{1cm}
A^{\dagger} + B^{\dagger} \in \{ (A + B )^- \}
 \Leftrightarrow r \left[ \begin{array}{cc} AB^*A & AA^*B + AB^*B \\
  BA^*A + BB^*A & BA^*B \end{array} \right]  = r(A) + r( B ) - r(N). 
 \hfill (10.2)
 \cr}
$$ } 
{\bf Proof.}\,  It follows by  (2.2) and block elementary  operations that 
$$\displaylines{
\hspace*{1cm}
 r[ \, N - N( \, A^{\dagger} + B^{\dagger} \, )N  \,]  \hfill
\cr
\hspace*{1cm}
 = r \left[ \begin{array}{ccc}
  A^*AA^*  & 0 & A^*N  \\ 0 & B^*BB^*  &  B^*N  \\ NA^* &  NB^* & N
  \end{array} \right] - r( A )  - r( B ) \hfill
\cr
\hspace*{1cm}
  =  r \left[ \begin{array}{ccc} -A^*BA^* &  - A^*AB^* - A^*BB^*  & 0 
\\ - B^*AA^* - B^*BA^* & -B^*AB^* & 0 \\
  0 & 0 & N  \end{array} \right] - r( A )  - r( B )  \hfill
\cr
\hspace*{1cm}
 = r \left[ \begin{array}{cc} AB^*A & AA^*B + AB^*B \\
  BA^*A + BB^*A & BA^*B \end{array} \right] + r(N) - r(A) - r( B ). \hfill
\cr}
$$
Thus we have  (10.1) and (10.2).   \qquad  $\Box$ 

\medskip

A general result is given below. 

\medskip

\noindent  {\bf Theorem 10.2.}\, {\em Let $ A_1, \ , A_2, \, \cdots, \,
 A_k \in {\cal C}^{ m \times n}$ be given and let
 $ A =  A_1 + A_2 + \cdots + A_k ,$
 $ X =  A_1^{\dagger} + A_2^{\dagger} + \cdots  +A_k^{\dagger}$. Then
$$\displaylines{
\hspace*{1.5cm}
 r( \, A - AXA \,) = r( \, DD^*D - PA^*Q \, )  - r(D) + r(A), \hfill (10.3)
\cr 
\hspace*{0cm}
where \hfill
\cr
\hspace*{1cm}
 D= {\rm diag}( \, A_1, \, A_2, \, \cdots ,\, A_k \,),
  \ \  \ P^* = [ \,  A_1^*, \ A_2^*, \ \cdots, \ A_k^* \,], \ \ \
  Q = [ \, A_1, \, A_2, \, \cdots , \, A_k \, ]. \hfill
 \cr
\hspace*{0cm}
In \ particular, \hfill
\cr
\hspace*{1cm}
 X \in \{ A^- \}  \Leftrightarrow  r( \, DD^*D - PA^*Q \, )
 = r(D) - r(A), \  \ i.e., \  \  PA^*Q \leq_{rs} DD^*D. \hfill (10.4)
 \cr}
$$}
{\bf Proof.}\, Let $ P_1 = [\, I_n, \, \cdots, \, I_n \,]$ and
$Q_1 = [\, I_m, \, \cdots, \, I_m \,]^*$. Then $ X = P_1D^{\dagger}Q_1$. In that case,
it follows by  (2.1) that 
\begin{eqnarray*}
 r( \, A - AXA \, ) &= & r( \, A - AP_1D^{\dagger}Q_1A  \,) \\
 & = & r \left[ \begin{array}{cc}
  D^*DD^*  & D^*Q_1A  \\ AP_1D^*  &  A
  \end{array} \right] - r( D ) \\
& = & r \left[ \begin{array}{cc}  D^*DD^* - D^*Q_1AP_1D^*  & 0  \\
 0 &  A \end{array} \right] - r( D ) \\
& = & r (\, D^*DD^* - D^*Q_1AP_1D^* \, ) + r(A) - r( D ) \\
& = & r (\, DD^*D - DP_1^*A^*Q_1^*D \, ) + r(A) - r( D )
\end{eqnarray*}
as required for (10.3).   \qquad  $\Box$ 

\medskip

\noindent  {\bf Theorem 10.3.}\, {\em Let $ A, \, B  \in
{\cal C}^{ m \times n}$ be given and let $ N = A + B $. Then
 $$\displaylines{
\hspace*{1.5cm}
  r( \, N^{\dagger} - A^{\dagger} - B^{\dagger} \, )
 =  r \left[ \begin{array}{cccc}
  -NN^*N  & 0 & 0 & N  \\ 0 & AA^*A  & 0 & A  \\ 0 & 0 & BB^*B & B \\
   N & A & B & 0 \end{array} \right]
  - r(N) - r(A) - r(B).   \hfill (10.5)
\cr
\hspace*{0cm}
In \ particular, \hfill
\cr
\hspace*{1.5cm}
  N^{\dagger}  = A^{\dagger} + B^{\dagger} \Leftrightarrow
  r \left[ \begin{array}{cccc}
  -NN^*N  & 0 & 0 & N  \\ 0 & AA^*A  & 0 & A  \\ 0 & 0 & BB^*B & B \\
  N & A & B & 0   \end{array} \right]  = r(N) + r(A) + r(B). \hfill (10.6)
\cr}
$$}
{\bf Proof.}\, Follows immediately from (2.7).
\qquad  $\Box$ 

\medskip

It is well known that for any two nonsingular matrices $ A $ and $B$, there
 always is $ A(\, A^{-1} + B^{-1} \, )B = A + B$. Now for Moore-Penrose
  inverses of matrices we have the following. 

\medskip

\noindent  {\bf Theorem 10.4.}\, {\em Let $ A, \, B  \in {\cal C}^{ m \times n}$ be
 given. Then 
 $$\displaylines{
\hspace*{1.5cm}
  r[ \, A + B  -  A( \, A^{\dagger} + B^{\dagger} \, )B \, ]
 = r \left[ \begin{array}{c} A \\  B  \end{array} \right] +
 r[ \, A,  \ B \, ] - r( A ) - r( B ),  \hfill (10.7)
\cr
\hspace*{0cm}
 and \hfill
\cr
\hspace*{1.5cm}
 r[  \,  A^{\dagger} + B^{\dagger} -
  A^{\dagger}( A +  B )B^{\dagger} \, ]
  = r \left[ \begin{array}{c} A \\  B  \end{array} \right] +
  r[ \, A,  \ B \, ] - r( A ) - r( B ).  \hfill (10.8)
 \cr
\hspace*{0cm}
In  \ particular, \hfill
\cr
\hspace*{0.5cm}
 A( \, A^{\dagger} + B^{\dagger} \, )B
 = A + B \Leftrightarrow  A^{\dagger}( A +  B )B^{\dagger}
 = A^{\dagger} + B^{\dagger}  \Leftrightarrow  R(A) =  R(B)  \ \ and
  \ \ R(A^*) = R(B^*).  \hfill(10.9)
\cr}
 $$}
{\bf Proof.}\, Writing
 $$
 A + B -  A( \, A^{\dagger} +  B^{\dagger} \, )B   = A + B - [\, A , \ A \, ]
 \left[ \begin{array}{cc} A  & 0  \\  0 & B  \end{array} \right]^{\dagger}
 \left[ \begin{array}{c} B \\ B \end{array} \right]
$$
 and then applying  (2.1) to it produce  (10.7). Replacing $ A $ and
 $ B $ in (10.7) respectively by $ A^{\dagger}$ and $ B^{\dagger}$ leads
 to (10.8). The equivalences in  (10.9) follow from  (10.7) and
 (10.8). \qquad $ \Box$ 

\medskip

\noindent {\bf Theorem 10.5.}\, {\em Let $ A, \, B  \in {\cal C}^{ m \times n}$ be
 given and let $ N = A + B $. Then
$$\displaylines{
\hspace*{0cm}
 r[ \, N - N( \, (E_BAF_B)^{\dagger} + (E_ABF_A)^{\dagger} \, )N  \,]
 = r( N ) + 2r(A) + 2r(B) - r \left[ \begin{array}{cc} A  & B \\  B & 0
   \end{array} \right]
  - r \left[ \begin{array}{cc} B  & A \\  A & 0  \end{array} \right].
  \hfill (10.10)
\cr
\hspace*{0cm}
In  \ particular, \hfill
\cr
\hspace*{1cm}
 (E_BAF_B)^{\dagger} + (E_ABF_A)^{\dagger} \in \{ (A + B )^- \}
 \Leftrightarrow  r( \, A + B  \,) = r(E_BAF_B) + r(E_ABF_A).
  \hfill (10.11)
\cr}
 $$}

\noindent {\bf Proof.}\,  Let $ P = E_BAF_B$ and $ Q = E_ABF_A$. Then it is easy to
 verify that
 $$
 P^*B = 0, \ \ \ BP^* = 0, \ \ \  Q^*A = 0, \ \ \  AQ^* = 0, \ \  \ P^*PP^* = P^*AP^*, \  \ \ Q^*QQ^* = Q^*BQ^*.
 $$
 Thus we find by  (2.2) that
 \begin{eqnarray*}
 r[ \, N - N( \, P^{\dagger} + Q^{\dagger} \, )N \,]  & = & r \left[ \begin{array}{ccc}
  P^*PP^*  & 0 & P^*N  \\ 0 & Q^*QQ^*  &  Q^*N  \\ NP^* &  NQ^* & N
   \end{array} \right] - r( P )  - r( Q ) \\
  & = & r \left[ \begin{array}{ccc} P^*AP^*  & 0 & P^*A  \\ 0 & Q^*BQ^*  &
   Q^*B  \\ AP^* &  BQ^* & A + B  \end{array} \right] - r( P )  - r( Q ) \\
 & = & r \left[ \begin{array}{ccc} 0  & 0 & 0  \\ 0 & 0  & 0  \\ 0 & 0 &
 A + B  \end{array} \right] - r( P )  - r( Q ) = r(N) - r( P )  - r( Q ),
 \end{eqnarray*}
 where
 $$ 
  r( P ) = r \left[ \begin{array}{cc} A  & B \\  B & 0 \end{array} \right]
  - 2r(B), \ \ \  \  r(Q) = r \left[ \begin{array}{cc} B  & A \\ A  & 0
   \end{array} \right]- 2r(A),
 $$
 by  (1.4). Hence we have (10.10) and (10.11). \qquad  $\Box$ 

\medskip

\noindent  {\bf Theorem 10.6.}\, {\em Let $ A, \, B  \in
{\cal C}^{ m \times n}$ be given. Then 

{\rm (a)} \  $ r( \, (\, A + B\,)(\, A + B\,)^{\dagger}
  -  [\, A, \ B \,] [\, A, \ B\, ] ^{\dagger} \, )
  = r[ \, A,  \ B \, ] - r( \, A  + B \, ). $

{\rm (b)} \ $  r\left( \, (\, A + B \, )^{\dagger}(\, A + B \, ) -
\left[ \begin{array}{c} A \\  B  \end{array} \right]^{\dagger}
 \left[ \begin{array}{c} A \\  B  \end{array} \right] \right)  =
  r \left[ \begin{array}{c} A \\  B  \end{array} \right] -
  r(\, A + B\, ). $ \\
In particular$,$

{\rm (c)} \  $ (\, A + B\,)(\, A + B\,)^{\dagger}
= [\, A, \ B \,] [\, A, \ B\, ] ^{\dagger} \Leftrightarrow r[ \, A,  \ B \, ] = r( \, A  + B \, ) 
\Leftrightarrow R(A) \subseteq R( \, A  + B \, )$  and $ R(B) \subseteq R( \, A  + B \, )$.

{\rm (d)} \ $ (\, A + B \, )^{\dagger}(\, A + B \, ) =
\left[ \begin{array}{c} A \\  B  \end{array} \right]^{\dagger}
  \left[ \begin{array}{c} A \\  B  \end{array} \right]  \Leftrightarrow  
r \left[ \begin{array}{c} A \\  B  \end{array} \right] =
  r(\, A + B\, ) \Leftrightarrow
  R(A^*) \subseteq R( \, A^*  + B^* \, )$  and $ R(B^*) \subseteq
  R( \, A^*  + B^* \, )$.
}  

\medskip

\noindent  {\bf Proof.}\, Let $ N = A + B $ and $ M = [\, A, \ B\, ]$. Then
it follows from Theorem 7.2(a) that
\begin{eqnarray*}
r(\, NN^{\dagger} - MM^{\dagger}  \, ) &  = &
2r[\, N, \ M \, ] -r(N ) - r(M ) \\
& = & 2 r[\, A + B, \ A, \ B \,] - r(\, A + B \, ) - r(A) - r(B) \\
& = & r[\, A, \ B \, ] - r(A) - r(B ),
\end{eqnarray*}
as required for Part (a). Similarly we have Part (b). \qquad $ \Box$ 

\medskip

In general we have the following. 

\medskip

\noindent  {\bf Theorem 10.7.}\, {\em Let $ A_1, \,  A_2, \, \cdots, \,
 A_k \in {\cal C}^{ m \times n}$ be given and let
 $ A =  A_1 + A_2 + \cdots + A_k ,$
 $ M = [\,  A_1, \  A_2, \ \cdots, \ A_k \,]$ and
 $ N^* = [\, A_1^*, \  A_2^*, \ \cdots, \ A_k^* \,]$  . Then

{\rm (a)} \ $ r( \, AA^{\dagger} - MM^{\dagger} \,) = r(M) - r(A).$

{\rm (b)} \ $ r( \, A^{\dagger}A - N^{\dagger}N \,) = r(N) - r(A).$

{\rm (c)} \ $ AA^{\dagger} = MM^{\dagger} \Leftrightarrow r(M) = r(A) 
 \Leftrightarrow  R(A_i) \subseteq R(M), \ i = 1, \ 2,\ \cdots, \ k.$

{\rm (d)}  $A^{\dagger}A = N^{\dagger}N  \Leftrightarrow r(N) = r(A) \Leftrightarrow
R(A_i^*) \subseteq R(N^*), \ i = 1, \ 2,\ \cdots, \ k.$ } 

\medskip

\noindent  {\bf Theorem 10.8.}\, {\em Let $ A, \, B  \in {\cal C}^{ m \times n}$ be
 given and let $ N = A + B $. Then

{\rm (a)} \  $ r(AN^{\dagger}B ) = r(NA^*) + r(B^*N ) - r(N).$

{\rm (b)} \ $ r(\, AN^{\dagger}B \, ) = r(A) + r(B) - r(N),$  \ \ if $ R(A^*)
  \subseteq R( N^*) \ \ and \ \  R(B) \subseteq R(N).$

 {\rm (c)} \ $ r(\, AN^{\dagger}B - BN^{\dagger}A \,)
 = r \left[ \begin{array}{c} N \\ NA^* \end{array} \right] +
  r[ \, N,  \ AN^* \,] - 2r(N).$  

 {\rm (d)} \  $ AN^{\dagger}B = 0 \Leftrightarrow
 r(NA^*) + r(B^*N) = r(N).$

{\rm (e)} \ $ AN^{\dagger}B =  BN^{\dagger}A  \Leftrightarrow
  R(AN^* ) \subseteq R( N ) \ \ and \ \  R( A^*N ) \subseteq R( N^* ).$

{\rm (f)} \  $ AN^{\dagger}B =  BN^{\dagger}A,$  \ \  if $ R(A ) \subseteq
 R( N ) \ \ and \ \  R( A^* ) \subseteq R( N^* ).$ } 

\medskip

\noindent  {\bf Proof.}\, It follows by (2.1) that 
 \begin{eqnarray*}
  r(AN^{\dagger}B ) & = & r \left[ \begin{array}{cc} N^*NN^*   & N^*B  \\ AN^* & 0
  \end{array} \right] 
 - r( N ) \\
 & = & r \left[ \begin{array}{cc} N^*AN^* + N^*BN^* & N^*B  \\ AN^* & 0
 \end{array} \right] - r( N ) \\
 & = & r \left[ \begin{array}{cc} 0  & N^*B \\ AN^* &  0 \end{array} \right]
 - r( N ) =  r(NA^*) + r(B^*N ) - r(N),
 \end{eqnarray*}
as required for Part (a). Under $ R(A^*) \subseteq R(N^*)$ and $ R(B) \subseteq R(N),$ 
we know that
 $ r(NA^*) = r(A)$ and $ r(B^*N) = r(B)$. Thus we have Part (b). Similarly it follows by (2.1) that
 \begin{eqnarray*}
 \lefteqn{r ( \, AN^{\dagger}B - BN^{\dagger}A \,) } \\
  & = & r \left( \, [ \, A, \ B \,] \left[ \begin{array}{cr} N
  &  0 \\ 0  & -N  \end{array} \right]^{\dagger}
  \left[ \begin{array}{c} B \\ A \end{array} \right]\, \right) \\
 & = & r \left[ \begin{array}{ccc} N^*NN^*   & 0 & N^*B  \\ 0 & -N^*NN^*
  & N^*A \\ AN^* & BN^* & 0 \end{array} \right]  - 2r( N ) \\
 & = & r \left[ \begin{array}{ccc} N^*AN^*   &  N^*BN^*  & N^*B  \\
  -N^*AN^* & -N^*BN^*  & N^*A \\ AN^* & BN^* & 0  \end{array} \right] -
  2r(N) \\
 & = & r \left[ \begin{array}{ccc} 0 & 0 & N^*B  \\ 0 & 0  & N^*A \\ AN^* &
  BN^* & 0 \end{array} \right]  - 2r( N ) \\
 & = & r \left[ \begin{array}{c} N^*B \\   N^*A  \end{array} \right] +
 r[ \, AN^*,  \ BN^* \,] - 2r( N ) \\
 & = & r \left[ \begin{array}{c} N^*N \\   N^*A  \end{array} \right] +
 r[ \, AN^*,  \ NN^* \, ] - 2r( N ) = r \left[ \begin{array}{c} N \\
  N^*A  \end{array} \right] + r[ \, AN^*, \ N \, ] - 2r( N ),
 \end{eqnarray*}
 as required for  Part (c).  \qquad $\Box$ 

\medskip

It is well known that if $ R(A^*) \subseteq R( N^* )$ and $ R(B )
  \subseteq R( N ),$ the product
 $ A( \, A + B \, )^{\dagger}B$ is called the parallel sum of $ A $ and $ B $
  and often denoted by $P(A, \, B)$. The results in Theorem 10.8(b) and (f)
  show that if $ A $ and $ B $ are parallel summable, then
 $$
   r[\,P(A, \, B) \,] = r( A ) + r( B) - r( \, A + B  \,)  \ \ {\rm and }  \ \
  P(A, \, B) = P(B, \, A).
  $$
 These two properties were obtained by Rao and Mitra \cite{RM} with a different method. 

\medskip

The following three theorems are derived directly from  (2.1). Their proofs 
are omitted here.

\medskip

\noindent {\bf Theorem 10.9.}\, {\em Let $ A, \, B
\in {\cal C}^{ m \times n}$ be given. Then
$$\displaylines{
\hspace*{2cm}
  r \left(  \, \left[ \begin{array}{cc} A  & 0 \\ 0 & B   \end{array} \right]  -
  \left[ \begin{array}{c} A \\ B  \end{array} \right] ( \, A + B \, )^{\dagger}
 [ \, A, \ B \,] \,   \right) =  r(A) + r(B) - r( \, A + B  \,).
 \hfill (10.12)
 \cr
\hspace*{0cm}
In  \ particular, \hfill
\cr
\hspace*{2cm}
 \left[ \begin{array}{c} A \\ B  \end{array} \right]
  (\, A + B  \, )^{\dagger} [\, A,  \ B \, ] =
  \left[ \begin{array}{cc} A  & 0 \\ 0 & B   \end{array} \right]
  \Leftrightarrow r( \, A + B  \,) = r(A) + r(B). \hfill (10.13)
\cr}
$$} 

The equivalence  in (10.13) was established by Marsaglia and Styan \cite{MaS2}. 

\medskip

 \noindent {\bf Theorem 10.10.}\, {\em Let $ A_1, \,  A_2, \, \cdots, \,  A_k \in
 {\cal C}^{ m \times n}$ be given and denote
 $$\displaylines{
\hspace*{2cm}
 A = {\rm diag}( \, A_1, \ A_2, \ \cdots , \ A_k \, ), \qquad
  N =  A_1 + A_2 + \cdots  + A_k. \hfill
 \cr
\hspace*{0cm}
Then \hfill
\cr
\hspace*{2cm}
  r \left( \, A  - \left[ \begin{array}{c} A_1 \\ \vdots \\ A_k \end{array}
  \right]N^{\dagger} [ \, A_1,  \, \cdots, \, A_k \, ]  \,  \right)
  = r(A_1) + \cdots +  r(A_k) - r(N).   \hfill (10.14)
\cr
\hspace*{0cm}
 In \ particular, \hfill
\cr
\hspace*{2cm}
 \left[ \begin{array}{c} A_1 \\ \vdots \\ A_k \end{array} \right]N^{\dagger}
 [ \, A_1,  \, \cdots, \, A_k \, ] = A  \Leftrightarrow
 r(N) = r(A_1) + \cdots + r(A_k).  \hfill (10.15)
\cr}
 $$} 

The equivalence  in (10.15) was established by Marsaglia and Styan \cite{MaS2}. 

\medskip

\noindent  {\bf Theorem 10.11.}\, {\em Let $ A, \ B  \in {\cal C}^{ m \times n}$ be given
  and let $ N = A + B $. Then

 {\rm (a)} \ $ r( \, A - AN^{\dagger}A ) = r(NB^*N) + r(A) - r(N).$

 {\rm (b)} \ $ r(\, A - AN^{\dagger}A \, ) = r(A) + r(B) - r(N),$  \ \  
if \ $ R(A ) \subseteq R( N ) \ \ and \ \  R( A^* ) \subseteq R( N^* ).$ 

 {\rm (c)} \ $ N^{\dagger} \in  \{ A^- \}  \ \Leftrightarrow \ r(NB^*N) = r(N) - r(A).$ 

 {\rm (d)} \ $ N^{\dagger} \in  \{ A^- \}$  \  if $ r(N) = r(A) +r(B).$ } 

\medskip

\noindent  {\bf Proof.}\, Immediate by  (2.1). \qquad  $\Box$

\medskip

Notice that the sum $ A + B $ can be expressed as 
$$ 
A + B = [\, A, \ I \,]\left[ \begin{array}{c} I \\ B \end{array} \right] 
= [\, AA^{\dagger}, \ B \,]\left[ \begin{array}{c} A \\ B^{\dagger}B \end{array} \right] 
= [\, A, \ B \,]\left[ \begin{array}{c} A^{\dagger}A \\ B^{\dagger}B \end{array} \right]
= [\, A, \ B \,]\left[ \begin{array}{cc} A^{\dagger} & 0 \\ 0 & B^{\dagger} \end{array} \right] 
\left[ \begin{array}{c} A \\ B \end{array} \right].    
$$
Through the results in Chapter 8, we can also establish some rank equalities for sums of matrices. 

\medskip 

\noindent {\bf Theorem 10.12.}\, {\em Let $ A, \, B \in {\cal C}^{ m \times n}$ be given. Then

{\rm (a)} \ $ r \left(  \,   (\, A + B\,)  -  (\, A + B \, )\left[ \begin{array}{cc} I_n  \\ B \end{array} \right]^{\dagger}
 [ \, A, \ I_m \,]^{\dagger}  (\, A + B\,) \right) =  r(\, I_n - A^*B \,) + r( \, A + B  \,) - n.
$

{\rm (b)} \ $ \left[ \begin{array}{cc} I_n  \\ B \end{array} \right]^{\dagger}
 [ \, A, \ I_m \,]^{\dagger}\subseteq \{\, ( \, A + B  \,)^- \,\}  \Leftrightarrow 
 r(\, I_n - A^*B \,) + r( \, A + B  \,) = n. $ } 

\medskip

\noindent  {\bf Proof.}\, Write $ A + B = [\, A, \ I_m \,]\left[ \begin{array}{c} I_n \\ B \end{array} \right] = PQ$.
Then we find by (8.1) that
\begin{eqnarray*}
r(\, PQ - PQQ^{\dagger}P^{\dagger}PQ \,) & = &  r[\, P^*, \ Q \,] + r(PQ) - r(P) -r(Q) \\
& = & r\left[ \begin{array}{cc} A^* & I_n  \\  I_m & B \end{array} \right] +  r(\, A + B\,) - r[\, A, \ I_m \,] - 
r\left[ \begin{array}{c} I_n \\ B \end{array} \right] \\
& = &r(\, I_n - A^*B \,) + r( \, A + B  \,) - n.
\end{eqnarray*}
Parts (a) and (b) follow from it. \qquad $\Box$ 

\medskip

\noindent {\bf Theorem 10.13.}\, {\em Let $ A, \, B \in {\cal C}^{ m \times n}$ be given. Then

{\rm (a)} \  $ r \left(  \, (\, A + B \,)  -  (\, A + B \, )\left[ \begin{array}{c} A  \\ 
B^{\dagger}B \end{array} \right]^{\dagger}[ \, AA^{\dagger}, \  B \,]^{\dagger}(\, A + B\,) \right) = 
r\left[ \begin{array}{cc} A^* & A^*A \\ BB^* & B \end{array} \right] + r(\, A + B \,) - r[\, A, \ B \,]
- r\left[ \begin{array}{cc} A  \\ B \end{array} \right].$

{\rm (b)} \ $ r \left(  \,   (\, A + B \,)  -  (\, A + B \, )\left[ \begin{array}{c} A 
 \\ B \end{array} \right]^{\dagger}[ \, AA^{\dagger}, \  BB^{\dagger} \,]^{\dagger}(\, A + B\,) \right) = 
r\left[ \begin{array}{cc} A & B \\ A^*A & B^*B \end{array} \right] + r(\, A + B \,) - r[\, A, \ B \,]
- r\left[ \begin{array}{cc} A  \\ B \end{array} \right].$

{\rm (c)} \ $ \left[ \begin{array}{c} A  \\ B^{\dagger}B \end{array} \right]^{\dagger}
 [ \, AA^{\dagger}, \  B \,]^{\dagger} \subseteq \{\, ( \, A + B  \,)^- \,\}  
\Leftrightarrow r\left[ \begin{array}{cc} A^* & A^*A \\ BB^* & B \end{array} \right] = 
r\left[ \begin{array}{cc} A  \\ B \end{array} \right] + r[\, A, \ B \,] - r(\, A + B \,).$

{\rm (d)} \ $ \left[ \begin{array}{c} A  \\ B \end{array} \right]^{\dagger}
 [ \, AA^{\dagger}, \  BB^{\dagger}\,]^{\dagger} \subseteq \{\, ( \, A + B  \,)^- \,\}  
\Leftrightarrow r\left[ \begin{array}{cc} A & B \\ A^*A & B^*B \end{array} \right] = 
r\left[ \begin{array}{cc} A  \\ B \end{array} \right] + r[\, A, \ B \,] - r(\, A + B \,).$ 
}

\medskip

\noindent  {\bf Proof.}\, Writing  $ A + B = [ \, AA^{\dagger}, \  B \,]\left[ \begin{array}{c} A  \\ 
B^{\dagger}B \end{array} \right] = [ \, AA^{\dagger}, \  BB^{\dagger} \,]\left[ \begin{array}{c} A  \\ 
B \end{array} \right]$, and then applying  (8.1) to them yields Parts (a) and (b). \qquad $\Box$ 

\medskip

\noindent {\bf Theorem 10.14.}\, {\em Let $ A, \, B \in {\cal C}^{ m \times n}$ be given. Then
$$\displaylines{
\hspace*{2cm}
( \, A + B  \,)^{\dagger} = \left[ \begin{array}{cc} I_n  \\ B \end{array} \right]^{\dagger}
 [ \, A, \ I_m \,]^{\dagger} \hfill (10.15)
\cr
holds \ if \ and \ only \ if \hfill
\cr
\hspace*{2cm}
(\, I_m - BA^* \,)(\, A+B \,) = (\, A+B \,)(\, I_n - B^*A\,) = 0. \hfill (10.16)
\cr}
$$ }
{\bf Proof.}\, Write $ A + B = [\, A, \ I_m \,] \left[ \begin{array}{c} I_n \\ B \end{array} \right] = PQ$.
Then we find by Theorem 8.2(e) that $ (PQ)^{\dagger} = Q^{\dagger}P^{\dagger}$ if and only if
$$\displaylines{
\hspace*{1.5cm}
r[\, P^*PQ, \ Q \,] = r(Q), \ \ {\rm and}  \ \ r\left[ \begin{array}{c} P \\ PQQ^* \end{array} \right] = r(P). 
 \hfill (10.17)
\cr}
$$
Notice that $r(Q) = n$ and $r(P) = m,$ and 
$$\displaylines{
\hspace*{2cm}
r[\, P^*PQ, \ Q \,]  = r\left[ \begin{array}{cc} A^*(\, A+B \,) & I_n  \\  
A + B & B \end{array} \right] = n + r[\,(\, I_m - BA^* \,)(\, A+B \,)\,],\hfill
\cr
\hspace*{2cm}
r\left[ \begin{array}{c} P \\ PQQ^* \end{array} \right] 
= r\left[ \begin{array}{cc} A & I_m  \\  A + B & (\, A+B \,)B^* \end{array} 
\right] = m + r[\,(\, A+B \,)(\, I_n - B^*A \,) \,]. \hfill
\cr}
$$
In that case, (10.17) reduces to  (10.16). \qquad $\Box$ 

\medskip

\noindent {\bf Theorem 10.15.}\, {\em Let $ A, \, B \in {\cal C}^{ m \times n}$ be given. Then
$$\displaylines{
\hspace*{1cm}
 r \left(  \,   (\, A + B\,)  -  (\, A + B \, )
\left[ \begin{array}{c} A \\ B \end{array} \right]^{\dagger}
\left[ \begin{array}{cc} A & 0  \\ 0  & B \end{array} \right]
 [ \, A, \ B \,]^{\dagger} (\, A + B\,) \right) \hfill
 \cr
\hspace*{1cm}
=  r \left[ \begin{array}{cc} A + B & AA^* + BB^* 
 \\ A^*A + B^*B & A^*AA^* + B^*BB^* \end{array} \right] + r( \, A + B  \,)  
- r[\, A , \ B\, ] - r \left[ \begin{array}{c} A \\ B \end{array} \right].
 \hfill(10.18)
\cr
In \ particular, \hfill
\cr
\hspace*{1cm}
\left[ \begin{array}{c} A \\ B \end{array} \right]^{\dagger}
\left[ \begin{array}{cc} A & 0  \\ 0  & B \end{array} \right]
 [ \, A, \ B \,]^{\dagger} \subseteq \{\, (\, A + B\,)^- \,\} \hfill
\cr
holds \ if \ and \ only \ if \hfill
\cr
\hspace*{1cm}
r \left[ \begin{array}{cc} A + B & AA^* + BB^* 
 \\ A^*A + B^*B & A^*AA^* + B^*BB^* \end{array} \right] 
= r \left[ \begin{array}{c} A \\ B \end{array} \right] + r[\, A , \ B\, ] -
 r( \, A + B  \,). \hfill (10.19)
\cr}
$$}
{\bf Proof.}\, Writing $ A + B = [ \, I, \  I \,]\left[ \begin{array}{cc} A & 0  \\ 
0 & B \end{array} \right]\left[ \begin{array}{c} I \\  I \end{array} \right] = PNQ$, 
and then applying  (8.2) to it yields  (10.18). \qquad $\Box$ 

\medskip

\noindent {\bf Theorem 10.16.}\, {\em Let $ A, \, B \in {\cal C}^{ m \times n}$ be given. Then
$$\displaylines{
\hspace*{2cm}
 r \left(   (\, A + B\,)^{\dagger}  -  \left[ \begin{array}{c} A \\ B \end{array} \right]^{\dagger}
\left[ \begin{array}{cc} A & 0  \\ 0  & B \end{array} \right]
 [ \, A, \ B \,]^{\dagger} \right) \hfill
\cr
\hspace*{2cm}
=  r \left( \left[ \begin{array}{cc} A  & B
 \\ B & A \end{array} \right]\left[ \begin{array}{cc} A^*  & 0
 \\ 0 & B^* \end{array} \right] \left[ \begin{array}{cc} A  & B \\ B & A \end{array} \right] 
\right) - r(\, A + B\,). \hfill (10.20)
\cr
In \ particular, \hfill 
\cr
\hspace*{0cm}
 (\, A + B\,)^{\dagger} = \left[ \begin{array}{c} A \\ B \end{array} \right]^{\dagger}
\left[ \begin{array}{cc} A & 0  \\ 0  & B \end{array} \right]
 [ \, A, \ B \,]^{\dagger}  \Leftrightarrow r \left( \left[ \begin{array}{cc} A  & B
 \\ B & A \end{array} \right]\left[ \begin{array}{cc} A^*  & 0
 \\ 0 & B^* \end{array} \right] \left[ \begin{array}{cc} A  & B \\ B & A \end{array} \right] 
\right) = r(\, A + B\,). \hfill (10.21)
\cr}
$$}
{\bf Proof.}\, Writing $ A + B = [ \, I, \  I \,]\left[ \begin{array}{cc} A & 0  \\ 
0 & B \end{array} \right]\left[ \begin{array}{c} I \\  I \end{array} \right] = PNQ$, 
and then applying  (8.4) to it yields (10.20). \qquad $\Box$ 

\medskip

The above several results can also be extended to  sums of $k $ matrices. In the remainder of this chapter, we present a set of results  related to expressions of Moore-Penrose
 inverses of Schur complements. These results have appeared in the author's recent paper \cite{Ti5}.

\medskip

\noindent  {\bf Theorem 10.17.}\, {\em Let $ A \in {\cal C}^{ m \times n},
\, B \in {\cal C}^{ m \times k}, \, C \in {\cal C}^{ l \times n}$  and
 $ D \in {\cal C}^{ l \times n}$ be given$,$ and satisfy the 
 rank additivity condition
$$
r \left[ \begin{array}{cc} A & B \\ C  & D \end{array} \right]
= r \left[ \begin{array}{c} A  \\ C  \end{array} \right] +
r \left[ \begin{array}{c} B \\ D \end{array} \right]  = r[\, A , \ B\, ]
+ r[\, C, \ D \,].   \eqno (10.22) 
$$
Then the following inversion formula  holds
$$
\displaylines{
\hspace*{2cm}
( E_{B_2}S_DF_{C_2} )^{\dagger} = A^{\dagger} + A^{\dagger}BJ^{\dagger}( D )CA^{\dagger} +  C_1^{\dagger}[ \,
S_AJ^{\dagger}(D)S_A - S_A \, ] B_1^{\dagger} \hfill
\cr
\hspace*{5cm}
 - \ A^{\dagger}B [ \, I- J^{\dagger}(D)S_A \, ]B_1^{\dagger} -
C_1^{\dagger}[ \, I - S_AJ^{\dagger}(D) \, ]CA^{\dagger},  \hfill (10.23)
\cr } 
$$
where 
$$ 
S_A = D - CA^{\dagger}B,  \qquad S_D = A- BD^{\dagger}C, \qquad J(D )
 =  E_{C_1}S_AF_{B_1},
$$
$$ 
  B_1 = E_AB, \qquad B_2 = BF_D, \qquad  C_1 = CF_A, \qquad C_2 = E_DC.
$$ } 
{\bf Proof.}\, Follows immediately from the two expressions of $M^{\dagger}$ in Theorem 9.8.
\qquad $\Box$

\medskip

The results given below are all the special cases of the general formula
 (10.23). 

\medskip

\noindent { \bf Corollary 10.18.}\,  {\em If $ A, \, B, \, C$ and $ D $ satisfy 
$$
 R \left[ \begin{array}{c}  A \cr 0 \end{array} \right]
 \subseteq R \left[ \begin{array}{cc}   A  & B \cr C & D
 \end{array} \right], \qquad and \qquad R \left[ \begin{array}{c}  A^* \cr 0 \end{array} \right]
\subseteq R \left[ \begin{array}{cc}  A^*  & C^* \cr B^* & D^*  \end{array}
\right], \eqno (10.24)
$$
and the following two conditions 
$$ 
R( CS^*_D) \subseteq R( D ),  \qquad  R( B^*S_D) \subseteq R( D^* ),
\eqno (10.25)
$$
or more specifically satisfy the four conditions 
$$ 
R( C) \subseteq R( D ), \ \ \  R(B^*) \subseteq R( D^* ),  \ \ \ R( B)
\subseteq R( S_D ),
  \ \ \  R(C^*) \subseteq R(S_D^* ), \eqno  (10.26) 
$$
then the Moore-Penrose inverse of the Schur complement
$S_D= A- BD^{\dagger}C $ satisfies the inversion formula
$$
\displaylines{
\hspace*{2cm}
( A- BD^{\dagger}C)^{\dagger} = A^{\dagger} +  A^{\dagger}BJ^{ \dagger}(
 D )CA^{\dagger}  + C_1^{\dagger}[ \, S_AJ^{\dagger}(D)S_A - S_A \, ]
 B_1^{\dagger} 
\hfill
\cr
\hspace*{5cm}
 - \ A^{\dagger}B[ \, I- J^{\dagger}(D)S_A \, ]B_1^{\dagger} -
 C_1^{\dagger}[ \, I - S_AJ^{\dagger}(D) \, ]CA^{\dagger}, \hfill(10.27)
\cr}
 $$
where $ S_A, \, B_1 , \, C_1 $ and $J( D)$  are defined in  {\rm (10.23)}. } 

\medskip

\noindent {\bf Proof.}\, It is obvious that (10.25) is equivalent to
$ ( E_DC)S_D^* = 0$ and $ S_D^*( BF_D ) = 0,$
or equivalently 
$$ 
S_D( E_DC)^{\dagger} = 0  \qquad  {\rm and } \qquad  ( BF_D )^{\dagger}S_D
= 0.  \eqno (10.28)
 $$
These two equalities clearly imply that $ S_D , \, E_DC $ and $ BF_D$
satisfy  (9.31). Hence by Lemma 9.7, we know that under (10.24)
and (10.25), $ A, \, B, \, C $ and $ D  $ naturally satisfy (10.22).
Now substituting  (10.28)
into the left-hand side of (10.23) yields $ J^{\dagger}( A ) =
( A - BD^{\dagger}C )^{\dagger} $. Hence (10.23) becomes (10.27).
Observe  that (10.28) is a  special case of (10.25), hence (10.27)
 is also true under  (10.26). \qquad$\Box $  

\medskip

\noindent {\bf Corollary 10.19.}\, {\em If $ A, \, B, \, C$ and $ D $ satisfy
  {\rm (10.27)}$,$ {\rm (10.25)} and the following two conditions
$$ \displaylines{
\hspace*{3cm}
R( CF_A) \cap R ( S_A) = \{ 0 \}  \ \ \  and \ \ \
 R[(E_AB)^*] \cap R(S_A^*) = \{ 0 \}, \hfill (10.29)
\cr
then \hfill
\cr
\hspace*{0.8cm}
( A- BD^{\dagger}C)^{\dagger} = A^{\dagger} +  A^{\dagger}BJ^{\dagger}( D )CA^{\dagger} -  A^{\dagger}B[ \, I- J^{\dagger}(D)S_A \,]B_1^{\dagger} -
 C_1^{\dagger}[ \, I - S_AJ^{\dagger}(D) \, ]CA^{\dagger}, \hfill (10.30)
\cr}
$$ 
where $ S_A, \, B_1 , \, C_1 $ and $J( D)$ are defined in  {\rm (10.23)}. } 

\medskip

\noindent {\bf Proof.}\, According to Theorem 7.8, the two conditions in  (10.27)
 imply that $ S_AJ^{\dagger}(D)S_A =  S_A.$  Hence (10.27) is
 simplified  to (10.30).  \qquad$ \Box $ 

\medskip

\noindent {\bf Corollary 10.20.}\, {\em  If $ A, \, B, \, C$ and $ D $ satisfy
 {\rm (10.24)}$,$ {\rm (10.25)} and the following two
conditions
$$\displaylines{
\hspace*{4cm}
 R( BS^*_A ) \subseteq R ( A ) \ \ \ and  \ \ \  R( C^*S_A ) \subseteq
 R( A^* ), \hfill (10.31)
\cr
then \hfill
\cr
\hspace*{3cm}
( A- BD^{\dagger}C)^{\dagger} = A^{\dagger}+
A^{\dagger}BS_A^{\dagger}CA^{\dagger} -A^{\dagger}B( E_AB)^{\dagger}
- ( CF_A)^{\dagger}CA^{\dagger}. \hfill (10.32)
\cr}
$$
where $ S_A = D - CA^{\dagger}B$.} 

\medskip

 \noindent {\bf Proof.}\,  Clearly, (10.31) is equivalent to $
 ( E_AB)S^*_A = 0 $ and $ S^*_A(CF_A) = 0, $  which can   also equivalently
  be expressed as $S_A ( E_AB)^{\dagger} = 0 $ and $ (CF_A)^{\dagger}
  S_A = 0. $ In that case, $ J( D ) = E_{C_1}S_AF_{B_1} = S_A.$ Hence
   (10.27) is simplified to (10.32).  \qquad   $ \Box $ 

\medskip

\noindent {\bf Corollary 10.21.}\, {\em If $ A, \, B, \, C$ and $ D $ satisfy 
 {\rm (10.24)}$,$ {\rm (10.25)} and the following two conditions
 $$\displaylines{
\hspace*{4cm}
 R( B ) \subseteq R ( A ) \ \ \  and \ \ \  R( C^*) \subseteq R( A^*),
 \hfill (10.33)
\cr
then \hfill 
\cr
\hspace*{4cm}
 ( A- BD^{\dagger}C)^{\dagger} = A^{\dagger}+A^{\dagger}B ( D -
 CA^{\dagger}B )^{\dagger}CA^{\dagger}.  \hfill (10.34)
\cr}
$$}
{\bf Proof.}\,  The two inclusions in  (10.33) are equivalent
to $ E_AB = 0$ and $ CF_A= 0 $. Substituting them into (10.27) yields
(10.34). \qquad $ \Box $ 

\medskip

\noindent {\bf Corollary 10.22.}\, {\it If $ A, \, B, \, C$ and $ D $ satisfy
the following four conditions
$$\displaylines{
\hspace*{1cm}
 R( A) \cap R ( B) = \{ 0 \}, \ \ \  R(A^*) \cap R( C^* ) = \{ 0 \}, \ \ \  R( C ) = R ( D ), \ \ \ 
 R( B^* ) = R( D^*), \hfill (10.35)
\cr
then  \hfill
\cr
\hspace*{2cm}
 ( A- BD^{\dagger}C)^{\dagger} = A^{\dagger} - A^{\dagger}B( E_AB)^{\dagger}
  - ( CF_A)^{\dagger}CA^{\dagger} + ( CF_A)^{\dagger}S_A( E_AB)^{\dagger}.
\hfill (10.36) 
\cr}
$$} 
{\bf Proof.}\, Under  (10.35),  $ A, \, B, \, C$
and $ D $ naturally satisfy the rank additivity condition in  (10.22).
Besides, from (10.35) and Theorem 7.2(c) and (d) we can derive 
$$
   B_1^{\dagger } B_1 = B^{\dagger}B, \quad  C_1C_1^{\dagger } =
   CC^{\dagger}, \quad B_2 = 0, \quad  C_2 = 0, \quad J( D ) = 0.
$$
Substituting them into (10.23) yields (10.36). \qquad$ \Box $ 

\medskip 

If $ D $ is invertible, or $ D = I$,  or $ B =  C = - D $, then the
inversion formula (10.23) can reduce to some other simpler forms. For simplicity, we do not 
list them here.  

\markboth{YONGGE  TIAN }
{11. MOORE-PENROSE INVERSES OF BLOCK CIRCULANT MATRICES}

\chapter{ Moore-Penrose inverses of block circulant matrices}

\noindent  Inverses or Moore-Penrose inverses of circulant matrices  and block circulant matrices is an attractive 
topic in matrix theory and lots of results can be find in the literature (see, e.g., 
 \cite{DP, DG, Sea,  Sm1, Sm2, Tr}). To find the general expression for inverses or Moore-Penrose inverses of 
circulant matrices and block circulant matrices, a best method is to use various well-known factorizations of 
circulant matrices and block circulant matrices, and then derive from them  general expressions of inverses or 
Moore-Penrose inverses the matrices. In this chapter we mainly consider Moore-Penrose inverses of block circulant 
matrices, and then then drive from them some interesting  consequences related to sums of matrices. In addition, 
we shall also consider some extension of the work to quaternion matrices. 
         
For a circulant matrix $ C $  over the complex number field $ { \cal C }$ with the form 
$$
C = \left[ \begin{array}{cccc} a_0 & a_1 & \cdots  & a_{k-1}  \\  a_{k-1} & a_0 & \cdots  & a_{k-2}  \\
 \vdots  & \vdots & \ddots  & \vdots \\ a_1 & a_2 & \cdots  & a_0 \end{array} \right], \eqno (11.1)
$$ 
the following factorization is well known (see, e.g., Davis \cite{DP})
$$ 
U^*CU = {\rm diag}( \, \lambda_1, \,  \lambda_2, \, \cdots, \,   \lambda_k \, ), \eqno (11.2)
$$ 
where $ U $ is a unitary matrix of the form
$$ 
U = ( u_{pq})_{k \times k} , \qquad  u_{pq} = \frac{1}{\sqrt{k}}\omega^{(p-1)(q-1)}, 
\qquad \omega^{k} = 1,  \ \ and  \ \  \omega \neq 1, \eqno (11.3) 
$$ 
and 
$$
 \lambda_t = a_0 + a_1 \omega^{(t-1)} + a_2 (\omega^{(t-1)})^2 + \cdots + a_{k-1}(\omega^{(t-1)})^{k-1}, \ \ \ t = 1, \, \cdots, \, k.  \eqno (11.4) 
  $$ 
It is evident that the entries in the first row and first column of $ U $ are all
 $1/\sqrt{k}$, and 
$$  
\lambda_1 = a_0 + a_1 + \cdots + a_{k-1}. \eqno (11.5) 
$$
Observe that $ U$ in (11.3) is independent of $ a_0$---$a_{k-1}$ in
 (11.1). Thus  (11.2) can directly be extended to block circulant matrix as
follows. 

\medskip

\noindent {\bf Lemma 11.1.}\,  {\em Let 
$$
 A = \left[ \begin{array}{cccc} A_1 & A_2 & \cdots  & A_{k}  \\  A_{k} & A_1 & \cdots  & A_{k-1}  \\
 \vdots  &  \vdots & \ddots  & \vdots  \\ A_2 & A_3 & \cdots  & A_1 \end{array} \right] \eqno (11.6)
$$ 
be a block circulant matrix over the complex number field $ { \cal C },$ where $ A_t \in 
 {\cal C}^{ m \times n}, $  $ t = 1, \, \cdots, \, k. $  Then  $ A $ satisfies the following factorization
 equality 
$$ 
U_m^*AU_n = \left[ \begin{array}{cccc} J_1 &   &  &     \\   & J_2 &  &    \\
  &   & \ddots &   \\    &   &  & J_k \end{array} \right], \eqno (11.7)  
$$ 
where $ U_r $  and $ U_s$ are  two block unitary matrices 
$$ 
U_m = ( u_{pq}I_m)_{k \times k} , \ \ \ \ U_n = ( u_{pq}I_n)_{k \times k}, \eqno (11.8)
$$  
$u_{pq}$  is as in $(11.3),$ meanwhile
$$ 
 J_t = A_1 + A_2 \omega^{(t-1)} + A_3 (\omega^{(t-1)})^2 + \cdots + A_{k}(\omega^{(t-1)})^{k-1},
 \ \ \ \ \ t = 1, \ \cdots, \ k.  \eqno (11.9)   
$$ 
Especially$,$ the block entries in the first block rows and  first block
columns of $ U_m $  and  $U_n $ are all scalar products of $ 1/\sqrt{k}$ with identity 
matrices$,$ and $ J_1$ is 
$$  
J_1 =  A_1 + A_2 + \cdots + A_k. \eqno (11.10) 
$$ } 

Observe that $ J_1$ in (11.7) is the sum of $ A_1, \, A_2, \, \cdots, \, A_k$.
Thus (11.7) implies that the sum $\sum_{t= 1}^kA_t$ is closely linked to
its corresponding block circulant matrix through a unitary
factorization equality.  Recall a fundamental fact in the theory of generalized  inverses of matrices
 (see, e.g., Rao and Mitra \cite{RM}) that 
$$ 
( PAQ )^{\dagger}  = Q^*A^{\dagger}P^*,   \ \ \ \ \ {\rm  if} \ P \
 {\rm and}  \  Q  \ { \rm are  \ unitary.}   \eqno (11.11)
$$     
Then from (11.7) we can directly find the following.  

\medskip

\noindent {\bf Lemma 11.2.}\,  {\em  Let $ A $ be given in {\rm (11.6)}$,$ $ U_r $
and $U_s$ be given in {\rm (11.8)}. Then the Moore-Penrose inverse of $ A $ satisfies 
$$  
A^{\dagger}= U_n{\rm diag}( \, J_1^{\dagger}, \  J_2^{\dagger}, \
\cdots, \  J_k^{\dagger} \, )U_m^*. \eqno (11.12)  
$$ }
{\bf Proof.}\,  Since $ U_m $ and $U_n$ in (11.7) are unitary,
we find by (11.11) that
$$ 
( U_m^*A U_n )^{\dagger} = U_n^*A^{\dagger}U_m.  
$$ 
On the other hand, it is easily seen that 
$$ 
[ \, {\rm diag}( \, J_1, \,  J_2, \, \cdots, \,  J_k \, ) \, ]^{\dagger} 
= {\rm diag}( \, J_1^{\dagger}, \, J_2^{\dagger}, \, \cdots, \, J_k^{\dagger}
 \, ). 
$$ 
Thus (11.12) follows.   \qquad  $\Box $  

\medskip

The expression shows that the Moore-Penrose inverse of $ A $  can be completely determined by the Moore-Penrose 
inverses of $J_1$---$J_k$. Moreover, $ A^{\dagger}$ is also a block circulant matrix, this fact was pointed out 
by Cline, Plemmons and Worm in  \cite{CPW} and Smith in \cite{Sm1}.   

The generalizations of circulants and block circulants have many forms (see, e.g. \cite{Chi, CPW, DP, Wan, Wat}), and  various factorizations of these kinds of matrices can also be established. In that case, one can use the rank formulas in Chapter 8 to those factorizations, and then find from them various expressions for 
Moore-Penroses inverses of these generalized circulants and  generalized block circulants. But we do not intend to 
go further along this direction. Instead, our next work is to consider some remarkable applications of (11.12) to
 Moore-Penrose inverses of sums of matrices.   

\medskip

\noindent {\bf Theorem 11.3}(Tian \cite{Ti2, Ti5}).\,  {\em  Let $ A_1, \, A_2, \, \cdots, \, A_k  \in {\cal C }^{ m \times n}.$  
Then the  Moore-Penrose inverse of their sum satisfies 
$$
 ( \,A_1 +  A_2 + \cdots + A_k \, )^{\dagger} 
= \frac{1}{k}[\, I_n, \,  I_n, \, \cdots, \, I_n \,] 
\left[ \begin{array}{cccc} A_1 & A_2 & \cdots  & A_{k}  \\  A_{k} & A_1 & \cdots  & A_{k-1}  \\
 \vdots  & \vdots & \ddots  & \vdots \\  A_2 & A_3 & \cdots  & A_1 \end{array} \right]^{\dagger}  
\left[ \begin{array}{c} I_m \\ I_m \\ \vdots 
 \\ I_m \end{array} \right]. \eqno (11.13) 
$$ 
In particular$,$ if the block circulant matrix in it is nonsingular$,$ then 
$$
 ( \,A_1 +  A_2 + \cdots + A_k \, )^{-1} 
= \frac{1}{k}[\, I_m, \,  I_m, \, \cdots, \, I_m \,] 
\left[ \begin{array}{cccc} A_1 & A_2 & \cdots  & A_{k}  \\  A_{k} & A_1 & \cdots  & A_{k-1}  \\
 \vdots  & \vdots & \ddots  & \vdots \\ A_2 & A_3 & \cdots  & A_1 \end{array} \right]^{-1} 
\left[ \begin{array}{c} I_m \\ I_m \\ \vdots 
 \\ I_m \end{array} \right]. \eqno (11.14) 
$$} 
{\bf Proof.}\, Pre-multiply $ [\, I_n, \,  0, \, \cdots, \, 0 \,] $ and 
 post-multiply $ [\, I_m, \,  0, \, \cdots, \, 0 \,]^T $ on the both sides 
of  (11.12) and observe the structure of $ U_m $ and $U_n$ to yield (11.13).  \qquad $\Box $ 

\medskip

In \cite{Ti2} and  \cite{Ti5}, the author proved (11.13) in some direct but tedious  methods. New (11.3) is just a 
simple consequence on the Moore-Penrose inverse of a block circulant matrix. The identity (11.13) manifests 
that the Moore-Penrose inverse of a sum of matrices can be completely determined through the Moore-Penrose inverse of 
the corresponding  block circulant matrix. In this case,  if we can find the expression of the Moore-Penrose inverse of 
the block circulant matrix by some other methods (not by (11.2)), then we can get the expression
for the Moore-Penrose inverse of the sum of matrices. In fact, we have presented many results in Chapter 9 for 
 Moore-Penrose inverses of block matrices. Applying some  of them to  the block circulant matrix in (11.13), one
 can derive many new conclusions on  Moore-Penrose inverses of  sums of matrices. Here we present some of them.

Let $ A$ and $ B $ be two  $ m \times n $  matrices. Then according to
(11.13) we have
$$ 
( \, A + B \, ) ^{\dagger} = \frac{1}{2}[ \, I_n , \ I_n \, ]
 \left[ \begin{array}{cc}  A & B \cr B & A  \end{array} \right]^{\dagger}
  \left[ \begin{array}{c}  I_m \cr I_m  \end{array} \right]. \eqno (11.15)
$$
As a special case of (11.15), if we replace $ A + B $ in (11.15) by a complex
matrix $ A + iB $, where  both $ A $ and $ B $ are  real matrices, then (11.15)
becomes the equality
$$
 ( \, A + iB \, )^{\dagger }  = \frac{1}{2}[ \, I_n , \ I_n \, ]
  \left[ \begin{array}{cc}  A & iB \cr iB & A  \end{array} \right]^{\dagger}
   \left[ \begin{array}{c}  I_m \cr I_m  \end{array} \right] =
   \frac{1}{2} [ \, I_n , \ iI_n \, ]
    \left[ \begin{array}{cc}  A & -B \cr B & A  \end{array} \right]^{\dagger}
    \left[ \begin{array}{c}  I_m \cr -iI_m  \end{array} \right]. \eqno (11.16)
$$ 

Now applying Theorems 9.8 and  9.9 to (11.15) and (11.16) we find the following
two results, which was presented by the author in \cite{Ti5}.

\medskip 
 
\noindent {\bf Theorem 11.4.}\, {\em Let $ A $ and $ B $ be two  $ m \times n
$ complex  matrices$,$ and suppose that  they satisfy the rank additivity condition
$$\displaylines{
\hspace*{2cm}
 r \left[ \begin{array}{cc}  A & B \cr B & A  \end{array} \right] =
 r \left[ \begin{array}{cc}   A \cr B  \end{array} \right] +
  r  \left[ \begin{array}{cc} B \cr A  \end{array} \right] = r[ \, A , \ B \, ] +
  r[ \, B , \ A \, ],
  \hfill (11.17)
\cr}
$$
or alternatively 
$$\displaylines{
\hspace*{2cm}
R( A ) \subseteq R( A \pm B ) \ \ and  \ \ R( A^* )
\subseteq R( A^* \pm B^* ). \hfill (11.18)
\cr}
$$ 
Then

{\rm (a)}\,  The Moore-Penrose inverse of $ A + B $ can be expressed as 
$$ 
\displaylines{
\hspace*{2cm}
( \, A + B  \, )^ {\dagger} = J^{\dagger}( A ) + J^{\dagger}( B ) =
( \, E_{B_2}S_AF_{B_1})^{\dagger} +( E_{A_2}S_BF_{A_1} \, )^{\dagger},
\hfill (11.19)
\cr}
$$
where $ J ( A ) $ and $J(B) $ are, respectively, the rank complements of
$A$ and $ B $ in $  \left[ \begin{array}{cc}  A & B \cr B & A
\end{array} \right], $  and
$$
S_A = A - BA^{\dagger}B, \ \ S_B = B - AB^{\dagger}A, \ \  A_1 = E_BA, \ \ 
 A_2 = AF_B, \ \ B_1 = E_AB, \ \ B_2= BF_A.
$$

{\rm (b)}\,  The matrices  $ A, \ B ,$ and the two terms $ G_1 = J^{\dagger}( A )$ and 
$G_2 = J^{\dagger}( B) $ in the right-hand side of {\rm (11.19)} satisfy the
following several  equalities
$$
r(G_1) = r( A ), \qquad  r( G_2 ) = r( B ),
 $$
$$
(\, A + B\, ) ( \,A + B \,)^{\dagger} = AG_1 + BG_2,  \qquad
 ( \,A + B \, )^{\dagger}( \, A + B \, ) = G_1A + G_2B,
$$
$$
 AG_2 + BG_1= 0, \qquad  G_2A + G_1B = 0. 
$$ }
{\bf Proof.}\,  The equivalence of  (11.17) and (1.18) is derived from
  (1.13). We know from Theorem 9.8 that under the condition
(11.17), the Moore-Penrose inverse of
 $ \left[ \begin{array}{cc}  A & B \cr B & A \end{array} \right] $ can be 
expressed as 
 $$
 \left[ \begin{array}{cc}  A & B \cr B & A \end{array} \right]^{\dagger} =
  \left[ \begin{array}{cc}  J^{\dagger}( A ) & J^{\dagger}( B ) \cr
  J^{\dagger}( B ) & J^{\dagger}( A ) \end{array} \right]  = 
\left[ \begin{array}{cc} ( \, E_{B_2}S_AF_{B_1})^{\dagger} & 
( E_{A_2}S_BF_{A_1} \, )^{\dagger} \\ 
( E_{A_2}S_BF_{A_1} \, )^{\dagger} &  ( \, E_{B_2}S_AF_{B_1})^{\dagger}
\end{array} \right].  
$$
Then putting it in (11.15) immediately yields (11.19). The results in 
Part (b) are derived from Theorem 9.9. \qquad $\Box $  

\medskip

\noindent{\bf Theorem 11.5.}\,  {\em Let $ A+iB $ be an $ m \times n $
complex matrix, where $ A  $ and $ B $ are two real matrices$,$ and suppose that  $ A  $
and $ B $ satisfy
$$
 r \left[ \begin{array}{cc}  A & -B \cr B & A \end{array} \right] =
  r \left[ \begin{array}{cc}  A \cr B \end{array} \right] +
 r \left[ \begin{array}{cc}  - B \cr A \end{array} \right]  =
 r[ \, A , \ -B \, ] + r[ \, B , \ A \, ], \eqno (11.20)
 $$
or equivalently 
$$
R( A ) \subseteq R ( A \pm iB ) \ \ and  \ \  R( A^*)
 \subseteq R( A^T \pm i B^T ). \eqno (11.21)
$$
Then the Moore-Penrose inverse of $ A + iB $ can be expressed as 
$$ 
( \, A +i B \, )^{\dagger} = G_1 -iG_2 = [ \, E_{B_2}( \,A +
BA^{\dagger}B \,)F_{B_1}\, ]^{\dagger} -i[\,  E_{A_2}( \, B + AB^{\dagger}A \,)
 F_{A_1} \,]^{\dagger}, \eqno (11.22)
$$
where $  A_1 = E_BA, \ A_2 = AF_B, \ B_1 = E_AB$ and $B_2= BF_A. $ } 

\medskip

\noindent {\bf Proof.}\, Follows directly from Theorem 11.4. \qquad $\Box $

\medskip

\noindent {\bf Corollary 11.6.}\, {\em Suppose that $ A + iB $  is a
nonsingular complex matrix, where $ A $ and $ B $  are real.

{\rm (a)}\,  If  both $ A $ and $ B $ are nonsingular$,$ then 
$$
\displaylines{
\hspace*{2cm}
( \,A + iB  \,)^{-1} = ( \, A + BA^{-1}B \, )^{-1}
-i( \, B + AB^{-1}A \,)^{-1}. \hfill
\cr}
$$

{\rm (b)}\,  If both $ R( A ) \cap R( B ) = \{ 0 \} $ and  $ R( A^*) \cap R( B^*)
= \{ 0 \}, $ then
$$
\displaylines{
\hspace*{2cm}
 ( \,A + iB \, )^{-1} = (E_BAF_B )^{\dagger} -i(E_ABF_A )^{\dagger}. \hfill 
\cr}
$$

{\rm (c)}\, Let $ A = \lambda I_m, $ where $ \lambda$ is a real number such
that $\lambda I_m + iB $ is nonsingular$,$ then
$$
\displaylines{
\hspace*{2cm}
( \, \lambda I_m  + iB \, )^{-1} = \lambda (\,  \lambda^2I_m + B^2 \, )^{-1}
-i( \, \lambda^2 B + B^{\dagger}B^3B^{\dagger} \, )^{\dagger}. 
\hfill
\cr}
$$ }
{\bf Proof.}\, Follows directly from Theorem 11.4. \qquad $\Box $   

\medskip
 
As a special case of Theorem 11.5, we have the following interesting result: Suppose $ M = A + iB $ is a nipotent 
matrix, i.e., $ M^2  = 0$. Then its Moore-Penrose inverse can be expressed as 
$$
\displaylines{
\hspace*{2cm}
 ( \,A + iB \, )^{\dagger} = (E_BAF_B )^{\dagger} -i(E_ABF_A )^{\dagger}. \hfill 
\cr}
$$    
we leave it as an exercise to the reader.

{\rm (b)}\,  If both $ R( A ) \cap R( B ) = \{ 0 \} $ and  $ R( A^*) \cap R( B^*)
= \{ 0 \}, $ then
$$
\displaylines{
\hspace*{2cm}
 ( \,A + iB \, )^{-1} = (E_BAF_B )^{\dagger} -i(E_ABF_A )^{\dagger}. \hfill 
\cr}
$$

We next turn our attention to the Moore-Penrose inverse of sum of $ k $
matrices, and give some general formulas.

\medskip

\noindent {\bf Theorem 11.7.}\, {\em  Let $ A_1, \, A_2, \, \cdots, \, A_k \in
{\cal C}^{ m \times n } $ be given. If they satisfy the following rank additivity
condition
$$
\displaylines{
\hspace*{2cm}
 r( A ) = kr[ \, A_1 , \,  \cdots, \, A_k \, ] = kr[ \, A_1^* , \, \cdots, \,
 A_k^* \, ], \hfill (11.23)
\cr}
$$
where $ A$ is the circulant block matrix defined in {\rm (11.6)}, then 

{\rm (a)}\,  The Moore-Penrose inverse of the sum $ \sum_{i=1}^k A_i $ can
be expressed as
$$ \displaylines{
\hspace*{2cm}
( \, A_1 + A_2 + \cdots + A_k \,)^{\dagger} = J^{\dagger}( A_1) +
J^{\dagger}( A_2)+ \cdots  +  J^{\dagger}( A_k), \hfill (11.24)
\cr}
$$
where $J( A_i)$ is the rank complement of $ A_i( 1 \leq i \leq k ) $ in
$ A .$

{\rm (b)}\, The rank of $ J( A_i) $ is 
$$\displaylines{
\hspace*{2cm} 
r[J( A_i)] = r[ \, A_1 , \, \cdots, \, A_k \,] + r[ \, A_1^*, \ \cdots, \
A_k^* \, ] - r( A ) + r( D_i ),
\hfill (11.25)
\cr}
$$
where $  1 \leq i \leq k, \,  D_i $ is the $ (k-1) \times (k-1) $ block
matrix resulting from the deletion of the first block row and $i$th block column of $ A.$

{\rm (c)} \  $ A_1, \, A_2, \, \cdots, \, A_k $ and $J^{\dagger}( A_1), \,
J^{\dagger}( A_2), \,  \cdots, \,  J^{\dagger}( A_k) $ satisfy the following
two equalities
$$ 
\displaylines{
\hspace*{2cm} 
( \, A_1 + \cdots + A_k  \, ) ( \, A_1 + \cdots + A_k  \, )^{\dagger} =
A_1J^{\dagger}( A_1) + \cdots + A_kJ^{\dagger}( A_k),  \hfill
\cr
\hspace*{2cm}
( \, A_1+ \cdots + A_k \, )^{\dagger} ( \,A_1 + \cdots + A_k \, ) 
= J^{\dagger}( A_1)A_1 + \cdots + J^{\dagger}(A_k) A_k.  \hfill 
\cr}
$$ } 
{\bf Proof.}\, Follows from combining Theorem 9.18 with the equality (11.13). \qquad$ \Box$ 

\medskip
 
\noindent {\bf  Corollary 11.8.}\, {\em  Let $ A_1, \, A_2, \, \cdots, \, A_k \in
{\cal C}^{ m \times n } $. If they satisfy the following rank additivity
condition
$$
 r( \, A_1+ A_2 + \cdots + A_k \, ) = r(A_1 )+ r(A_2) + \cdots +r( A_k),
 \eqno (11.26)
$$ 
then the  Moore-Penrose inverse of the sum $ \sum_{i=1}^k A_i $ can
be expressed as
$$ 
( \, A_1 + A_2 + \cdots + A_k \, )^{\dagger} =
(E_{\alpha_1}A_1F_{\beta_1})^{\dagger}  +
( E_{\alpha_2}A_2F_{\beta_2})^{\dagger} + \cdots +
( E_{\alpha_k}A_kF_{\beta_k})^{\dagger},  \eqno (11.27)
$$
where $ \alpha_i $ and $ \beta_i $ are 
$$
\alpha_i = [\, A_1, \,  \cdots \, A_{i-1}, \, A_{i+1}, \, \cdots, \ A_k  \,], 
 \ \ \ \beta_i = \left[ \begin{array}{c}  A_1 \\ \vdots \\ A_{i-1} 
\\ A_{i+1} \\ \vdots \\ A_k \end{array} \right], \ \ \  i = 1, \, 2, \, \cdots, \, k.
$$  }
{\bf Proof.}\, We first show that under the condition (11.17) the
rank of the circulant matrix $ A $ in (11.6) is 
 $$ 
r( A )= k[ \, r(A_1 )+ r(A_2) + \cdots +r( A_k) \, ].  \eqno (11.28) 
$$      
According to (11.7), we see that
$$ 
 r( A )= r(J_1 )+ r(J_2) + \cdots +r( J_k).
$$
Under Eq. (11.26), the ranks of all $ J_i$ are the same, that is,
$$
r( J_i )= r(A_1 )+ r(A_2) + \cdots +r( A_k), \  \ \ \  i = 1, \, 2, \, \cdots, \, k. 
$$
Thus we have (11.28). In that case, applying the result in Corollary 9.19 to 
the circulant block matrix $ A $
 in  (11.13) produces the equality (11.27). \qquad $ \Box$   

\medskip
 
It is worth to point out that the formulas for
Moore-Penrose inverses of sums of matrices given in this chapter and those for 
 Moore-Penrose inverses of block matrices given in Chapter 9 are, in
fact, a group of dual results. That is to say, not only can we derive 
Moore-Penrose inverses of sums of matrices from Moore-Penrose inverses of
block matrices, but also we can make a contrary derivation. For simplicity,
 here we  only illustrate this assertion by a $ 2 \times 2$  block matrix.
 In fact,  for any  $ 2 \times 2 $  block matrix can factor as 
$$
 M = \left[ \begin{array}{cc}  A & 0 \cr 0 & D \end{array} \right] =
  \left[ \begin{array}{cc}  A & 0 \cr 0 & D \end{array} \right]  +
  \left[ \begin{array}{cc}  0 & B \cr C & 0  \end{array} \right]
  = N_1 + N_2.
$$
If $ M$ satisfies the rank additivity condition (9.22), then $ N_1 $ and
$ N_2 $ satisfy
$$
r \left[ \begin{array}{cc}  N_1 & N_2 \cr N_2 & N_1 \end{array} \right] =
r \left[ \begin{array}{cccc}  A & 0 & 0 & B \cr 0 & D & C & 0 \cr 0  & B & A
& 0 \cr C & 0 & 0 & D \end{array} \right] =
2r  \left[ \begin{array}{cc}  A & B \cr C & D \end{array} \right]
= 2r \left[ \begin{array}{cc}  N_1 \cr N_2 \end{array} \right] = 2r[\, N_1, \ N_2 \, ]. 
$$
Hence by Theorem 11.4, we have 
$$
 M_1^{\dagger } = ( \, N_1 + N_2  \,)^{\dagger} =
 J^{\dagger}( N_1 ) + J^{\dagger}( N_2 ), \eqno (11.29) 
$$ 
where $ J( N_1 ) $ and $ J(N_2 ) $ are, respectively, the rank complements
of $ N_1 $ and $ N_2 $  in  $ \left[ \begin{array}{cc} N_1 & N_2 \cr N_2 & N_1 
\end{array} \right]  $.
Written in an explicit form, (11.29) is exactly the formula (9.33). 

\medskip

Besides (11.13), some other identities  between  Moore-Penrose inverses of sums of matrices and 
Moore-Penrose inverses of block matrices can also be established. Here we present a result for the sum of four 
matrices.   
$$ 
 ( \,  A_0 +  A_1 + A_2 +  A_3 \, )^{\dagger }  = \frac{1}{4}[ \, I_n , \,  I_n , \, I_n, \,  
I_n \, ] \left[ \begin{array}{rrrr}  A_0 & A_1 & A_2 &  A_3   \\  A_1 &  A_0 &  A_3 &  A_2 
\\  A_2 &  A_3 &  A_0  & A_1 \\ A_3 &  A_2 &  A_1  & A_0 \end{array} \right]^{\dagger} 
\left[ \begin{array}{c}  I_m \\ I_m \\ I_m  \\ I_m  \end{array} \right], \eqno(11.30) 
$$ 
Clearly, the block matrix in (9.30) is not  $ 4 \times 4$  block circulant, but
 $ 2 \times 2$  block circulant with two  $ 2 \times 2$  block circulants in it.

In addition, we mention another interesting fact that (11.16) can be extended to 
any real quaternion matrix of the form $ A = A_0 +  iA_1 + jA_2 +  kA_3,$  
where $A_0$---$A_3$ are real $ m \times n$ matrices and 
$i^2 = j^2 = k^2  = -1,  \, ijk = -1$, as follows:
$$ 
 ( \,  A_0 +  iA_1 + iA_2 +  kA_3 \, )^{\dagger }  = \frac{1}{2}[ \, I_n , \ jI_n \, ]
  \left[ \begin{array}{cc}  A_0 +  iA_1 &  -(\, A_2 +  iA_3 \, )  \\  A_2 - iA_3&  A_0 - iA_1 
\end{array} \right]^{\dagger} \left[ \begin{array}{c}  I_m \\ -jI_m  \end{array} \right], \eqno (11.31) 
$$ 
and 
$$ 
 ( \,  A_0 +  iA_1 + iA_2 +  kA_3 \, )^{\dagger }  = \frac{1}{4}[ \, I_n , \  iI_n , \ jI_n, \  
kI_n \, ] \left[ \begin{array}{rrrr}  A_0 & -A_1 & -A_2 &  - A_3   \\  A_1 &  A_0 &  A_3 &  - A_2 
\\  A_2 &  -A_3 &  A_0  & A_1 \\ A_3 &  A_2 &  -A_1  & A_0 \end{array} \right]^{\dagger} 
\left[ \begin{array}{c}  I_m \\ -iI_m \\ -jI_m  \\ -kI_m  \end{array} \right]. \eqno (11.32) 
$$ 

Moreover denote $ ( \,  A_0 +  iA_1 + iA_2 +  kA_3 \, )^{\dagger} 
= G_0 +  iG_1 + iG_2 +  kG_3.$ Then 
$$
\left[ \begin{array}{rrrr}  A_0 & -A_1 & -A_2 &  - A_3   \\  A_1 &  A_0 &  A_3 &  - A_2 
\\  A_2 &  -A_3 &  A_0  & A_1 \\ A_3 &  A_2 &  -A_1  & A_0 \end{array} \right]^{\dagger} = 
\left[ \begin{array}{rrrr}  G_0 & -G_1 & -G_2 &  - G_3   \\  G_1 &  G_0 &  G_3 &  - G_2 
\\  G_2 &  -G_3 &  G_0  & G_1 \\ G_3 &  G_2 &  -G_1  & G_0 \end{array} \right]. \eqno (11.33) 
$$ 
These equalities are in fact derived from the following two universal factorization equalities (see \cite {Ti6}) 
$$
 P_{2m} \left[ \begin{array}{rr}  A  & 0  \\ 0 &  A \end{array} 
\right] P_{2n}^*= \left[ \begin{array}{cc}  A_0 +A_1i &  -(A_2 +A_3i) \\ A_2 - A_3i &  A_0 - A_1i \end{array} \right],  
\eqno (11.34)
 $$  
where $P_{2m} $ and $ P_{2n}^* $  are the following two unitary quaternion  matrices  
$$
 P_{2m} = \frac{1}{ \sqrt{2}} \left[ \begin{array}{cc}  \ I_m  &  - iI_m \\ -jI_m &  kI_m  \end{array} \right], 
\qquad
 P_{2n}^* = \frac{1}{ \sqrt{2}} \left[ \begin{array}{cc}  \ I_n  & jI_n \\ iI_n &  -kI_n  \end{array} \right], 
$$
and
$$ 
Q_{4m}  \left[ \begin{array}{cccc} A &  &  &  \\  & A &  & \\
 &  &  A &    \\  & & & A  \end{array} \right] Q_{4n}^*= \left[ \begin{array}{rrrr} A_0 & -A_1 &  -A_2 & - A_3  \\ A_1 & A_0 & - A_3 & A_2 \\
A_2 & A_3 &  A_0 &  -A_1  \\ A_3 & -A_2 &  A_1 & A_0  \end{array} \right], \eqno (11.35)
$$ 
where $ Q_{4t}$ is the following  unitary quaternion matrix    
$$
Q_{4t} =  Q_{4t}^* =  \frac{1}{2} \left[ \begin{array}{cccc} I_t & iI_t & jI_t & kI_t  \\ -iI_t  & I_t  
& kI_t  & -jI_t  \\ -jI_t  & -kI_t  & I_t  &  iI_t   \\  -kI_t  & jI_t  & -iI_t  & I_t
   \end{array} \right], \ \ \ \  t = m, \ n.  
$$ 
Based on (11.31)---(11.33), one find easily determine expressions of Moore-Penrose inverses of 
 any real quaternion matrices, especially the inverses of nonsingular matrices.

Furthermore, it should be pointed out that the above work can extend to matrices over any $2^n$-dimensional 
real and complex Clifford algebras through a set of universal similarity factorization equalities 
established in the author's recent papers \cite {Ti6} and  \cite {Ti9}.

\markboth{YONGGE  TIAN }
{12. RANK EQUALITIES FOR SUBMATRICES IN  MOORE-PENROSE INVERSES}

\chapter{Rank equalities for submatrices in Moore-Penrose inverses}

\noindent Let 
 $$ 
 M = \left[ \begin{array}{cc} A  & B  \\ C & D  \end{array} \right]
  \eqno (12.1)
 $$ 
 be a $ 2 \times2 $ block matrix over ${\cal C}$, where  $ A \in
 {\cal C}^{m \times n}, \ B \in {\cal C}^{m \times k}, \
 C \in {\cal C}^{l \times n}, \ D \in {\cal C}^{l \times k}$, and let
 $$ 
 V_1 = \left[ \begin{array}{c} A \\ C  \end{array} \right], \qquad
 V_2 = \left[ \begin{array}{c} B \\ D  \end{array} \right], \qquad
  W_1 = [ \, A, \ B \,], \qquad W_2 = [ \, C, \ D \, ].\eqno (12.2)
 $$ 
 Moreover, partition the Moore-Penrose inverse of $ M $ as 
 $$ 
 M^{\dagger} = \left[ \begin{array}{cc} G_1  & G_2  \\ G_3 & G_4
 \end{array} \right], \eqno (12.3)
 $$ 
 where $G_1 \in {\cal C}^{n \times m}$. As is well known, the expressions
  of the submatrices  $ G_1$---$G_4$ in (12.3) are quite complicated if there are no 
restrictions  on the blocks in $ M $ (see, e.g., Hung and Markham \cite{HM1},  
 Miao \cite{M2}). In that case, it is hard to find properties of submatrices in 
$ M^{\dagger}$.  In the present chapter,  we consider a simpler problem---what is the ranks of submatrices in $ M^{\dagger}$, when $ M$ is
 arbitrarily given? This  problem was examined by Robinson \cite{Rob} and 
Tian \cite{Ti3}. In this chapter, we shall give this problem a new discussion.

\medskip

\noindent  {\bf Theorem 12.1.}\, {\em Let $ M$  and $M^{\dagger}$ be
given by  {\rm (12.1)} and {\rm (12.3)}. Then
 $$ \displaylines{
\hspace*{2cm}
 r(G_1) =  r \left[ \begin{array}{cc} V_2D^*W_2  & V_1  \\ W_1 & 0
 \end{array} \right] -r( M ), \ \ \
 r(G_2) =  r \left[ \begin{array}{cc} V_2B^*W_1  & V_1  \\ W_2 & 0
 \end{array} \right] -r( M ),
 \hfill (12.4) 
 \cr
\hspace*{2cm}
 r(G_3) =  r \left[ \begin{array}{cc} V_1C^*W_2  & V_2  \\ W_1 & 0
 \end{array} \right] - r( M ), \ \ \
 r(G_4) =  r \left[ \begin{array}{cc} V_1A^*W_1  & V_2  \\ W_2 & 0
 \end{array} \right] - r( M ),
  \hfill (12.5)
\cr}
 $$ 
 where $ V_1, \ V_2, \ W_1$ and $ W_2 $ are defined in {\rm (12.2)}. }

\medskip

\noindent {\bf Proof.}\, We only show the first equality in  (12.4). In fact
 $ G_1$ in (12.3) can be written as 
 $$ 
 G_1 = [\, I_n, \ 0 \,] M^{\dagger} \left[ \begin{array}{c} I_m  \\ 0
  \end{array} \right] = PM^{\dagger}Q.
  \eqno (12.6)  
 $$ 
 Then applying  (2.1) to it we find
 \begin{eqnarray*}
 r(G_1)  &= & r \left[ \begin{array}{cc} M^*MM^*  & M^*Q   \\ PM^* & 0
 \end{array} \right] -r( M ) \\
 & = &  r \left[ \begin{array}{cc} MM^*M  & MP^*   \\ Q^*M & 0 \end{array}
 \right] -r( M ) \\
 & = & r \left[ \begin{array}{cc} [\, V_1, \ V_2 \,]M^*
 \left[ \begin{array}{c} W_1  \\ W_2
 \end{array} \right]  & V_1  \\ W_1 & 0 \end{array} \right] -r( M ) \\
  & = & r \left[ \begin{array}{cc} [\, 0, \ V_2 \,]M^*
   \left[ \begin{array}{c} 0 \\ W_2
 \end{array} \right]  & V_1  \\ W_1 & 0  \end{array} \right] -r( M ) 
=  r\left[ \begin{array}{cc} V_2D^*W_2  & V_1  \\ W_1 & 0  \end{array}
 \right] -r( M ),
 \end{eqnarray*}
 establishing the first equality in (12.4). \qquad $\Box$ 

\medskip

\noindent  {\bf Corollary 12.2.}\, {\em Let $ M$  and $M^{\dagger}$ be
given by {\rm (12.1)} and {\rm (12.3)}. If 
$$
\displaylines{
\hspace*{2cm}
 r( M ) = r(V_1) + r(V_2),  \ \ i.e., \ \ R( V_1) \cap  R(V_2) = \{ 0 \},
   \hfill (12.7)
\cr
then \hfill
\cr
\hspace*{2cm}
 r(G_1) =  r \left[ \begin{array}{cc} A  & B  \\ D^*C & D^*D  \end{array}
 \right] -
 r\left[ \begin{array}{c} B  \\ D \end{array} \right], \ \ \ 
 r(G_2) =  r \left[ \begin{array}{cc} B^*A  & B^*B  \\ C & D  \end{array}
 \right] -
 r\left[ \begin{array}{c} B  \\ D \end{array} \right], \hfill (12.8)
\cr
\hspace*{2cm}
 r(G_3) =  r \left[ \begin{array}{cc} A  & B  \\ C^*C & C^*D  \end{array}
 \right] -
 r\left[ \begin{array}{c} A  \\ C \end{array} \right], \ \ \ 
 r(G_4) =  r \left[ \begin{array}{cc} A^*A  & A^*B  \\ C & D  \end{array}
 \right] -
 r\left[ \begin{array}{c} A  \\ C \end{array} \right]. \hfill (12.9)
\cr}
 $$}
{\bf Proof.}\, Under  (12.7),  we also know that  $ R( V_1)
  \cap  R( V_2D^*W_2) = \{ 0 \}$. Thus the first
 equality in (12.4) becomes 
 \begin{eqnarray*}
 r(G_1) & =  & r \left[ \begin{array}{cc} V_2D^*W_2  & V_1  \\ W_1 & 0
 \end{array} \right] -r( M ) \\
  & =  & r \left[ \begin{array}{c} V_2D^*W_2 \\ W_1  \end{array} \right]
  + r(V_1) -r( M ) =  r \left[ \begin{array}{c} D^*W_2 \\ W_1  \end{array} \right]
 -r( V_2)   =  r \left[ \begin{array}{c} W_1 \\ D^*W_2  \end{array} \right] -
 r( V_2 ),
 \end{eqnarray*}
establishing  the first one in  (12.8). Similarly, we can show the
 other three in  (12.8) and (12.9). \qquad $\Box$ 

\medskip

 Similarly, we have the following. 

\medskip

 \noindent {\bf Corollary 12.3.}\, {\em Let $ M$  and $M^{\dagger}$ be given
  by  {\rm (12.1)} and {\rm (12.3)}. If
 $$ \displaylines{
\hspace*{2cm}
 r( M ) = r(W_1) + r(W_2),  \ \ i.e., \ \ R( W_1^*) \cap  R( W^*_2)
 = \{ 0 \},   \hfill (12.10)
 \cr
then \hfill
\cr
\hspace*{2cm}
 r(G_1) =  r \left[ \begin{array}{cc} A  & BD^*  \\ C & DD^*  \end{array}
 \right] - r[\, C, \ D \, ], \ \ \
 r(G_2) =  r \left[ \begin{array}{cc} A  & BB^*  \\ C & DB^*  \end{array}
  \right] -  r[\, A, \ B \, ],
 \hfill (12.11)
 \cr
\hspace*{2cm}
 r(G_3) =  r \left[ \begin{array}{cc} AC^*  & B  \\ CC^* & D  \end{array}
 \right] - r[\, C, \ D \, ],\ \ \
 r(G_4) =  r \left[ \begin{array}{cc} AA^*  & B  \\ CA^* & D  \end{array}
 \right] - r[\, A, \ B \, ].
 \hfill (12.12)
 \cr}
$$}
\hspace*{0.3cm} Combining the above two corollaries, we obtain the following, 
 which is previously shown in Corollary 9.9.

\medskip

\noindent  {\bf Corollary 12.4.}\, {\em Let $ M$  and $M^{\dagger}$ be given
 by  {\rm (12.1)} and {\rm (12.3)}. If $ M $
 satisfies the rank additivity condition
 $$ \displaylines{
\hspace*{3cm}
 r( M ) = r(V_1) + r(V_2) = r(W_1) + r(W_2),  \hfill (12.13)
 \cr
 then \hfill
 \cr
\hspace*{3cm}
 r(G_1) = r( D ) + r(V_1) + r( W_1) - r( M ),  \hfill (12.14)
 \cr
\hspace*{3cm}
 r(G_2) = r( B ) + r(V_1) + r( W_2) - r( M ),  \hfill(12.15)
 \cr
\hspace*{3cm}
 r(G_3) = r( C ) + r(V_2) + r( W_1) - r( M ),  \hfill (12.16)
 \cr
\hspace*{3cm}
 r(G_4) = r( A ) + r(V_2) + r( W_2) - r( M ).  \hfill (12.17)
\cr}
 $$ 
}
{\bf Proof.}\, We only show  (12.14). Under  (12.13), we find that 
 $$ 
 r \left[ \begin{array}{cc} V_2D^*W_2  & V_1  \\ W_1 & 0  \end{array}
 \right] = r(V_2D^*W_2 ) + r(V_1) +
 r(W_1), 
 $$ 
 where 
 $$
  r(V_2D^*W_2 ) = r \left[ \begin{array}{cc} BD^*C  & BD^*D  \\
   DD^*C & DD^*D \end{array} \right] = r(D).
 $$ 
 Thus the first equality in  (12.4) reduces to (12.14).
 \qquad $\Box$ 

\medskip

\noindent  {\bf Corollary 12.5.}\, {\em Let $ M$  and $M^{\dagger}$
be given by {\rm (12.1)} and {\rm (12.3)}. If $ M $
 satisfies the rank additivity condition
 $$
 r( M ) = r(A) + r(B) +r(C) + r(D),  \eqno (12.18)
 $$
 then
 $$ 
 r(G_1) = r(A), \ \ \ \ r(G_2) = r(C), \ \ \ \ r(G_3) = r(B), \ \ \ \
  r(G_4) = r(D). \eqno (12.19)
 $$ }
{\bf Proof.}\, Follows directly from (12.14)---(12.17).
 \qquad $\Box$ 

\medskip

\noindent {\bf Corollary 12.6.}\, {\em Let $ M$  and $M^{\dagger}$ be given by
   {\rm (12.1)} and {\rm (12.3)}. If
 $$ 
 r( M ) = r(V_1),   \ \ \ i.e.,  \ \ \ R( V_2)  \subseteq R(V_1),
  \eqno (12.20)
 $$ 
 then
 $$ 
 r(G_1) = r( A ),  \qquad r(G_2) = r( C ), \eqno (12.21) 
 $$ 
 $$ 
 r(G_3) =  r \left[ \begin{array}{cc} V_1C^*C  & V_2  \\ A & 0
 \end{array} \right] - r(V_1), \ \ \ \  r(G_4) =
 r \left[ \begin{array}{cc} V_1A^*A  & V_2  \\ C & 0
  \end{array} \right] - r( V_1).\eqno (12.22)
 $$  }
{\bf Proof.}\, The inclusion in  (12.20)
implies that
 $$
  R(V_2D^*W_2) \subseteq  R( V_1), \ \ \ R( V_2D^*W_1) \subseteq
   R( V_1),  \ \ \
  R(B) \subseteq  R(A),  \ \ \  R(D) \subseteq  R(C).
 $$ 
 Thus the two rank equalities in  (12.4) become 
 $$ 
 r(G_1) =  r \left[ \begin{array}{cc} V_2D^*W_2  & V_1  \\ W_1 & 0
  \end{array} \right] -r( M ) = r(V_1) +
 r( W_1) - r( M ) = r(W_1) = r( A ), 
 $$ 
 $$ 
 r(G_2) =  r \left[ \begin{array}{cc} V_2B^*W_1  & V_1  \\ W_2 & 0
  \end{array} \right] -r( M ) = r(V_1) +
 r( W_2) - r( M ) = r(W_2) = r(C),   
 $$ 
 and the two rank equalities in  (12.5) become 
 \begin{eqnarray*}
 r(G_3) & =  & r \left[ \begin{array}{cc} V_1C^*W_2  & V_2  \\ W_1 & 0
 \end{array} \right] -r( M ) \\
  & =  & r \left[ \begin{array}{ccc} V_1C^*C  & V_1C^*D & V_2 \\
  A & B & 0 \end{array} \right] -r( M )
 = r \left[ \begin{array}{cc} V_1C^*C  & V_2 \\ A & 0 \end{array}
  \right] -r( M ),
 \end{eqnarray*}
 \begin{eqnarray*}
 r(G_4) & =  & r \left[ \begin{array}{cc} V_1A^*W_1  & V_2  \\ W_2 & 0
  \end{array} \right] -r( M ) \\
  & =  & r \left[ \begin{array}{ccc} V_1A^*A  & V_1A^*B & V_2 \\
  C & D & 0 \end{array} \right] -r( M )
 = r \left[ \begin{array}{cc} V_1A^*A  & V_2 \\ C & 0 \end{array}
 \right] -r( M ).
 \end{eqnarray*}
 Hence we have  (12.21) and (12.22). \qquad $\Box$ 

\medskip

 Similarly, we have the following. 

\medskip

\noindent  {\bf Corollary 12.7.}\, {\em Let $ M$ and $M^{\dagger}$
be given by {\rm (12.1)} and {\rm (12.3)}. If
 $$ \displaylines{
\hspace*{2cm}
 r( M ) = r(W_1),  \ \ \ i.e., \ \ \ R( W_2^*)  \subseteq R(W_1^*),
  \hfill (12.23) 
\cr
 then \hfill
\cr
\hspace*{2cm}
r(G_1) = r( A ),  \qquad r(G_3) = r(B), \hfill (12.24) 
\cr
\hspace*{2cm}
r(G_2) =  r \left[ \begin{array}{cc} BB^*W_1 & A \\ W_2 & 0  \end{array}
  \right] - r(W_1), \ \ \ \
 r(G_4) =  r \left[ \begin{array}{cc} AA^*W_1  & B  \\ W_2 & 0
 \end{array} \right] - r( W_1).
 \hfill (12.25)  
 \cr}
$$ }
\hspace*{0.3cm}Combining the above two corollaries, we obtain the following.

\medskip

\noindent  {\bf Corollary 12.8.}\, {\em Let $ M$ and $M^{\dagger}$ be
given by {\rm (12.1)} and {\rm (12.3)}. If 
 $$ 
 r( M ) = r(A),   \eqno (12.26)
 $$ 
 then
 $$ 
 r(G_1) = r( A ),  \qquad r(G_2) = r(C),   \qquad r(G_3) = r(B),
  \eqno (12.27)
 $$ 
and 
 $$ 
 r(G_4) =  r \left[ \begin{array}{cc} AA^*A & B  \\ C & 0  \end{array}
 \right] - r(A). \eqno (12.28)
 $$  }
{\bf Proof.}\, Clearly   (12.26) implies that $ r(M)
= r( V_1) = r( W_1)$. Thus we have  (12.27) by
 Corollaries 12.6 and 12.7. On the other hand, (12.26) is also
 equivalent to $ AA^{\dagger}B = B, \
   CA^{\dagger}A= C$ and $ D = CA^{\dagger}B$ by (1.5). Hence 
 \begin{eqnarray*}
 r(G_4) &= & r \left[ \begin{array}{cc} V_1A^*W_1  & V_2 \\ W_2 & 0
  \end{array} \right] - r(M) \\
  &= & r \left[ \begin{array}{cc} AA^*W_1  & B  \\ CA^*W_1 & D \\
 W_2 & 0 \end{array} \right] - r(A) \\
&= & r \left[ \begin{array}{cc} AA^*W_1  & B  \\ W_2 & 0
  \end{array} \right] - r( W_1) \\
  &= & r \left[ \begin{array}{ccc} AA^*A  & AA^*B & B  \\ C & D & 0
   \end{array} \right] - r(A)
 = r \left[ \begin{array}{cc} AA^*A  & B  \\ C & 0  \end{array} \right]
 - r(A),
 \end{eqnarray*}
 which is (12.28). \qquad $\Box$ 

\medskip

 Next we list a group of rank inequalities derived from (12.4) and (12.5). 

\medskip

\noindent {\bf Corollary 12.9.}\, {\em Let $ M$ and $M^{\dagger}$ be given by  
 {\rm (12.1)} and {\rm (12.2)}. Then 
the rank of $ G_1$ in $ M^{\dagger}$ satisfies the rank inequalities
$$\displaylines{
\hspace*{2cm}
r(G_1) \leq r( D ) + r[ \, A , \ B \, ] + r \left[ \begin{array}{c} A
\\ C  \end{array} \right] - r( M ), \hfill (12.29) 
 \cr
\hspace*{2cm}
r(G_1) \geq  r[ \, A , \ B \, ] + r \left[ \begin{array}{c} A  \\ C
\end{array} \right] - r( M ),
\hfill (12.30) 
 \cr
\hspace*{2cm}
r(G_1) \geq r( D ) - r[ \, C , \ D \, ] - r \left[ \begin{array}{c} B
\\ D  \end{array} \right] + r( M ).
\hfill (12.31) 
\cr}
$$ } 
{\bf Proof.}\, Observe that 
$$  
r(V_1) + r(W_1) \leq r \left[ \begin{array}{cc} V_2D^*W_2  & V_1
 \\ W_1 & 0  \end{array} \right]
 \leq  r(D) +  r(V_1) + r(W_1). 
$$ 
Putting them in the first rank equality in   (12.4), we obtain
(12.29) and (12.30).
To show  (12.31), we need the following rank equality
$$ 
r(CA^{\dagger}B) \geq  r \left[ \begin{array}{cc} A & B  \\ C  & 0
\end{array} \right]
- r \left[ \begin{array}{c} A  \\ C  \end{array} \right] -
 r[ \, A , \ B \, ] + r(A), \eqno (12.32)
$$ 
which is derived from (1.6). Now applying (12.32) to
 $ PM^{\dagger}Q$ in (12.6), we obtain
\begin{eqnarray*}
r(G_1) = r( PM^{\dagger}Q) & \geq & r \left[ \begin{array}{cc} M
 & Q  \\ P & 0  \end{array} \right]
- r \left[ \begin{array}{c} M  \\ P  \end{array} \right] -
r[ \, M , \ Q \,] + r(M) \\
 &= & r( D ) - r[\, C , \ D \,] - r \left[ \begin{array}{c} B
 \\ D  \end{array} \right] + r( M ),
\end{eqnarray*}
which is (12.31). \qquad  $\Box$ 

\medskip
 
Rank inequalities for  the block entries $ G_2, \ G_3 $ and $ G_4$ in (12.3)
 can also be derived in the similar way shown above. Finally let
 $ D = 0 $ in (12.1). Then the results in (12.4) and (12.5) can be simplified to
  the following. 

\medskip

\noindent {\bf Thoerem 12.10.}\, {\em Let 
$$ 
M_1 = \left[ \begin{array}{cc} A  & B  \\ C & 0 \end{array} \right],
 \eqno (12.33)
$$ 
and denote the Moore-Penrose inverse of $ M $ as 
$$ 
M_1^{\dagger} = \left[ \begin{array}{cc} G_1  & G_2  \\ G_3 & G_4
 \end{array} \right], \eqno (12.34)
$$ 
where $G_1 \in {\cal C}^{n \times m}$. Then
$$ 
r(G_1) = r \left[ \begin{array}{c} A  \\ C  \end{array} \right]
 + r[ \, A , \ B \, ] - r( M_1 ),  \ \ \ \
r(G_2) = r(C), \ \ \ \ r(G_3) = r(B), \eqno (12.35)
$$ 
$$ 
r(G_4) = r \left[ \begin{array}{ccc} AA^*A  & AA^*B  & B \\  CA^*A
 & CA^*B  & 0  \\ C & 0 & 0
\end{array} \right] - r( M_1 ).    \eqno (12.36)
$$ } 
\hspace*{0.4cm} Various consequences of  (12.35) and (12.36) can also be derived.
But we omit them here.

\markboth{YONGGE  TIAN }
{13. RANK EQUALITIES FOR DRAZIN INVERSES}

\chapter{Rank equalities for Drazin inverses }

\noindent As one of the important types of generalized  inverses of matrices, the Drazin inverses 
and their applications have well been examined in the literature. Having established so many rank 
equalities in the preceding chapters, one  might naturally consider how to extend that work from 
Moore-Penrose inverses to Drazin inverses. To do this, we only need to use a basic identity on 
the Drazin inverse of a matrix $ A^D = A^k( A^{2k+1})^{\dagger}A^k$ (see, e.g., Campbell and Meyer \cite{CM2}). 
In that case,  the rank formulas obtained in the preceding chapters can all apply to 
establish various rank equalities for matrix expressions involving Drazin inverses of matrices. 
  
\medskip

 \noindent {\bf Theorem 13.1.}\, {\em Let $ A \in {\cal C}^{m \times m}$
with ${\rm Ind}(A) = k.$ Then

{\rm (a)} \ $ r(\,I_m \pm A^D \,) = r( \, I_m \pm A \, ).$

{\rm (b)} \ $ r[\,I_m - (A^D)^2 \,] = r( \, I_m  - A^2 \, ).$ 
}
\medskip

\noindent {\bf Proof.}\, Observe that $R( A^k) = R(A^{2k+1})$ and 
 $R[( A^k)^*] = R[(A^{2k+1})^*]$. Thus applying (1.7) and then (1.16) to 
$ I_m - A^D = I_m - A^k( A^{2k+1})^{\dagger}A^k $
yields
\begin{eqnarray*} 
r(\,I_m - A^D \,) &= & r[\,I_m - A^k( A^{2k+1})^{\dagger}A^k \,]  \\
&= &  r \left[ \begin{array}{cc} A^{2k+1} & A^k \\ A^k & I_m
\end{array} \right] -r(A^{2k+1}) \\
&= &  r \left[ \begin{array}{cc} A^{2k+1}- A^{2k} &  0 \\ 0 &
I_m \end{array} \right] -r(A^k) \\
& =&  r( \, A^{2k+1}- A^{2k} \, ) + m - r(A^k) =  r(A^{2k}) +  r( \, I_m - A \,) - r(A^k)
= r( \, I_m - A \,).
\end{eqnarray*} 
Similarly we can find $ r(\,I_m + A^D \,) = r( \,  I_m + A \, )$. Note  by  (1.12) that 
$$
r[\,I_m  - (A^D)^2 \,]  = r(\,I_m + A^D \,) +  r(\,I_m - A^D \,) - m = 
r(\,I_m + A \,) +  r(\,I_m - A \,) - m = r( \, I_m  - A^2 \, ),
$$
establishing Part (b).     \qquad $ \Box$ 

\medskip
   
\noindent {\bf Theorem 13.2.}\, {\em Let $ A \in {\cal C}^{m \times m}$
with ${\rm Ind}(A) = k.$ Then

{\rm (a)} \ $ r(\,A - AA^D \,) = r(\,A - A^DA \,) = r(\, A -A^2 \,).$

{\rm (b)} \ $ r(\,A - AA^DA \,) = r(A) - r(A^D),$ i.e.$,$  $AA^DA \leq_{rs} A$. 

{\rm (c)} \ $ AA^D = A^DA = A \Leftrightarrow A^2 = A.$  } 

\medskip

The results in Theorem 13.2(d) is well known, see, e.g., Campbell and Meyer \cite{CM2}.  

\medskip

\noindent {\bf Proof.}\, Applying (1.6)  and (1.16) to $A - AA^D$
yields
\begin{eqnarray*} 
r(\, A - AA^D \,) &= & r[\,A  - A^{k+1}( A^{2k+1})^{\dagger}A^k \,]  \\
&= &  r \left[ \begin{array}{cc} A^{2k+1} & A^{k+1} \\ A^k & A \end{array}
\right] -r(A^{2k+1}) \\
&= &  r \left[ \begin{array}{cc} A^{2k+1}- A^{2k} &  0  \\ 0 & A \end{array}
 \right] -r(A^k) \\
&=&  r( \, A^{2k+1}- A^{2k} \, ) + r(A) - r(A^k) \\
&=&  r(A^{2k}) + r( \, I_m  - A \, ) - m  + r(A) - r(A^k). \\
&=&  r( \, A  - A^2 \, ). \\
\end{eqnarray*} 
as required for Part (a).  Notice that $A^D$ is an outer inverse of $ A $. Thus it follow by (5.6) 
that
\begin{eqnarray*} 
r(\, A - AA^DA \,) = r(A) - r(A^D) = r(A) - r(A^k),
\end{eqnarray*} 
as required for Part (b). The results in Parts (c) and (d) follow
from Part (a). \qquad  $ \Box$ 

\medskip

\noindent {\bf Theorem 13.3.}\, {\em Let $ A \in {\cal C}^{m \times m}$
with ${\rm Ind}(A) = k.$

{\rm (a)} \  $  r(\,A - A^D \,)= r( \, A^{k+2} - A^k \, ) + r(A) -
r( A^k ) = r( \, A  - A^{3} \, ).$

{\rm (b)} \ $  r(\,A - A^D \,)=  r(A) - r(A^D), \  i.e., \ A^D \leq_{rs} A  \Leftrightarrow  
A^{k+2} = A^k.$

{\rm (c)} \ $  r(\,A - A^{\#} \,)= r( \, A^3 - A \, ),$ if ${\rm Ind}(A) = 1.$ 

{\rm (d)} \  $  A^{\#} = A  \Leftrightarrow  A^3 = A.$ } 

\medskip

\noindent {\bf Proof.}\,  Applying (1.6) and (1.16) to $A - A^D $  yields
\begin{eqnarray*} 
r(\, A - A^D \,) &= & r[\,A  - A^{k}( A^{2k+1})^{\dagger}A^k \,]  \\
&= &  r \left[ \begin{array}{cc} A^{2k+1} & A^{k} \\ A^k & A \end{array}
\right] -r(A^{2k+1}) \\
&= &  r \left[ \begin{array}{cc} A^{2k+1}- A^{2k-1} &  0  \\ 0 & A
\end{array} \right] -r(A^k) \\
&= & r( \, A^{k+2}- A^{k} \, ) + r(A) - r(A^k) \\
&= & r(A^k) + r( \, I_m  - A^2 \, ) - m  + r(A) - r(A^k) \\
&= & r( \, A  - A^3 \, ),
\end{eqnarray*} 
as required for Part(a). The results in Parts (b), (c) and (d) follow 
immediately from it, where the result in Part (d) is well known.  
. \qquad $ \Box$  

\medskip

Similarly, we can establish the following two. 

\medskip

\noindent {\bf Theorem 13.4.}\, {\em Let $ A \in {\cal C}^{m \times m}$
with ${\rm Ind}(A) = k.$

{\rm (a)}\, If $ k \geq 2,$ then $r(\,A^D - A^2 \,) =  r( \, A^5 - A^2 \, ).$

{\rm (b)}\,  If $ k = 1,$ then $r(\,A^{\#} - A^2 \,) = r( \, A^4 - A \, ).$

{\rm (c)}\, $ A^2 = A^D \Leftrightarrow A^5 = A^2$.

{\rm (d)}\, $ A^2 = A^{\#} \Leftrightarrow A^4 = A$ when $ k = 1$. 
} 

\medskip

The two equivalence relations in Theorem 13.4(c) and (d) were obtained by Grass and Trenkler \cite{GT} when they considered generalized and hypergeneralized projectors.     

\medskip

\noindent {\bf Theorem 13.5.}\, {\em Let $ A \in {\cal C}^{m \times m}$
with ${\rm Ind}(A) = k.$

{\rm (a)}\  $r(\,A^D - A^3 \,) = r( \, A^7 -  A^3 \, ) .$

{\rm (b)}\,  If $ k = 2,$ then $r(\,A^D - A^3 \,) = r( \, A^6 - A^2 \, ).$

{\rm (c)}\,  If $ k = 1,$ then $r(\,A^{\#} - A^3 \,) = r( \, A^5 - A \, ).$

{\rm (d)} \ $ A^3 = A^D \Leftrightarrow A^7 = A^3.$ 

{\rm (e)} \  $ A^3 = A^D \Leftrightarrow A^6 = A^2$ when $ k = 2$. 

{\rm (f)} \ $ A^3 = A^{\#} \Leftrightarrow A^5 = A$ when $ k = 1$.

{\rm (g)}\, In general$,$  $r(\,A^D - A^t \,) = r( \, A^{2t+1} - A^t \, ).$
} 

\medskip

\noindent {\bf Theorem 13.6.}\, {\em Let $ A \in {\cal C}^{m \times m}$
with ${\rm Ind}(A) = k.$

{\rm (a)} \  $  r(\,A^D - AA^DA \,) = r(\,  A^k - A^{k+2}\,).$

{\rm (b)} \ $A^D = AA^DA \Leftrightarrow  A^{k+2} = A^k.$ } 

\medskip

\noindent {\bf Proof.}\, Observe that 
$$
A^{k+1}(\,A^D - AA^DA \,) = A^k - A^{k+2}  \ \ {\rm and} \ \  
({A^D})^{k+1}(\, A^k - A^{k+2} \,) = A^D - AA^DA. 
$$
Thus Part (a) follows. \qquad $ \Box$  

\medskip

A square matrix $ A $ is said to be quasi-idempotent if $ A^{k+1} = A^k$ for some positive integer 
$k$. In a recent paper by Mitra \cite{Mi7}, quasi-idempotent matrices and the related topics are 
well examined. The result given below reveals a new aspect on quasi-idempotent matrix.     

\medskip

\noindent {\bf Theorem 13.7.}\, {\em Let $ A \in {\cal C}^{m \times m}$ with  ${\rm Ind}(A) = k.$ Then

{\rm (a)} \ $ r[\,A^D \pm (A^D)^2 \,] = r( \, A^{k+1} \pm A^k \, ).$ 

{\rm (b)} \ $  (A^D)^2 = A^D, \ i.e., \  A^D \ is \ idempotent \Leftrightarrow  A^{k+1} = A^k,$ i.e.$,$  $ A $ is quasi-idempotent.

{\rm (c)} \ $ r[\,A^{\#} \pm (A^{\#})^2 \,] = r( \, A^2 \pm A \,).$ 

{\rm (d)} \ $A^{\#}$ is idempotent $\Leftrightarrow$ $ A $ is idempotent.
}

\medskip

\noindent {\bf Proof.}\, By (2.3) and $(A^D)^2 = (A^2)^D $ we find that 
\begin{eqnarray*}
\lefteqn{r[\, A^D - (A^D)^2 \,] } \\
&= & r[\, A^{k}( A^{2k+1})^{\dagger}A^k  - A^{2k}( A^{4k+2})^{\dagger}A^{2k}
\,]  \\
&= &  r \left[ \begin{array}{ccc} -A^{2k+1} & 0  & A^{k} \\  0 & A^{4k+2}
 & A^{2k} \\ A^k & A^{2k} & 0  \end{array} \right] -r(A^{2k+1}) -
 r(A^{4k+2}) \\
&= &  r \left[ \begin{array}{ccc} -A^{2k+1} & 0  & A^{k} \\  0 & A^{2k+2}
& A^{k} \\ A^{k} & A^{k} & 0 \end{array} \right] -2r(A^{k})  \\
&= &  r \left[ \begin{array}{ccc} 0  & 0  & A^{k} \\  0 &
A^{2k+2}- A^{2k+1} & 0 \\ A^{k} & 0 & 0 \end{array} \right] -2r(A^{k}) 
=  r( \, A^{2k+2}- A^{2k+1} \, ) =  r( \, A^{k+1}- A^{k} \, ).
\end{eqnarray*} 
Similarly, we can obtain $ r[\,A^D + (A^D)^2 \,] = r( \, A^{k+1} + A^k \,).$ 
The results in Parts (b)---(d) follow immediately  from Part (a). 
\qquad $ \Box$

\medskip

\noindent {\bf Theorem 13.8.}\, {\em Let $ A \in {\cal C}^{m \times m}$
with ${\rm Ind}(A) = k.$ Then

{\rm (a)} \ $ r[\,A^D - (A^D)^3 \,] = r(\,  A^k - A^{k+2}\,).$

{\rm (b)} \ $ r[\,A^{\#} - (A^{\#})^3 \,] = r( \, A - A^3 \,),$ if ${\rm Ind}(A) = 1.$

{\rm (c)}\, $ (A^D)^3 = A^D  \Leftrightarrow A^D = AA^DA \Leftrightarrow A^{k+2} = A^k.$ 

{\rm (d)} \ $A^{\#}$ is tripotent if and only if  $A$  is tripotent.}
   
\medskip
 
\noindent {\bf Proof.}\, Notice that
$$
   A[\, A^D - (A^D)^3 \,]A  = AA^DA - A^D \ \ \ {\rm and}  \ \ \
 A^D( \, AA^DA - A^D \,)A^D = A^D - (A^D)^3.  
$$
Thus we have $ r[\, A^D - (A^D)^3 \,] = r(\, AA^DA - A^D \, )$. In that case  we have 
Part (a) by Theorem 13.6(a). 
yields
\begin{eqnarray*}
r[\, A^D - (A^D)^3 \,] & = & r[ \, A^D + (A^D)^2 \, ] + r[ \, A^D -
(A^D)^2 \,] -r(A^D) \\
& = & r( \, A^{k+1} + A^k \, ) +  r( \, A^{k+1} - A^k \,) - r(A^k) 
= r( \, A^{k+1} + A^k \, ) +  r( \, A^{k+1} - A^k \,) - r(A^k) 
,  
\end{eqnarray*}
as required for  Part (a). The equivalence in Part (c) follows directly from Part (a).   \qquad $ \Box$

\medskip

\noindent {\bf Theorem 13.9.}\, {\em Let $ A \in {\cal C}^{m \times m}$
 with ${\rm Ind}(A) = k.$ Then

{\rm (a)} \ $ r[\,AA^D - (AA^D)^* \,] = 2r[ \, A^{k}, \  (A^k)^* \,] - 2r(A^k).$

{\rm (b)} \ $ r[\,(AA^D)(AA^D)^* - (AA^D)^* (AA^D)\,] = 2r[ \, A^{k}, \  (A^k)^* \,] - 2r(A^k).$

{\rm (c)} \ $ r[\,AA^{\#} - (AA^{\#})^* \,] = r[\,(AA^{\#})(AA^{\#})^* - (AA^{\#})^*(AA^{\#}) \,] =  2r[ \, A, \  A^* \,] - 2r(A),$ if ${\rm Ind}(A) = 1$.

{\rm (d)} \ $AA^D = (AA^D)^*  \Leftrightarrow (AA^D)(AA^D)^* = (AA^D)^*(AA^D) \Leftrightarrow R(A^{k}) = R[(A^k)^*],$ i.e.$,$ $A^k $ is EP. 

{\rm (e)} \ $AA^{\#} = (AA^{\#})^* \Leftrightarrow (AA^{\#})(AA^{\#})^* = 
(AA^{\#})^*(AA^{\#}) \Leftrightarrow R(A^*) = R(A), \ i.e., \ A  \ is \ EP.$}

\medskip

\noindent {\bf Proof.}\,  Note that both $ AA^D$ and  $(AA^D)^*$ are idempotent. It 
follows from (3.1) that
\begin{eqnarray*} 
r[\, AA^D - (AA^D)^* \,] & = & r \left[ \begin{array}{c} AA^D \\ (AA^D)^*
\end{array} \right] + r[\, AA^D, \ (AA^D)^* \,] - r(AA^D) - r[(AA^D)^*] \\
& = & 2r[\, AA^D, \ (AA^D)^* \,] - 2r(A^D) \\
&= & 2r[ \, A^{k}, \ (A^k)^* \, ] - 2r(A^k),
\end{eqnarray*} 
as required for Part (a).  Part (b) follows from  Part (a) and Corollary 3.26(d). 
The results in Parts (c)---(e) follow immediately from Part (a). \qquad $ \Box$ 

\medskip

\noindent {\bf Theorem 13.10.}\, {\em Let $ A \in {\cal C}^{m \times m}$
 with ${\rm Ind}(A) = k.$ Then

{\rm (a)}  \ $ r(\, A^{\dagger} - A^D \,) =  r \left[ \begin{array}{c} A^{k}
\\ A^* \end{array} \right]
+ r[ \, A^{k}, \ A^* \,] - r(A^k) - r(A).$ 

{\rm (b)}  \ $ r(\,A^{\dagger} - A^D \,) = r(A^{\dagger}) - r(A^D), \ i.e., \  A^D \leq_{rs} A^{\dagger} 
 \Leftrightarrow R(A^k) \subseteq R(A^*)  \ and  \ R[(A^k)^*] \subseteq r(A), \ i.e., \  A \ is \ 
power-EP$.

{\rm (c)} \ $ r(\, A^{\dagger} - A^{\#} \,) = 2r[\, A, \ A^* \,] - 2r(A),$ If \ ${\rm Ind\,}(A) = 1.$  

{\rm (d)\cite{BeG}} \  $ A^{\dagger} = A^{\#} \Leftrightarrow R(A^*) = R(A), \ i.e.,  \ A \ is \ EP.$ } 

\medskip

\noindent {\bf Proof.}\, Since both $ A^{\dagger}$ and $ A^D$ are 
outer inverses of $ A $, it follows from (5.1) that
\begin{eqnarray*}
r(\,A^{\dagger} - A^D \,) & = &
r \left[ \begin{array}{c} A^{\dagger} \\ A^D \end{array} \right]
+ r[\, A^{\dagger}, \ A^D \,]  - r(A^{\dagger}) - r(A^{D}) \\
& = & r \left[ \begin{array}{c} A^* \\ A^{k} \\ \end{array} \right] +
r[ \, A^*, \ A^{k} \,] - r(A) - r(A^k),
\end{eqnarray*} 
as required for Part (a). The results in Part (b)---(d) follows 
immediately from Part (a). \qquad $ \Box$ 

\medskip

\noindent {\bf Theorem 13.11.}\, {\em Let $ A \in {\cal C}^{m \times m}$
with ${\rm Ind}(A) = k.$ Then

{\rm (a)} \ $ r(\, AA^{\dagger} - AA^D \,) =  r \left[ \begin{array}{c}
A^{k} \\ A^* \end{array} \right] - r(A^k).$

{\rm (b)} \ $ r(\, A^{\dagger}A - A^DA \,) = r[ \, A^{k}, \ A^* \,] -
r(A^k).$

{\rm (c)} \  $ r[\, A^k(A^k)^{\dagger} - AA^D \,] =
r \left[ \begin{array}{c} A^{k} \\ (A^*)^k \end{array} \right] - r(A^k).$

{\rm (d)} \ $ r[\, (A^k)^{\dagger}A^k - A^DA \,] = r[ \, A^{k},
\ (A^*)^k \,] - r(A^k).$ \\
In particular$,$

{\rm (e)} \ $ r(\, AA^{\dagger} - AA^{\#} \,) =  r(\, A^{\dagger}A - A^{\#}A \,) =
r[ \, A, \ A^* \,] - r(A),$ if \ ${\rm Ind\,}(A) = 1.$  

{\rm (f)} \ $ r(\, AA^{\dagger} - AA^D \,) =
 r(AA^{\dagger}) -r( AA^D ), \ i.e., \ AA^D \leq_{rs} AA^{\dagger}  \Leftrightarrow R[(A^k)^*] 
\subseteq R(A)$.

{\rm (g)} \ $ r(\, A^{\dagger}A - A^DA \,) =
 r(A^{\dagger}A) -r( A^DA), \ i.e., \ A^DA \leq_{rs} A^{\dagger}A  
\Leftrightarrow R(A^k) \subseteq R(A^*)$.

{\rm (h)} \ $  A^k(A^k)^{\dagger} = AA^D \Leftrightarrow \
(A^k)^{\dagger}A^k = A^DA  \Leftrightarrow  A^k $ is EP.

{\rm (i)} \  $  AA^{\dagger} = AA^{\#} \Leftrightarrow  A^{\dagger}A = A^{\#}A 
 \Leftrightarrow A$ is EP.} 

\medskip

\noindent {\bf Proof.}\, Note that both $ AA^{\dagger}$ and $AA^D$ are idempotent.
Then it follows by (3.1) that
\begin{eqnarray*} 
r(\,AA^{\dagger} - AA^D \,) &= &
r \left[ \begin{array}{c} AA^{\dagger} \\ AA^D  \end{array} \right]
+ r[\, AA^{\dagger}, \ AA^D \,] - r(AA^{\dagger}) - r(AA^D) \\
& = & r \left[ \begin{array}{c} A^* \\ A^k \end{array} \right]  +
r[\, A, \ A^k\,]  -r(A) - r(A^k) = r \left[ \begin{array}{c} A^* \\ A^k \end{array} \right] - r(A^k),
\end{eqnarray*}
as required for Part (a). Similarly we can establish Parts (b)---(d).  Part (e) is a special case of
 (b)---(d). Based on them we easily get Parts (f)---(i).  \qquad
 $ \Box$ 

\medskip

\noindent {\bf Theorem 13.12.}\, {\em Let $ A \in {\cal C}^{m \times m}$
with ${\rm Ind}(A) = k.$ Then

{\rm (a)}  \ $ r(\, A^{\dagger}A^D - A^DA^{\dagger} \,) =
r \left[ \begin{array}{c} A^{k} \\ A^* \end{array} \right] +
r[ \, A^{k}, \  A^* \, ]  - 2r(A) = r(\, A^{\dagger}A^k - A^kA^{\dagger} \,).$

{\rm (b)} \ $ r(\, A^{\dagger}AA^DA - AA^DAA^{\dagger} \,) =
r \left[ \begin{array}{c} A^{k} \\ A^* \end{array} \right] +
r[ \, A^{k}, \  A^* \, ]  - 2r(A).$

{\rm (c)} \ $ r(\, A^{\dagger}AA^D - A^DAA^{\dagger} \,) =
r \left[ \begin{array}{c} A^{k} \\ A^* \end{array} \right] +
r[ \, A^{k}, \  A^* \, ]  - 2r(A).$

{\rm (d)} \ $ r(\, A^{\dagger}A^{\#} - A^{\#}A^{\dagger} \,) = 2r[ \, A, \  A^* \, ] 
- 2r(A),$ if ${\rm Ind}(A) = 1.$ \\

{\rm (e)}  \ $ r(\, A^{\dagger}AA^{\#} - A^{\#}AA^{\dagger} \,) = 2r[ \, A, \  A^* \, ] - 2r(A),$ if ${\rm Ind}(A) = 1.$ \\
In particular$,$

{\rm (f)} \ $ A^{\dagger}$ commutes with $A^D$ $\Leftrightarrow$ 
$ A^{\dagger}$ commutes with $AA^DA$  $\Leftrightarrow$  
$ A^{\dagger}$ commutes with $AA^D$ $\Leftrightarrow$
$ A^{\dagger}$ commutes with $A^k$ $\Leftrightarrow$ $ AA^{\dagger}= AA^D  \ 
and \ A^DA = A^{\dagger}A$  $\Leftrightarrow$ 
$ R(A^k) \subseteq R(A^*) $ and $[(A^k)^*] \subseteq R(A),$ i.e.$,$  A is power-EP.  

{\rm (g)} \  $ A^{\dagger}$ commutes with $A^{\#}$ $\Leftrightarrow$ 
$ A^{\dagger}$ commutes with $AA^{\#}$ $\Leftrightarrow$
 A is EP.}

\medskip

\noindent {\bf Proof.}\, Applying (2.2) to $ A^{\dagger}A^D - A^DA^{\dagger}$
yields
\begin{eqnarray*} 
r(\, A^{\dagger}A^D - A^DA^{\dagger} \,) & = & r \left[ \begin{array}{ccc}
A^*AA^* & 0  & A^*A^D \\  0 & -A^*AA^* & A^* \\ A^* & A^DA^* & 0  \end{array} \right] - 2r(A)  \\
&= &  r \left[ \begin{array}{ccc} A^*AA^* & A^*A^DAA^*  & A^*A^D \\
 0 & 0 & A^* \\ A^* & A^DA^* & 0  \end{array} \right] - 2r(A) \\
 &= &  r \left[ \begin{array}{ccc} 0  & 0  & A^*A^D \\
 0 & 0 & A^* \\ A^* & A^DA^* & 0  \end{array} \right] - 2r(A) \\
 &=& r \left[ \begin{array}{c} A^D \\ A^* \end{array} \right] + r[ \, A^D, \
   A^* \, ]  - 2r(A) = r \left[ \begin{array}{c} A^k \\ A^* \end{array} \right] + 
r[ \, A^k, \  A^* \, ] - 2r(A).
\end{eqnarray*} 
The second one in Part (a) is from Theorem 13.10(a) and (b). The two  formulas in Part (b) and (c) are derived from (4.1) by noting that 
$AA^D =A^DA$ is idempotent.  Parts (d)---(g) are direct  consequences of Parts (a)---(c). \qquad$ \Box$ 

\medskip

\noindent {\bf Theorem 13.13.}\, {\em Let $ A \in {\cal C}^{m \times m}$
with ${\rm Ind}(A) = k.$ Then

{\rm (a)} \ $ r[\, (AA^{\dagger})A^D - A^D(AA^{\dagger}) \,] =  r \left[ \begin{array}{c}
A^{k} \\ A^* \end{array} \right] - r(A).$

{\rm (b)} \ $ r[\, (A^{\dagger}A)A^D - A^D(A^{\dagger}A) \,]  = r[ \, A^{k}, \ A^* \,] -
r(A).$  

{\rm (c)} \  $ r(\, A^{\dagger}A^D - A^DA^{\dagger} \,) =  
r[\, (AA^{\dagger})A^D - A^D(AA^{\dagger}) \,] +  r[\, (A^{\dagger}A)A^D - A^D(A^{\dagger}A) \,].$ 

{\rm (d)} \ $ r[\, (AA^{\dagger})A^{\#} - A^{\#}(AA^{\dagger}) \,] =  r[\, (A^{\dagger}A)A^{\#}
 - A^{\#}(A^{\dagger}A) \,] =  r[ \, A, \ A^* \,] - r(A).$

{\rm (e)} \ $ A^D$ commutes with $AA^{\dagger}  \Leftrightarrow R[(A^k)^*] \subseteq  R(A).$   

{\rm (f)} \ $ A^D$ commutes with $A^{\dagger}A \Leftrightarrow R(A^k) \subseteq  R(A^*).$   

{\rm (g)}  \ $ A^{\dagger}A^D = A^DA^{\dagger}$ $ \Leftrightarrow$  $A^D$ commutes  with 
$A^{\dagger}A$  and $A^D$  commutes with  $A^{\dagger}A$  $\Leftrightarrow$  $R(A^k) \subseteq  R(A^*)$ 
and  $R(A^k) \subseteq  R(A^*)$ $\Leftrightarrow$ $A$ is power-EP.    

{\rm (h)} \ $ A^{\dagger}A^{\#} = A^{\#}A^{\dagger}$ $\Leftrightarrow$  $A^{\#}$ commutes with 
$A^{\dagger}A$ $\Leftrightarrow$  $A^{\#}$ commutes with $A^{\dagger}A$  $\Leftrightarrow$  $A$ is EP.}   

\medskip

\noindent {\bf Proof.}\, Note that both $ AA^{\dagger}$ and $ A^{\dagger}A$ are idempotent.
Thus  Parts (a) and (b) can easily be established through (4.1). Contrasting Parts (a) and (b) with 
Theorem 13.12(a) yields Part (c). Parts (d)---(g) are direct consequences of Parts (a)---(b).  
\qquad$ \Box$

\medskip

\noindent {\bf Theorem 13.14.}\, {\em Let $ A \in {\cal C}^{m \times m}$
with ${\rm Ind}(A) = k.$ Then

{\rm (a)} \ $ r(\, A^*A^D - A^DA^* \,) =  r\left[ \begin{array}{ccc}
A^k( \, AA^* - A^*A \, )A^k  & 0  &
 A^kA^*   \\  0 & 0 & A^k  \\ A^*A^k  & A^k  & 0   \end{array} \right]
  - 2r(A^k).$

{\rm (b)} \ $ r(\, A^*A^{\#} - A^{\#}A^* \,) =  r\left[ \begin{array}{ccc}
A( \, AA^* - A^*A \, )A  & 0  &
 AA^*   \\  0 & 0 & A  \\ A^*A  & A  & 0   \end{array} \right]
  - 2r(A).$

{\rm (c)} \ $ r(\, A^*A^D - A^DA^* \,) =  r( \, A^{k+1}A^*A^k -
A^kA^*A^{k+1} \,),$ if $R(A^*A^k) \subseteq R(A^k)$ and $ R[A(A^k)^*]
\subseteq R[(A^k)^*].$

{\rm (d)} \  $ r(\, A^*A^D - A^DA^* \,) =  r\left[ \begin{array}{c}
A^kA^* \\ A^k  \end{array} \right]
 + r[ \, A^k,  \  A^*A^k \,] - 2r(A^k),$  if  $ A^{k+1}A^*A^k =
  A^kA^*A^{k+1}.$ 

{\rm (e)} \ $ A^*A^D = A^DA^*  \Leftrightarrow  R(A^*A^k)
\subseteq R(A^k), \ R[A(A^k)^*] \subseteq R[(A^k)^*]$ and $ A^{k+1}A^*A^k
= A^kA^*A^{k+1}.$ 

{\rm (f)} \ $ r(\, A^*A^{\#} - A^{\#}A^* \,) =  r( \, A^2A^*A - AA^*A^2\,),$ if $A$ is EP.

{\rm (g)} \ $ A^*A^{\#} = A^{\#}A^* \Leftrightarrow  A^2A^*A = AA^*A^2 $ and $A$ is EP
 $ \Leftrightarrow$  $ A$ is both EP and star-dagger.}
 
\medskip

\noindent {\bf Proof.}\, Applying (2.3) to $ A^*A^D - A^DA^*$ yields 
\begin{eqnarray*} 
r(\, A^*A^D - A^DA^* \,) & = & r[\, A^*A^{k}( A^{2k+1})^{\dagger}A^k  
 -  A^{k}( A^{2k+1})^{\dagger}A^kA^* \,] \\
 &= & r \left[ \begin{array}{ccc} -A^{2k+1} & 0
& A^k \\  0 & A^{2k+1} & A^kA^* \\ A^*A^k & A^k & 0  \end{array} \right] -
  2r(A^{2k+1})  \\
&= & r \left[ \begin{array}{ccc} -A^{2k+1} & 0 & A^k \\
 -A^{k+1}A^*A^k & 0 & A^kA^* \\ A^*A^k & A^k & 0  \end{array} \right] -
 2r(A^{k})  \\
&= & r \left[ \begin{array}{ccc} 0  & 0  & A^k \\
A^{k}A^*A^{k+1} -A^{k+1}A^*A^k  & 0 & A^kA^* \\ A^*A^k & A^k & 0
\end{array} \right] - 2r(A^{k}) \\
& = & r \left[ \begin{array}{ccc} A^{k}( \, AA^* -A^*A \, )A^k & 0
& A^kA^* \\
 0 & 0 & A^k \\ A^*A^k & A^k &  0 \end{array} \right] - 2r(A^{k}). 
\end{eqnarray*} 
as required for Part (a). Parts (b)---(d) are special cases of  Part (a). 
Applying Lemma 1.2(f) to the rank equality in Part (a) yields 
Part (e).  Parts (f) and (g) follow from Part (b).  \qquad  $ \Box$ 

\medskip

Similarly we can also establish  the following four theorems, which proofs are omitted. 

\medskip

\noindent {\bf Theorem 13.15.}\, {\em Let $ A \in {\cal C}^{m \times m}$
with ${\rm Ind}(A) = k.$ Then

{\rm (a)} \ $ r(\, AA^*A^D - A^DA^*A \,) =  r\left[ \begin{array}{ccc}
A^k( \, A^2A^* - A^*A^2 \, )A^k  & 0  &
 A^kA^*A   \\  0 & 0 & A^k  \\ AA^*A^k  & A^k  & 0   \end{array} \right]
  - 2r(A^k).$

{\rm (b)} \  $ r(\, A^kA^*A^{D} - A^{D}A^*A^k \,) = r( \, A^{2k+1}A^*A^k
 - A^kA^*A^{2k+1} \,).$ 

{\rm (c)} \ $ r(\, AA^*A^{\#} - A^{\#}A^*A \,) = r( \, A^3A^*A - AA^*A^3 \,),$ if 
  ${\rm Ind}(A) = 1.$ 

{\rm (d)} \ $ AA^*A^D = A^DA^*A  \Leftrightarrow  R(AA^*A^k) =  R(A^k), \ R[(A^kA^*A)^*] = R[(A^k)^*]$ and $ A^{k+2}A^*A^k = A^kA^*A^{k+2}.$ 

{\rm (e)}  \ $A^kA^*A^{D} = A^{D}A^*A^k \Leftrightarrow A^{2k+1}A^*A^k 
= A^kA^*A^{2k+1}.$ 

{\rm (f)\cite{HaSp2}}\, If $A$ is star-dagger$,$  then  $AA^*A^D = A^DA^*A.$

{\rm (g)} \ $ AA^*A^{\#} = A^{\#}A^*A \Leftrightarrow  A^3A^*A = AA^*A^3.$ }

\medskip

\noindent {\bf Theorem 13.16.}\, {\em Let $ A \in {\cal C}^{m \times m}$
with ${\rm Ind}(A) = k.$ Then

{\rm (a)} \ $ r(\,  AA^DA^* - A^*A^DA  \,) =  r\left[ \begin{array}{c}
A^k \\  A^kA^* \end{array} \right] + r[\,A^k,  \ A^*A^k \,] - 2r(A^k) = 
r[\,  A(A^k)^{\dagger} - (A^k)^{\dagger}A  \,].$  

{\rm (b)} \  $ r(\, AA^{\#}A^* - A^*A^{\#}A \,) = 2r[\,A,  \ A^* \,] - 2r(A) ,$ if 
  ${\rm Ind}(A) = 1.$ 

{\rm (c)} \ $ AA^DA^* = A^*A^DA \Leftrightarrow 
A(A^k)^{\dagger} = (A^k)^{\dagger}A \Leftrightarrow R(A^*A^k) =  R(A^k)$ and $R[A(A^k)^*]  =  R[(A^k)^*].$ 

{\rm (d)} \ $ AA^{\#}A^* = A^*A^{\#}A \Leftrightarrow$  A is EP. }

\medskip

\noindent {\bf Theorem 13.17.}\, {\em Let $ A \in {\cal C}^{m \times m}$
with ${\rm Ind}(A) = k.$ Then 

{\rm (a)} \ $ r[\,  AA^D(A^*)^k - (A^*)^kA^DA \,] = 2 r[\,A^k,  \ (A^k)^* \,] - 2r(A^k).$ 

{\rm (b)} \ $ AA^D(A^*)^k = (A^*)^kA^DA \Leftrightarrow $  $A^k$ is EP.}

\medskip

\noindent {\bf Theorem 13.18.}\, {\em Let $ A \in {\cal C}^{m \times m}$
with ${\rm Ind}(A) = 1$ and  $ \lambda$ is a nonzero complex number. Then

{\rm (a)}  \ $ r[\, AA^{\#}(\, AA^* +  \lambda A^*A \,) - (\, AA^* +  
\lambda A^*A \,)A^{\#}A \,] = 2r[\,A,  \ A^* \,] - 2r(A). $

{\rm (b)}  \ $AA^{\#}$ commutes with $ AA^* + \lambda A^*A \Leftrightarrow $ 
   A is EP.}

\medskip

\noindent {\bf Theorem 13.19.}\, {\em Let $ A \in {\cal C}^{m \times m}$
with ${\rm Ind}(A) = k$. Then

{\rm (a)}  \ $ r[\, (AA^{D})^*A^{\dagger} -A^{\dagger}(AA^{D})^* \,] 
= r\left[ \begin{array}{c} A^kA^*A \\  A^k \end{array} \right] + r[\,AA^*A^k,  \ A^k \,] - 2r(A^k).$

{\rm (b)} \ $(AA^{D})^*A^{\dagger} = A^{\dagger}(AA^{D})^*$ $ \Leftrightarrow$ 
$ R(AA^*A^k) = R(A^k)$ and $ R[(A^kA^*A)^*] = R[(A^k)^*].$

{\rm (c)} \ $(AA^{\#})^*A^{\dagger} = A^{\dagger}(AA^{\#})^*$, if ${\rm Ind}(A) = 1.$ }

\medskip

\noindent {\bf Proof.}\, Apply (4.1) to $(AA^{D})^*A^{\dagger} - A^{\dagger}(AA^{D})^*$ to yield 
\begin{eqnarray*} 
 r[\, (AA^{D})^*A^{\dagger} -A^{\dagger}(AA^{D})^* \,] & = & r\left[ \begin{array}{c}  
(AA^{D})^*A^{\dagger} \\ (AA^{D})^*  \end{array} \right] + 
r[\, A^{\dagger}(AA^{D})^*, \ (AA^{D})^* \,] - 2r(AA^{D}) \\
 & = & r\left[ \begin{array}{c}  (A^k)^*A^{\dagger} \\ (A^k)^* \end{array} \right] + 
r[\, A^{\dagger}(A^k)^*, \ (A^k)^* \,] - 2r(A^k).
\end{eqnarray*} 
Next applying (2.1), we find that 
$$ 
r\left[ \begin{array}{c}  (A^k)^*A^{\dagger} \\ (A^k)^*  \end{array} \right] 
= r[\,AA^*A^k,  \ A^k \,]  \ \ {\rm and}  \ \  r[\, A^{\dagger}(A^k)^*, \ (A^k)^* \,] 
= r\left[ \begin{array}{c} A^kA^*A \\  A^k \end{array} \right].
$$
Thus we get Part (a).  \qquad  $ \Box$ 

\medskip

\noindent {\bf Theorem 13.20.}\, {\em Let $ A \in {\cal C}^{m \times m}$
with ${\rm Ind}(A) = k.$ Then

{\rm (a)}  \ $ r[\, A - (A^D)^{\dagger} \,] =  r\left[ \begin{array}{ccc}
A  & A^k &
(A^k)^*  \\  A^k & 0 & 0  \\ (A^k)^*  & 0 & 0 \end{array} \right] -
2r(A^k).$

{\rm (b)} \ $  r[\, A - (A^D)^{\dagger} \,] = r(A) -r(A^k),$ \ if $ A^k $ is EP.

{\rm (c)} \ $ (A^{\#})^{\dagger} = A  \Leftrightarrow  A $ is EP. }

\medskip

\noindent {\bf Proof.}\, According to Cline's identity 
$ (A^D)^{\dagger} = (A^k)^{\dagger}A^{2k+1}(A^k)^{\dagger}$(see \cite{BeG} and \cite{Cl3}), we find  by (2.8) that
\begin{eqnarray*} 
r[\, A - (A^D)^{\dagger} \,] & = &r[\, A -
(A^k)^{\dagger}A^{2k+1}(A^k)^{\dagger} \,] \\
& = & r \left[ \begin{array}{ccc} (A^k)^*A^{2k+1}(A^k)^* & (A^k)^*A^k(A^k)^*
 & 0 \\
  (A^k)^*A^k(A^k)^* & 0  & (A^k)^* \\ 0 & (A^k)^* & -A \end{array} \right]
  - 2r(A^{k})  \\
& = & r \left[ \begin{array}{ccc} A^{2k+1} & A^k(A^k)^*  & 0 \\
  (A^k)^*A^k & 0 & (A^k)^* \\ 0 & (A^k)^* & -A \end{array} \right] -
  2r(A^{k})  \\
& = & r \left[ \begin{array}{ccc} A^{2k+1} & 0  &  A^{k+1} \\
  (A^k)^*A^k & 0 & (A^k)^* \\ 0 & (A^k)^* & -A \end{array} \right] -
  2r(A^{k})  \\
& = & r \left[ \begin{array}{ccc} 0 & 0  &  A^{k+1} \\ 0 & 0 & (A^k)^*
\\ A^{k+1} & (A^k)^* & -A \end{array} \right] - 2r(A^{k})  \\ 
& = & r\left[
\begin{array}{ccc} A  & A^k  & (A^k)^*  \\ A^k & 0 & 0  \\ (A^k)^*
& 0 & 0\end{array} \right] - 2r(A^k).
\end{eqnarray*} 
as required for Part (a).  The results in Part (b) and (c) follows
 immediately from Part (a). \qquad$ \Box$ 

\medskip

\noindent {\bf Theorem 13.21.}\, {\em Let $ A \in {\cal C}^{m \times m}$
with ${\rm Ind}(A) = k.$ Then

{\rm (a)}  \ $ r[\, AA^DA - (A^D)^{\dagger} \,] = 2r[ \, A^k, \ (A^k)^* \,] - 2r(A^k).$

{\rm (b)}  \ $ (A^D)^{\dagger} = AA^DA  \Leftrightarrow  A^k $ is EP. } 

\medskip

\noindent {\bf Proof.}\, It is easy to verify that both $ AA^DA$ and $ (A^D)^{\dagger}$
 are  outer inverses of $ A^D$. In that case it follows from (5.1) that
\begin{eqnarray*}
r[\, AA^DA - (A^D)^{\dagger} \,] &= &
r \left[ \begin{array}{c}  AA^DA  \\ (A^D)^{\dagger}  \end{array} \right]
+ r[\, AA^DA, \  (A^D)^{\dagger} \,]  - r(AA^DA) - r[(A^D)^{\dagger}] \\
& = & r \left[ \begin{array}{c} A^D \\ (A^D)^* \\ \end{array} \right] +
r[ \, A^D, \ (A^D)^* \,] - 2r(A^k) \\
& = & r \left[ \begin{array}{c} A^k \\ (A^k)^* \\ \end{array} \right] +
r[ \, A^k, \ (A^k)^* \,] - 2r(A^k),
\end{eqnarray*}
as required for Part (a). The result in Part (b) follows immediately from
Part (a). \qquad $ \Box$ 

\medskip

\noindent {\bf Theorem 13.22.}\, {\em Let $ A, \, B \in {\cal C}^{m \times m}$
with ${\rm Ind}(A) = k$ and $ {\rm Ind}(B) = l.$ Then

{\rm (a)} \ $ r(\, AA^D - BB^D \,) =  r\left[ \begin{array}{c} A^k \\ B^l
 \end{array} \right] + r[\ A^k, \ B^l \,] - r(A^k) -r(B^l).$

{\rm (b)} \  $ r(\, AA^{\#} - BB^{\#} \,) =  r\left[ \begin{array}{c} A \\ B
 \end{array} \right] + r[\ A, \ B \,] - r(A) -r(B),$ If ${\rm Ind}(A) =
 {\rm Ind}(B) = l.$

{\rm (c)} \ $  AA^D = BB^D \Leftrightarrow   R(A^k) = R(B^l)$
and $ R[(A^k)^*] = R[(B^l)^*].$

{\rm (d)}  \ $ r(\, AA^D - BB^D \,)$ is nonsingular  $ \Leftrightarrow 
r\left[ \begin{array}{c} A^k \\ B^l   \end{array} \right] = r[ \ A^k,
\ B^l \,] = r(A^k)  + r(B^l) = m  \Leftrightarrow R(A^k)
\oplus R(B^l) =  \ R[(A^k)^*] \oplus R[(B^l)^*] = {\cal C}^m. $

{\rm (e)}\, In particular$,$ if $ {\rm Ind}\left[ \begin{array}{cc} A & B \\ 0 & D \end{array} \right] = 1.$ Then
$$
r \left(\,  \left[ \begin{array}{cc} A & B \\ 0 & D \end{array} \right]
\left[ \begin{array}{cc} A & B \\ 0 & D \end{array} \right]^{\#} - 
\left[ \begin{array}{cc} AA^{\#} & 0 \\ 0 & DD^{\#} \end{array} \right] \, \right)
= r[\, A, \ B \, ] + r\left[ \begin{array}{cc} B \\ D \end{array} \right] - 
r\left[ \begin{array}{cc} A & B \\ 0 & D \end{array} \right].
$$ } 

\medskip

\noindent {\bf Proof.}\, Note that both $ AA^D$ and $BB^D$ are idempotent.
Then it follows from (3.1) that
\begin{eqnarray*}
r(\, AA^D - BB^D \,)
&= & r \left[ \begin{array}{c} AA^D \\ BB^D \end{array} \right]
+ r[\, AA^D, \ BB^D \,] - r(AA^D) - r(BB^D) \\
&= & r \left[ \begin{array}{c} A^D \\ B^D \end{array} \right]
+ r[\, A^D, \ B^D \,] - r(A^D) - r(B^D) \\
& = & r\left[ \begin{array}{c} A^k \\ B^l   \end{array} \right] + r[\, A^k, \ B^l \,]
- r(A^k) -r(B^l),
\end{eqnarray*} 
as required for  Part (a). The results in Parts (b)---(e) follow
immediately from Part (a). \qquad$ \Box$

\medskip

\noindent {\bf Theorem 13.23.}\, {\em Let $ A, \, B \in {
\cal C}^{m \times m}$ with ${\rm Ind}(A) =k.$ Then

{\rm (a)} \ $ r(\, AA^DB - BA^DA \,) =  r\left[ \begin{array}{c} A^k \\
A^kB   \end{array} \right]
+ r[ \, A^k, \ BA^k \, ] - 2r(A^k).$

{\rm (b)} \ $ r(\, A^DAA^{\dagger} - A^{\dagger}AA^D \,) =  r\left[
\begin{array}{c} A^k \\ A^* \end{array} \right] + r[ \, A^k, \ A^* \,]
- 2r(A) = r(\, A^DA^{\dagger} - A^{\dagger}A^D \,).$ \\
In particular$,$

{\rm (c)} \ $  AA^DB = BA^DA  \Leftrightarrow  R(BA^k) = R(A^k)$
and $ R[(A^kB)^*] = R[(A^k)^*].$

{\rm (d)} \ $ A^DAA^{\dagger} = A^{\dagger}AA^D   \Leftrightarrow
  A^DA^{\dagger} = A^{\dagger}A^D   \Leftrightarrow   R(A^k)
 \subseteq R(A^*)$ and $ R[(A^k)^*] \subseteq R(A).$ } 

\medskip

\noindent {\bf Proof.}\, Note that  $ AA^D =A^DA$ is idempotent. It follows
 by (4.1) that
\begin{eqnarray*} 
r(\, AA^DB - BA^DA \,)
 &= & r \left[ \begin{array}{c} AA^DB \\ A^DA \end{array} \right]
  + r[\, BA^DA, \ AA^D \,] -r(AA^D) -r(A^DA) \\
 &= & r \left[ \begin{array}{c} A^DB \\ A^D \end{array} \right]
  + r[\, BA^D, \ A^D \,] - 2r(A^D) \\
 &= & r \left[ \begin{array}{c} A^kB \\ A^k \end{array} \right]
  + r[\, BA^k, \ A^k \,] - 2r(A^k).
\end{eqnarray*} 
Thus we have Parts (a). Replacing $ B$ by $ A^{\dagger}$ in Part (a) and
simplifying it yields the first equality in Part (b). The second equality in
Part (b) follows from Theorem 13.22(a). \qquad $ \Box$ 

\medskip

\noindent {\bf Theorem 13.24.}\, {\em Let $ A \in {\cal C}^{m \times n},
\, B \in {\cal C}^{m \times m}$ with ${\rm Ind}(B) =k$ and $ C \in
{\cal C}^{n \times n}$ with ${\rm Ind}(C) = l.$ Then

{\rm (a)} \  $ r(\, BB^DA - AC^DC \,) =  r\left[ \begin{array}{c} B^kA
\\ C^l  \end{array} \right]
+ r[ \, AC^l, \ B^k \,] - r(B^k) -r(C^l).$

{\rm (b)} \  $  BB^DA = AC^DC   \Leftrightarrow   R(AC^l)
\subseteq R(B^k)$ and
$ R[(B^kA)^*] \subseteq R[(C^l)^*].$ } 

\medskip

\noindent {\bf Proof.}\, Note that both $BB^D$ and $C^DC$ are idempotent.
Then it follows from  (4.1) that
\begin{eqnarray*} 
r(\, BB^DA - AC^DC \,)
&= &  r \left[ \begin{array}{c} BB^DA \\ C^DC \end{array} \right]
+ r[\, AC^DC, \ BB^D \,] - r(BB^D) -  r(C^DC) \\
&= &  r \left[ \begin{array}{c} B^DA \\ C^D \end{array} \right]
+ r[\, AC^D, \ B^D \,] - r(B^D) -  r(C^D) \\
&= &  r \left[ \begin{array}{c} B^kA \\ C^l \end{array} \right]
+ r[\, AC^l, \ B^k \,] - r(B^k) -  r(C^l),
\end{eqnarray*}
as required for Part (a).  \qquad $ \Box$

\medskip

\noindent {\bf Theorem 13.25.}\, {\em Let $ A \in {\cal C}^{m \times n},
 \, B \in {\cal C}^{m \times m} $ with
${\rm Ind}(B) = k $ and  $ C \in {\cal C}^{n \times n}$ with ${\rm Ind}(C)
= l.$ Then

{\rm (a)} \ $ r[ \, A, \ B^k \, ] = r(B^k) + r( \, A - BB^DA \,).$

{\rm (b)} \ $ r \left[ \begin{array}{c}  A  \\ C^l  \end{array} \right]
 = r(C^l) + r(\, A - AC^DC \,).$

{\rm (c)} \ $  r\left[ \begin{array}{cc}  A  & B^k  \\ C^l &  0 \end{array}
\right]
= r(B^k) + r(C^l) + r [ \,( \ I_m - BB^D \, )A(\, I_n - C^DC \, ) \, ]. $ }
 
\medskip

\noindent {\bf Proof.}\, Applying  (1.7) to $ A - BB^DA$ yields 
\begin{eqnarray*} 
r(\, A - BB^DA \,) &= & r[\,A  - B^{k+1}( B^{2k+1})^{\dagger}B^kA \,]  \\
&= &  r \left[ \begin{array}{cc} B^{2k+1} & B^{k}A \\ B^{k+1} & A
\end{array} \right] -r(B^{2k+1}) \\
 & = &r \left[ \begin{array}{cc} 0 & 0
\\ B^{k+1} & A \end{array} \right] -r(B^{k})  =  r[ \, A, \ B^k \, ] - r(B^{k}),
\end{eqnarray*} 
as required for Part (a). Similarly we can show Parts (b) and (c).
\qquad $ \Box$ 

\medskip

\noindent {\bf Theorem 13.26.}\, {\em Let $ A, \, B \in {\cal C}^{m \times m}$
with
${\rm Ind}(A) = k$ and ${\rm Ind}(B) = l.$ Then

{\rm (a)} \  $ r(\, AB - ABB^DA^DAB \,) =  r\left[ \begin{array}{cc} A^{2k}
& A^kB^l   \\ B^lA^k
 & B^{2l}    \end{array} \right] + r(AB) - r(A^k) - r(B^l).$ 

{\rm (b)} \ $ r(\, AB - ABB^{\#}A^{\#}AB \,) =  r\left[ \begin{array}{cc} A^{2}
& AB   \\ BA & B^{2} \end{array} \right] + r(AB) - r(A) - r(B),$ if 
 ${\rm Ind}(A) = {\rm Ind}(B) = 1.$

{\rm (c)} \ $ B^DA^D \in \{ (AB)^- \}  \Leftrightarrow 
 r\left[ \begin{array}{cc} A^{2k} & A^kB^l \\
 B^lA^k  & B^{2l} \end{array} \right] = r(A^k) + r(B^l) - r(AB).$  

{\rm (d)} \  $ B^{\#}A^{\#} \in \{ (AB)^- \}  \Leftrightarrow 
 r\left[ \begin{array}{cc} A^{2} & AB \\
 BA  & B^{2} \end{array} \right] = r(A) + r(B) - r(AB).$ }

\medskip 

\noindent {\bf Proof.}\, It follows by (2.9) that 
\begin{eqnarray*} 
r(\, AB - ABB^DA^DAB \,) &
= & r[\,  AB - AB^{k+1}( B^{2k+1})^{\dagger}B^lA^{k}( A^{2k+1})^{\dagger}
A^{k+1}B \,] \\
 &= & r \left[ \begin{array}{ccc} B^{l}A^k & B^{2l+1} & 0 \\ A^{2k+1}  & 0
 & A^{k+1}B \\
  0 & AB^{l+1} &  -AB \end{array} \right] - r(A^{2k+1}) -  r(B^{2l+1}) \\
 &= & r \left[ \begin{array}{ccc} B^{l}A^k & B^{2l+1} & 0 \\  A^{2k+1}
  &  A^{k+1}B^{l+1}  & 0 \\ 0 & 0 &  -AB \end{array} \right] - r(A^{k}) -  r(B^{l}) \\ 
 &= & r \left[ \begin{array}{cc} B^{l}A^k & B^{2l+l}  \\ A^{2k +1}  & 
 A^{k}B^{l+l} \end{array} \right]  + r(AB) - r(A^{k}) - r(B^{l}) \\
& = & r \left[ \begin{array}{cc} B^{l}A^k & B^{2l}  \\ A^{2k}  &  A^{k}B^{l}
\end{array} \right]  + r(AB) - r(A^{k}) - r(B^{l}).
\end{eqnarray*} 
Thus we have  Parts (a).  \qquad $ \Box$

\medskip

\noindent {\bf Theorem 13.27.}\, {\em Let $ A, \ B \in
{\cal C}^{m \times m}$ with ${\rm Ind}(A) = k$ and ${\rm Ind}(B) = l.$ Then

{\rm (a)} \ $ r(\, AA^DB^DB -  BB^DA^DA \,) =
  r\left[ \begin{array}{c} A^{k} \\ B^l \end{array} \right] +
 r[\, A^k, \ B^l \,] + r(A^kB^l)  + r(B^lA^k) - 2r(A^k) - 2r(B^l).$

{\rm (b)}\  $ r(\, AA^{\#}B^{\#}B -  BB^{\#}A^{\#}A \,) =
  r\left[ \begin{array}{c} A \\ B \end{array} \right] +
 r[\, A, \ B \,] + r(AB)  + r(BA) - 2r(A) - 2r(B),$ if 
${\rm Ind}(A) = {\rm Ind}(B) = 1.$

{\rm (c)} \ $ AA^DB^DB = BB^DA^DA  \Leftrightarrow
r\left[ \begin{array}{c} A^{k}  \\ B^l \end{array} \right]
= r(A^k) + r(B^l) - r(A^kB^l) \ and  \ 
r[ \, A^k, \ B^l \,] = r(A^k) + r(B^l) - r(B^lA^k).$ 

{\rm (d)} \ $AA^{\#}B^{\#}B = BB^{\#}A^{\#}A \Leftrightarrow
r\left[ \begin{array}{c} A  \\ B \end{array} \right]
= r(A) + r(B) - r(AB) \ and  \ 
r[ \, A, \ B \,] = r(A) + r(B) - r(BA).$ } 

\medskip

\noindent {\bf Proof.}\, Note that both $ AA^D = A^DA$ and $BB^D = B^DB$ are
idempotent. Then it follows by  (3.26) that
\begin{eqnarray*} 
\lefteqn{r(\, AA^DB^DB -  BB^DA^DA  \,) } \\
& = & r\left[ \begin{array}{c} AA^D \\ B B^D \end{array} \right] +
r[ \, A^DA, \ B^DB \,] + r(AA^DB^DB)  + r(BB^DA^DA) - 2r(AA^D) - 2r(BB^D)\\
& = & r\left[ \begin{array}{c} A^k \\ B^l \end{array} \right] +
r[ \, A^k, \ B^l \,] + r(A^kB^l) + r(B^lA^k) - 2r(A^k) - 2r(B^l),
\end{eqnarray*}
as required for Part (a).  \qquad $ \Box$ 

\medskip

\noindent {\bf Theorem 13.28.}\, {\em Let $ A, \, B  \in
{\cal C}^{ m \times m}$ with
 ${\rm Ind}( \, A+B \,) = k$ and
 denote $ N = A + B $. Then

{\rm (a)} \ $ r(AN^DB ) = r(AN^k) + r(N^kB ) - r(N^k).$
 
{\rm (b)} \ $ r(\, AN^DB \, ) = r(A) + r(B) - r(N^k),$  \ \ if $ R(B)
 \subseteq R(N^k) \ \ and \ \  R(A^*) \subseteq R[(N^k)^* ].$

{\rm (c)} \ $ r(\, AN^DB - BN^DA \,) = r \left[ \begin{array}{c} N^k \\
 N^kB \end{array} \right]
 + r[\, N^k,  \ BN^k \,] - 2r(N^k).$   

{\rm (d)} \ $ AN^DB =  BN^DA  \Leftrightarrow  R(BN^k) \subseteq R(N^k)
 \ \ and \ \
 R[(N^kB)^*] \subseteq R[(N^k)^*].$  } 

\medskip

\noindent {\bf Proof.}\, It follows by (1.7), that 
\begin{eqnarray*}
r(AN^DB )& = &r[ \, A N^{k}( N^{2k+1})^{\dagger}N^kB \, ] \\
& = &  r \left[ \begin{array}{cc} N^{2k+1}   & N^kB  \\ AN^k & 0
  \end{array} \right] - r(N^{2k+1} ) \\
& = &  r \left[ \begin{array}{cc} 0  & N^kB  \\ AN^k & 0 \end{array} \right]
 - r(N^{k} )= r(AN^k) + r(N^kB ) - r(N^k),
\end{eqnarray*}
which is the first equality in Part (a). The second equality in Part (a)
follows from $  r(AN^k) = r(AN^D),$  $r(N^kB ) = r(N^DB)$ and $ r(N^D) = r(N^k).$  Under $ R(B) \subseteq R(N^k)$
 and $ R(A^*) \subseteq
 R[(N^k)^*],$ it follows  that $ r(AN^k) = r(A)$ and $ r(N^kB) = r(B)$. Thus
 Part (a) becomes Part (b). Next applying (2.3) to  $AN^DB - BN^DA$ yields
\begin{eqnarray*}
\lefteqn{r ( \, AN^DB - BN^DA \,) }\\
& = &  r[ \, AN^{k}( N^{2k+1})^{\dagger}N^kB -
BN^{k}( N^{2k+1})^{\dagger}N^kA \, ]  \\
&= &  r \left[ \begin{array}{ccc} -N^{2k+1}  & 0 & N^{k}B \\ 0 &  N^{2k+1} &
 N^kA \\ AN^{k} &  BN^{k}  & 0 \end{array} \right] - 2r(N^{2k+1})\\
&= &  r \left[ \begin{array}{ccc} -N^kAN^k  & -N^kBN^k & N^kB \\  N^kAN^k &
 N^kBN^k & N^kA \\ AN^{k} &  BN^{k}  & 0 \end{array} \right] - 2r(N^{k})\\
&= &  r \left[ \begin{array}{ccc} 0 & 0 & N^kB \\  0 & 0 & N^kA \\ AN^{k} &
 BN^{k}  & 0 \end{array} \right] - 2r(N^{k})\\
&= & r \left[ \begin{array}{c} N^kA \\ N^kB \end{array} \right] +
r[ \, AN^k,  \ BN^k \,] - 2r(N^k)
= r \left[ \begin{array}{c} N^k \\ N^kB \end{array} \right] +
r[ \, N^k,  \ BN^k \,] - 2r(N^k).
\end{eqnarray*}
Thus we have Parts (c) and (d).  \qquad $\Box$

\medskip

\noindent {\bf Theorem 13.29.}\, {\em Let $ A, \, B  \in
 {\cal C}^{ m \times m}$ be given$,$ and let $ N = A + B $ with 
${\rm Ind}( \, A+B \,) = k$. Then

{\rm (a)} \ $ r \left(  \, \left[ \begin{array}{cc} A  & 0 \\ 0 & B
 \end{array} \right]  - \left[ \begin{array}{c} A \\ B  \end{array} \right]
 ( \, A + B \, )^D[ \, A, \ B \,] \,
\right) =  r(A) + r(B) - r(N^k).$ 

{\rm (b)} \ $  \left[ \begin{array}{c} A \\ B  \end{array} \right]( \, A + B
 \, )^D [\, A,  \ B \, ] = \left[ \begin{array}{cc} A  & 0 \\ 0 & B
  \end{array} \right] \Leftrightarrow 
 {\rm Ind}( \, A+B \,) \leq 1 $ and $ r( \, A + B  \,) = r(A) + r(B).$ }
  
\medskip

\noindent {\bf Proof.}\, It follows by  (1.7) that 
\begin{eqnarray*} 
\lefteqn { r \left(\, \left[ \begin{array}{cc} A  & 0 \\ 0 & B \end{array} \right] -
 \left[ \begin{array}{c} A \\ B  \end{array} \right] N^D[ \, A, \ B \, ] \,
 \right) } \\
& = & r \left(\, \left[ \begin{array}{cc} A  & 0 \\ 0 & B \end{array}
\right] - \left[ \begin{array}{c} A \\ B  \end{array} \right]
 N^k ( N^{2k+1})^{\dagger} N^k[ \, A, \ B \, ] \, \right) \\
&= &  r \left[ \begin{array}{ccc} N^{2k+1} & N^{k}A  & N^{k}B \\ AN^k & A &
 0 \\ BN^{k} & 0 & B \end{array} \right] - r(N^k) \\
&= &  r \left[ \begin{array}{ccc} 0 & 0 & 0  \\ 0 & A  &  0 \\ 0 & 0 & B
\end{array} \right] - r(N^k) = r(A) + r(B) - r(N^k),
\end{eqnarray*} 
which is exactly Part (a). Note that $ r(N^k) \leq r(N) = r(\, A + B \,)
\leq r(A) + r(B).$ Thus $ r(N^k) = r(A) + r(B)$ is equivalent to
$ {\rm Ind}(N) \leq 1 $ and $ r(N) = r(A) + r(B).$ \qquad $ \Box$ 

\medskip

In general we have the following. 

\medskip

\noindent {\bf Theorem 13.30.}\, {\em Let $ A_1, \,  A_2, \, \cdots, \,  A_k
\in {\cal C}^{ m \times m}$ with
$ {\rm Ind}(N) = k,$ where $ N =  A_1 + A_2 + \cdots  +A_k,$ and let
$ A = {\rm diag}( \, A_1, \, A_2, \, \cdots , \, A_k \, ).$ Then

{\rm (a)} \ $ r \left( \, A  - \left[ \begin{array}{c} A_1 \\ \vdots \\
A_k \end{array} \right]N^D[ \, A_1,  \ \cdots, \ A_k \, ]  \,  \right) =
r(A_1) + \cdots +  r(A_k) - r(N^k). $

{\rm (b)} \ $ \left[ \begin{array}{c} A_1 \\ \vdots \\ A_k \end{array}
\right]N^D[ \, A_1,  \, \cdots, \, A_k \,] = A  \Leftrightarrow 
{\rm Ind}(N) \leq 1 $ and  $r(N) = r(A_1) + \cdots + r(A_k). $ } 

\medskip

\noindent {\bf Theorem 13.31.}\,  {\em  Let $ A_1, \, A_2, \, \cdots, \ ,A_k 
 \in {\cal C }^{ m \times m}.$  
Then the  Drazin inverse of their sum satisfies the following equality
$$
 ( \,A_1 +  A_2 + \cdots + A_k \, )^D 
= \frac{1}{k}[\, I_m, \,  I_m, \, \cdots, \, I_m \,] 
\left[ \begin{array}{cccc} A_1 & A_2 & \cdots  & A_{k}  \\  A_{k} & A_1 & \cdots  & A_{k-1}  \\
 \vdots  & \vdots & \ddots  & \vdots \\ A_2 & A_3 & \cdots  & A_1 \end{array} \right]^D 
\left[ \begin{array}{c} I_m \\ I_m \\ \vdots 
 \\ I_m \end{array} \right]. \eqno (13.1) 
$$} 
{\bf Proof.}\, Since the given matrices are  square, (11.7) can be written as 
$$
U^*_m A U_m = {\rm diag}(\, J_1, \, J_2, \, \cdots, \, J_k \ ).
$$
In that case, it is easy to verify that
$$
(U^*_m A U_m)^D = U^*_m A^DU_m,
$$  
and 
$$
[\, {\rm diag}(\, J_1, \, J_2, \, \cdots, \, J_k \, ) \,]^D 
= {\rm diag}(\, J_1^D, \, J_2^D, \, \cdots, \, J_k^D \, )
$$
Thus we have 
\begin{eqnarray*}
J_1^D  = [\, I_m, \, 0, \, \cdots, \, 0  \, ] U^*_m A^DU_m 
[\, I_m, \ 0, \ \cdots, \ 0  \, ]^T  =   \frac{1}{k}[\, I_m, \, I_m, \, \cdots, \, I_m  \, ]A^D [\, I_m, \, I_m, \, \cdots, \, I_m \, ]^T,
\end{eqnarray*}
which is (13.1). \qquad $ \Box$

\medskip

\noindent {\bf Theorem 13.32.}\,  {\em  Let $ A + iB \in {\cal C }^{ m \times m},$  where $ A $ and $ B $ are real. 
Then the  Drazin inverse of $ A + iB $ satisfies the identity
$$
 ( \,A + iB \, )^D 
= \frac{1}{2}[\, I_m, \ iI_m \,]  \left[ \begin{array}{rr} A & -B \\ B & A  \end{array} \right]^D 
\left[ \begin{array}{c} I_m \\ -iI_m \\\end{array} \right]. \eqno (13.2) 
$$} 
{\bf Proof.}\, Observe that 
$$
\left[ \begin{array}{rr} A & iB \\ iB  & A  \end{array} \right] = \left[ \begin{array}{cc} I_m & 0  \\ 0 & iI_m \end{array} 
\right]\left[ \begin{array}{rr} A & -B \\ B & A  \end{array} \right] \left[ \begin{array}{rr} I_m & 0  \\ 0 & iI_m 
 \end{array} \right]^{-1}.
$$ 
Thus 
$$ 
\left[ \begin{array}{rr} A & iB \\ iB  & A  \end{array} \right]^D = \left[ \begin{array}{cc} I_m & 0  \\ 0 & iI_m \end{array} 
\right]\left[ \begin{array}{rr} A & -B \\ B & A  \end{array} \right]^D \left[ \begin{array}{cc} I_m & 0  \\ 0 & iI_m 
 \end{array} \right]^{-1}.  \eqno (13.3) 
$$ 
In that case, applying (13.1) and then (13.3) to $ A +iB$ yields (13.2). \qquad $ \Box$ 

\medskip 

The identities in (11.31)---(11.33) can also be extended to the Drazin inverse of a real quaternion matrix.  

\medskip

\noindent {\bf Theorem 13.33.}\,  {\em  Let  $A_0 +  iA_1 + iA_2 +  kA_3$ be an $ m \times m $ real quaternion 
matrix.  Then its Drazin inverse satisfies    
$$ 
 ( \,  A_0 +  iA_1 + iA_2 +  kA_3 \, )^D  = \frac{1}{2}[ \, I_m , \ jI_m \, ]
  \left[ \begin{array}{cc}  A_0 +  iA_1 &  -(\, A_2 +  iA_3 \, )  \\  A_2 - iA_3&  A_0 - iA_1 
\end{array} \right]^D\left[ \begin{array}{c}  I_m \\ -jI_m  \end{array} \right], \eqno (13.4) 
$$ 
and 
$$ 
 ( \,  A_0 +  iA_1 + iA_2 +  kA_3 \, )^{D}  = \frac{1}{4}[ \, I_m , \  iI_m , \ jI_m, \  
kI_m \, ] \left[ \begin{array}{rrrr}  A_0 & -A_1 & -A_2 &  - A_3   \\  A_1 &  A_0 &  A_3 &  - A_2 
\\  A_2 &  -A_3 &  A_0  & A_1 \\ A_3 &  A_2 &  -A_1  & A_0 \end{array} \right]^D 
\left[ \begin{array}{c}  I_m \\ -iI_m \\ -jI_m  \\ -kI_m  \end{array} \right]. \eqno (13.5) 
$$ 
Moreover denote $ ( \,  A_0 +  iA_1 + iA_2 +  kA_3 \, )^{D} = G_0 +  iG_1 + iG_2 +  kG_3.$ Then 
$$
\left[ \begin{array}{rrrr}  A_0 & -A_1 & -A_2 &  - A_3   \\  A_1 &  A_0 &  A_3 &  - A_2 
\\  A_2 &  -A_3 &  A_0  & A_1 \\ A_3 &  A_2 &  -A_1  & A_0 \end{array} \right]^D = 
\left[ \begin{array}{rrrr}  G_0 & -G_1 & -G_2 &  - G_3   \\  G_1 &  G_0 &  G_3 &  - G_2 
\\  G_2 &  -G_3 &  G_0  & G_1 \\ G_3 &  G_2 &  -G_1  & G_0 \end{array} \right]^D. \eqno (13.6) 
$$}

As is well known that Drazin inverses of block matrices are quite difficult to determine in general. However, if a 
block matrix has some special pattern, its Drazin inverse can still be presented. Motivated by the expressions (9.78) 
and (8.82)---(9.85),  we can find the following. 

Let 
$$
M =  \left[ \begin{array}{cccc} A & B & \cdots  & B \\  B & A & \cdots  & B \\  
\vdots & \vdots & \ddots & \vdots \\ B & B & \cdots & A  \end{array} \right]_{ k\times k}, \eqno (13.7)
$$
where both $A$ and $B$ are $ m \times m$ matrices. Then 
$$
M^{D}=  \left[ \begin{array}{cccc} S & T & \cdots  & T \\  T & S & \cdots  & T \\  
\vdots & \vdots & \ddots & \vdots \\ T & T & \cdots & S  \end{array} \right]_{ k\times k}, \eqno (13.8)
$$ 
where 
$$
S = \frac{1}{k}[\, A + (k-1)B \,]^{D} + \frac{k-1}{k}(\, A - B \,)^D , \qquad  
T = \frac{1}{k}[\, A + (k-1)B \,]^{D} - \frac{1}{k}(\, A - B \,)^D. \eqno (13.9)
$$   
In fact we see from (9.77) that $ M = P_mNP_m^{-1}$, when  $m = n$. In that case,   $ M^D = P_mN^DP_m^{-1}$ holds.
Written in an explicit form, it is (13.7).  
The expression (9.78) illustrates that $M^{D}$ also has the same pattern as $M$.

\markboth{YONGGE  TIAN }
{14. RANK EQUALITIES FOR SUBMATRICES IN  DRAZIN INVERSES}

\chapter{Rank equalities for submatrices in Drazin inverses}

\noindent Let 
$$ 
M = \left[ \begin{array}{cc} A  & B  \\ C & D  \end{array} \right]
 \eqno (14.1)
$$ 
be a square block matrix over ${\cal C}$, where
 $ A \in {\cal C}^{m \times m} $ and  $D \in {\cal C}^{n \times n}$,
$$ 
V_1 = \left[ \begin{array}{c} A \\ C  \end{array} \right], \ \ \ \
 V_2 = \left[ \begin{array}{c} B \\ D  \end{array} \right], \ \ \ \
  W_1 = [ \, A, \ B \, ], \ \ \ \ W_2 = [\, C, \ D \,], \eqno (14.2)
$$ 
and partition the Drazin inverse of $ M $ as 
$$ 
M^D = \left[ \begin{array}{cc} G_1  & G_2  \\ G_3 & G_4  \end{array}
\right], \eqno (14.3)
$$ 
where $G_1 \in {\cal C}^{m \times m}$. It is, in general, quite difficult to 
give the expression of $ G_1$---$G_4$.  In this chapter we consider a simpler problem---the ranks 
of the submatrices $ G_1$---$G_4$ in  (14.3). 

\medskip

\noindent {\bf Theorem 14.1.}\, {\em  Let $ M$ and $M^D$ be given
by {\rm (14.1)} and {\rm (14.3)}  with ${\rm Ind}(M) \geq 1.$ Then the ranks of
$ G_1$---$G_4$ in  {\rm (14.3)} can be determined by the following
formulas
$$ \displaylines{
\hspace*{2cm}
r(G_1) =  r \left[ \begin{array}{cc} M^kJ_1M^k  & M^{k-1}V_1  \\
W_1M^{k-1} & 0  \end{array} \right]
- r( M^k ),  \hfill (14.4)
\cr
\hspace*{2cm}
r(G_2) =  r \left[ \begin{array}{cc} M^kJ_2M^k  & M^{k-1}V_2  \\
 W_1M^{k-1} & 0  \end{array} \right] - r( M^k ), \hfill(14.5)
\cr
\hspace*{2cm}
r(G_3) =  r \left[ \begin{array}{cc} M^kJ_3M^k  & M^{k-1}V_1  \\
W_2M^{k-1} & 0 \end{array} \right]
- r( M^k ),  \hfill(14.6)
\cr
\hspace*{2cm} 
 r(G_4) =  r \left[ \begin{array}{cc} M^kJ_4M^k  & M^{k-1}V_2  \\
 W_2M^{k-1} & 0  \end{array} \right] - r( M^k ), \hfill(14.7)
\cr}
$$ 
where $ V_1, \ V_2, \ W_1$ and $ W_2 $ are defined in {\rm (14.2)}$,$
and
$$
J_1 = \left[ \begin{array}{cc} -A & 0  \\ 0 & D \end{array} \right] ,
 \ \ \ \
J_2 = \left[ \begin{array}{cc} 0 & B  \\ -C &  0 \end{array} \right] ,
\ \ \ \
J_3 = \left[ \begin{array}{cr} 0 & -B  \\ C &  0 \end{array} \right] ,
 \ \ \ \
J_4 = \left[ \begin{array}{cr} A & 0  \\ 0 & -D \end{array} \right].
\eqno (14.8)
$$ } 
{\bf Proof.}\, We only show (14.4). In fact $ G_1$ in  (14.3)
can be written as
$$ 
G_1 = [\, I_m, \ 0 \,] M^D \left[ \begin{array}{c} I_m  \\ 0  \end{array}
 \right] = P_1M^DQ_1 =
P_1M^k( M^{2k+1})^{\dagger}M^kQ_1,
$$ 
where $ P_1 = [\, I_m, \ 0 \,]$ and $ Q_1 = \left[ \begin{array}{c} I_m
 \\ 0  \end{array} \right]$. Then  it follows by Eq.(1.6) and block 
elementary operations that
\begin{eqnarray*}
r(G_1) & = & r \left[ \begin{array}{cc} M^{2k+1}  & M^kQ_1   \\ P_1M^k & 0
  \end{array} \right] -
r( M^{2k+1} ) \\
& = & r \left[ \begin{array}{cc} M^{2k+1}-M^kQ_1P_1MM^k - M^kMQ_1P_1M^k &
 M^kQ_1 \\ P_1M^k & 0  \end{array} \right] - r(M^k) \\
& = & r \left[ \begin{array}{cc} M^{k}( \, M - Q_1P_1M - MQ_1P_1 \,)M^k &
 M^kQ_1 \\ P_1M^k & 0  \end{array} \right] - r(M^k) \\
& = & r \left[ \begin{array}{cc} M^kJ_1M^k  & M^{k-1}V_1  \\ W_1M^{k-1} & 0
 \end{array} \right]
- r( M^k ),
\end{eqnarray*}
which is exactly the equality (14.4). \qquad$\Box$ 

\medskip

The further simplification of  (14.4)---(14.7) is quite difficult,
because the powers of $ M$ occur in them.  However if $ M $ in  (14.1)
satisfies some additional conditions, the four rank equalities in  (14.4)---(14.7) can reduce to simpler forms.
We next present some of them. The first one is related to the well-known
result on the Drazin inverse of an upper triangular block
matrix (see Campbell and Meyer \cite{CM2}).
$$\displaylines{
\hspace*{4cm}
 \left[ \begin{array}{cc} A  & B  \\ 0 & N  \end{array} \right]^D
 = \left[ \begin{array}{cc} A^D  & X  \\ 0 & N^D  \end{array} \right],
 \hfill (14.9)
\cr
where  \hfill
\cr
\hspace*{0cm}
X = (A^D)^2 \left[ \, \sum_{i=0}^{l-1}(A^D)^iBN^i \,\right]
( \, I_n - N^DN \,) + ( \, I_m - AA^D \,)
\left[ \, \sum_{i=0}^{k-1}A^iB (N^D)^i \, \right] (N^D)^2 - A^DBN^D,
 \hfill (14.10)
\cr}
$$ 
and $ {\rm Ind}(A) = k, \  {\rm Ind}(N) = l.$ 

\medskip

\noindent {\bf Theorem 14.2.}\, {\em The rank of the submatrix $ X $ in {\rm (14.9)}
 is
$$ 
r(X) =  r \left[ \begin{array}{ccc} A^k & P_t(B)  & 0  \\ 0 & A^tBN^t &
P_t(B) \\  0 & 0 & N^l    \end{array} \right] -
r \left[ \begin{array}{cc} A^k & P_t(B)   \\ 0 & N^l  \end{array} \right],
\eqno (14.11)
$$ 
where $ t = {\rm Ind}\left[ \begin{array}{cc} A & B  \\ 0 & N  \end{array}
\right], \
 P_t(B) = \sum_{i=0}^{t-1}A^{t-i-1}BN^i$. In particular if $ A^kBN^l = 0,$
 then
$$ 
r(X) = r[\,  A^k, \  P_t(B) \,] +  r \left[ \begin{array}{c} P_t(B) \\ N^l
\end{array} \right] - r \left[ \begin{array}{cc} A^k & P_t(B)   \\ 0 & N^l
 \end{array} \right]. \eqno (14.12)
$$ 
In particular if $ R[P_t(B)] \subseteq R(A^k)$ and $ R[( P_t(B))^*] \subseteq R[(N^l)^*],$
 then
$ r(X) = r( A^kBC^l).$ }

\medskip
  
\noindent {\bf Proof.}\, It is easy to verify that 
$$ 
 M^t = \left[ \begin{array}{cc} A & B  \\ 0 & N \end{array} \right]^t
 = \left[ \begin{array}{cc} A^t & P_t(B)   \\ 0 & N^t  \end{array} \right],
   \ \ \ {\rm  and } \ \ \ P_{2t+1}(B) = A^{t+1} P_t(B) + P_t(B)N^{t+1} + A^tBN^t.
$$ 
Then applying (1.7) to $ X = [\, I_m, \ 0 \,] \left[ \begin{array}{cc}
 A^D & X  \\ 0 & N^D  \end{array} \right] \left[ \begin{array}{c} 0  \\ I_n
  \end{array} \right] = P_1 M^t(M^{2t+1})^{\dagger}M^tQ_2$, we find that
\begin{eqnarray*}
r(X) & = & r \left[ \begin{array}{cc} M^{2t+1}  & M^tQ_2   \\ P_1M^t & 0
\end{array} \right]
- r( M^{2k+1} ) \\ 
& = & r \left[ \begin{array}{ccc} A^{2t+1} & P_{2t+1}(B) & P_t(B) \\ 0
 & N^{2t+1} &  N^t  \\
A^t & P_t(B) & 0 \end{array} \right] - r( M^{k} ) \\
& = & r \left[ \begin{array}{ccc} 0 & A^tBN^t & P_t(B) \\ 0  & 0 &  N^t
 \\
A^t & P_t(B) & 0 \end{array} \right] - r(M^{k}) \\
&  =  & r \left[ \begin{array}{ccc} A^k & P_t(B)  & 0  \\ 0 & A^tBD^t & P_t(B)
  \\  0 & 0 & N^l    \end{array} \right] -
   r \left[ \begin{array}{cc} A^k & P_t(B)   \\ 0 & N^l  \end{array} \right].
\end{eqnarray*}
Thus we have  the desired results.  \qquad $\Box$ 

\medskip

\noindent {\bf Theorem 14.3.}\, {\em Let $ M$ be given by  {\rm (14.1)}
with ${\rm Ind}(M)= 1.$ Then the ranks of
 $G_1$---$G_4$ in the group inverse of $ M $ in  {\rm (14.3)} can be
 expressed as
$$ 
\displaylines{
\hspace*{2cm} 
r(G_1) =  r \left[ \begin{array}{cc} V_2DW_2  & V_1  \\ W_1 & 0
 \end{array} \right] - r(M), \qquad
r(G_2) =  r \left[ \begin{array}{cc} V_1BW_2  & V_2  \\ W_1 & 0
 \end{array} \right] - r(M), \hfill (14.13)
\cr
\hspace*{2cm} 
r(G_3) =  r \left[ \begin{array}{cc} V_2CW_1  & V_1  \\ W_2 & 0
  \end{array} \right] - r(M), \qquad
r(G_4) =  r \left[ \begin{array}{cc} V_1AW_1  & V_2  \\ W_2 & 0
\end{array} \right] - r(M), \hfill (14.14)
\cr}
$$  
where $ V_1, \ V_2, \ W_1$ and $ W_2 $ are defined in  {\rm (14.2)}. } 

\medskip 

\noindent {\bf Proof.}\,  Note that $ M^{\#} = M(M^3)^{\dagger}M$ when
 ${\rm Ind}(M)= 1$. Thus $ G_1$ in  (14.13) can be written as 
$ G_1= W_1(M^3)^{\dagger}V_1$.  In that case it follows by (1.7) that 
\begin{eqnarray*}
r(G_1) &= & r \left[ \begin{array}{cc} M^3  & V_1  \\ W_1 & 0  \end{array}
\right] -r(M^3) \\
 & =& r \left[ \begin{array}{cc} [\, 0, \ V_2 \, ]M \left[ \begin{array}{cc}
 0 \\ W_2  \end{array} \right]  & V_1  \\ W_1 & 0 \end{array} \right] -r(M)
  =  r \left[ \begin{array}{cc} V_2DW_2  & V_1  \\ W_1 & 0  \end{array}
   \right] - r(M).
\end{eqnarray*}
In the same manner we can show the other three in  (14.13) and (14.14).
\qquad $\Box$

\medskip

\noindent {\bf Corollary 14.4.}\, {\em  Let $ M$ be given by  {\rm (14.1)}
with ${\rm Ind}(M)= 1.$

{\rm (a)}\, If $ M $ satisfies the rank additivity condition
$$\displaylines{
\hspace*{3cm} 
r(M) = r(V_1) + r(V_2) = r(W_1) + r(W_2), \hfill
\cr}
$$
then the ranks of  $G_1$---$G_4$ in the group inverse of $ M $ in
  {\rm (14.3)} can be expressed as
$$\displaylines{
\hspace*{3cm} 
r(G_1) = r(V_1) + r( W_1) + r( V_2DW_2) - r(M), \hfill
\cr
\hspace*{3cm} 
r(G_2) = r(V_2) + r( W_1) + r( V_1BW_2) - r(M),  \hfill
\cr
\hspace*{3cm} 
r(G_3) = r(V_1) + r( W_2) + r( V_2CW_1) - r(M), \hfill
\cr
\hspace*{3cm} 
r(G_4) = r(V_2) + r( W_2) + r( V_1AW_1) - r(M). \hfill
\cr}
$$

{\rm (b)}\, If $ M $ satisfies the rank additivity condition
$$\displaylines{
\hspace*{2cm} 
r(M) = r(A) + r(B) + r(C) + r(D), \hfill
\cr}
$$
then the ranks of $G_1$---$G_4$ in the group inverse of $ M $ in
 {\rm (14.3)} satisfy 
$$\displaylines{
\hspace*{2cm} 
r(G_1) = r(A) - r(D) + r( V_2DW_2),  \qquad r(G_2) = r(B) - r(C) +
r(V_1BW_2), \hfill
\cr
\hspace*{2cm} 
r(G_3) = r(C) - r(B) + r( V_2CW_1),  \qquad r(G_2) = r(D) - r(A) +
 r(V_1AW_1). \hfill
\cr}
$$
where $ V_1, \ V_2, \ W_1$ and $ W_2 $ are defined in {\rm (14.2)}. } 

\medskip

In addition, we have some inequalities on ranks of submatrices
in the group inverse of a block matrix.

\medskip

\noindent {\bf  Corollary 14.5.}\, {\em  Let $ M$ be given by  {\rm (14.1)}
with
${\rm Ind}(M)= 1.$  Then the ranks of the matrices $ G_1$---$G_4$ in {\rm (14.3)} satisfy the following rank inequalities

{\rm (a)} \ $r(G_1) \geq  r(V_1) + r( W_1) - r(M).$ 

{\rm (b)} \ $r(G_1) \leq  r(V_1) + r( W_1) + r(D) - r(M).$  
 
{\rm (c)} \ $r(G_2) \geq  r(V_2) + r( W_1) - r(M).$ 

{\rm (d)} \ $r(G_2) \leq  r(V_2) + r( W_1) + r(B) - r(M).$  
 
{\rm (e)} \ $r(G_3) \geq  r(V_1) + r( W_2) - r(M).$ 

{\rm (f)} \ $r(G_3) \leq  r(V_1) + r( W_2) + r(C) - r(M).$  
 
{\rm (g)} \ $r(G_4) \geq  r(V_2) + r( W_2) - r(M).$ 

{\rm (h)} \ $r(G_4) \leq  r(V_2) + r( W_2) + r(A) - r(M).$ } 

\medskip

\noindent {\bf Proof.} \, Follows from  (14.13) and (14.14). \qquad $ \Box$.

\markboth{YONGGE  TIAN }
{15. REVERSE ORDER LAWS FOR DRAZIN INVERSES OF PRODUCTS OF MATRICES  }

\chapter{Reverse order laws for Drazin inverses}

\noindent In this chapter we consider reverse order laws for  Drazin inverses of products of matrices.
 We will give necessary and sufficient conditions for $( ABC)^D = C^DB^DA^D$ to hold and 
then present some of its consequences.

\medskip

\noindent {\bf Lemma 15.1.}\, {\em Let $ A, \, X \in {\cal C}^{m \times m}$
with
${\rm Ind}(A) = k.$ Then $ X = A^D$ if and only if
$$ 
A^{k+1}X = A^k,  \qquad  XA^{k+1} = A^k, \qquad and \qquad r(X) = r(A^k). \eqno (15.1)
$$}
{\bf Proof.}\, Follows from the definition of the Drazin inverse of a matrix. \qquad $ \Box$ 

\medskip

\noindent {\bf Lemma 15.2.}\, {\em  Let $ A, \, B, \, C \in
{\cal C}^{m \times m}$ with ${\rm Ind}(A) = k_1, \
{\rm Ind}(B) = k_2 $ and $ {\rm Ind}(C) = k_3$. Then the product
$ C^DB^DA^D$ of the Drazin inverses of $ A, \ B,$ and $C$ can be expressed in the form 
$$\displaylines{
\hspace*{1.5cm}
 C^DB^DA^D = [\, C^{k_3}, \ 0, \ 0  \, ]  \left[ \begin{array}{ccc} 0 & 0
 & A^{2k_1+1} \\ 0 & B^{2k_2+1} & B^{k_2}A^{k_1}  \\ C^{2k_3+1}  &
 C^{k_3}B^{k_2} & 0 \end{array} \right]^{\dagger}
\left[ \begin{array}{c} A^{k_1} \\ 0  \\ 0  \end{array} \right]
 := PN^{\dagger}Q, \hfill (15.2)
\cr}
$$ 
where $ P,  \ N$ and $ Q$ satisfy the three properties
$$\displaylines{
\hspace*{1.0cm} 
R(Q) \subseteq R(N), \ \ \ \  R(P^*) \subseteq R(N^*), \ \ \ \ \ 
r(N) = r(A^{k_1}) + r(B^{k_2}) + r(C^{k_3}). \hfill (15.3)
\cr}
$$ }   
{\bf Proof.}\, It is easy to verify that the $ 3 \times 3$ block
matrix $ N $ in  (15.2) satisfies the conditions in Lemma 8.8.
Hence  it follows by (8.8) that
$$ \displaylines{
\hspace*{1.0cm}
 N^{\dagger} = \left[ \begin{array}{cccc}
& (C^{2k_3+1})^{\dagger}C^{k_3}B^{k_2}(B^{2k_2+1})^{\dagger}B^{k_2}
 A^{k_1}(A^{2k_1+1})^{\dagger}A^{k_1} & * & *  \\ & *  & * & 0 \\
  & * & 0 & 0 \end{array} \right].
\hfill (15.4)
\cr}
$$  
Thus we have (15.2). The three properties in  (15.3) follows from the
structure of $ N $.   \qquad $ \Box$

\medskip

The main results of the chapter are the following two. 

\medskip

\noindent {\bf Theorem 15.3.}\, {\em  Let $ A, \, B, \, C
\in {\cal C}^{m \times m}$ with
${\rm Ind}(A) = k_1, \ {\rm Ind}(B) = k_2 $ and $ {\rm Ind}(C) = k_3,$ and
 denote $ M = ABC$ with ${\rm Ind}(M) = t$. Then the reverse order law
 $ ( ABC )^D = C^DB^DA^D$ holds if and only if $ A, \ B$ and $ C $ satisfy
 the three rank equalities
$$
\displaylines{
\hspace*{1.0cm}
 r \left[ \begin{array}{cccc}  0  &  0  & A^{2k_1+1}  &  A^{k_1} \\
0 & B^{2k_2+1} & B^{k_2}A^{k_1} &  0  \\ C^{2k_3+1}  & C^{k_3}B^{k_2} & 0
 & 0 \\  M^{t+1}C^{k_3} & & 0 & M^t  \end{array} \right]
=  r(A^{k_1}) + r(B^{k_2}) + r(C^{k_3}), \hfill (15.5)
\cr
\hspace*{1.0cm}
r \left[ \begin{array}{cccc}  0  & 0 & A^{2k_1+1}  &  A^{k_1} M^{t+1} \\ 0
 & B^{2k_2+1} & B^{k_2}A^{k_1} &  0  \\ C^{2k_3+1}  & C^{k_3}B^{k_2} & 0 & 0
   \\  C^{k_3} & & 0 & M^t  \end{array} \right]
=  r(A^{k_1}) + r(B^{k_2}) + r(C^{k_3}), \hfill (15.6)
\cr
\hspace*{1.0cm}
r \left[ \begin{array}{cc} B^{2k_2+1} & B^{k_2}A^{k_1} \\  C^{k_3}B^{k_2} &
0  \end{array} \right]
=  r(B^{k_2}) + r(M^t). \hfill (15.7)
\cr }
$$ } 
{\bf Proof.}\, Let $ X = C^DB^DA^D$. Then by definition of the
Drazin
inverse, $X = M^D$ if and only if $ M^{t+1}X = M^t$, $ XM^{t+1} = M^t $ and
 $r(X) = r(M^t)$, which are equivalent to
$$ \displaylines{
\hspace*{1.5cm}
 r( \, M^k - M^{k+1}X \,) = 0, \qquad r( \, M^k - XM^{k+1}\,)
  = 0 \ \  {\rm and}  \ \  r(X) = r(M^t).
\hfill (15.8)
\cr}
$$ 
Replacing $ X $ in (15.8) by $ X = PN^{\dagger}Q $ in  (15.2) and
applying (1.7) them, we find that
$$\displaylines{
\hspace*{1.5cm}
 r( \, M^t - M^{t+1}X \,) =  r( \, M^t - M^{t+1}PN^{\dagger}Q \,) =
 r \left[ \begin{array}{cc} N & Q \\ M^{t+1}P & M^t  \end{array} \right]
 - r(N), \hfill
\cr
\hspace*{1.5cm}
 r( \, M^t - XM^{t+1} \,) =  r( \, M^t - PN^{\dagger}QM^{t+1} \,) =
 r \left[ \begin{array}{cc} N & Q M^{t+1} \\ P & M^t  \end{array} \right] -
 r(N),\hfill
\cr
\hspace*{1.5cm}
 r(X) =  r( PN^{\dagger}Q) 
=  r \left[ \begin{array}{cc} N & Q  \\ P & 0 \end{array} \right] - r(N).\hfill
\cr}
$$ 
Putting them in  (15.8), we obtain  (15.5)---(15.7). \qquad $ \Box$ 

\medskip

\noindent {\bf Theorem 15.4.}\, {\em  Let $ A, \, B, \, C \in {\cal C}^{m \times m}$ with
${\rm Ind}(A) = k_1, \
{\rm Ind}(B) = k_2 $ and $ {\rm Ind}(C) = k_3,$ and let $ M = ABC$
with ${\rm Ind}(M) = t$. Then the reverse order law
$ ( ABC )^D = C^DB^DA^D$ holds if and only if $ A, \ B$ and $ C $ satisfy
the following  rank equality
$$
r \left[ \begin{array}{ccccc}  0  &  0  & A^{2k_1+1}  &  A^{k_1}  & 0 \\ 0
 & B^{2k_2+1} & B^{k_2}A^{k_1} &  0  & 0  \\ C^{2k_3+1}  & C^{k_3}B^{k_2} &
  0 & 0 &  0  \\ C^{k_3}  & 0 & 0 & 0 & M^t  \\ 0  & 0 & 0 & M^t
 & M^{2t+1}  \end{array} \right]
=  r(A^{k_1}) + r(B^{k_2}) + r(C^{k_3}) - r(M^t). \eqno (15.9)
$$ }
{\bf Proof.}\, Applying (2.3) to $ ( ABC )^D - C^DB^DA^D
= M^t(M^{2t+1})^{\dagger}M^t - PN^{\dagger}Q$,
we find that
\begin{eqnarray*}
r[ \, ( ABC )^D - C^DB^DA^D \, ] & = & r[ \, PN^{\dagger}Q  -
 M^t(M^{2t+1})^{\dagger}M^t \,] \\
& = & r \left[ \begin{array}{ccc} N & 0 & Q \\  0  & -M^{2t+1}  & M^t \\ P &
M^t & 0 \end{array} \right] - r(N) - r(M^t) \\
& = & r \left[ \begin{array}{ccc} N & 0 & Q \\  0  & 0  & M^t \\ P & M^t &
-M^{2t+1}
 \end{array} \right] - r(N) - r(M^t).
\end{eqnarray*} 
Thus  (15.9) follows  by putting $ P, \ N $ and $Q $ in it. \qquad $
\Box$ 

\medskip

We next give some particular cases of the above two theorems. 

\medskip

\noindent {\bf Corollary 15.5.}\, {\em  Let $ A, \, B, \, C \in
{\cal C}^{m \times m}$ with ${\rm Ind}(B) = k$ and ${\rm Ind}(ABC) = t,$
where $A$ and $C$ are nonsingular. Then

{\rm (a)} \ $  r[ \, (ABC)^D - C^{-1}B^DA^{-1} \, ]
= r\left[ \begin{array}{c} B^k \\ (ABC)^tA \end{array} \right] +
 r[ \,  B^k, \ C(ABC)^t \, ] - r(B^k) -r[(ABC)^t].$

{\rm (b)} \ $ (ABC)^D = C^{-1}B^DA^{-1}  \Leftrightarrow
  R[C(ABC)^t] = R(B^k) \ \ and \ \
 R\{[(ABC)^tA]^*\} = R[(B^k)^*] .$ } 

\medskip

\noindent {\bf Proof.}\, It is easy to verify that both
 $(ABC)^D $ and $ C^{-1}B^DA^{-1}$ are outer inverses of $ABC$. Thus
it follows from  (5.1) that
$$
\displaylines{
\hspace*{1cm}
 r[ \, (ABC)^D - C^{-1}B^DA^{-1} \, ]\hfill
\cr
\hspace*{1cm}  =  r\left[ \begin{array}{c} (ABC)^D \\  
C^{-1}B^DA^{-1}  \end{array} \right] +
 r[ \, (ABC)^D, \  C^{-1}B^DA^{-1} \, ] - r[(ABC)^D] - r(B^D) \hfill
\cr
\hspace*{1cm}
 =   r\left[ \begin{array}{c} (ABC)^tA \\  B^k  \end{array} \right] +
 r[ \, C(ABC)^t, \ B^k \, ] - r[(ABC)^t] - r(B^k), \hfill
\cr}
$$
as required for Part (a). Notice that 
$$
 r\left[ \begin{array}{c} B^k \\ (ABC)^tA \end{array} \right] \geq  r(B^k),  \ \ \ 
r\left[ \begin{array}{c} B^k \\ (ABC)^tA \end{array} \right] \geq  r[(ABC)^t], 
$$
and
$$
r[ \,  B^l, \ C(ABC)^t \, ] \geq  r(B^k), \ \ \ r[ \,  B^l, \ C(ABC)^t \, ] \geq  r[(ABC)^t].
$$
Then Part (b) follows from Part (a). \qquad $\Box$ 

\medskip

\noindent {\bf Corollary 15.6.}\, {\em  Let $ A, \, B, \, C
\in {\cal C}^{m \times m}$ with ${\rm Ind}(A) = k_1, \,
{\rm Ind}(B) = k_2 $ and $ {\rm Ind}(C) = k_3,$ and let $ M = ABC$ with
${\rm Ind}(M) = t.$  Moreover suppose that
 $$ 
AB= BA, \qquad AC= CA, \qquad BC = CB.  \eqno (15.10)
$$ 
Then the reverse order law $ ( ABC )^D = C^DB^DA^D$ holds if and only if
$ A, \ B$ and $ C $ satisfy {\rm (15.7)}. } 

\medskip

\noindent {\bf Proof.}\, It is not difficult to verify that under (15.10),
the two rank equalities in  (15.5) and (15.6)
become two identities. Thus,  (15.7) becomes a necessary and
sufficient condition for $ ( ABC )^D = C^DB^DA^D$  to hold. \qquad
$ \Box$ 

\medskip

\noindent {\bf  Corollary 15.7.}\, {\em  Let $ A, \, B \in {\cal C}^{m \times
m}$ with ${\rm Ind}(A) = k, \, {\rm Ind}(B) = l $ and ${\rm Ind}(AB) = t$.
Then the following three are equivalent$:$

{\rm (a)} \ $ ( AB)^D = B^DA^D.$

{\rm (b)}  \ $ r \left[ \begin{array}{cccc} 0 & A^{2k+1}  &  A^{k}  & 0 \\
 B^{2l+1} & B^{l}A^{k} &  0  & 0  \\  B^{l}  & 0 & 0  & (AB)^t  \\ 0  & 0
   & (AB)^t & (AB)^{2t+1}  \end{array} \right]
=  r(A^{k}) + r(B^{l})  - r[(AB)^t].$

{\rm (c)}\, The following three rank equalities are all satisfied
$$
\displaylines{
\hspace*{1.5cm}
r[(AB)^t] = r(B^lA^k), \hfill
\cr
\hspace*{1.5cm}
 r \left[ \begin{array}{ccc}  0  & A^{2k+1}  &  A^{k} \\
B^{2l+1} & B^{l}A^{k} &  0  \\ (AB)^{t+1}B^{l}  & 0 & -(AB)^t  \end{array}
\right] =
r(A^{k}) + r(B^{l}), \hfill
\cr
\hspace*{1.5cm}
 r \left[ \begin{array}{ccc} 0  & A^{2k+1}  &  A^{k}(AB)^{t+1} \\
B^{2l+1} & B^{l}A^{k} &  0  \\ B^{l} & 0 & -(AB)^t  \end{array} \right] =  r(A^{k}) + r(B^{l}).\hfill
\cr}
$$ }  
{\bf Proof.}\, Letting $ C = I_m$ in  (15.9) results in
Part (b), and letting $ B = I_m $ and replacing $ C $
 by $ B $  in Theorem 15.4 result in Part (c). \qquad $ \Box$ 

\markboth{YONGGE  TIAN }
{16. RANK EQUALITIES FOR WEIGHTED MOORE-PENROSE INVERSES OF MATRICES  }

\chapter{Ranks equalities for weighted Moore-Penrose inverses}

\noindent The weighted Moore-Penrose inverse of a matrix $ A \in {\cal C}^{m \times n} $ with
 respect to two positive definite  matrices $ M \in {\cal C}^{m \times m} $ and 
$N \in {\cal C}^{n \times n}$ is defined to be the unique solution of the following four matrix 
equations
$$ 
AXA = A, \qquad  XAX = X, \qquad  (MAX)^* = MAX, \qquad (NXA)^* = NXA,
 \eqno (16.1)
$$ 
and this $ X $ is often denoted by $ X = A^{\dagger}_{M,N}$. In particular, when
 $ M = I_m$ and $ N = I_n$,
 $A^{\dagger}_{M,N}$ is the standard Moore-Penrose inverse $ A^{\dagger}$
 of $ A$. As is well known (see, e.g., Rao and Mitra \cite{RM}), the weighted Moore-Penrose inverse
$A^{\dagger}_{M,N}$ of $A$ can be written as a  matrix
 expressions involving a standard Moore-Penrose inverse as follows 
$$
A^{\dagger}_{M,N} =
N^{-\frac{1}{2}}(M^{\frac{1}{2}}AN^{-\frac{1}{2}})^{\dagger}M^{\frac{1}{2}},
\eqno (16.2)
$$ 
where $ M^{\frac{1}{2}}$ and $ N^{\frac{1}{2}}$ are the positive definite
square roots of $M$ and $ N $, respectively. According to (16.2), it is easy to verify that
$$
R( A^{\dagger}_{M,N}) = R(N^{-1}A^*), \ \  {\rm and}  \ \
R[( A^{\dagger}_{M,N})^*] = R(MA). \eqno (16.2)
$$
Based on these basic facts and the rank formulas in Chapters 2---5, we now can establish
 various rank equalities related to weighted Moore-Penrose inverses of 
matrices, and the consider their various consequences.

\medskip

\noindent {\bf Theorem 16.1.}\, {\em Let $ A \in {\cal C}^{m \times n}$ be
given$,$
$ M \in {\cal C}^{m \times m}$ and $ N \in {\cal C}^{n \times n}$ be two
positive definite matrices. Then

{\rm (a)} \ $ r(  \, A^{\dagger} -  A^{\dagger}_{M,N} \, )
= r \left[ \begin{array}{c} A \\  AN  \end{array} \right] +
r[ \, A, \ MA \, ] - 2r ( A ).$

{\rm (b)} \  $ r(  \, A^{\dagger} -  A^{\dagger}_{M,I} \, )
= r[ \, A, \ MA \, ] - r ( A ). $

{\rm (c)} \ $  r(  \, A^{\dagger} -  A^{\dagger}_{I,N} \, )
= r \left[ \begin{array}{c} A \\  AN  \end{array} \right] - r ( A ).$

{\rm (d)} \ $ A^{\dagger}_{M,N} = A^{\dagger}  \Leftrightarrow
 R(MA) = R(A)$ and $ R[(AN)^*] = R(A^*).$ }

\medskip

\noindent {\bf Proof.}\, Note that $ A^{\dagger}$ and $ A^{\dagger}_{M,N}$
are outer inverses of $ A $. Thus it follows from (5.1) that
\begin{eqnarray*}
 r( \, A^{\dagger} -  A^{\dagger}_{M,N} \,) 
& = & r \left[ \begin{array}{c} A^{\dagger} \\ A^{\dagger}_{M,N} \end{array}
\right] + r[\, A^{\dagger}, \  A^{\dagger}_{M,N} \,] - r(A^{\dagger}) -
 r( A^{\dagger}_{M,N}) \\
 & = & r \left[ \begin{array}{c} A^* \\ (MA)^* \end{array} \right]
 + r[\, A^*, \ N^{-1}A^* \,] - 2r(A) \\
& = &r \left[ \begin{array}{c} A \\  AN \end{array} \right] +
 r[ \, A, \ MA \,] - 2r(A).
\end{eqnarray*}
Parts (a)---(c) follow immediately from it.  \qquad $ \Box$ 

\medskip

\noindent {\bf Theorem 16.2.}\, {\em Let $ A \in {\cal C}^{m \times n}$ be
given$,$
 $ M \in {\cal C}^{m \times m}$ and $ N \in {\cal C}^{n \times n}$ be two
positive definite matrices. Then

{\rm (a)} \ $ r( \,  AA^{\dagger}_{M,N} - AA^{\dagger} \,)
= r[\, A, \ MA \,] - r (A).$

{\rm (b)} \  $r(\, A^{\dagger}_{M,N}A -   A^{\dagger}A \, ) =
r \left[ \begin{array}{c} A \\  AN  \end{array} \right] - r ( A ).$

{\rm (c)} \ $ AA^{\dagger}_{M,N} = AA^{\dagger}  \Leftrightarrow
 R(MA) = R(A).$

{\rm (d)} \ $A^{\dagger}_{M,N}A = A^{\dagger}A   \Leftrightarrow
 R[(AN)^*] = R(A^*).$ } 

\medskip

\noindent {\bf Proof.}\, Note that both $ AA^{\dagger}$  and
$AA^{\dagger}_{M,N}$ are idempotent. It follows from  (3.1)
that
\begin{eqnarray*}
 r( \, AA^{\dagger} -  AA^{\dagger}_{M,N} \,) 
& = & r \left[ \begin{array}{c} AA^{\dagger} \\ AA^{\dagger}_{M,N}
 \end{array} \right] + r[\, AA^{\dagger}, \  AA^{\dagger}_{M,N} \,]
 - r(AA^{\dagger}) -  r(AA^{\dagger}_{M,N})  \\
 & = & r \left[ \begin{array}{c} A^* \\ (MA)^* \end{array} \right]
 + r[\, A, \ A \,] - 2r(A)  \\
& =  & r[ \, A, \ MA \,] - r(A),
\end{eqnarray*}
as required for Part (a). Similarly we can show Part (b).  \qquad $ \Box$ 

\medskip

\noindent {\bf Theorem 16.3.}\, {\em Let $ A \in {\cal C}^{m \times m}$ be
given$,$ and $ M, \, N  \in {\cal C}^{m \times m}$
be  two positive definite matrices. Then 

{\rm (a)} \ $ r( \,  AA^{\dagger}_{M,N} - A^{\dagger}_{M,N}A \,) =
 r[ \, A^*, \ MA \, ] +
 r[ \, A^*, \ NA \,]  - 2r(A).$

{\rm (b)} \ $ AA^{\dagger}_{M,N} = A^{\dagger}_{M,N}A
\Leftrightarrow  R(MA) = R(NA) = R(A^*)  \Leftrightarrow
 $ both MA and NA are EP. } 

\medskip

\noindent {\bf Proof.}\, Note that both $ AA^{\dagger}$  and
$AA^{\dagger}_{M,N}$ are idempotent. It follows by  (3.1)
that
\begin{eqnarray*}
 r( \, AA^{\dagger}_{M,N} -  A^{\dagger}_{M,N}A \,) 
& = & r \left[ \begin{array}{c} AA^{\dagger}_{M,N} \\
A^{\dagger}_{M,N}A \end{array} \right] +
r[\, AA^{\dagger}_{M,N}, \  A^{\dagger}_{M,N}A \,]
 - r(AA^{\dagger}_{M,N}) -  r(A^{\dagger}_{M,N}A)  \\
& = & r \left[ \begin{array}{c} A^{\dagger}_{M,N} \\
A\end{array} \right] + r[\, A, \  A^{\dagger}_{M,N} \,] - 2r(A) \\
& = & r \left[ \begin{array}{c} (MA)^* \\
A \end{array} \right] + r[\, A, \ N^{-1}A^* \,] - 2r(A),
\end{eqnarray*}
as required for Part (a). Part(b) follows immediately from Part (a).  \qquad $ \Box$ 

\medskip

Based on the result in Theorem 16.3(b), we can extend the concept of EP matrix to weighted case: 
A square matrix $ A $ is said to be {\em weighted EP} if both $ MA$ and
$NA$ are EP, where both $ M $ and $ N $ are two positive definite matrices. It is expected that
weighted EP matrix would have some nice properties. But we do not intend to go 
further along this direction 
in the thesis.

\medskip

\noindent {\bf Theorem 16.4.}\, {\em Let $ A \in {\cal C}^{m \times m}$ be
given$,$ and $ M, \, N  \in {\cal C}^{m \times m}$
be  two positive definite matrices. Then 

{\rm (a)} \ $ r( \,  AA^{\dagger}_{M,N} -\overline{ A^{\dagger}_{M,N}A} \,) =
 r[ \, A^T, \ MA \, ] + r[ \, A^T, \ N^TA \,]  - 2r(A).$

{\rm (b)} \  $ AA^{\dagger}_{M,N} =\overline{ A^{\dagger}_{M,N}A}
\Leftrightarrow \ R(MA) = R(N^TA) = R(A^T)  \Leftrightarrow
 $ both $MA$ and $N^TA$ are EP. } 

\medskip

\noindent {\bf Proof.}\, Follows from (3.1) by  noting that both $ AA^{\dagger}_{M,N}$  and
$ \overline{ AA^{\dagger}_{M,N}}$ are idempotent. \qquad $ \Box$ 

\medskip

\noindent {\bf Theorem 16.5.}\, {\em Let $ A \in {\cal C}^{m \times m}$ be
given
with {\rm Ind}$(A) = 1,$ and  $ M, \, N  \in {\cal C}^{m \times m}$ be two
positive definite matrices. Then

{\rm (a)}  \ $r( \, A^{\dagger}_{M,N} - A^{\#} \,)
= r[\, A^*, \ MA \,] + r[\, A^*, \ NA \,] - 2r(A).$

{\rm (b)} \ $A^{\dagger}_{M,N} = A^{\#}  \Leftrightarrow  R(MA) =
  R(NA) = R(A^*),$ i.e., $ A $ is weighted EP. }

\medskip

\noindent {\bf Proof.}\, Note that both $ A^{\dagger}$  and
$A^{\#}$ are outer inverses of $A$. It follows by  (5.1)
that
\begin{eqnarray*}
 r( \, A^{\dagger}_{M,N} -  A^{\#} \,) 
& = & r \left[ \begin{array}{c} A^{\dagger}_{M,N} \\ A^{\#} \end{array}
\right] + r[\, A^{\dagger}_{M,N}, \  A^{\#} \,]
 - r(A^{\dagger}_{M,N}) -  r(A^{\#})  \\
& = & r \left[ \begin{array}{c} (MA)^* \\ A \end{array} \right]
+ r[\, N^{-1}A^*, \ A \,] - 2r(A) \\
& = & r[\, A^*, \ MA \,] + r[\, A^*, \ NA \,] - 2r(A),
\end{eqnarray*}
as required for Part (a).  \qquad $ \Box$

\medskip

\noindent {\bf Theorem 16.6.}\, {\em Let $ A \in {\cal C}^{m \times m}$ be
given
with ${\rm Ind}(A) = 1,$ and $ M, \, N  \in {\cal C}^{m \times m}$ be two
positive definite matrices. Then

{\rm (a)} \ $ r( \, AA^{\dagger}_{M,N} - AA^{\#} \,) = r[\, A^*, \ MA \,]
 - r (A).$

{\rm (b)} \ $r( \, A^{\dagger}_{M,N}A - A^{\#}A \,) = r[\, A^*, \ NA \,] -
r (A).$

{\rm (c)}  \ $r( \, A^{\dagger}_{M,N} - A^{\#} \,) = r( \, AA^{\dagger}_{M,N}
- AA^{\#} \,) + r( \, A^{\dagger}_{M,N}A - A^{\#}A \,).$ \\
In particular$,$

{\rm (d)}\  $ AA^{\dagger}_{M,N} = AA^{\#}  \Leftrightarrow  R(MA)
= R(A^*),$ i.e.$,$ $ MA $ is EP.

{\rm (e)} \  $A^{\dagger}_{M,N}A = A^{\#}A   \Leftrightarrow
 R(NA) = R(A^*),$ i.e.$,$ $ NA $ is EP.

{\rm (f)} \ $A^{\dagger}_{M,N} = A^{\#} \Leftrightarrow
 AA^{\dagger}_{M,N} = AA^{\#}$  \ and \ $A^{\dagger}_{M,N}A = A^{\#}A$. }

\medskip

\noindent {\bf Proof.}\, Note that both $ AA^{\dagger}$  and
$AA^{\#}$ are idempotent. It follows from (5.1)
that
\begin{eqnarray*}
 r( \, AA^{\dagger}_{M,N} -  AA^{\#} \,) 
& = & r \left[ \begin{array}{c} AA^{\dagger}_{M,N} \\ AA^{\#} \end{array}
\right] + r[\, AA^{\dagger}_{M,N}, \  AA^{\#} \,]
 - r(AA^{\dagger}_{M,N}) -  r(AA^{\#})  \\
& = & r \left[ \begin{array}{c} A^{\dagger}_{M,N} \\ A^{\#} \end{array}
\right] + r[\, A, \  A \,] - 2r(A)  \\
& = & r \left[ \begin{array}{c} (MA)^* \\ A^* \end{array} \right]  = r[\, A^*, \ MA \,] - r (A),
\end{eqnarray*}
as required for Part (a).  \qquad $ \Box$ 

\medskip

\noindent {\bf Theorem 16.7.}\, {\em Let $ A \in {\cal C}^{m \times m}$ be
given with
${\rm Ind}(A) = k,$  and $ M, \, N  \in {\cal C}^{m \times m}$ be two positive
definite matrices. Then

{\rm (a)} \ $r( \, A^{\dagger}_{M,N} - A^{D} \,) =
 r \left[ \begin{array}{c} A^kM^{-1} \\  A^* \end{array} \right] +
 r[\, NA^k, \ A^* \,] - r(A) - r (A^k). $

{\rm (b)} \ $r( \, A^{\dagger}_{M,N} - A^{D} \,) = r(A^{\dagger}_{M,N}) -
 r(A^{D})  \Leftrightarrow  R(NA^k) \subseteq r(A^*) $ and
 $ R[(A^kM^{-1})^*] \subseteq r(A).$  }

\medskip

\noindent {\bf Proof.}\, Note that both $ A^{\dagger}_{M,N}$  and
$A^D$ are outer inverses of $A$. It follows by  (5.1)
that
\begin{eqnarray*}
 r( \, A^{\dagger}_{M,N} -  A^D \,) 
& = & r \left[ \begin{array}{c} A^{\dagger}_{M,N} \\ A^D \end{array}
\right] + r[\, A^{\dagger}_{M,N}, \  A^D \,] - r(A^{\dagger}_{M,N}) -
r(A^D)  \\
& = & r \left[ \begin{array}{c} (MA)^* \\ A^k \end{array} \right]
+ r[\, N^{-1}A^*, \ A^k \,] - 2r(A) \\
& = & r \left[ \begin{array}{c} A^* \\ A^kM^{-1} \end{array} \right]
+ r[\, A^*, \ NA^k \,] - 2r(A),
\end{eqnarray*}
as required for Part (a).  \qquad $ \Box$ 

\medskip

\noindent {\bf Theorem 16.8.}\, {\em Let $ A \in {\cal C}^{m \times m}$ be
given
with ${\rm Ind}(A) = k,$ and $ M, \, N  \in {\cal C}^{m \times m}$ be two
positive definite matrices. Then

{\rm (a)}  \ $r( \, AA^{\dagger}_{M,N} - AA^{D} \,) =
r \left[ \begin{array}{c} A^kM^{-1} \\  A^* \end{array} \right] -
r (A^k). $

{\rm (b)} \ $r( \, A^{\dagger}_{M,N}A - A^{D}A \,) =
r[\, NA^k, \ A^* \,] - r (A^k).$

{\rm (c)} \ $r( \, A^{\dagger}_{M,N} - A^{D} \,)
=r( \, AA^{\dagger}_{M,N} - AA^{D} \,) + ( \, A^{\dagger}_{M,N}A - A^{D}A \,)
 +  r (A^k) - r(A). $ } 

\medskip

\noindent {\bf Proof.}\, Note that both $ AA^{\dagger}_{M,N}$  and
$AA^D$ are idempotent. It follows from  (5.1)
that
\begin{eqnarray*}
 r( \, AA^{\dagger}_{M,N} -  AA^D \,) 
& = & r \left[ \begin{array}{c} AA^{\dagger}_{M,N} \\ AA^D \end{array}
\right] + r[\, AA^{\dagger}_{M,N}, \  AA^D \,] -
 r(AA^{\dagger}_{M,N}) - r(AA^D)  \\
& = & r \left[ \begin{array}{c} A^{\dagger}_{M,N} \\ A^D \end{array}
\right] + r[\, A, \  A^D \,] -
 r(A) - r(A^k)  \\
& = & r \left[ \begin{array}{c} (MA)^* \\ A^k \end{array} \right] - r(A)  
=  r \left[ \begin{array}{c}  A^* \\ A^kM^{-1} \end{array} \right] -
r(A),
\end{eqnarray*}
as required for Part (a). Similarly we can show Part (b). Combining
Theorem 16.6(a) and Theorem 16.7(a) yields Part (c).  \qquad $ \Box$ 

\medskip

\noindent {\bf Theorem 16.9.}\, {\em Let $ A \in {\cal C}^{m \times n}$ be
given, $M, \, N \in {\cal C}^{m \times m}$ be two positive definite matrices.
Then

{\rm (a)} \ $r( \, A^{\dagger}_{M,N}A^k -A^kA^{\dagger}_{M,N}  \,)
= r\left[ \begin{array}{c} A^k \\  A^*M \end{array} \right]  +
 r[\, A^k, \ N^{-1}A^* \,] - 2r(A).$

{\rm (b)} \ $A^{\dagger}_{M,N}A^k  = A^kA^{\dagger}_{M,N}
 \Leftrightarrow  R(A^k) \subseteq
 R(N^{-1}A^*) \ \ and  \ \  R[(A^k)^*]  \subseteq R(MA). $ } 

 \medskip

\noindent {\bf Proof.}\, Follows from  (4.1). \qquad $ \Box$

\medskip

Based on the result in Theorem 16.9(b), we can extend the concept of power-EP matrix to weighted case:   A square matrix $ A $ is said to
be {\em weighted power-EP} if both $R(A^k) \subseteq R(N^{-1}A^*)$ and $R[(A^k)^*]
\subseteq R(MA)$ hold, where both $ M $ and $ N $ are  positive definite
matrices. 

\medskip

\noindent {\bf Theorem 16.10.}\, {\em Let $ A \in {\cal C}^{m \times m}$ be
given
with ${\rm Ind}(A) = k,$ and $ M, \, N  \in {\cal C}^{m \times m}$ be two
positive definite matrices. Then

{\rm (a)} \ $r( \, A^{\dagger}_{M,N}A^{D} - A^{D}A^{\dagger}_{M,N} \,)
= r \left[ \begin{array}{c} A^k \\  A^*M \end{array} \right] +
r[\, A^k, \ N^{-1}A^* \,] - 2r(A) = r( \, A^{\dagger}_{M,N}A^k -A^kA^{\dagger}_{M,N}  \,).$

{\rm (b)} \ $A^{\dagger}_{M,N}A^{D} = A^{D}A^{\dagger}_{M,N}
 \Leftrightarrow  R(A^k) \subseteq
 R(N^{-1}A^*) \ and  \   R[(A^k)^*] \subseteq R(MA),$ i.e.$,$ $ A $ is weighted power-EP.}

\medskip
  
\noindent {\bf Proof.}\, Follows from (4.1).   \qquad $ \Box$

\medskip

\noindent {\bf Theorem 16.11.}\, {\em Let $ A \in {\cal C}^{m \times n}$ be
given, $ M, \, S \in {\cal C}^{m \times m}$ and  $ N, \, T  \in {\cal C}^{n \times n} $
 be four positive definite matrices. Then

{\rm (a)} \ $r( \, A^{\dagger}_{M,N} - A^{\dagger}_{S,T}  \,) =
r \left[ \begin{array}{c} AN^{-1} \\  AT^{-1} \end{array} \right] +
 r[\, MA, \ SA \,] - 2r(A).$

{\rm (b)} \ $A^{\dagger}_{M,N} = A^{\dagger}_{S,T}   \Leftrightarrow
R(MA) =  R(SA)$ and $
 R[(AN^{-1})^*] = R[(AT^{-1})^*].$  }

\medskip

\noindent {\bf Proof.}\, Note that both $ A^{\dagger}_{M,N}$ and
$A^{\dagger}_{P,Q}$ are outer inverses of $ A $. Thus it follows by
Eq.(5.1) that
\begin{eqnarray*}
 r( \, A^{\dagger}_{M,N} -  A^{\dagger}_{S,T} \,) 
& = & r \left[ \begin{array}{c} A^{\dagger}_{M,N} \\ A^{\dagger}_{S,T}
 \end{array} \right] + r[\, A^{\dagger}_{M,N}, \  A^{\dagger}_{S,T} \,]
  - r(A^{\dagger}_{M,N}) - r( A^{\dagger}_{S,T}) \\
 & = & r \left[ \begin{array}{c} (MA)^* \\ (SA)^* \end{array} \right]
 + r[\, N^{-1}A^*, \ T^{-1}A^* \,] - 2r(A) \\
& = & r \left[ \begin{array}{c} AN^{-1} \\  AT^{-1} \end{array} \right] +
 r[ \, MA, \ SA \,] - 2r(A),
\end{eqnarray*}
establishing Part (a).  \qquad $ \Box$ 

\medskip

\noindent {\bf Theorem 16.12.}\, {\em Let $ A \in {\cal C}^{m \times n}$ be
given$,$
$ M, \, S  \in {\cal C}^{m \times m},$ $ N, \, T  \in {\cal C}^{n \times n} $
 be four positive definite matrices. Then

{\rm (a)} \  $r( \, AA^{\dagger}_{M,N} - AA^{\dagger}_{S,T}  \,)
= r[\, MA, \ SA \,] - 2r(A).$

{\rm (b)} \  $r( \, A^{\dagger}_{M,N}A - A^{\dagger}_{S,T}A  \,)
= r \left[ \begin{array}{c} AN^{-1} \\  AT^{-1} \end{array} \right]
- r(A).$

{\rm (c)} \ $r( \, A^{\dagger}_{M,N} - A^{\dagger}_{S,T}  \,)
= r( \, AA^{\dagger}_{M,N} - AA^{\dagger}_{S,T}  \,) +
r( \, A^{\dagger}_{M,N}A - A^{\dagger}_{S,T}A  \,)$. }

\medskip

\noindent {\bf Proof.}\, Follows from  (3.1) by noticing that
 $ AA^{\dagger}_{M,N}$, $A^{\dagger}_{M,N}A,$ $ AA^{\dagger}_{S,T}$ and 
$ A^{\dagger}_{S,T}A$ are idempotent matrices. \qquad $ \Box$

\medskip

\noindent {\bf Theorem 16.13.}\, {\em Let $ A \in {\cal C}^{m \times m}$ be
an idempotent or tripotent matrix$,$  and $ M, \, N  \in {\cal C}^{m \times m}$ 
be two positive definite matrices. Then

{\rm (a)} \ $r( \,  A - A^{\dagger}_{M,N}  \,)
= r[\, A^*, \ MA \,] + r[\, A^*, \ NA \,] - 2r(A).$

{\rm (b)} \ $ A = A^{\dagger}_{M,N} \Leftrightarrow  R(MA) =
  R(NA) = R(A^*),$ i.e.$,$ $ A $ is weighted EP. }

\medskip

\noindent {\bf Proof.}\, Note that $ A,  \, A^{\dagger}_{M,N} \in A\{2\}$ when
 $ A $ is idempotent or tripotent.  It follows by (5.1)  that
\begin{eqnarray*}
 r( \,  A -  A^{\dagger}_{M,N}\,) 
& = & r \left[ \begin{array}{c} A \\  A^{\dagger}_{M,N}  \end{array}
\right] + r[\, A, \  A^{\dagger}_{M,N} \,]
 -  r(A ) - r(A^{\dagger}_{M,N})  \\
& = & r \left[ \begin{array}{c}  A  \\ (MA)^* \end{array} \right]
+ r[\, A, \  N^{-1}A^* \,] - 2r(A) \\
& = & r[\, A^*, \ MA \,] + r[\, A^*, \ NA \,] - r(A) -r(A^k),
\end{eqnarray*}
as required for Part (a).  \qquad $ \Box$ 

\medskip

\noindent {\bf Theorem 16.14.}\, {\em Let $ A, \, B \in {\cal C}^{m \times m}$
be given, $ M, \, N \in {\cal C}^{m \times m}$ be two positive definite
 matrices. Then

{\rm (a)} \ $r( \, AA^{\dagger}_{M,N}B - BA^{\dagger}_{M,N}A  \,)
= r\left[ \begin{array}{c} A \\  A^*MB \end{array} \right]  +
 r[\, A, \ BN^{-1}A^* \,] - 2r(A).$

{\rm (b)} \ $r( \, A^{\dagger}_{M,N}AB - BAA^{\dagger}_{M,N}  \,) =
 r\left[ \begin{array}{c} AB \\  A^*M \end{array} \right]  +
 r[\, BA, \ N^{-1}A^* \,] - 2r(A).$

{\rm (c)} \ $AA^{\dagger}_{M,N}B  = BA^{\dagger}_{M,N}A
 \Leftrightarrow  R(BN^{-1}A^*) \subseteq
 R(A) \ \ and  \ \    R(B^*MA)  \subseteq R(A^*). $ 

{\rm (d)} \ $A^{\dagger}_{M,N}AB  = BAA^{\dagger}_{M,N}
 \Leftrightarrow  R(BA) \subseteq
 R(N^{-1}A^*) \ \ and  \ \  R[(AB)^*]  \subseteq R(MA). $ } 

\medskip

\noindent {\bf Proof.}\, Parts (a) and (b) Follow  from (4.1) by noticing 
that both $ AA^{\dagger}_{M,N}$ and $ A^{\dagger}_{M,N}A$ are idempotent. \qquad $ \Box$

\medskip

\noindent {\bf Theorem 16.15.}\, {\em Let $ A \in {\cal C}^{m \times m}$
be given with ${\rm Ind}(A) =1,$ and $ P, \, Q  \in {\cal C}^{m \times m}$ be two nonsingular matrices.
Then

{\rm (a)} \ $r[ \, (PAQ)^{\dagger} - Q^{-1} A^{\#}P^{-1} \,] =
r \left[ \begin{array}{c} A \\  A^*P^*P \end{array} \right] +
 r[\, A , \ QQ^*A^* \,] - 2r(A).$

{\rm (b)} \ $   (PAQ)^{\dagger} = Q^{-1} A^{\#}P^{-1} \Leftrightarrow
R(QQ^*A^*) =  R(A)$ and $  R(P^*PA) = R(A^*).$  }

\medskip

\noindent {\bf Proof.}\, It is easy to verify that both
$ (PAQ)^{\dagger}$ and $ Q^{-1}A^{\#}P^{-1}$ are outer inverses of $ PAQ $.
Thus it follows by (5.1) that
\begin{eqnarray*}
\lefteqn{ r[ \, (PAQ)^{\dagger} - Q^{-1} A^{\#}P^{-1} \,] }\\ 
& = & r \left[ \begin{array}{c} (PAQ)^{\dagger} \\ Q^{-1} A^{\#}P^{-1}
 \end{array} \right] + r[ \, (PAQ)^{\dagger}, \  Q^{-1} A^{\#}P^{-1} \,]
  - r[(PAQ)^{\dagger}] - r[ Q^{-1} A^{\#}P^{-1}]\\
& = & r \left[ \begin{array}{c} (PAQ)^*P \\ A^{\#}
 \end{array} \right] + r[ \, Q(PAQ)^*, \  A^{\#} \,]
  - 2r(A)\\
& = & r \left[ \begin{array}{c} A^*P^*P \\ A
 \end{array} \right] + r[ \, QQ^*A^*, \  A \,]
  - 2r(A),
\end{eqnarray*}
establishing Part (a) and then Part (a).  \qquad $ \Box$ 

\medskip

\noindent {\bf Theorem 16.16.}\, {\em Let $ A \in {\cal C}^{m \times m}$ be
given$,$ $ M, \ N  \in {\cal C}^{m \times m}$
be  two positive definite matrices$,$  and $ P, \, Q  \in {\cal C}^{m \times m}$ be two nonsingular matrices. Then

{\rm (a)} \ $r[ \, (PAQ)^{\dagger} - Q^{-1} A^{\dagger}_{M,N}P^{-1} \,] =
r \left[ \begin{array}{c} A \\  AQQ^*N \end{array} \right] +
 r[\, A, \ M^{-1}P^*PA^* \,] - 2r(A).$

{\rm (b)} \  $ (PAQ)^{\dagger} = Q^{-1} A^{\dagger}_{M,N}P^{-1}
\Leftrightarrow R( M^{-1}P^*PA^*) =  R(A)$ and $  R(NQQ^*A^*)
 = R(A^*).$  }

\medskip

\noindent {\bf Proof.}\, It is easy to verify that both
$ (PAQ)^{\dagger}$ and $ Q^{-1}A^{\dagger}_{M,N}P^{-1}$ are outer inverses
of $ PAQ $. Thus it follows by (5.1) that
\begin{eqnarray*}
\lefteqn{ r[ \, (PAQ)^{\dagger} - Q^{-1} A^{\dagger}_{M,N}P^{-1} \,] } \\
& = & r \left[ \begin{array}{c} (PAQ)^{\dagger} \\ Q^{-1}
A^{\dagger}_{M,N}P^{-1}  \end{array} \right] +
r[ \, (PAQ)^{\dagger}, \  Q^{-1} A^{\dagger}_{M,N}P^{-1} \,]
  - r[(PAQ)^{\dagger}] - r[ Q^{-1} A^{\dagger}_{M,N}P^{-1}]\\
& = & r \left[ \begin{array}{c} (PAQ)^*P \\ A^{\dagger}_{M,N}
 \end{array} \right] + r[ \, Q(PAQ)^*, \  A^{\dagger}_{M,N} \,]
  - 2r(A)\\
& = & r \left[ \begin{array}{c} A^*P^*P \\ (MA)^*
 \end{array} \right] + r[ \, QQ^*A^*, \  N^{-1}A^* \,]
  - 2r(A)\\
& = & r \left[ \begin{array}{c} A \\  AQQ^*N \end{array} \right] +
 r[\, A, \ M^{-1}P^*PA \,] - 2r(A),
\end{eqnarray*}
establishing Part (a).  \qquad $ \Box$

\markboth{YONGGE  TIAN }
{17. REVERSE ORDER LAWS FOR WEIGHTED MOORE-PENROSE INVERSES}

\chapter{Reverse order laws for  weighted  Moore-Penrose inverses}

\noindent Just as for  Moore-Penrose inverses and  Drazin inverses of products of matrices, we 
can also consider reverse order laws for weighted Moore-Penrose inverses of products of matrices. Noticing the basic fact in (16.2), we can easily extend 
the results in Chapter 8 to weighted  Moore-Penrose inverses of  products of matrices. 

\medskip

\noindent {\bf Theorem 17.1.}\, {\em Let $ A \in {\cal C}^{m \times n}, \, B \in
{\cal C}^{n \times k},$ and $C \in {\cal C}^{k \times l }$ be given and
let $ J = ABC$.
 Let $ M \in {\cal C}^{m \times m},$
$N \in {\cal C}^{ l \times l}, $ $ P \in {\cal C}^{n \times n},$
 and $Q \in {\cal C}^{k \times k}$
be four positive  definite  matrices. Then the following
 three statements are equivalent$:$

{\rm (a)}  \  $(ABC)^{\dagger}_{M,N} =
C^{\dagger}_{Q,N}B^{\dagger}_{P,Q}A^{\dagger}_{M,P}$.

{\rm (b)} \ $( M^{\frac{1}{2}}ABCN^{-\frac{1}{2}})^{\dagger} 
= (Q^{\frac{1}{2}} C N^{-\frac{1}{2}})^{\dagger}
 ( P^{\frac{1}{2}}BQ^{-\frac{1}{2}})^{\dagger}
( M^{\frac{1}{2}}AP^{-\frac{1}{2}})^{\dagger}. $  
 
{\rm (c)} \ $ r \left[ \begin{array}{ccc}
 BQ^{-1}B^*PB  & 0  & BC \\ 0 & -JN^{-1}J^*MJ &
JN^{-1}C^*QC  \\ AB & AP^{-1}A^*MJ &  0  \end{array} \right]
= r(B) + r(J). $ } 

\medskip

\noindent {\bf  Proof.}\, The equivalence of Part (a) and Part (b) follows directly
from applying (16.2) to the both sides of $(ABC)^{\dagger}_{M,N} =
C^{\dagger}_{Q,N}B^{\dagger}_{P,Q}A^{\dagger}_{M,P}$ and simplifying. 
 Observe that the left-hand side of Part (b) can also be written as
$$ 
( M^{\frac{1}{2}}ABCN^{-\frac{1}{2}})^{\dagger} 
=[ \, (M^{\frac{1}{2}}AP^{-\frac{1}{2}}) (P^{\frac{1}{2}}B Q^{-\frac{1}{2}}) 
(Q^{\frac{1}{2}}CN^{-\frac{1}{2}}) \,]^{\dagger}.
$$
In that case, we see by Theorem 8.11 that Part (b) holds if and only if 
$$  
r \left[ \begin{array}{ccc}  B_1B_1^*B_1  & 0  & B_1C_1 \\
0 & -J_1J_1^*J_1 & J_1C_1^*C_1  \\
A_1B_1 & A_1A_1^* J_1 &  0  \end{array} \right]
=  r(B_1) + r(J_1),  
$$
where 
$$ 
A_1 = M^{\frac{1}{2}}AP^{-\frac{1}{2}}, \ \ \
 B_1 = P^{\frac{1}{2}}B Q^{-\frac{1}{2}},\ \ \
C_1 = Q^{\frac{1}{2}}CN^{-\frac{1}{2}},  \ \ \ J_1 =
 M^{\frac{1}{2}}ABCN^{-\frac{1}{2}}.
$$ 
Simplifying this rank equality by the given condition
that $ M, \ N, \ P $ and $ Q $ are positive definite,
we obtain the rank equality in Part (c). \qquad $\Box$ 

\medskip

\noindent {\bf Corollary 17.2.}\, {\em Let $ A \in {\cal C}^{m \times n},
\, B \in {\cal C}^{n \times k}$ and $C \in {\cal C}^{k \times l }$
be given and let $ J = ABC$.
 Let $P \in {\cal C}^{ n \times n }$ and
$Q \in {\cal C}^{ k \times k}$ be two positive  definite matrices.
Then the following three statements are equivalent:

{\rm (a)} \  $(ABC)^{\dagger}
 = C^{\dagger}_{Q,I}B^{\dagger}_{P,Q}A^{\dagger}_{I,P}$.

{\rm (b)} \  $(ABC)^{\dagger}
= ( Q^{\frac{1}{2}}C)^{\dagger}
 ( P^{\frac{1}{2}}BQ^{-\frac{1}{2}})^{\dagger}
(AP^{-\frac{1}{2}})^{\dagger}. $
 
{\rm (c)} \  $r \left[ \begin{array}{ccc}  BQ^{-1}B^*PB  & 0  & BC \\ 0 & -JJ^*J &
JC^*QC  \\ AB & AP^{-1}A^*J &  0  \end{array} \right]
 = r(B) + r(J). $  } 

\medskip

\noindent {\bf Proof.}\, Follows from Theorem 17.1 by setting $ M $ and $ N$ as
identity matrices. \qquad $\Box$ 

\medskip

\noindent {\bf Corollary 17.3.}\, {\em Let $ A \in {\cal C}^{m \times n},
\, B \in {\cal C}^{n \times k}$, and $C \in {\cal C}^{k \times l }$
be given and denote $ J = ABC$. Let $ M \in {\cal C}^{m \times m},$
$N \in {\cal C}^{ l \times l}$ be two positive definite  matrices.
Then the following three statements are equivalent$:$

{\rm (a)} \  $(ABC)^{\dagger}_{M,N}
= C^{\dagger}_{I,N}B^{\dagger}A^{\dagger}_{M,I}$.

{\rm (b)} \  $( M^{\frac{1}{2}}ABCN^{-\frac{1}{2}})^{\dagger}
= (CN^{-\frac{1}{2}})^{\dagger}B^{\dagger}
( M^{\frac{1}{2}}A)^{\dagger}. $  
 
{\rm (c)}  \ $r \left[ \begin{array}{ccc}
 BB^*B  & 0   BC \\ 0 & -JN^{-1}J^*MJ &
JN^{-1}C^*C  \\ AB & AA^*MJ &  0  \end{array} \right]
= r(B) + r(J). $ } 

\medskip

\noindent {\bf Proof.}\, Follows from Theorem 17.1 by setting $ P $
and $ Q$ as identity matrices.  \qquad $\Box$ 

\medskip

\noindent {\bf Corollary 17.4.}\, {\em Let $ A \in {\cal C}^{m \times n}, \,
B \in {\cal C}^{n \times k}$,   and 
$C \in {\cal C}^{k \times l }$ be given with $r(A)
= n $ and $r( C ) = k$.
 Let $ M \in {\cal C}^{m \times m},$ $N \in {\cal C}^{ l \times l},$
 $P \in {\cal C}^{ n \times n },$  and $Q \in {\cal C}^{ k \times k}$
be four positive  definite matrices. Then the following
two statements are equivalent$:$

{\rm (a)} \ $(ABC)^{\dagger}_{M,N}
= C^{\dagger}_{Q,N}B^{\dagger}_{P,Q}A^{\dagger}_{M,P}$.

{\rm (b)} \  $ R(\,P^{-1}A^*MAB\,) \subseteq  R(B) \ \ and
\ \ R[ \,(BCN^{-1}C^*Q)^*\, ] \subseteq  R(B^*). $ }

\medskip

\noindent {\bf  Proof.}\, The given condition $r( A) = n $
and $r(C)= k$ is equivalent to
$ A^{\dagger}A = I_n$, $ CC^{\dagger} = I_k$,  and $r( ABC ) = r(B)$. In that case,
we can show by block elementary operations that 
$$
 \left[ \begin{array}{ccc}  BQ^{-1}B^*PB  & 0  & BC \\
  0 & -JN^{-1}J^*MJ &
JN^{-1}C^*QC  \\ AB & AP^{-1}A^*MJ &  0  \end{array} \right]
 \ \  {\rm and } \ \
\left[ \begin{array}{ccc}   0 & 0  & B \\ 0 & 0 & BCN^{-1}C^*Q
\\ B & P^{-1}A^*MAB &  0  \end{array} \right] 
$$
are equivalent, the detailed is omitted here. This result implies that  
\begin{eqnarray*}
\lefteqn { r \left[ \begin{array}{ccc}  BQ^{-1}B^*PB
 & 0  & BC \\ 0 & -JN^{-1}J^*MJ &
JN^{-1}C^*QC  \\ AB & AP^{-1}A^*MJ &  0  \end{array} \right] } \\
& = & r \left[ \begin{array}{ccc} 0 & 0  & B \\ 0 & 0
 & BCN^{-1}C^*Q  \\ B & P^{-1}A^*MAB &  0  \end{array} \right]  \\
& = & r \left[ \begin{array}{c} B  \\ BCN^{-1}C^*Q
 \end{array} \right] + r [\,  B, \ P^{-1}A^*MAB \,].
\end{eqnarray*}
Thus  under the given condition of this corollary,  Part (c)
of Theorem 17.1 reduces to
$$
r\left[ \begin{array}{c} B  \\ BCN^{-1}C^*Q
\end{array} \right] +
r [\,  B, \ P^{-1}A^*MAB \,] = 2r(B),
$$ 
which is obviously equivalent to Part (c) of this corollary.
  \qquad $\Box$ 

\medskip

\noindent {\bf Corollary 17.5.}\, {\em Let $ A \in {\cal C}^{m \times m}, \, B
 \in {\cal C}^{m \times n},$and  $C \in {\cal C}^{n \times n}$ be given with
 $A $ and $ C $ nonsingular.   Let $ M, \, P \in {\cal C}^{m \times m} $
  and $N, \ Q  \in {\cal C}^{ n \times n}$ be four positive definite
   Hermitian
 matrices. Then 

{\rm (a)} \ $(ABC)^{\dagger}_{M,N} = C^{-1}B^{\dagger}_{P,Q}A^{-1}
\Leftrightarrow
R( \,P^{-1}A^*MAB \,) = R(B) $ and  $ R[\, (BCN^{-1}C^*Q)^* \,] = R(B^*). $

{\rm (b)} \ $(ABC)^{\dagger}_{M,N} = C^{-1}B^{\dagger} A^{-1}
 \Leftrightarrow
R(\,A^*MAB\, ) = R(B)$ and $  \ \ R[\, (BCN^{-1}C^*)^*\, ] = R(B^*).$

{\rm (c)}  \ $(ABC)^{\dagger} = C^{-1}B^{\dagger}_{P,Q}A^{-1}
 \Leftrightarrow
R( \,P^{-1}A^*AB \,) = R(B)  $ and  $ R[\, (BCC^*Q)^*\, ] = R(B^*). $ \\
In particular,  the following two identities hold 
$$\displaylines{
\hspace*{1.5cm}
(ABC)^{\dagger}_{M,N} = C^{-1}B^{\dagger}_{(A^*MA),\:
(CN^{-1}C^*)^{-1}} A^{-1}, \hfill (17.1)
\cr
\hspace*{1.5cm}
(ABC)^{\dagger} = C^{-1}B^{\dagger}_{(A^*A), \: (CC^*)^{-1}} A^{-1}. \hfill (17.2)
\cr}
$$} 
{\bf Proof.}\, Let $ A $ and $ C $ be nonsingular matrices in Corollary 17.4.
 We can obtain Part (a) of this corollary. Parts (a) and (b) are special
 cases of Part (a). The equality (17.1) follows from Part (a)
by setting $ P = A^*MA$ and $ Q = (CN^{-1}C^*)^{-1}$. \qquad $\Box$ 

\medskip

\noindent {\bf Theorem 17.6.}\, {\em Let $ A \in {\cal C}^{m \times n}, \, B \in
 {\cal C}^{n \times k}$, and
$C \in {\cal C}^{k \times l }$ be given and denote $ J = ABC$.
Let $ M \in {\cal C}^{m \times m}, \,N \in {\cal C}^{ l \times l}, \, P \in {\cal C}^{n \times n},$
and $Q \in {\cal C}^{k \times k}$
be four positive  definite matrices. Then the following two
statements are equivalent$:$

{\rm (a)} \ $(ABC)^{\dagger}_{M,N}
 = (BC)^{\dagger}_{P,N}B(AB)^{\dagger}_{M,Q}$.

{\rm (b)} \ $r  \left[ \begin{array}{ccc}
 B^*  & (AB)^*MJ  & B^*PBC \\
 JN^{-1}(BC)^* & 0 & 0 \\ ABQ^{-1}B^* &  0 &  0  \end{array} \right] 
= r(B) + r(J). $ } 

\medskip

\noindent {\bf  Proof.} \, Write $ABC $ as $ ABC = (AB)B^{\dagger}_{P,Q}(BC)$ and notice that
$ ( B^{\dagger}_{P,Q})^{\dagger}_{Q,P} = B$.  Then by Theorem 17.1, we know
that
\begin{eqnarray*}
(ABC)^{\dagger}_{M,N} = [ \,(AB)B^{\dagger}_{P,Q}(BC) \,]^{\dagger}_{M,N} 
=  (BC)^{\dagger}_{P,N} ( B^{\dagger}_{P,Q})^{\dagger}_{Q,P}(AB)^{\dagger}_{M,Q} = (BC)^{\dagger}_{P,N}B(AB)^{\dagger}_{M,Q}
\end{eqnarray*}
holds if and only if 
$$
r \left[ \begin{array}{ccc}
B^{\dagger}_{P,Q}P^{-1}(B^{\dagger}_{P,Q})^*Q B^{\dagger}_{P,Q}
 & 0  & B^{\dagger}_{P,Q}BC \\ 0 & -JN^{-1}J^*MJ & JN^{-1}(BC)^*P(BC)
 \\ ABB^{\dagger}_{P,Q}
 & ABQ^{-1}(AB)^*MJ & 0 \end{array} \right] = r(
 B^{\dagger}_{P,Q}) + r(J). \eqno (17.3)
$$ 
Note by (1.5) that
\begin{eqnarray*}
B^{\dagger}_{P,Q}P^{-1}(B^{\dagger}_{P,Q})^*Q B^{\dagger}_{P,Q}
 = Q^{-\frac{1}{2}}( P^{\frac{1}{2}} BQ^{-\frac{1}{2}})^{\dagger} 
[\,(P^{\frac{1}{2}} BQ^{-\frac{1}{2}})^{\dagger} \,]^*
( P^{\frac{1}{2}} AQ^{-\frac{1}{2}})^{\dagger}
P^{\frac{1}{2}}. 
\end{eqnarray*}
Thus by block elementary operations, we can deduce that
 (17.3) is equivalent
 to Part (c) of the theorem. The details are omitted.  \qquad $\Box$ 

\medskip
 
\noindent {\bf Corollary 17.7.}\, {\em Let $ A \in {\cal C}^{m \times n}, \, B
\in {\cal C}^{n \times k}$, and
$C \in {\cal C}^{k \times l }$ be given and denote $ J = ABC$.
 Let $ M \in {\cal C}^{m \times m},$
$N \in {\cal C}^{ l \times l}, $ $ P \in {\cal C}^{n \times n},$
and $Q \in {\cal C}^{k \times k}$
be four positive  definite matrices. If 
$$\displaylines{
\hspace*{2cm}
r(ABC )= r(B), \hfill (17.4)
\cr}
$$   
then the weighted Moore-Penrose inverse of the product $ABC$ satisfies the
following two equalities
$$\displaylines{
\hspace*{2cm}
(ABC)^{\dagger}_{M,N} = (BC)^{\dagger}_{P,N}B(AB)^{\dagger}_{M,Q}, \hfill (17.5)
\cr
\hspace*{0cm}
and \hfill
\cr
\hspace*{2cm}
(ABC)^{\dagger}_{M,N} =
(B^{\dagger}_{P,Q}BC)^{\dagger}_{Q,N}B^{\dagger}_{P,Q}
(ABB^{\dagger}_{P,Q}  )^{\dagger}_{M,Q}. \hfill (17.6)
\cr}
$$}
{\bf Proof.}\, Under (17.4), we know that
$$\displaylines{
\hspace*{2cm}
r(AB) = r(BC) = r(B),\hfill
\cr}
$$ 
which is equivalent to 
$$ \displaylines{
\hspace*{2cm}
R(BC) = R(B), \ \ \ {\rm and} \ \ \ R(B^*A^*) = R(B^*).\hfill
\cr}
$$ 
Based on them we further obtain
$$ \displaylines{
\hspace*{2cm}
R(B^*PBC) = R(B^*PB) =
R[ \,(B^*P^{\frac{1}{2}})(B^*P^{\frac{1}{2}})^* \,]
 =  R(B^*P^{\frac{1}{2}})
= R(B^*),\hfill
\cr
\hspace*{0cm}
and \hfill
\cr
\hspace*{2cm}
R(BQ^{-1}B^*A^*) = R(BQ^{-1}B^*) = R[\, (BQ^{-\frac{1}{2}})
 (BQ^{-\frac{1}{2}})^* \,]
= R(BQ^{-\frac{1}{2}}) = R(B). \hfill
\cr}
$$
Under these two conditions, the left-hand side of Part (b)
 in Theorem 17.6 reduces to $2r(B)$.
Thus  Part (b) in Theorem 17.6 is indentity under  (17.4).
Therefore we have  (17.5) under (17.4). Consequently
writing $ ABC $ as $ ABC = (AB)B^{\dagger}_{P,Q}(BC)$ and applying
 (17.5) to it yields (17.6). \qquad $\Box$ 

\medskip

Some applications of Corollary 17.7 are given below. 

\medskip

\noindent {\bf Corollary 17.8.}\, {\em Let $ A, \, B \in {\cal C}^{m \times n}$ be given$,$
$ M \in {\cal C}^{m \times m},$ $N \in {\cal C}^{ n \times n},$ $P \in {\cal C}^{ 2m \times 2m},$ and
 $Q \in {\cal C}^{2n \times 2n}$ be
four positive definite matrices. If $ A $ and $ B $ satisfy
the rank additivity condition
$$ 
r(\, A + B \,)  = r(A)  + r(B), \eqno (17.7)
$$   
then the weighted Moore-Penrose of  $ A + B $ satisfies the two
equalities
$$
(\, A + B \,)^{\dagger}_{M,N}
 = \left[ \begin{array}{c} A  \\ B \end{array} \right]^{\dagger}_{P,N} 
\left[ \begin{array}{cc} A  & 0  \\  0 & B \end{array} \right]
[\, A, \ B \,]^{\dagger}_{M,Q}
, \eqno (17.8)
$$
$$
(\, A + B \,)^{\dagger}_{M,N}
= \left[ \begin{array}{c} A^{\dagger}_{M,N}A  \\ B^{\dagger}_{M,N}B
 \end{array} \right]^{\dagger}_{Q,N} \left[ \begin{array}{cc}
 A^{\dagger}_{M,N} & 0  \\  0 & B^{\dagger}_{M,N} \end{array} \right]
 [\, AA^{\dagger}_{M,N}, \ BB^{\dagger}_{M,N}\,]^{\dagger}_{M,P}
. \eqno (17.9)
$$ }
{\bf Proof.}\, Write $ A + B $ as
$$ 
A + B = [\, I_m, \  I_m \,] \left[ \begin{array}{cc} A  & 0  \\  0 & B
\end{array} \right]\left[ \begin{array}{c} I_n \\ I_n \end{array} \right]
 := UDV.
$$
Then the condition  (17.7) is equivalent to $r(UDV) =
 r(D)$. Thus it turns out that  
$$
 (UDV)^{\dagger}_{M,N} =
(DV)^{\dagger}_{P,N}D(UD)^{\dagger}_{M,Q},
$$ 
which is exactly (17.8). Next write $ A + B $ as
$$ 
A + B = [\, A, \  B \,] \left[ \begin{array}{cc} A^{\dagger}_{M,N}
 & 0  \\ 0 & B^{\dagger}_{M,N} \end{array} \right]
 \left[ \begin{array}{c} A \\ B \end{array} \right]
 := U_1D_1V_1.
$$
Then the condition  (17.7) is also  equivalent to
 $ r(U_1D_1V_1) = r(D_1)$. Thus it follows by  (17.5) that 
$ (U_1D_1V_1)^{\dagger}_{M,N} =
(D_1V_1)^{\dagger}_{Q,N}D_1(U_1D_1)^{\dagger}_{M,P}$,
which is exactly (17.9). \qquad $\Box$ 

\medskip

A generalization of Corollary 17.8 is presented below, the proof is omitted. 

\medskip

\noindent {\bf Corollary 17.9.}\, {\em Let $ A_1, \,  \cdots, \,  A_k \in
{\cal C}^{m \times n}$ be given$,$
and let $ M \in {\cal C}^{m \times m},$ $N \in {\cal C}^{ n \times n},$
$P \in {\cal C}^{ km \times km}, $ and
 $Q \in {\cal C}^{kn \times kn}$ be four  positive  definite Hermitian
 matrices. If
$$ 
r(\, A_1 +  \cdots + A_k \,)  = r(A_1)  + \cdots +
r(A_k), \eqno (17.10)
$$   
then the weighted Moore-Penrose inverse of the sum satisfies the following
 two equalities
$$\displaylines{
\hspace*{1cm}
(\, A_1 +  \cdots + A_k \,)^{\dagger}_{M,N} 
= \left[ \begin{array}{c} A_1 \\ \vdots \\ A_k \end{array}
\right]^{\dagger}_{P,N}
\left[ \begin{array}{ccc} A_1  & &   \\   & \ddots &  \\ & & A_k \end{array}
 \right]
[\, A_1, \ \cdots, \  A_k \,]^{\dagger}_{M,Q}, \hfill (17.11)
\cr}
$$
$$
\displaylines{
\hspace*{1cm}
\left( \sum_{t=1}^{k} A_t \right)^{\dagger}_{M,N} =
 \hfill
\cr
\hspace*{1cm}
\left[ \begin{array}{c} (A_1)_{M,N}^{\dagger}A_1 \\ \vdots \\
(A_k)_{M,N}^{\dagger}A_k \end{array} \right]^{\dagger}_{Q,N}
\left[ \begin{array}{ccc} (A_1)_{M,N}^{\dagger}  & &   \\   & \ddots &
 \\ & & (A_k)_{M,N}^{\dagger} \end{array} \right]
[\, A_1(A_1)_{M,N}^{\dagger}, \ \cdots, \  A_k(A_k)_{M,N}^{\dagger}
 \,]^{\dagger}_{M,P}. \hfill (17.12)
\cr}
$$ }


\noindent {\bf Corollary 17.10.}\, {\em Let $ A \in {\cal C}^{m \times n},$ $ B\in
{\cal C}^{m \times k},$ $ C \in {\cal C}^{l \times n},$ $ A
\in {\cal C}^{l \times k}$  be given$,$
$ M, \, P \in {\cal C}^{(m+l) \times (m+l)},$ $N, \ Q \in
{\cal C}^{ (n + k)\times (n +k) }$ be four positive  definite
matrices. If
$$ \displaylines{
\hspace*{2cm}
r \left[ \begin{array}{cc} A  & B \\  C & D \end{array} \right]
=r(A), \hfill (17.13)
\cr}
$$
or equivalently $ AA^{\dagger}B = B,$ $CA^{\dagger}A = C$ and
$ D = CA^{\dagger}B,$ then
$$\displaylines{
\hspace*{2cm}
\left[ \begin{array}{cc} A  & B \\  C & D \end{array} \right]^{\dagger}_{M,N} 
= \left[ \begin{array}{cc} A  & B \\  0 & 0 \end{array} \right]^{\dagger}_{P,N}
\left[ \begin{array}{cc} A  & 0 \\  0 & 0 \end{array} \right] 
\left[ \begin{array}{cc} A  & 0 \\  C & 0 \end{array}
\right]^{\dagger}_{M,Q}. \hfill (17.14)
\cr}
$$ 
In particular$,$ 
$$\displaylines{
\hspace*{2cm}
\left[ \begin{array}{cc} A  & B \\  C & D \end{array}
\right]^{\dagger}_{M,N}
= [\, A, \ B \,]^{\dagger}_{I,N} A \left[ \begin{array}{c} A \\ C \end{array} \right]^{\dagger}_{M,I}. 
\hfill (17.15)
\cr}
$$}       
{\bf Proof.}\, Under (17.13), we see that 
$$ \displaylines{
\hspace*{2cm}
\left[ \begin{array}{cc} A  & B \\  C & D \end{array} \right] 
=\left[ \begin{array}{cc} I_m  & 0 \\  CA^{\dagger} & I_l \end{array} \right]
\left[ \begin{array}{cc} A  & 0 \\  0 & 0 \end{array} \right]
\left[ \begin{array}{cc} I_n  & A^{\dagger}B \\  0 & I_k \end{array} \right]
  := ULV.\hfill
\cr}
$$ 
Thus by Corollary 17.6, we obtain 
$$\displaylines{
\hspace*{2cm}
\left[ \begin{array}{cc} A  & B \\  C & D \end{array} \right]^{\dagger}_{M,N} 
= (LV)^{\dagger}_{P,N}L(UL)^{\dagger}_{M,Q},  \hfill
\cr}
$$ 
which  is exactly (17.14). When $ P = I_{m + l}$ and  $ Q = I_{n + k}$,
we have 
\begin{eqnarray*} 
\left[ \begin{array}{cc} A  & B \\  0 & 0 \end{array}
\right]^{\dagger}_{I,N}
\left[ \begin{array}{c} I_m \\ 0 \end{array} \right]
& = & N^{-\frac{1}{2}}\left( \, \left[ \begin{array}{cc} A  & B \\  0 & 0
\end{array} \right]
N^{-\frac{1}{2}} \, \right)^{\dagger} \left[ \begin{array}{c} I_m \\ 0
\end{array} \right] \\
& = & N^{-\frac{1}{2}} \left[ \begin{array}{c} ( \,
\left[\, A,\ B \,\right] N^{-\frac{1}{2}} \, )^{\dagger} \\
\left[\, 0 , \ 0 \, \right]^{\dagger} \end{array} \right]
\left[ \begin{array}{c} I_m \\ 0 \end{array} \right] \\
& = &  N^{-\frac{1}{2}} ( [ \, A, \ B \, ] N^{-\frac{1}{2}})^{\dagger}
= [\, A, \ B \,]^{\dagger}_{I,N}.
\end{eqnarray*} 
Similarly we can deduce 
$$\displaylines{
\hspace*{2cm} 
[\, I_n, \ 0 \,]\left[ \begin{array}{cc} A  & 0 \\  C & 0 \end{array} \right]^{\dagger}_{M,I}
 =\left[ \begin{array}{c} A \\ C \end{array} \right]^{\dagger}_{M,I}. \hfill
\cr}
$$
Putting both of them in (17.14) yields  (17.15).  \qquad $\Box$

\markboth{YONGGE  TIAN }
{18. EXTREME RANKS OF $ A-BXC$ }

\chapter{Extreme  ranks of $ A- BXC$}

\noindent The basic tool for establishing the whole work in the monograph is the rank formula (2.1) for
 the Schur complement $ D - CA^{\dagger}B$. Motivated by (2.1),  one might naturally consider 
the rank of a generalized Schur complement $ D - CA^-B$, where $ A^-$ is an inner inverse of $ A$. 
Since $ A^-$ is not unique in general, the rank of $ D - CA^-B$ will depend on the choice of $ A^-$. 
Thus a fundamental problem related to a generalized Schur complement $ D - CA^-B$  is to
 find its maximal and minimal possible ranks with respect to the choice of $ A^-$. Notice that 
the general expression of  $ A^-$ is  $ A^- = A^{\dagger} + F_AV +  WE_A
,$  where both $ V $ and $ W $ are arbitrary matrices. As a consequence, 
$$
 D - CA^-B =  D - CA^{\dagger}B - CF_AVB  - CWE_AB. 
$$
This expression implies that $ D - CA^-B$ is in fact a matrix expression with two independent 
variant matrices. This fact motivates us to consider another basic problem in matrix 
theory---maximal and minimal possible  ranks of linear matrix expressions with variant 
matrices.  In this chapter, we consider the simplest case--- the maximal and the minimal ranks of the matrix expression $ A - BXC$ with respect to the variant matrix $X$ and then discuss
  some related topics. Throughout the  symbols $ E_A $ and
 $ F_A $ stand for the two oblique projectors $ E_A  = I - AA^-$ and
 $F_A = I - A^-A$ induced by $ A$. 

\medskip


The following result is well known (see, e.g., \cite{RM}).  

\medskip

\noindent {\bf Lemma 18.1.}\, {\em Suppose $BXC = A$ is a linear matrix equation 
over  an arbitrary field $ {\cal F},$ where  $ A \in {\cal F}^{ m \times n }, \, B \in {\cal F}^{ m \times k }, \,
 C \in {\cal F}^{l \times n}$  are given. Then it is consistent if and only if $  
R( A ) \subseteq  R( B )$ and $R( A^T )
\subseteq R( C^T ),$ or equivalently $ BB^- AC^- C = A.$  In that case$,$  the 
general solution  of $BXC = A$ can be expressed as
$$ 
X = B^- AC^- + U - B^- BUCC^- , \ \  or \ \ X = B^- AC^- + F_BV + WE_A,      
$$ 
where $ U, \ V $ and $ W $ are arbitrary matrices. In particular$,$ the solution of 
 $BXC = A$  is unique if and only if  $B$ has
full column rank and $C$ has full row rank.}

\medskip

In order to determine the maximal and minimal ranks of  $ A- BXC$ with respect to $ X $, we first
  establish two rank identities for $ A- BXC$ through (1.4) and (1.5). 

\medskip

\noindent {\bf  Theorem 18.2.}\, {\em The matrix expression  $ A- BXC$ satisfies
 the rank identity
$$
r( \, A - BXC \, ) = r \left[ \begin{array}{c}  A  \\ C  \end{array} \right]
 +  r[ \, A, \ B \, ] -r(M) +
r[ \,E_{T_1}(\, X + TM^-S \, )F_{S_1} \, ], \eqno (18.1)
$$ 
where $ M = \left[ \begin{array}{cc}  A  & B \\ C  & 0 \end{array} \right], \
 T = [ \, 0, \ I_k \,]$ and $S = \left[ \begin{array}{c}  0  \\ I_l
  \end{array} \right],$ $ T_1 = TF_M,$ and $ S_1 = E_MS$.}

\medskip

\noindent {\bf Proof.}\, It is easy to verify by block elementary operations
of matrix that 
$$
r( \, A - BXC \, ) = r \left[ \begin{array}{ccc}  A & B & 0  \\ C & 0 & I_l
\\ 0 & I_k & -X  \end{array} \right] - k - l  =
r\left[ \begin{array}{cr}  M & S  \\ T & -X   \end{array} \right] - k - l.
\eqno (18.2)
$$ 
Applying (1.6) to the block matrix in it, we find  that
\begin{eqnarray*}
 r\left[ \begin{array}{cr}  M & S  \\ T & -X   \end{array} \right] &= & r \left[ \begin{array}{c}  M  \\ T  \end{array} \right] +  r[ \, M, \ S \, ] -r(M)+ r[ \,E_{T_1}(\, X + TM^-S \, )F_{S_1} \, ] \\
&= & r \left[ \begin{array}{c}  A  \\ C  \end{array} \right] +  r[ \, A, \ B \, ] + k + l -r(M)+ 
 r[ \,E_{T_1}(\, X + TM^-S \, )F_{S_1} \, ].
 \end{eqnarray*}
Putting it in (18.2) yields (18.1). \qquad  $\Box$.

\medskip

\medskip

\noindent {\bf Theorem  18.3.}\, {\em The matrix expression  $ A- BXC$ satisfies
 the rank identity
$$ 
 r( A - BXC ) = r[ \, A, \  B  \,] +   r \left[ \begin{array}{c}
 A \\ C \end{array} \right] -
r \left[ \begin{array}{cc}  A  & B  \\ C &  0 \end{array} \right]
+ r( E_{A_2}AF_{A_1} - E_{A_2}BXCF_{A_1}),
 \eqno (18.3)
$$ 
where $ A_1 = E_BA, \ A_2 = AF_C,$  and the matrix equation 
$E_{A_2}BXCF_{A_1} = E_{A_2}AF_{A_1}$ is consistent.}

\medskip

\noindent {\bf Proof}.\,  We first establish the following rank equality 
$$
\displaylines{
\hspace*{2cm}
r \left[ \begin{array}{cc}  A  & B  \\ C &  0 \end{array} \right]
= r \left[ \begin{array}{c}  A  \\ C  \end{array} \right]  +
 r[ \, A, \  B \, ] - r( A ) + r( E_{A_2}A F_{A_1}).  \hfill (18.4)
\cr}
$$
Observe that
$$\displaylines{
\hspace*{2cm}
r( E_BAF_C ) = r \left[ \begin{array}{cc}  A  & AF_C  \\ E_BA  &  0
\end{array} \right] - r( A ), \hfill
\cr}
$$
and also observe from (1.4) that 
$$
\displaylines{
\hspace*{2cm}
 r \left[ \begin{array}{cc}  A  & AF_C  \\ E_BA  &  0 \end{array} \right]
 = r(E_BA ) +  r(AF_C) +
 r(E_{A_2}A F_{A_1}). \hfill
\cr}
$$
We obtain $r( E_BAF_C ) = r(  E_BA ) +  r( AF_C) - r( A ) +  r( E_{A_2}A F_{A_1}).$
Putting it in (1.4) and applying  (1.2) and (1.3), we get  (18.4). Next
replace the matrix $ A $ in  (18.4) by $ p(X) = A - BXC$  and notice that
$$
\displaylines{
\hspace*{2cm}
r \left[ \begin{array}{cc}  A - BXC  & B  \\ C &  0 \end{array} \right]
 = r \left[ \begin{array}{cc} A  & B  \\ C &  0 \end{array} \right],
  \qquad  r \left[ \begin{array}{c}  A - BXC   \\ C  \end{array} \right]
  = r \left[ \begin{array}{c}  A  \\ C  \end{array} \right], \hfill
\cr
\hspace*{1cm}
r[  \,A - BXC, \  B \, ] = r[ \, A, \ B \, ],   \qquad   E_B( A - BXC ) = E_BA, \qquad 
 ( A - BXC )F_C = AF_C. \hfill
\cr}
$$
Then (18.4) becomes
$$\displaylines{
\hspace*{2cm}
r \left[ \begin{array}{cc}  A  & B  \\ C &  0 \end{array} \right]
= r \left[ \begin{array}{c}  A  \\ C  \end{array} \right]  +
 r[ \, A, \  B \, ] - r( A - BXC) + r( E_{A_2}A F_{A_1} -
 E_{A_2}BXC F_{A_1} ), \hfill
\cr}
$$
establishing (18.3). On the other hand, from
$ E_{A_2}A_2 = 0$ and $ A_1F_{A_1} = 0$ we can deduce that
$E_{A_2}AC^-C = E_{A_2}A $ and $ BB^-AF_{A_1} = AF_{A_1}.$ Thus 
$R( E_{A_2}A F_{A_1}) = R( E_{A_2}BB^-AF_{A_1} )
\subseteq R( E_{A_2}B )$ and $R[ \, ( E_{A_2}A F_{A_1})^T \, ]$ $ 
= R[\, ( E_{A_2}AC^-CF_{A_1})^T ]
 \subseteq R[ \, (CF_{A_1})^T \, ].$ Both of them imply that the matrix equation $ E_{A_2}BXCF_{A_1} =
 E_{A_2}AF_{A_1}$ is consistent. \qquad   $\Box $

\medskip


On the basis of (18.3), we establish the main result of the chapter.   

\medskip

\noindent {\bf  Theorem 18.4.}\, {\em \  Let $ A \in {\cal F}^{ m \times n }, \, B 
\in {\cal F}^{ m \times k }$ and $ C \in {\cal F}^{ l \times n }$ be given. Then

{\rm (a)}\, The maximal rank of $A - BXC $ with respect to $ X$ is   
$$
\displaylines{
\hspace*{2cm}
\max_{X}r( \, A - BXC \, ) = \min  \left\{ \, r[ \, A, \  B \, ], \ \ \ \
r \left[ \begin{array}{c}  A  \\ C  \end{array} \right]  \, \right\}.
\hfill (18.5)
\cr}
$$ 

{\rm (b)}\, The minimal rank of $A - BXC $ with respect to $ X$ is 
$$
\displaylines{
\hspace*{2cm}
\min_{X}r( \, A - BXC \, )  =  r[  \, A, \ B \, ] +
r \left[ \begin{array}{c}  A  \\ C  \end{array} \right] -
r \left[ \begin{array}{cc}  A  & B  \\ C &  0 \end{array} \right].
\hfill (18.6)
\cr}
$$ 

{\rm (c)}\, The general expression of $X$ satisfying {\rm (18.5)} can be written as
$$\displaylines{
\hspace*{2cm}
X = - TM^-S + U,  \hfill (18.7) 
\cr}
$$ 
where $U$ is chosen such that $ r(E_{T_1}UF_{S_1}) = \min \{ \, r(E_{T_1}), \ r(F_{S_1})  \, \},$
where $M, \ S,  \ T, \ S_1$ and $T_1$ are defined in {\rm (18.1)}. 

{\rm (d)}\, The matrix $X$ satisfying 
 {\rm (18.6)} is determined by the matrix equation $E_{T_1}(\, X + TM^-S \, )F_{S_1} = 0, $
and can be written as 
$$\displaylines{
\hspace*{2cm}
X = - TM^-S + T_1V + WS_1, \hfill (18.8)
\cr}
$$
where $ V $ and $W$ are arbitrary.}

\medskip

\noindent {\bf Proof.}\, Eq.\,(18.2) implies that 
$$
\displaylines{
\hspace*{1cm}
\max_Xr( \, A - BXC \, ) = r \left[ \begin{array}{c}  A  \\ C  \end{array}
\right] +  r[ \, A, \ B \, ] -r(M) + \max_Xr[ \,E_{T_1}(\, X + TM^-S \,
)F_{S_1} \, ], \hfill (18.9)
\cr
\hspace*{1cm}
\min_Xr( \, A - BXC \, ) = r \left[ \begin{array}{c}  A  \\ C  \end{array}
\right] +  r[ \, A, \ B \, ] -r(M) + \min_Xr[ \,E_{T_1}(\, X + TM^-S \, )F_{S_1} \, ]. \hfill (18.10)
\cr}
$$ 
It is obvious that
$$
\max_X r[ \, E_{T_1}(\, X + TM^-S \, )F_{S_1} \, ] = \max_Yr( \, E_{T_1}YF_{S_1} \, ) = \min \{ \ r(E_{T_1}),  \ \ \ r(F_{S_1})  \ \}, \eqno (18.11)
$$ 
and the matrix satisfying it can be written as (18.7). According to (1.2) and (1.3), 
we find that
$$
 r(E_{T_1}) = k - r(T_1) = k - r(TF_M) = k -  r \left[ \begin{array}{c} M \\ T \end{array} \right] + r(M )  =  r(M ) -r \left[ \begin{array}{c} A \\ C \end{array} \right], 
$$ 
$$
 r(F_{S_1}) = l - r(S_1) = k - r(E_MS) = k -  r[\, M, \ S \,] + r(M) = r(M) - r[\, A, \ B \,].
$$ 
Putting both of them in (18.11) and then  (18.11) in  (18.9) yields (18.5). The 
results in  (18.6) and (18.8) are direct consequences of (18.10). \qquad$ \Box$. 
    
\medskip

The maximal and the minimal ranks of $A - BXC $ with respect to $ X$ can also be determined through 
the rank identity (18.3). In that case, the matrix $ X $ satisfying (18.5) and (18.6) can 
respectively be determined by the expression matrix $ E_{A_2}AF_{A_1} - E_{A_2}BXCF_{A_1},$
where the corresponding  matrix equation $E_{A_2}BXCF_{A_1} =E_{A_2}AF_{A_1}$ is consistent.    

\medskip

\noindent {\bf Corollary 18.5.}\, {\em Let $ p( X ) = A - BXC$ be given over ${\cal F}$ with 
$ B \neq 0$ and $ C \neq 0.$  Then 

{\rm (a)}\, The rank of $ A- BXC $ is invariant with respect to the choice of $ X $ if and only if 
$$
R \left[ \begin{array}{c} B \\ 0  \end{array} \right] \subseteq R \left[ \begin{array}{c}  A  \\ C  \end{array} \right]  \ \  \ or  \  \ \ 
 R([ \, C, \, 0 \,]^T) \subseteq R([ \, A,  \ B \,]^T).  \eqno (18.12)
$$

{\rm (b)}\, The range $ R( \, A - BXC  \, ) $ is invariant  with respect to the choice of $ X $ if and only if
$$ 
R \left[ \begin{array}{c}  B  \\ 0  \end{array} \right] \subseteq  R\left[ \begin{array}{c}  A  \\ C  \end{array} \right].  \eqno (18.13)
$$ 

{\rm (c)}\, The range  $R[ \,( \,A- BXC \,)^T\,]$  is invariant with respect to the choice of $ X $ if and only if 
$$ 
R(\,[  \, C ,\ 0 \, ]^T \,) \subseteq R(\,[ \, A , \ B \, ]^T \,). \eqno (18.14)
$$

{\rm (d)}\, The rank of $ A- BXC $ is invariant with respect to the choice of $ X $ if and only if the range $ R( \, A - BXC  \, ) $ is invariant  with respect to the choice of $ X $  or 
the range $R[ \,( \,A- BXC \,)^T\,]$ is invariant with respect to the choice of $ X $.} 

\medskip

\noindent {\bf Proof.}\, From  (18.5) and (18.6), we obtain 
$$  
\max_Xr(\, A - BXC \,) - \min_X r(\, A - BXC\, ) 
= \min \left\{ \, r \left[ \begin{array}{cc}  A  & B  \\ C &  0
\end{array} \right] - r\left[ \begin{array}{c}  A  \\ C
\end{array} \right], \ \  r \left[ \begin{array}{cc}  A  & B  \\ C &
0 \end{array} \right] - r[\, A, \ B \, ]  \, \right \}.
$$
Let its right-hand side be zero, then we get  (18.12). To show Part (b), we use
a fundamental fact that two matrices $P $ and $Q$ have the same range,
i.e., $ R(P) = R(Q),$ if and only if $ r[ \, P, \  Q
\, ]$ $= r(P) = r(Q)$.  Applying this fact to $ A- BXC$,
we know that the range $ R( A - BXC ) $ is invariant with respect
to the choice of $ X $ if and only if
$$ 
 r[ \,  A - BXC, \  A - BYC  \, ] =  r(\, A - BXC \, ) =  r( \, A - BYC \, )
 \eqno (18.15)
$$ 
holds for all $ X $ and $ Y $. Obviously this equality holds for all $ X$
and $ Y $ if and only if
$$ 
 r( \, A - BXC \, ) = r( A ), \eqno (18.16) 
$$   
holds for all $ X $, and  
$$ 
r[ \, A - BXC, \  A - BYC \, ] =  r \left(   [  \, A, \ A \, ] - B[ \,  X, \
Y \, ]\left[ \begin{array}{cc} C  & 0  \\ 0 &  C \end{array} \right]
\ \right) = r( A )  \eqno (18.17)
$$ 
holds for all $ X $ and $ Y $. According to Part (a), the equality
(18.16) holds for all $ X $ if and only if  (18.12) holds, and the equality
(18.17) holds if and only if
$$ 
r \left[ \begin{array}{ccc}  A  & A & B  \\ C & 0 & 0 \\ 0 & C & 0 \end{array} \right]
=  r \left[ \begin{array}{cc}  A & A  \\ C  & 0 \\ 0 & C \end{array} \right] \ \ \ {\rm or}  \ \ \
  r \left[ \begin{array}{ccc}  A  & A & B  \\ C & 0 & 0 \\ 0 & C & 0 \end{array} \right]
  = r[ \, A, \ A, \  B \, ],
$$
that is, 
$$ 
r \left[ \begin{array}{cc}  A  & B  \\ C &  0 \end{array} \right]
=  r \left[ \begin{array}{c}  A  \\ C  \end{array} \right] \ \ \ {\rm or}  \ \ \
  r \left[ \begin{array}{cc}  A  & B  \\ C &  0 \end{array} \right]
  = r[ \, A, \  B \, ]
- r(C).   \eqno (18.18)  
$$ 
Note that $ B \neq 0$ and $  C \neq 0$. Thus combining  (18.12) with (18.18),
we know that  (18.15) holds if and only if 
$$ 
r \left[ \begin{array}{cc}  A  & B  \\ C &  0\end{array} \right] =
 r \left[ \begin{array}{c}  A  \\ C  \end{array} \right],
$$ 
which is equivalent to (18.13). Similarly we can show Part (c). Contrasting 
Parts (a)---(c) yields Part (b). \qquad $\Box $    

\medskip

\noindent {\bf Corollary 18.6.}\, {\em The matrix satisfying  {\rm (18.6)} is 
unique if and only if
$$ 
 r( B ) = k,   \ \ \  r( C ) = l ,  \ \ \  and \ \ \ r \left[ \begin{array}{cc}  A  & B  \\ C &  0 \end{array} \right]
= r \left[ \begin{array}{c}  A  \\ C  \end{array} \right] + r( B )
 =  r[ \, A, \ B \, ] + r( C ). \eqno (18.19)
$$ 
In that  case$,$ the unique matrix satisfying {\rm (18.6)} is
$$ 
X = - [\, 0, \ I_k \, ] \left[ \begin{array}{cc}  A  & B   \\ C  & 0  \end{array} \right]^-
 \left[ \begin{array}{c}  0  \\ I_l  \end{array} \right].  \eqno  (18.20) 
$$}
\noindent {\bf Proof.}\, The  matrix satisfying (18.6) is unique
if and only if the solution to the equation $E_{T_1}(\, X + TM^-S \, )F_{S_1} = 0$ 
 is unique, which is equivalent to
$$ 
r(E_{T_1}) = k  \ \ \ {\rm and}   \ \ \  r(F_{S_1}) = l.  \eqno (18.21)
$$ 
Recall that
$$
r(E_{T_1}) =  r(M ) - r \left[ \begin{array}{c} A \\ C \end{array} \right],  \ \ \ 
{\rm and}  \ \ \  r(F_{S_1}) = r(M) - r[\, A, \ B \,],
$$ 
and $r( B ) \leq k $ and $ r( C ) \leq l.$ Thus (18.21) is equivalent to
  (18.19), and the unique matrix is  (18.20). \qquad $\Box$.

\medskip

\noindent {\bf Corollary 18.7.}\,  {\em  The following four statements are equivalent$:$

{\rm (a)} \  $ \min_X r ( A - BXC ) = r ( A ). $ 
   
{\rm (b)} \ $ r \left[ \begin{array}{cc}  A  & B  \\ C &  0 \end{array}
\right] = r \left[ \begin{array}{c}  A  \\ C  \end{array} \right] +
 r[ \, A, \ B \, ] - r(A).$

{\rm (c)} \ $E_{T_1}TM^-S F_{S_1} = 0.$

{\rm (d)} \ $E_{A_2}AF_{A_1}  = 0.$

{\rm (e)} \  $ E_{C_1}CA^-BF_{B_1}  = 0,$  where $ A_1 = E_BA, \
A_2 = AF_C, \  B_1 = E_AB, \ C_1 = CF_A.$}

\medskip

\noindent {\bf Proof.} \ Follows immediately from the combination of
 (18.6), (18.1), (18.2) and (1.6).
 \qquad $\Box $ 

\medskip

In the remainder of this section, we present some equivalent statements for
the results in Theorem 18.4.

Suppose that $  B \in {\cal F}^{ m \times k}, \,
 C \in {\cal F}^{ l \times n },\,  P \in {\cal F}^{ s \times m}, \,
 Q \in {\cal F}^{ n \times t},$ and let $ \Theta $
  be the  matrix set
$$ 
\Theta = \{ \, Z \in {\cal F}^{ m \times n} \ | \ R(Z)
\subseteq R(B) \ \ {\rm and } \ \  R(Z^T) \subseteq R
(C^T)\, \}.  \eqno (18.22)
$$ 

Then we have the following results. 

\medskip

\noindent {\bf Theorem 18.8.}\, {\em  Let $ A \in {\cal F}^{ m \times n}$
be given and $ \Theta $ be defined in {\rm (18.21)}.  Then  

{\rm (a)}\, The maximal rank of $ A - Z $ subject to $ Z \in \Theta $ is  
$$\displaylines{
\hspace*{2cm}
\max_{Z \in \Theta} r( A - Z ) = \min \left\{ \, r[ \, A, \  B  \,],
\qquad r \left[ \begin{array}{c}  A  \\ C  \end{array} \right] \, \right \},
\hfill (18.23)
\cr}
$$ 
and the matrix $Z$ satisfying  {\rm (18.23)} can be written in the form 
$$
\displaylines{
\hspace*{2cm} 
Z  = -[ \,0 , \ B \,]\left[ \begin{array}{cc}  A  & B \\ C & 0 \end{array} 
\right]^-\left[ \begin{array}{c} 0  \\ C  \end{array} \right] -  BUC,  \hfill (18.24)
\cr}
$$  
where $U$ is chosen such that $r(E_{T_1}UF_{S_1}) = \min \{ \, r(E_{T_1}),  \ \ \ r(F_{S_1})  \, \},$
where $M, \ S,  \ T, \ S_1$ and $T_1$ are as in {\rm (18.1)}.

{\rm (b)}\, The minimal rank of $ A - Z $ subject to $ Z \in \Theta $ is  
$$ \displaylines{
\hspace*{2cm} 
\min_{Z \in \Theta} r( A - Z ) =  r[ \, A, \ B \, ] +
r \left[ \begin{array}{c}  A  \\ C  \end{array} \right] -
r \left[ \begin{array}{cc}  A  & B  \\ C &  0 \end{array} \right],
\hfill (18.25)
\cr}
$$ 
and the general expression of the matrix $Z$ satisfying  {\rm (18.24)}
can be written as
$$\displaylines{
\hspace*{2cm} 
Z  = -[ \,0 , \ B \,]\left[ \begin{array}{cc}  A  & B \\ C & 0 \end{array} 
\right]^-\left[ \begin{array}{c} 0  \\ C  \end{array} \right] + 
BT_1VC + BWS_1C, \hfill (18.26)
\cr}
$$
where $ V $ and $W$ are arbitrary.}

\medskip

\noindent {\bf Proof.}\, From the structure  of $ \Theta $ in (18.22)
we easily see that  $ \Theta $ can equivalently
be expressed as 
$$ 
\Theta = \{ \, Z = BXC \ | \ X \in {\cal F}^{k \times l} \, \}.
$$ 
Thus the rank of $ A - Z $ subject to $ Z \in \Theta $ can be written as 
$$
r( A - Z ) =  r( A - BXC )= r[ \, A, \  B \, ] +  r \left[ \begin{array}{c}  A \\ C \end{array}
\right] - r \left[ \begin{array}{cc}  A  & B  \\ C &  0 \end{array} \right]
+ r[ \,E_{T_1}(\, X + TM^-S \, )F_{S_1} \, ]. \eqno (18.27)
$$
In this case, applying Theorem 18.4 to this equality, we
 obtain the desired results in the
theorem.  \qquad  $\Box $    

\medskip

The matrix $ Z \in \Theta$ satisfying (18.25) is well known as a {\it shorted matrix} of
$ A $ relative to $ \Theta$. Thus (18.26) is in fact the general
 expression  of shorted matrices  of $ A $ relative to $ \Theta $.
 One  of the most important aspects on shorted matrices is concerning  
their uniqueness, which has been well examined by lots of
 authors (see, e.g., \cite{Ande}, \cite{Car}, \cite{Mi5}, \cite{MP}). Now from the
 general result in Theorem 18.6 and 18.8(b) and we easily get the following known
 result.  

\medskip
  
\noindent {\bf  Theorem  18.9}\cite{Mi5}.\, {\em  Let
$ A \in {\cal F}^{m \times n}$ be given and $ \Theta $ be defined
in  {\rm (18.22)}.  Then the shorted matrix of $ A $ relative to $ \Theta $
is unique if and only if
$$ \displaylines{
\hspace*{2cm}
r \left[ \begin{array}{cc}  A  & B  \\ C &  0 \end{array} \right]
= r \left[ \begin{array}{c}  A  \\ C  \end{array} \right] + r( B )  =
 r[ \, A, \ B \, ] + r( C ). \hfill (18.28)
\cr}
$$ 
In that case$,$ the unique shorted matrix is 
$$\displaylines{
\hspace*{2cm} 
Z  = -[ \,0 , \ B \,]\left[ \begin{array}{cc}  A  & B \\ C & 0 \end{array} 
\right]^-\left[ \begin{array}{c} 0  \\ C  \end{array} \right],  \hfill (18.29)
\cr}
$$
and this matrix is invariant
with respect to the choice of the inner inverse in it. }

\markboth{YONGGE  TIAN }
{19. EXTREME RANKS OF $A - B_1 X_1 C_1 - B_2 X_2 C_2$ }

\chapter{Extreme  ranks of $A - B_1 X_1 C_1 - B_2 X_2 C_2$}

\noindent In order to find the maximal and the minimal ranks of $ D - CA^-B$ with respect to $ A^-$, we need to 
know  maximal and minimal ranks of 
$$\displaylines{
\hspace*{2cm} 
 p( X_1, \,  X_2 ) = A - B_1 X_1 C_1 - B_2 X_2 C_2 \hfill (19.1)  
\cr}
$$  
under the two conditions 
$$\displaylines{
\hspace*{2cm} 
 R(B_1) \subseteq R(B_2)  \ \ {\rm and}  \ \
 R(C_2^T) \subseteq R(C_1^T),  \hfill (19.2)
\cr}
$$ 
where $ A,  \, B_1,  \, B_2, \,  C_1 $ and $ C_2 $ are given,  $ X_1 $ and $ X_2 $ are two independent variant matrices  over ${\cal F}$. 

\medskip

\noindent {\bf Theorem 19.1.}\, {\em Let $ p( X_1, \, X_2 )$ be given by
 {\rm (19.1)} and {\rm (19.2)}. Then the maximal rank of $ p( X_1, \, X_2 )$
with respect to $ X_1 $  and $ X_2 $ is
$$ \displaylines{
\hspace*{2cm} 
\max_{X_1, \, X_2} r[ \, p( X_1, \, X_2 )\,] = \min  \left\{ \, r[ \, A, \  B_2 \, ], \ \ 
r \left[ \begin{array}{c}  A  \\ C_1  \end{array} \right], \ \ r \left[ \begin{array}{cc}  A   
& B_1\\ C_2   & 0 \end{array} \right] \, \right\}. \hfill (19.3)
\cr}
$$ }
{\bf Proof.}\, Applying  (18.5) to  $ p(X_1, \, X_2)$ in (19.1)
we first obtain
\begin{eqnarray*}
\max_{X_2} r[\, p( X_1, \, X_2 )\,] & = & \min  \left\{ \, r( \, A - B_1X_1C_1, \  B_2 \, ), \ \ 
r \left[ \begin{array}{c}  A - B_1X_1C_1 \\ C_2  \end{array} \right] \, \right\} \\
& = & \min  \left\{ \, r[\, A,  \  B_2 \, ], \ \ \ r \left[ \begin{array}{c}  A - B_1X_1C_1 \\ C_2  \end{array} \right] \, \right\}. 
\end{eqnarray*} 
Next applying  (18.5) to $ \left[ \begin{array}{c}  A - B_1X_1C_1 \\ C_2
\end{array} \right], $ we obtain
\begin{eqnarray*} 
\max_{X_1} r \left[ \begin{array}{c}  A - B_1X_1C_1 \\ C_2  \end{array} \right] =  \max_{X_1} r \left( \,  \left[ \begin{array}{c}  A \\ C_2  \end{array} \right] -
 \left[ \begin{array}{c}  B_1 \\ 0  \end{array} \right] X_1C_1  \, \right) =  \min  \left\{ \, r \left[ \begin{array}{c}  A  \\ C_1  \end{array} \right], \ \ \   r \left[ \begin{array}{cc}  A   & B_1\\ C_2   & 0 \end{array} \right] \, \right\}.
\end{eqnarray*} 
Combining the above two results  yields (19.3). \qquad $ \Box$ 

\medskip

\noindent {\bf Theorem 19.2.}\, {\em Let $ p( X_1, \, X_2 )$ be given by  {\rm (19.1)} and 
{\rm (19.2)}. Then the minimal rank of $ p( X_1, \, X_2 )$ with respect to $ X_1 $  and $ X_2 $ is 
$$ 
\min_{X_1, \, X_2} r[ \, p( X_1, \, X_2 )\,] = r[ \, A, \  B_2 \, ] +
 r \left[ \begin{array}{c}  A  \\ C_1  \end{array} \right] + r \left[\begin{array}{cc}  A  
 & B_1 \\ C_2 & 0 \end{array} \right]- r \left[ \begin{array}{cc}  A   & B_1 \\ C_1   & 0 \end{array} \right]
- r \left[ \begin{array}{cc}  A   & B_2 \\ C_2   & 0 \end{array} \right]. 
   \eqno (19.4)
$$ }
{\bf Proof.}\,  Applying  (18.6) to $ p( X_1, \, X_2 )$ in (19.1) we
first obtain
\begin{eqnarray*}
\min_{X_2} r[\, p( X_1, \, X_2 )\,] & = & r[ \, A - B_1X_1C_1, \  B_2 \, ] + r \left[ \begin{array}{c}  A - B_1X_1C_1 \\ C_2  \end{array} \right]  - r \left[ \begin{array}{cc}  A - B_1X_1C_1  & B_2  \\ C_2  & 0  \end{array} \right] \\
& = &  r[ \, A, \ B_2 \, ] + r \left[ \begin{array}{c} A - B_1X_1C_1 \\ C_2  \end{array} \right]  - r \left[ \begin{array}{cc}  A & B_2  \\ C_2  & 0  \end{array} \right]. 
\end{eqnarray*} 
Next applying (18.6) to $ \left[ \begin{array}{c}  A - B_1X_1C_1 \\ C_2
\end{array} \right],$  we find
$$
\min_{X_1} r \left[ \begin{array}{c}  A - B_1X_1C_1 \\ C_2  \end{array} \right] 
 =  \min_{X_1} r \left(  \left[ \begin{array}{c}  A \\ C_2  \end{array} \right] -
 \left[ \begin{array}{c}  B_1 \\ 0  \end{array} \right] X_1C_1  \right) =  r \left[ \begin{array}{cc}  A  & B_1 \\ C_2  & 0 \end{array} \right] +  r \left[ \begin{array}{c}  A  \\ C_1  \end{array} \right]  - r \left[ \begin{array}{cc}  A   & B_1 \\ C_1   & 0 \end{array} \right].
$$
Combining the above two results yields (19.4). \qquad $ \Box$  

\medskip

The matrices $X_1$ and $ X_2$ satisfying  (19.3) and (19.4) can also be derived through 
the two expressions in  (18.7) and (18.8). But 
their expressions are somewhat complicated in form and are omitted them here.

Eq.\,(19.4) can also be written as 
\begin{eqnarray*}  
\min_{X_1, \, X_2} r[\, p( X_1 \, X_2 )\,] & = & \left( \, r[\, A, \  B_2 \, ] + r \left[ \begin{array}{c}  A  \\ C_1  \end{array} \right] -  r \left[ \begin{array}{cc}  A   & B_2 \\ C_1   & 0 \end{array} \right] \
 \right) \\
 & & + \left(\, r \left[ \begin{array}{cc}  A   & B_2 \\ C_1   & 0 \end{array} \right] + r \left[ \begin{array}{cc}  A  & B_1 \\ C_2  & 0 \end{array} \right] - r \left[ \begin{array}{cc}  A   & B_1 \\ C_1   & 0 \end{array} \right] - r\left[ \begin{array}{cc}  A   & B_2 \\ C_2   & 0 \end{array} \right] \, \right).
\end{eqnarray*}
It is easy to verify that under (19.2) the  two quantities in the parentheses on the right hand-side of the above equality are nonnegative. Thus the right hand-side of (19.4) is also nonnegative, although this is not evident from its expression. 

Some direct consequences of Theorems 19.1 and 19.2 are given below.  

\medskip

\noindent {\bf Corollary 19.3.}\, {\em  Let $ p(X_1, \, X_2)$ be given by
 {\rm (19.1)} and {\rm (19.2)}. Then the rank of $ p(X_1, \, X_2)$ is invariant with
respect to the choice of $ X_1 $ and $ X_2$ if and only if
$$ \displaylines{
\hspace*{2cm}
r \left[ \begin{array}{cc}  A  & B_1  \\ C_1 &  0 \end{array} \right] =  r \left[ \begin{array}{c}  A 
 \\ C_1  \end{array} \right] \ \ \  \  and  \ \ \ r \left[ \begin{array}{cc}  A  & B_2  \\ C_2 &  0 \end{array} \right] = r \left[ \begin{array}{cc}  A  & B_1  \\ C_2 &  0 \end{array} \right],  \hfill (19.5) 
\cr
\hspace*{0cm}
or \hfill
\cr
\hspace*{2cm}
r \left[ \begin{array}{cc}  A  & B_2  \\ C_2 &  0 \end{array} \right] = r[  \, A, \  B_2 \, ]  \ \ \  \  and  \ \ \  r \left[ \begin{array}{cc}  A  & B_1  \\ C_1 &  0 \end{array} \right] = r \left[ \begin{array}{cc}  A  & B_1  \\ C_2 &  0 \end{array} \right] , 
\hfill (19.6) 
\cr
\hspace*{0cm}
or\hfill
\cr
\hspace*{2cm} 
r \left[ \begin{array}{cc}  A  & B_1  \\ C_1 &  0 \end{array} \right] =  r \left[ \begin{array}{c}  A 
 \\ C_1  \end{array} \right] \ \ \  \  and  \ \ \ r \left[ \begin{array}{cc}  A  & B_2  \\ C_2 &  0 \end{array} \right] = r[  \, A, \  B_2 \,]. \hfill (19.7)
\cr}
$$ } 
{\bf Proof.}\, Combining (19.3) and (19.4), we obtain the following
$$
\displaylines{
\hspace*{2cm} 
\max_{X_1,\, X_2} r[\, p(X_1, \,X_2)\,] - \min_{X_1,\, X_2} r[\, p(X_1, \,X_2)\,] =
 \min \{ \, s_1, \ \ s_2, \ \ s_3 \, \}, \hfill
\cr}
$$
where
$$
\displaylines{
\hspace*{2cm} 
s_1 = r \left[ \begin{array}{cc}  A  & B_1  \\ C_1 &  0 \end{array} \right]  
+ r \left[ \begin{array}{cc} A  & B_2  \\ C_2 &  0 \end{array} \right]
- r \left[ \begin{array}{c}  A  \\ C_1  \end{array} \right] - r \left[ \begin{array}{cc}  A  & B_1  \\ C_2 &  0 \end{array} \right], \hfill
\cr
\hspace*{2cm} 
s_2 =  r \left[ \begin{array}{cc}  A  & B_1  \\ C_1 &  0 \end{array} \right]  
+ r \left[ \begin{array}{cc}  A  & B_2  \\ C_2 &  0 \end{array} \right]
-  r[ \, A, \  B_2 \, ] - r \left[ \begin{array}{cc}  A  & B_1  \\ C_2 &  0 \end{array} \right],  \hfill
\cr
\hspace*{2cm} 
s_3 = r \left[ \begin{array}{cc}  A  & B_1  \\ C_1 &  0 \end{array} \right]  
+ r \left[ \begin{array}{cc}  A  & B_2  \\ C_2 & 0 \end{array} \right]
- r \left[ \begin{array}{c}  A  \\ C_1  \end{array} \right] - r[ \, A, \ B_2 \,]. \hfill
\cr}
$$
Let the right-hand side of the above equality be zero. Then we obtain (19.5)---(19.7). \qquad  $\Box $ 

\medskip

\noindent {\bf Corollary 19.4.}\, {\em  Let $ p(X_1, \, X_2)$ be given by
 {\rm (19.1)} and {\rm (19.2)} with $ B_1 \neq 0$ and $ C_2 \neq 0$. Then 

{\rm (a)}\, The range $R[\, p(X_1, \, X_2)\,]$ is invariant  with respect to the choice of $ X_1 $ and $X_2$ if and only if
$$ 
r \left[ \begin{array}{cc}  A  & B_1  \\ C_1 &  0 \end{array} \right] =  r \left[ \begin{array}{c}  A 
 \\ C_1  \end{array} \right] \ \ \  \  and  \ \ \ r \left[ \begin{array}{cc}  A  & B_2  \\ C_2 &  0 \end{array} \right] = r \left[ \begin{array}{cc}  A  & B_1  \\ C_2 &  0 \end{array} \right].  \eqno (19.8) 
$$  

{\rm (b)}\, The range $R[\, p^T(X_1, \, X_2)\,]$ is invariant  with respect to the choice of 
$ X_1 $ and $X_2$ if and only if 
$$ 
r \left[ \begin{array}{cc}  A  & B_2  \\ C_2 &  0 \end{array} \right] = r[  \, A, \  B_2 \, ]  \ \ \  \  and  \ \ \ 
r \left[ \begin{array}{cc}  A  & B_1  \\ C_1 &  0 \end{array} \right] = r \left[ \begin{array}{cc}  A  & B_1  \\ C_2 &  0 \end{array} \right]. \eqno (19.9) 
$$} 
{\bf Proof.}\, It is obvious that the range $R[\, p(X_1,
\,X_2)\,]$ is invariant  with respect to the choice of $ X_1 $ and $X_2$ if
and only if
$$ \displaylines{
\hspace*{2cm} 
 r[ \,  p(X_1, \, X_2), \  p(Y_1, \, Y_2) \, ] =  r[\, p(X_1, \, X_2)\, ] =  r[\, p(Y_1, \, Y_2)\, ] =r(A) 
\hfill (19.10) 
\cr}
$$ 
holds for all $ X_1, \  X_2, \ Y_1$ and $ Y_2 $. By Corollary 19.3,
$ r[\, p(X_1, \,X_2)\, ] = r(A)$ holds  for all $ X_1, \  X_2 $ if and only
if one of (19.5)---(19.7) holds. On the other hand,
$$ \displaylines{
\hspace*{2cm} 
[ \,  p(X_1, \, X_2), \  p(Y_1, \, Y_2) \, ] = [\, A, \ A \, ] - B_1[\, X_1, \ Y_1 \, ] 
\left[ \begin{array}{cc}  C_1  & 0 \\ 0 & C_1 \end{array} \right] -  B_2[\, X_2, \ Y_2 \,] 
\left[ \begin{array}{cc}  C_2  & 0 \\ 0 & C_2 \end{array} \right]. \hfill
\cr}
$$   
Then according to Corollary 19.3, this expression satisfies (19.10) if and
only if
$$ \displaylines{
\hspace*{1cm}
r \left[ \begin{array}{cc}  A  & B_1  \\ C_1 &  0 \end{array} \right] =  r \left[ \begin{array}{c}  A 
 \\ C_1  \end{array} \right] \ \ \  \  {\rm and} \ \ \  r \left[ \begin{array}{cc}  A  & B_1  \\ C_2 &  0 \end{array} \right] = r \left[ \begin{array}{cc}  A  & B_2  \\ C_2 &  0 \end{array} \right],  \hfill (19.11) 
\cr}
$$ 
or 
$$\displaylines{
\hspace*{1cm} 
r \left[ \begin{array}{cc}  A  & B_2  \\ C_2 &  0 \end{array} \right] + r(C_2) = r[  \, A, \  B_2 \, ]  \ \ \  \ {\rm and}  \ \ \  r \left[ \begin{array}{cc}  A  & B_1  \\ C_2 &  0 \end{array} \right] + r(C_2) = 
r \left[ \begin{array}{cc}  A  & B_1  \\ C_1 &  0 \end{array} \right], \hfill (19.12) 
\cr}
$$ 
or 
$$ \displaylines{
\hspace*{1cm}
r \left[ \begin{array}{cc}  A  & B_1  \\ C_1 &  0 \end{array} \right] =  r \left[ \begin{array}{c}  A 
 \\ C_1  \end{array} \right] \ \ \  \  {\rm and}  \ \ \ r \left[ \begin{array}{cc}  A  & B_2  \\ C_2 &  0 \end{array} \right]  + r(C_2) = r[  \, A, \  B_2 \, ]. \hfill (19.13)
\cr}
$$
Contrasting (19.11)---(19.13) with (19.5)---(19.7) and noticing the
condition $ B_1 \neq 0$ and $ C_2 \neq 0$, we find that (19.10) holds if
and only if (19.11), i.e., (19.8) holds. Similarly we can show Part (b).
\qquad  $\Box $

\medskip

If one of  $ B_1, \, B_2, \, C_1 $ and $ C_2$ in  (19.1) is a null matrix,
then $ p(X_1, \, X_2)$ becomes an expression with a single variant matrix in it. 
In that case, the range invariance criterion is listed in Corollary 18.5(b) and (c). 

\medskip

\noindent {\bf Corollary 19.5.}\, {\em Let $ A \in {\cal F}^{m \times n},  \, B \in {\cal F}^{m \times k}$ and 
$ C \in {\cal F}^{l \times n}$ be given. Then    
$$\displaylines{
\hspace*{1.5cm}
\max_{X, \, Y} r( \, A - BX - YC \,)  = \min \left\{  m , \ \ \ n , \ \ \ r \left[ \begin{array}{cc}  A   & B \\ C  & 0 \end{array} \right] \right\}, \hfill (19.14)
\cr
\hspace*{1.5cm}
\min_{X, \, Y} r( \, A - BX - YC \,)  = r \left[ \begin{array}{cc}  A   & B  \\ C & 0 \end{array} \right]
- r(B) - r(C).  \hfill (19.15) 
\cr}
$$
A pair of matrices $ X $ and $ Y $ satisfying {\rm (19.5)} can be written as 
$$ \displaylines{
\hspace*{1.5cm}
X = B^-A + UC + (\, I_k - B^-B \,)V_1, \ \ \ Y = (\, I_m - BB^- \,)AC^- - BU + V_2(\, I_l - CC^- \,), \hfill (19.16)
\cr}
$$ 
where $ U, \, V_1 $ and $ V_2$ are arbitrary. } 

\medskip

\noindent {\bf Proof.}\, Eqs.\,(19.14) and (19.15) follow immediately from
 (19.3) and (19.4). Putting  (19.16) in $A - BX - YC$ yields 
$$ \displaylines{
\hspace*{1.5cm} 
A - BX - YC  = (\, I_m - BB^- \,)A(\, I_n - C^-C \,).  \hfill
\cr}
$$
Thus we have (19.15) by (1.4). \qquad  $ \Box$      

\markboth{YONGGE  TIAN }
{20. EXTREME RANKS OF $A - B_1XC_1$ SUBJECT TO $ B_2 XC_2 = A_2$ }

\chapter{Extreme ranks of $A_1 - B_1XC_1$ subject to $ B_2XC_2 = A_2$ }

\noindent  Based on the results in Chapter 19, we are now able to find the maximal and
the minimal ranks of $A_1 - B_1XC_1$ subject to a consistent linear matrix
equation $B_2XC_2 = A_2$. The corresponding results will widely be  used in the sequel.

\medskip

\noindent {\bf Theorem 20.1.}\, {\em Suppose that the matrix equation
 $B_2XC_2 = A_2$ is a consistent.  Then

{\rm (a)}\, The maximal rank of $ p(X) = A_1 - B_1XC_1$ subject to
$B_2XC_2 = A_2$ is
$$\displaylines{
\hspace*{1cm} 
\max_{B_2XC_2 = A_2}r[\,p( X)\,] = \min  \left\{ \,  r \left[ \begin{array}{ccc}  A_1  & 0  & B_1 \\  0 & -A_2   & B_2  
\\  C_1 & C_2   & 0 \end{array} \right]-r(B_2) -r(C_2),  \ \ r \left[ \begin{array}{c}  A_1  \\ C_1  \end{array} \right],  \ \   r[ \, A_1, \  B_1 \, ] \, \right\}. \hfill (20.1)
\cr}
$$  

{\rm (b)}\, The minimal rank of $ p(X) = A_1 - B_1XC_1$ subject to
$B_2XC_2 = A_2$ is
$$\displaylines{
\hspace*{1cm}
\min_{B_2XC_2 = A_2}r[\,p(X)\,]  \hfill
\cr
\hspace*{0.5cm} = r[\, A_1, \  B_1 \, ] + r \left[ \begin{array}{c}  A_1  \\
C_1  \end{array} \right]  - r \left[ \begin{array}{ccc}  A_1  &  B_1  & 0 \\
C_1 & 0  & C_2  \end{array} \right]  - r \left[ \begin{array}{cc}  A_1  &
B_1  \\ C_1 & 0  \\ 0 & B_2  \end{array} \right] + r \left[
\begin{array}{ccc}  A_1  & 0  & B_1 \\  0 & -A_2   & B_2  \\  C_1 & C_2
& 0 \end{array} \right].  \hfill (20.2)
\cr}
$$
} 
{\bf Proof.}\, Note from Lemma 18.1 that the general solution of the consistent linear matrix equation
 $ B_2XC_2 = A_2$ can be written as $ X = X_0 + F_{B_2}V + W E_{C_2}$, where
 $ X_0 = B_2^-A_2C^-_2,$  $ V$  and $ W $ are arbitrary. Putting it in
 $ p(X) = A_1 - C_1XB_1$ yields
$$ \displaylines{
\hspace*{2cm} 
p(X) = A - B_1 F_{B_2}VC_1 - B_1WE_{C_2}C_1, \hfill
\cr}
 $$ 
where $ A = A_1 - B_1X_0C_1$. Observe that $R(B_1 F_{B_2}) \subseteq R(B_1)$ and 
$R[(E_{C_2}C_1)^T] \subseteq R[(C_1)^T]$. Thus it follows by (19.3) and (19.4) that
$$ \displaylines{
\hspace*{1cm} 
\max_{B_2XC_2 = A_2}r[\,p( X)\,] \hfill
\cr
\hspace*{1cm} 
 =  \max_{V, \, W} r( \, A - B_1 F_{B_2}VC_1 - B_1WE_{C_2}C_1 \, ) \hfill
\cr
\hspace*{1cm} 
 =   \min  \left\{ \ r[ \, A, \  B_1 \, ], \ \ \ 
r \left[ \begin{array}{c}  A  \\ C_1  \end{array} \right],  \ \ \ r \left[ \begin{array}{cc}  A   & B_1F_{B_2} \\ E_{C_2} C_1   & 0 \end{array} \right] \ \right\}, \hfill
\cr
\hspace*{0cm} 
and \hfill
\cr
\hspace*{1cm}
\min_{B_2XC_2 = A_2}r[\,p( X)\,]  \hfill
\cr
\hspace*{1cm}
 = \min_{V, \, W} r( \, A - B_1F_{B_2}VC_1 - B_1WE_{C_2}C_1 \, ) \hfill
\cr
\hspace*{1cm}
 = r[ \, A, \  B_1 \, ] + \ r \left[ \begin{array}{c}  A  \\ C_1  \end{array} \right] +
 r \left[ \begin{array}{cc}  A   & B_1 F_{B_2} \\ E_{C_2} C_1   & 0
 \end{array} \right]  -  r \left[ \begin{array}{cc}  A   & B_1F_{B_2} \\
 C_1   & 0 \end{array} \right] - r \left[ \begin{array}{cc}  A   & B_1 \\
 E_{C_2}C_1   & 0 \end{array} \right]. \hfill
\cr}
$$
Simplifying the ranks of the block matrices by Lemma 1.1, we see that
$$
\displaylines{
\hspace*{1cm}
r[ \, A, \  B_1 \, ] = r[\, A_1 - B_1X_0C_1, \  B_1 \,] = r[ \, A_1, \  B_1 \, ],  \ \ \ 
r \left[ \begin{array}{c}  A  \\ C_1  \end{array} \right] = r \left[ \begin{array}{c} A_1 - B_1X_0C_1  \\ C_1  \end{array} \right] = r \left[ \begin{array}{c}  A_1  \\ C_1  \end{array} 
\right],  \hfill
\cr
\hspace*{1cm}
\begin{array}{ rcl}
r \left[ \begin{array}{cc}  A   & B_1F_{B_2} \\ E_{C_2} C_1   & 0 \end{array} \right] 
& = & r \left[ \begin{array}{ccc}  A_1 - B_1X_0C_1  & B_1 & 0 \\  C_1 & 0  & C_2 \\ 0 & B_2  & 0 \end{array} \right]  - r(B_2) -r(C_2)  \\ 
& = & r \left[ \begin{array}{ccc}  A_1  & B_1 & 0 \\  C_1 & 0  & C_2 \\ 0 & B_2  & -A_2 \end{array} \right]  - r(B_2) -r(C_2),   
\end{array} \hfill
\cr
\hspace*{1cm}
r \left[ \begin{array}{cc}  A   & B_1F_{B_2} \\ C_1   & 0 \end{array} \right] = r \left[ \begin{array}{cc}  A_1 - B_1X_0C_1  & B_1  \\  C_1 & 0  \\ 0 & B_2  \end{array} \right]  - r(B_2) = r \left[ \begin{array}{cc}  A_1  & B_1 \\  C_1 & 0   \\  0 & B_2  \end{array} \right]  - r(B_2),  
\hfill
\cr
\hspace*{1cm}
 r \left[ \begin{array}{cc}  A   & B_1 \\ E_{C_2} C_1   & 0 \end{array} \right] = r \left[ \begin{array}{ccc}  A_1 - B_1X_0C_1  & B_1 & 0 \\  C_1 & 0  & C_2  \end{array} \right]-r(C_2)  =r \left[ \begin{array}{ccc}  A_1  & B_1 & 0 \\  C_1 & 0  & C_2  \end{array} \right] - r(C_2). \hfill
\cr}
$$
Putting them in the above two rank equalities yields  (20.1) and (20.2).
\qquad  $ \Box$

\medskip

Eq.\,(20.2) can also be written as
$$
\displaylines{
\hspace*{1cm}
\min_{B_2XC_2 = A_2}r( \, A_1 - B_1XC_1\,) = \left( \, r[ \, A_1, \  B_1 \, ] + r \left[ \begin{array}{c}  A_1  \\ C_1  \end{array} \right] -  r \left[ \begin{array}{cc}  A_1   & B_1 \\ C_1   & 0 \end{array} \right] \, \right) \hfill
\cr
\hspace*{1.5cm}
+ \left( \, r \left[ \begin{array}{ccc}  A_1  & 0  & B_1 \\  0 & -A_2   & B_2  \\  C_1 & C_2   & 0 \end{array} \right] + \   r \left[ \begin{array}{cc}  A_1   & B_1 \\ C_1   & 0 \end{array} \right] - \ r \left[ \begin{array}{ccc}  A_1  &  B_1  & 0 \\ C_1 & 0  & C_2  \end{array} \right] - r \left[ \begin{array}{cc}  A_1  &  B_1  \\ C_1 & 0  \\ 0 & B_2  \end{array} \right]  \, \right), 
\hfill
\cr}
$$ 
and the two quantities in the parentheses on the right hand-side of the
above equality are nonnegative.

\medskip

Some direct consequences are given below.

\medskip
  
\noindent {\bf Corollary 20.2.}\, {\em Suppose that $B_1XC_1 = A_1$ and
$B_2XC_2 = A_2$ are consistent$,$ respectively. Then
$$\displaylines{
\hspace*{0cm}
 \max_{B_2XC_2 = A_2}r( \, A - B_1XC_1 \,)  =  \min \left\{ 
 r \left[ \begin{array}{ccc} A_1  & 0  & B_1 \\  0 & -A_2   & B_2  \\
 C_1 & C_2   &  0 \end{array} \right]-r(B_2) -r(C_2),  \ \ r( C_1 ), \ \
 r( B_1 ) \right\}, \hfill (20.3)
\cr
\hspace*{0cm}
and \hfill
\cr
\hspace*{0cm}
 \min_{B_2XC_2 = A_2}r(\, A - B_1XC_1 \,)  = r \left[ \begin{array}{ccc}
 A_1  & 0  & B_1 \\  0 & -A_2   & B_2  \\  C_1 & C_2   & 0 \end{array}
 \right] - r \left[ \begin{array}{c}  B_1  \\ B_2  \end{array} \right]
 - r[ \, C_1, \ C_2 \, ].  \hfill (20.4)
\cr}
$$}
{\bf Proof.}\, The consistency of $B_1XC_1 = A_1$ implies that 
$ R(A_1)\subseteq R(B_1)$ and  $ R(A_1^T)\subseteq R(C_1^T),$  the consistency
 of $B_2XC_2 = A_2$ implies that  $ R(A_2) \subseteq R(B_2),$  and  $R(A_2^T) 
\subseteq R(C_2^T).$ In that
 case, (20.1) and (20.2) simplify to (20.3) and (20.4).   \qquad  $ \Box$

\medskip
 
Notice a simple fact that the pair of matrix equations $B_1XC_1 = A_1$ and
$B_2XC_2 = A_2$ have a common solution if and only if $B_1XC_1 = A_1$  and
$B_2XC_2 = A_2$ are consistent, respectively, and
$$ \displaylines{
\hspace*{2cm} 
\min_{B_2XC_2 = A_2}r( \, A_1 - B_1XC_1 \,) =
\min_{B_1XC_1 = A_1}r( \,  A_2 - B_2XC_2 \, ) = 0. \hfill
\cr}
$$ 
We immediately find from (20.4) the following well-known results. 

\medskip

\noindent {\bf Corollary 20.3}\cite{Mi6}\cite{Van}.\, {\em The pair of matrix
equations  $B_1XC_1 = A_1$ and $B_2XC_2 = A_2$ have a common solution if
and only if $B_1XC_1 = A_1$ and $B_2XC_2 = A_2$ are consistent$,$
respectively$,$ and
$$ \displaylines{
\hspace*{2cm} 
r \left[ \begin{array}{ccc} A_1  & 0  & B_1 \\  0 & -A_2   & B_2  \\
C_1 & C_2   & 0 \end{array} \right]
= r \left[ \begin{array}{c}  B_1  \\ B_2  \end{array} \right] +
r[ \, C_1, \ C_2 \, ]. \hfill
\cr}
$$ }
{\bf Corollary 20.4.}\, {\em Suppose that the pair of matrix
equations  $B_1XC_1 = A_1$ and  $B_2XC_2 = A_2$ are consistent$,$
respectively$,$ and denote their solution sets by
$$ \displaylines{
\hspace*{2cm} 
\Omega_1 = \{ \, X  \, | \, B_1XC_1 = A_1 \,  \} \ \ and  \ \  \Omega_2 = \{ \, X  \, | \, B_2XC_2 = A_2 \,  \}. \hfill
\cr}
$$ 
Then 

{\rm (a)}\, $ \Omega_2 \subseteq  \Omega_1$ holds if and only if $ B_1 = 0$ or $C_1 = 0$ or $ 
r \left[ \begin{array}{ccc} A_1  & 0  & B_1 \\  0 & -A_2   & B_2  \\  C_1 & C_2   & 0 \end{array} \right] 
= r(B_2) + r(C_2). $  

{\rm (b)} Under $ B_i \neq 0 $ and $ C_i \neq 0, \, i = 1, \, 2, $  the two equations  $B_1XC_1 = A_1$ and  $B_2XC_2 = A_2$ have the same solution set$,$ i.e.$,$ $\Omega_1 = \Omega_2,$ if and only if 
$$\displaylines{
\hspace*{2cm} 
r \left[ \begin{array}{ccc} A_1  & 0 & B_1 \\  0 & -A_2   & B_2  \\  C_1 & C_2   & 0 \end{array} \right] = r(B_1) + r(C_1) =  r(B_2) + r(C_2). \hfill
\cr}
$$} 
\hspace*{0.4cm} Another result related to the pair of matrix equations $B_1XC_1 = A_1$ and
$B_2XC_2 = A_2$ is given  below, which was presented by the  author in \cite{Ti7}. 

\medskip
 
\noindent {\bf Corollary 20.5.}\, {\em Suppose that $B_1X_1C_1 = A_1$ and  $B_2X_2C_2 = A_2$ 
are consistent$,$ respectively$,$ where $ X_1$ and $ X_2 $ have the same size. Then  
$$
 \min_{ \begin{array}{c}B_1X_1C_1 = A_1 \\ B_2X_2C_2 = A_2 \end{array}} r(\, X_1 - X_2 \, ) 
= r \left[ \begin{array}{ccc} A_1  & 0  & B_1 \\ 0 & -A_2   & B_2  \\
C_1 & C_2   & 0 \end{array} \right] - r \left[ \begin{array}{c}  B_1  \\ B_2
  \end{array} \right] - r[ \, C_1, \ C_2 \, ].  \eqno (20.5)
$$ } 
\hspace*{0.4cm} Finally we present two results on  rank invariance and range invariance of
$ A_1 - B_1XC_1$ subject to  $B_2XC_2 = A_2$.

\medskip

\noindent {\bf Theorem 20.6.}\, {\em Suppose that $B_2XC_2 = A_2$ is  consistent. Then the rank of $ A_1 - B_1XC_1$
 is invariant subject to  $B_2XC_2 = A_2$ if and only if 
$$\displaylines{
\hspace*{1cm}
r \left[ \begin{array}{cc}  A_1  &  B_1  \\ C_1 & 0 \\ 0 & B_2  \end{array} \right] = 
 r \left[ \begin{array}{c}  A_1  \\ C_1  \end{array} \right] + r(B_2) \ \ and  \ \
 r \left[ \begin{array}{ccc}  A_1  & 0  & B_1 \\  0 & -A_2   & B_2  \\  C_1 & C_2   & 0 \end{array} \right] 
 = r \left[ \begin{array}{ccc}  A_1  &  B_1  & 0 \\ C_1 & 0  & C_2  \end{array} \right] +
 r(B_2), \hfill
\cr
\hspace*{0cm}
or \hfill
\cr
\hspace*{1cm}
r \left[ \begin{array}{ccc}  A_1  &  B_1  & 0 \\ C_1 & 0 & C_2  \end{array} \right] =r[\, A_1, \ B_1 \, ] + r(C_2) \ \  and \ \  r \left[ \begin{array}{ccc}  A_1  & 0 & B_1 \\  0 & -A_2   & B_2  \\  C_1 & C_2   & 0 \end{array} \right]  = r \left[ \begin{array}{cc}  A_1  &  B_1  \\ C_1 & 0  \\ 0 & B_2  \end{array} \right] + r(C_2),  \hfill
\cr
\hspace*{0cm}
or \hfill
\cr
\hspace*{1cm} 
r \left[ \begin{array}{cc}  A_1  &  B_1  \\ C_1 & 0  \\ 0 & B_2  \end{array} \right] = 
 r \left[ \begin{array}{c}  A_1  \\ C_1  \end{array} \right] + r(B_2) \ \ and \ \ 
r \left[ \begin{array}{ccc}  A_1  &  B_1  & 0 \\ C_1 & 0  & C_2  \end{array} \right] =r[ \, A_1, \ B_1 \,] + r(C_2). \hfill
\cr}
$$ }  
{\bf Proof.}\, It is obvious that the rank of $ A_1 - B_1XC_1$ is invariant subject to $ B_2XC_2 = A_2 $ if and only if 
$$ \displaylines{
\hspace*{2cm} 
\max_{B_2XC_2 = A_2}r( \, A_1 - B_1XC_1 \,) = \min_{B_2XC_2 = A_2}r( \, A_1 - B_1XC_1 \,). \hfill
\cr}
$$
Applying Theorem 20.1 to it produces the desired result in the 
theorem.  \qquad  $ \Box$ 

\medskip

\noindent {\bf Theorem 20.7.}\, {\em Suppose that $B_2XC_2 = A_2$ is  consistent with $ B_1F_{B_2} \neq 0 $ and
 $ C_1E_{C_2} \neq 0.$   Then

{\rm (a)}\, The range $ R(A_1 - B_1XC_1)$ is invariant subject to  $B_2XC_2 = A_2$ if and only if
$$
r \left[ \begin{array}{cc}  A_1  &  B_1  \\ C_1 & 0  \\  0 & B_2  \end{array} \right] = 
 r \left[ \begin{array}{c}  A_1  \\ C_1  \end{array} \right] + r(B_2) \ \ and  \ \ r \left[ \begin{array}{ccc}  A_1  & 0   & B_1 \\  0  & -A_2   & B_2  \\  C_1 & C_2   & 0  \end{array} \right] 
 = r \left[ \begin{array}{ccc}  A_1  &  B_1  & 0  \\ C_1 & 0   & C_2  \end{array} \right] + r(B_2).
$$ 

{\rm (b)}\, The range $R[(A_1 - B_1XC_1)^T]$ is invariant subject to  $B_2XC_2 = A_2$ if and only if
$$ 
r \left[ \begin{array}{ccc}  A_1  &  B_1  & 0  \\ C_1 & 0   & C_2  \end{array} \right] =r[ \, A_1, \ B_1 \,] + r(C_2)  \ and  \  r \left[ \begin{array}{ccc}  A_1  & 0   & B_1 \\  0  & -A_2   & B_2  \\  C_1 & C_2   & 0  \end{array} \right]  = r \left[ \begin{array}{cc}  A_1  &  B_1  \\ C_1 & 0   \\ 0  & B_2  \end{array} \right] + r(C_2).
$$  
}
{\bf Proof.}\, Follows from  Theorem 19.4.  \qquad  $ \Box$  

\markboth{YONGGE  TIAN }
{21. EXTREME RANKS OF THE SCHUR COMPLEMENT $D - CA^- B$ }

\chapter{Extreme ranks of the Schur complement $ D - CA^- B$}

\noindent With the proper background of rank formulas presented in the previous chapter, we are now 
able to systematically deal with ranks of generalized Schur complements and  various related  topics. 
As is well known, for a given block matrix $ M = \left[ \begin{array}{cc}  A  &  B  \\ C & D  \end{array} \right] $ 
 over an arbitrary field ${\cal F }$, where $ A, \ B, \ C $ and $ D $ are $ m\times n$,  $ m\times k$,  $ l\times n$ and
 $ l\times k$ matrices, respectively, a generalized  Schur complement of $ A $ in $ M$ is defined to be  
$$ 
 S_A = D - CA^-B, \eqno(21.1) 
$$
where $ A^-$ is an inner inverse of $A$, i.e.,  $ A^- \in \{ \, X \, | \, AXA = A  \, \}.$ 
As one of the most important matrix expressions in matrix theory, there have been many results in 
the literature on generalized Schur complements and their applications (see, e.g., \cite{Ando, Bh, BCHM, Car, 
CHM, Fi, MaS2, MP, Ou, St}). Some of the work focused on equalities and inequalities for ranks of 
generalized Schur complements. The two rank well-known inequalities (see \cite{BCHM, Car, MaS2}) related to 
the Schur complement $ S_A$ are given by  
$$\displaylines{
\hspace*{2cm}
 r \left[ \begin{array}{cc}  A  & B  \\ C & D  \end{array} \right] \geq r(A) + r( \, D - CA^-B \, ), 
\hfill (21.2)
\cr
\hspace*{2cm}
 r \left[ \begin{array}{cc}  A  & B  \\ C & D  \end{array} \right] \leq r \left[ \begin{array}{cc} 
 A  & B  \\ C & CA^-B  \end{array} \right]+ r( \, D - CA^-B \, ). 
\hfill (21.3)
\cr}
$$    
Both of them in fact give upper and lower bounds for the rank of the 
the Schur complement $ D - CA^-B$, but they are not, in general,  the maximal and the minimal ranks of 
$D - CA^-B$ with respect to $ A^-$. Note that $A^-$ is in fact a solution of the
matrix equation $ AXA = A$.  Thus  Schur complement $ D - CA^-B$ may be regarded as a matrix expression 
$ D - CXB$, where $ X $ is a solution of the  matrix equation $ AXA = A$. In that case, applying the rank formulas
in Chapter 20,  we can simply establish the following.

\medskip

\noindent {\bf Theorem 21.2.}\, {\em Let $ S_A =  D - CA^-B$ be given by
  {\rm (21.1)}. Then

{\rm (a)}\, The maximal rank of $ S_A $ with respect to $ A^-$ is
$$ \displaylines{
\hspace*{1cm}
\max_{A^-}r( \,D - CA^-B \, )  = \min  \left\{ r[ \, C, \ D \,] ,
\ \ \ r \left[ \begin{array}{c}
 B  \\ D  \end{array} \right],  \ \ \  r \left[ \begin{array}{cc}  A  & B \\  C & D   \end{array} \right]
- r(A)  \right\}. \hfill (21.4)
\cr}
$$  

{\rm (b)}\, The minimal rank of $ S_A $ with respect to $ A^-$ is
$$
\min_{A^-}r( \,D - CA^-B \, ) = r(A) + r[\, C, \ D \,] + r \left[ \begin{array}{c} B  \\ D  \end{array} \right] + r\left[ \begin{array}{cc}  A  &  B  \\ C & D  \end{array} \right]  - \  r \left[ \begin{array}{ccc}  A  & 0  & B \\ 0   & C  & D  \end{array} \right] - 
 r \left[ \begin{array}{cc}  A  &  0   \\ 0  & B  \\ C & D  \end{array}
 \right]. \eqno (21.5)
$$ 
} 
{\bf Proof.}\, It is quite obvious that
$$
 \max_{A^-}r( \,D - CA^-B \, ) = \max_{AXA= A}r( \, D - CXB \, ),  \ \ \ \ 
\min_{A^-}r( \,D - CA^-B \, )  = \min_{AXA= A}r( \,D - CXB \, ). 
$$ 
Thus we obtain (21.4) and (21.5) by Theorem 20.1.  \qquad  $ \Box$ 

\medskip 

Eq.\,(21.5) can also be written as 
$$
\displaylines{
\hspace*{0.5cm}
\min_{A^-}r( \,D - CA^-B \, ) = \left( r[ \,C, \  D \, ] + r \left[ \begin{array}{c} B \\ D  \end{array} \right] -  r \left[ \begin{array}{cc} 0   & B \\ C   & D \end{array} \right] \right) \hfill
\cr
\hspace*{3.6cm}
 + \left( r \left[ \begin{array}{cc}  A  & B \\ C & D  \end{array} \right] + 
  r \left[ \begin{array}{cc} 0    & B \\ C & D \end{array} \right] + r( A )  - \ r \left[ \begin{array}{ccc}  A  &   0  & B \\ 0  & C  & D  \end{array} \right] -  r \left[ \begin{array}{cc}  A  &  0   \\ 0  & B  \\ C & D  \end{array} \right]  \right), \hfill 
\cr}
$$
and the two quantities in the parentheses on the right hand-side of the above equality are nonnegative. 

The two formulas in (21.4) and (21.5) can further simplify when $ A,$ 
$ B,$ $ C $ and $ D $ satisfy some conditions, such as, $R(D)
 \subseteq R(C)$ and $ R(D^T) \subseteq R(B^T);$ 
 $ R(D) \cap R(C) = \{ 0 \}$ and $ R(D^T) \cap R(B^T)
 = \{ 0 \};$ $R(C) \subseteq R(D)$ and $ R(B^T)
 \subseteq R(D^T).$ The reader can easily list the corresponding  results.

\medskip
  
\noindent {\bf Corollary 21.3.}\, {\em The rank of $ D - CA^-B $ is invariant  with respect to the choice of
 $ A^-$  if and only if 
$$\displaylines{
\hspace*{2cm}
r \left[ \begin{array}{ccc}  A  & 0  & B \\ 0  & C  & D  \end{array} \right] = r\left[ \begin{array}{cc}  A  &  B  \\ C & D  \end{array} \right]  \ \ and  \ \ r \left[ \begin{array}{cc}  A  &  0  \\ 0 & B  \\ C & D  \end{array} \right] = r\left[ \begin{array}{c} B  \\ D  \end{array} \right] + r(A), \hfill
\cr
\hspace*{0cm}
or \hfill
\cr
\hspace*{2cm}
r \left[ \begin{array}{ccc}  A  & 0  & B \\ 0  & C  & D  \end{array} \right] = r[ \, C, \ D \, ] + r(A) \ \ and \ \   r \left[ \begin{array}{cc}  A  &  0  \\ 0 & B  \\ C & D  \end{array} \right] = r\left[ \begin{array}{cc}  A  &  B  \\ C & D  \end{array} \right], \hfill
\cr
\hspace*{0cm}
or \hfill
\cr
\hspace*{2cm}
r \left[ \begin{array}{ccc}  A  & 0  & B \\ 0  & C  & D  \end{array} \right] = r[ \, C, \ D \, ] + r(A) \ \ and \ \   r \left[ \begin{array}{cc}  A  &  0  \\ 0 & B  \\ C & D  \end{array} \right] =r\left[ \begin{array}{c} B  \\ D  \end{array} \right] + r(A). \hfill
\cr}
$$ } 
{\bf Proof.}\, It is obvious that the rank of $ D - CA^-B $ is invariant with respect to the choice of $ A^-$ if and only if 
$$ \displaylines{
\hspace*{2cm}
\max_{A^-}r( \, D - CA^-B \,) = \min_{A^-}r( \, D - CA^-B \,). \hfill
\cr}
$$
Applying Theorem 20.6 to it leads to the desired result in the corollary.  \qquad  $ \Box$  

\medskip

\noindent {\bf Corollary 21.4.}\, {\em  Let $ S_A $ be given by  {\rm (21.1)} with $ E_AB \neq 0$ and 
$ CF_A \neq 0 $. 

 {\rm (a)}\, The range  $R( \,D - CA^-B \,) $ is invariant  with respect to the choice of
 $ A^-$  if and only if 
$$\displaylines{
\hspace*{2cm}
r \left[ \begin{array}{ccc}  A  & 0 & B \\ 0  & C  & D  \end{array} \right] = r\left[ \begin{array}{cc}  A  &  B  \\ C & D  \end{array} \right]  \ \ and  \ \ r \left[ \begin{array}{cc}  A  &  0  \\ 0 & B  \\ C & D  \end{array} \right] = r\left[ \begin{array}{c} B  \\ D  \end{array} \right] + r(A). \hfill
\cr}
$$

{\rm (b)}\, The range $R[ \,( \,D - CA^-B \,)^T \,] $ is invariant  with respect to the choice of
 $ A^-$  if and only if 
$$\displaylines{
\hspace*{2cm}
r \left[ \begin{array}{ccc}  A  & 0  & B \\ 0  & C  & D  \end{array} \right] = [\, C,  \ D \, ] + r(A ) \ \ and  \ \ r \left[ \begin{array}{cc}  A  &  0  \\ 0 & B  \\ C & D  \end{array} \right] = r\left[ \begin{array}{cc}  A  &  B  \\ C & D  \end{array} \right].  \hfill
\cr}
$$ 
} 
{\bf Proof.} \ Follows from  Theorem 20.7.  \qquad  $ \Box$

\medskip

Combining (21.2) and (21.4), (21.3) and (21.5), we derive the following several 
results. 

\medskip
   
\noindent {\bf Theorem 21.5.}\, {\em  Let $ S_A $ be given by  {\rm (21.1)}.  Then

{\rm (a)}\, There is an  $A^- \in \{ A^- \} $ such that
$$\displaylines{
\hspace*{2cm}
 r \left[ \begin{array}{cc}  A  & B  \\ C & D  \end{array} \right] = r(A) + r( \, D - CA^-B \, ), 
\hfill (21.6)
\cr
\hspace*{0cm} 
if \ and  \ only \ if \hfill
\cr
\hspace*{2cm}
 r \left[ \begin{array}{cc}  A  & B  \\ C & D  \end{array} \right] \leq r(A) + \min \left\{ \  r[ \, C, \ D \,],  \ \ \  r \left[ \begin{array}{c} B \\  D   \end{array} \right]  \  \right\}. \hfill (21.7)
\cr}
$$

{\rm (b)}\, The equality {\rm (21.6)} holds for all
$A^- \in \{ A^- \}$ if and only if
$$\displaylines{
\hspace*{2cm}
r\left[ \begin{array}{ccc}  A  &  0  & B \\ 0 & C  & D  \end{array} \right] = r(A) + r[ \, C , \ D \, ] \ \ and  \ \ \  r \left[ \begin{array}{cc}  A  &  0  \\ 0 & B  \\ C & D  \end{array} \right]
 = r( A ) +  r \left[ \begin{array}{c} B  \\ D  \end{array} \right].  \hfill (21.8) 
\cr}
$$ 
}  
{\bf Proof.}\, Note from  (21.2) that $ r \left[ \begin{array}{cc}  A  & B  \\ C & D  
\end{array} \right]  
- r(A )$ is an upper bound for  $ r( \, D - CA^-B \, )$. Thus there is an
$A^- \in \{ A^- \} $ such that (21.6) holds if and only if
$$\displaylines{
\hspace*{2cm}
 \max_{A^-}r( \,D - CA^-B \, ) = r \left[ \begin{array}{cc}  A  & B  \\ C & D  \end{array} \right]  
- r(A ). \hfill
\cr}
$$ 
Putting (21.4) in it  immediately yields (21.7). On the other hand, (21.6)
holds for all $A^- \in \{ A^- \}$ if and only if  
$$\displaylines{
\hspace*{2cm}
 \min_{A^-}r( \,D - CA^-B \, ) = r \left[ \begin{array}{cc}  A  & B  \\ C & D  \end{array} \right]  
- r(A ). \hfill
\cr}
$$ 
Putting (21.5) in it  yields (21.8).  \qquad  $ \Box$  

\medskip

The rank equality (21.6) was  examined by Carlson in \cite{Car} and Marsaglia and Styan in \cite{MaS2}. 
Their conclusion is that (21.6) holds if and only if
$ ( \, I - AA^- \, )B( \, I - S_A^-S_A \, ) = 0,$ $ ( \, I - S_A^-S_A \, )C( \, I - A^-A \, ) 
= 0$ and  $ ( \, I - AA^- \, )BS_A^-C( \, I - A^-A \, ) = 0$. In comparison,  (21.7) has no
 inner inverses in it, thus it is simpler and  is easier to verify. 

\medskip

\noindent {\bf Theorem 21.6.}\, {\em  Let $ S_A $ be given by {\rm (21.1)}.
Then

{\rm (a)}\, There is an  $A^- \in \{ A^- \} $ such that
$$\displaylines{
\hspace*{1.5cm}
 r \left[ \begin{array}{cc}  A  & B  \\ C & D  \end{array} \right]
 = r \left[ \begin{array}{cc}
 A  & B  \\ C & CA^-B  \end{array} \right]+ r( \, D - CA^-B \, ), 
\hfill (21.9)
\cr
\hspace*{0cm} 
if \ and  \ only \ if \hfill
\cr
\hspace*{1.5cm}
r\left[ \begin{array}{ccc}  A  &  0  & B \\ 0 & C  & D  \end{array} \right] = r[\, A, \ B  \,] + r[ \, C , \ D \, ] \ \ and  \ \ \  r \left[ \begin{array}{cc}  A  &  0  \\ 0 & B  \\ C & D  \end{array} \right]
 = r \left[ \begin{array}{c} A  \\ C \end{array} \right] +
 r \left[ \begin{array}{c} B  \\ D  \end{array} \right].  \hfill (21.10)
\cr}
$$

{\rm (b)}\, The equality  {\rm (21.9)} holds for all  $A^- \in \{ A^- \}$ if and only if 
$$\displaylines{
\hspace*{1.5cm}
R(B) \subseteq R(A) \ \ \ and   \ \ \ R(C^T) \subseteq
 R(A^T), \hfill (21.11)
\cr
\hspace*{0cm} 
or  \hfill
\cr
\hspace*{1.5cm}
r \left[ \begin{array}{cc}  A  & B \\ C & D  \end{array} \right] = r \left[ \begin{array}{c} A  \\ C \end{array} \right] +  r \left[ \begin{array}{c} B  \\ D  \end{array} \right]  \ \ and
 \ \ R(B) \subseteq R(A),   \hfill (21.12) 
\cr
\hspace*{0cm} 
or  \hfill
\cr
\hspace*{1.5cm}
r\left[ \begin{array}{cc}  A  &  B \\ C  & D  \end{array} \right] = r[\, A, \ B \,] + r[ \, C, \ D \, ] \ \ and  \ \ \ R(C^T) \subseteq R(A^T). \hfill (21.13)
\cr}
$$ 
}
{\bf Proof.}\, Note from  (21.3) that $ r \left[ \begin{array}{cc}  A  & B  \\ C & D  \end{array} \right]  
- r \left[ \begin{array}{cc}  A  & B  \\ C & CA^-B  \end{array} \right]$ is a lower bound for  
$ r( \, D - CA^-B \, )$. Thus (21.9) holds if and only if 
$$
 \min_{A^-}r( \,D - CA^-B \, ) = r \left[ \begin{array}{cc}  A  & B  \\ C & D  \end{array} \right] -
 r \left[ \begin{array}{cc}  A  & B  \\ C & CA^-B  \end{array} \right] 
= r \left[ \begin{array}{cc}  A  & B  \\ C & D  \end{array} \right] -
 r \left[ \begin{array}{c}  A \\ C \end{array} \right] - r[\,  A, \ B \,] + r(A). 
$$ 
Combining it with  (21.5) yields 
$$\displaylines{
\hspace*{1cm} 
r\left[ \begin{array}{ccc}  A  &  0  & B \\ 0 & C  & D  \end{array} \right] + 
r \left[ \begin{array}{cc}  A  &  0  \\ 0 & B  \\ C & D  \end{array} \right] = r \left[ \begin{array}{c} A  \\ C \end{array} \right] +  r \left[ \begin{array}{c} B  \\ D  \end{array} \right] + r[\, A, \ B \,] + r[ \, C, \ D \, ], \hfill
\cr}
$$ 
which is obviously equivalent to (21.10). On the other hand,  (21.9) holds 
for all $A^- \in \{ A^- \}$ if and only if  
$$\displaylines{
\hspace*{1cm}
 \max_{A^-}r( \,D - CA^-B \, ) = r \left[ \begin{array}{cc}  A  & B  \\ C & D  \end{array} \right] - r \left[ \begin{array}{cc}  A  & B  \\ C & CA^-B  \end{array} \right]. \hfill 
\cr}
$$ 
Combining it with (21.4) yields  (21.11)---(21.13). \qquad  $ \Box$  

\medskip

As a special case of  Schur complements, the rank and the range of the product 
$ CA^-B$ and their applications were examined by Baksalary and Kala in \cite{BK2}, 
 Baksalary and Mathew in\cite{BM} and Gross in \cite{Gro1}. Based
on the previous several theorems and corollaries, we now have the following three corollaries.

\medskip

\noindent {\bf Corollary 21.7.}\, {\em Let $A \in {\cal F}^{m \times n}, \,
B \in {\cal F}^{m \times k}$ and $ C \in {\cal F}^{l \times n}$ be given. Then

{\rm (a)}\, The maximal rank of $ CA^-B $ with respect to $ A^-$ is
$$ \displaylines{
\hspace*{2cm}
\max_{A^-}r( CA^-B)  = \min \left\{ \, r(B) , \ \ \ r(C),  \ \ \ r \left[ \begin{array}{cc}  A  & B \\  C & 0   \end{array} \right] - r(A) \,  \right\}. \hfill (21.14)
\cr}
$$  

{\rm (b)}\, The minimal rank of $ CA^-B $ with respect to $ A^-$ is
$$ \displaylines{
\hspace*{2cm}
\min_{A^-}r(CA^-B) = r\left[ \begin{array}{cc}  A  &  B  \\ C & 0  \end{array} \right] - r \left[ \begin{array}{c} A  \\ C  \end{array} \right] - r[ \, A, \ B \, ] + r(A). \hfill (21.15)
\cr}
$$ 

{\rm (c)}\, There is an $ A^- \in \{ A^-\}$ such that $ CA^-B = 0$ if and only if
$$\displaylines{
\hspace*{2cm}
 r\left[ \begin{array}{cc}  A  &  B  \\ C & 0  \end{array} \right] = r \left[ \begin{array}{c} A  \\ C  \end{array} \right] + r[ \, A, \ B \, ] - r(A). \hfill (21.16)
\cr}
$$

{\rm (d)}\, $ CA^-B = 0$ holds for all $ A^- \in \{ A^-\}$ if and only if $ B = 0 $ or $C = 0$ or
$ r\left[ \begin{array}{cc}  A  &  B  \\ C & 0  \end{array} \right] = r(A).$ 

{\rm (e)\cite{Gro1}}\,  The rank of $ CA^-B $ is invariant  with respect to the choice of
 $ A^-$  if and only if 
$$\displaylines{
\hspace*{2cm}
R(B) \subseteq  R(A)  \ \  and  \ \  R(C^T) \subseteq  R(A^T), \hfill
\cr
\hspace*{0cm} 
or  \hfill (21.17)
\cr
\hspace*{2cm}
r \left[ \begin{array}{cc}  A  & B \\ C  & 0 \end{array} \right] = r[\, A,  \ B \,] + r(C) 
\ \ \ and  \ \ \ R(C^T) \subseteq  R(A^T), \hfill (21.18)
\cr
\hspace*{0cm} 
or  \hfill
\cr
\hspace*{2cm}   
 r \left[ \begin{array}{cc}  A & B  \\ C & 0 \end{array} \right] = 
r\left[ \begin{array}{c} A  \\ C  \end{array} \right] + r(B) \ \ \ and  \ \ \  
R(B) \subseteq R(A). \hfill (21.19)
\cr}
$$
}
{\bf Corollary 21.8.}\, {\em Let $A \in {\cal F}^{m \times n}, \, B \in {\cal F}^{m \times k}$ and 
$ C \in {\cal F}^{l \times n}$ be given with $ B\neq 0$ and $ C\neq 0$. Then

{\rm (a)\cite{Gro1}}\,  The range  $ R(CA^-B) $ is invariant  with respect to the choice of
 $ A^-$  if and only if  $R(B) \subseteq R(A)$ and  $R(C^T) \subseteq  R(A^T),$ or
$$
\displaylines{
\hspace*{2cm}
r \left[ \begin{array}{cc}  A  & B \\ C  & 0 \end{array} \right] = r[\, A,  \ B \,] + r(C) 
\ \ \ and  \ \ \ R(C^T) \subseteq R(A^T). \hfill  (21.20)
\cr}
$$
 
{\rm (b)\cite{Gro1}}\, The range $ R[(CA^-B)^T] $ is invariant  with respect to the choice of
 $ A^-$  if and only if $R(B) \subseteq  R(A)$ and $R(C^T) \subseteq  R(A^T),$ or
$$ 
\displaylines{
\hspace*{2cm}   
r \left[ \begin{array}{cc}  A & B  \\ C & 0 \end{array} \right] = r\left[ \begin{array}{c} A  \\ C  \end{array} \right] + r(B) \ \ \ and  \ \ \  R(B) \subseteq R(A). \hfill (21.21)
\cr}
$$

{\rm(c)\cite{BM}}\, The rank of $ CA^-B $ is invariant  with respect to the choice of
 $ A^-$  if and only if $R(CA^-B)$ or $R[(CA^-B)^T] $ is
  invariant  with respect to the choice of $ A^-$. } 

\medskip 

\noindent {\bf Corollary 21.9.}\, {\em Let $A \in {\cal F}^{m \times n}, \, B \in {\cal F}^{m \times k}$ and 
 $ C \in {\cal F}^{l \times n}$. Then 
$$ 
\displaylines{
\hspace*{0.5cm}  
 \max_{A^-}r(A^-B) = r(B), \qquad   \max_{B^-}r(B^-A) = r(A),  \hfill (21.22)
\cr
\hspace*{0.5cm}  
 \max_{A^-}r(AA^-B) = \max_{B^-}r(BB^-A) = \min \{ \, r(A), \ \ r(B) \, \},  \hfill (21.23)
\cr
\hspace*{0.5cm}   
 \max_{A^-}r(CA^-) = r(C), \qquad   \max_{C^-}r(AC^-) = r(A),  \hfill (21.24)
\cr
\hspace*{0.5cm}  
 \max_{A^-}r(CA^-A) = \max_{C^-}r(AC^-C) = \min \{ \, r(A), \ \ r(C) \, \},  \hfill (21.25)
\cr
\hspace*{0.5cm}  
 \min_{A^-}r(AA^-B) = \min_{B^-}r(BB^-A)= \min_{A^-}r(A^-B) = \min_{B^-}r(B^-A) = r(A) + r(B) - r[ \, A, \ B \, ],  \hfill (21.26)
\cr
\hspace*{0.5cm}  
\min_{A^-}r(CA^-A) = \min_{C^-}r(AC^-C) = \min_{A^-}r(CA^-) = \min_{C^-}r(AC^-) = r(A) + r(C) - r \left[ \begin{array}{c} A  \\ C  \end{array} \right]. \hfill(21.27)
\cr}
$$ 
In particular$,$

{\rm (a)}\, There are $ A^-$  and $ B^-$ such that $ A^-B = 0$ and $ B^-A = 0$ if and only if  $R(A) \cap R(B) = \{ 0 \}.$ 

{\rm (b)}\, There are  $ A^-$ and $ C^-$ such that $ CA^- = 0$ and $ AC^- = 0$  if and only if $R(A^T) \cap R(C^T) = \{ 0 \}.$ 
} 

\medskip

The two formulas in (21.4) and (22.5) can help to
establish various rank equalities for matrix expressions that involve inner inverses of matrices,
and then to derive from them various consequences. Below are some of them. 

\medskip

\noindent {\bf Theorem 21.10.}\, {\em Let $A, \, B \in {\cal F}^{m \times n} $ be given. Then 
$$\displaylines{ 
\hspace*{0cm}
\max_{B^-}r( \, A - AB^-A \,)  = \min \left\{ \ r(A),  \ \ \
r(\, B - A \,)  - r(B) + r(A) \ \right\}, \hfill (21.28)
\cr 
\hspace*{0cm}
\min_{B^-}r( \, A - AB^-A \,) =  \min_{A^-, \, B^-}r( \, A^- - B^- \,) =
r(\, A - B \,) + r(A) + r(B) -  r[\, A, \ B \, ] - r
\left[ \begin{array}{c} A  \\ B  \end{array} \right]. \hfill (21.29)
\cr}
$$
In particular$,$

{\rm (a)}\, $ A $ and $ B $ have a common inner inverse if and only if 
$ r(\, A - B \,) = r \left[ \begin{array}{c} A  \\ B  \end{array} \right]  + r[\, A, \ B \,]  
-  r(A) - r(B).$

{\rm (b)}\, The inclusion $\{ B^- \} \subseteq \{ A^- \}$ holds if and only if $ A = 0$ or $ r(\, B - A \,) =
r(B) - r(A).$

{\rm (c)\cite{Mi6}}\, $\{ A^- \} = \{ B^- \}$ holds if and only if $ A = B.$

{\rm (d)} \ $\{ A^- \} \cap \{ B^- \} = \O$ holds if and only if $ 
 r(\, A - B \,) > r \left[ \begin{array}{c} A  \\ B  \end{array} \right]  + r[ \, A, \ B \,]  
-  r(A) - r(B).$

{\rm (e)}\, If $R(A) \cap R(B) = \{ 0\}$ and $R(A^T) \cap R(B^T) = \{ 0\}, $ 
then there must exist $ A^- \in \{ A^- \}$ and $ B^- \in \{ B^- \}$ such
that $  A^- = B^-$. }

\medskip

\noindent {\bf Proof.}\, Eq.\,(21.28) follows from (21.4);  (21.29) follows
from (21.5) and (20.5). The results in Parts (a)---(e) are direct consequences
 of (21.28) and (21.29).  \qquad $ \Box$

\medskip

A lot of consequences can be derived from Theorem 21.10. For example, let $ B = A^k$ in (21.29). Then we get 
$$
\displaylines{
\hspace*{1.5cm}
\min_{A^-, \, (A^k)^-}r[ \, A^- - (A^k)^- \,] =
r(\, A - A^k \,) + r(A^k) -  r(A). \hfill (21.30)
\cr}
$$   
Thus $ A $ and $A^k$ have a common inner inverse if and only if  $r(\, A - A^k \,) = r(A) - r(A^k)$. In that case, 
  $\{ A^- \} \subseteq \{ (A^k)^- \}$ holds by Theorem 21.10(b).

Replacing $ A $ and $ B $ in (21.29) by $ I_m - A $ and $ A $, respectively, we can get by (1.16) 
$$
\min_{(I_m - A)^-,\, A^-}r[ \, (I_m - A)^-  -  A^- \,] = r(I_m - 2A) + r(I_m - A)+  r(A) - 2m = 
r[\, A(I_m - A)(I_m - 2A)\,]. \eqno (21.31)
$$  
Thus $ I_m - A$ and $A $ have a common inner inverse if and only if $A(I_m - A)(I_m - 2A) = 0$.

Replacing $ A $ and $ B $ in (21.29) by $ A - I_m  $ and $ A $, respectively, we can get  
$$
\min_{( A - I_m)^-,\, A^-}r[ \, ( A - I_m)^-  -  A^- \,] = r(I_m - A)+  r(A) - m = 
r( \, A^2 - A\,). \eqno (21.32)
$$  
Thus $A$ is idempotent if and only if $ A - I_m $ and $A $ have a common inner inverse, this fact could be regarded as a new characterization of idempotent matrix.

Replacing $ A $ and $ B $ in (21.29) by $ I_m + A $ and $ A $, respectively, we can get by (1.11)  
$$
\min_{(I_m +  A)^-,\, A^-}r[ \, (I_m + A)^-  -  A^- \,] =  r(I_m + A)+  r(A) - m = 
r(A + A^2). \eqno (21.33)
$$  
Thus $ I_m + A$ and $A $ have a common inner inverse if and only if $A^2 = -A$.

Replacing $ A $ and $ B $ in (21.29) by $ I_m + A $ and $ I_m - A $, respectively, we can get by (1.15) 
$$
\min_{(I_m + A)^-,\, (I_m - A)^-}r[ \, (I_m + A)^-  -  (I_m - A)^- \,] = r(A) + r(I_m + A)+  r(I_m - A) - 2m 
= r(\, A^3 - A \,). \eqno (21.34)
$$  
In particular, $I_m + A$ and $ I_m - A$ have a common inner inverse if and only if $A $ is tripotent.

Replacing $ A $ and $ B $ in (21.15) by $ A + I_m$ and $  A - I_m$, respectively, we can get by (1.12) 
$$
\min_{( A + I_m )^-,\, ( A - I_m)^-}r[ \, ( A + I_m )^-  -  ( A - I_m)^- \,] =  r(  A  + I_m )+  r( A - I_m) - m 
= r(\, A^2 - I_m \,). \eqno (21.35)
$$  
This implies that $A $ is involutory if and only if $ A + I_m $ and $ A - I_m$ have a common inner inverse, this fact could be regarded as a new characterization of involutory matrix.

Now  suppose $ \lambda_1 \neq  \lambda_2$ are two scalars. Then it is easy to show by (21.29) and (1.16) the following 
two rank equalities
$$
\displaylines{
\hspace*{0.5cm}
\min_{(\lambda_1  I_m  - A )^-,\, ( \lambda_2I_m - A)^-}r[ \, (\lambda_1  I_m  - A )^-  - 
( \lambda_2I_m - A)^- \,] = r[ \, ( \lambda_1I_m - A)( \lambda_2I_m - A) \,], \hfill (21.36)
\cr
\hspace*{0.5cm}
\min_{( I_m  - \lambda_1A )^-,\, (  I_m - \lambda_2A )^-}r[ \, ( I_m  - \lambda_1A )^- - 
( I_m  - \lambda_2A )^-\,] = r[ \, A( I_m  - \lambda_1A )( I_m  - \lambda_2A ) \,]. \hfill (21.37)
\cr}
$$ 
Thus the two matrices  $\lambda_1 I_m  - A$  and  $ \lambda_2I_m - A$  have a common inner inverse if and only if 
$( \lambda_1I_m - A)( \lambda_2I_m - A) = 0$.  The two matrices $ I_m  - \lambda_1A$ and  $I_m  - \lambda_2A$ 
have  have a common inner inverse if and only if  $A( I_m  - \lambda_1A )( I_m  - \lambda_2A ) = 0$.

Again replacing $ A $ and $ B $ in (21.15) by $ A^k + A  $ and $ A^k - A $, respectively, we can get by (1.14)  
 $$
\min_{(A^k + A)^-,\, (A^k - A)^-}r[ \, (A^k + A)^-  -  (A^k - A)^- \,] = r(A^k + A)+  
r(A^k - A) - r(A) = r(\, A^{2k-1}- A \,). \eqno  (21.38)
$$  
In particular,  $A^{2k-1}= A$  if and only if $A^k + A$ and $A^k - A$ have a common inner inverse.

In general, suppose that $p(x)$ and $ q(x)$ are two polynomials without common roots. Then there is  
$$
\min_{p^-(A),\, q^-(A)}r[ \, p^-(A)  - q^-(A)  \,] = r[\, p^2(A)q(A) -  p(A)q^2(A) \, ]. 
\eqno  (21.39)
$$ 
Thus $ p(A)$ and $ q(A)$  have a common inner inverse if and only if $p^2(A)q(A) = p(A)q^2(A)$. 

>From (21.29) we also get 
$$
\displaylines{ 
\hspace*{1.5cm}
\min_{( A + B)^-, \, A^-}r[ \, ( \, A  + B \,)^- - A^-  \,] \hfill (21.40)
\cr
\hspace*{1.5cm}
=  \min_{( A + B)^-, \, B^-}r[ \, ( \, A  + B \,)^- - B^-\, ] =  r(\, A + B \,) + r(A) + r(B) -  r[\, A, \ B \, ] - r
\left[ \begin{array}{c} A  \\ B  \end{array} \right]. \hfill (21.41)
\cr}
$$ 
Hence we see that if  $R(A) \cap R(B) = \{ 0\}$ and $R(A^T) \cap R(B^T) = \{ 0\}, $ 
then $ A + B $ and $A$  must have a common inner inverse, meanwhile 
 then $ A + B $ and $B$  must have a common inner inverse.  

\bigskip

Now let $ M = \left[ \begin{array}{cc} A & B   \\ C  & 0 \end{array} \right]$ and  
$ N = \left[ \begin{array}{cc} 0 & B   \\ C  & 0 \end{array} \right]$. Then we get from (21.29) that 
$$
\displaylines{ 
\hspace*{1.5cm}
\min_{M^-, \, N^-}r( \, M^- - N^- \,) =  r\left[ \begin{array}{cc} A & B \\ C  & 0 \end{array} \right] -  
r\left[ \begin{array}{c} A  \\ C  \end{array} \right] - r[\, A, \ B \, ] + r(A). \hfill (21.42)
\cr}
$$
Hence $ M $ and $N$ have a common inner inverse if and only if 
$$
\displaylines{ 
\hspace*{1.5cm}
r \left[ \begin{array}{cc} A & B   \\ C  & 0 \end{array} \right] = r\left[ \begin{array}{c} A  \\ C 
 \end{array} \right] + r[\, A, \ B \, ] - r(A). \hfill (21.43)
\cr}
$$  
\hspace*{0.3cm} Next let $ M = \left[ \begin{array}{cc} A & B   \\ C  & 0 \end{array} \right]$ and  
$ N = \left[ \begin{array}{cc} A & 0  \\ 0  & 0 \end{array} \right]$. Then we can also get from (21.29) that 
$$
\displaylines{ 
\hspace*{1.5cm}
\min_{M^-, \, N^-}r( \, M^- - N^- \,) =  r\left[ \begin{array}{cc} A & B   \\ C  & 0 \end{array} \right] -  
r\left[ \begin{array}{c} A  \\ C  \end{array} \right] - r[\, A, \ B \, ] + r(A). \hfill (21.44)
\cr}
$$
Hence $ M $ and $N$ have a common inner inverse also if and only if (21.43) holds.  

Furthermore let $ M = \left[ \begin{array}{cc} A & B   \\ C  & D \end{array} \right]$ and  
$ N = \left[ \begin{array}{cc} A & 0  \\ 0  & D \end{array} \right]$. Then we can also derive from (21.29) that 
$$
\displaylines{ 
\hspace*{0.5cm}
\min_{M^-, \, N^-}r( \, M^- - N^- \,) \hfill
\cr
\hspace*{0.5cm}
 =  r(A) + r(B) + r(C) + r(D) + r\left[ \begin{array}{cc} A & B   
\\ C  & D \end{array} \right] -   r[\, A, \ B \, ] - r[\, C, \ D \, ]
 - r\left[ \begin{array}{c} A  \\ C  \end{array} \right]  - r\left[ \begin{array}{c} B  \\ D \end{array} \right]. 
\hfill (21.45)
\cr}
$$
Hence $ M $ and $N$ have a common inner inverse if and only if   
$$
\displaylines{ 
\hspace*{0.5cm}
r\left[ \begin{array}{cc} A & B   
\\ C  & D \end{array} \right] = r\left[ \begin{array}{c} A  \\ C  \end{array} \right] +
 r\left[ \begin{array}{c} B  \\ D \end{array} \right] +
 r[\, A, \ B \, ] + r[\, C, \ D \, ] -  r(A) - r(B) - r(C) - r(D) 
\hfill (21.46)
\cr}
$$
holds. In particular, if $r(M) =  r(A) + r(B) + r(C) + r(D)$, then $M $ and $ N$ must have a common inner inverse.   

Finally let  $ M = \left[ \begin{array}{cc} A & B   \\ C  & D \end{array} \right]$ and  
$ N = \left[ \begin{array}{cc} A & 0  \\ 0  & 0 \end{array} \right]$. Then we derive from (21.29) that 
$$
\displaylines{ 
\hspace*{0cm}
\min_{M^-, \, N^-}r( \, M^- - N^- \,)  =   r\left[ \begin{array}{cc} A & B   
\\ C  & D \end{array} \right] + r\left[ \begin{array}{cc} 0 & B   
\\ C  & D \end{array} \right] + r(A) - r\left[ \begin{array}{ccc} A & 0 & B   
\\ 0 & C  & D \end{array} \right] - r\left[ \begin{array}{cc} A & 0 \\ 0 & B   
\\ C  & D \end{array} \right]. \hfill (21.47)
\cr}
$$
Hence $ M $ and $N$ have a common inner inverse if and only if   
$$
\displaylines{ 
\hspace*{1.5cm}
r\left[ \begin{array}{cc} A & B \\ C  & D \end{array} \right] 
=r\left[ \begin{array}{cc} A & 0 \\ 0 & B \\ C  & D \end{array} \right] + 
r\left[ \begin{array}{ccc} A & 0 & B \\ 0 & C  & D \end{array} \right] -  r\left[ \begin{array}{cc} 0 & B   
\\ C  & D \end{array} \right] - r(A). \hfill (21.48)
\cr}
$$
{\bf Theorem 21.11.}\, {\em Let $A, \, B \in {\cal F}^{m \times n} $ be given. Then 
$$
\displaylines{ 
\hspace*{1cm}
\max_{(A + B)^-}r[ \, A( \, A + B \, )^-B \,]\hfill
\cr
\hspace*{1cm}
= \max_{(A + B)^-}r[ \, B( \, A + B \, )^-A \,] = \min \left\{\, r(A),  \ r(B), \  r(A) + r(B) - r( A + B) \, \right\}, 
\hfill (21.49)
\cr
and \hfill
\cr
\hspace*{1cm}
\min_{(A + B)^-}r[ \, A( \, A + B \, )^-B \,] \hfill
\cr
\hspace*{1cm} 
= \min_{(A + B)^-}r[ \, B( \, A + B \, )^-A \,] =  r(\, A + B \,) + r(A) + r(B) -  r[ \, A, \ B \, ] - r \left[ \begin{array}{c} A  \\ B  \end{array} \right]. \hfill   (21.50)
\cr}
$$
In particular$,$

{\rm (a)}\, There is an $ ( \, A + B \, )^-$ such that 
$  A( \, A + B \, )^-B  = 0$ if and only if $r(\, A + B \,) = $  $ r \left[ \begin{array}{c} A  \\ B  \end{array} \right] $ $ + r[\, A, \ B \, ] -  r(A) - r(B).$

{\rm (b)}  $  A( \, A + B \, )^-B  = 0$ holds  for all $ ( \, A + B \, )^- $ if and only if $ A = 0 $ or $ B =0 $ or $r(\, A + B \,) = r(A) + r(B).$

{\rm (c)\cite{PSS}}\,  The rank of  $A( \, A + B \, )^-B$ is invariant with respect to the choice of 
 $( \, A + B \, )^-$ if and only if
 $R(B) \subseteq R( \, A +B \,)$  and $R(A^T) \subseteq R( \, A^T +B^T \,),$ that is$,$ $A$ 
and $B$ are parallel summable.}

\medskip

A parallel result to  Theorem 21.11 is   

\medskip

\noindent {\bf Theorem 21.12.}\, {\em Let $A, \, B \in {\cal F}^{m \times n} $ be given. Then 
$$
\displaylines{ 
\hspace*{1cm}
\max_{A^-,\, B^-}r[ \, A^-( \, A + B \, )B^- \,] = \max_{A^-,\, B^-}r[ \, B^-( \, A + B \, )A^- \,] = r(\, A + B \, ), 
\hfill (21.51)
\cr
and \hfill
\cr
\hspace*{1cm}
\min_{A^-,\, B^-}r[ \, A^-( \, A + B \, )B^- \,] \hfill
\cr
\hspace*{1cm}
=
\min_{A^-,\, B^-}r[ \, B^-( \, A + B \, )A^- \,] =  r(\, A + B \,) + r(A) + r(B) -  r[ \, A, \ B \, ] -
 r \left[ \begin{array}{c} A  \\ B  \end{array} \right]. \hfill   (21.52)
\cr}
$$
In particular$,$

{\rm (a)}\, There are  $A^-, \, B^-$ such that 
$  A^-( \, A + B \, )A^-  = 0$ if and only if $r(\, A + B \,) = r \left[ \begin{array}{c} A  \\ B  \end{array} \right] $ $ + r[\, A, \ B \, ] -  r(A) - r(B).$

{\rm (b)}\, The rank of  $A^-( \, A + B \, )B^-$ is invariant with respect to the choice of 
 $A^-, \,  B^-$ if and only if $R(A) = R(B)$  and $R(A^T)= R(B^T)$.}

\medskip

\noindent {\bf Proof.}\, According to (21.22) and (21.26) we first get 
$$
\max_{A^-}r[ \, A^-( \, A + B \, )B^- \,] = r[\, ( \, A + B \, )B^-\, ], 
$$ 
\begin{eqnarray*}
\min_{A^-}r[ \, A^-( \, A + B \, )B^- \,] & = & r(A) + r[\, ( \, A + B \, )B^- \,] -
 r[\, A, \ ( \, A + B \, )B^- \,]  \\
 & = &r(A) + r[\, ( \, A + B \, )B^- \,] -
 r[\, A, \ B \,] 
\end{eqnarray*}
Next by (21.24) and (21.27), we find 
$$
\max_{B^-}r[ \, ( \, A + B \, )B^- \,] = r( \, A + B \, ), 
$$
$$ 
\min_{B^-}r[ \, ( \, A + B \, )B^- \,] =  r(B) + r( \, A + B \, ) - 
 r\left[ \begin{array}{c} A +B \\ B  \end{array} \right] =  r(B) + r( \, A + B \, ) - 
 r\left[ \begin{array}{c} A \\ B  \end{array} \right].
$$
Combining them yields (21.51) and  (21.52). \qquad $ \Box$.

\medskip

From (21.52) we can also find some interesting consequences. For example, let $ B = I_m - A $ in (21.52), we can get 
 $$
\min_{A^-,\, (I_m - A)^-}r[ \, A^-(I_m - A)^- \,] =
\min_{A^-,\, (I_m - A)^-}r[ \, A^-(I_m - A)^- \,] =  r( \,A - A^2\,). \eqno (21.53)
$$ 
Thus $ A$ is idempotent if and only if there are $ A^- $  and $ (I_m - A)^-$ such that
 $A^-(I_m - A)^- = 0$, which could be regard as a new characterization for 
idempotent matrix.    

Replacing $ A $ and $ B $  in (21.52) by  $I_m + A $  and $I_m - A $, respectively,  we can get 
 $$
\min_{(I_m + A)^-,\, (I_m - A)^-}r[ \, (I_m + A)^-(I_m - A)^- \,] =
\min_{(I_m + A)^-,\, (I_m - A)^-}r[ \, (I_m - A)^-(I_m + A)^- \,] =  r( \,I_m - A^2\,). \eqno (21.54)
$$ 
Thus $ A$ is involutory if and only if there are $  (I_m + A)^- $  and $ (I_m - A)^-$ such that
 $ (I_m + A)^-(I_m - A)^- = 0$, which could be regard as a new characterization for  involutory matrix.    

\medskip  

In general replacing $ A $ and $ B $  in (21.52) by  $ \lambda_1 I_m - A $  and $ -(\lambda_2 I_m - A)$,
 respectively, where  $\lambda_1 \neq  \lambda_2$,  we can get 
$$
\displaylines{ 
\hspace*{1.5cm}
\min_{(\lambda_1I_m - A)^-,\, (\lambda_2I_m - A)^-}r[ \, ( \lambda_1 I_m - A)^-( \lambda_2 I_m - A)^- \,]
 =r[ \, ( \lambda_1 I_m - A)( \lambda_2 I_m - A) \,]. \hfill (21.55)
\cr}
$$ 
Thus there are $ ( \lambda_1 I_m - A)^- $ and $( \lambda_2 I_m - A)^- = 0$ such that 
 $( \lambda_1 I_m - A)^-( \lambda_2 I_m - A)^- = 0$ if and only if 
 $( \lambda_1 I_m - A)( \lambda_2 I_m - A)= 0$. 

Motivated by  (21.55), we find the the following general result. 

\medskip  

\noindent {\bf Theorem 21.13.}\, {\em Let $A \in {\cal F}^{m \times n}$ and  
$ \lambda_1, \, \cdots, \, \lambda_k \in {\cal F}$ with  $ \lambda_i \neq  \lambda_j$ for all $ i \neq j $.  Then
$$
\min_{(\lambda_1  I_m  - A )^-, \,\cdots, \, (\lambda_k I_m  - A )^-} 
r[\, (\lambda_1  I_m  - A )^- \cdots 
 (\lambda_k I_m  - A )^-  \, ] = r[\, (\lambda_1  I_m  - A )\cdots 
(\lambda_k I_m  - A )  \, ].  \eqno (21.56) 
$$}
\noindent {\bf Proof.}\, According to (21.26) we first get 
$$\displaylines{ 
\hspace*{0cm}
\min_{(\lambda_1  I_m  - A )^-} r[\, (\lambda_1  I_m  - A )^- \cdots 
 (\lambda_k I_m  - A )^-  \, ]  \hfill
\cr
\hspace*{0cm}
= r(\lambda_1  I_m  - A ) + r[\, (\lambda_2 I_m  - A )^- \cdots (\lambda_k I_m  - A )^-\,]
  - r[\, (\lambda_1 I_m  - A ), \ (\lambda_2 I_m  - A )^-  \cdots (\lambda_k I_m  - A )^- \,].
  \hfill (21.57) 
\cr}
$$
Notice that $ \lambda_i \neq  \lambda_j$  for $ i \neq j$.  Then there must be 
$$ 
r[\, (\lambda_1 I_m  - A ), \ (\lambda_2 I_m  - A )^-  \cdots (\lambda_k I_m  - A )^- \,] = m
  \eqno (21.58) 
$$
for all $ (\lambda_2 I_m  - A )^-,   \cdots, \  (\lambda_k I_m  - A )^-$. To show this fact, we need 
the following two rank formulas
$$
\displaylines{ 
\hspace*{1.5cm}
\min_{B^-} r[\, A, \ B^- \, ] =  r(A) + r(B) - r(BA),  \hfill (21.59) 
\cr
\hspace*{1.5cm}
\min_{B^-} r[\, A, \ B^-C \, ] =  r(A) + r(B) - r(BA) + r[\, BA, \ C \, ] - 
r[\, B, \ C \, ] ,  \hfill (21.60)
\cr}
$$
We see first  by (21.59) and (1.16) that for $ \lambda_t \neq \lambda_j, \ j = 1, \, \cdots, \, i$   there is
$$ 
\displaylines{ 
\hspace*{0.5cm}
\min_{(\lambda_t I_m  - A )^-} r[\, (\lambda_1 I_m  - A ) \cdots (\lambda_{i} I_m  - A ) , \ 
 (\lambda_t I_m  - A )^- \,]  \hfill
\cr
\hspace*{0.5cm}
= r[\, (\lambda_1 I_m  - A ) \cdots (\lambda_{i} I_m  - A ) \,] + r(\lambda_t I_m  - A ) 
- r[\, (\lambda_1 I_m  - A ) \cdots (\lambda_{i} I_m  - A )(\lambda_t I_m  - A ) \,]  = m. 
 \hfill
\cr}
$$
That is to say,  
$$
r[\, (\lambda_1 I_m  - A ) \cdots (\lambda_{i} I_m  - A ) , \  (\lambda_t I_m  - A )^- \,] =m 
$$ 
holds for any $(\lambda_t I_m  - A )^- $ with $ \lambda_t \neq \lambda_j, \ j = 1, \, \cdots, \, i$.  
Now suppose that 
$$
r[\, (\lambda_1 I_m  - A ) \cdots (\lambda_{i} I_m  - A ) , \  (\lambda_{t+1} I_m  - A )^- 
\cdots (\lambda_k I_m  - A )^- \,]  = m 
$$
holds for all $ (\lambda_{t+1} I_m  - A )^-,  \, \cdots, \, (\lambda_k I_m  - A )^-$ and $ 1\leq i < 
t < k$. Then we can obtain by (21.60), (1.16) and induction hypothesis that 
$$
\displaylines{ 
\hspace*{0cm}
\min_{(\lambda_{t} I_m  - A )^-}r[\, (\lambda_1 I_m  - A ) \cdots (\lambda_{i} I_m  - A ) , \ 
 (\lambda_{t} I_m  - A )^-\cdots (\lambda_k I_m  - A )^- \,]  \hfill
\cr
\hspace*{0cm}
= r[\, (\lambda_1 I_m  - A ) \cdots (\lambda_{i} I_m  - A ) \, ]  + r(\lambda_t I_m  - A ) 
-  r[\, (\lambda_1 I_m  - A ) \cdots (\lambda_{i} I_m  - A )(\lambda_t I_m  - A )\,]  \hfill
\cr
\hspace*{0.5cm}
+ \ r[\, (\lambda_1 I_m  - A ) \cdots (\lambda_{i} I_m  - A )(\lambda_t I_m  - A ), \ 
 (\lambda_{t+1} I_m  - A )^- \cdots (\lambda_k I_m  - A )^- \,]  \hfill
 \cr
\hspace*{0.5cm}
- \  r[\, (\lambda_t I_m  - A ), \   (\lambda_{t+1} I_m  - A )^-\cdots (\lambda_k I_m  - A )^- \,]  \hfill
\cr
\hspace*{0cm}
= m + m - m = m,  \hfill
\cr}
$$
that is,
$$
 r[\, (\lambda_1 I_m  - A ) \cdots (\lambda_{i} I_m  - A ), \ 
 (\lambda_{t} I_m  - A )^-\cdots (\lambda_k I_m  - A )^- \,] =m  \eqno (21.61) 
$$
holds for all $ (\lambda_{t} I_m  - A )^-,  \, \cdots, \, (\lambda_k I_m  - A )^-$ and $ 1\leq i < 
t < k$.
When $ i = 1$ and $ t = 2$, (21.61) becomes (21.58). In that case, (21.57) reduces to  
$$
\min_{(\lambda_1  I_m  - A )^-} r[\, (\lambda_1  I_m  - A )^- \cdots 
 (\lambda_k I_m  - A )^-  \, ] = r(\lambda_1  I_m  - A ) + r[\, (\lambda_2 I_m  - A )^- \cdots (\lambda_k I_m  - A )^-\,] -m. \eqno (21.62)
$$ 
Repeatedly applying (21.62) for the product $(\lambda_2 I_m  - A )^- \cdots (\lambda_k I_m  - A )^-$ in 
(21.62), we eventually  get 
$$
\min_{(\lambda_1  I_m  - A )^-, \,\cdots, \, (\lambda_k I_m  - A )^-} 
r[\, (\lambda_1I_m  - A )^- \cdots 
 (\lambda_kI_m  - A )^-  \, ] = r(\lambda_1 I_m - A ) + \cdots  + r(\lambda_kI_m  - A )
 - m(k-1),
$$
which, by (1.16), is the desired formula (21.56).   \qquad $\Box$ 
 
\medskip

When two square matrices $A $ and $ B$ of the same size are nonsingular, it is well known 
that $ A^{-1} + B^{-1} = A^{-1}(A + B)B^{-1}$. This fact  motivates us to consider the 
relationship between  $ A^- + B^-$ and $ A^-(A  + B)B^-$ in general case.  Using 
the rank formula (21.5) we can simply find that
$$
\min_{A^-, \, B^-} r[\,  A^- + B^- - A^-(A  + B)B^- \,]  = 0, 
$$
which implies the following. 

\medskip
 
\noindent {\bf Theorem 21.14.}\, {\em Let $A, \, B \in {\cal F}^{m \times n} $ be given. Then there must
exist $ A^-$ and $ B^-$ such that 
$$ 
A^- + B^- = A^-(A + B)B^- \eqno (21.63)
$$ 
holds.} 
 
\medskip

As applications, we can simply get from (21.63) that there must exist $ A^-$ and $(I_m - A)^- $ such that 
$$ 
A^- + (I_m - A)^- = A^-(I_m - A)^-.  \eqno (21.64)
$$ 
and there must exist $ A(I_m + A)^- $ and $(I_m - A)^- $ such that 
$$ 
(I_m + A)^- + (I_m - A)^- = 2(I_m + A)^- (I_m - A)^-.  \eqno (21.65)
$$  
This result leads to the following conjecture .

\medskip
 
\noindent {\bf Conjecture 21.15.}\, {\em Let  $A \in {\cal F}^{m \times n}$ and  
$ \lambda_1, \, \cdots, \, \lambda_k \in {\cal F}$ with  $ \lambda_i \neq  \lambda_j.$ Then  there
 exist $(\lambda_1  I_m  - A )^-, \,\cdots, \, (\lambda_k I_m  - A )^- $ such that 
$$
\frac{1}{p_1}(\lambda_1  I_m  - A )^- + \cdots +  \frac{1}{p_k}(\lambda_k I_m  - A )^- = 
 (\lambda_1  I_m  - A )^- \cdots (\lambda_k I_m  - A )^-.  \eqno (21.66)
$$
where
$$
p_i = (\lambda_1 - \lambda_i)\cdots(\lambda_{i-1} - \lambda_i)
(\lambda_{i+1} - \lambda_i)\cdots(\lambda_k - \lambda_i),  \ \ \ i = 1, \, \cdots, \, k.  
$$}

\noindent {\bf Theorem 21.16.}\, {\em Let $A, \, B \in {\cal F}^{m \times n} $ be given and let $ M  
= {\rm diag}(\, A, \ B \, )$ and  $ N = A + B$ . Then 
$$
\max_{N^-}r \left(   M  -
\left[ \begin{array}{c} A  \\ B  \end{array} \right] N^- [ \, A, \ B \, ]  \right) 
=  \min_{N^-}r \left( M -
\left[ \begin{array}{c} A  \\ B  \end{array} \right] N^- [ \, A, \ B \, ]  \right) 
=  r(A) + r(B) - r(N). \eqno (21.67) 
$$
That is$,$ the rank of  $M -\left[ \begin{array}{c} A  \\ B  \end{array} \right] N^- [ \, A, \ B \, ]$  is 
invariant with respect to the choice of $N^-$. In general$,$ for $A_1, \, A_2, \, \cdots, \, A_k 
\in {\cal F}^{m \times n},$ there is 
\begin{eqnarray*}
\max_{N^-}r \left(   M - \left[ \begin{array}{c} A_1  \\ \vdots\\ A_k  \end{array} \right] N^- 
[ \, A_1, \ \cdots, \  A_k \, ]  \right) & = & \min_{N^-}r \left(  M - \left[ \begin{array}{c} A_1 
 \\ \vdots \\ A_k  \end{array} \right] N^- [ \, A_1, \ \cdots, \  A_k \, ]  \right) \\ 
& = & r(A_1) + \cdots + r(A_k) - r(N),   
\end{eqnarray*} 
where $ M ={\rm diag}( \, A_1, \, \cdots, \,  A_k \, )$ and $ N =A_1 + \cdots + A_k.$ In particular$,$ the 
equality  
$$ 
 \left[ \begin{array}{c} A_1  \\ \vdots \\ A_k  \end{array} \right] ( \,A_1 + \cdots + A_k \,)^-[ \, A_1, \ \cdots, \  A_k \, ] = \left[ \begin{array}{ccc} A_1  & & \\ &  \ddots &  \\ & &  A_k \end{array} \right] 
 \eqno (21.68) 
$$ 
holds for all $( \,A_1 + \cdots + A_k \,)^- $ if and only if $ r(\,A_1 + \cdots + A_k \, ) = r(A_1) + \cdots + r(A_k)$.} 

\medskip

\noindent {\bf Theorem 21.17.}\, {\em \ Let $ M = \left[ \begin{array}{cc} A  & B \\ C & D \end{array} \right]$ 
 be a partitioned matrix over ${\cal F}.$ Then 
$$
 \max_{A^-}r \left( M - \left[ \begin{array}{c} A  \\ C  \end{array} \right] A^- [ \, A, \ B \, ] \right) = \min_{A^-}r \left(  M -
\left[ \begin{array}{c} A  \\ C  \end{array} \right] A^- [\, A, \ B \, ] \right) =  r(M) -r(A). 
$$
That is$,$ 
$$
 r\left[ \begin{array}{cc} A  & B \\ C & D \end{array} \right]  = r(A) +
r \left( M - \left[ \begin{array}{c} A  \\ C  \end{array} \right] A^- [ \, A, \ B \, ] \right),
$$
which is exactly the formula {\rm (1.5)}.}

\medskip

\noindent {\bf Theorem  21.18.}\, {\em Let $ A \in  {\cal F}^{m \times m}$ be
given. Then
$$\displaylines{ 
\hspace*{1.5cm}
\max_{A^-} r( \, AA^- - A^-A \, ) = \min \{ \, 2m - 2r(A),  \ \ 2r(A) \, \}, \hfill (21.69)
\cr
\hspace*{1.5cm}
\min_{A^-} r( \, AA^- - A^-A \, ) = 2r(A)- 2r(A^2).  \hfill (21.70)
\cr}
$$ } 
{\bf Proof.}\, Since both  $ AA^- $ and $ A^-A $  are idempotent, we see by  (3.1) that  the rank of  
$ AA^- - A^-A$ can be written as 
$$\displaylines{
\hspace*{2cm}
r( \, AA^- - A^-A \, ) = r \left[ \begin{array}{c}  AA^- \\ A^-A \end{array} \right] + 
r[ \, AA^-,  \ A^-A \, ] - r(AA^-) -  r(A^-A). \hfill
\cr}
$$
Note that $ r(AA^-) = r(A^-A) = r(A)$,  $r \left[ \begin{array}{c}  AA^- \\ A^-A \end{array} \right]
= r \left[ \begin{array}{c}  AA^- \\ A \end{array} \right]$ and $r[ \, AA^-,  \ A^-A \, ] = r[ \, A, \ A^-A \, ]$. Then  
$$
\displaylines{
\hspace*{2cm} 
r( \, AA^- - A^-A \, ) = r \left[ \begin{array}{c}  AA^- \\ A \end{array} \right] + r[ \, A,  \ A^-A \, ] - 2r(A). 
\hfill
\cr}
$$
On the other hand, form  the general expression of $ A^- = A^{\sim} + F_AV + WE_A,$ we also know that 
$ AA^- = AA^{\sim} + AWE_A,$ and  $ A^-A = A^{\sim}A + F_AVA.$ Thus $ AA^- $ and $ A^-A$ are in fact two independent matrix
expressions. In that case, we see that  
$$ 
\displaylines{
\hspace*{2cm}
\max_{A^-} r( \, AA^- - A^-A \, ) = \max_{A^-} r\left[ \begin{array}{c}  AA^- \\ A \end{array} \right] + 
\max_{A^-} r[ \, A,  \ A^-A  \, ] - 2r(A), \hfill (21.71)
\cr
\hspace*{2cm}
\min_{A^-} r( \, AA^- - A^-A \, ) = \min_{A^-} r\left[ \begin{array}{c}  AA^- \\ A \end{array} \right] + 
\min_{A^-} r[ \, A,  \ A^-A  \, ] - 2r(A). \hfill (21.72)
\cr}
$$
According to  (21.4) and (21.5), we easily find that 
$$\displaylines{
\hspace*{2cm}
\max_{A^-} r \left[ \begin{array}{c} AA^- \\ A \end{array} \right] =  \max_{A^-}
 r \left( \left[ \begin{array}{c} 0  \\ A \end{array} \right]  
+ \left[ \begin{array}{c} A \\ 0  \end{array} \right]A^- \right) =  \min \{ \, 2r(A), \ \ m  \, \}, \hfill
\cr
\hspace*{2cm}
\min_{A^-} r \left[ \begin{array}{c} AA^- \\ A \end{array} \right]  =  \min_{A^-}r \left( 
\left[ \begin{array}{c} 0  \\ A \end{array} \right]  
+ \left[ \begin{array}{c} A \\ 0  \end{array} \right]A^- \right)  =  2r(A) - r(A^2), \hfill
\cr
\hspace*{2cm}
\max_{A^-} r[ \, A,  \ A^-A \, ]  =  \max_{A^-} r( \, [ \, A,  \ 0 \, ]  +  A^-[ \, 0, \ A \, ] \, ) 
=  \min \{ \, 2r(A), \ \ m  \, \}, \hfill
\cr
\hspace*{2cm}
\min_{A^-} r[ \, A,  \ A^-A \, ]  =  \min_{A^-} r( \, [ \, A,  \ 0  \, ]  + A^-V[ \, 0, \ A \, ] \, ) = 2r(A) - r(A^2). 
\hfill
\cr
}  
$$
Putting the above four results in  (21.71) and (21.72) yields (21.69) and (21.70). \qquad$ \Box$ 

\medskip
  
\noindent {\bf Corollary 21.19.}\, {\em Let $ A \in  {\cal F}^{m \times m}$ be given. 

{\rm (a)}\, There is an $A^-$ such that $ AA^- - A^-A$ is nonsingular if
and only if $ m$ is even and $r(A) = m/2$.

{\rm (b)}\, There is an $A^-$ such that $ AA^- = A^-A$ if and only if
$r(A^2) = r(A)$.

{\rm (c)}\, The rank of  $ AA^- - A^-A$ is invariant with respect to the choice of $A^- $ if and only if 
$A^2 = 0 $ or  $r(A^2) = 2r(A) - m$. } 

\medskip

The two rank formulas (21.69) and (21.70) manifest that the maximal and minimal ranks of 
$ AA^- - A^-A$ are even. Recall from (6.1) that the rank  of $ AA^{\dagger}  - A^{\dagger}A$  is even, too.
Thus we have the following conjecture. 

\medskip

\noindent {\bf Conjecture 21.20.}\, {\em Let $ A \in  {\cal F}^{m \times m}$ be given. Then the rank of 
the  matrix expression $ AA^- - A^-A$ is even for any $A^-$.} 

\medskip
  
In the  same way we can  establish  the following. The details are omitted. 

\medskip
  
\noindent {\bf Theorem  21.21.}\, {\em Let $ A \in  {\cal F}^{m \times k}$ and  $ B \in  {\cal F}^{l \times m}$ be given. Then

 {\rm (a)}\, The maximal and the minimal ranks of $ AA^- - B^-B$  with
 respect to $ A^-$  and $ B^-$ are 
$$ \displaylines{ 
\hspace*{1.5cm}
\max_{A^-,\, B^-} r( \, AA^- - B^-B \, ) =
\min \{ \, 2m - r(A) - r(B),  \ \ r(A) + r(B) \, \}, \hfill (21.73)
\cr
\hspace*{1.5cm}
\min_{A^-,\, B^-} r( \, AA^- - B^-B \, ) = r(A) + r(B) - 2r(BA). \hfill (21.74)
\cr}
$$ 

{\rm (b)}\, There are $A^-$  and $B^-$ such that $ AA^- - B^-B$ is
nonsingular if and only if$r(A) + r(B) = m$.  

{\rm (c)}\, There are $A^-$  and $B^-$ such that $ AA^- = B^-B$ if and
only if $r(A) + r(B) = 2r(BA)$.

{\rm (d)} \, The rank of $ AA^- - B^-B$ is invariant with respect to the
choice of $A^- $  and $ B^-$ if and only if $BA = 0$ or $r(BA) = r(A) + r(B) - m$. } 

\medskip

As for extreme ranks of $ AA^- + B^-B$ we shall present them in Chapter 27. Moreover, we  can also determine 
the maximal and the minimal ranks of $ BB^-A - AC^-C$.

\medskip

\noindent {\bf Theorem  21.22.}\, {\em Let $ A \in  {\cal F}^{m \times n}, \, B \in {\cal F}^{m \times k }$  and
 $ C \in  {\cal F}^{l \times n}$  be given. Then 
$$ \displaylines{ 
\hspace*{1cm}
\max_{B^-,\, C^-} r( \, BB^-A - AC^-C \, ) \hfill
\cr
\hspace*{1cm}
 = \min \left\{ \, r[\, A, \ B \,], \ \ \  r\left[ \begin{array}{c} A \\ C \end{array} \right], \ \ \
 r(B) + r(C), \ \ \ r\left[ \begin{array}{c} A \\ C \end{array} \right]  + r[\, A, \ B \,] 
- r(B) - r(C) \,  \right\}, \hfill (21.75)
\cr
and \hfill
\cr
\hspace*{1cm}
\min_{B^-,\, C^-} r( \, BB^-A - AC^-C \, ) = r\left[ \begin{array}{c} A  \\ C  \end{array} \right]
 + r[\, A, \ B \,] + r(B) + r(C) - 2r\left[ \begin{array}{cc} A & B \\ C &  0 \end{array} 
\right]. \hfill (21.76)
\cr}
$$ }
{\bf Proof.}\,  According to (4.1), the rank of  $ BB^-A - AC^-C$ can be written as 
\begin{eqnarray*} 
r( \, BB^-A - AC^-C \, ) & = & r \left[ \begin{array}{c}  BB^-A \\ C^-C \end{array} \right] + 
r[ \, AC^-C,  \ BB^- \, ] - r(BB^-) -  r(C^-C) \\
&  = &  r \left[ \begin{array}{c}  BB^-A \\ C \end{array} \right] + r[ \, AC^-C,  \ B \, ] - r(B) - r(C). 
\end{eqnarray*} 
Hence
$$  \displaylines{ 
\hspace*{1cm}
\max_{B^-,\, C^-} r( \, BB^-A - AC^-C \, )= \max_{B^-} r\left[ \begin{array}{c}  BB^-A \\ C \end{array} \right] + \max_{C^-} r[ \, AC^-C,  \ B \, ] - r(B) - r(C), \hfill (21.77)
\cr 
\hspace*{1cm}
\min_{B^-,\, C^-} r( \, BB^-A - AC^-C \, )= \min_{B^-} r\left[ \begin{array}{c}  BB^-A \\ C \end{array} \right] + \min_{C^-} r[ \, AC^-C,  \ B \, ] - r(B) - r(C). \hfill (21.78) 
\cr}
$$
According to  (21.4) and (21.5),  we easily find that 
$$ \displaylines{ 
\hspace*{1.5cm} 
\max_{B^-} r\left[ \begin{array}{c}  BB^-A \\ C \end{array} \right] = \min \left\{ \, r(B) + r(C),  \ \ 
 r\left[ \begin{array}{c} A \\ C \end{array} \right]  \, \right\}, \hfill
\cr
\hspace*{1.5cm}
\min_{B^-} r\left[ \begin{array}{c}  BB^-A \\ C \end{array} \right] = r(B) + r(C) + r\left[ \begin{array}{c} A \\ C \end{array} \right] -  r\left[ \begin{array}{cc} A & B  \\ C & 0 \end{array} \right], \hfill 
\cr
\hspace*{1.5cm}
 \max_{C^-} r[ \, AC^-C,  \ B \, ]  =  \min \{ \, r[ \, A,  \ B \, ],  \ \ \   r(B) + r(C) \, \}, \hfill
\cr 
\hspace*{1.5cm}
 \min_{C^-} r[ \, AC^-C,  \ B \, ] = r(B) + r(C) +  r[ \, A,  \ B \, ] - r\left[ \begin{array}{cc} A & B  \\ C & 0 \end{array} \right].  \hfill
\cr}
$$
Putting them in  (21.77) and (21.78) yields  (21.75) and (21.76). \qquad $ \Box$ 

\medskip
 
\noindent {\bf Corollary 21.23.}\, {\em Let $ A \in  {\cal F}^{m \times n}, \,B \in {\cal F}^{m \times k }$  and
 $ C \in  {\cal F}^{l \times m}$  be given.

{\rm (a)}\, Assume  $ A $ is square.  Then there are  $B^-$  and $C^-$ such that $ BB^-A - AC^-C$ is
nonsingular if and only if $ A,$ $B$ and $C$ satisfy the following rank equality  
$$
\displaylines{
\hspace*{1.5cm}
r\left[ \begin{array}{c} A \\ C \end{array} \right]
= r[ \, A,  \ B \, ] = r(B) + r(C) =m. \hfill (21.79)
\cr}
$$ 

{\rm (b)}\, There are  $B^-$  and $C^-$ such that $ BB^-A = AC^-C$ if
and only if  $ A,$ $B$ and $C$ satisfy the rank additivity
condition
$$\displaylines{
\hspace*{1.5cm}
r\left[ \begin{array}{cc} A & B  \\ C & 0 \end{array} \right] =r\left[ \begin{array}{c} A \\ C \end{array} \right] + r(B) =  r[ \, A,  \ B \, ] + r(C). \hfill (21.80)
\cr}
$$

{\rm (c)}\, The rank of  $ BB^-A - AC^-C$ is invariant with respect to the choice of $B^- $ and $ C^-$ 
 if and only if 
$$\displaylines{
\hspace*{1.5cm}
r\left[ \begin{array}{cc} A & B  \\ C & 0 \end{array} \right] = r(B) + r(C), \ \ \ or \ \ \ 
r\left[ \begin{array}{cc} A & B  \\ C & 0 \end{array} \right] =r\left[ \begin{array}{c} A \\ C \end{array} \right] =  r[ \, A,  \ B \, ]. \hfill  (21.81)
\cr}
$$ }
{\bf Proof.}\, Follows from (21.75) and (21.76).  \qquad$ \Box$ 

\medskip
 
\noindent {\bf Theorem  21.24.}\, {\em Let $ A \in  {\cal F}^{m \times n},\, B \in {\cal F}^{k \times m}$  and
 $ C \in {\cal F}^{n \times l}$  be given.

{\rm(a)}\, The maximal and minimal ranks of $ B^-BA - ACC^-$ with
respect to $ B^-$ and $C^-$ are given by
$$  \displaylines{ 
\hspace*{1.5cm}
\max_{B^-,\, C^-} r( \, B^-BA - ACC^- \, ) = \min \{ \, r(BA) + r(AC),   \ \ m + n - r(B) - r(C)  \, \}, 
\hfill (21.82)
\cr
\hspace*{1.5cm}
\min_{B^-,\, C^-} r( \, B^-BA - ACC^- \, ) =  r(BA) + r(AC) - 2r(BAC). \hfill (21.83)
\cr}
$$ 

{\rm (b)}\, There are  $B^-$  and $C^-$ such that $ B^-BA = ACC^-$ if and only if $ A,$ $B$ and $C$ satisfy 
the rank equality $ r(BAC) = r(BA) = r(AC)$. 

{\rm (c)}\, The rank of  $ B^-BA - ACC^-$ is invariant with respect to the choice of $B^- $ and $ C^-$ 
 if and only if 
$$\displaylines{
\hspace*{1.5cm}
 BAC = 0 \ \ \  or  \ \ \  r(BAC) = r(B) + r(AC) - m = r(BA) + r(C) - n. \hfill (21.84)
\cr}
$$} 
\hspace*{0.3cm} The proof of  Theorem 21.24 is similar to that of Theorem 21.22 and is, therefore, omitted. Replacing $ A $ in  (21.82) and (21.83) by $A^{k-1}$,  and $ B $ and $ C $ in  (21.82) and (21.83) by
 $ A $, we directly obtain the following. 

\medskip

\noindent {\bf Corollary 21.25.}\, {\em Let $ A \in  {\cal F}^{m \times m}$ be given.

{\rm (a)}\, The maximal and  the minimal ranks of $ A^kA^- - A^-A^k$  with
respect to $ A^-$ are
$$ \displaylines{ 
\hspace*{2cm} 
\max_{A^-} r( \, A^kA^- - A^-A^k \, ) = \min \{ \, 2m - 2r(A),  \ \ 2r(A^k) \, \}, \hfill (21.85)
\cr
\hspace*{2cm} 
\min_{A^-} r( \, A^kA^- - A^-A^k \, ) = 2r(A^k)- 2r(A^{k+1}).  \hfill (21.86)
\cr}
$$ 

{\rm (b)}\, There is an $A^-$ such that $ A^kA^- - A^-A^k$ is nonsingular if and only if $ m$
 is even and  $r(A^k) = r(A) = m/2$.  

{\rm (c)}\, There is an $A^-$ such that $ A^kA^- = A^-A^k$ if and only if $r(A^{k+1}) = r(A^k).$   

{\rm (d)}\, The rank of  $ A^kA^- - A^-A^k$ is invariant with respect to the choice of $A^-$ if and only if 
$A^{k+1} = 0 $ or $r(A^{k+1}) = r(A^k) + r(A) - m$. }

\medskip

\noindent {\bf Theorem 21.26.}\, {\em Let $ A \in  {\cal F}^{m \times n},\,  B \in {\cal F}^{m \times p}$ and 
$ C \in  {\cal F}^{l \times n}$ be given. Then

{\rm (a)}\, The maximal and minimal ranks of $ I_m - AA^- - BB^-$ and  $ I_m - A^-A - C^-C$ with
respect to $ A^-, \ B^-$ and $C^-$ are
$$
\displaylines{ 
\hspace*{1.5cm} 
\max_{A^-, B^-} r( \, I_m - AA^- - BB^- \, ) = m - |\, r(A) - r(B) \, |, \hfill (21.87)
\cr
\hspace*{1.5cm} 
\min_{A^-, B^-} r( \, I_m -  AA^- - BB^- \, ) =  m + r(A) + r(B) - 2r[\, A, \ B \,].  \hfill (21.88)
\cr
\hspace*{1.5cm} 
\max_{A^-, C^-} r( \, I_n - A^-A - C^-C \, ) = n - |\, r(A) - r(C) \, |, \hfill (21.89)
\cr
\hspace*{1.5cm} 
\min_{A^-, B^-} r( \, I_n -  A^-A - C^-C \, ) =  n + r(A) + r(C) - 2r\left[ \begin{array}{c} A \\ C \end{array} \right].  \hfill (21.90)
\cr}
$$ 

{\rm (b)}\, There are $A^-$ and $B^-$ such that $ I_m - AA^- - BB^-$ is nonsingular if and only if 
$r(A) = r(B)$.  

{\rm (c)}\, There are  $A^-$ and $B^-$ such that $ AA^- + BB^- = I_m$ if and only if $r[\, A, \ B \,] = 
r(A) + r(B) = m$. 

{\rm (d)}\, There are $A^-$ and $C^-$ such that $ I_n - A^-A - C^-C$ is nonsingular if and only if 
$r(A) = r(C)$.  

{\rm (e)}\, There are  $A^-$ and $C^-$ such that $ A^-A + C^-C = I_n$ if and only if 
$r\left[ \begin{array}{c} A \\ C \end{array} \right] = r(A) + r(C) = n$. 
}

\medskip

\noindent {\bf Proof.}\,  According to (3.8), the rank of  $ I_m - AA^- - BB^-$ can be written as 
\begin{eqnarray*} 
r( \, I_m -  AA^- - BB^- \, ) & = & r(AA^-BB^-) + r(BB^-AA^-) - r(AA^-) - r(BB^-) + m \\
  & = & r(AA^-B) + r(BB^-A) - r(A) - r(B) + m.
\end{eqnarray*} 
Hence
$$ \displaylines{
\hspace*{1cm}
\max_{A^-,\, B^-} r( \,I_m -  AA^-  - BB^- \, )= \max_{A^-} r(AA^-B) + \max_{B^-} r(BB^-A) - r(A) - r(B) + m, 
\hfill (21.91)  
\cr
\hspace*{1cm}
\min_{A^-,\, B^-} r( \,I_m -  AA^-  - BB^- \, )= \min_{A^-} r(AA^-B) + \min_{B^-} r(BB^-A) - r(A) - r(B) + m. 
\hfill (21.92)  
\cr}
$$
It follows from Corollaries 21.7 and 21.9 that 
$$\displaylines{
\hspace*{2cm}
 \max_{A^-} r(AA^-B) = \max_{B^-} r(BB^-A) = \min \{ \,  r(A),  \ \ r(B) \, \}, \hfill
\cr
\hspace*{2cm}
 \min_{A^-} r(AA^-B) = \min_{B^-} r(BB^-A) = r(A) +  r(B) - r[\, A, \ B \,]. \hfill
\cr}
$$ 
Putting them in (21.91) and (21.92) yields (21.87) and (21.88). By the similar approach, we can get 
 (21.88) and (21.89).  \qquad $\Box$ 

\medskip

Replace $ B $ in (21.88) and (21.90) by $ I_m - A $, we get by noticing (1.11)  
 $$
\displaylines{ 
\hspace*{1.5cm} 
\min_{A^-, (I_m - A)^-} r[ \, I_m - AA^- - (I_m - A)(I_m - A)^- \, ] \hfill  
\cr
\hspace*{1.5cm} =
\min_{A^-, (I_m - A)^-} r[ \, I_m - A^-A - (I_m - A)^-(I_m - A) \, ] \hfill 
\cr
\hspace*{1.5cm}
= r(\, A - A^2 \,). \hfill (21.93)
\cr}
$$ 
Thus the following three statements are equivalent:
 
(a)  \ There are $ A^- $ and  $ (I_m - A)^-$  such that $ AA^- + (I_m - A)(I_m - A)^- = I_m$.

(b) \ There are $ A^- $ and  $ (I_m - A)^-$  such that $A^-A  + (I_m - A)^-(I_m - A) = I_m$.

(c)\,  $ A$ is idempotent. 

The two two statements in (a)  and (b) could be regarded as new characterizations of idempotent matrix.
 
Replace $ A $  and $ B $ in  (21.88) and (21.90) by $ I_m + A $ and  $ I_m - A $, respectively,  we get
 by noticing (1.12)  
$$
\displaylines{ 
\hspace*{1.5cm} 
\min_{(I_m + A)^-, (I_m - A)^-} r[ \, I_m - (I_m + A)(I_m + A)^-  - (I_m - A)(I_m - A)^- \, ] \hfill
\cr
\hspace*{1.5cm}
=\min_{(I_m + A)^-, (I_m - A)^-} r[ \, I_m -  (I_m + A)^- (I_m + A) -(I_m - A)^-(I_m - A) \, ]  \hfill
\cr
\hspace*{1.5cm}
 = r(\, I_m - A^2 \,). \hfill (21.94)
\cr}
$$ 
Thus the following three statements are equivalent:
 
(a)  \ There are $ (I_m + A)^-$ and  $ (I_m - A)^-$  such that $(I_m + A)(I_m + A)^- + (I_m - A)(I_m - A)^- = I_m$.

(b) \ There are $(I_m + A)^- $ and  $ (I_m - A)^-$  such that  $(I_m + A)^-(I_m + A)  + (I_m - A)^-(I_m - A) = I_m$.

(c)\,  $ A$ is involutory. 

The two two statements in (a) and (b) could be regarded as new characterizations of involutory matrix.

In general, suppose that $p(x)$ and $ q(x)$ are two polynomials without common roots. Then there is  
 by  (21.88), (21.90) and (1.17) the following  
$$
\min_{p^-(A),\, q^-(A)}r[ \, I_m - p(A)p^-(A)  -q(A)q^-(A)  \,] =
\min_{p^-(A),\, q^-(A)}r[ \, I_m - p^-(A)p(A)  -q^-(A)q(A)  \,]  = r[\, p(A)q(A) \, ]. \eqno (21.95)
$$ 
Thus there exist $p^-(x)$ and $ q^-(x)$  such that $p(A)p^-(A)  + q(A)q^-(A) = I_m$ if and only if
 $p(A)q(A) = 0$.  We leave its verification to the reader.

\medskip

\noindent {\bf Theorem 21.27.}\, {\em Let $ A \in  {\cal F}^{m \times m}$ be given. Then
$$\displaylines{
\hspace*{1cm}
\max_{A^-} r( \, I_m \pm  A^k - AA^- \, ) = \min_{A^-} r( \, I_m \pm  A^k - AA^- \, ) = r(A^{k+1}) - r(A) + m, \hfill (21.96)
\cr
\hspace*{1cm}
\max_{A^-} r( \, I_m \pm  A^k - A^-A \, ) = \min_{A^-} r( \, I_m \pm  A^k - A^-A \, ) = r(A^{k+1}) - r(A) + m. \hfill (21.97)
\cr}
$$ 
That is$,$ the equalities   
$$\displaylines{ 
\hspace*{1.5cm} 
r(A^{k+1}) =  r(A) - m  + r( \, I_m \pm  A^k - AA^- \, ) = r(A) - m  + r( \, I_m \pm  A^k -  A^-A \, )
 \hfill (21.98)
\cr}
$$
hold for any $ A^-$.
}

\medskip

\noindent {\bf Proof.}\,  Applying (21.4) and (21.5) to  $ I_m \pm  A^k - AA^-$ and  $I_m \pm  A^k - A^-A$  yields 
the desired results.  \qquad $\Box$ 

\medskip

We leave the proof of the following result to the reader.

\medskip

\noindent {\bf Theorem 21.28.}\, {\em Let $ A \in  {\cal F}^{m \times m}$ be given. Then
$$\displaylines{
\hspace*{1cm}
\max_{(I_m - A)^-} r\left[ \, (I_m - A)^- -  \sum_{i = 0}^{k -1} A^i \, \right] = \min \{ \,  m,  \ \ \ m + r( \, I_m - A^k \, ) - r(A) \, 
\},  \hfill (21.99)
\cr
and \hfill
\cr
\hspace*{1cm}
\min_{(I_m - A)^-} r\left[ \, (I_m - A)^- -  \sum_{i = 0}^{k -1} A^i \, \right] \hfill
\cr
\hspace*{1cm}
=\min_{(I_m - A)^-} r[ \, (I_m - A)(I_m - A)^- -  (\, I_m - A^k \,) \, ] \hfill
\cr
\hspace*{1cm}
=\min_{(I_m - A)^-} r[ \, (I_m - A)^-(I_m - A) -  (\, I_m - A^k \,)\, ] = r(\, A^k - A^{k+1} \,). \hfill (21.100)
\cr}
$$ 
In particular$,$ the following four statements are equivalent$:$

{\rm (a)}\, $ \sum_{i = 0}^{k -1} A^i \in \{ (I_m - A)^- \}.$
 
{\rm (b)}\, $ (I_m - A)(I_m - A)^- = I_m - A^k.$

{\rm (c)}\, $ (I_m - A)^-(I_m - A) = I_m - A^k.$
 
{\rm (d)}\, $ A^{k+1} = A^k,$ i.e., $ A $ is quasi-idempotent.
}

\medskip

A parallel result to (21.100) is 
$$
\displaylines{
\hspace*{1cm}
\min_{(\sum_{i = 0}^{k -1} A^i)^-} r\left[ \,  \left( \sum_{i = 0}^{k -1} A^i \right)^- - (I_m - A) \, \right]
= r(\, A^k + A^{k+1} + \cdots + A^{2k} \,). \hfill (21.101)
\cr}
$$ 
It  implies that $I_m - A \in \{ \, ( \sum_{i = 0}^{k -1} A^i)^- \}$ if and only if 
$ A^k + A^{k+1} + \cdots + A^{2k} = 0.$

It is expected that  one can further establish numerous rank equalities among $ A^-$, $(I - A)^-$ and 
 polynomials of $A$, and then derive from them various conclusions related $ A^-$ and $(I - A)^-$. we leave
this work to the reader.

In the remainder of this chapter,  we consider the rank  of the difference $  A^- - PN^-Q$ and  then 
present some of their consequences.     

\medskip

\noindent {\bf Theorem 21.28.}\, {\em Let $A \in {\cal F}^{m \times n}, \, P \in {\cal F}^{n \times p} \,  N \in {\cal F}^{q \times p}$ and $Q \in {\cal F}^{q \times n} $ be given. Then 
$$
\min_{A^-, \, N^-}r( \, A^- - PN^-Q \,) = r(\, N - QAP \,) + r(A) + r(N) -  r[\, A, \ QAP  \, ] - r
\left[ \begin{array}{c} A  \\ QAP \end{array} \right]. \eqno (21.102)
$$
In particular$,$  there are $ A^- $ and $ N^- $ such that $ A^- = PN^-Q$  if and only if 
$$
  r(\, N - QAP \,) = r \left[ \begin{array}{c} N  \\ QAP \end{array} \right]  + r[\, N, \  QAP \,]  
-  r(A) - r(N). \eqno (21.103)
$$}

\medskip

\noindent {\bf Proof.}\,  From (21.5) we can get the following simple result 
 $$
\min_{A^-}r( \, A^- - D \,) = r(\, A - ADA \,).
$$
Applying it to  $ A^- - PN^-Q$ we first  get 
$$
\min_{A^-}r( \, A^- - PN^-Q \,) = r(\, A - A PN^-Q A \,).
$$
Next applying (21.5) to its right hand side and simplifying, we then have 
$$
\min_{N^-}r(\, A - A PN^-Q A \,) = r(\, N - QAP \,) + r(A) + r(N) -  r[\, A, \ QAP  \, ] - r
\left[ \begin{array}{c} A  \\ QAP \end{array} \right]. 
$$
Combining the above  two equalities results in (21.102). \qquad $\Box$ 

\medskip

Clearly (21.29) could be regarded as a special case of (21.102).  Now applying 
(21.102) to the block matrix $ M = \left[ \begin{array}{cc} A & B  \\ C & 0 \end{array} \right]$, we can simply 
get following 
$$
\min_{A^-, \, M^-} r\left( A^-  - [\, I_n, \ 0 \,]M^- \left[\begin{array}{c} I_m  \\ 0 
\end{array} \right]  \right)
 = r(A) + r(M) + r\left[ \begin{array}{cc} 0  &  B \\ C & D \end{array} \right] - 
r \left[ \begin{array}{ccc} A & 0 & B \\ 0 & C & D \end{array} \right] -
r \left[ \begin{array}{cc} A & 0  \\ 0 & B \\  C & D \end{array}
\right]. \eqno (21.104)
$$
This result implies that there is an $ M^-$ which upper left block is an inner inverse of $ A $ if and only if 
$$
 r \left[ \begin{array}{ccc} A & 0 & B \\ 0 & C & D \end{array} \right] +
r \left[ \begin{array}{cc} A & 0  \\ 0 & B \\  C & D \end{array}
\right] =  r(A) + r(M) + r\left[ \begin{array}{cc} 0  &  B \\ C & D \end{array} \right]. \eqno (21.105)
$$
Applying (21.102) to a block circulant matrix $M$ generated by $ k $ matrices $ A_1, \, A_2, \, \cdots, \, 
 A_k$ and their sum $ A =  A_1 + A_2 + \cdots +A_k$, we can also get  
$$
\min_{A^-, \, M^-} r\left( A^-  -  \frac{1}{k}[\, I, \,\cdots \, I \,]M^-
\left[\begin{array}{c} I  \\ \vdots \\ I \end{array} \right]  \right) = 0, \eqno (21.106)
$$
that is, there must exist $ A^- $ and $M^-$ such that  
$$
A^-  = [\, I, \,\cdots,  \, I \,]M^-\left[\begin{array}{c} I  \\ \vdots \\ I \end{array} \right]. 
\eqno (21.107)
$$

Besides the result in Theorem 21.28, we can also determine the relationship
 between the two matrix sets $\{ PN^-Q \}$ and $\{ A^- \}$. Here we only list the 
 main results without proofs.   
  
\medskip
 
\noindent {\bf Theorem 21.29.}\, {\em Let $ A \in {\cal F}^{m \times n}, \,  
 N \in {\cal F}^{k \times l}, \, P \in {\cal F}^{n \times l}$  and $ Q \in  
{\cal F}^{k \times m}$ be given with $r(P) = n$ and $r(A) = m.$  Then  
$$  
\displaylines{ 
\hspace*{1.5cm}
\max_{N^-} r( \, A - APN^-QA \, ) = \min \{ \, r(A),   \ \ r(N - QAP)  + r(A) - r(N)
  \, \}, \hfill (21.108)
\cr}
$$ 
In particular$,$ the following set inclusion 
$$
 \{ PN^-Q \} \subseteq  \{ A^- \},  \eqno (21.109) 
$$     
holds if and only if    
$$
  r( N - QAP ) = r( N ) - r(A).  \eqno(21.110)
$$}
\noindent {\bf Theorem 21.30.}\, {\em Let $ A \in {\cal F}^{m \times n}, \,  
 N \in {\cal F}^{k \times l}, \, P \in {\cal F}^{n \times l}$  and $ Q \in  
{\cal F}^{k \times m}$ be given with $r(P) = n$ and $r(A) = m.$  Then  
$$
\displaylines{
\hspace*{1cm}
\max_{A^-}\min_{N^-}r( \, A^- - PN^-Q \, ) = \min \left \{ \, n + r(N) - 
r\left[ \begin{array}{c} N \\ P \end{array} \right], \ \ m + r(N) + r[\, N, \, Q \,], \right.  \hfill
\cr
\hspace*{4cm}
 \left. m + n + r(N - QAP) + r(N) - r(A) - r[\, N, \, Q \,] - 
r\left[ \begin{array}{c} N \\ P \end{array} \right] \right \}. \hfill (21.111)
\cr}
$$
In particular$,$ the following set inclusion 
$$
 \{ A^- \} \subseteq  \{ PN^- Q\},  \eqno (21.112) 
$$     
holds if and only if $ R(N) \cap R(Q) \neq \{ 0 \},$ or $ R(N^T) \cap R(P^T) 
\neq \{ 0 \},$ 
or   
$$
r( N - QAP ) = r\left[ \begin{array}{c} N \\ P \end{array} \right]  + r[\, N, \, Q \,] -    r( N ) + r(A) - m  - n.  \eqno(21.113)
$$}
\noindent {\bf Theorem 21.31.}\, {\em Let $ A \in {\cal F}^{m \times n}, \,  
 N \in {\cal F}^{k \times l}, \, P \in {\cal F}^{n \times l}$  and $ Q \in  
{\cal F}^{k \times m}$ be given with $r(P) = n$ and $r(A) = m.$  Under the 
condition $ R(N) \cap R(Q) = \{ 0 \} $ and $ R(N^T) \cap R(P^T) = \{ 0 \},$  the following equality 
$$
 \{ PN^-Q \} = \{ A^- \}  \eqno (21.114) 
$$     
holds if and only if 
$$
r\left[ \begin{array}{cc} N  & Q \\ P & 0  \end{array} \right] 
= r\left[ \begin{array}{c} N \\ P \end{array} \right]  + r(Q) 
 = r[\, N, \, Q \,]  + r(P), \eqno (21.115)    
$$
and 
$$ 
A = -[\, 0, \ I_m \,] \left[ \begin{array}{cc} N  & Q \\ P & 0  \end{array} \right]^-
 \left[ \begin{array}{c} 0 \\ I_n \end{array} \right]. \eqno (21.116)
$$}

Based on the above three theorems, we can establish the following several results.

\medskip

\noindent {\bf Theorem 21.32.}\,  {\em  Let $ A_1, \, A_2, \, \cdots, \, 
 A_k$ be given matrices of the same size$,$ and let $A =  A_1 + A_2 + \cdots +A_k$. Denote by $ M $
 the block circulant matrix generated by $ A_1, \, A_2, \, \cdots, \,  A_k$. Then 
 $ A $ and $ M $ satisfy 
$$ 
\{\, A^- \, \}  = \left\{  \frac{1}{k} [\, I, \,\cdots, \, I \,]M^-
\left[\begin{array}{c} I  \\ \vdots \\ I \end{array} \right]  \right\} 
. \eqno (21.117)
$$}
\noindent {\bf Theorem 21.33.}\,  {\em  Let $ A + iB $ be an $ m \times n $ complex matrix. Then 
$$ 
\{\, (A + iB)^- \, \}  = \left\{  \frac{1}{2} [\, I_n,  \, iI_n \,]
\left[\begin{array}{cr} A  & -B  \\ B & A  \end{array} \right]^-
\left[\begin{array}{c} I_m  \\  -iI_m \end{array} \right]  \right\} . \eqno (21.118)
$$}
\noindent {\bf Theorem 13.34.}\,  {\em  Let  $A_0 +  iA_1 + iA_2 +  kA_3$ be an $ m \times n $ real quaternion 
matrix.  Then 
$$ 
  \left\{  \, ( \,  A_0 +  iA_1 + iA_2 +  kA_3 \, )^-  \, \right\}   
= \left\{  \frac{1}{4}[ \, I_n , \  iI_n , \ jI_n, \  
kI_n \, ] \left[ \begin{array}{rrrr}  A_0 & -A_1 & -A_2 &  - A_3   \\  A_1 &  A_0 &  A_3 &  - A_2 
\\  A_2 &  -A_3 &  A_0  & A_1 \\ A_3 &  A_2 &  -A_1  & A_0 \end{array} \right]^- 
\left[ \begin{array}{c}  I_m \\ -iI_m \\ -jI_m  \\ -kI_m  \end{array} \right]  \right\} . \eqno (21.119) 
$$ 
}

\markboth{YONGGE  TIAN }
{22. GENERALIZED INVERSES OF MULTIPLE MATRIX PRODUCTS }

\chapter{ Generalized inverses of multiple matrix products}

\noindent Generalized inverses of products of matrices have been an attractive topic in the theory of generalized inverses matrices. 
Various results related to reverse order laws for g-inverses, reflexive g-inverses, and the Moore-Penrose inverses of matrix products can be found in the literature. Generally speaking,  this work has two main directions according the classification of generalized inverses of matrices, one of which is concerned with reverse order laws for the Moore-Penrose inverses of matrix products. Up till now, necessary and sufficient conditions for the reverse order laws $ (AB)^{\dagger} = B^{\dagger}A^{\dagger},$ \ $ (ABC)^{\dagger} = (BC)^{\dagger}B(AB)^{\dagger},$ \ $(ABC)^{\dagger} = C^{\dagger}B^{\dagger}A^{\dagger},$ \ $(A_1A_2 \cdots A_k)^{\dagger} = A_k^{\dagger} \cdots A_2^{\dagger}A_1^{\dagger}$ to hold  have well been established. The other direction of this work is concerned with reverse order 
laws for inner inverses, reflexive inner inverses, as well as several other types of generalized inverses of matrix products. 
Some earlier and recent work gives a complete consideration for the reverse order laws $ (AB)^- = B^- A^-$, \ $ (AB)_r^- = B_r^- A_r^- $, \ $(AB)^{\dagger} = B_{mr}^-A_{lr}^-$, and $(AB)_{MN}^{\dagger} = B^-A^-$ (see, e.g., \cite{BG1, DW,  RM, SS1, We1, We2, WHG}). In this chapter we first present some rank equalities 
related to $(BC)^-B(AB)^-$  and $C^-B^-A^-$, and then apply them to establish the relationships 
between $(BC)^-B(AB)^- $ and $ (ABC)^-,$  $ C^-B^-A^-$ and $ (ABC)^-.$

The following lemma comes directly from  (22.4) and (22.5), which will be used in the sequel. 

\medskip

\noindent {\bf Lemma 22.1.}\, {\em Let $ A \in {\cal F}^{ m \times n}, \, B \in 
{\cal F}^{ m \times k}, \, C \in {\cal F}^{ l \times n}$ and $ D \in {\cal F}^{ l \times k}$
 be given.

{\rm (a)}\,  If $R(D) \subseteq R(C)$ and $ R(D^T) \subseteq R(B^T),$ then 
$$ 
\displaylines{
\hspace*{2cm}
\max_{A^-}r( \,D - CA^-B \, )  = \min \left\{  r(B) , \ \ \ r(C),  \ \ \  r \left[ \begin{array}{cc}  A  & B \\  C & D   \end{array} \right] - r(A)   \right\}, \hfill (22.1)
\cr
\hspace*{2cm} 
\min_{A^-}r( \,D - CA^-B \, ) = r (A ) - r[ \, A, \ B \, ] - r \left[ \begin{array}{c} A  \\ C  \end{array} \right] + r\left[ \begin{array}{cc}  A  &  B  \\ C & D  \end{array} \right]. \hfill (22.2)
\cr}
$$
 
{\rm (b)}\, In particular$,$ 
$$ \displaylines{
\hspace*{2cm}
\max_{A^-}r( \,D - CA^- \, )  = \min \left\{ \, m,  \ \ \ r[\, C, \ D \, ] , \ \ \ r( \, C- DA \, ) - r(A) + m  \, \right\}, \hfill (22.3)
\cr
\hspace*{2cm}
\max_{A^-}r( \,D - A^-B \, )  = \min \left\{ \,  n, \ \ \  r \left[ \begin{array}{c} B  \\ D  \end{array} \right] , \ \ \  
r( \, B- AD \, ) - r(A) + n  \, \right\}. \hfill (22.4)
\cr}
$$ }
{\bf Proof.}\, Follows immediately from (21.5) and (21.6). \qquad $\Box$

\medskip

It is quite obvious  that there are 
$(AB)^- \in \{ \, (AB)^- \}$ and $(BC)^- \in \{ \, (BC)^- \}$ such that $ (BC)^-B(AB)^-  
$ $\subseteq \{ \, (ABC)^- \}$ if and only if  
$$ \displaylines{
\hspace*{2cm}
\min_{(AB)^-, \,(BC)^-} r[ \, ABC - (ABC)(BC)^-B(AB)^-(ABC) \, ] = 0,  \hfill (22.5) 
\cr}
$$ 
and the set inclusion $ \{ (BC)^-B(AB)^- \} \subseteq \{ (ABC)^- \}$ holds if and only if  
$$ \displaylines{
\hspace*{2cm}
\max_{(AB)^-,\, (BC)^-} r[ \, ABC - (ABC)(BC)^-B(AB)^-(ABC) \, ] = 0.  \hfill (22.6) 
\cr}
$$ 
These two equivalence statements clearly show that the relationship between $(BC)^-B(AB)^-$  and 
$ (ABC)^-$ can be characterized by extreme ranks of matrix expressions involving 
$(BC)^-,  (AB)^-$,  and $ (ABC)^-$. According to the rank formulas in Lemma 22.1, we  
easily find the following special result. 

\medskip

\noindent {\bf Theorem 22.2.}\, {\em Let $ A \in {\cal F}^{ m \times n}, \, B \in {\cal F}^{ n \times p},$ and
 $ C \in {\cal F}^{ p \times q}$ be given$,$ and let $ M = ABC$.  Then 
$$ 
\min_{(AB)^-} r[ \, M- M(BC)^-B(AB)^-M \, ] = 0,  \ \ \ \ \min_{(BC)^-} r[ \, M - M(BC)^-B(AB)^-M \, ] = 0,  \eqno (22.7) 
$$ 
and
$$ 
\max_{(AB)^-, \, (BC)^-} r[ \, M - M(BC)^-B(AB)^-M \, ] = \min \left\{ \, r(M),  \ \ 
\ r(M)- r(AB) - r(BC) + r(B) \, \right\}. \eqno (22.8)
$$ } 
{\bf Proof.}\, Applying  (22.2) to $M - M(BC)^-B(AB)^-M$, we find 
\begin{eqnarray*}
\lefteqn{ \min_{(AB)^-} r[ \, M - M(BC)^-B(AB)^-M \, ] } \\ 
& = & r(AB) - r[\, AB, \ M \,] - r \left[ \begin{array}{c} AB  \\
 M(BC)^-B  \end{array} \right] +  r \left[ \begin{array}{cc} AB  & M  \\ M(BC)^-B & M \end{array} \right] \\
 & = & r(AB) - r[ \, AB, \ 0 \, ] - r \left[ \begin{array}{c} AB  \\ M(BC)^-B  \end{array} \right] +  
r \left[ \begin{array}{cc} AB  & 0 \\ M(BC)^-B & 0 \end{array} \right]  = 0, 
\end{eqnarray*}
\begin{eqnarray*}
\lefteqn{ \min_{(BC)^-} r[ \,M - M(BC)^-B(AB)^-M \, ] }\\
& = & r(BC) - r[\, BC, \ B(AB)^-M \,] - 
r\left[ \begin{array}{c} BC  \\ M \end{array} \right] +  r \left[ \begin{array}{cc} BC  & B(AB)^-M  \\ M & M \end{array} \right] \\ 
& = & r(BC) - r[ \, BC, \  B(AB)^-M  \, ] - r \left[ \begin{array}{c} BC  \\ 0 \end{array} \right] +  
r \left[ \begin{array}{cc} BC  &  B(AB)^-M \\ 0 & 0 \end{array} \right]  = 0. 
\end{eqnarray*}
Both of them are (22.7).  Next by (22.1), we find
\begin{eqnarray*}
\max_{(AB)^-} r[ \, M - M(BC)^-B(AB)^-M \, ] & = &\min \left\{   r[M(BC)^-B], \ \  r(M),  \ \  r \left[ \begin{array}{cc} AB & M \\
 M(BC)^-B & M \end{array} \right] - r(AB)  \right\} \\
& = &\min \left\{ r(M),  \ \ \ r \left[ \begin{array}{c} AB \\ M(BC)^-B \end{array} \right] -r(AB) \right\}. 
\end{eqnarray*}
Accroding to  (21.4) and $ r(M) \leq r(BC)$, we also find that 
\begin{eqnarray*}
\max_{(BC)^-} r \left[ \begin{array}{c} AB \\ M(BC)^-B \end{array} \right] & = & 
\max_{(BC)^-} r \left( \left[ \begin{array}{c} AB \\ 0 \end{array} \right] - \left[ \begin{array}{c} 0  \\ -M \end{array} \right](BC)^-B \right) \\
& = & \min  r\left\{ r\left[ \begin{array}{cc} AB  & 0  \\ 0 & -M \end{array} \right], \ \ 
 r\left[ \begin{array}{c} AB \\ 0 \\ B \end{array} \right], \ \ r \left[ \begin{array}{cc} BC  & B  \\ 0 & AB \\ -M & 0 \end{array} \right] -r(BC) \right\} \\
& = & \min \{  \, r(AB) + r(M), \ \ r(B),  \ \ r(B) + r(M) - r(BC) \, \}  \\  
& = & \min \{  \, r(AB) + r(M), \ \ r(B) + r(M) - r(BC) \, \}.    
\end{eqnarray*}
Combining the above two equalities, we obtain 
\begin{eqnarray*}
\max_{(AB)^-, \, (BC)^-} r[ \, M - M(BC)^-B(AB)^-M \, ] & = & \min \left\{ \, r(M), \ \ \max_{(BC)^-}r \left[ \begin{array}{c} AB \\ M(BC)^-B \end{array} \right] -r(AB) \, \right\} \\
& = & \min \left\{ \, r(M), \ \  r(M)- r(AB) - r(BC) + r(B) \, \right\},
\end{eqnarray*}
which is exactly (22.8). \qquad  $ \Box$.  

\medskip

Combining  (22.5) and (22.6) with  (22.7) and (22.8), we obtain the main result in this chapter. 

\medskip
 
\noindent {\bf Theorem 22.3.}\, {\em Let $ A \in {\cal F}^{ m \times n}, \, B \in {\cal F}^{ n \times p},$ and
 $ C \in {\cal F}^{ p \times q}$ be given.   

{\rm (a)}\, For every $ (AB)^- \in \{ \, (AB)^- \},$  there must be a $ (BC)^- \in \{ \,(BC)^-\} $ such that $ (BC)^-B(AB)^-  \in \{\, (ABC)^- \}$ holds.   

{\rm (b)}\, For every $ (BC)^- \in \{\, (BC)^- \},$  there must be an $ (AB)^- \in \{ \,(AB)^-\} $ such that $ (BC)^-B(AB)^-  \in \{ \, (ABC)^- \}$ holds. 

{\rm (c)}\, The set inclusion  $ \{ \,(BC)^-B(AB)^- \}  \subseteq \{ \,(ABC)^- \}$ holds if and only if 
$$\displaylines{
\hspace*{2cm}
 ABC = 0  \ \  or  \ \  r( ABC) = r(AB) + r(BC ) - r(B). \hfill (22.9)
\cr}
$$  

{\rm (d)}\, In particular$,$ if $ r(ABC) = r(B),$  then $ \{ \,(BC)^-B(AB)^- \}  \subseteq \{ \,(ABC)^- \}$ holds. } 

\medskip

As a direct consequence by setting $ B = I$ in Theorem 22.3, we obtain the following.

\medskip
 
\noindent {\bf Corollary 22.4.}\, {\em Let $ A \in {\cal F}^{ m \times n}$ and $ B \in {\cal F}^{ n \times p}$  be  given.   

{\rm (a)}\, For every $ A^- \in \{ \,A^- \},$  there must be a $ B^- \in \{ \, B^-\} $ such that $ B^-A^-  \in \{  \,(AB)^- \}$ holds.   

{\rm (b)}\, For every $ B^- \in \{ \,B^- \},$  there must be an $ A^- \in \{ \, A^-\} $ such that $ B^-A^-  \in \{ \, (AB)^- \}$ holds. 

{\rm (c)}\, The set inclusion  $ \{ \, B^-A^- \}  \subseteq  \{ \,(AB)^- \}$ holds if and only if 
$$\displaylines{
\hspace*{2cm}
 AB = 0  \ \ or  \ \ r( AB) = r(A) + r(B) - n. \hfill (22.10)
\cr}
$$ }  
\hspace*{0.3cm} Necessary and sufficient conditions for $ \{ \,B^-A^- \} \subseteq \{ \, (AB)^- \}$ were previously examined by 
Gross in \cite{Gro2}, Werner in  \cite{We1} and \cite{We2}.  The results given there are in fact equivalent to (22.10)  
The results in  Theorem 22.3 and Corollary 22.4 can help us to establish various relationship between generalized inverses of matrices. We next present one of them. 

\medskip

\noindent {\bf Corollary 22.5.}\, {\em Let $ A, \, B  \in {\cal F}^{ m \times n}$ be  given and let
 $ M = {\rm diag}(\,  A, \,  B \, ),$ $ N = A + B,$ $  S =  \left[ \begin{array}{c} A  \\ B  \end{array} \right],$ and
 $ T = [ \, A, \ B \, ].$ Then 
$$ 
\max_{ S^-, \, T^-} r \left(  N - N\left[ \begin{array}{c} A  \\ B  \end{array} \right]^- M[ \, A, \ B \, ]^-N \right) = \min \left\{ r(N),  \ \ \ r(N)- r\left[ \begin{array}{c} A  \\ B  \end{array} \right] 
- r [ \, A, \ B \, ] + r(A) + r(B)  \right\}. \eqno (22.11)
$$
In particular$,$  the set inclusion 
$$ \displaylines{
\hspace*{2cm}
\left\{ \left[ \begin{array}{c} A  \\ B  \end{array} \right]^- \left[ \begin{array}{cc} A  & 0 \\  0 & B  \end{array} \right][ \, A, \ B \, ]^-  \right\} 
\subseteq \{  \, (\, A + B \, )^-  \, \} \hfill (22.12) 
\cr}
$$ 
holds if and only if 
$$\displaylines{
\hspace*{2cm}
 A + B  = 0 \ \  or \ \ r( \, A + B \, ) = r \left[ \begin{array}{c} A  \\  B \end{array} \right] + r [ \, A, \ B \, ] - r(A) - r(B). 
\hfill (22.13) 
\cr}
$$} 
{\bf Proof.}\, Writing $A + B $ as the product 
$$ 
\displaylines{
\hspace*{2cm}
A + B = [\, I_m , \ I_m \, ] \left[ \begin{array}{cc} A  & 0 \\  0 & B  \end{array} \right] \left[ \begin{array}{c} I_n  \\ I_n  \end{array} \right] = PDQ. \hfill
\cr}
$$
Then (22.11) follows from  (22.8). Consequently (22.12) and (22.13) 
follow form  (22.11).  \qquad $\Box$ 

\medskip

Eq.\,( 22.12) has some interesting consequences. For example, in the case of $ N = A + ( I_m - A) = I_m$, then 
$$ 
\max_{ S^-, \, T^-} r \left(  I_m  - \left[ \begin{array}{c} A  \\ I_m - A   \end{array} \right]^- 
\left[ \begin{array}{cc} A  & 0  \\ 0   & I_m - A  \end{array} \right] [ \, A, \ I_m - A  \, ]^- \right) 
= r(A) + r( I_m - A)  - m = r( A - A^2). 
$$
This implies that $ A $ is idempotent if and only if 
$$
\left[ \begin{array}{c} A  \\ I_m - A   \end{array} \right]^- 
\left[ \begin{array}{cc} A  & 0  \\ 0   & I_m - A  \end{array} \right] [ \, A, \ I_m - A  \, ]^- = I_m
$$
holds for any inner inverses in it. This fact could be regarded as a characterization for idempotent matrix.
In addition for any  two idempotent matrices $ A $ and $ B$, there is 
$$ 
\max_{ S^-, \, T^-} r \left(  (A - B)   - (A - B) \left[ \begin{array}{c} A  \\  -B    \end{array} \right]^- 
\left[ \begin{array}{cc} A  & 0  \\ 0   &  -B  \end{array} \right] [ \, A, \ -B  \, ]^-( A -B)  \right) 
= 0. 
$$  
Thus    
$$
\left\{ \left[ \begin{array}{c} A  \\  -B  \end{array} \right]^- 
\left[ \begin{array}{cc} A  & 0 \\  0 & -B  \end{array} \right][ \, A, \ -B \, ]^-  \right\} 
\subseteq \{  \, (\, A - B \, )^-  \, \}  
$$  
holds for any two idempotent matrices $ A $ and $ B$.

An extension of Corollary 22.5 is given below. 

\medskip

\noindent {\bf Corollary 22.6.}\, {\em Let $ A_1, \, A_2, \, \cdots, \, A_k \in {\cal F}^{ m \times n}$ be  given.  The  set inclusion 
$$ 
\left\{ \, \left[ \begin{array}{c} A_1  \\ \vdots \\ A_k \end{array} \right]^- \left[ \begin{array}{ccc} A_1  & & \\  
 & \ddots  & \\ & & A_k \end{array} \right][ \, A_1, \, \cdots, \,  A_k \, ]^-  \ \right\}  \subseteq  \{  \ (\, A_1 + \cdots  + A_k \, )^-  \ \} 
$$ 
holds if and only if 
$$
A_1 + \cdots  + A_k  = 0 \ \ \  or \ \ \ r(\, A_1 + \cdots  + A_k \,) = r \left[ \begin{array}{c} A_1  \\  \vdots \\ A_k \end{array} \right] +
 r [ \, A_1, \, \cdots, \, A_k  \, ] - r(A_1) - \cdots - r(A_k).
$$ } 
\hspace*{0.3cm} For some products of matrices,  the two equalities in (22.10) and (22.11) are satisfied, for example, (1.10)---(1.12), (1.14), (3.1), (4.1) and so on could be regarded as the  special 
cases of (22.10) and (22.11).  Thus based on them and Corollary 22.4(c) and Corollary 
22.5(c), one can establish various set inclusions for inner inverses of products of matrices. We leave them to the reader.  

Without much effort, we can also find a necessary and sufficient condition for the set inclusion $ \{ \,C^-B^- A^-\}$
 $\subseteq \{ \, (ABC)^- \}$ to hold. A rank formula related to this inclusion can be established using the rank formulas in (22.3) and (22.4). 

\medskip

\noindent {\bf Theorem 22.7.}\, {\em Let $ A \in {\cal F}^{ m \times n}, \, B \in {\cal F}^{ n \times p}$ and $ C \in {\cal F}^{ p \times q}$ be given, and let $ M = ABC$.  Then 
$$ \displaylines{
\hspace*{0.5cm}
\max_{A^-, \, B^-, \, C^-} r(\, M - MC^-B^-A^-M \, ) = \min \left\{ \, r(M),  \ \ \ r(M)- r(A) - r(B) - r(C) + n + p  \, \right\}. \hfill (22.14)
\cr}
$$ } 
{\bf Proof.}\, We determine the maximal rank of $  M - MC^-B^-A^-M $ subject to  $ A^-, \ B^- $,  and $ C^-$ through the following step   
$$ 
\displaylines{
\hspace*{1cm}
\max_{A^-, \, B^-, \, C^-} r( \, M - MC^-B^-A^-M \, ) =\max_{C^-}\max_{B^-}\max_{A^-}r( \, M - MC^-B^-A^-M \, ). 
\hfill (22.15)
\cr}
$$ 
According to  (21.4) we first find 
$$
\displaylines{
\hspace*{1.5cm}
\max_{A^-}r(\, M - MC^-B^-A^-M \, ) \hfill
\cr
\hspace*{1.5cm}
= \min \left\{ \, r[ \, M, \ MC^-B^- \,] , \ \ \ r \left[ \begin{array}{c}  M  \\ M  \end{array} \right],  \ \ \  r \left[ \begin{array}{cc}  A  & M \\  MC^-B^- & M  \end{array} \right] - r(A) \ \right\} \hfill
\cr
\hspace*{1.5cm}
=  \min \left\{ \ r(M) ,  \ \ \  r(\, A -  MC^-B^- \, )  + r(M) - r(A) \, \right\}.
\hfill (22.16)
\cr}
$$
Next applying (22.3) to $ A - MC^-B^-$ and noticing that  $ r(A) \leq n$, we obtain
\begin{eqnarray*}
\max_{B^-}r( \, A - MC^-B^- \,) & = & \min \left\{ \, r[ \, A, \ MC^- \,] , \ \ \ r(\, MC^- - AB \,) - r(B) + n , \ \ \  n  \ \right\} \\
& = & \min \left\{ \ r(A) , \ \ \ r(\, MC^- - AB \,) - r(B) + n , \ \ \  n  \, \right\} \\
& = & \min \left\{ \ r(A) , \ \ \ r(\,  AB - MC^- \,) - r(B) + n \ \right\}.
\end{eqnarray*}
Consequently applying  (22.3) to $  AB - MC^-$ and the noticing that $r(AB) \leq p$, we further find 
\begin{eqnarray*}
\max_{C^-}r( \, AB - MC^- \,) & = & \min \left\{ \, r[ \, AB, \ M \,] , \ \ \ r(\, M - ABC \,) - r(C) +  p , \ \ \  p \ \right\} \\
& = & \min \left\{ \, r(AB) , \ \ \ p - r(C) \, \right\}. 
\end{eqnarray*}
Putting the above three results in (22.16) and noticing the Sylvester's law $ r(AB) \geq r(A) + r(B) - n,$  we eventually obtain 
\begin{eqnarray*}
\lefteqn{\max_{A^-, \, B^-,\, C^-} r[ \, M - MC^-B^-A^-M \, ] }\\
& = & \min \left\{ \, r(M), \ \ \max_{C^-} \max_{B^-} r( \, A - MC^-B^- \,) + r(M) - r(A) \, \right\} \\
& = & \min \left\{ \, r(M), \ \ \max_{C^-} r( \, AB - MC^- \,) + r(M) - r(A) -r(B) + n  \, \right\} \\
& = & \min \left\{ \, r(M), \ \ r(M) + r(AB) - r(A) -r(B) + n,  \ \  r(M) - r(A) -r(B) -r(C) + n +p \, \right\} \\
& = & \min \left\{ \, r(M), \ \ \  r(M) - r(A) -r(B) -r(C) + n +p \, \right\},
\end{eqnarray*}
establishing (22.14).  \qquad  $\Box$ 

\medskip

It is quite obvious that the set inclusion $\{ \, C^-B^-A^- \} \subseteq \{ \, (ABC)^-\}$ holds if and only if 
$$ \displaylines{
\hspace*{1.5cm}
\max_{A^-, \, B^-, \, C^-} r( \, M - MC^-B^-A^-M \, ) =0. \hfill
\cr}
$$ 
Thus from Theorem 3.1, we immediately obtain the following.

\medskip
 
\noindent {\bf Theorem 22.8.}\, {\em Let $ A \in {\cal F}^{ m \times n}, \, B \in {\cal F}^{ n \times p},$ and $ C \in {\cal F}^{ p \times q}$ be given.  Then the set inclusion $\{\, C^-B^-A^- \} \subseteq \{ \, (ABC)^-\}$ holds if and only if 
$$ 
\displaylines{
\hspace*{1.5cm}
ABC = 0 \ \ or  \ \ r(ABC) =  r(A) + r(B) + r(C) - n - p. \hfill (22.17)
\cr}
$$ }  
\hspace*{0.3cm} The results in Theorems 22.7 and 22.8 can easily be extended to inner inverse of  multiple matrix products and its proof is omitted.

\medskip

\noindent {\bf Theorem 22.9.}\, {\em Let $ A_1 \in {\cal F}^{ n_1 \times n_2}, \, A_2 \in {\cal F}^{ n_2 \times n_3}, \, \cdots, \,A_k \in {\cal F}^{ n_k \times n_{k+1}}$ be given$,$ and denote $ M = A_1A_2 \cdots A_k$.  Then 
$$ 
\max_{A^-_1,\,\cdots,\,A_k^-} r( \, M - MA_k^- \cdots A_1^-M \,) = \min \left\{ \, r(M),  \ \ \ r(M)- r(A_1) - \cdots - r(A_k) + n_2 + \cdots + n_k  \, \right\}.  \eqno (22.18)
$$ }    
{\bf Theorem 22.10.}\, {\em Let  $ A_1 \in {\cal F}^{ n_1 \times n_2}, \, 
A_2 \in {\cal F}^{ n_2 \times n_3}, \, \cdots, \, A_k \in {\cal F}^{ n_k \times n_{k+1}}$ be given. Then the set inclusion 
$$
\{ \, A_k^- \cdots A_2^-A_1^-  \, \} \subseteq \{ \, (\, A_1A_2 \cdots A_k \,)^- \} \eqno (22.19)
$$
holds if and only if 
$$ 
A_1A_2 \cdots A_k = 0 \ \ or \ \ r( \, A_1A_2 \cdots A_k \,) = r(A_1) + r(A_2) + \cdots  +r(A_k) - n_2 - n_3  -\cdots - n_k.  \eqno (22.20)
$$} 
\hspace*{0.3cm} Combining Theorem 22.10 and the rank equality (1.16), we then get the following interesting result. 

\medskip
\noindent {\bf Theorem 22.11.}\, {\em Let $ A \in {\cal F}^{m \times m}$ be given$,$ $\lambda_1,  \, \lambda_2, \, 
 \cdots, \, \lambda_k \in {\cal F}$ with $ \lambda_i \neq  \lambda_j$  for $ i \neq j, $ 
and denote $M = (\lambda_1I - A)^{t_1} ( \lambda_2I - A)^{t_2} \cdots 
( \lambda_kI - A)^{t_k},$ where $t_1, \, t_2,  \, \cdots, \, t_k$ are  any positive integers. Then
 the following set inclusion holds
$$
\left\{ \, [\,( \lambda_kI - A)^{t_k}]^-  \cdots  [\, ( \lambda_2I - A)^{t_2}]^- [\,( \lambda_1I - A)^{t_1} \,]^-  \, \right\} 
\subseteq \left\{ \, M^- \, \right\}. 
\eqno (22.21)
$$
In general$,$ the following set inclusion
$$
\left\{ \, [\,( \lambda_{i_1} I - A)^{t_{i_1}}]^- [\, ( \lambda_{i_2}I - A)^{t_{i_2}}]^- 
\cdots [\,( \lambda_{i_k} I - A)^{t_{i_k}} \,]^-  \, \right\} 
\subseteq \left\{ \, M^- \, \right\}. 
\eqno (22.22)
$$
also holds$,$ where $ i_1, \,  i_2, \, \cdots, \, i_k$ are any permutation of $ 
1, \,  2, \, \cdots, \, k$.     
}

\medskip
As a special consequence, we see from (22.22) that the following six set inclusions all hold
$$
\left\{ \, A^-(I_m - A)^- \, \right\} \subseteq \left\{ \, (A - A^2)^- \, \right\},  \qquad 
\left\{ \, (I_m - A)^-A^-   \, \right\} \subseteq \left\{ \, (A - A^2)^- \, \right\}, \eqno (22.23)
$$
$$
\left\{ \, (I_m - A)^-(I_m + A)^-   \, \right\} \subseteq \left\{ \, (I_m - A^2)^- \, \right\},  \  \ 
\left\{ \, (I_m + A)^-(I_m - A)^-   \, \right\} \subseteq \left\{ \, (I_m - A^2)^- \, \right\}, \eqno (22.24)
$$
$$
\left\{ \, A^-(I_m - A)^-(I_m + A)^-   \, \right\} \subseteq \left\{ \, (A - A^3)^- \, \right\},  \ \ 
\left\{ \, A^-(I_m + A)^-(I_m - A)^-   \, \right\} \subseteq \left\{ \, (A - A^3)^- \, \right\}, \eqno (22.25)
$$
$$
\left\{ \, (I_m - A)^- A^-(I_m + A)^-   \, \right\} \subseteq \left\{ \, (A - A^3)^- \, \right\},  \ \
\left\{ \, (I_m + A)^-A^-(I_m - A)^-   \, \right\} \subseteq \left\{ \, (A - A^3)^- \, \right\}, \eqno (22.26)
$$
$$
\left\{ \, (I_m - A)^-(I_m + A)^-A^-  \, \right\} \subseteq \left\{ \, (A - A^3)^- \, \right\},  \ \
\left\{ \, (I_m + A)^-(I_m - A)^-A^-   \, \right\} \subseteq \left\{ \, (A - A^3)^- \, \right\}. \eqno (22.27)
$$

\markboth{YONGGE  TIAN }
{23. GENERALIZED INVERSES OF SUMS OF  MATRICES}

\chapter{Generalized inverses of sums of matrices}

\noindent In this chapter we  establish some rank equalities related for
 sums of inner inverses of matrices and then use them to deal with
  the following several problems:

\begin{enumerate}
\item[(I)]  The relationship between  $ A^- + B^- $ and $ (\, A + B\, )^-$.

\item[(II)]  The relationship between  $A_1^- + A_2^-  + \cdots + A_k^- $ and $ ( \, A_1 + A_2  + \cdots + A_k \, )^-$.

\item[(III)] The relationship between   $ \{ \, A^- + B^- \, \} $ and $ \{ \, C^- \, \}$.

\item[(IV)] The relationship between  $ \{ \, A_1^- + A_2^-  + \cdots + A_k^- \, \} $ and $ \{ \, C^- \, \}$. 
\end{enumerate}

\medskip

We first present a formula for the dimension of the intersection of $k$
matrices, which will be applied in the sequel.

\medskip

\noindent {\bf  Lemma 23.1}\cite{Ti8}.\, {\em Let $ [\, A_1,  \, A_2, \
\cdots \, A_k \, ] \in {\cal F}^{m \times n}$. Then
$$ \displaylines{
\hspace*{2cm}
{\rm dim} [ \, R(A_1) \cap R(A_2) \cap \cdots \cap  R(A_k) \, ]
= r(N) + r(Q) - r[\, N, \ Q \, ],  \hfill (23.1)
\cr}
$$
where $ N = {\rm diag}( \, A_1,  \, A_2, \, \cdots \, A_k \, ), \
 Q = [\, I_m,  \ I_m, \ \cdots \ I_m \, ]^T.$
In particular$,$ 
$$\displaylines{
\hspace*{2cm}
R(A_1) \cap R(A_2) \cap \cdots \cap R(A_k) = \{ 0 \} \Leftrightarrow  R(N) \cap R(Q) = \{ 0 \}. \hfill (23.2) 
\cr}
$$ } 
{\bf Proof.}\, Let $ X \in {\cal F}^{m \times t}$ be a matrix
satisfying $ R (X) = \cap_{i = 1}^k R(A_i)$. The this $ X $  can be
written as $X = A_1X_1 =  A_2X_2 =  \cdots = A_kX_k.$ Consider it as a
system of matrix equations. It  can equivalently be written as
$$
\displaylines{
\hspace*{2cm}
\left[ \begin{array}{ccccc} I_m & - A_1 & & &  \\ I_m & & - A_2 & &  \\ \vdots & & & \ddots & \\
 I_m &  & & & -A_m  \end{array} \right] \left[ \begin{array}{c} X \\ X_1 \\ \vdots \\ X_k \end{array} \right]  = 0,
\hfill
\cr}
$$ 
or briefly $ [ \, Q, \ -N \,] Y = 0$. Solving for $ X$, we obtain  its general solution is
$$ \displaylines{
\hspace*{2cm}
X = [ \, I_m, \ 0 \,]( \, I -  [ \, Q, \ -N \,]^-[ \, Q, \ -N \,] \, )V, \hfill
\cr}
$$     
where $ V $ is arbitrary. The maximal rank of $X$, according to (1.3) is 
\begin{eqnarray*}
r(X ) & = &  r (\, [ \, I_m, \ 0 \,]( \, I -  [ \, Q, \ -N \,]^-[ \, Q, \
-N \,] \, ) \, ) \\
& = & r \left[ \begin{array}{cc} I_m & 0 \\ Q & -N \end{array} \right] - r[ \, Q, \ -N \,]  \\
& = &  m + r(N) - r[ \, Q, \ -N \,] =  r(N) + r(Q) - r[\, N, \ Q \, ], 
 \end{eqnarray*} 
which is the dimension of $ \cap_{i = 1}^k R (A_i)$. \qquad  $\Box$ \\

\noindent{ \Large 23.1. \ {\bf The relationships between $ A^- + B^- $ and  $(\, A + B \,)^-$} } \\

\noindent We first establish some rank equalities related to $  A + B $ and
$ A^- + B^- $.

\medskip  
 
\noindent {\bf Theorem 23.2.}\, {\em  Let $A, \, B \in {\cal F}^{m \times n}$
 be given$,$ and let $ M = A + B.$ Then 
$$ 
\displaylines{
\hspace*{0.5cm}
\max_{A^-,\, B^-}r[ \, M -  M( \, A^- + B^- \, )M \,] = \min \left\{ r(M), \ \ \ 
r \left[ \begin{array}{cc} M   & A \\ B  & M  \end{array} \right]  + r(M) - r(A) - r(B)  \right\},
 \hfill (23.3) 
\cr
\hspace*{0cm}
and \hfill
\cr
\hspace*{0.5cm}
\min_{A^-,\, B^-}r[ \, M -  M( \, A^- + B^- \, )M \,] \ \hfill
\cr
\hspace*{1cm}
 = r(A) + r(B) + r(M)  +  \ r \left[ \begin{array}{cc} M   & A \\ B  & M  \end{array} \right]  
- r \left[ \begin{array}{ccc} A  & 0 & B \\ 0 & B & A    \end{array} \right] - 
 r \left[ \begin{array}{cc} A  & 0 \\  0  & B \\ B & A   \end{array} \right]. \hfill (23.4) 
\cr}
$$}
{\bf Proof.}\, We first show that 
$$ 
\displaylines{
\hspace*{1.5cm}
\{ \, A^- + B^- \, \}  = \left\{  [ \, I_n, \ I_n \,]N^- \left[ \begin{array}{c} I_m \\ I_m
 \end{array} \right]  \right\}, \ \ {\rm where}  \ \   N = \left[ \begin{array}{cc} A  & 0 \\ 0  & B  \end{array} \right]. 
 \hfill (23.5) 
\cr}
$$
In fact, the general expression of $ A^- + B^-$ can be written as 
$$\displaylines{
\hspace*{1.5cm} 
A^- + B^- = A^{\sim} + B^{\sim} + F_AV_1 + V_2E_A + F_BW_1 + W_2E_B, \hfill (23.6)
\cr}
$$ 
where $ A^{\sim}$ and $ B^{\sim}$ are two special inner inverses of $ A $
and $ B $, $ V_1,$ $ V_2,$ $ W_1,$ and  $ W_2$ are arbitrary. The general expression of $ N^- $ is
\begin{eqnarray*} 
N^- &  =  & N^{\sim} + F_NS + TE_N \\ 
& = & \left[ \begin{array}{cc} A^{\sim}   & 0 \\ 0  & B^{\sim}   \end{array} \right]+ \left[ \begin{array}{cc} F_A   & 0 \\ 0  & F_B \end{array} \right] \left[ \begin{array}{cc} S_1   & S_2 \\ S_3  & S_4 \end{array} \right] + \left[ \begin{array}{cc} T_1   & T_2 \\ T_3  & T_4 \end{array} \right] \left[ \begin{array}{cc} E_A   & 0 \\ 0  & E_B \end{array} \right] \\
& = & \left[ \begin{array}{cc} A^{\sim} + F_AS_1 + T_1E_A    &  F_AS_2 + T_2E_B  \\  F_BS_3 + T_3E_A  & 
 B^{\sim} + F_BS_4 + T_4E_B    \end{array} \right], 
\end{eqnarray*} 
where $ S_1$---$S_4$ and  $ T_1$---$T_4$ are arbitrary. In that case, we
have the general expression
\begin{eqnarray*}  
[ \, I_n, \ I_n \,]N^- \left[ \begin{array}{c} I_m  \\ I_m  \end{array} \right] & = & A^{\sim} + B^{\sim} + F_AS_1 + T_1E_A  + F_AS_2 + T_2E_B  + F_BS_3 + T_3E_A  + F_BS_4 + T_4E_B \\
& = & A^{\sim} + F_A(S_1 + S_2) + (T_1+ T_3)E_A  + B^{\sim} + F_B(S_3 + S_4) + (T_2 +T_4)E_B.
\end{eqnarray*} 
This expression is the same as (23.6). Thus (23.5) holds. This fact implies
that
$$\displaylines{
\hspace*{1.5cm} 
\max_{A^-,\, B^-}r[ \, M -  M( \, A^- + B^- \, )M \,] = \max_{N^-}r \left( M -  [ \, M, \  M \,]N^- \left[ \begin{array}{c} M  \\ M  \end{array} \right]  \right),\hfill
\cr
\hspace*{1.5cm} 
\min_{A^-,\, B^-}r[ \, M -  M( \, A^- + B^- \, )M \,] = \min_{N^-}r \left( M -  [ \, M, \  M \,]N^- \left[ \begin{array}{c} M  \\ M  \end{array} \right] \right). \hfill
\cr}
$$  
Applying (22.1) and (22.2) to the right-hand sides of the above two
equalities, we obtain
\begin{eqnarray*}  
\lefteqn{  \max_{N^-}r \left( M -  [ \, M, \  M \,]N^- \left[ \begin{array}{c} M  \\ M  \end{array} \right] \right) } \\
& = & \min \left\{ r(M), \ \ r \left[ \begin{array}{ccc} A  & 0 & M \\ 0 & B & M \\ M & M & M \end{array} \right] -r(N)  \right\} \\
& = & \min \left\{ r(M), \ \ r \left[ \begin{array}{ccc} A - M  & - M \\ -M & B - M  \end{array} \right]  + r(M) -r(A) - r(B)  \right\} \\
& = & \min \left\{ r(M), \ \ r \left[ \begin{array}{ccc} M  & A \\  B  & M  \end{array} \right]  + r(M) -r(A) - r(B)  \right\}, 
\end{eqnarray*} 
which is exactly (23.3), and 
$$
\displaylines{
\hspace*{1.5cm}
\min_{N^-}r \left( M -  [ \, M, \  M \,]N^- \left[ \begin{array}{c} M  \\ M
 \end{array} \right] \right) \hfill
\cr
\hspace*{1.5cm} 
=  r(N) - r \left[ \begin{array}{ccc} A  & 0 & M \\ 0 & B & M \end{array} \right] - r \left[ \begin{array}{cc} A  & 0 \\ 0 & B  \\  M & M \end{array} \right] + r \left[ \begin{array}{ccc} A  & 0 & M \\ 0 & B & M \\ M & M & M \end{array} \right] \hfill
\cr
\hspace*{1.5cm}
=  r(A) + r(B) + r(M) + r \left[ \begin{array}{ccc} M  & A \\  B  & M  \end{array} \right] - r \left[ \begin{array}{ccc} A  & 0 & B \\ 0 & B & A \end{array} \right] - r \left[ \begin{array}{cc} A  & 0 \\ 0 & B  \\  B & A \end{array} \right], \hfill
\cr}
$$
which is exactly (23.4). \qquad $\Box$  

\medskip

Two direct consequences can be derived from (23.3) and (23.4). 
\medskip

\noindent {\bf Theorem 23.3.}\, {\em Let $A, \, B \in {\cal F}^{m \times n}$  be given, and let $ M = A + B.$ Then there exist $ A^- \in \{ A^- \}$ and $B^- \in \{ B^- \}$ such that $  A^- + B^- \in \{ \,( \,A + B \,)^- \,\}$ holds if and only if
$$\displaylines{
\hspace*{2cm}
r \left[ \begin{array}{cc} M & A \\ B  & M \end{array} \right]  = r \left[ \begin{array}{cc} A  & 0 \\  0  & B \\ B & A   \end{array} \right]  + r \left[ \begin{array}{ccc} A  & 0 & B \\ 0 & B & A  \end{array} \right] - r(M) -r(A) -r(B). 
\hfill
\cr}
$$ }   
{\bf Theorem 23.4.}\, {\em  Let $A, \, B \in {\cal F}^{m \times n}$
 be given$,$ and let $ M = A + B \neq 0.$ The the following four statements
 are equivalent:

{\rm (a)}\, $ \{ \, A^- + B^- \,\} \subseteq \{ \,( \,A + B \,)^- \,\}$.

{\rm (b)}\, $ r\left[ \begin{array}{cc} M   & A \\ B  & M  \end{array}
\right] = r(A) + r(B) - r(\, A + B \, ).$

{\rm (c)}\, $ r \left[ \begin{array}{ccc} A  & 0 & M \\ 0 & B & M  \\ M & M & M   \end{array} \right] = r \left[ \begin{array}{ccc} A  & 0 & B \\ 0 & B & A  \\ B & A &  -2M   \end{array} \right] = r\left[ \begin{array}{cc} A & 0 \\ 0 & B  \end{array} \right].$ 

{\rm (d)}\, $ R(A) = R(B), \ R(A^T) = R(B^T)$ and $ A + B = - \frac{1}{2}( \, AB^-A + BA^-B \, )$.  } 

\medskip
   
\noindent {\bf Proof.}\, The equivalence of Parts (a) and (b) follows immediately from 
(22.3). The equivalence of Parts (b) and (c) follows from the rank equality 
$$\displaylines{
\hspace*{2cm}
 r \left[ \begin{array}{ccc} A  & 0 & M \\ 0 & B & M  \\ M & M & M   \end{array} \right] = 
r\left[ \begin{array}{cc} M   & A \\ B  & M  \end{array} \right] + r(M). \hfill
\cr}
$$ 
The equivalence of Parts (c) and (d) follows from (1.5). \qquad  $\Box$ 

\medskip

\noindent {\bf Theorem 23.5.}\, {\em Let $A, \, B \in {\cal F}^{m \times n}$
 be given$,$ and let $ M = A + B.$ Then
$$ \displaylines{
\hspace*{0cm}
\min_{A^-,\, B^-}r( \, M^- - A^-  - B^-  \,) = r( \, M - AM^-B \, ) - r \left[ \begin{array}{c} A \\ B \end{array} \right] - r[\, A, \ B \, ] + r(A) + r(B),  \hfill (23.7) 
\cr
\hspace*{0cm}
\min_{M^-, \, A^-,\, B^-}r(\, M^- - A^-  - B^- \,) \hfill
\hspace*{0cm}
\cr
 = r(M) + r(A) + r(B)  +
 r\left[ \begin{array}{cc} M  &  A  \\ B & M  \end{array} \right] -
 r\left[ \begin{array}{ccc} A & 0 & B \\ 0 & B & A  \end{array} \right] -
 r\left[ \begin{array}{cc} A  &  0  \\ 0 & B \\ B & A  \end{array} \right], 
\hfill (23.8)
\cr
\hspace*{0cm}
\max_{M^-}\min_{A^-,\, B^-}r( \, M^- - A^-  - B^-  \,) = \min \left\{  \, r(A) + r(B) - r \left[ \begin{array}{c} A \\ B \end{array} \right], \ \ \ \ \   r(A) +  r(B) - r[\, A, \ B \, ], \right. \hfill
\cr
\hspace*{5.5cm} 
\left.    r \left[ \begin{array}{cc} M  & A \\ B  & M  \end{array} \right] - r \left[ \begin{array}{c} A \\ B \end{array} \right]- r[\, A, \ B \, ] - r(M)  + r(A) +  r(B) \, \right\}.  \hfill (23.9) 
\cr}
$$ }
{\bf Proof.}\, According to (23.5) and (22.2) we first find 
\begin{eqnarray*}  
\lefteqn {\min_{A^-,\, B^-}r( \, M^- - A^-  - B^-  \,) } \\
& = & \min_{N^-} \left( \, M^- -  [ \, I_n, \ I_n \,]N^- \left[ \begin{array}{c} I_m \\ I_m  \end{array} \right]  \right)  \\
& = & r(N) - r \left[ \begin{array}{ccc} A  & 0 & I_m \\ 0 & B & I_m \end{array} \right] 
- r \left[ \begin{array}{cc} A  & 0 \\ 0 & B \\ I_n & I_n   \end{array} \right]  + r \left[ \begin{array}{ccc} A  & 0 & I_m \\ 0 & B & I_m  \\ I_n & I_n & M^-  \end{array} \right]  \\
&  = & r(N) - r \left[ \begin{array}{ccc} 0  & 0 & I_m \\ -A & B & 0 \end{array} \right] 
- r \left[ \begin{array}{cc} 0  & -A \\  0  & B \\ I_n & 0   \end{array} \right]  + r \left[ \begin{array}{ccc} 0  &  - M + AM^-B & 0 \\ 0 & 0 & I_m  \\ I_n & 0 &  0  \end{array} \right]  \\
& = & r(\,  M - AM^-B \, ) - r \left[ \begin{array}{c} A \\ B \end{array} \right]- r[\, A, \ B \, ] + r(A)+ r(B), 
\end{eqnarray*}
which is eaxactly (23.7). Next applying (21.4) and (21.5) to $  M - AM^-B$,
we obtain
$$
\displaylines{
\hspace*{1cm}
\min_{M^-} r(\,  M - AM^-B \, ) \hfill
\cr
\hspace*{1cm}
 =  r(M) + r[\, A, \ M \, ] +  r \left[ \begin{array}{c} B \\ M \end{array}
  \right] + r \left[ \begin{array}{cc} M  & B \\ A & M  \end{array}
  \right] - r \left[ \begin{array}{ccc} M  & 0 & B \\ 0 & A & M \end{array}
  \right] - r \left[ \begin{array}{cc} M  & 0 \\  0  & B \\ A & M
  \end{array} \right] \hfill
\cr
\hspace*{1cm}
=  r(M) + r[\, A, \ B \, ] +  r \left[ \begin{array}{c} A \\ B
\end{array} \right] + r \left[ \begin{array}{cc} M  & A \\ B & M
\end{array} \right] - r \left[ \begin{array}{ccc} A  & 0 & B \\ 0 & B & A
\end{array} \right] - r \left[ \begin{array}{cc} A  & 0 \\  0 & B \\ B & A
\end{array} \right],\hfill
\cr
\hspace*{1cm}
\max_{M^-} r(\,  M - AM^-B \, ) = \min \left\{ \,  r[\, A, \ B \, ], \ \ \
r \left[ \begin{array}{c} A \\ B \end{array} \right], \ \ \
r \left[ \begin{array}{cc} M  & B \\ A & M  \end{array} \right] - r(M) \,
 \right\}.  \hfill
 \cr}
$$
Putting them in (23.7) resprctively yields (23.8) and (23.9). \qquad
$\Box$

\medskip
  
Two direct consequences of Theorem 23.5 are given below.  

\medskip

\noindent {\bf Theorem 23.6.}\, {\em  Let $A, \, B \in {\cal F}^{m \times n}$  be given. Then for a given 
$( \, A + B \, )^-,$  there exist $ A^- \in \{ A^-\} $ and  $ B^- \in \{ B^-\} $ such that $ A^-  + B^-  
= ( \, A + B \, )^-$ if and only if   
$$ \displaylines{
\hspace*{2cm}
r[ \, A + B  - A( \, A + B\,)^-B \, ] =  r \left[ \begin{array}{c} A \\ B \end{array} \right] + r[\, A, \ B \, ] - r(A) - r(B).  \hfill
\cr}
$$  } 
{\bf Theorem 23.7.}\, {\em  Let $A, \ B \in {\cal F}^{m \times n}$ be given. Then $ \{ ( \, A + B \, )^- \} \subseteq \{ \, A^- + B^- \, \}$ holds if and only if 
$$ 
\displaylines{
\hspace*{2cm}
R(A) \cap R(B) = \{ 0\},  \ \ \  or \ \ \  R(A^T) \cap R(B^T) = \{ 0\},  \hfill
\cr
\hspace*{0cm}
or \hfill
\cr
\hspace*{2cm}
r \left[ \begin{array}{cc} A+B  & A \\ B  & A + B  \end{array} \right] = r \left[ \begin{array}{c} A \\ B \end{array} \right] + r[\, A, \ B \, ] +  r( \, A+B \, ) - r(A) - r(B). \hfill
\cr}
$$ }
\hspace*{0.4cm} Combining Theorems 23.4 and 23.7, one can easily establish a necessary and
sufficient condition for $ \{ ( \, A + B \, )^- \} = \{ \, A^- + B^- \, \}$
to hold. we shall, however, to present it in Section 23.3 as a special
case of a general result. \\

\noindent{ \Large 23.2. \ {\bf The relationships between $ A_1^- + A_2^-  + \cdots + A^-_k $ and
$(A_1 + A_2  + \cdots + A_k \,)^-$ }} \\

\noindent  The results in the preceding section can directly be extended to sums of
 $k$ matrices. We present them below without detailed proofs. 

\medskip

\noindent {\bf Theorem 23.8.}\, {\em  Let $ A_1, \,  A_2, \,  \cdots,\,  A_k
\in {\cal F}^{m \times n}$ be given$,$ and let $ M = A_1 + A_2  + \cdots +
A_k.$ Then
$$  
\displaylines{
\hspace*{0cm}
\max_{A_1^-,\,\cdots,\,A_k^-}r[ \, M -  M( \, A_1^-  + \cdots + A^-_k  \, )M \,] = \min \left\{ \, r(M), \ \ \ 
r( \, N - QMP \, ) + r(M) - r(N) \, \right\},
 \hfill (23.10) 
\cr
\hspace*{0cm}
and \hfill
\cr
\hspace*{0cm}
\min_{A_1^-,\,\cdots,\,A_k^-}r[ \, M -  M( \, A_1^-  + \cdots + A^-_k  \, )M \,] = r(M) + r(N) + r( \, N - QMP \, ) - r[\, N, \ QM \, ] - r \left[ \begin{array}{c} N \\ MP  \end{array} \right], 
 \hfill (23.11) 
\cr}
$$
where $ N = {\rm diag}( \, A_1,  \, A_2, \, \cdots, \, A_k \, ), \ P = [\, I_n,
 \, I_n, \, \cdots \, I_n \, ], \ Q = [\, I_m,  \, I_m, \, \cdots \, I_m \, ]^T.$ }

 \medskip

\noindent {\bf Proof.}\, It is easy to verify that  
$$ 
\displaylines{
\hspace*{2cm}
\{ \, A_1^-  + \cdots + A^-_k  \, \}  = \{ \, PN^-Q \, \}.  \hfill (23.12) 
\cr}
$$
In that case, it follows by (22.1) that 
$$  
\displaylines{
\hspace*{1cm}
\max_{A_1^-,\,\cdots,\,A_k^-}r[ \, M -  M( \, A_1^-  + \cdots + A^-_k  \, )M \,] \hfill 
\cr
\hspace*{1cm}
=  \max_{N^-}r( \, M -  MPN^-QM \, ) \hfill
\cr
\hspace*{1cm}
= \min \left\{ \, r(M), \ \ \  r \left[ \begin{array}{cc} N  & QM \\  MP & M \end{array} \right] -r(N) \  \right\} \hfill
\cr
\hspace*{1cm}
= \min \left\{  \, r(M), \ \ \ r( \, N - QMP \, ) + r(M) - r(N) \, \right\}, \hfill
\cr}
$$
which is exactly (3.1). Applying (22.2), we also obtain
$$  \displaylines{
\hspace*{1cm} 
\min_{A_1^-,\,\cdots,\,A_k^-}r[ \, M -  M( \, A_1^-  + \cdots + A^-_k  \,
 )M \,] \hfill
\cr
\hspace*{1cm}
= \min_{N^-}r( \, M -  MPN^-QM \, ) \hfill
\cr
\hspace*{1cm}
= r (N) - r[ \, N, \ QM \, ] - r \left[ \begin{array}{c} N \\ MP \end{array} \right] + r\left[ \begin{array}{cc}  N  &  QM  \\ MP & M  \end{array} \right] \hfill
\cr
\hspace*{1cm}
= r (N) - r[ \, N, \ QM \, ] - r \left[ \begin{array}{c} N \\ MP \end{array} \right] + r( \, N - QMP \, ) + r(M), \hfill
\cr}
$$  
which is exactly (23.11). \qquad  $ \Box$  

\medskip

Two direct consequences of (23.10) and (23.11) are listed below. 

\medskip

\noindent {\bf Theorem 23.9.}\, {\em Let $ A_1, \,  A_2, \,  \cdots,\,  A_k
\in {\cal F}^{m \times n}$ be given$,$ and denote $ M = A_1 + A_2  + \cdots
+ A_k.$ Then there exist $ A^-_i \in \{ A^-_i \},$ $ i = 1,\, 2, \, \cdots, \,
k $  such that $A_1^-  + A_2^- + \cdots + A^-_k  \in \{ M^-\}$ if and only
if
$$\displaylines{
\hspace*{2cm}
r( \, N - QMP \, ) =  r \left[ \begin{array}{c} N \\ MP  \end{array} \right]
+ r[\, N, \ QM \, ] - r(M) - r(N), \hfill
\cr}
$$
where $ N, \ P $ and $ Q $ are defined in Theorem 23.8. } 

\medskip

\noindent {\bf Theorem 23.10.}\, {\em Let $ A_1, \,  A_2, \,  \cdots,\,
A_k \in {\cal F}^{m \times n}$ be given$,$ and let $ M = A_1 + A_2  + \cdots
+ A_k.$ The the following four statements are equivalent$:$

{\rm (a)} \ $ \{ \, A_1^-  + A_2^- \cdots + A^-_k \,\} \subseteq \{ \,(
\,A_1 + A_2  + \cdots + A_k \,)^- \,\}.$

{\rm (b)} \ $ r( \, N - QMP \, ) =  r(N) - r(QMP).$ 

{\rm (c)} \ $ r \left[ \begin{array}{cc} N & QM \\ MP & M \end{array}
\right] = r(N).$

{\rm (d)} \ $ R(M) = R(A_i), \ R(M^T) = R(A_i^T),$ $ i = 1, \, 2, \, \cdots,
\ k,$ and $M = MPN^-QM,$ where $ N, \ P$ and $ Q $ are defined in
Theorem 23.8.  }

\medskip

\noindent {\bf Proof.}\, The equivalence of Parts (a) and (b) follows
immediately from (23.10). The equivalence of Parts (b) and (c) is
evident. The equivalence of Parts (c) and (d) follows from (1.5).
\qquad $\Box$

\medskip

\noindent {\bf Theorem 23.11.}\, {\em  Let $ A_1, \,  A_2, \,  \cdots,\,
A_k \in {\cal F}^{m \times n}$ be given$,$ and denote $ M = A_1 + A_2  +
\cdots + A_k.$  Then
$$  \displaylines{
\hspace*{0cm}
\min_{A_1^-,\,\cdots,\,A_k^-}r( \, M^- - A_1^-  -\cdots - A^-_k \,)
= r(N) - r[\, N, \ Q \, ] - r \left[ \begin{array}{c} N \\ P \end{array}
\right] + r\left[ \begin{array}{cc} N  &  Q  \\ P & M^-  \end{array}
\right], \hfill (23.13)
\cr
\hspace*{0cm} 
\min_{M^-,\, A_1^-,\,\cdots,\,A_k^-}r( \, M^- - A_1^-  -\cdots - A^-_k \,)
= r(M) + r(N) + r(\, N - QMP \,) - r[\, N, \ QM \, ] -  r \left[
\begin{array}{c} N \\ MP \end{array} \right], \hfill (23.14)
\cr
\hspace*{0cm}
\max_{M^-}\min_{A_1^-,\,\cdots,\,A_k^-}r( \, M^- - A_1^-  -\cdots - A^-_k
\,) = \min \left\{  \,  r(N)  + m - r[\, N, \ Q \,], \ \ \ r(N)  + n -
r \left[ \begin{array}{c} N \\ P \end{array} \right],  \right. \hfill
\cr
\hspace*{4cm}
\left.  m + n + r( \, N- QMP \, ) + r(N) - r(M) - r[\, N, \ Q \, ] -
r \left[ \begin{array}{c} N \\ P \end{array} \right] \ \right\},
\hfill (23.15)
\cr}
$$ 
where $ N, \ P$ and $ Q $ are defined in Theorem 23.8.} 

\medskip

\noindent {\bf Proof.}\, According to (23.12) and (22.2) we first find
that
$$\displaylines{
\hspace*{1cm}  
\min_{A^-_1,\,\cdots,\,A_k^-} r( \, M^- - A_1^- - \cdots - A_k^- \,) \hfill
\cr
\hspace*{1cm}
 =  \min_{N^-} ( \, M^- - PN^-Q \,) = r(N) - r[\, N, \ Q \, ] - r \left[ \begin{array}{c} N \\ P \end{array} \right] + r\left[ \begin{array}{cc} N  &  Q  \\ P & M^-  \end{array} \right],  \hfill
\cr}
$$
which is eaxactly (23.13). Applying (21.5) to the $ 2 \times 2 $ block
matrix in the above equality, we further obtain
\begin{eqnarray*}  
\lefteqn{\min_{M^-}r\left[ \begin{array}{cc} N  &  Q  \\ P & M^-  \end{array} \right] } \\
& = & \min_{M^-}r\left( \, \left[ \begin{array}{cc} N  &  Q  \\ P & 0  \end{array} \right] +
 \left[ \begin{array}{c} 0 \\ I_n \end{array} \right]M^-[\, 0, \ I_m \,]\,
 \right)  \\
& = & r(M) + r\left[ \begin{array}{ccc} N  &  Q  & 0  \\ P & 0  & I_n
\end{array} \right] + r\left[ \begin{array}{cc} N  &  Q  \\ P & 0 \\ 0 &
I_m \end{array} \right] +
r\left[ \begin{array}{ccc} -M  &  0  & I_m  \\ 0 & N & Q \\I_n & P & 0 \end{array} \right] \\
& &  - \ \ r\left[ \begin{array}{cccc} -M  &  0  & 0 & I_m  \\ 0 & 0 & N & Q \\ 0 &  I_n & P & 0 \end{array} \right] - r\left[ \begin{array}{ccc} -M & 0 & 0 \\  0 & 0 & I_m  \\ 0 & N & Q \\I_n & P & 0 \end{array} \right] \\
& = & r(M) +  r[\, N, \ Q \, ] +  r \left[ \begin{array}{c} N \\ P \end{array} \right] +
r\left[ \begin{array}{ccc} 0  &  0  & I_m  \\ 0 & N - QMP  & 0  \\  I_n & 0 & 0 \end{array} \right]
 - r\left[ \begin{array}{ccc} M  & 0 & I_m  \\ 0 & N  & Q \end{array} \right] - 
r\left[ \begin{array}{ccc} M  & 0  \\ 0  & N \\ I_n & P \end{array} \right] \\
& = & r(M) +  r[\, N, \ Q \, ] + r \left[ \begin{array}{c} N \\ P \end{array} \right] + r(\,  N - QMP\, ) 
- r[\, N, \ QM \,] - r \left[ \begin{array}{c} N \\ MP \end{array} \right].
\end{eqnarray*}  
Putting it in (23.13) yields (23.14). Next applying (21.5), we obtain the
following
\begin{eqnarray*}  
\max_{M^-}r\left[ \begin{array}{cc} N & Q \\ P & M^- \end{array} \right]
 & = & \max_{M^-}r\left( \, \left[ \begin{array}{cc} N  &  Q  \\ P & 0  \end{array} \right] +
 \left[ \begin{array}{c} 0 \\ I_n \end{array} \right]M^-[\, 0, \ I_m \,]\, \right)  \\
& = & \min \left\{ r\left[ \begin{array}{ccc} N  &  Q  & 0  \\ P & 0  & I_n \end{array} \right], \ \ \  
r\left[ \begin{array}{cc} N  &  Q  \\ P & 0 \\ 0 & I_n \end{array} \right], \ \ \  r\left[ \begin{array}{ccc} -M  &  0  & I_m  \\ 0 & N & Q \\I_n & P & 0 \end{array} \right] - r(M) \right\} \\
& = & \min \left\{  n + r[\, N, \ Q \, ], \ \ \ m +  r\left[ \begin{array}{c} N \\ P \end{array} \right], \ \ \  m + n + r(\, N - QMP \, ) - r(M)  \right\}. 
\end{eqnarray*}  
Putting it in (23.13) yields (23.15). \qquad  $\Box$

\medskip

Two direct consequences of Theorem 23.11 are given below. 

\medskip

\noindent {\bf Theorem 23.12.}\, {\em Let $ A_1, \,  A_2, \,  \cdots,\,  A_k
\in {\cal F}^{m \times n}$ be given$,$ and let $ M = A_1 + A_2  + \cdots +
A_k.$  Then  for a given $M^-,$  there exist $ A^-_i \in \{ A^-_i \},$ $
i = 1, \, 2, \, \cdots, \, k $  such that $ A_1^-  + A_2^-  + \cdots + A_k^-
= M^-$ if and only if $ M^- $  satisfies
$$ 
\displaylines{
\hspace*{2cm}
r\left[ \begin{array}{cc} N  &  Q  \\ P & M^-  \end{array} \right] = r\left[ \begin{array}{c} N \\ P \end{array} \right] +  r[\, N, \ Q \, ] - r(N),  \hfill
\cr}
$$ 
where $ N, \ P$ and $ Q $ are defined in Theorem 23.8. } 

\medskip

\noindent {\bf Theorem 23.13.}\, {\em Let $ A_1, \,  A_2, \,  \cdots,\,  A_k
\in {\cal F}^{m \times n}$ be given$,$ and let $ M = A_1 + A_2  + \cdots +
A_k.$  Then the set inclusion
$$  \displaylines{
\hspace*{2cm}
\{ \, (\, A_1 + A_2  + \cdots + A_k \, )^- \, \} \subseteq  \{ \, A_1^-  + A_2^-  + \cdots + A_k^- \, \}  \hfill(23.16) 
\cr}
$$ 
holds  if and only if 
$$ \displaylines{
\hspace*{2cm}
R(A_1) \cap R(A_2) \cap \cdots \cap R(A_k)  = \{ 0 \}, \hfill (23.17)   
\cr
\hspace*{0cm}
or \hfill
\cr
\hspace*{2cm}
R(A_1^T) \cap R(A_2^T) \cap \cdots \cap R(A_k^T)  = \{ 0 \}, \hfill (23.18)   
\cr
\hspace*{0cm}
or \hfill
\cr
\hspace*{2cm}
r( \, N- QMP \, ) = r \left[ \begin{array}{c} N \\ P \end{array} \right] + r[\, N, \ Q \, ] - r(N) + r(M) - m - n, \hfill (23.19)  
\cr}
$$
where $ N, \ P$ and $ Q $ are defined in Theorem 23.8. } 

\medskip

\noindent {\bf Proof.}\, It is easy to see that the set inclusion in (23.16) hold if and only if 
$$\displaylines{
\hspace*{2cm}
\max_{M^-}\min_{A_1^-,\,\cdots,\,A_k^-}r( \, M^- - A_1^-  -\cdots - A^-_k \,) = 0.  \hfill
\cr}
$$
In light of (23.15), the above equality is equivalent to 
$$\displaylines{
\hspace*{2cm}
r[\, N, \ Q \,] = r(N) + r(Q),  \ \ \ {\rm or} \ \ \  r \left[ \begin{array}{c} N \\ P \end{array} \right] 
= r(N) + r(P) , \hfill (23.20) 
\cr}
$$
or (23.18) holds. The two rank equalities in (23.20) are
equivalent to (23.17) and (23.18) according to (23.1).   \qquad $\Box$ \\

\noindent{ \Large 23.3. \ {\bf The relationships between $ \{ \, A^- + B^- \, \} $ and 
 $ \{ \, C^- \, \}$, and parallel sum of two matrices}}  \\
 
\noindent It is well known (see \cite{MO}, \cite{RM}) that the parallel sum of two
matrices $ A $ and $ B $ of the same size is defined to be $
C:= A(\, A + B \, )^-B,$ whenever this product is invariant with respect
to the choice of $(\, A + B \,)^-$. One of the  well-known nice
properties on parallel sum of two matrices is $\{ \, A^- + B^- \,\}
= \{ C^- \}$. This set equality motivates us to consider the
relationship between the two sets $\{ \, A^- + B^- \, \}$ and
$ \{ C^- \}$ in general cases, where $ A, \ B$, and $ C$ are  any
three given matrices of the same size.
Just as what we do in Section 23.1, we first establish several
basic rank equalities related to generalized inverses of  $ A, \ B$,
and $ C$, and then deduce from them various relationships between
$\{ \, A^- + B^- \,\}$ and $ \{ \, C^- \}$.

\medskip
       
\noindent {\bf Theorem 23.14.}\, {\em  Let $A, \, B, \, C \in {\cal F}^{m \times
n}$ be given. Then
$$
\displaylines{
\hspace*{0cm}
\max_{A^-,\, B^-}r[ \, C -  C( \, A^- + B^- \, )C \,] = \min \left\{ r(C), \ \ \ 
r \left[ \begin{array}{cc}  A - C  & C \\ C & B - C \end{array} \right]  + r(C) - r(A) - r(B)  \right\},
 \hfill (23.21) 
\cr}
$$
$$
\min_{A^-,\, B^-}r[ \, C -  C( \, A^- + B^- \, )C \,] = r(A) + r(B) + r(C) + \  r \left[ \begin{array}{cc} A - C  & C \\ C & B - C  \end{array} \right]  
- r \left[ \begin{array}{ccc} A  & 0 & C \\ 0 & B & C    \end{array} \right] - 
 r \left[ \begin{array}{cc} A  & 0 \\ 0 & B \\ C & C \end{array} \right].
 \eqno (23.22)
$$} 
{\bf Proof.}\, By (23.5) and (22.1), we easily find
\begin{eqnarray*}  
\max_{A^-,\, B^-}r[ \, C -  C( \, A^- + B^- \, )C \,]& = & \max_{N^-}r \left( C -  [ \, C, \  C \,]N^- \left[ \begin{array}{c} C  \\ C  \end{array} \right]  \right) \\
& = & \min \left\{  r(C), \ \ r \left[ \begin{array}{ccc} A  & 0 & C \\ 0 & B & C \\ C & C & C \end{array} \right] -r(N)   \right\} \\
& = & \min \left\{  r(C), \ \ r \left[ \begin{array}{cc} A - C  & - C \\ -C & B - C \end{array} \right]  + r(C) -r(A) - r(B)   \right\}, 
\end{eqnarray*} 
which is (23.21). Next applying (23.5) and (22.2), we obtain
$$
\displaylines{
\hspace*{1cm}
\lefteqn{ \min_{A^-,\, B^-}r[ \, C -  C( \, A^- + B^- \, )C \,] } \hfill
\cr
\hspace*{1cm}
=  \min_{N^-} r \left(  \, C -  [ \, C, \  C \,]N^- \left[ \begin{array}{c} C \\ C \end{array} \right] \, \right)  
\hfill
\cr
\hspace*{1cm}
 =  r(N) - r \left[ \begin{array}{ccc} A  & 0 & C \\ 0 & B & C \end{array} \right] - r \left[ \begin{array}{cc} A  & 0 \\0 & B  \\  C & C \end{array} \right] + r \left[ \begin{array}{ccc} A  & 0 & C \\ 0 & B & C \\ C & C & C \end{array} \right] \hfill
\cr
\hspace*{1cm}
 =  r(A) + r(B) + r(C) + r \left[ \begin{array}{ccc} A - C & C \\  C  & B - C  \end{array} \right] - r \left[ \begin{array}{ccc} A  & 0 & C \\ 0 & B & C \end{array} \right] - r \left[ \begin{array}{cc} A  & 0 \\ 0 & B  \\  C & C \end{array} \right], \hfill
\cr}
$$ 
which is exactly (23.22). \qquad $\Box$  

\medskip

Two consequences can directly be derived from (23.21) and (23.22). 

\medskip

\noindent {\bf Theorem 23.15.}\, {\em Let $A, \, B, \, C \in {\cal F}^{m \times n}$  be given. Then there exist $ A^- \in \{ A^- \}$ and $B^- \in \{ B^- \}$ such that $  A^- + B^- \in \{ \, C^- \,\}$ holds if and only if
$$
 r \left[ \begin{array}{cc}A - C  & C \\ C & B - C \end{array} \right]
  = r \left[ \begin{array}{cc}  C & B \\ A  & A + B  \end{array} \right]
   = r \left[ \begin{array}{cc} A  & 0 \\  0  & B \\ B & A
   \end{array} \right]  + r \left[ \begin{array}{ccc} A  & 0 & B \\ 0 & B
   & A  \end{array} \right] - r(N) -r(A) -r(B). \eqno (23.23)
$$ }   
{\bf Theorem 23.16.}\, {\em  Let $A, \, B, \, C  \in {\cal F}^{m \times n}$  be given with $ C \neq 0.$ Then 
the following four statements
are equivalent$:$

{\rm (a)} \ $ \{ \, A^- + B^- \,\} \subseteq \{ \, C^-\}$.

{\rm (b)} \ $  r \left[ \begin{array}{cc} A - C  & C \\ C & B - C \end{array} \right] = r \left[ \begin{array}{cc}  C & B \\ A  & A + B  \end{array} \right]  = r(A) + r(B) - r(C).$ 

{\rm (c)} \ $ r \left[ \begin{array}{ccc} A  & 0 & C \\ 0 & B & C  \\ C & C & C \end{array} \right] = r\left[ \begin{array}{cc} A & 0 \\ 0 & B  \end{array} \right].$ 

{\rm (d)} \ $ R(C) \subseteq R(A), \ R(C) \subseteq R(B), \
R(C^T) \subseteq R(A^T), \ R(C^T) \subseteq r(B^T)$ and
$ C = CA^-C + CB^-C$.}

\medskip

\noindent {\bf Proof.}\, The equivalence of Parts (a) and (b) follows from
(23.21). The equivalence of Parts (b) and (c) is evident. The equivalence
of Parts (c) and (d) follows from (1.5). \qquad  $\Box$

\medskip

\noindent {\bf Theorem 23.17.}\, {\em Let $A, \, B, \, C \in {\cal F}^{m \times
n}$ be given. Then
$$ \displaylines{
\hspace*{0cm}
\min_{A^-,\, B^-}r( \, C^- - A^-  - B^-  \,) = r( \, A + B  - AC^-B \, )
- r \left[ \begin{array}{c} A \\ B \end{array} \right] - r[\, A, \ B \, ]
+ r(A) + r(B), \hfill (23.24)
\cr}
$$
$$ 
\min_{C^-, \, A^-,\, B^-}r(\, C^- - A^-  - B^- \,) = r(C) + r(A) + r(B)  +
 r\left[ \begin{array}{cc} C  &  B  \\ A & A +B  \end{array} \right] -
 r\left[ \begin{array}{ccc} A & 0 & C \\ 0 & B & C  \end{array} \right] -
 r\left[ \begin{array}{cc} A  &  0  \\ 0 & B \\ C & C \end{array} \right], 
 \eqno (23.25) 
$$
$$ \displaylines{
\hspace*{0cm}
\max_{C^-}\min_{A^-,\, B^-}r( \, C^- - A^-  - B^-  \,) = \min \left\{
\, r(A) + r(B) - r \left[ \begin{array}{c} A \\ B \end{array} \right],
\ \ \ r(A) +  r(B) - r[\, A, \ B \, ], \right. \hfill
\cr
\hspace*{4.5cm}
\left.  r \left[ \begin{array}{cc} C  & B \\ A  & A + B  \end{array}
\right] - r \left[ \begin{array}{c} A \\ B \end{array} \right]- r[\, A, \
B \, ] - r(C)  + r(A) +  r(B) \, \right\}. \hfill (23.26)
\cr}
$$ }
\hspace*{0.4cm} The proof of this theorem is much similar to that of Theorem 23.5 and is,
therefore, omitted. Two direct consequences of Theorem 23.17 are given
below.

\medskip

\noindent {\bf Theorem 23.18.}\, {\em  Let $A, \, B, \, C \in {\cal F}^{m \times
n}$ be given. Then for a given
$C^-,$  there exist $ A^- \in \{ A^-\} $ and  $ B^- \in \{ B^-\} $ such that
 $ A^-  + B^- = C^-$ if and only if $ A, \ B $ and $C^- $ satisfies    
$$ 
\displaylines{
\hspace*{2cm}
r[ \, A + B  - AC^-B \, ] =  r \left[ \begin{array}{c} A \\ B \end{array}
\right] + r[\, A, \ B \, ] - r(A) - r(B). \hfill (23.27)
\cr}
$$  } 
\noindent {\bf Theorem 23.19.}\, {\em  Let $A, \ B, \ C \in {\cal F}^{m
\times n}$ be given. Then $ \{ C^- \} \subseteq \{ \, A^- + B^- \, \}$
holds if and only if
$$ 
\displaylines{
\hspace*{2cm}
R(A) \cap R(B) = \{ 0\},  \ \ \  or \ \ \  R(A^T) \cap R(B^T) = \{ 0\},
\hfill (23.28)
\cr
\hspace*{0cm}
or \hfill
\cr
\hspace*{2cm}
r \left[ \begin{array}{cc} C & B \\ A  & A + B  \end{array} \right] =
r \left[ \begin{array}{c} A \\ B \end{array} \right] + r[\, A, \ B \, ]
+ r(C) - r(A) - r(B).  \hfill (23.29)
\cr}
$$}
\hspace*{0.3cm} The two conditions in (23.28) have no relation with the matrix $ C $,
which implies that under (23.28),
$ \{ C^- \} \subseteq \{ \, A^- + B^- \, \}$ holds for any choice of $ C $.
Now setting $ C = 0 $, then $ \{ \, 0^- \} = {\cal F}^{ n \times m}$. Thus
under (23.28), there is
$$ \displaylines{
\hspace*{2cm}
\{ \, A^- + B^- \, \} =  \{ \, 0^- \} = {\cal F}^{ n \times m}.
\hfill (23.30)
\cr}
$$ 
Conversely, if (23.30) holds, then it is easy to see  from Theorem 23.19
that (23.28) holds. Thus (23.28) is a necessary and sufficient condition
for (23.30) to hold.

\medskip

Now combining Theorems 23.16 and 23.19, we obtain the following result. 

\medskip

 \noindent {\bf Theorem 23.20.}\, {\em  Let $A, \, B, \,  C \in {\cal F}^{m
 \times n}$  be given and suppose that
$$ \displaylines{
\hspace*{2cm}
R(A) \cap R(B) \neq \{ 0\},  \ \ \  and  \ \ \  R(A^T) \cap R(B^T) \neq
\{ 0\}. \hfill (23.31)
\cr}
$$
Then the equality 
$$ \displaylines{
\hspace*{2cm}
\{ \ A^- + B^- \, \} =  \{ C^- \} \hfill (23.32) 
\cr}
$$ 
holds if and only if 
$$\displaylines{
\hspace*{2cm} 
R(B) \subseteq R( \, A + B  \,), \ \   R(A^T) \subseteq R( \,A^T + B^T \, )
\ \ and  \ \ C = A( \, A + B  \,)^- B, \hfill (23.33)
\cr}
$$ 
that is, $ A $ and $B$ are parallel summable and $C$ is the parallel sum
of $ A $and $ B $. In that case, the rank of $ C $ satisfies the following
rank equality
$$ 
\displaylines{
\hspace*{2cm}
r(C) = r(A ) + r(B) - r( \, A + B  \,). \hfill (23.34) 
\cr}
$$}
 {\bf Proof.}\, Assume first that (23.32) holds. By Theorem 23.19,
 we know that  $A, \ B$ and $ C $
satisfy (23.29). On the other hand, (23.32) implies that 
$$ \displaylines{
\hspace*{2cm}
\min_{A^-,\, B^-}r( \, A^-  + B^-  \,)  = \min_{C^-}r( C^-). \hfill (23.35)    
\cr}
$$ 
It is well-known that the minimal rank of $ C^-$ is  $r(C)$. On the other
hand, it follows from (23.24) that
$$ \displaylines{
\hspace*{2cm} 
\min_{A^-,\, B^-}r( \, A^-  + B^-  \,) =  r( \, A + B \,)+ r(A) + r(B) -
r[\, A, \ B \, ]  - r \left[ \begin{array}{c} A \\ B \end{array}
\right]. \hfill
\cr}
$$
Thus (23.35) is equivalent to 
$$\displaylines{
\hspace*{2cm}
r(C) = r( \, A + B \,)+ r(A) + r(B) - r[\, A, \ B \, ]  -
r \left[ \begin{array}{c} A \\ B \end{array} \right]. \hfill (23.36)
\cr}
$$
Putting it (23.29) yields 
$$\displaylines{
\hspace*{2cm} 
r \left[ \begin{array}{cc} C & B \\ A  & A + B  \end{array} \right]
= r( \, A + B \,), \hfill (23.37)
\cr}
$$
which, by (1.5), is equivalent to (23.33), meanwhile
\begin{eqnarray*} 
r(C) &= & \min_{( A + B)^-}r[ \, A( \, A + B \,)^-B \,]  \\
& = & r( \, A + B \,) -  r[\, A + B, \ B \, ]  - r \left[ \begin{array}{c}
A + B \\ B \end{array} \right] +  r \left[ \begin{array}{cc } A + B  &  B
\\ A & 0 \end{array} \right] \\
& = & r(A ) + r(B) - r( \, A + B  \,). 
\end{eqnarray*}  
Conversely, if (23.33) holds, then (23.34) and (23.37) also hold.
Combining both of them shows that the two rank equalities in Theorem
23.16(b) and (23.29) are satisfied.  Therefore (23.32) holds.  \qquad
$ \Box$  

\medskip

The equivalence of (23.32) and (23.33) was previously proved by Mitra and
Odell in \cite{MO}. But the assumption (23.31) is neglected there.
As shown in (23.30), if $ R(A) \cap R(B) = \{ 0 \}$ and $ R(A^T) \cap
 R(B^T) \neq \{ 0 \}$, then the equality $ \{ \, A^- + B^- \, \} =
 \{ \, 0^- \}$. This, however, does not imply that
 $A$ and $ B $ are parallel summable and $ A(\, A + B \, )^-B = 0 $,
 in general. An example is
$$ \displaylines{
\hspace*{2cm}
   A = \left[ \begin{array}{c} 1 \\ 0 \end{array} \right],  \qquad       
 B = \left[ \begin{array}{c} 0 \\ 1 \end{array} \right],  \hfill     
\cr}
$$
both of which satisfy $ \{ \, A^- + B^- \, \} =  \{ \, 0^- \} $, since
 $ R(A) \cap R(B) = \{ 0 \}$. If we let  
$$\displaylines{
\hspace*{2cm}
( \, A + B \, )^- = \left[ \begin{array}{c} 1 \\ 1 \end{array} \right]^- = \frac{1}{2}[\, 1, \ 1 \,], \hfill
\cr\hspace*{0cm}
then \hfill
\cr
\hspace*{2cm}
A( \, A + B \, )^-B  = \frac{1}{2} \left[ \begin{array}{c} 1 \\ 0 \end{array} \right][\, 1, \ 1 \,]\left[ \begin{array}{c} 0 \\ 1 \end{array} \right] 
= \frac{1}{2} \left[ \begin{array}{c} 1 \\ 0 \end{array} \right] \neq 0. \hfill
\cr
\hspace*{0cm}
If \ we  \ let  \hfill
\cr
\hspace*{2cm}
( \, A + B \, )^- = \left[ \begin{array}{c} 1 \\ 1 \end{array} \right]^- =[\, 1, \ 0 \,], \hfill
\cr
\hspace*{0cm}
then \hfill
\cr
\hspace*{2cm}
A( \, A + B \, )^-B  = \left[ \begin{array}{c} 1 \\ 0 \end{array} \right][\, 1, \ 0 \,]\left[ \begin{array}{c} 0 \\ 1 \end{array} \right] 
= \left[ \begin{array}{c} 0 \\ 0 \end{array} \right]. \hfill
\cr}
$$
Therefore $ A$ and $B$ are not parallel summable. 

Some interesting consequences can be derived from the above results. For example, let $ B = I_m  - A $ and  $ C= I_m$ in 
(23.24), we then get 
$$
\min_{A^-,\, (I_m  - A)^-}r[ \, I_m - A^-  - (I_m  - A)^-  \,] = r[\, (A -A^2) -  (A -A^2)^2 \, ].
$$  
Thus there are $ A^- $ and  $ (I_m - A)^-$  such that $ A^- + (I_m - A)^- = I_m$ if and only if $A - A^2 $ is 
idempotent.

Replace  $ B = I_m  - A $ and  $ C= A - A^2$ in Theorem 23.20. Then it is easy to verify that 
 these $ A,\ B$ and $ C$ satisfy the condition (23.33). Thus the set equality 
$$ \displaylines{
\hspace*{2cm}
 \{  \,( A - A^2)^- \, \} = \{ \, A^- + ( I_m - A)^- \, \} \hfill  
\cr}
$$ 
holds for any $ A$. In other words, the matrices $ A $ and $ I_m - A $ are always parallel summable, and  $  A - A^2$ 
is their parallel sum. Recall (22.23), we then get the following  
$$
\{  \, A^-( I_m - A)^- \, \} \subseteq \{ \, A^- + ( I_m - A)^- \, \}  \ \ \  {\rm and } \ \ \
 \{  \, ( I_m - A)^-A^- \, \} \subseteq \{ \, A^- + ( I_m - A)^- \, \}. 
$$
When one of  $ A $ and $ I_m - A$ is nonsingular, say, $A$,  there is  
$$
\{  \, A^{-1}( I_m - A)^- \, \} \subseteq \{ \, A^{-1} + ( I_m - A)^- \, \}  \ \ \  {\rm and } \ \ \
 \{  \, ( I_m - A)^-A^{-1} \, \} \subseteq \{ \, A^{-1} + ( I_m - A)^- \, \}. 
$$
When both $ A $ and $ I_m - A$  nonsingular, the above becomes a trivial result 
$A^{-1}(I_m - A)^{-1} =  A^{-1} + ( I_m - A)^{-1}$.  

Replace $  A, \ B$ and $ C$ in Theorem 23.20 by $ I_m  +  A, \ I_m  - A $ and  $(I_m - A^2)/2$, respectively. Then it is 
easy to verify that they  satisfy the condition (23.33). Thus the set equality 
$$
 \{  \,( I_m - A^2)^- \, \} = \left\{ \, \frac{1}{2}( I_m + A)^-  +  \frac{1}{2}( I_m - A)^- \, \right\}
$$ 
holds for any $ A$. In other words, the matrices $ I_m +  A $ and $ I_m - A $ are always parallel summable, 
and the matrix $ (I_m - A^2)/2$ is their parallel sum. Recall (22.24), then we also get the following  
$$
\{  \, ( I_m + A)^-( I_m - A)^- \, \} \subseteq \left\{ \, \frac{1}{2}( I_m + A)^-  + 
 \frac{1}{2}( I_m - A)^- \,  \right\}, 
$$
$$
\{  \, ( I_m - A)^-( I_m + A)^- \, \} \subseteq \left\{ \, \frac{1}{2}( I_m + A)^-  + 
 \frac{1}{2}( I_m - A)^- \,  \right\}. 
$$
When one of  $I_m + A$ and $ I_m - A$ is nonsingular, say, $ I_m + A$,  there is  
$$
\{  \, (I_m + A)^{-1}( I_m - A)^- \, \} \subseteq  \left\{ \,  \frac{1}{2} ( I_m + A)^{-1} +  \frac{1}{2}( I_m - A)^- \,  \right\},
$$
$$
 \{  \, ( I_m - A)^- ( I_m + A)^{-1} \, \} \subseteq  \left\{ \,  \frac{1}{2} ( I_m + A)^{-1}
 + \frac{1}{2} ( I_m - A)^- \,  \right\}. 
$$
When both $I_m + A $ and $ I_m - A$ are  nonsingular, the above becomes a trivial result 
$2(I_m + A)^{-1}(I_m - A)^{-1} =  (I_m + A)^{-1} + ( I_m - A)^{-1}$.  

In general, suppose that $ \lambda_1 \neq  \lambda_2$ are two scalars,  and replacing 
$A, \ B$ and $ C$ in Theorem 23.20 by $\lambda_1 I_m  - A, \, \lambda_2 I_m  - A $ and
  $(\lambda_1  I_m  - A )(\lambda_2  I_m  - A ) / (\lambda_1 - \lambda_2)$, respectively. Then it is 
easy to verify that they  satisfy the condition (23.33). Thus the set equality 
$$ 
 \{  \, [\, (\lambda_1  I_m  - A )(\lambda_2  I_m  - A ) \,]^- \, \} = 
\left\{ \, \frac{1}{\lambda_1 - \lambda_2}(\lambda_1I_m - A)^-  +  \frac{1}{\lambda_2 -
 \lambda_1}( \lambda_2I_m - A)^- \, \right\}
$$ 
holds for any $ A$. In other words, the matrices $\lambda_1  I_m  - A $ and $\lambda_2 I_m  - A $ 
are parallel summable, and the matrix 
$(\lambda_1  I_m  - A )(\lambda_2  I_m  - A )/(\lambda_1 - \lambda_2) $ is their parallel sum.  Recall (22.22), 
we then get the following  
$$
 \{  \,(\lambda_1  I_m  - A )^-(\lambda_2  I_m  - A )^- \, \} \subseteq 
\left\{ \, \frac{1}{\lambda_1 - \lambda_2}(\lambda_1 - I_m A)^-  +  \frac{1}{\lambda_2 -
 \lambda_1}( \lambda_2I_m - A)^- \, \right\}.
$$ 
This result motivates us to guess that for $ \lambda_1, \, \cdots, \, \lambda_k$ with  $ \lambda_i \neq  \lambda_j$
 for all $ i \neq j$,there is 
$$
 \{ (\lambda_1  I_m  - A )^-  \cdots  (\lambda_k I_m  - A )^-  \} \subseteq \left\{  \frac{1}{p_1}(\lambda_1 -I_m A)^-   + \cdots  + \frac{1}{p_k}(\lambda_k I_m - A)^-  \right\},
$$ 
where
$$ 
p_i = (\lambda_1 - \lambda_i)\cdots(\lambda_{i-1} - \lambda_i)
(\lambda_{i+1} - \lambda_i)\cdots(\lambda_k - \lambda_i),  \ \ \ i = 1, \, 2, \, \cdots, \, k.
$$
We leave it as an open problem to the reader. \\

\noindent{ \Large 23.4. \ {\bf The relationships between  $ \{ \, A_1^- + A_2^-  + \cdots +
A^-_k \, \} $ and  $ \{ C^- \}$, and parallel sum of $k$ matrices}} \\
 
\noindent  The results in  Section 23.3 can easily be generalized to sums of $k$
matrices, which can help to extend the concept of the parallel sum of
two matrices  to $k$ matrices, and establish a set of results on
 parallel sums of $k$ matrices.

 \medskip

\noindent {\bf Theorem 23.21.}\, {\em  Let $ A_1, \,  A_2, \,  \cdots, \,  A_k , \, C \in {\cal F}^{m \times n}$ be given. Then 
$$ \displaylines{
\hspace*{0cm}
\max_{A_1^-,\,\cdots,\,A_k^-}r[ \, C -  C( \, A_1^-  + \cdots + A^-_k  \, )C \,] = \min \left\{ \, r(C), \ \ \ 
r( \, N - QCP \, ) + r(C) - r(N) \, \right\},
 \hfill (23.38) 
\cr
\hspace*{0cm}
\min_{A_1^-,\,\cdots,\,A_k^-}r[ \, C -  C( \, A_1^-  + \cdots + A^-_k  \, )C \,] = r(C) + r(N) + r( \, N - QCP \, ) - r[\, N, \ QC \, ] - r \left[ \begin{array}{c} N \\ CP  \end{array} \right], 
 \hfill (23.39) 
\cr}
$$
where $N = {\rm diag}( \, A_1,  \, A_2, \, \cdots, \, A_k \, ), \
P = [\, I_n,  \, I_n, \, \cdots, \, I_n \, ], $ and $ Q = [\, I_m,  \, I_m, \,
\cdots, \, I_m \, ]^T.$ }

\medskip

\noindent {\bf Proof.}\, According to (23.12), (22.1) and (22.2), we easily
find that
\begin{eqnarray*}
\max_{A_1^-,\,\cdots,\,A_k^-}r[ \, C -  C( \, A_1^-  + \cdots + A^-_k  \, )C \,] 
& = & \max_{N^-}r( \, C -  CPN^-QC \, ) \\
& = & \min \left\{ \, r(C), \ \ \  r \left[ \begin{array}{cc} N  & QC \\  CP & C \end{array} \right] -r(N) \  \right\} \\
& = & \min \left\{  \, r(C), \ \ \ r( \, N - QCP \, ) + r(C) - r(N) \, \right\},
\end{eqnarray*} 
\begin{eqnarray*}  
\min_{A_1^-,\,\cdots,\,A_k^-}r[ \, C -  C( \, A_1^-  + \cdots + A^-_k  \,
)C \,] 
& = & \min_{N^-}r( \, C -  CPN^-QC \, ) \\ & = & r (N) - r[ \, N, \ QC \, ] - r \left[ \begin{array}{c} N \\ CP \end{array} \right] + r\left[ \begin{array}{cc}  N  &  QC  \\ CP & C  \end{array} \right] \\
& = & r (N) - r[ \, N, \ QC \, ] - r \left[ \begin{array}{c} N \\ CP \end{array} \right] + r( \, N - QCP \, ) + r(C), 
\end{eqnarray*}  
establishing (23.38) and (23.39). \qquad  $ \Box$  \medskip

Two consequences can directly be derived from (23.38) and (23.39). 

\medskip

\noindent {\bf Theorem 23.22.}\, {\em Let $ A_1, \,    A_2, \,  \cdots, \, A_k,
\, C \in {\cal F}^{m \times n}$ be given. Then there exist $ A^-_i \in
\{ A^-_i \},$ $ i = 1, \, 2, \, \cdots, \, k $  such that $ A_1^-  + A_2^- +
\cdots + A^-_k  \in \{ C^-\},$ if and only if 
$$
r( \, N - QCP \, ) =  r \left[ \begin{array}{c} N \\ CP  \end{array}
\right] + r[\, N, \ QC \,] - r(N) - r(C),
$$
where $ N, \ P $ and $ Q $ are defined in Theorem 23.21. } 

\medskip

\noindent {\bf Theorem 23.23.}\, {\em Let $ A_1, \,  A_2, \,  \cdots,\,  A_k,
\, C \in {\cal F}^{m \times n}$ be given.  The the following four
statements are equivalent$:$

{\rm (a)} \ $ \{ \, A_1^-  + A_2^- \cdots + A^-_k \,\} \subseteq
\{ C^- \}.$

{\rm (b)} \ $ r( \, N - QCP \, ) =  r(N) - r(QCP).$ 

{\rm (c)} \ $ r \left[ \begin{array}{cc} N & QC \\ CP & C \end{array}
\right] = r(N).$

{\rm (d)} \ $ R(C) = r(A_i), \ R(C^T) = r(A_i^T),$ $ i = 1$---$k,$ and
$C = CPN^-QC$, where $ N, \ P$ and $ Q $ are defined in Theorem 23.21.  } 

\medskip

\noindent {\bf Theorem 23.24.}\, {\em  Let $ A_1, \,  A_2, \,  \cdots, \, A_k,
\, C \in {\cal F}^{m \times n}$ be given. Then
$$ \displaylines{
\hspace*{0.5cm}
\min_{A_1^-,\,\cdots,\,A_k^-}r( \, C^- - A_1^-  -\cdots - A^-_k \,)
= r(N) - r[\, N, \ Q \, ] - r \left[ \begin{array}{c} N \\ P \end{array}
\right] + r\left[ \begin{array}{cc} N  &  Q  \\ P & C^-  \end{array}
\right], \hfill
\cr
\hspace*{0.5cm}
\min_{C^-,\, A_1^-,\,\cdots,\,A_k^-}r( \, C^- - A_1^-  -\cdots - A^-_k \,)
= r( C ) + r(N) + r(\, N - QCP \,) - r[\, N, \ QC \, ] -
r \left[ \begin{array}{c} N \\ CP \end{array} \right], \hfill
\cr
\hspace*{0.5cm}
\max_{C^-}\min_{A_1^-,\,\cdots,\,A_k^-}r( \, C^- - A_1^-  -\cdots - A^-_k
\,) = \min \left\{  \,  r(N)  + m - r[\, N, \ Q \,], \ \ \ r(N)  + n -
r \left[ \begin{array}{c} N \\ P \end{array} \right], \right.  \hfill
\cr
\hspace*{5cm}
\left.  m + n + r( \, N- QCP \, ) +  r(N) - r(C) - r[\, N, \ Q \, ] -
r \left[ \begin{array}{c} N \\ P \end{array} \right] \ \right\},
\hfill
\cr}
$$
where $ N, \ P$ and $ Q $ are defined in Theorem 23.21. } 

\medskip

The proof of this theorem is much like that of Theorem 23.17 and is,
therefore, omitted. Two direct consequences of Theorem 23.24 are given below.  

\medskip

\noindent {\bf Theorem 23.25.}\, {\em Let $ A_1, \,  A_2, \,  \cdots,\,  A_k,
\, C \in {\cal F}^{m \times n}$ be given.  Then  for a given $C^-,$
there exist $ A^-_i \in \{ A^-_i \},$ $ i = 1, \, 2, \, \cdots, \, k $
such that
$$\displaylines{
\hspace*{2cm}
 A_1^-  + A_2^-  + \cdots + A_k^- = C^-, \hfill
\cr}
$$
if and only if $ C^- $  satisfies
$$ \displaylines{
\hspace*{2cm}
r\left[ \begin{array}{cc} N  &  Q  \\ P & C^-  \end{array} \right] =
r\left[ \begin{array}{c} N \\ P \end{array} \right] +  r[\, N, \ Q \, ] -
r(N),  \hfill
\cr}
$$ 
where $ N, \ P$ and $ Q $ are defined in Theorem 23.21. } 

\medskip

\noindent {\bf Theorem 23.26.}\, {\em Let $ A_1, \,  A_2, \,  \cdots,\,  A_k,
\ C \in {\cal F}^{m \times n}$ be given. Then the set inclusion
$$ \displaylines{
\hspace*{1.5cm}
  \{ C^-\} \subseteq \{ \, A_1^-  + A_2^-  + \cdots + A_k^- \, \}
  \hfill (23.40)
\cr
\hspace*{0cm}
holds  \ if \ and \ only  \ if  \hfill
\cr
\hspace*{1.5cm}
R(A_1) \cap R(A_2) \cap \cdots \cap R(A_k)  = \{ 0 \}, \ \ \ or \ \ \ 
 R(A_1^T) \cap R(A_2^T) \cap \cdots \cap R(A_k^T)  = \{ 0 \}, \hfill (23.41)   
\cr
\hspace*{0cm}
or \hfill
\cr
\hspace*{1.5cm}
r( \, N- QCP \, ) = r \left[ \begin{array}{c} N \\ P \end{array} \right]
+ r[\, N, \ Q \, ] - r(N) + r(C) - m - n, \hfill (23.42)
\cr}
$$
where $ N, \ P$ and $ Q $ are defined in Theorem 23.21. }

\medskip

The result in Theorem 23.26 implies the following special case.

\medskip

\noindent {\bf Corollary 23.27.}\, {\em Let $ A_1, \,  A_2, \,  \cdots,\,
A_k \in {\cal F}^{m \times n}$ be given. Then the equality
$$ \displaylines{
\hspace*{2cm}
  \{ \, A_1^-  + A_2^-  + \cdots + A_k^- \, \} = \{ 0^-\} =
  {\cal F}^{n \times m} \hfill (23.43)
\cr}
$$ 
holds if and only if $ A_1, \,  A_2, \,  \cdots,\,  A_k$ satisfy
{\rm (23.41)}.}

\medskip

\noindent {\bf Theorem 23.28.}\, {\em Let $ A_1, \,  A_2, \,  \cdots,\,  A_k,
\ C \in {\cal F}^{m \times n}$ be given with
$$ \displaylines{
\hspace*{2cm}
\cap_{i=1}^k R(A_i) \neq \{ 0 \} \ \ \ and \ \ \ \cap_{i=1}^k R(A_i^T)
\neq \{ 0 \}. \hfill (23.44)
\cr
\hspace*{0cm}
Then \ the \ equality \hfill
\cr
\hspace*{2cm}
  \{ \, A_1^-  + A_2^-  + \cdots + A_k^- \, \} = \{ C^-\}  \hfill (23.45) 
\cr}
$$ 
holds if and only if they satisfy the following rank additivity condition  
$$ \displaylines{
\hspace*{2cm}
r\left[ \begin{array}{cc} N  &  Q  \\ P & 0  \end{array} \right] =
r\left[ \begin{array}{c} N \\ P \end{array} \right] + r(Q) =
r[\, N, \ Q \, ] + r(P),  \hfill (23.46)
\cr
\hspace*{0cm}
and \hfill
\cr
\hspace*{2cm}
 C = - [\, 0, \ I_m \, ]\left[ \begin{array}{cc} N  &  Q  \\ P & 0
 \end{array} \right]^- \left[ \begin{array}{c}  0 \\ I_n \end{array}
 \right], \hfill (23.47)
\cr}
$$  
where $ N, \ P$ and $ Q $ are defined in Theorem 23.21. In that case$,$
the rank of $ C $ satisfies the equality
$$ \displaylines{
\hspace*{2cm} 
r(C) = r(N) + r(Q) -  r[\, N, \ Q \, ] = r(N) + r(P) - r\left[ \begin{array}{c} N \\ P \end{array} \right], 
\hfill (23.48)
\cr
\hspace*{0cm} 
or \ more  \ precisely \hfill
\cr 
\hspace*{2cm}
r(C) = {\rm dim}[ \,\cap_{i=1}^k R(A_i) \,] = {\rm dim}[ \,\cap_{i=1}^k
R(A_i^T) \,]. \hfill (23.49)
\cr}
$$ }
{\bf Proof.}\, Assume first that (23.45) holds. Then it follows
from Theorem 23.26 and (23.44) that $ A_1,$ $A_2,$ $ \cdots,$  $ A_k,$  
 and $C$ satisfy (23.42). On the other hand, (23.45) implies that
$$ \displaylines{
\hspace*{2cm}
\min_{A_1^-,\,\cdots,\,A_k^-}r( \, A_1^-  + \cdots + A^-_k \,)
= \min_{C^-}r( C^-), \hfill (23.50)
\cr}
$$ 
which, by the first equality in Theorem 23.24, is equivalent to 
$$\displaylines{
\hspace*{2cm}
r(C) = r \left[ \begin{array}{cc} N & Q \\ Q  & 0  \end{array} \right] -
r \left[ \begin{array}{c} N \\ P \end{array} \right] - r[\, N, \ Q \, ]
+ r(N). \hfill (23.51)
\cr}
$$
Putting it (23.42) yields 
$$\displaylines{
\hspace*{2cm}
r(\,  N - QCP \, ) = r \left[ \begin{array}{cc} N & Q \\ P  & 0
\end{array} \right] - m - n. \hfill (23.52)
\cr}
$$ 
On the other hand, it is easy to verify that 
$$\displaylines{
\hspace*{2cm} 
r \left[ \begin{array}{ccc} N & Q & 0  \\ P  & 0 & I_n \\  0 & I_m & -C
\end{array} \right] = r \left[ \begin{array}{ccc} N - QCP & 0 & 0  \\ 0
& 0 & I_n \\  0 & I_m &  0  \end{array} \right] =  m + n +r(\,  N - QCP
\, ). \hfill
\cr}
$$ 
Thus (23.52) is equivalent to  
$$\displaylines{
\hspace*{2cm} 
r \left[ \begin{array}{ccc} N & Q & 0  \\ P  & 0 & I_n \\  0 & I_m & -C
\end{array} \right] =
 r \left[ \begin{array}{cc} N & Q \\ P  & 0 \end{array} \right]. \hfill
\cr} 
$$ 
In light of (1.5), the rank equality is further equivalent to 
$$ \displaylines{
\hspace*{2cm}
r \left[ \begin{array}{ccc} N & Q & 0  \\ P  & 0 & I_n \end{array}
 \right] = r \left[ \begin{array}{cc} N & Q \\ P  & 0 \end{array}
 \right],   \ \ \ \ r \left[ \begin{array}{cc} N & Q \\ P & 0 \\  0 & I_m
 \end{array} \right] = r \left[ \begin{array}{cc} N & Q \\ P  & 0
 \end{array} \right],  \hfill (23.53)
\cr
\hspace*{0cm}
and \hfill
\cr
\hspace*{2cm}
 C = - [\, 0, \ I_m \, ]\left[ \begin{array}{cc} N  &  Q  \\ P & 0
 \end{array} \right]^- \left[ \begin{array}{c}  0 \\ I_n \end{array}
 \right], \hfill
\cr}
$$
which are exactly (23.46) and (23.47). Consequently  combining
(23.51) with (23.46) yields (23.48), and then yields (23.49) by Lemma
23.1. Conversely if (23.46)---(23.48) hold, then it is easy to verify
that (23.42)  and Theorem 23.23(c) are all satisfied, both of which
imply that (23.45) holds. \qquad  $\Box$

\medskip
 
On the basis of Theorem 23.28, we now can reasonably extend the concept
of parallel sums of two matrices to $k$ matrices.

\medskip

\noindent {\bf Defination.}\, The $k $ matrices $ A_1, \,  A_2, \,  \cdots,\,
 A_k \in {\cal F}^{m \times n}$ are said to be
 parallel summable, if the matrix product 
$$ \displaylines{
\hspace*{2cm}
- [\, 0, \ I_m \, ]\left[ \begin{array}{cc} N  &  Q  \\ P & 0  \end{array}
\right]^- \left[ \begin{array}{c}  0 \\ I_n \end{array} \right] \hfill (23.54)
\cr}
$$
is invariant with respect to the choice of the inner inverse in it,
where $ N, \ P$ and $ Q $ are defined in Theorem 23.21. In that case,
the matrix product in (23.54) is called  the  parallel sum of
$  A_1, \,  A_2, \,  \cdots,\,  A_k $ and denoted by $ p(\, A_1, \,  A_2, \,
\cdots,\,  A_k \, ).$
\medskip

Various properties on parallel sums of $k$ matrices can easily be
derived from Theorem 23.28. Below are some of them, which are quite
 analogous to those for  parallel sums of two matrices.

\medskip
     
\noindent {\bf Theorem 23.29.}\, {\em A null matrix is parallel
summable with any other matrices of the same size$,$ and their
parallel sum is also a null matrix. }

\medskip
 
\noindent {\bf Proof.}\, Let $ A = \left[ \begin{array}{cc} N & Q \\ P & 0
\end{array} \right]$ in (23.54).
 If one of  $ A_1, \,  A_2, \,  \cdots,\,  A_k$ is null, then it is easy to
 verify that
$$\displaylines{
\hspace*{2cm}
r\left[ \begin{array}{cc} N  &  Q  \\ P & 0  \end{array} \right]
= r(N) + m + n. \hfill
\cr}
$$ 
In that case, applying Corollary 21.7(a) to (23.54),
we obtain
\begin{eqnarray*}
\max_{A^-} r\left( \, [\, 0, \ I_m \, ]A^- \left[ \begin{array}{c}  0 \\ I_n \end{array} \right] \, \right) & = & \min \left\{ \ m , \ \ \ n, \ \ \ r \left[ \begin{array}{ccc} N & Q & 0 \\ P & 0 & I_n \\ 0 & I_m & 0 \end{array} \right] -r \left[ \begin{array}{cc} N & Q \\ P & 0 \end{array} \right] \ \right\}\\
 &= & \min \left\{ m , \ \ \ n, \ \ \ m + n +  r(N) -
 \left[ \begin{array}{cc} N & Q \\ P & 0 \end{array} \right] 
 \right\} = 0.
\end{eqnarray*} 
This result implies that (23.54) is always null with respect to the
choice of $ A ^- $. Thus $ A_1, \, A_2, \, \cdots,\, A_k$  are
parallel summable.  \qquad  $\Box $

\medskip

\noindent {\bf Theorem 23.30.}\, {\em Let  $ A_1, \, A_2, \,
\cdots,\, A_k \in {\cal F}^{m \times n}$ be nonnull matrices.  Then they
are  parallel summable if and  only if
$$\displaylines{
\hspace*{2cm} 
R \left[ \begin{array}{c} 0 \\ I_n \end{array} \right] \subseteq R \left[ \begin{array}{cc} N  & Q \\ 
 P & 0  \end{array} \right] \ \ \   and \ \ \ R( \, [\, 0, \ I_m \,]^T \, )
 \subseteq R\left( \, \left[ \begin{array}{cc} N  & Q \\ P & 0  \end{array}
 \right]^T \, \right), \hfill (23.55)
\cr
\hspace*{0cm} 
or \ equivalently \hfill
\cr
\hspace*{2cm} 
r\left[ \begin{array}{cc} N  &  Q  \\ P & 0  \end{array} \right] =
r\left[ \begin{array}{c} N \\ P \end{array} \right] + r(Q) =
r[\, N, \ Q \, ] + r(P),  \hfill (23.56)
\cr}
$$ 
where $ N, \ P$ and $ Q $ are defined in Theorem 23.1. } 

\medskip

\noindent {\bf Proof.}\, It is well-known that (see \cite{RM}, \cite{Mi6}) that a
product $AB^-C $ is invariant with
respect to the choice of $ A^- $ if and only if $ R(A^T) \subseteq R(B^T)$
and  $ R(C) \subseteq R(B)$. Applying this assertion to (23.54)
immediately leads to (23.55). The equivalence of (23.55) and (23.56) is
obvious. \qquad $\Box$

\medskip

\noindent {\bf Theorem 23.31.}\, {\em Let $ A_1, \, A_2, \, \cdots,\, A_k \in
{\cal F}^{m \times n}$  be given. If they are parallel summable$,$ then

{\rm (a)} \ $ \{ \, [\, p(\, A_1, \, A_2, \, \cdots,\,  A_k \,) \, ]^- \, \}
=  \{ \, A_1^-  + A_2^-  + \cdots + A_k^- \, \}.$

{\rm (b)} \ $ p(\, A_1, \,  A_2, \, \cdots,\ A_k \,) =
p(\, A_{i_1}, \,  A_{i_2}, \, \cdots,\,  A_{i_k} \,),$ \  where $i_1, \, i_2, \,
\cdots, \ i_k$ are any permutation of $ 1, \ 2, \ \cdots, \ k$.

{\rm (c)} \ $ p(\, A_1^T, \,  A_2^T, \, \cdots, \,  A_k^T \,) =
[ \, p(\, A_1, \,  A_2, \, \cdots, \, A_k \,) \, ]^T$. }

\medskip

\noindent {\bf Proof.}\,  If any one of $ A_1, \, A_2, \, \cdots,\,  A_k $
is null, then Parts (a)---(c) are naturally valid by Theorem 23.29. Now
suppose that $  A_1, \, A_2, \, \cdots,\,  A_k $ are nonnull and parallel
summable. Then by Theorems 23.30 and 23.28, we immediately see that the
equality in Part (a) holds. The equality in Part (b) comes
from a trivial equality $  \{ \, A_1^-  + A_2^-  + \cdots + A_k^- \, \} =
 \{ \, A_{i_1}^- + A_{i_2}^- + \cdots + A_{i_k}^- \,\}, $ and
 Theorem 21.10(c). By Theorem 23.30, we also know that if
 $ A_1, \, A_2, ,\ \cdots,\,  A_k $ satisfy (23.56), then
$ A_1^T, \, A_2^T, \, \cdots,\,  A_k^T$ naturally satisfy 
$$ \displaylines{
\hspace*{2cm}
r\left[ \begin{array}{cc} N^T  &  P^T  \\ Q^T & 0  \end{array} \right] = r\left[ \begin{array}{c} N^T \\ Q^T \end{array} \right] + r(P^T) =  r[\, N^T, \ P^T \, ] + r(Q^T), \hfill
\cr}
$$ 
where $ N, \ P$ and $ Q $ are defined in Theorem 23.21. Thus
$ A_1^T, \, A_2^T, \, \cdots,\,  A_k^T$ are also parallel
summable.  In that case, it follows from (23.54) that 
\begin{eqnarray*}
p(\, A_1^T, \,  A_2^T, \, \cdots, \,  A_k^T \,) & = & - [\, 0, \ I_n \, ]\left[ \begin{array}{cc} N^T  &  P^T \\ Q^T & 0  \end{array} \right]^- \left[ \begin{array}{c}  0 \\ I_m \end{array} \right] \\
& = & - [\, 0, \ I_n \, ]\left( \, \left[ \begin{array}{cc} N  & Q \\ P & 0  \end{array} \right]^- \, \right)^T \left[ \begin{array}{c}  0 \\ I_m \end{array} \right] \\
& = & - \left( \, [\, 0, \ I_m \, ]\left[ \begin{array}{cc} N  & Q \\ P & 0
  \end{array} \right]^- \left[ \begin{array}{c}  0 \\ I_n \end{array}
  \right] \, \right)^T =  [ \, p(\, A_1, \,  A_2, \, \cdots, \, A_k \,) \,
  ]^T,
\end{eqnarray*}
which is the result in Part (c). \qquad   $\Box$ 

\medskip

\noindent {\bf Theorem 23.32.}\, {\em Let $ A_1, \, A_2,  \, \cdots,\, A_k \in
{\cal F}^{m \times n}$ be given$,$ and  $ B  \in {\cal F}^{m \times m},$ $ C \in {\cal F}^{n\times n}$ are two nonsingular matrices. Then
 $ A_1, \, A_2, \, \cdots,\, A_k $ are parallel summable if and only if
 $ BA_1C, \, BA_2C, \, \cdots,\, BA_kC $ are are parallel summable. In that
 case$,$
$$ \displaylines{
\hspace*{2cm}
p(\, BA_1C, \, BA_2C, \, \cdots,\, BA_kC \, )  =
Bp(\, A_1, \, A_2, \, \cdots,\,  A_k \,)C. \hfill (23.57)
\cr}
$$  }
{\bf Proof.}\, If any one of $ A_1, \, A_2, \, \cdots,\,  A_k $
is null, then (23.57) is a trivial result by Theorem 23.30. Now suppose
that $  A_1, \, A_2, \, \cdots,\,  A_k $ are nonnull and denote
$$ \displaylines{
\hspace*{2cm}
\widehat{B} = {\rm diag}( \, B, \, B, \, \cdots,\, B \, ),  \ \ \ \widehat{C} = {\rm diag}( \, C, \, C, \, \cdots,\, C \, ).    
 \hfill 
\cr}
$$ 
Since $ B $ and $ C $ are nonsingular, $\widehat{B} $ and $  \widehat{C}$ are nonsingular, too. In that case, it is easy to verify that 
$$
\displaylines{
\hspace*{2cm}
 r\left[ \begin{array}{cc} \widehat{B}N \widehat{C} &  Q  \\ P & 0  \end{array} \right]
= r\left[ \begin{array}{cc} N &  \widehat{B}^{-1}Q  \\ P\widehat{C}^{-1} & 0  \end{array} \right] 
= r\left[ \begin{array}{cc} N &  \widehat{B}^{-1}QB  \\ CP\widehat{C}^{-1} & 0  \end{array} \right] = 
r\left[ \begin{array}{cc} N &  Q  \\ P & 0  \end{array} \right], \hfill
\cr
\hspace*{2cm} 
r\left[ \begin{array}{c} \widehat{B}N \widehat{C} \\ P \end{array} \right]
= r\left[ \begin{array}{c} N  \\ P\widehat{C}^{-1} \end{array} \right] 
= r\left[ \begin{array}{c} N \\ CP\widehat{C}^{-1} \end{array} \right] = 
r\left[ \begin{array}{c} N \\ P  \end{array} \right],  \hfill
\cr
\hspace*{2cm}
r[\, \widehat{B}N \widehat{C}, \ Q \,] = r[\, N, \ \widehat{B}^{-1}Q \,] = r[\, N, \ \widehat{B}^{-1}QB\,] = r[\, N, \ Q \,]. \hfill
\cr}
$$
Combining them with (23.56) clearly shows that $ A_1,  \, A_2, \, \cdots, \,
A_k $ are parallel summable if and only if $ BA_1C, \, BA_2C, \, \cdots,\,
BA_kC $ are parallel summable.  From the nonsingularity of  $ B $ and
$ C $, we also see that
$$\displaylines{
\hspace*{2cm}
 \left[ \begin{array}{cc} \widehat{B}N \widehat{C} &  Q  \\ P & 0  \end{array} \right]^- 
= \left[ \begin{array}{cc} \widehat{C}^{-1} &  0 \\ 0 & B \end{array} \right]
\left[ \begin{array}{cc} N &  Q  \\ P & 0  \end{array} \right]^-
\left[ \begin{array}{cc} \widehat{B}^{-1} &  0 \\ 0 & C\end{array} \right]. \hfill
 \cr}
$$    
Thus it follows from (23.54) that 
\begin{eqnarray*}
p(\, BA_1C, \ BA_2C, \, \cdots,\, BA_kC  \,) & = & - [\, 0, \ I_m \, ]\left[ \begin{array}{cc} \widehat{B}N \widehat{C} &  Q  \\ P & 0  \end{array} \right]^- \left[ \begin{array}{c}  0 \\ I_n \end{array} \right] \\
& = & - [\, 0, \ I_m \, ]\left[ \begin{array}{cc} \widehat{C}^{-1} &  0
\\ 0 & B \end{array} \right]\left[ \begin{array}{cc} N &  Q  \\ P & 0
\end{array} \right]^- \left[ \begin{array}{cc} \widehat{B}^{-1} &  0 \\ 0 &
C \end{array} \right] \left[ \begin{array}{c}  0 \\ I_n \end{array} \right]
\\ & = & - B[\, 0, \ I_m \, ]\left[ \begin{array}{cc} N  &  Q  \\ P & 0  \end{array} \right]^- \left[ \begin{array}{c}  0 \\ I_n \end{array} \right]C \\ 
& = & Bp(\, A_1, \, A_2, \, \cdots,\,  A_k \,)C,
\end{eqnarray*}
which is (23.57). \qquad $ \Box$

\markboth{YONGGE  TIAN }
{24. RANKS AND INDEPENDENCE OF SOLUTIONS TO $ BXC = A$}

\chapter{Ranks and independence of submatrices in solutions to $ BXC = A$}

\noindent Suppose that $ BXC = A $ is a consistent matrix equation over an arbitrary
field ${\cal F}$, where $A \in {\cal F}^{m \times n}, \ B \in
{\cal F}^{m \times k}$ and $ C \in {\cal F}^{l \times n}$ are given. Then
it can factor in the form
$$\displaylines{
\hspace*{2cm}
[ \, B_1, \ B_2 \, ] \left[ \begin{array}{cc}  X_1  & X_2 \\ X_3 & X_4
 \end{array} \right]
\left[ \begin{array}{c} C_1 \\ C_2 \end{array} \right] = A,   \hfill (24.1) 
\cr}
$$  
where $ X_1 \in {\cal F}^{k_1 \times l_1}, \, X_2 \in {\cal F}^{k_1
 \times l_2}, \, X_3 \in {\cal F}^{k_2 \times l_1} $  and $ X_4 \in
  {\cal F}^{k_2 \times l_2},$ $k_1  + k_2 = k, \, l_1 + l_2 = l$. 
In this chapter, we determine  maximal and minimal possible ranks of
submatrices $ X_1$---$X_4$ in a solution to (24.1).
   
Possible ranks of solutions of linear matrix equations and  various related topics  have been considered
 previously by several authors. For example, Mitra in \cite{Mi1} examined solutions with fixed ranks for 
the matrix equations $AX = B$ and  $ AXB = C$;  Mitra in \cite{Mi2} 
 gave common solutions of minimal rank of the pair of matrix equations $ AX = C$, $XB = D $; Uhlig in 
\cite{Uh} presented maximal and minimal possible ranks of solutions of the equation $ AX= B$;
 Mitra \cite{Mi6} described common solutions with the minimal rank to the pair of matrix equations 
$ A_1XB_1 = C_1$ and $ A_2XB_2 = C_2$.  Besides the work in the chapter, we shall also consider in the next two 
chapters possible ranks of the  two real matrices $ X_0$ and $X_1$ in solutions to the complex matrix equation 
 $ B(\, X_0 + iX_1 \,)C = A$, as well as possible ranks and independence of solutions to the matrix equation 
 $ B_1XC_1  + B_2YC_2 = A$. 

For convenience of representation, we adopt the notation for the collections
of the submatrices $X_1$---$X_4$ in (24.1)
$$ \displaylines{
\hspace*{2cm}
S_i = \left\{ \ X_i \ \left| \ [ \, B_1, \ B_2 \, ]
\left[ \begin{array}{cc}  X_1  & X_2 \\ X_3 & X_4 \end{array} \right]
\left[ \begin{array}{c} C_1 \\ C_2 \end{array} \right] = A \right.
\ \right\}, \qquad  i = 1, \, 2, \, 3, \, 4 . \hfill (24.2)
\cr}
$$  

It is easily seen that $ X_1$---$X_4$ in (24.1) can be written as 
$$\displaylines{
\hspace*{2cm}
X_1 =  [ \, I_{k_1}, \ 0 \, ] X \left[ \begin{array}{c} I_{l_1} \\  0
\end{array} \right] = P_1XQ_1, \qquad
 X_2 =  [ \, I_{k_1}, \ 0  \, ] X \left[ \begin{array}{c} 0 \\  I_{l_2}
  \end{array} \right] = P_1XQ_2, \hfill (24.3)
\cr
\hspace*{2cm}
X_3 =  [ \, 0, \ I_{k_2} \, ] X \left[ \begin{array}{c} I_{l_1} \\  0
 \end{array} \right] = P_2XQ_1, \ \ \ \
 X_4 =  [ \, 0, \ I_{k_2} \, ] X \left[ \begin{array}{c} 0 \\  I_{l_2}
  \end{array} \right] = P_2XQ_2.
\hfill (24.4)
\cr}
$$ 

Since $ BXC = A$ is consistent, its general solution can be written
as $ X = B^-AC^- + F_BV + WE_C$. Putting it in (24.3) and (24.4) yields the general expressions of
 $ X_1$---$X_4$ as follows 
$$\displaylines{
\hspace*{2cm}
X_1 = P_1X_0Q_1 +P_1F_BV_1  + W_1E_CQ_1, \ \ \
  X_2 = P_1X_0Q_2 +P_1F_BV_2 + W_1E_CQ_2,
\hfill (24.5)
\cr
\hspace*{2cm}
X_3 = P_2X_0Q_1 +P_2F_BV_1  + W_2E_CQ_1, \ \ \   X_4 = P_2X_0Q_2
+P_2F_BV_2  + W_2E_CQ_2,
\hfill (24.6)
\cr}
$$ 
where $ X_0 = B^-AC^-, \ V =[\, V_1, \ V_2 \, ]$ and $ W = \left[
\begin{array}{c} W_1 \\ W_2  \end{array} \right]$.

\medskip

\noindent {\bf Theorem 24.1.}\, {\em Suppose that the  matrix equation
 {\rm (24.1)} is  consistent. Then 
$$
\displaylines{
\hspace*{2cm}
\max_{X_1 \in S_1}r(X_1) = \min \left\{ \ k_1, \ \ \ l_1, \ \ \
r \left[ \begin{array}{cc}  A  & B_2 \\
 C_2 & 0  \end{array} \right] - r(B) - r(C) + k_1 + l_1  \, \right\}, \hfill
 (24.7)
 \cr
\hspace*{2cm}
\min_{X_1 \in S_1}r(X_1) =  r \left[ \begin{array}{cc}  A  & B_2 \\
C_2 & 0  \end{array} \right] - r(B_2) - r(C_2). \hfill (24.8)
\cr}
$$ } 

\noindent {\bf Proof.}\, It is quite obvious that to determine maximal and
minimal ranks of $X_1$ in (24.1) is in fact to determine maximal and
minimal ranks of $P_1XQ_1$ subject to the consistent equation $ BXC = A $.
Thus applying (20.3) and (20.4) to $ X_1 = P_1XQ_1$  produces the
following two expressions
$$
\displaylines{
\hspace*{1.5cm}
\max_{BXC = A}r( P_1XQ_1) = \min \left\{ \,  r( P_1 ), \ \ \  r( Q_1 ),
 \ \ \
 r \left[ \begin{array}{ccc} 0  & 0 & P_1 \\  0 & A   & B  \\  Q_1 & C
  & 0 \end{array} \right]
-r(B) -r(C) \, \right\},   \hfill
\cr
\hspace*{1.5cm}
\min_{BXC = A}r( P_1XQ_1) = r \left[ \begin{array}{ccc} 0  & 0  & P_1 \\
 0 & A   & B  \\  Q_1 & C   & 0 \end{array} \right] -
 r \left[ \begin{array}{c} P_1  \\ B  \end{array} \right] -
 r[ \, Q_1, \  C \, ].  \hfill
\cr}
$$
Putting  the given matrices $ B = [ \, B_1, \  B_2 \, ], \
C = \left[ \begin{array}{c} C_1 \\ C_2  \end{array} \right], \ P_1$ and
$Q_1$ in them and simplifying yields the desired formulas (24.7) and
(24.8). The details are omitted. \qquad  $ \Box$

\medskip

Maximal and minimal ranks of the submatrices $X_2, \, X_3,$  and $ X_4$ in
 (24.1) can also be derived in the same manner. We omit them here
 for simplicity. The two formulas in (24.7) and (24.8) can help to
 characterize structure of solutions to (24.1).  Next are some of them.

\medskip
 
\noindent {\bf Corollary 24.2.}\, {\em Suppose that the  matrix equation
{\rm (24.1)} is consistent. Then 

{\rm (a)}\, Eq.\,{\rm (24.1)} has a solution with the form
$ X =  \left[\begin{array}{cc} 0  & X_2 \\
 X_3 & X_4 \end{array} \right] $ if and only if $ r \left[ \begin{array}{cc}  A  & B_2 \\ C_2 & 0  \end{array} \right] =  r(B_2) + r(C_2).$

{\rm (b)}\,  All the solutions of {\rm (24.1)} have the form
$ X =  \left[\begin{array}{cc} 0  & X_2 \\
 X_3 & X_4 \end{array} \right] $ if and only if 
$$\displaylines{
\hspace*{2cm}
 \left[ \begin{array}{cc}  A  & B_2 \\ C_2 & 0 \end{array} \right] = r(B)
+ r(C) - k_1 - l_1, \hfill (24.9)
\cr}
$$
or equivalently
$$
r \left[ \begin{array}{cc}  A  & B_2 \\ C_2 &  0 \end{array} \right] =
 r(B_2) + r(C_2), \ r(B_1) = k_1, \ \   r(C_1) = l_1, \ \
R(B_1) \cap R(B_2) = \{ 0\} \ \ and \ \
R(C_1^T) \cap R(C_2^T) = \{ 0\}. \eqno (24.10)
$$}
{\bf Proof.}\, Part (a) and (24.9) follows directly from
 (24.7) and (24.8). On the other hand, observe that 
\begin{eqnarray*}
\lefteqn {r \left[ \begin{array}{cc}  A  & B_2 \\ C_2 & 0  \end{array}
 \right] - r(B)- r(C) + k_1 + l_1 }\\
& = &  \left( \, r \left[ \begin{array}{cc}  A  & B_2 \\ C_2 & 0
\end{array} \right] - r(B_2) - r(C_2) \, \right) +  [ \, k_1 + r(B_2) -
r(B) \, ] + [ \, l_1  + r(C_2) - r(C) \, ].
\end{eqnarray*} 
Thus (24.9) is equivalent to (24.10).  \qquad  $ \Box$

\medskip
 
\noindent {\bf Theorem 24.3.}\, {\em Suppose that the  matrix equation
{\rm (24.1)} is  consistent.  Then

{\rm (a)}\, Eq.\,{\rm (24.1)} has a solution with the form  $ X =
 \left[ \begin{array}{cc} X_1  &  0 \\
 X_3 &  0 \end{array} \right] $ if and only if $ R(A^T) \subseteq
 R(C^T_1).$

{\rm (b)}\, Eq.\,{\rm (24.1)} has a solution with the form  $ X =
 \left[ \begin{array}{cc} X_1  &  X_2 \\
 0 & 0 \end{array} \right] $ if and only if $ R(A)
 \subseteq R(B_1).$

{\rm (c)}\, Eq.\,{\rm (24.1)} has a solution with the form  $ X =
 \left[ \begin{array}{cc} X_1  & 0 \\
 0 & 0 \end{array} \right] $ if and only if $ R(A)
 \subseteq R(B_1)$ and $ R(A^T) \subseteq R(C^T_1).$ } 

\medskip

\noindent {\bf Proof.}\, According to (20.3) and (20.4), we find that 
$$ \displaylines{
\hspace*{2cm}
\min_{BXC = A} r\left[ \begin{array}{c}  X_2  \\ X_4  \end{array}
 \right]  = \min_{BXC = A} r( XQ_2)
= r\left[ \begin{array}{c}  A  \\ C_1  \end{array} \right] - r(C_1), \hfill
\cr
\hspace*{2cm}
\min_{BXC = A} r[ \,  X_3, \  X_4 \,] = \min_{BXC = A} r(P_2X) =
r[ \, A, \ B_1 \, ] - r(B_1). \hfill
\cr}
$$ 
Thus we have Parts (a) and (b). The result in Part (c) is evident.
\qquad $ \Box$ 

\medskip

Note from (24.5) and (24.6) that $ X_1$ and $X_4,$ $ X_2$ and $X_3$
are independent in their expressions, i.e., both of them do not involve
 the same  variant matrices, thus we have the following.

\medskip

 \noindent {\bf Theorem 24.4.}\, {\em Suppose that the  matrix equation
 {\rm (24.1)} is  consistent.  Then

{\rm (a)}\, Eq.\,{\rm (24.1)} must have two solutions with the forms  
$$ \displaylines{
\hspace*{2cm}
X =  \left[ \begin{array}{cc} \widehat{X_1} & X_2 \\ X_3 & \widehat{X_4}
\end{array} \right], \qquad X =  \left[ \begin{array}{cc} X_1 & \widehat{X_2} \\ \widehat{X_3}
 & X_4 \end{array} \right], \hfill
\cr}
$$
where $ \widehat{X_1}$---$\widehat{X_4}$ in them with the ranks 
$$ \displaylines{
\hspace*{2cm}
r(\widehat{X_1} ) = \min_{X_1 \in S_1}r(X_1) =  r \left[ \begin{array}{cc}
  A  & B_2 \\ C_2 & 0  \end{array} \right] - r(B_2)
- r(C_2), \hfill
\cr
\hspace*{2cm}
 r(\widehat{X_2} ) = \min_{X_2 \in S_2}r(X_2) =  r \left[ \begin{array}{cc}
 A  & B_2 \\ C_1 & 0  \end{array} \right] - r(B_2) - r(C_1), \hfill
\cr
\hspace*{2cm}
r(\widehat{X_3}) = \min_{X_3 \in S_3}r(X_3) =  r \left[ \begin{array}{cc}
 A  & B_1 \\ C_2 & 0  \end{array} \right] - r(B_1) - r(C_2), \hfill
\cr
\hspace*{2cm}
r(\widehat{X_4}) = \min_{X_4 \in S_4}r(X_4) =  r \left[ \begin{array}{cc}
 A  & B_1 \\ C_1 & 0  \end{array} \right] - r(B_1) - r(C_1). \hfill
\cr}
$$

{\rm (b)}\, Eq.\,{\rm (24.1)} has a solution with the form  $ X =  \left[
 \begin{array}{cc} 0  &  X_2 \\
 X_3 & 0 \end{array} \right], $ if and only if
$$ \displaylines{
\hspace*{2cm}
 r \left[ \begin{array}{cc}  A  & B_1 \\ C_1 & 0  \end{array} \right] =
  r(B_1) + r(C_1), \ \  and \ \
  r \left[ \begin{array}{cc}  A  & B_2 \\ C_2 & 0  \end{array} \right] =
   r(B_2) + r(C_2). \hfill
\cr}
$$ 

{\rm (c)}\, Eq.\,{\rm (24.1)} has a solution with the form  $ X =
\left[ \begin{array}{cc} X_1 &  0 \\
 0  & X_4 \end{array} \right], $ if and only if
$$ \displaylines{
\hspace*{2cm}
 r \left[ \begin{array}{cc}  A  & B_1 \\ C_2 & 0  \end{array} \right] =
  r(B_1) + r(C_2), \ \  and \ \
  r \left[ \begin{array}{cc}  A  & B_2 \\ C_1 & 0  \end{array} \right]
  = r(B_2) + r(C_1). \hfill
\cr}
$$ } 
\hspace*{0.4cm} The result in Theorem 24.4(c) in fact implies a necessary and sufficient
condition for the matrix equation
$ B_1X_1C_1 + B_2X_4C_2  = A$ to be solvable, which was first established
 by \"{O}zg\"{u}ler in \cite{Oz}.

The uniqueness of the submatrices $ X_1$---$X_4$ in (24.1) can be
determined by (24.5) and (24.6).

\medskip

\noindent {\bf Theorem 24.5.}\, {\em Suppose that the  matrix equation
{\rm (24.1)} is  consistent. The submatrix  $ X_1$  in {\rm (24.1)} is
unique if and only if {\rm (24.1)} satisfies the
following four conditions 
$$\displaylines{
\hspace*{2cm}
  r(B_1) = k_1, \ \ \ r(C_1) = l_1, \ \ \  R(B_1) \cap R(B_2)
   = \{ 0\},  \ \ \ R(C_1^T) \cap R(C_2^T) = \{ 0\}.
   \hfill (24.11)
\cr}
$$}
{\bf Proof.}\ It is easy to see from (24.5) that $ X_1$ is
unique if and only if $P,_1F_B  = 0$ and $ E_CQ_1 = 0,$
where we find by (1.2) and (1.3) that
$$\displaylines{
\hspace*{0.5cm}
P_1F_B  = 0  \Rightarrow r\left[ \begin{array}{cc} P_1 \\ B  \end{array}
\right] = r(B) \Rightarrow  k_1 + r(B_2) = r(B) \Rightarrow r(B_1) = k_1
\ {\rm and}  \ R(B_1) \cap R(B_2) = \{ 0\}, \hfill
\cr
\hspace*{0.5cm}
E_CQ_1 = 0 \Rightarrow r[\, Q_1, \ C\, ] = r(C) \Rightarrow l_1 + r(C_2)
 = r(C) \Rightarrow  r(C_1) = l_1 \ {\rm and} \ R(C_1^T) \cap R(C_2^T) =
 \{ 0\}. 
\cr}
$$ 
Thus we have (24.11).  \qquad $ \Box$  

\medskip

The following result is concerning the independence of submatrices in
solutions to (24.1).

\medskip

\noindent {\bf Theorem 24.6.}\, {\em Suppose that the  matrix equation
{\rm (24.1)} is  consistent with $ B \neq 0 $ and $ C \neq 0 $.

{\rm (a)}\, Consider $S_1$---$S_4$ in {\rm (24.1)} as four independent
 matrix sets.  Then
$$
\displaylines{
\hspace*{1cm}
\max_{X_i \in S_i} r \left(  A - [ \, B_1, \ B_2 \, ]
 \left[ \begin{array}{cc}  X_1  & X_2 \\ X_3 & X_4 \end{array} \right]
  \left[ \begin{array}{c} C_1 \\ C_2 \end{array} \right] \right) \hfill
  \cr
\hspace*{1cm}
 =  \min \left\{ r(B), \ \ \ r(C), \ \ \ r(B_1) + r(B_2) - r(B)
 + r(C_1) + r(C_2) - r(C) \right\}. \hfill  (24.12)
\cr}
$$

{\rm (b)}\, The four submatrices $ X_1$---$ X_4$  in {\rm (24.1)} are
independent$,$ that is$,$ for any choice of $ X_i \in S_i (i = 1, \, 2, \, 3,
\ 4),$ the corresponding matrix $ X =  \left[ \begin{array}{cc}  X_1  & X_2 \\ X_3 &
X_4 \end{array} \right] $ is a solution of {\rm (24.1)}$,$  if and only if
$$\displaylines{
\hspace*{1cm}
 R(B_1) \cap R(B_2) = \{ 0\} \ \ and  \ \  R(C_1^T)
 \cap R(C_2^T) = \{ 0\}.
\hfill (24.13)
\cr}
$$ }
{\bf Proof.}\, According to (24.5) and (24.6), the general
expressions of  $ X_1$---$ X_4$ in $S_1$---$S_4$ can independently be
written as
$$
X_1 = P_1X_0Q_1 +P_1F_BV_1  + W_1E_CQ_1, \ \ \   X_2 = P_1X_0Q_2 +
P_1F_BV_2 + W_2E_CQ_2,
$$ 
$$
X_3 = P_2X_0Q_1 +P_2F_BV_3  + W_3E_CQ_1, \ \ \
X_4 = P_2X_0Q_2 +P_2F_BV_4  + W_4E_CQ_2,
$$ 
where $ X_0 = B^-AC^-,$  $V_1$---$V_4$ and $W_1$---$W_4$ are arbitrary.
Putting them in $ X $ yields
$$
\displaylines{
\hspace*{1cm}
\left[ \begin{array}{cc}  X_1  & X_2 \\ X_3 & X_4 \end{array}
 \right] \hfill
\cr
\hspace*{1cm}
=  \left[ \begin{array}{c}  P_1  \\ P_2 \end{array} \right]X_0[ \, Q_1, \ Q_2 \, ] +
\left[ \begin{array}{cc}  P_1F_B  &  0 \\  0 &  P_2F_B \end{array}
\right]\left[ \begin{array}{cc} V_1  &
 V_2 \\ V_3 & V_4 \end{array} \right] + \left[ \begin{array}{cc} W_1  & W_2
  \\ W_3 & W_4 \end{array}
 \right] \left[ \begin{array}{cc} E_CQ_1  & 0 \\ 0 & E_CQ_2  \end{array}
 \right] \hfill
\cr
\hspace*{1cm}
= X_0 + GV + WH,  \hfill
\cr}
$$
where $ G = {\rm diag}( \, P_1F_B, \  P_2F_B \, ),  \, H ={\rm diag}( \,
 E_CQ_1, \,  E_CQ_2 \, ).$  Applying (1.6) to it, we find
$$\displaylines{
\hspace*{1cm}
 \max_{X_i \in S_i} r \left(  \,  A - [\, B_1, \ B_2 \, ]
\left[ \begin{array}{cc}  X_1  & X_2 \\ X_3 & X_4 \end{array} \right]
\left[ \begin{array}{c} C_1 \\ C_2 \end{array} \right] \,  \right)
\hfill
\cr
\hspace*{1cm}
= \max_{V, \, W}r( \, BGVC + BWHC \, ) = \min \left\{ \, r(B), \ \
 r(C),  \ \ r(BG) + r(HC) \, \right\}. \hfill (24.14)
\cr}
$$
According to (1.2) and (1.3), we see that
\begin{eqnarray*}
r(BG) & = & r[ \, B_1P_1F_B, \  B_2P_2F_B \, ] \\
& = & r \left[ \begin{array}{cc}  B_1P_1 & B_2P_2  \\ B & 0 \\ 0 & B
\end{array} \right] - 2r(B) = r\left[ \begin{array}{cccc}  B_1 & 0 & 0 & B_2  \\ B_1 & B_2 & 0 & 0
 \\  0 & 0 & B_1 & B_2 \end{array} \right] - 2r(B) = r(B_1) +r(B_2) -r(B),
\end{eqnarray*}  
$$
\displaylines{
\hspace*{0cm} r(HC)  \hfill
\cr
\hspace*{0cm}
 = r\left[ \begin{array}{c} E_CQ_1C_1 \\ E_CQ_2C_2 \end{array} \right]
=r\left[ \begin{array}{ccc}  Q_1C_1 & C & 0   \\ Q_2C_2 & 0 & C \end{array}
\right] - 2r(C)  = r\left[ \begin{array}{ccc}  C_1 & C_1 & 0  \\ 0 & C_2 & 0  \\  0 & 0 & C_1 \\  C_2 & 0 & C_2 \end{array} \right] - 2r(C)
 = r(C_1) +r(C_2) -r(C). \hfill
\cr}
$$
Putting them in (24.14), we obtain (24.12).  The result in Part (b) is a
direct consequence of (24.12). \qquad $ \Box$

\medskip

Let 
$$\displaylines{
\hspace*{2cm}
 M = \left[ \begin{array}{cc}  A  & B \\ C & D \end{array} \right] \hfill
 (24.15)
\cr}
$$  
be a partitioned matrix over ${\cal C}$, where $A \in {\cal F}^{m \times n},
 \, B \in {\cal F}^{m \times k}, \,  C \in {\cal F}^{l \times n}$ and $ D \in
  {\cal F}^{l \times k}$, and  write its inner inverse in the block  form
$$
\displaylines{
\hspace*{2cm}
 M^- = \left[ \begin{array}{cc}  G_1  & G_2 \\ G_3 & G_4 \end{array}
 \right], \hfill (24.16)
\cr}
$$ 
where $ G_1 \in {\cal F}^{n \times m}.$ In this section, we determine
 maximal and minimal ranks of the blocks $ G_1$---$G_4$ in (24.16) and
 consider their relationship with $A, \ B, \
C$ and $ D $.

For convenience of representation, we adopt the notation 
$$ \displaylines{
\hspace*{2cm}
T_i = \left\{ \, G_i \ \left| \ \left[ \begin{array}{cc}  G_1  & G_2 \\ G_3
& G_4 \end{array} \right] \in
\{ M^- \}  \right.  \, \right\}, \ \ \ \  i = 1, \, 2, \, 3, \, 4 .
\hfill (24.17)
\cr}
$$ 
\hspace*{0.4cm} Notice that $M^-$ is in fact a solution to the matrix equation $ MXM = M$.
Thus applying the results in Theorem 24.1 to (24.15) and (24.16), we
find the following.

\medskip

\noindent {\bf Theorem 24.7.}\, {\em Let $ M $ and $ M^-$ be given
by {\rm (24.15)} and {\rm (24.16)}. Then 
$$
\displaylines{
\hspace*{2cm}
\max_{G_1 \in T_1}r(G_1) = \min \left\{ \, m , \ \ \ n, \ \ \ m + n + r(D)
- r(M) \, \right\}, \hfill (24.18)
\cr
\hspace*{2cm}
\min_{G_1 \in T_1}r(G_1) = r(M) + r(D) - r[\, C, \ D \, ] -
r \left[ \begin{array}{c} B \\ D  \end{array} \right]. \hfill  (24.19)
\cr}
$$}
{\bf Proof.}\, Follows from (24.7) and (24.8). \qquad
 $ \Box$

\medskip

\noindent {\bf Corollary 24.8.}\, {\em Let $ M $ and $ M^-$ be given
by {\rm (24.15)} and {\rm (24.16)}.  Then

{\rm (a)}\,  $M$ has a g-inverse  with the form  $ M^- =
\left[ \begin{array}{cc} 0  & G_2 \\
 G_3 & G_4 \end{array} \right] $ if and only if \ $ r(M) =
  r \left[ \begin{array}{c} B \\ D  \end{array} \right] +
  r[\, C, \ D \, ] - r(D).$

{\rm (b)}\, All the g-inverses of $M$ have the form $ M^- =
 \left[ \begin{array}{cc} 0  & G_2 \\
 G_3 & G_4 \end{array} \right] $ if and only if \  $ r(M) =
 m + n -r(D).$ }

 \medskip

\noindent {\bf Proof.}\, Follows from Theorem 24.7.   \qquad $
\Box$

\medskip

\noindent {\bf Corollary 24.9.}\, {\em Let $ M $ and $ M^-$ be given by
 {\rm (24.15)} and {\rm (24.16)}.

{\rm (a)}\,  $M$ has a g-inverse  with the form  $ M^- =
\left[ \begin{array}{cc} G_1  & 0 \\
 G_3 &  0 \end{array} \right] $ if and only if \  $ R( \,
 [\, C, \ D\, )^T \, ]
\subseteq R( \, [\, A, \ B \, ]^T \, ).$

{\rm (b)}\,  $M$ has a g-inverse  with the form  $ M^- =  \left[
\begin{array}{cc} G_1  & G_2 \\
 0 &  0 \end{array} \right] $ if and only if \  $ R\left[
 \begin{array}{c} B \\ D  \end{array} \right]
 \subseteq R \left[ \begin{array}{c} A \\ C  \end{array}
 \right].$

{\rm (c)}\,  $M$ has a g-inverse  with the form  $ M^- =  \left[
\begin{array}{cc} G_1  &  0 \\
 0 &  0 \end{array} \right] $ if and only if \  $r(M) = r(A)$.  }

\medskip

\noindent {\bf Proof.}\, Follows from Corollary 24.3.   \qquad
  $ \Box$ 

\medskip

\noindent {\bf Corollary 24.10.}\, {\em Let $ M $ and $ M^-$ be given by
 {\rm (24.15)} and {\rm (24.16)}. Then

{\rm (a)}\,  $M$ has a g-inverse  with the form  $ M^- =
\left[ \begin{array}{cc} G_1  & 0 \\
 0  & G_4 \end{array} \right], $ if and only if 
$$\displaylines{
\hspace*{2cm}
 r(M) =  r \left[ \begin{array}{c} A \\ C \end{array} \right] +
 r[\, C, \ D \, ] - r(C)
= r \left[ \begin{array}{c} B \\ D \end{array} \right] + r[\, A, \ B \, ]
 - r(B). \hfill
\cr}
$$

{\rm (b)}\,  $M$ has a g-inverse  with the form  $ M^- =  \left[
\begin{array}{cc} 0  & G_2 \\
 G_3  & 0 \end{array} \right],$ if and only if 
$$\displaylines{
\hspace*{2cm}
 r(M) =  r \left[ \begin{array}{c} B \\ D \end{array} \right] + r[\, C,
 \ D \, ] - r(D)
= r \left[ \begin{array}{c} A \\ C \end{array} \right] + r[\, A, \ B \, ]
 - r(A). \hfill
\cr}
$$ }
{\bf Proof.}\, Follows from Theorem 24.4 (b) and (c).  \qquad
$ \Box$

\medskip

\noindent {\bf Corollary  24.11} (Rao and Yanai \cite{RY}).\, {\em Let $ M $ and $ M^-$ be given
by {\rm (24.15)} and {\rm (24.16)}.  Then the submatrix  $ G_1$
in {\rm (24.16)}
is unique if and only if $M$ satisfies the following three conditions
$$\displaylines{
\hspace*{2cm}
 r[\, A, \ B \, ] = m, \ \ \  r \left[ \begin{array}{c} A \\ C \end{array}
 \right] = n,  \ \ \
 r(M) = n + r \left[ \begin{array}{c} B \\ D \end{array} \right] = m + r[\,
  C, \ D \, ]. \hfill
\cr}
$$ } 
{\bf Proof.}\, Follows from Theorem 24.5.   \qquad  $ \Box$

 \medskip

\noindent {\bf Theorem 24.12.}\, {\em Let $ M $ and $ M^-$ be given by
 {\rm (24.15)} and {\rm (24.16)}.  

{\rm (a)}\, Consider $T_1$---$T_4$ in {\rm (24.17)} as
 four independent
matrix sets.  Then
$$
\max_{G_i \in T_i} r \left( M - M \left[ \begin{array}{cc}
  G_1  & G_2 \\ G_3 & G_4 \end{array} \right] M  \right)  
 = \min \left\{ \ r(M),  \ \ r \left[ \begin{array}{c} A \\ C
 \end{array} \right] +
 r \left[ \begin{array}{c} B \\ D \end{array} \right] + r[\, A, \ B \, ]
 + r[\, C, \ D \,] - 2r(M) \, \right\}.
 $$

{\rm (b) (Rao and Yanai \cite{RY})}\, The four submatrices $ G_1$---$ G_4$
in {\rm (24.16)}
are independent if and only if $M $ satisfies the following rank
additivity condition
$$\displaylines{
\hspace*{2cm}
r(M) = r \left[ \begin{array}{c} A \\ C \end{array} \right] + r \left[
 \begin{array}{c} B \\ D \end{array} \right]  = r[\, A, \ B \, ] + r[\, C, \
  D \, ]. \hfill
\cr}
$$ }
{\bf Proof.}\, Follows from Theorem 24.6.  \qquad  $ \Box$

\medskip

In the remainder of this section, we consider the relationship between
 $ \{A^- \} $ and $ T_1$,
 $ \{B^- \} $ and $ T_3$,   $ \{C^- \} $ and $ T_2$, $ \{D^- \} $ and $
 T_4$, where $ T_1$---$T_4$ are defined in (24.17). 

 \medskip

\noindent {\bf Theorem 24.13.}\, {\em Let $ M $ and $ M^-$ be given by
 {\rm (24.15)} and {\rm (24.16)}. Then 
$$
\displaylines{
\hspace*{1cm}
\max_{G_1 \in T_1}r( \, A - AG_1A \,) = \min \left\{ \,  r(A), \ \ \ \
r(A) +
 r\left[ \begin{array}{cc} 0  &  B \\ C & D \end{array} \right] - r(M)
  \, \right\}, \hfill (24.20)
\cr
\hspace*{1cm}
\min_{G_1 \in T_1}r( \, A - AG_1A \,) = r(A) + r(M) + r\left[
\begin{array}{cc} 0  &  B \\ C & D \end{array} \right] - r \left[ \begin{array}{ccc} A & 0 & B \\ 0 & C & D \end{array} \right] -
r \left[ \begin{array}{cc} A & 0  \\ 0 & B \\  C & D \end{array}
\right].  \hfill  (24.21)
\cr}
$$}
{\bf Proof.}\, Let $ P = [\, I_n, \ 0 \,]$ and
$Q = [\, I_m, \ 0 \,]^T$. Then  according to (22.1) and (22.2),
we find that
\begin{eqnarray*}
\max_{G_1 \in T_1}r( \, A - AG_1A \,) &= & \max_{M^-}r( \, A - APM^-QA \,) \\
& = & \min \left\{ \,  r(AP), \ \ r(QA),  \ \ r\left[\begin{array}{cc} M
 & QA \\ AP & A \end{array} \right] - r(M) \, \right\} \\
 & = & \min \left\{ \  r(A), \ \ \  r(\, M - QAP \,) + r(A) - r(M) \,
 \right\} \\
& = & \min \left\{ \,  r(A), \ \ \  r\left[\begin{array}{cc} 0  &  B \\ C & D \end{array} \right] 
+ r(A) - r(M) \, \right\}, 
\end{eqnarray*}
\begin{eqnarray*}
\min_{G_1 \in T_1}r( \, A - AG_1A \,) &= & \min_{M^-}r( \, A - APM^-QA \,) \\
& = & r(M) -  r[\, M, \ QA \,] -
r\left[ \begin{array}{c} M \\ AP \end{array} \right] +
r\left[ \begin{array}{cc} M & QA \\ AP & A  \end{array} \right]\\
& = & r(A) + r(M) + r(\, M - QAP \,)  -  r[\, M, \ QA \,] -
r\left[ \begin{array}{c} M \\ AP \end{array} \right] \\
& = &  r(A) + r(M) + r\left[ \begin{array}{cc} 0  &  B \\ C & D \end{array}
\right] - r \left[ \begin{array}{ccc} A & 0 & B \\ 0 & C & D \end{array}
\right] - r \left[ \begin{array}{cc} A & 0  \\ 0 & B \\  C & D \end{array}
\right],
\end{eqnarray*}  
establishing (24.20) and (24.21). \qquad  $ \Box$  

\medskip

A similar result to (24.21) was presented in (21.104).  

\medskip

\noindent {\bf Corollary 24.14.}\, {\em Let $ M $ and $ M^-$ be given by
{\rm (24.15)} and {\rm (24.16)}.  Then

{\rm (a)}\,  $M$ has a g-inverse  with the form  $ M^- =
\left[ \begin{array}{cc} A^-  & G_2 \\ G_3 & G_4 \end{array} \right] $
if and only if
$$\displaylines{
\hspace*{2cm}
 r \left[ \begin{array}{cc} A & 0  \\ 0 & B \\  C & D \end{array} \right] + r \left[ \begin{array}{ccc} A & 0 & B \\ 0 & C & D \end{array} \right] = r(A) + r(M) + r\left[ \begin{array}{cc} 0  &  B \\ C & D \end{array} \right]. \hfill
\cr}
$$

{\rm (b)}\, $ T_1 \subseteq \{A^-\},$ i.e.$,$ any $ G_1$ in $T_1$ is a g-inverse of $ A$ if and only if 
 $ r(M) = r(A) + r \left[ \begin{array}{cc} 0 & B \\  C & D \end{array}
 \right].$  }

 \medskip

\noindent {\bf Proof.}\, Follows immediately from Theorem 24.13.
 \qquad  $ \Box$

 \medskip

\noindent {\bf Corollary 24.15.}\, {\em Let $ M $ and $ M^-$ be given by
 {\rm (24.15)} and {\rm (24.16)}, and $ T_1$---$T_4$ are given by
 {\rm (24.17)}.
Then
$$\displaylines{
\hspace*{2cm}
 T_1 \subseteq \{A^-\},  \qquad T_2 \subseteq \{C^-\}, \qquad T_3 \subseteq \{B^-\},   \ \ \ \ \ 
T_4 \subseteq \{D^-\} \hfill (24.22)
\cr}
$$
are all satisfied if and only if
$$ \displaylines{
\hspace*{2cm}
r(M) = r(A) + r(B) + r(C) + r(D). \hfill (24.23)
\cr}
$$
}
{\bf Proof.}\, If (24.22) holds, then 
\begin{eqnarray*}
r(M) = r(MM^-) = tr(MM^-)  & = & tr \left[ \begin{array}{cc} AG_1 + BG_3 &   AG_2 + BG_4
 \\ CG_1 + DG_3  & CG_2 + DG_4 \end{array} \right] \\
& = & tr(AG_1) + tr (BG_3) + tr(CG_2) + tr(DG_4) \\
& = & tr(AA^-) + tr (BB^-) + tr(CC^-) + tr(DD^-) \\
& = &r(A) + r(B) + r(C) + r(D).
\end{eqnarray*}
Conversely, if (24.23) is satisfied, then (24.22) naturally holds by
Corollary 24.14(b). \qquad  $ \Box$

\medskip

When $ D = 0$ in the above theorems and corollaries, the corresponding
results can further simplify. We leave them to the reader.

\markboth{YONGGE  TIAN }
{25. EXTREME RANKS OF  $X_0 $ AND $X_1 $ TO $B(\, X_0 + iX_1 \,)C = A$} 

\chapter{Extreme ranks of $ X_0$ and $X_1$ in solutions to $B(\, X_0 + iX_1 \,)C = A$}

\noindent Suppose $ BXC = A $ is a complex matrix equation. Then it can be written as 
$$ \displaylines{
\hspace*{2cm}
( \, B_0 + iB_1 \,) ( \, X_0 + iX_1 \, )( \, C_0 + iC_1 \,) = (\, A_0 + iA_1 
\,), \hfill (25.1) 
\cr}
$$ 
where $ A_0, \,  A_1 \in {\cal R}^{m \times n}, \,B_0, \,  B_1 \in
{\cal R}^{m \times k}, \,C_0, \,  C_1 \in {\cal R}^{l \times n},$  and
 $ X_0, \,  X_1 \in {\cal R}^{k \times l}.$
In this chapter we determine maximal and minimal ranks of two real matrices
$ X_0$ and $X_1$ in solutions to the complex matrix equation in (25.1), 
and then present some consequences. To do so, we need the following result.

\medskip

\noindent {\bf  Lemma 25.1.}\, {\em The complex matrix equation
 {\rm (25.1)} is consistent if and only of the following real matrix
equation
$$\displaylines{
\hspace*{2cm}
 \left[ \begin{array}{cr}  B_0   & -B_1  \\ B_1 & B_0 \end{array} \right]
\left[ \begin{array}{cc}  Y_1   & Y_2  \\ Y_3 & Y_4 \end{array} \right]
\left[ \begin{array}{cr}  C_0   & -C_1  \\ C_1 & C_0 \end{array} \right]
= \left[ \begin{array}{cr}  A_0   & -A_1  \\ A_1 & A_0 \end{array} \right],
\hfill (25.2)
\cr}
$$ 
is consistent over the real number field ${\cal R}$. In that case the
general solution of {\rm (25.1)} can be written as
$$\displaylines{
\hspace*{2cm} 
X  = X_0 + iX_1 = \frac{1}{2} (\, Y_1 + Y_4 \, ) +  \frac{i}{2}
(\, Y_3 - Y_2 \, ), \hfill (25.3)
\cr}
$$ 
where $ Y_1$---$Y_4$ are the general solutions of {\rm (25.2)} over
${\cal R}.$  Written in an explicit form$,$ $ X_0$ and $ X_1$ in
 {\rm (25.3)}  are
$$
\displaylines{
\hspace*{2cm}
X_0 =  \frac{1}{2}P_1\phi^-(B)\phi(A) \phi^-(C)Q_1 +
\frac{1}{2}P_2\phi^-(B)\phi(A) \phi^-(C)Q_2 \hfill
\cr
\hspace*{4cm}
  + \ [\, P_1F_{\phi(B)}, \  P_2F_{\phi(B)} 
 \,] \left[ \begin{array}{c} V_1 \\ V_2  \end{array} \right]  +
 [\, W_1, \ W_2 \,] \left[ \begin{array}{c} E_{\phi(C)}Q_1 \\
 E_{\phi(C)}Q_2  \end{array} \right],  \hfill
\cr
\hspace*{2cm}
X_1   =  \frac{1}{2}P_2\phi^-(B)\phi(A) \phi^-(C)Q_1 - \frac{1}{2}P_1
\phi^-(B)\phi(A) \phi^-(C)Q_2  \hfill
\cr
\hspace*{4cm}
 + \  [\, P_2F_{\phi(B)}, \  -P_1F_{\phi(B)} \,] \left[ \begin{array}{c}
 V_1 \\ V_2  \end{array} \right]  + [\, W_1, \ W_2 \,]
 \left[ \begin{array}{c} -E_{\phi(C)}Q_2 \\
 E_{\phi(C)}Q_1  \end{array} \right], \hfill
\cr}
$$
where $ \phi(M) = \phi(M_0 + iM_1) =\left[ \begin{array}{cr}  M_0   & -M_1
\\ M_1 & M_0 \end{array} \right],$ $ P_1 = [\, I_k, \ 0 \,], $ $ P_2 = [\, 0
, \ I_k \,],$ $ Q_1 = \left[ \begin{array}{c} I_l \\ 0  \end{array}
\right],$ $ Q_2 =
 \left[ \begin{array}{c} 0 \\ I_l  \end{array} \right],$ $ V_1, \, V_2, \,
 W_1 $ and $ W_2$ are arbitrary over ${\cal R}$.  }

\medskip
   
\noindent {\bf Proof.}\, It is well known that for any $ M = M_0 +iM_1 \in
{\cal C}^{ m \times n}$, there  is
$$ 
 \frac{1}{2} \left[ \begin{array}{cr}  I_m   &  iI_m  \\ -iI_m & -I_m \end{array} \right]\left[ \begin{array}{cc}  M_0 +iM_1  &  0\\ 0 & M_0  - iM_1
 \end{array} \right] \left[ \begin{array}{cr}  I_m   &  iI_m  \\ -iI_m
  & -I_m \end{array} \right] = \left[ \begin{array}{cr}  M_0   & -M_1
  \\ M_1 & M_0 \end{array} \right] = \phi(M), \eqno (25.4)
$$ 
where $\phi(\cdot)$ satisfies the following operation properties 
 
(i)\, $ M  = N  \Leftrightarrow  \phi(M) = \phi(N).$ 

(ii)\,  $ \phi( \, M + N \, ) = \phi(M) + \phi(N), \ \
\phi( MN) = \phi(M)\phi(N), \ \ \    \phi( kM) = k\phi(M),
 \ k \in {\cal R}. $   

(iii)\,  $ \phi(M) = K_{2m}\phi(M)K_{2n}^{-1}$, where $ K_{2t} = \left[ \begin{array}{cr} 0 &  I_t  \\ 
 -I_t & 0 \end{array} \right], \ \ t = m, \ n.$   

(iv)\,  $r[ \phi(M)] = 2r(M).$ \\
Suppose now that (25.1) has a solution $ X $ over ${\cal C}.$ Applying
the above properties (i) and  (ii) to it yields
$$ \displaylines{
\hspace*{2cm}
\phi(B)\phi(X)\phi(C) = \phi(A), \hfill (25.5) 
\cr}
$$ 
which shows that $\phi(X)$ is a solution to (25.2). Conversely suppose
that (25.2) has a solution
 $ \widehat{Y} = \left[ \begin{array}{cc}  Y_1   & Y_2  \\ Y_3 & Y_4
 \end{array} \right]
\in {\cal R}^{2k \times 2l}$, i.e.,  $\phi(B) \widehat{Y} \phi(C) =
\phi(A)$. Then applying the above property (iii) to it yields
$$\displaylines{
\hspace*{1.5cm} 
 K_{2m}\phi(B)K_{2k}^{-1}\widehat{Y} K_{2l}\phi(C)K_{2n}^{-1} =
 K_{2m} \phi(A)K_{2n}^{-1}, \hfill 
 \cr}
$$ 
consequently
$$ \displaylines{
\hspace*{1.5cm}
\phi(B)( \, K_{2k}^{-1} \widehat{Y} K_{2l} \,)\phi(C) = \phi(A), \hfill
\cr} 
$$ 
which shows that $ K_{2k}^{-1}\widehat{Y} K_{2l}$ is a solution of
 (25.5), too. Thus
$ \frac{1}{2} ( \, \widehat{Y} + K_{2k}\widehat{Y} K_{2l} \, )$ is a
solution of (25.2), and this solution  has the form
$$ \displaylines{
\hspace*{1.5cm}
\frac{1}{2}(\,  \widehat{Y} + \frac{1}{2} K_{2k}^{-1}\widehat{Y} K_{2l} \,) = 
\frac{1}{2}\left[ \begin{array}{cc}  Y_1 & Y_2  \\ Y_3 & Y_4 \end{array} \right] 
+ \frac{1}{2}\left[ \begin{array}{cr} Y_4 & -Y_3  \\ -Y_2 & Y_1 \end{array}
 \right] = \frac{1}{2}\left[ \begin{array}{cc}  Y_1 + Y_4  &  -( \,Y_3 - Y_2 \, ) \\ Y_3 - Y_2  &  Y_1 + Y_4 \end{array} \right].  \hfill
\cr}
$$ 
Let $ \widehat{X} =  \frac{1}{2} (\, Y_1 + Y_4  \,) + \frac{i}{2}
( \,Y_3 - Y_2 \, )$. Then
$ \phi( \widehat{X}) = \frac{1}{2} ( \, \widehat{Y} + K_{2k}^{-1}\widehat{Y}
 K_{2l} \, )$ is a solution of (25.5).
Thus by the above property (i), we know that $ \widehat{X}$ is a solution
of (25.1). The above derivation shows that
the two equations (25.1) and (25.2) have the same consistency condition
and their solutions satisfy the equality (25.3). Observe that
$ Y_1$ ---$Y_4$ in (25.2) can be written as
$$ \displaylines{
\hspace*{1.5cm}
Y_1 = P_1YQ_1, \ \ \  \ Y_2 = P_1YQ_2, \ \ \ \  Y_3 = P_2YQ_1, \ \ \ \
Y_4 = P_2YQ_2, \hfill
\cr}
$$
where $ Y = \left[ \begin{array}{cc}  Y_1   & Y_2  \\ Y_3 & Y_4 \end{array}
\right]$, and the general solution of (25.2) can be written as
$$ \displaylines{
\hspace*{1.5cm}
Y = \phi^-(B)\phi(A) \phi^-(C) + 2F_{\phi(B)}[\, V_1, \ V_2 \,] +
2 \left[ \begin{array}{c} W_1 \\ W_2  \end{array} \right]E_{\phi(C)}. \hfill
\cr
Hence \hfill
\cr
\hspace*{1.5cm}
Y_1 = P_1YQ_1 =P_1\phi^-(B)\phi(A) \phi^-(C)Q_1 + 2P_1F_{\phi(B)}V_1 +
2 W_1E_{\phi(C)}Q_1, \hfill
\cr
\hspace*{1.5cm}
Y_2 = P_1YQ_2 =P_1\phi^-(B)\phi(A) \phi^-(C)Q_2 + 2P_1F_{\phi(B)}V_2 +
2W_1E_{\phi(C)}Q_2, \hfill
\cr
\hspace*{1.5cm}
Y_3 = P_2YQ_1 =P_2\phi^-(B)\phi(A) \phi^-(C)Q_1 + 2P_2F_{\phi(B)}V_1 +
2 W_2E_{\phi(C)}Q_1, \hfill
\cr
\hspace*{1.5cm} 
Y_4 = P_2YQ_2 =P_2\phi^-(B)\phi(A) \phi^-(C)Q_2 + 2P_2F_{\phi(B)}V_2 +
2W_2E_{\phi(C)}Q_2. \hfill
\cr}
$$  
Putting them in (25.3) yields the general expressions of the
two real matrices $ X_0$ and $ X_1$.  \qquad $\Box$

\medskip

\noindent {\bf Theorem 25.2.}\, {\em  Suppose the matrix equation
 {\rm (25.1)} is consistent$,$ and denote
$$  \displaylines{
\hspace*{1.5cm}
S_0 = \{ \, X_0 \in {\cal R}^{k \times l} \ | \ B(X_0 + i X_1)C = A \, \}, \ \ \ 
S_1 = \{ \, X_1 \in {\cal R}^{k \times l} \ | \ B(X_0 + i X_1)C = A\, \}.
\hfill (25.6)
\cr}
$$

Then

{\rm (a)}\, The maximal and the minimal ranks of $X_0$ are given by 
$$
\displaylines{
\hspace*{1.5cm}
\max_{X_0 \in S_0}r(X_0) = \min  \left\{ \, k, \ \  l, \ \ k + l + r \left[
\begin{array}{rrc} A_0  & -A_1  & B_0 \\ A_1 & A_0 & B_1 \\
 C_0 & -C_1 & 0 \end{array} \right] - 2r(B) - 2r(C) \, \right\}, \hfill (25.7)
\cr
\hspace*{1.5cm}
\min_{X_0 \in S_0}r(X_0) =  r \left[ \begin{array}{rrc} A_0  & -A_1  & B_0
\\ A_1 & A_0 & B_1 \\ C_0 & -C_1 & 0  \end{array} \right] -
r \left[ \begin{array}{c}  B_0 \\
 B_1  \end{array} \right] - r[\, C_0, \ C_1\, ]. \hfill (25.8)
 \cr}
$$

{\rm (b)}\, The maximal and the minimal ranks of $ X_1$ are given by  
$$
\displaylines{
\hspace*{1.5cm}
\max_{X_1 \in S_1}r(X_1) = \min  \left\{ \, k, \ \  l, \ \ k + l +
r \left[ \begin{array}{rrc} A_0  & -A_1
 & B_0 \\ A_1 & A_0 & B_1 \\
 C_0 & C_1 & 0  \end{array} \right]  - 2r(B) - 2r(C) \, \right\}, \hfill (25.9)
\cr
\hspace*{1.5cm}
\min_{X_1 \in S_1}r(X_1)   =  r \left[ \begin{array}{rrc} A_0  & -A_1
 & B_0 \\ A_1 & A_0 & B_1 \\
 C_0 & C_1 & 0  \end{array} \right] - r \left[ \begin{array}{c}  B_0 \\ B_1
 \end{array} \right] - r[\, C_0, \ C_1\,]. \hfill (25.10)
\cr}
$$ }
{\bf Proof.}\, Applying (19.14) and (19.15) to $ X_0$ in
 (25.3) yields
$$ \displaylines{
\hspace*{2cm}
\max_{X_0 \in S_0}r(X_0) = \min \{ \, k, \ \ \   l,  \ \ \  r(M) \,  \},
\hfill (25.11)
\cr
\hspace*{2cm}
\min_{X_0 \in S_0}r(X_0) = r(M) - r[\, P_1F_{\phi(B)}, \  P_2F_{\phi(B)}
\, ] - r\left[ \begin{array}{c}
 E_{\phi(C)}Q_1 \\  E_{\phi(C)}Q_2  \end{array} \right],  \hfill (25.12)
\cr
\hspace*{0cm} 
where \hfill
\cr
\hspace*{2cm} 
M = \left[ \begin{array}{ccc} \frac{1}{2}P_1\phi^-(B)\phi(A) \phi^-(C)Q_1 +
\frac{1}{2}P_2\phi^-(B)\phi(A) \phi^-(C)Q_2 &  P_1F_{\phi(B)} & P_2F_{\phi(B)
}  \\ E_{\phi(C)}Q_1  & 0 & 0 \\  E_{\phi(C)}Q_2 & 0 & 0 \end{array}
\right]. \hfill
\cr}
$$ 
Note that $ \phi(B) \phi^-(B)\phi(C) = \phi(A)$ and $ \phi(A) \phi^-(C)
\phi(C) =\phi(A)$. By (1.3), (1.4) and (1.5), it is not difficult but tedious to find
that
\begin{eqnarray*} 
\lefteqn{r(M)} \\
& = & r\left[ \begin{array}{ccccc} \frac{1}{2}P_1\phi^-(B)\phi(A)
\phi^-(C)Q_1 +  \frac{1}{2}P_2\phi^-(B)\phi(A) \phi^-(C)Q_2 &  P_1 & P_2
& 0 & 0 \\ Q_1 & 0 & 0 & \phi(C) & 0 \\ Q_2 & 0 & 0 & 0 & \phi(C)
 \\  0 & \phi(B) & 0 & 0 & 0 \\ 0 & 0 & \phi(B) & 0 & 0 \end{array} \right] \\
  & & \ \ - 2r[\phi(B)] - 2r[\phi(C)] \\
& = & r\left[ \begin{array}{ccccc} 0 &  P_1 & P_2  & 0 & 0 \\ Q_1 & 0 & 0 &
\phi(C) & 0 \\ Q_2 & 0 & 0 & 0 & \phi(C)  \\  0 & \phi(B) & 0 & 0 & 0 \\
0 & 0 & \phi(B) & 0 & \phi(A) \end{array} \right] - 2r[\phi(B)] -
2r[\phi(C)] \\
& = & r \left[ \begin{array}{rrr} 0  & C_0  & -C_1 \\ B_0 & A_0 & -A_1 \\
 C_1 & A_1 & A_0  \end{array} \right] - r[\phi(B)]  - r[\phi(C)] + k + l \\ 
& = &  r \left[ \begin{array}{rrc} A_0  & -A_1  & B_0 \\ A_1 & A_0 & B_1 \\
 C_0 & -C_1 & 0  \end{array} \right] - 2r(B) - 2r(C) + k  + l, 
\end{eqnarray*} 
\begin{eqnarray*} 
r[\, P_1F_{\phi(B)}, \  P_2F_{\phi(B)} \, ] & = & r\left[ \begin{array}{cc}
P_1 & P_2 \\ \phi(B) & 0
\\ 0 & \phi(B)  \end{array} \right] - 2r[\phi(B)]  \\ 
 & = & r\left[ \begin{array}{ccr} -A_1 & 0 & -A_0 \\ A_0 & 0  & -A_1  \\  0 & A_0 & -A_1 \\ 
0 & A_1 & A_0   \end{array} \right] - 2r[\phi(B)] + k \\
& = & r \left[ \begin{array}{c}  B_0 \\ B_1  \end{array} \right] -
r[\phi(B)] + k  = r \left[ \begin{array}{c}  B_0 \\ B_1  \end{array}
\right] - 2r(B) + k,
\end{eqnarray*} 
\begin{eqnarray*} 
r \left[ \begin{array}{c} E_{\phi(C)}Q_1 \\ E_{\phi(C)}Q_2 \end{array}
\right] & = & r \left[ \begin{array}{ccc} \phi(C) & 0 & Q_1 \\ 0 & \phi(C)
& Q_2 \end{array} \right] - 2r[\phi(C)] \\
& = & r \left[ \begin{array}{cccc} C_1 & C_0 & 0 & 0  \\  
0 & 0  & C_0 & -C_1 \\ -C_0 & C_1 & C_1 & C_0 \\  \end{array} \right] -
2r[\phi(C)] + l \\
& = & r[\, C_0, \ C_1 \,] - r[\phi(C)] + l = r[\, C_0, \ C_1 \,] - 2r(C)
+ l.
\end{eqnarray*}  
Putting them in (25.10) and (25.12) yields (25.7) and (25.8). Similarly
 we can establish (25.9) and (25.10).  \qquad $ \Box$

Below is a direct consequence of Theorem 25.2. 

\medskip

\noindent {\bf Corollary 25.3.}\, {\em  Suppose the matrix equation {\rm (25.1)}
 is consistent. Then

{\rm (a)}\, Eq.\,{\rm (25.1)} has a real solution $ X \in {\cal R}^{k \times l}$ if
and only if
$$ \displaylines{
\hspace*{2cm}
r \left[ \begin{array}{crc} A_0  & -A_1  & B_0 \\ A_1 & A_0 & B_1 \\
 C_0 & C_1 & 0  \end{array} \right] = r \left[ \begin{array}{c}  B_0 \\ B_1
  \end{array} \right] + r[\, C_0, \ C_1\,]. \hfill (25.13)
\cr}
$$ 
 
 {\rm (b)}\, All the solutions of {\rm (25.1)} are real if and only if   
$$\displaylines{
\hspace*{2cm} 
r \left[ \begin{array}{rrc} A_0  & -A_1  & B_0 \\ A_1 & A_0 & B_1 \\
 C_0 & C_1 & 0  \end{array} \right] =  2r(B) +2r(C) - k - l. \hfill (25.14)
 \cr}
$$

{\rm (c)}\, Eq.\,{\rm (25.1)} has a pure imaginary solution $ X = iX_1,$
where $ X_1 \in {\cal R}^{k \times l}, $ if and only if
$$ \displaylines{
\hspace*{2cm}
r \left[ \begin{array}{rrc} A_0  & -A_1  & B_0 \\ A_1 & A_0 & B_1 \\
 C_0 & -C_1 & 0 \end{array} \right] = r \left[ \begin{array}{c}  B_0 \\ B_1
  \end{array} \right] + r[\, C_0, \ C_1  \, ]. \hfill (25.15)
  \cr}
$$
 
 {\rm (d)}\, All the solutions of {\rm (25.1)} are  pure imaginary if and only if
$$ \displaylines{
\hspace*{2cm}
r \left[ \begin{array}{rrc} A_0  & -A_1  & B_0 \\ A_1 & A_0 & B_1 \\ C_0 &
-C_1 & 0  \end{array} \right]  = 2r(B) +2r(C) - k - l. \hfill (25.16)
\cr}
$$}

We next consider extreme ranks of $ A_0 - B_0X_0C_0 $  and  $ A_1 - B_1X_1C_1 $ with respect to  the real matrices 
$ X_0$ and $ X_1$ in  solution of $ ( \, B_0 + iB_1 \, ) ( \, X_0 + iX_1 \, )( \, C_0 + iC_1 \, ) =  A_0 + iA_1$, 
the corresponding results can  be used  in the next section to determine the relationships of  $A$ and $C$, $B$ and 
$D$ in generalized inverse $ (\, A + iB \, )^- = C + iD.$ 

\medskip

\noindent {\bf Theorem 25.4.}\, {\em  Suppose the matrix equation {\rm (25.1)} is consistent$,$ and $ S_0$ is defined by 
{\rm (25.6)}. Then 

{\rm (a)}\, The maximal rank of $ A_0 - B_0X_0C_0 $ is  
\begin{eqnarray*}  
\lefteqn{ \max_{X_0 \in S_0}r(\, A_0 - B_0X_0C_0 \,) } \\
& = & \min  \left\{ \ r[ \, A_0,  \ B_0 \, ], \ \  
 r \left[ \begin{array}{c}  A_0 \\ C_0  \end{array} \right],  \ \ r \left[ \begin{array}{ccccc}
 -A_0 & B_0 & 0 & 0 & 0 \\  C_0 & 0 & 0 & C_0 & -C_1 \\  0 & 0 & 0 & C_1 & C_0  \\ 0 & B_0 & -B_1 & A_0 & -A_1  \\ 0 &
 B_1 & B_0 & A_1 & A_0  \end{array} \right] -2r(B)- 2r(C) 
\ \right\}.    
\end{eqnarray*} 

{\rm (b)}\, The minimal  rank of $ A_0 - B_0X_0C_0 $ is  
\begin{eqnarray*} 
\min_{X_0 \in S_0}r(\, A_0 - B_0X_0C_0 \,) & = & r[ \, A_0, \ B_0 \, ] + r \left[ \begin{array}{c} A_0 \\ 
C_0  \end{array} \right]  + r \left[ \begin{array}{ccccc}  -A_0 & B_0 & 0 & 0 & 0 \\  C_0 & 0 & 0 & C_0 & -C_1 \\ 0 & 
0 & 0 & C_1 & C_0  \\  0 & B_0 & -B_1 & A_0 & -A_1  \\ 0 & B_1 & B_0 & A_1 & A_0  \end{array} \right] \\
 & & \  - \ r \left[ \begin{array}{ccc} A_0 & B_0 & 0 \\  C_0 & 0 & 0  \\ 0 & B_0 & -B_1  \\ 0 & B_1 & B_0 \end{array} \right] - r \left[ \begin{array}{cccc}  C_0 & A_0 & 0 & 0 \\ B_0  & 0 & B_0 & -B_1 \\
 0 & 0 & B_1 & B_0   \end{array} \right]. 
\end{eqnarray*} }

{\bf Proof.}\, Putting the general expression $ X_0$ in (25.3) in  $A_0 - B_0X_0C_0$ yields 
\begin{eqnarray*}  
A_0 - B_0X_0C_0 & = & A_0 - \frac{1}{2}B_0P_1\phi^-(B)\phi(A) \phi^-(C)Q_1C_0 -  \frac{1}{2}B_0P_2\phi^-(B)\phi(A) 
\phi^-(C)Q_2B_0 \\
& &   - \ B_0[\, P_1F_{\phi(B)}, \  P_2F_{\phi(B)}  \,] \left[ \begin{array}{c} V_1 \\ V_2  \end{array} \right]C_0
  - B_0[\, W_1, \ W_2 \,] \left[ \begin{array}{c} E_{\phi(C)}Q_1 \\  E_{\phi(C)}Q_2  \end{array} \right]C_0 \\
 &= & N -  B_0GVC_0 -  B_0 WHC_0. 
\end{eqnarray*}
Then according to (19.3)  and (19.4), we get 
\begin{eqnarray*} 
\max_{X_0 \in S_0}r(\, A_0 - B_0X_0C_0 \,)  & = & \max_{V, \,W} r( \, N -  B_0GVC_0 -  B_0 WHC_0 \, ) \\
 & = & \min  \left\{ \  r[ \, N,  \ B_0 \, ], \ \  r \left[ \begin{array}{c} N \\ C_0  \end{array} \right], \ \ \ 
r \left[ \begin{array}{cc} N & B_0G \\ HC_0 & 0 \end{array} \right]  \ \right\}, 
\end{eqnarray*} 
\begin{eqnarray*}  
\lefteqn{\min_{X_0 \in S_0}r(\, A_0 - B_0X_0C_0 \,)} \\
 & = & \min_{V, \,W} r( \, N -  B_0GVC_0 -  B_0 WHC_0 \, ) \\
& = & r[\, N,  \ B_0 \, ] +  r \left[ \begin{array}{c} N \\ C_0  \end{array} \right] + 
r \left[ \begin{array}{cc} N & B_0G \\ HC_0  & 0  \end{array} \right] - \ r \left[ \begin{array}{cc} N & B_0G \\ C_0  & 0 \end{array} \right] -  r \left[ \begin{array}{cc}
 N & B_0 \\ HC_0  & 0 \end{array} \right].     
\end{eqnarray*} 
Simplyfying the rank equalities by Lemma 1.1 may eventually results in the two equalities in parts (a) and (b). 
The processes, however, are quite tedious and are therefore omitted them here.  \qquad $ \Box$  

\medskip

\noindent {\bf Theorem 25.5.}\, {\em Suppose the matrix equation {\rm (25.1)} is consistent, and $ S_1$ is defined by
{\rm (25.6)}. Then 

{\rm (a)}\, The maximal rank of $ A_1 + B_1X_1C_1 $ is  
\begin{eqnarray*}  
\lefteqn{\max_{X_1 \in S_1}r(\, A_1 + B_1X_1C_1 \,) } \\
 &= & \min  \left\{ \ r[ \, A_1,  \ B_1 \, ], \ \  
 r \left[ \begin{array}{c}  A_1 \\ C_1  \end{array} \right],   \ \ r \left[ \begin{array}{ccccc}
 A_1 & B_1 & 0 & 0 & 0  \\  C_1 & 0  & 0 & C_1 & C_0 \\ 0 & 0 & 0 & -C_0 & C_1  \\ 0 & B_1 & B_0 & -A_1 & -A_0  \\ 0 & 
-B_0 & B_1 & A_0 & -A_1  \end{array} \right] -2r(B)- 2r(C) 
\ \right\}.    
\end{eqnarray*}  

{\rm (b)}\, The minimal rank of $ A_1 + B_1X_1C_1 $ is 
\begin{eqnarray*} 
\min_{X_1 \in S_1}r(\, A_1 + B_1X_1C_1 \,) & = & r[ \, A_1,  \ B_1 \,] + r \left[ \begin{array}{c} A_1 \\ C_1  \end{array} \right]  +  r \left[ \begin{array}{ccccc}
 A_1 & B_1 & 0 & 0 & 0  \\  C_1 & 0  & 0 & C_1 & C_0 \\ 0 & 0 & 0 & -C_0 & C_1  \\ 
 0 & B_1 & B_0 & -A_1 & -A_0  \\ 0 & -B_0 & B_1 & A_0 & -A_1  \end{array} \right] \\
 & & \ \  - r \left[ \begin{array}{ccc} A_1 & B_1 & 0 \\  C_1 & 0 & 0  \\ 0 & B_1 & B_0  \\ 0 & -B_0 & 
B_1   \end{array} \right] - r \left[ \begin{array}{cccc}  A_1 & B_1 & 0 & 0 \\ C_1  & 0 & C_1 & C_0 \\
  0 & 0 & -C_0 & C_1   \end{array} \right]. 
\end{eqnarray*}}

{\bf Proof.}\, Writing (25.1) in the following equivalent form 
$$
 ( \, A_1 - iA_0 \, ) ( \, X_1 - iX_0 \, )( \, B_1 - iB_0 \, ) =  -C_1 + iC_0,
$$ 
and then applying Theorem 2.4 to it  yields (a) and (b).  \qquad  $ \Box$  

\medskip

\noindent {\bf Corollary 25.6.}\, {\em  Suppose the matrix equation {\rm (25.1)} is consistent$,$  and the two sets
 $ S_0$ and  $ S_1$ are defined by {\rm (25.6)}. 

{\rm (a)}\, If $ B_0X_0 C_0 = A_0$ is consistent over ${\cal R},$ then  
\begin{eqnarray*} 
\lefteqn{ \min_{X_0 \in S_0}r(\, A_0 - B_0X_0C_0 \,) }  \\
& = & r \left[ \begin{array}{ccccc}  -A_0 & B_0 & 0 & 0 & 0  \\  C_0 & 0  & 0 & C_0 & -C_1 \\ 0 & 0 & 0 & C_1 & C_0  \\ 
0 & B_0 & -B_1 & A_0 & -A_1  \\ 0 & B_1 & B_0 & A_1 & A_0  \end{array} \right] -
 r \left[ \begin{array}{cc} B_0 & 0  \\ 0 & B_1  \\ B_1 & B_0  \end{array} \right] - 
r \left[ \begin{array}{ccc} C_0 & 
0  & C_1 \\  0 & C_1 & C_0   \end{array} \right]. 
\end{eqnarray*} 

{\rm (b)}\, If $ B_1X_1 C_1 = -A_1$ is consistent over ${\cal R},$  then  
\begin{eqnarray*} 
\lefteqn{ \min_{X_1 \in S_1}r(\, A_1 + B_1X_1C_1 \,) }  \\
& = &  r \left[ \begin{array}{ccccc}
 A_1 & B_1 & 0 & 0 & 0 \\  C_1 & 0  & 0 & C_1 & C_0 \\ 0 & 0 & 0 & -C_0 & C_1  \\ 
0 & B_1 & B_0 & -A_1 & -A_0  \\ 0 & -B_0 & B_1 & A_0 & -A_1  \end{array} \right] - 
r \left[ \begin{array}{cc} B_0 & 0  \\ 0 & B_1  \\ B_1 & B_0  \end{array} \right] - 
r \left[ \begin{array}{ccc} C_0 & 0 & C_1 \\ 0 & C_1 & C_0   \end{array} \right]. 
\end{eqnarray*} 

{\rm (c)}\, The two linear matrix equations 
$$ 
 ( \, B_0 + iB_1 \, ) ( \, X_0 + iX_1 \, )( \, C_0 + iC_1 \, ) =  A_0 + iA_1  \ \ \ and  \ \ \ 
 B_0X_0 C_0 = A_0
$$   
have a common solution for $ X_0 \in {\cal R}^{p \times q}$ if and only if 
$$ 
r \left[ \begin{array}{ccccc}  -A_0 & B_0 & 0 & 0 & 0  \\  C_0 & 0  & 0 & C_0 & -C_1 \\ 
 0 & 0 & 0 & C_1 & C_0  \\ 0 & B_0 & -B_1 & A_0 & -A_1  \\ 
0 & B_1 & B_0 & A_1 & A_0  \end{array} \right] =
 r \left[ \begin{array}{cc} B_0 & 0  \\ 0 & B_1  \\ B_1 & B_0  \end{array} \right] +  
r \left[ \begin{array}{ccc} C_0 & 
0 & C_1 \\ 0 & C_1 & C_0   \end{array} \right]. 
$$ 

{\rm (d)}\, The two linear matrix equations 
$$ 
 ( \, B_0 + iB_1 \, ) ( \, X_0 + iX_1 \, )( \, C_0 + iC_1 \, ) =  A_0 + iA_1  \ \ \ and  \ \ \ 
 B_1X_1 C_1 = -A_1
$$   
have a common solution for $ X_1 \in {\cal R}^{p \times q}$ if and only if 
$$ 
r \left[ \begin{array}{ccccc}
A_1 & B_1 & 0 & 0 & 0  \\ C_1 & 0  & 0 & C_1 & C_0 \\ 0 & 0 & 0 & -C_0 & C_1  \\ 0 & B_1 & B_0 & -A_1 & -A_0  \\ 0 & 
-B_0 & B_1 & A_0 & -A_1  \end{array} \right] =
r \left[ \begin{array}{cc} B_0 & 0  \\ 0 & B_1  \\ B_1 & B_0  \end{array} \right] + 
r \left[ \begin{array}{ccc} C_0 & 0 & C_1 \\ 0 & C_1 & C_0 \end{array} \right].
$$}
{\bf Proof.}\, Follows from simplification of Theorems 25.4 and 25.5  under the assumptions of the corollary. 
\qquad  $ \Box$  

\medskip

Applying Theorem 25.2 to the real matrices $ C$ and  $D$ in
in the inner inverse $(\, A + iB \, )^- = C + iD$, we get the following.

\medskip

\noindent {\bf Theorem 25.7.}\, {\em  Let $N = A + iB \in {\cal C}^{m \times n}$ be
given$,$  and denote
$$ \displaylines{
\hspace*{1cm}
T_1 = \{ \ C \in {\cal R}^{n \times m} \ | \  C + iD \in \{ N^- \} \ \},
\ \ \ T_2 = \{ \, D \in {\cal R}^{n \times m} \ | \  C + iD \in \{ N^- \} \
\}. \hfill (25.17)
\cr}
$$ 
  
{\rm (a)}\, The maximal and the minimal ranks of $ C $ in {\rm (25.17)} are 
$$ 
\displaylines{
\hspace*{2cm}
\max_{C \in T_1}r(C) = \min  \left\{ \ m , \ \   n,  \ \  m + n + r(A) -
r \left[ \begin{array}{cr}  A & -B  \\ B   & A \end{array} \right] \
\right\}, \hfill
\cr
\hspace*{2cm}
\min_{C \in T_1}r(C) = r \left[ \begin{array}{cr}  A  & -B  \\ B   & A \end{array} \right] -
  r \left[ \begin{array}{cc}  A  \\ B  \end{array} \right] - r[ \, A, \ B
  \, ] + r(A).  \hfill
\cr}
$$ 

{\rm (b)}\, The maximal and the minimal ranks of $ D $ in {\rm (25.17)} are 
$$\displaylines{
\hspace*{2cm}
\max_{D \in T_2}r(D) = \min  \left\{ \ m ,  \ \  n, \ \  m + n + r(B) -
r \left[ \begin{array}{cr}  A & -B  \\ B   & A \end{array} \right] \
\right\},  \hfill
\cr
\hspace*{2cm}
\min_{D \in T_2}r(D) = r \left[ \begin{array}{cr}  A  & -B  \\ B   & A
\end{array} \right] - r \left[ \begin{array}{cc}  A  \\ B  \end{array}
\right] - r[\, A, \ B \, ] + r(B). \hfill
\cr}
$$ 

{\rm (c)}\, If $ A = 0,$ then 
$$\displaylines{
\hspace*{2cm} 
\max_{C \in T_1}r(C) = \min  \left\{ \ m ,  \ \  n,  \ \  m + n - 2r(B) \
 \right\}. \hfill
\cr}
$$ 

{\rm (d)}\, If $ B = 0,$ then 
$$ \displaylines{
\hspace*{2cm}
\max_{D \in T_2}r(D) = \min \left\{ \, m , \ \ n, \ \  m + n - 2r(A) \,
\right\}. \hfill
\cr}
$$ 

{\rm (e)}\, If $ R(A) \subseteq  R(N)$ and  $ R(A^*) \subseteq  R(N^*),$
then
$$\displaylines{
 \hspace*{2cm}
\min_{C \in T_1}r(C) = r(A), \qquad  \min_{D \in T_2}r(D) = r(B). \hfill  
\cr}
$$ }
{\bf Proof.}\, Follows from replacing $A, \ B $ and $ C $ all by
$ N = A + iB$ in  Theorem 25.2. \qquad  $ \Box$

\medskip

\noindent {\bf Corollary 25.8.}\, {\em  Let $N = A + iB \in
{\cal C}^{m \times n}$ be given$,$  $T_1 $ and $ T_2$ be
defined  in {\rm (25.17)}. 
 
{\rm (a)}\,  $ N $ has a real generalized inverse if and only if 
$$ \displaylines{
\hspace*{2cm}
r \left[ \begin{array}{cr}  A  & -B  \\ B   & A \end{array} \right] =
r \left[ \begin{array}{cc}  A  \\ B  \end{array} \right] + r[\, A, \ B \,]
- r(B). \hfill
\cr}
$$ 

{\rm (b)} \, $ N $ has a pure imaginary generalized inverse if and only if 
$$ \displaylines{
\hspace*{2cm}
r \left[ \begin{array}{cr}  A  & -B  \\ B   & A \end{array} \right] =
r \left[ \begin{array}{cc}  A  \\ B  \end{array} \right] + r[ \, A, \ B \, ] - r(A). \hfill
\cr}
$$ }
{\bf Proof.}\, Follows directly from Theorem 25.4(a) and (b).
 \qquad $ \Box$

\medskip

\noindent {\bf Theorem 25.9.}\, {\em  Let $N = A + iB \in {\cal C}^{m \times n}$ be given,  $T_1 $ and $ T_2$ be 
defined  by {\rm (25.17)}. 
 
{\rm (a)}\, The maximal and  the minimal ranks of $ A - ACA $  are 
$$ 
\max_{C \in T_1}r(\,  A - ACA \,) = \min  \left\{ \ r(A), \ \ \ r(A) + r \left[ \begin{array}{cr} 0 & B  \\ B   & A \end{array} \right] -r \left[ \begin{array}{cr} A & -B  \\ B   & A \end{array} \right] \  \right\},  
$$ 
$$ 
\min_{C \in T_1}r( \,  A - ACA \, ) = r(A) + r \left[ \begin{array}{cc} 0 & B  \\ B   & A \end{array} \right]  + r \left[ \begin{array}{cr} A & -B  \\ B   & A \end{array} \right] -   r \left[ \begin{array}{ccc}  A  & 0  & B \\ 0 & B  & A  \end{array} \right] -  r \left[ \begin{array}{cc}  A  & 0 \\ 0 & B  \\ B & A  \end{array} \right]. 
$$

{\rm (b)}\, The maximal and the minimal ranks of $ B + BDB $ are 
$$ 
\max_{D \in T_2}r(\,  B + BDB \,) = \min \left\{ \ r(B), \ \ \ r(B) + r \left[ \begin{array}{cc} 0 & A  \\ A   & B \end{array} \right] - r \left[ \begin{array}{cr} B & -A  \\ A   & B \end{array} \right] \  \right\},  
$$ 
$$
 \min_{D \in T_2}r( \,  B + BDB  \, )=  r(B) + r \left[ \begin{array}{cc} 0 & A  \\ A   & B \end{array} \right] + r \left[ \begin{array}{cr} B & -A  \\ A   & B \end{array} \right] - r \left[ \begin{array}{ccc}  A  & 0 & B \\ 0  & B  & A  \end{array} \right] -  r \left[ \begin{array}{cc}  A  & 0 \\ 0 & B  \\ B & A  \end{array} \right]. 
$$}
{\bf Proof.}\, Follows from replacing $A, \ B $ and $ C $ all by $N = A + iB$ in Theorems 25.4 and 25.5. \qquad $ \Box$  

\medskip

\noindent {\bf Corollary 25.10.}\, {\em  Let $N = A + iB \in {\cal C}^{m \times n}$ be given$,$  $T_1 $ and $ T_2$ be 
defined by {\rm (25.17)}. 
 
{\rm (a)}\,  $ N $ has a generalized inverse with the form $N^- = A^- + iD$ if and only if 
$$ 
r \left[ \begin{array}{ccc} A & 0  & B \\ 0  & B  & A  \end{array} \right] + r \left[ \begin{array}{cc}  A  & 0 \\ 0 & B  \\ B & A  \end{array} \right] = r(A) + r \left[ \begin{array}{cc} 0 & B  \\ B  & A \end{array} \right] + r \left[ \begin{array}{cr} A & -B  \\ B   & A \end{array} \right].
$$ 

{\rm (b)}\, $T_1 \subseteq \{ A^- \}$, i.e., all the  generalized inverses of $ N $ have the form
 $N^- = A^- + iD$ if and only if 
$$ 
r \left[ \begin{array}{cr} A & -B  \\ B   & A \end{array} \right]  = r(A) +  r \left[ \begin{array}{cc} 0 & B  \\ B & A \end{array} \right].
$$ 

{\rm (c)} \ $ N $ has a generalized inverse  with the form $N^- = C - iB^-$ if and only if 
$$ 
r \left[ \begin{array}{ccc} A & 0  & B \\ 0 & B  & A \end{array} \right] + r \left[ \begin{array}{cc}  A  &  0 \\ 0 & B  \\ B & A  \end{array} \right] = r(B) + r \left[ \begin{array}{cc} 0 & A  \\ A & B \end{array} \right] + r \left[ \begin{array}{cr} A & -B  \\ B   & A \end{array} \right].
$$ 

{\rm (d)} \ $T_2 \subseteq \{ -B^- \}$, i.e., all the  generalized inverses of $ N $ have the form $N^- = C - iB^-$ if and only if 
$$ 
r \left[ \begin{array}{cr} A & -B  \\ B   & A \end{array} \right]  = r(B) + 
r \left[ \begin{array}{cc} 0 & A \\ A   & B \end{array} \right].
$$}
{\bf Proof.} \ Follows directly from Theorem 25.9.  \qquad $\Box$

\markboth{YONGGE  TIAN}
{26. RANKS AND INDEPENDENCE OF SOLUTIONS TO $B_1XC_1 + B_2YC_2 = A $} 

\chapter{Ranks and independence of solutions of the matrix equation $ B_1XC_1 + B_2YC_2 = A $}

We consider in the chapter possible  ranks of solutions $ X $ and 
$ Y $  of the matrix equation

Suppose 
$$
 B_1XC_1 + B_2YC_2= A,  \eqno (26.1)
$$
is a consistent linear matrix  over an arbitrary field ${\cal F}$, where $ B_1, \ C_1,  \ B_2, \ C_2$  and $A$ are $ m \times p,$ $q \times n$, $m \times s$, $t \times n $ and $m \times n$ matrices,
respectively. In this chapter, We consider the maximal and the minimal  ranks of solutions $ X $ and 
$ Y $  of (26.1), as well as independence of solutions $ X $ and $Y $ of (26.1).

As one of  basic linear matrix equations, (26.1) has been well examined
in matrix theory and its applications (see, e.g., \cite{BK1, ChP, He, Hu, Oz, Per, Ti0, XWZ}). 
Its solvability conditions and general
solutions for $ X $ and $ Y $ are completely  established by using ranks and generalized inverse of matrices. On the basis of those results and the rank formulas in the previous chapters, we now can give complete solutions to the above 
two problems.

The basic tools for investigating the above problems are the following
several known results on ranks and generalized inverses of matrices.

\medskip

\noindent {\bf  Lemma 26.1.}\, {\em Let $A \in {\cal F}^{m \times n}, \, B_1 \in
{\cal F}^{m \times k_1},     \, B_2 \in {\cal F}^{m \times k_2}, \, C_1 \in
{\cal F}^{l_1 \times n}$ and $C_2 \in {\cal F}^{l_2 \times n}$ be given. Then 
$$\displaylines{
\hspace*{0.1cm}
\max_{X, \, Y, \, Z} r( \, A - B_1XC_1 - B_2Y - ZC_2 \,)  =
\min \left\{  \  m , \ \ n , \ \
r \left[ \begin{array}{ccc} A & B_1 & B_2 \\ C_2  & 0 & 0 \end{array} \right], \ \  
r \left[ \begin{array}{cc} A & B_2  \\ C_1 & 0 \\ C_2 & 0 \end{array} \right]
 \ \right\}, \hfill (26.2)
\cr
\hspace*{0.1cm}
 \min_{X, \, Y, \, Z} r( \, A - B_1XC_1 - B_2Y - ZC_2 \,) =
 r \left[ \begin{array}{ccc} A & B_1 & B_2 \\ C_2  & 0 & 0 \end{array}
 \right] + r \left[ \begin{array}{cc} A & B_2  \\ C_1 & 0 \\ C_2 & 0
  \end{array} \right] \hfill
\cr
\hspace*{6cm}
  - r \left[ \begin{array}{ccc} A & B_1 & B_2 \\ C_1  & 0 & 0  \\  C_2  & 0 & 0 \end{array}
 \right] - r(B_2) - r(C_2).  \hfill (26.3)
\cr}
$$ }
\hspace*{0.3cm} This lemma can be simply derived from the rank formulas in Theorem 18.4 and Corollary 19.5, and its proof is omitted here.  

Concerning the general solution of (26.1), the following is well known.

\medskip

\noindent {\bf  Lemma 26.2.}\, {\em Suppose the matrix equation is given by
{\rm (26.1)}. Then

{\rm (a)\cite{Ti0}}\, The general solution of the homogeneous equation $ B_1XC_1 + B_2YC_2 = 0$ can factor as
$$ \displaylines{
\hspace*{2cm}
X = X_1X_2 + X_3, \qquad   Y = Y_1Y_2 + Y_3, \hfill
\cr}
$$ 
where $ X_1$---$X_3$ and $ Y_1$---$Y_3$ are the general solutions of the following four simple   homogeneous matrix 
equations
$$ \displaylines{
\hspace*{2cm}
B_1X_1 = - B_2Y_1, \qquad   X_2C_1 = Y_2C_2, \qquad  B_1X_3C_1 = 0, \qquad B_2Y_3C_2 =0, \hfill (26.4)
\cr}
$$ 
Solving these four equations and putting their general solutions in $ X $ and $ Y $ yields
$$\displaylines{
\hspace*{2cm}
X = S_1F_G UE_HT_1 + F_{B_1}V_1 + V_2E_{C_1},  \qquad    Y =  S_2F_G UE_HT_2 +  F_{B_2}W_1 + 
W_2E_{C_2}, \hfill
\cr}
$$
where $S_1 = [\, I_p, \ 0 \, ],$ $ S_2=[\, 0, \ I_s \, ],$ $ T_1 = \left[ \begin{array}{c} I_q \\ 0
\end{array} \right],$ $T_2 = \left[ \begin{array}{c} 0 \\ I_t  \end{array}
\right]$, $ G = [\, B_1, \ B_2 \, ]$ and $ H = \left[ \begin{array}{r} C_1 \\ - C_2
 \end{array} \right];$ the matrices $ U,$ $ V_1,$ $ V_2,$ $ W_1$ and $ W_2$ are arbitrary.

{\rm (b)\cite{Ti0}}\, Suppose the matrix equation {\rm (26.1)} is consistent. Then its general solution can 
factor as
$$ \displaylines{
\hspace*{2cm}
X = X_0  + X_1X_2 + X_3, \qquad   Y = Y_0 + Y_1Y_2 + Y_3, \hfill
\cr}
$$ 
where $X_0$ and  $ Y_0$ are a pair of particular solutions to {\rm (26.1)}$,$ $ X_1$---$X_3$ and 
$ Y_1$---$Y_3$ are the general solutions of the four simple matrix equations in {\rm (26.4)}. Written in an explicit form$,$ the general solution of { \rm (26.1)} is
$$
\displaylines{ 
\hspace*{2cm}
X =  X_0 + S_1F_G UE_HT_1 + F_{B_1}V_1 + V_2E_{C_1}, \hfill (26.5)
\cr
\hspace*{2cm}
 Y =  Y_0 + S_2F_G UE_HT_2 +  F_{B_2}W_1 + W_2E_{C_2}. \hfill (26.6)
 \cr}
$$
}
\hspace*{0.4cm} Various expressions of a pair of particular solutions of (26.1) can be found in
\cite{ChP}, \cite{He}, \cite{Hu}, \cite{Oz} and \cite{XWZ}. However we only use $X_0$ and $Y_0$ in form when 
determining possible ranks of solutions to (26.1),  we do not intend to present their explicit expressions  in 
(25.5) and (25.6). 

For convenience of representation, we adopt the following notation
$$ \displaylines{
\hspace*{0.5cm}
J_1 = \{ \, X \in {\cal F}^{p \times q} \ | \ B_1XC_1 + B_2YC_2= A  \, \}, \ \ \ \
J_2 = \{ \, Y \in {\cal F}^{s \times t} \ | \  B_1XC_1 + B_2YC_2= A \, \}. \hfill (26.7)
\cr}
$$
The two expressions in (26.5) and (26.6) clearly show that the general solution 
$X $ and $ Y$ of (26.1) 
are in fact two linear matrix expressions, each of which involves three independent variant
matrices.  In that case, apply Lemma 26.1 to obtain the following.

\medskip
   
\noindent {\bf Theorem 26.3.}\, {\em  Suppose that the matrix equation {\rm (26.1)} is
consistent$,$ and $ J_1$ and $J_2$ are defined  in {\rm (26.7)}. Then 

{\rm (a)}\, The maximal and  the minimal ranks of  solution $X$ of
{\rm (26.1)} are
$$
\max_{X \in J_1}r(X) = \min \left\{ p,  \  q,  \ p + q + r[\, A, \ B_2 \,] - r[\, B_1, \ B_2
 \,] - r(C_1),  \  p + q + r\left[ \begin{array}{c} A \\ C_2  \end{array} \right] - 
r\left[ \begin{array}{c} C_1 \\ C_2  \end{array} \right] -r(B_1) \right\}, \eqno (26.8)
$$
$$
\displaylines{
\hspace*{0cm}
\min_{X \in J_1}r(X) = r[\, A, \ B_2 \,] + r\left[ \begin{array}{c} A \\ C_2
\end{array} \right] -
r\left[ \begin{array}{cc} A & B_2 \\ C_2 & 0 \end{array} \right]. \hfill (26.9)
\cr}
$$

{\rm (b)}\, The maximal and  the minimal ranks of solution $Y$ of {\rm (26.1)}
are
$$
\max_{Y \in J_2}r(Y) = \min  \left\{ s,  \  t,  \  s + t + r[\, A, \ B_1 \,] - r[\, B_2, \ B_1 \,]
 - r(C_2),  \   s + t + r\left[ \begin{array}{c} A \\ C_1  \end{array} \right]
 - r\left[ \begin{array}{c} C_2 \\ C_1  \end{array} \right] -r(B_2) \right\},
 \eqno (26.10)
$$
$$
\displaylines{
\hspace*{0cm}
\min_{Y \in J_2}r(Y) = r[\, A, \ B_1 \,] + r\left[ \begin{array}{c} A \\ C_1
\end{array} \right] - r\left[ \begin{array}{cc} A &  B_1 \\ C_1 & 0 \end{array}
\right].  \hfill (26.11)
\cr}
$$}
{\bf Proof.}\, Applying (26.2) and (26.3) to  (26.5) yields
\begin{eqnarray*} 
\max_{X \in J_1}r(X) & = & \max_{U,\, V_1, \, V_2}r( \,  X_0 +
S_1F_G UE_HT_1  + F_{B_1}V_1 + V_2E_{C_1}\, ) \\
  & = & \min \left\{ \ p, \ \  q, \ \ r \left[ \begin{array}{ccc} X_0 & F_{B_1}
  & S_1F_G \\ E_{C_1}  & 0 & 0 \end{array} \right], \ \  r \left[
  \begin{array}{cc} X_0 & F_{B_1}  \\ E_{C_1} & 0 \\ E_HT_1 & 0 \end{array}
  \right] \ \right\},
\end{eqnarray*} 
$$
\displaylines{
\hspace*{0cm}
 \min_{X \in J_1}r(X) \hfill
\cr
\hspace*{0cm}  
= \min_{U,\, V_1, \, V_2}r( \,  X_0 + S_1F_G UE_HT_1  + F_{B_1}V_1 + V_2E_{C_1} \, ) \hfill
\cr
\hspace*{0cm}
=  r \left[ \begin{array}{ccc}  X_0 & F_{B_1}  & S_1F_G \\ E_{C_1}  & 0 & 0 \end{array} \right] + r \left[ \begin{array}{cc} X_0 & F_{B_1}  \\ E_{C_1} & 0 \\ E_HT_1 & 0
\end{array} \right] - r \left[ \begin{array}{ccc} X_0 & F_{B_1} & S_1F_G  \\ E_{C_1} & 0 & 0  \\ E_HT_1 & 0 & 0 
 \end{array} \right]  - r(F_{B_1})  - r(E_{C_1}). \hfill
\cr}
$$ 
By Lemma 1.1 and $ B_1X_0C_1+ B_2Y_0C_2 = A $,  we find that
$r(F_{B_1}) = p - r(B_1), \ r(E_{C_1}) = q - r(C_1)$,  and
$$
\displaylines{
\hspace*{0cm}
r \left[ \begin{array}{ccc} X_0 & F_{B_1}  & S_1F_G \\ E_{C_1}  & 0 & 0 \end{array}
\right] \hfill
\cr
\hspace*{0cm}
 =  r \left[ \begin{array}{cccc} X_0 & I_p & S_1 & 0  \\ I_q  & 0 & 0 & C_1
\\ 0 & B_1 & 0 & 0  \\ 0 & 0 & G & 0 \end{array} \right] - r(B_1)- r(C_1) - r(G) \hfill
\cr
\hspace*{0cm}
 =  r \left[ \begin{array}{cccc} 0  & I_p & 0  & 0  \\ 0  & 0 & 0 & 0  \\ 0 & 0  & -B_1S_1  & B_1X_0C_1  
\\ 0 & 0 & G & 0 \end{array} \right] - r(B_1)- r(C_1) - r(G) \hfill
\cr
\hspace*{0cm}
=  r \left[ \begin{array}{ccc} - B_1 & 0 & B_1X_0C_1  \\ B_1  & B_2 & 0 \end{array}
\right] + p + q - r(B_1)- r(C_1) - r(G) \hfill
\cr
\hspace*{0cm}
=  r[\, B_1, \  B_1X_0C_1 \, ] + p + q - r(C_1) - r(G) =  r[\, B_2, \  A \, ] + p + q - r(C_1) - r(G),   \hfill
\cr}
$$
$$
\displaylines{
\hspace*{0cm}
r \left[ \begin{array}{cc} X_0 & F_{B_1}  \\ E_{C_1} & 0 \\ E_HT_1 & 0\end{array}
\right] \hfill
\cr
\hspace*{0cm}
 =  r \left[ \begin{array}{cccc} X_0 & I_p & 0 & 0  \\ I_q  & 0 & C_1 & 0 \\
T_1 & 0 & 0 & H  \\ 0 & B_1 & 0 & 0 \end{array} \right] - r(B_1)- r(C_1) - r(H) 
 \hfill
\cr
\hspace*{0cm}
= r \left[ \begin{array}{cccc} 0  & I_p & 0 & 0  \\  0  & 0 & 0 & 0 \\
0 & 0 & - T_1C & H  \\ 0 & B_1X_0C  & 0 & 0 \end{array} \right] - r(B_1)- r(C_1) - r(H) \hfill
\cr
\hspace*{0cm}
 =  r \left[ \begin{array}{cc} -C_1 & C_1  \\ 0 & -C_2 \\ B_1X_0C_1 & 0 \end{array}
\right] + p + q - r(B_1) -  r(C_1) - r(H) \hfill
\cr
\hspace*{0cm}
= r \left[ \begin{array}{c} C_2 \\ B_1X_0C_1  \end{array} \right] + p + q -
r(B_1) - r(H) = r \left[ \begin{array}{c} C_2 \\ A  \end{array} \right] + p + q - 
r(B_1) - r(H), \hfill 
\cr}
$$
$$
\displaylines{
\hspace*{0cm}
r \left[ \begin{array}{ccc}  X_0 & F_{B_1} & S_1F_G  \\ E_{C_1} & 0 & 0  \\ E_HT_1 & 0 & 0 \end{array} 
\right] \hfill
\cr
\hspace*{0cm}
 =  r \left[ \begin{array}{ccccc} X_0 & I_p & S_1 & 0 & 0  \\ I_q  & 0 & 0
& C_1 & 0 \\
 T_1 & 0 & 0 & 0 & H  \\  0 & B_1 & 0 & 0 & 0 \\ 0 & 0 & G & 0 & 0  \end{array}
  \right] - r(B_1)- r(C_1) - r(G)
- r(H)    \hfill
\cr
\hspace*{0cm}
= r \left[ \begin{array}{ccccc} 0 & I_p & 0 & 0 & 0  \\ I_q  & 0 & 0
& 0 & 0 \\  0 & 0 & 0 & -T_1C & H  \\  0 & 0 & -B_1S_1 & B_1X_0C_1 & 0 \\ 0 & 0 & G & 0 & 0  \end{array} 
\right] - r(B_1)- r(C_1) - r(G)
- r(H)    \hfill
\cr
\hspace*{0cm}
 =  r \left[ \begin{array}{cccc} 0 & 0 & -C_1 & C_1   \\ 0 & 0 & 0 & -C_2 \\  -B_1
& 0 & B_1X_0C_1 & 0 \\ B_1 & B_2  & 0 & 0  \end{array} \right] + p + q  - r(B_1)- r(C_1) - r(G) - r(H) 
 \hfill
\cr
\hspace*{0cm}
 =  r \left[ \begin{array}{cccc} 0 & 0 & -C_1 & 0   \\ 0 & 0 & 0 & -C_2 \\
-B_1 & 0 & 0 & 0 \\ 0 & B_2  & 0 & A  \end{array} \right] + p + q  - r(B_1)- r(C_1) - r(G) - r(H)  \hfill
\cr
\hspace*{0cm} 
=  r\left[ \begin{array}{cc} A &  B_2 \\ C_2 & 0 \end{array} \right] + p + q
- r(G) - r(H).  
\hfill
\cr}
$$
Thus we have the two formulas in Part (a). By the  similar approach, we can establish
Part (b). \qquad  $\Box$

\medskip

Furthermore, we can also find the  maximal and the minimal ranks of $ B_1XC_1 $ and
$ B_2YC_2$ in (26.1) when it is consistent.

\medskip

\noindent {\bf Theorem 26.4.}\, {\em  Suppose that the matrix equation {\rm(26.1)} is consistent$,$ 
and $ J_1$ and $J_2$ are defined by  {\rm (26.7)}. Then 
$$
\displaylines{
\hspace*{0cm}
\max_{X\in J_1}r(B_1XC_1) =  \min  \left\{ \  r[\, A, \ B_2 \,] - r[\, B_1, \ B_2 \,]
+ r(B_1),  \ \ \  r\left[ \begin{array}{c} A \\ C_2  \end{array} \right] -
r\left[ \begin{array}{c} C_1 \\ C_2  \end{array} \right] + r(C_1)  \ \right\},
 \hfill (26.12)
\cr
\hspace*{0cm}
\min_{X \in J_1}r(B_1XC_1) = r[\, A, \ B_2 \,] + r\left[ \begin{array}{c} A \\ C_2
\end{array} \right] -
r\left[ \begin{array}{cc} A & B_2 \\ C_2 & 0 \end{array} \right],  \hfill (26.13)
\cr
\hspace*{0cm}
\max_{Y \in J_2}r(B_2YC_2) = \min  \left\{ \ r[\, A, \ B_1 \,] - r[\, B_2, \ B_1 \,] +
r(B_2), \ \  \ r\left[ \begin{array}{c} A \\ C_1  \end{array} \right] -
r\left[ \begin{array}{c} C_2 \\ C_1  \end{array} \right]  + r(C_2)  \ \right\},
\hfill (26.14)
\cr
\hspace*{0cm}
\min_{Y \in J_2}r(B_2YC_2) = r[\, A, \ B_1 \,] + r\left[ \begin{array}{c} A \\ C_1
 \end{array} \right] -
r\left[ \begin{array}{cc} A &  B_1\\ C_1 & 0 \end{array} \right]. \hfill (26.15)
\cr}
$$}
{\bf Proof.}\, Putting (26.4) to $B_1XC_1$ and then applying (18.5) and
(18.6) , we find that
$$
\displaylines{
\hspace*{1cm}
\max_{X \in J_1}r(B_1XC_1) \hfill
\cr
\hspace*{1cm}
 =  \max_{U}r( \,  B_1X_0C_1 + B_1S_1F_G UE_HT_1C_1 \, ) = \min \left\{ r[ \,  B_1X_0C_1, \  B_1S_1F_G \, ], \ \ r \left[ \begin{array}{c} B_1X_0C_1 \\
 E_HT_1C_1 \end{array} \right]  \right\}, \hfill 
\cr
\hspace*{1cm}
\min_{X \in J_1}r(B_1XC_1) \hfill
\cr
\hspace*{1cm}
 = \min_{U}r( \,  B_1X_0C_1 + B_1S_1F_G UE_HT_1C_1 \, ) \hfill
\cr
\hspace*{1cm}
=  r[ \,  B_1X_0C_1, \  B_1S_1F_G \, ] + r \left[ \begin{array}{c} B_1X_0C_1 \\ E_HT_1C_1 \end{array} \right] - r \left[ \begin{array}{cc} B_1X_0C_1  & B_1S_1F_G \\ E_HT_1C_1  & 0 \end{array}
 \right].  \hfill
\cr}
$$
Simplifying the ranks of the block matrices in them by Lemma 1.1 and
$ B_1X_0C_1 + B_2Y_0C_2 = A $, we get that
\begin{eqnarray*}  
r[ \,  B_1X_0C_1, \  B_1S_1F_G \, ] & = & r \left[ \begin{array}{cc} B_1X_0C_1  & B_1S_1
\\ 0 & G \end{array} \right]- r(G) \\
 & = & r \left[ \begin{array}{ccc} B_1X_0C_1  & B_1 & 0 \\ 0 & B_1  & B_2 \end{array}
 \right] - r(G)\\
& = & r[ \,  B_1X_0C_1, \ B_2 \, ] + r(B_1) - r(G) =  r[ \, A, \ B_2 \, ] + r(B_1) -
r(G),
\end{eqnarray*} 
\begin{eqnarray*}  
r \left[ \begin{array}{c} B_1X_0C_1 \\ E_HT_1C_1 \end{array} \right] 
& = & r \left[ \begin{array}{cc} B_1X_0C_1  & 0 \\ T_1C_1 & H \end{array} \right]
 -r(H)  \\
& = & r \left[ \begin{array}{cr} B_1X_0C_1  & 0 \\ C_1 & C_1 \\ 0 & -C_2 \end{array}
\right] -r(H) \\
& = & r \left[ \begin{array}{c} B_1X_0C_1  \\ C_2 \end{array} \right] + r(C_1) - r(H) =  
r \left[ \begin{array}{c} A \\ C_2 \end{array} \right] + r(C_1) - r(H), 
\end{eqnarray*} 
\begin{eqnarray*} 
 r \left[ \begin{array}{cc} A_1X_0C_1  & A_1S_1F_G \\ E_HT_1C_1  & 0 \end{array}
 \right] & = &
  r \left[ \begin{array}{ccc} B_1X_0C_1  & B_1S_1 & 0  \\ T_1C_1  & 0  & H  \\ 0 & G
  & 0 \end{array} \right]
- r(G) - r(H) \\
 & = & r \left[ \begin{array}{cccc} B_1X_0C_1  & B_1 & 0 & 0  \\ C_1  & 0  & 0 & C_1
  \\ 0 & 0  & 0 & -C_2  \\ 0 & B_1 & B_2 & 0  \end{array} \right] - r(G) - r(H) \\
 & = & r \left[ \begin{array}{cccc} 0 & B_1 & 0 & 0  \\ C_1  & 0  & 0 & 0  \\
  0& 0  & 0 & C_2  \\ 0 & 0 & B_2 & B_1X_0C_1  \end{array} \right] - r(G) - r(H) \\
& = & r \left[ \begin{array}{cc} A & B_2 \\ C_2 & 0  \end{array} \right] + r(B_1)
+ r(C_1)- r(G) - r(H).
\end{eqnarray*}
Therefore we have (26.12) and (26.13).  Similarly we can show (26.14) and (26.15).
\qquad $\Box$

\medskip

Contrasting (26.9), (26.11), (26.13), (26.13) with (18.5),  we find the following relations  
$$\displaylines{
\hspace*{2cm}
\min_{X \in J_1}r(X) =  \min_{X \in J_1}r(B_1XC_1) = 
\min_{Y}r( \, A - B_2YC_2\, ),  \hfill
\cr
\hspace*{2cm}
\min_{Y \in J_2}r(Y)  = \min_{Y \in J_2}r(B_2YC_2) = \min_{X}r( \, A - B_1XC_1 \, ) . \hfill
\cr}
$$
{\bf Theorem 26.5.}\, {\em  Suppose that the matrix equation {\rm(26.1)} is
consistent$,$ and consider $ J_1$ and $J_2$ in {\rm (26.7)} as two independent
matrix sets. Then
$$ 
\max_{X \in J_1, \, Y \in J_2 }r(\, A - B_1XC_1 - B_2YC_2 \, ) =  \min  \left\{ r(B_1) + r(B_2) - r[\, B_1, \ B_2 \,],
  \ \  r(C_1) + r(C_2) - r\left[ \begin{array}{c} C_1 \\ C_2  \end{array} \right]
  \right\}. \eqno (26.16)
$$
In particular$,$

{\rm (a)}\, Solutions $ X$ and $ Y $ of {\rm(26.1)} are independent$,$ that
is$,$ for any $ X \in J_1$ and $ Y \in J_2$ the pair $ X$ and $ Y $ satisfy
{\rm(26.1)}, if and only if
$$ 
 R(B_1) \cap R(B_2) = \{ 0 \}, \ \ \  or \ \ \  R(C_1^T) \cap R(C_2^T) = \{ 0 \},
  \eqno (26.17)
$$ 
where $ R( \cdot) $ denotes the column space of a matrix. 

{\rm (b)}\, Under {\rm (26.17)}$,$ the general solution of {\rm(26.1)} can be
written as the two independent forms
$$ 
X =  X_0 + S_1Q_G U_1P_HT_1 + F_{B_1}V_1 + V_2E_{C_1}
,   \qquad Y =  Y_0  + S_2Q_G U_2P_HT_2 + F_{B_2}W_1 + 
W_2E_{C_2}, \eqno (26.18)
$$  
where $X_0$ and $ Y_0$ are a pair of special solutions of {\rm (26.1)},
$U_1$, $U_2$, $ V_1$, $ V_2$, $ W_1$ and $ W_2$ are arbitrary. }

\medskip

\noindent {\bf Proof.}\,  Writing (26.5) and (26.6) as two independent matrix expressions, that is, replacing 
$U$ in (26.5) and (26.6)  by $U_1$ and $U_2$ respectively, and then  putting them in 
$ A- B_1XC_1 - B_2YC_2$ yields
\begin{eqnarray*} 
A- B_1XC_1 - B_2YC_2  & = & A- B_1X_0C_1 - B_2Y_0C_2 -  B_1S_1F_G U_1E_HT_1C_1 - B_2S_2F_G U_2E_HT_2C_2  \\
& = & - B_1S_1F_G U_1E_HT_1C_1 - B_2S_2F_G U_2E_HT_2C_2  \\
& = & - B_1S_1F_G U_1E_HT_1C_1 + B_1S_1F_G U_2E_HT_1C_1 \\
& = &  B_1S_1F_G ( \, - U_1 + U_2 \, )E_HT_1C_1, 
\end{eqnarray*} 
where $ U_1$ and $U_2$ are arbitrary. Then by (18.5), it follows that 
\begin{eqnarray*} 
\max_{X \in J_1, \, Y \in J_2 }r(\, A- B_1XC_1 - B_2YC_2\, ) & = &
\max_{U_1, \, U_2}r[ \,  B_1S_1F_G ( \, - U_1 + U_2 \, )E_HT_1C_1 \,] \\
& = &  \min  \left\{ \ r( B_1S_1F_G), \ \ \ r(E_HT_1C_1) \ \right\},  
\end{eqnarray*} 
where 
$$ 
r( B_1S_1F_G) = r\left[ \begin{array}{c} B_1S_1 \\ G  \end{array} \right] - r(G)
= r\left[ \begin{array}{cc} B_1 & 0  \\ B_1  & B_2 \end{array} \right] -
r(G) = r(B_1) + r(B_2) - r(G),
$$ 
$$ 
r(E_HT_1C_1) = r[\, T_1C_1, \ H \,] - r(H) = r\left[ \begin{array}{cr} C_1 & C_1
\\ 0  & - C_2 \end{array} \right]
 - r(H) = r(C_1) + r(C_2) - r(H).
$$ 
Therefore, we have (26.16). The result in (26.17) follows directly from
(26.16) and the solutions in (26.18)
follow from (26.5) and (26.7). \qquad  $\Box$

\markboth{YONGGE  TIAN }
{27. MORE ON EXTREME RANKS OF $ A - B_1 X_1 C_1 - B_2 X_2 C_2$} 

\chapter{ More on extreme ranks of $ A - B_1 X_1 C_1 - B_2 X_2 C_2$ and  related topics }

\noindent In Chapter 19 we have presented  extreme ranks of a matrix expression  
 $$ \displaylines{
\hspace*{2cm}
p(X_1, \, X_2) = A - B_1X_1C_1 - B_2X_2C_2, \hfill (27.1) 
\cr}
$$ 
with respect to $ X_1$ and $X_2$ under some restrictions on the given 
matrices in it. In this chapter we get rid of the restrictions 
to determine the maximal and the minimal ranks of 
$ p(X_1, \, X_2)$ with respect to $ X_1$ and $ X_2$. 

In the two papers \cite{CJRW} by Johnson and \cite{Wo3} by  Woerdeman,  maximal and minimal rank
completions of partial banded block matrices were well examined. Two
general methods for finding maximal and minimal ranks of partial banded
block matrices were established in these two papers. According to the
 general methods, we can simply find the following two special results
 for a $ 3 \times 3 $ partial banded block matrix.

\medskip

\noindent {\bf Lemma 27.1}\cite{CJRW}\cite{Wo3}.\, {\em Let 
$$\displaylines{ 
\hspace*{2cm}
 M  = r \left[ \begin{array}{ccc}  A_{11}  & A_{12} & X \\ A_{21} & A_{22}
 & A_{23}  \\
 Y & A_{32} & A_{33} \end{array} \right], \hfill (27.2)
\cr}
$$
where $ A_{ij} \in {\cal F}^{m_i \times n_j} \ (1 \leq i, \ j \leq 3)$ are
given$,$  $ X \in {\cal F}^{m_1 \times n_3}$ and
$ Y \in {\cal F}^{m_3 \times n_1}$ are two variant matrices. Then
$$
\displaylines{ 
\hspace*{0cm}
\max_{X,  \ Y} r(M) = \min \left\{ m_3 + n_3 + r \left[ \begin{array}{cc}
 A_{11}  & A_{12} \\ A_{21}  & A_{22}
 \end{array} \right], \ \  \
 m_1 + n_1 + r \left[ \begin{array}{cc}  A_{22}  & A_{23} \\ A_{32}  &
 A_{33}  \end{array} \right], \right.    \hfill
\cr
\hspace*{5cm}  
\left.  m_1 + m_3 +  r [ \, A_{21}, \   A_{22}, \  A_{23} \, ],  \ \ \
 n_1 + n_3 + r \left[ \begin{array}{cc}  A_{12} \\ A_{22} \\ A_{32}
 \end{array} \right]  \right\},  \hfill (27.3)
\cr
\hspace*{0cm} 
and \hfill
\cr
\hspace*{0cm}
\min_{X,  \ Y} r(M) = r [ \, A_{21}, \   A_{22}, \  A_{23} \, ] + 
r \left[ \begin{array}{c}  A_{12} \\ A_{22} \\ A_{32} 
\end{array} \right]
 + \max  \left\{  r \left[ \begin{array}{cc}  A_{11}  & A_{12} \\  A_{21}
  & A_{22}  \end{array} \right] -
 r \left[ \begin{array}{cc}  A_{12} \\ A_{22} \end{array} \right] - 
r [ \, A_{21}, \ A_{22} \,] \right. , \hfill 
\cr
\hspace*{6.5cm}  
\left.  r \left[ \begin{array}{cc}  A_{22}
 & A_{23} \\
A_{32}  & A_{33}  \end{array} \right] - r \left[ \begin{array}{cc}  A_{22}
 \\ A_{32} \end{array} \right] -
 r[ \, A_{22}, \ A_{23} \,] \right\}. \hfill (27.4) 
\cr}
$$ }
\hspace*{0.3cm} Notice that the block matrix $M$ in (27.2) and the matrix expression in
(27.1) have  two independent variant matrices, respectively.  This fact
motivates us to express the rank of (27.1) as the rank of a block matrix,
and then apply (27.3) and (27.4) to determine extreme ranks of
(27.1) with respect to $ X_1$ and $X_2$.  

It is easy to verify by block elementary
operations of matrices that the rank  of $p(X_1, \, X_2)$ in (27.1) satisfies
the equality
$$\displaylines{
\hspace*{2cm}
 r[\, p(X_1, \, X_2)\, ] = 
r \left[ \begin{array}{ccccc}  0 & 0 & 0 & I_{p_2} & - X_2  \\  0 & 0 & C_2
& 0  & I_{q_2} \\
  0 & B_1 & A & B_2 & 0 \\ I_{q_1} & 0 & C_1 & 0 & 0\\  -X_1  & I_{p_1} &
  0 & 0 & 0
\end{array} \right] - p_1 - p_2 - q_1 - q_2. \hfill (27.5) 
\cr}
$$
Applying Lemma 1.1 to the block matrix in (27.5) and simplying, we obtain
the main result of the chapter.

\medskip

\noindent  {\bf Theorem 27.2.}\, {\em Let $ p(X_1, \, X_2)$ be given by
{\rm (27.1)}.  Then
$$
\displaylines{
\hspace*{0cm} 
\max_{X_1, \, X_2} r[ \, p( X_1, \, X_2 )\,] = \min  \left\{ \ r[ \, A, \
B_1, \  B_2 \, ], \ \
r \left[ \begin{array}{c}  A  \\ C_1 \\ C_2  \end{array} \right], \ \
r \left[ \begin{array}{cc}  A
& B_1\\ C_2   & 0 \end{array} \right], \ \ r \left[ \begin{array}{cc}  A   
& B_2 \\ C_1   & 0 \end{array} \right]  \ \right\}, \hfill (27.6)
\cr
\hspace*{0cm} 
and \hfill
\cr
\hspace*{0cm} 
\min_{X_1, \, X_2} r[ \, p( X_1, \, X_2 )\,] \hfill
\cr
\hspace*{1cm} 
= r \left[ \begin{array}{c}  A  \\ C_1 \\ C_2  \end{array} \right] + 
 r[ \, A, \ B_1, \  B_2 \, ] + \max \left\{ r\left[ \begin{array}{cc} A & B_1\\ C_2 & 0 \end{array} \right] - r\left[ \begin{array}{ccc} A   
& B_1 & B_2  \\ C_2   & 0 & 0  \end{array} \right] - r\left[ \begin{array}{cc} A   
& B_1 \\ C_1   & 0 \\ C_2 & 0  \end{array} \right] \right., \hfill
\cr
\hspace*{5.5cm} 
 \left. r \left[ \begin{array}{cc}  A   
& B_2 \\ C_1 & 0 \end{array} \right]  - r\left[ \begin{array}{ccc} A   
& B_1 & B_2  \\ C_1   & 0 & 0  \end{array} \right] - r\left[ \begin{array}{cc} A   
& B_2 \\ C_1   & 0 \\ C_2 & 0  \end{array} \right] \ \right\}. \hfill (27.7)
\cr}
$$
} 
\hspace*{0.3cm}The two rank equalities in (27.6) and (27.7) can help to reveal some
fundamental properties of the matrix expression $ p( X_1, \, X_2 )$ in
(27.1).  For example, let
$$\displaylines{
\hspace*{2cm} 
\max_{X_1, \, X_2} r[ \, p( X_1, \, X_2 )\,] = \min_{X_1, \, X_2} r[ \, p( X_1, \, X_2 )\,]  = r(A), \hfill
\cr}
$$
one can immediately establish  a necessary and sufficient condition for the
rank of $p( X_1, \, X_2 )$ to
be invariant with respect to $  X_1 $ and  $X_2$. Notice that two matrices
$M$ and $ N$ have the same column space if and only if $ r[ \, M,  \ N \,]
= r(M ) = r(N)$. Thus the  column space of  $ p( X_1, \, X_2 )$ is invariant
with respect to $  X_1 $ and  $X_2$ if and only if 
$$\displaylines{
\hspace*{2cm} 
r[ \, p( X_1, \, X_2 ), \ p( Y_1, \, Y_2 ) \,] =  
r[ \, p( X_1, \, X_2 )\,] =
r[\, p( Y_1, \, Y_2 ) \,] = r(A) \hfill
\cr}
$$
holds for all $ X_1, \ X_2 , \ Y_1, \ Y_2$,  where
\begin{eqnarray*}
[\, p( X_1, \, X_2 ), \ p( Y_1, \, Y_2 ) \,] & =
& [\,  A - B_1X_1C_1 - B_2X_2C_2,  \  A - B_1Y_1C_1 - B_2Y_2C_2 \,] \\
& =  & [\, A, \  A \,] - B_1[\, X_1, \  Y_1\,]\left[ \begin{array}{cc} C_1 &  0 \\ 0  & C_1 \end{array} \right] 
 -  B_2[\, X_2, \  Y_2 \,]\left[ \begin{array}{cc} C_2 &  0 \\ 0  & C_2
 \end{array} \right].
\end{eqnarray*} 
Thus applying Theorem 27.1 to the equality, one can also establish a
necessary and sufficient condition for the column space
of $p( X_1, \, X_2 )$ to be invariant with respect to $  X_1 $ and
 $X_2$.  Moreover let (27.7) be zero, we can trivialy obtain a solvability
 condition for the matrix equation $ B_1X_1C_1 + B_2X_2C_2 = A$, which has
been established previously  by  \"{O}zg\"{u}ler in \cite{Oz}.
 
\medskip

\noindent  {\bf Corollary 27.3.}\, {\em There exist 
$ X_1$ and $X_2$ such that $ B_1X_1C_1 + B_2X_2C_2 = A$ if and only if 
$$ \displaylines{
\hspace*{2cm} 
r[ \, A, \ B_1, \  B_2 \, ] = r[ \,B_1, \  B_2 \, ] ,  \qquad   r \left[
\begin{array}{c}  A  \\ C_1 \\ C_2  \end{array} \right] =
r \left[ \begin{array}{c} C_1 \\ C_2  \end{array} \right], \hfill
\cr
\hspace*{2cm} 
 r \left[ \begin{array}{cc}  A & B_1 \\ C_2 & 0 \end{array}
 \right] = r(B_1) + r(C_2),  \qquad  r \left[ \begin{array}{cc}  A & B_2 \\ C_1 & 0 \end{array}
 \right] = r(B_2) + r(C_1).   \hfill
\cr}
$$}

Combining the two formulas (27.6), (26.7) and those in Chapter 19, we can also
 establish some more general results for linear  matrix expressions with four
 two-sided  independent variant matrices, which, in turn, will apply to determine 
extreme  ranks for some more general matrix expressions.    

\medskip

\noindent {\bf Theorem 27.4.}\, {\em  Let 
$$ \displaylines{
\hspace*{2cm} 
p(X_1, \, X_2, \  X_3, \, X_4) = A - B_1X_1C_1 - B_2X_2C_2 - B_3X_3C_3 - B_4X_4C_4, \hfill (27.8) 
\cr}
$$ 
be a linear  matrix expression with four two-sided terms over an arbitrary field ${\cal F},$  and suppose that the given 
matrices satisfying the conditions 
$$\displaylines{
\hspace*{2cm} 
R(B_i) \subseteq R(B_2),  \ \ and   \ \  R(C_j^T) \subseteq R(C_1^T), \ \ \ \ i = 1, \ 3, \ 4, \ \ j = 2, \ 3, \ 4. 
\hfill (27.9)
\cr
\hspace*{0cm} 
Then   \hfill
\cr
\hspace*{0cm} 
\max_{X_i} r[ \, p(X_1, \, X_2, \  X_3, \, X_4)\,] = \min  \left\{ r[ \, A, \  B_2 \, ],  \
r \left[ \begin{array}{c}  A  \\ C_1 \end{array} \right],  \
r \left[ \begin{array}{cc}  A  & B_1 \\ C_2 & 0\\ C_3 & 0 \\ C_4 & 0  \end{array} \right],  \ r \left[ \begin{array}{cccc} A   & B_1 & B_3 & B_4 \\ C_2  & 0 & 0 & 0 \end{array} \right], \right. \hfill 
\cr
\hspace*{7cm} 
\left. r \left[ \begin{array}{ccc}  A  & B_1 & B_3 \\ C_2 & 0 & 0 \\ 
C_4 & 0 & 0  \end{array} \right], \ \   r \left[ \begin{array}{ccc}  A  & B_1 & B_4 \\ C_2 & 0 & 0 \\ 
C_3 & 0 & 0  \end{array} \right]  \right\}, \hfill (27.10)
\cr}
$$
and 
$$
\displaylines{ 
\hspace*{0cm} 
\min_{X_i} r[ \, p(X_1, \, X_2, \  X_3, \, X_4)\,]  \hfill
\cr
\hspace*{0cm} 
= r \left[ \begin{array}{cc}  A  & B_1 \\ C_2 & 0 \\ C_3 & 0 \\ C_4 & 0 \end{array} \right] + r \left[ \begin{array}{cccc} A   
& B_1 & B_3 & B_4 \\ C_2   & 0 & 0 & 0 \end{array} \right] + r \left[ \begin{array}{c}  A  \\ C_1 \end{array} \right] 
+  r[ \, A, \  B_2 \, ] - r \left[ \begin{array}{cc}  A  & B_1 \\ C_1 & 0
\end{array} \right] - r \left[ \begin{array}{cc} A
& B_2 \\ C_2   & 0 \end{array} \right] \hfill
\cr
\hspace*{0.5cm} 
+  \max \left\{ r\left[ \begin{array}{ccc} A & B_1 & B_3 \\  C_2 & 0 & 0
\\  C_4 & 0 & 0 \end{array} \right]
- r\left[ \begin{array}{cccc} A & B_1 & B_3  & B_4 \\ C_2   & 0 & 0 & 0
\\ C_4 & 0 & 0 & 0  \end{array}
\right] - r\left[ \begin{array}{ccc} A & B_1 & B_3  \\ C_2 & 0 & 0
\\ C_3 & 0 & 0 \\ C_4 & 0 & 0 \end{array}
\right] \right., \hfill
\cr
\hspace*{2.5cm} 
 \left. r \left[ \begin{array}{ccc}  A   
& B_1  & B_4 \\ C_2 & 0 & 0  \\ C_3 & 0 & 0  \end{array} \right] -
r\left[ \begin{array}{cccc} A & B_1 & B_3 & B_4 \\ C_2   & 0 & 0 & 0 \\
C_3 & 0 & 0 & 0 \end{array} \right] - r\left[ \begin{array}{ccc} A
& B_1 & B_4 \\ C_2 & 0 & 0 \\ C_3 & 0 & 0 \\  C_4 & 0 & 0 \end{array}
\right]  \right\}. \hfill (27.11)
\cr}
$$
}
{\bf Proof.}\, We only show (27.11). Under (27.9), we apply (19.4) 
to the two variant matrices $ X_1$ and $X_2$ in (27.8) to yield
$$
\displaylines{ 
\hspace*{0cm} 
\min_{X_1, \, X_2} r[ \, p( X_1, \, X_2, \ X_3, \, X_4 )\,] \hfill
\cr
\hspace*{0cm} 
= r[ \, A - B_3X_3C_3 - B_4X_4C_4, \ B_2 \, ] + r \left[ \begin{array}{c}  A - B_3X_3C_3 - B_4X_4C_4  \\ C_1 \end{array} 
\right] + r\left[ \begin{array}{cc} A - B_3X_3C_3 - B_4X_4C_4  & B_1 \\ C_2 & 0 \end{array} \right] \hfill
\cr
\hspace*{0.5cm} 
- \ r \left[ \begin{array}{cc}  A - B_3X_3C_3 - B_4X_4C_4  & B_1 \\ C_1   & 0 \end{array} \right]
- r \left[ \begin{array}{cc}  A- B_3X_3C_3 - B_4X_4C_4 & B_2 \\ C_2 & 0 \end{array} \right] \hfill
\cr
\hspace*{0cm}
= r[ \, A, \ B_2 \, ] + r \left[ \begin{array}{c}  A \\ C_1 \end{array} 
\right] - r\left[ \begin{array}{cc} A   & B_1 \\ C_1 & 0 \end{array} \right] -
 r \left[ \begin{array}{cc}  A & B_2 \\ C_2   & 0 \end{array} \right]
+ r \left[ \begin{array}{cc}  A- B_3X_3C_3 - B_4X_4C_4 & B_1 \\ C_2 & 0 \end{array} \right]. \hfill(27.12)
\cr}
$$
Notice that 
$$
 \left[ \begin{array}{cc}  A- B_3X_3C_3 - B_4X_4C_4 & B_1 \\ C_2 & 0 \end{array} \right]
= \left[ \begin{array}{cc}  A & B_1 \\ C_2 & 0 \end{array} \right] -  \left[ \begin{array}{c} B_3 
\\  0 \end{array} \right]X_3[\, C_3, \ 0 \,]  -  \left[ \begin{array}{c} B_4 \\  0 \end{array} \right]X_4[\, C_4, \ 0 \,].  $$
In that case, applying (27.7) to it and then putting the corresponding
 result in (27.12) yields (27.11). \qquad $\Box$  

\medskip

\noindent {\bf Corollary  27.5.}\, {\em  Let  
$$ \displaylines{
\hspace*{2cm} 
p(X_1, \, X_2, \  X_3, \, X_4) = A - B_1X_1 - X_2C_2 - B_3X_3C_3 - B_4X_4C_4 \hfill (27.13) 
\cr}
$$ 
be a linear matrix expression over an arbitrary field ${\cal F}$ with two one-sided terms and 
two two-sided terms. Then  
$$
\displaylines{
\hspace*{0cm} 
\max_{ \{X_i\}} r[ \, p(X_1, \, X_2, \  X_3, \, X_4)\,] = \min  \left\{ m, \ \  n ,  \ \
r \left[ \begin{array}{c}  A  \\ C_1 \end{array} \right],  \
r \left[ \begin{array}{cc}  A  & B_1 \\ C_2 & 0\\ C_3 & 0 \\ C_4 & 0  \end{array} \right],  \ \ 
r \left[ \begin{array}{cccc} A   & B_1 & B_3 & B_4 \\ C_2  & 0 & 0 & 0 \end{array} \right], \right. \hfill 
\cr
\hspace*{7cm} 
\left. r \left[ \begin{array}{ccc}  A  & B_1 & B_3 \\ C_2 & 0 & 0 \\ 
C_4 & 0 & 0  \end{array} \right], \ \   r \left[ \begin{array}{ccc}  A  & B_1 & B_4 \\ C_2 & 0 & 0 \\ 
C_3 & 0 & 0  \end{array} \right]  \right\}, \hfill (27.14)
\cr
\hspace*{0cm} 
and \hfill
\cr
\hspace*{0cm} 
\min_{\{X_i\}} r[ \, p(X_1, \, X_2, \  X_3, \, X_4)\,] = r \left[ \begin{array}{cc}  A  & B_1 \\ C_2 & 0 \\ C_3 & 0 \\ C_4 & 0 \end{array} \right] + r \left[ \begin{array}{cccc} A   
& B_1 & B_3 & B_4 \\ C_2   & 0 & 0 & 0 \end{array} \right] - r(B_1) - r(C_2) \hfill
\cr
\hspace*{2.3cm} 
+  \max \left\{ r\left[ \begin{array}{ccc} A & B_1 & B_3 \\  C_2 & 0 & 0
\\  C_4 & 0 & 0 \end{array} \right]
- r\left[ \begin{array}{cccc} A & B_1 & B_3  & B_4 \\ C_2   & 0 & 0 & 0
\\ C_4 & 0 & 0 & 0  \end{array}
\right] - r\left[ \begin{array}{ccc} A & B_1 & B_3  \\ C_2 & 0 & 0
\\ C_3 & 0 & 0 \\ C_4 & 0 & 0 \end{array}
\right] \right., \hfill
\cr
\hspace*{3.2cm} 
 \left. r \left[ \begin{array}{ccc}  A   
& B_1  & B_4 \\ C_2 & 0 & 0  \\ C_3 & 0 & 0  \end{array} \right] -
r\left[ \begin{array}{cccc} A & B_1 & B_3 & B_4 \\ C_2   & 0 & 0 & 0 \\
C_3 & 0 & 0 & 0 \end{array} \right] - r\left[ \begin{array}{ccc} A
& B_1 & B_4 \\ C_2 & 0 & 0 \\ C_3 & 0 & 0 \\  C_4 & 0 & 0 \end{array}
\right]  \right\}. \hfill (27.15)
\cr}
$$
In particular$,$ the matrix equation 
$$ \displaylines{
\hspace*{2cm} 
B_1X_1 + X_2C_2 + B_3X_3C_3 + B_4X_4C_4 = A \hfill (27.16) 
\cr}
$$ 
is consistent if and only if
$$
\displaylines{
\hspace*{0.5cm} 
r \left[ \begin{array}{cc}  A  & B_1 \\ C_2 & 0\\ C_3 & 0 \\ C_4 & 0  \end{array} \right] = 
r \left[ \begin{array}{cc}  0  & B_1 \\ C_2 & 0\\ C_3 & 0 \\ C_4 & 0  \end{array} \right], 
\ \ \
r \left[ \begin{array}{cccc} A   & B_1 & B_3 & B_4 \\ C_2  & 0 & 0 & 0 \end{array} \right] = r \left[ \begin{array}{cccc} 0  & B_1 & B_3 & B_4 \\ C_2  & 0 & 0 & 0 \end{array} \right], 
 \hfill (27.17)
\cr
\hspace*{0.5cm} 
 r \left[ \begin{array}{ccc}  A  & B_1 & B_3 \\ C_2 & 0 & 0 \\ 
C_4 & 0 & 0  \end{array} \right] = 
r \left[ \begin{array}{ccc}  0  & B_1 & B_3 \\ C_2 & 0 & 0 \\ C_4 & 0 & 0  \end{array} \right],  \ \ \  r \left[ \begin{array}{ccc}  A  & B_1 & B_4 \\ C_2 & 0 & 0 \\ 
C_3 & 0 & 0  \end{array} \right] = r \left[ \begin{array}{ccc}  0  & B_1 & B_4 \\ C_2 & 0 & 0 \\ 
C_3 & 0 & 0  \end{array} \right]. \hfill (27.18)
\cr}
$$
}\hspace*{0.3cm} Based on (27.10) and (27.11), we can determine extreme  ranks
 of $ A_1 - B_1XC_1$ subject to a pair of 
consistent matrix equations $ B_2XC_2 = A_2$ and $ B_3XC_3 = A_3$. From 
them we can find lots of valuable results related to solvability and 
solutions of some matrix equations. We shall present them in the next chapter. 

Some more general work than those in the chapter is to determine extreme ranks of a linear  matrix expression
$ A - B_1X_1C_1 -  B_2X_2C_2 - B_3X_3C_3$  with respect to $X_1, \ X_2$ and $X_3$, as well as 
$ A - B_1X_1C_1 - \cdots - B_kX_kC_k, \ k > 3$  with respect to $X_1$---$X_k$ without any restrictions to the given 
matrices in them. According to the method presented by Johnson in \cite{Jo}, the maximal rank of 
$ A - B_1X_1C_1 - \cdots - B_kX_kC_k$  can completely be determined. Here we only list the case for $ k = 3$ without 
its tedious proof.

\medskip

\noindent  {\bf Theorem 27.6.}\, {\em  Let 
$$ 
p (\, X_1, \, X_2, \  X_3 \, ) = A - B_1X_1C_1 - B_2X_2C_2 - B_3X_3C_3
$$ 
be a matrix expression over an arbitrary field ${\cal F}.$ Then  
$$
\displaylines{
\hspace*{0cm} 
\max_{\{X_i\}} r[ \, p(X_1, \, X_2, \  X_3 \,)\,] = \min  \left\{r \left[ \begin{array}{c}  A  \\ C_1 \\ C_2 \\ C_3  \end{array} \right],  \ \ r \left[ \begin{array}{cc}  A  & B_1 \\ C_2 & 0\\ C_3 & 0 \end{array} \right],  \ \ 
r \left[ \begin{array}{cc}  A  & B_2 \\ C_1 & 0\\ C_3 & 0 \end{array} \right],  \ \ 
r \left[ \begin{array}{cc}  A  & B_3 \\ C_1 & 0\\ C_2 & 0 \end{array} \right],  \right. \hfill 
\cr
\hspace*{1cm} 
\left.
r \left[ \begin{array}{ccc}  A  & B_1 & B_2 \\  C_3 & 0 & 0  \end{array} \right],   \ \ 
r \left[ \begin{array}{ccc}  A  & B_1 & B_3 \\  C_2 & 0 & 0  \end{array} \right],   \ \ 
r \left[ \begin{array}{ccc}  A  & B_2 & B_3 \\  C_1 & 0 & 0  \end{array} \right], \ \ r[ \, A, \  B_1,  \ B_2, \ B_3 \, ]  \right\}. \hfill (27.19)
\cr}
$$}

>From the eight block matrices in (27.19), the reader can easily infer the maximal rank of
 $ A - B_1X_1C_1 - \cdots - B_kX_kC_k,$  in which, the $ 2^k$ block matrices  are much similar to those in (27.19). 

As to the minimal rank of $ A - B_1X_1C_1 - \cdots - B_kX_kC_k$ when $ k \geq 3$, the process to find it becomes quite 
complicated. We do not find at present a general method to solve this challenging problem.  However, as we have seen in 
Theorem 27.4 and Corollary 27.5, if the given matrices in a  matrix expression satisfy some restrictions, then we 
can still find its minimal rank. Here we list two simple results.  

\medskip

\noindent {\bf Theorem 27.7.}\, {\em Suppose $A_{ij}( 1 \leq i, \ j \leq 3)$ are all
 nonsingular matrices of order $m$. Then
$$\displaylines{
\hspace*{3cm}
\min_{X, \, Y, \, Z} r\left[ \begin{array}{ccc}  X &  A_{12}  &  A_{13} \\
A_{21} & Y & A_{23} \\ A_{31} & A_{32} &  Z \end{array} \right]
 =\left\{ \begin{array}{ll} r(A)/2   &  if \ r(A) \ is \ 
even \\ \left[ \, r(A) + 1 \, \right]/2 & if \ r(A) \ is \ odd
\end{array} \right.. \hfill (27.20)
\cr
where \hfill
\cr
\hspace*{3cm}
A = \left[ \begin{array}{ccc}  0 &  -A_{12}  &  -A_{13} \\ A_{21} & 0 &
 -A_{23} \\ A_{31} & A_{32} &  0 \end{array} \right].  \hfill (27.21)
\cr}
$$
}
{\bf Proof.}\, According to (27.4), we first find that 
$$
\displaylines{
\hspace*{0cm}
\min_{X, \, Z} r\left[ \begin{array}{ccc}  X &  A_{12}  &  A_{13} \\
A_{21} & Y & A_{23} \\ A_{31} & A_{32} &  Z \end{array} \right] = r [ \, A_{21}, \   Y, \  A_{23} \, ] + 
r \left[ \begin{array}{c}  A_{12} \\ Y \\ A_{32} \end{array} \right] \hfill
\cr
\hspace*{1cm}
 + \max  \left\{  r \left[ \begin{array}{cc}  A_{12}  & A_{13} \\  Y
  & A_{23}  \end{array} \right] -
 r \left[ \begin{array}{cc}  A_{12} \\ Y \end{array} \right] - 
r [ \, Y, \ A_{23} \,], \ \   r \left[ \begin{array}{cc}  A_{21}
 & Y \\ A_{31}  & A_{32}  \end{array} \right] - r \left[ \begin{array}{cc}  Y
 \\ A_{32} \end{array} \right] -
 r[ \, A_{21}, \ Y \,] \right\} \hfill 
\cr
\hspace*{0cm}
= 2m + \max  \left\{  r \left[ \begin{array}{cc}  0  & A_{13} \\  Y - A_{23}A^{-1}_{13}A_{12}  & 0 
  \end{array} \right] - 2m, \ \   r \left[ \begin{array}{cc}  0
 & Y - A_{21}A^{-1}_{31}A_{32} \\ A_{31}  & 0  \end{array} \right] - 2m  \right\} \hfill 
\cr
\hspace*{0cm}
=  m + \max  \left\{ \  r(\, Y - A_{23}A^{-1}_{13}A_{12} \,), \ \ \   r (\,  Y - A_{21}A^{-1}_{31}A_{32} \,) \,
 \right\} \hfill 
\cr
\hspace*{0cm}
= m +  \max  \left\{ \  r( \widehat{Y}), \ \ \   r (\,  \widehat{Y} - M \,) \, \right\}, \hfill 
\cr}
$$ 
where $  \widehat{Y} =  Y - A_{23}A^{-1}_{13}A_{12}$ and $ M =  A_{21}A^{-1}_{31}A_{32} - A_{23}A^{-1}_{13}A_{12}$. 
Thus
$$
\min_{X, \, Y, \, Z} r\left[ \begin{array}{ccc}  X &  A_{12}  &  A_{13} \\
A_{21} & Y & A_{23} \\ A_{31} & A_{32} &  Z \end{array} \right] 
=  m + \min_{\widehat{Y}} \max  \left\{ \  r( \widehat{Y}), \ \ \   r (\,  M - \widehat{Y} \,) \, \right\}. 
$$ 
Notice that $  r (\,  M -  \widehat{Y}\,) \geq  r(M) -  r( \widehat{Y})$ for all $\widehat{Y}$. We see that 
$$
\displaylines{
\hspace*{2cm}
 \max  \left\{ \  r( \widehat{Y}), \ \ \   r (\,  M  - \widehat{Y} \,) \, \right\} \geq 
\max  \left\{ \  r( \widehat{Y}), \ \ \    r(M) - r(\widehat{Y}) \, \right\},  \hfill
\cr
and \hfill
\cr
\hspace*{2cm} 
\min_{\widehat{Y}} \max  \left\{ \  r( \widehat{Y}), \ \ \   r (\,  M - \widehat{Y} \,) \, \right\}  \geq
\min_{\widehat{Y}} \max  \left\{ \  r( \widehat{Y}), \ \ \    r(M) - r(\widehat{Y}) \, \right\}. \hfill
\cr}
$$
Since $ \widehat{Y}$ is  arbitrary, we easily get  that 
$$
\displaylines{
\hspace*{2cm}
\min_{\widehat{Y}} \max  \left\{ \  r( \widehat{Y}), \ \ \  r(M) -  r(\widehat{Y})  \ \right\} 
= \left\{ \begin{array}{ll} r(M)/2   &  if \ r(M) \ is \ 
even \\ \left[ \, r(M) + 1 \, \right]/2 & if \ r(M) \ is \ odd
\end{array} \right.. \hfill
\cr}
$$
Consequently,
$$
\displaylines{
\hspace*{2cm}
\min_{\widehat{Y}} \max  \left\{ \  r( \widehat{Y}), \ \ \   r (\,  M - \widehat{Y} \,) \, \right\}  \geq
 \left\{ \begin{array}{ll} r(M)/2   &  if \ r(M) \ is \ 
even \\ \left[ \, r(M) + 1 \, \right]/2 & if \ r(M) \ is \ odd
\end{array} \right.. \hfill (27.22)
\cr}
$$
We next show that  the lower bound in the right side pf (27.22) can be reached by the left hand side of (27.22) by 
choosing some $\widehat{Y}$. In fact, suppose $ M $ can factor as 
$ M =  P \left[ \begin{array}{cc} I_k  & 0  \\ 0  & 0  \end{array} \right]Q,$  where $ P$ and $ Q $ are nonsingular.
If $ k = r(M)$ is even, we take $ \widehat{Y}=  P \left[ \begin{array}{cc} I_{k/2}  & 0  \\ 0  & 0  \end{array} 
\right]Q.$  In that case, $ r(\, M - \widehat{Y} \, )  = r( \widehat{Y}) = k/2$. 
If $ k = r(M)$ is odd, then we take $ \widehat{Y}=  P \left[ \begin{array}{cc} I_{(k+1)/2}  & 0  \\ 0  & 0  \end{array} 
\right]Q.$ In that case, $ r(\, M - \widehat{Y} \, )  = (k-1)/2$ and $r( \widehat{Y}) = (k+1)/2$. There two cases show 
that  the left hand side of (27.22) can be reached by the right hand side of (27.22). Hence we have 
$$
\displaylines{
\hspace*{2cm}
\min_{\widehat{Y}} \max  \left\{ \,  r( \widehat{Y}), \ \ \   r (\,  M - \widehat{Y} \,) \, \right\} =
 \left\{ \begin{array}{ll} r(M)/2   &  if \ r(M) \ is \ 
even \\ \left[ \, r(M) + 1 \, \right]/2 & if \ r(M) \ is \ odd
\end{array} \right.. \hfill
\cr}
$$
Consequently, 
$$
\displaylines{
\hspace*{2cm}
\min_{X, \, Y, \, Z} r\left[ \begin{array}{ccc}  X &  A_{12}  &  A_{13} \\
A_{21} & Y & A_{23} \\ A_{31} & A_{32} &  Z \end{array} \right] 
 =  m + \left\{ \begin{array}{ll} r(M)/2   &  if \ r(M) \ is \ 
even \\ \left[ \, r(M) + 1 \, \right]/2 & if \ r(M) \ is \ odd
\end{array} \right.. \hfill (27.23)
\cr}
$$
On the other hand, it is easy to verify that 
$$
r\left[ \begin{array}{ccc} 0 &  -A_{12}  &  -A_{13} \\
A_{21} & 0 & -A_{23} \\ A_{31} & A_{32} &  0 \end{array} \right]  = 
r\left[ \begin{array}{ccc} 0 &  0  &  -A_{13} \\
0  & A_{23}A^{-1}_{13}A_{12} - A_{21}A^{-1}_{31}A_{32} & 0  \\ A_{31} &  0  & 0 \end{array}
 \right]  = 2m + r(M). 
$$ 
Hence $ r(M) = r(A) - 2m$. Putting it in (27.23) yields (27.20).  \qquad $ \Box$ 

\medskip

Clearly the $ 3 \times 3 $ block matrix in (27.20) is a special case of 
$A - B_1X_1C_1 -  B_2X_2C_2 - B_3X_3C_3$. If $ A_{ij}'$s are singular or 
are not square, the formula (27.20) is not valid. But we guess that its minimal rank
can be expressed through  the ranks of $ A $ and its submatrices.   

\medskip

\noindent  {\bf Theorem 27.8.}\, {\em  Let 
$$\displaylines{
\hspace*{1cm}  
p(X_1, \ \cdots, \  X_5 \, ) = A - B_1X_1C_1 - B_2X_2C_2 - B_3X_3C_3 - B_4X_4C_4 - B_5X_5C_5 \hfill
\cr}
$$ 
be a matrix expression over an arbitrary field ${\cal F}$ with the given matrices satisfying the conditions 
$$\displaylines{
\hspace*{1cm} 
R[\, B_2, \ B_3 \, ] \subseteq R(B_1),  \ \ and   \ \  R[\, C_4^T, \ C_5^T \,] \subseteq R(C_1^T). \hfill
\cr
\hspace*{0cm} 
Then  \hfill
\cr
\hspace*{0cm} 
\min_{X_i} r[ \, p(X_1, \ \cdots , \  X_5\,)\,] \hfill
\cr
\hspace*{0cm} 
= \min_{X_4, \, X_5}r[\,  A - B_4X_4C_4 - B_5X_5C_5, \ B_1\,] + \min_{X_2, \, X_3}r \left[ \begin{array}{c} 
 A - B_2X_2C_2 - B_3X_3C_3  \\ C_1  \end{array} \right] -  r \left[ \begin{array}{cc}  A  & B_1 \\ C_1 & 0  \end{array} \right]. \ \hfill (27.24)
\cr}
$$}

Clearly the two minimal ranks in (27.24) can further be determined by (27.7). We leave it to the reader. 

It should be mentioned  that although  failing  to give the minimal rank of the matrix expression $ A $ $ -B_1X_1C_1 - 
B_2X_2C_2 - B_3X_3C_3$ in general cases,  we can still express  consistency condition using rank equalities 
 for the corresponding linear matrix equation   
$$
B_1X_1C_1 + B_2X_2C_2 + B_3X_3C_3 = A.  \eqno (27.25) 
$$ 
Here we only list the result, its proof will presented in chapter 28. 

\medskip

\noindent  {\bf Theorem 27.9.}\, {\em The matrix equation {\rm (27.25)} is consistent if and only if the following 
nine rank equalities hold  
$$  
\begin{array}{ll}
r \left[ \, A, \ B_1, \ B_2, \  B_3 \, \right] = r  \left[ \,B_1, \  B_2,  \ B_3 \, \right], &  r \left[
\begin{array}{c}  A  \\ C_1 \\ C_2 \\ C_3 \end{array} \right] =
r \left[ \begin{array}{c} C_1 \\ C_2 \\ C_3  \end{array} \right],   \\
 r\left[ \begin{array}{cc}  A  & B_1 \\ C_2 & 0 \\ C_3 & 0 \end{array} \right] 
= r \left[ \begin{array}{c}  C_2 \\ C_3 \end{array}\right] + r(B_1), &  r \left[ \begin{array}{cc} A & B_2 \\ C_1 & 0\\ C_3 & 0 \end{array} \right] = r \left[ \begin{array}{c}
  C_1 \\ C_3 \end{array} \right] + r(B_2), \\
 r \left[ \begin{array}{cc}  A  & B_3 \\ C_1& 0\\ C_2 & 0 \end{array} \right] 
= r \left[ \begin{array}{c}  C_1 \\ C_2 \end{array}\right] + r(B_3), &  r \left[ \begin{array}{ccc}  A  & B_1 & B_2 \\  C_3 & 0 & 0  \end{array} \right]
 = r\left[\, B_1, \ B_2\, \right] + r(C_3),  \\
r \left[ \begin{array}{ccc}  A  & B_1 & B_3 \\  C_2 & 0 & 0  \end{array} \right] 
= r\left[\, B_1, \ B_3 \, \right] + r(C_2), &
r \left[ \begin{array}{ccc}  A  & B_2 & B_3 \\  C_1 & 0 & 0  \end{array} \right] = r\left[\, B_2, \ B_3 \,\right] + 
r(C_1),
\end{array}
$$
$$
r \left[ \begin{array}{ccccc}  A  &  0 & B_1 & 0  &  B_3  \\
                             0  & -A & 0   & B_2&  B_3\\ 
                             C_2& 0  & 0   & 0  &  0  \\
                              0 & C_1& 0  & 0 & 0  \\
                            C_3 & C_3& 0  & 0 & 0 \end{array} \right] 
= r\left[ \begin{array}{cc} C_2  & 0  \\ 0  & C_1 \\ C_3 & C_3 \end{array} \right] + 
 r \left[ \begin{array}{ccc} B_1  & 0 & B_3 \\ 0   & B_2  & B_3  \end{array} \right]. 
$$}

Of course, those rank equalities can also  equivalently be expressed by equivalence of matrices, column or row spaces of matrices, generalized inverses of matrices, and decompositions of matrices, and so on. The reader can easily
 list them according to Lemma 1.2.    

Solvability and solutions of linear matrix equations have been one of principal topics in matrix 
theory and its applications. Based on the well-known Kronecker product of matrices, one can simply transform 
any kind of linear matrix equations to a standard form $ Mx= b$, and then solve through it. 
Nearly all characteristics of the original equations, however, are lost in this kind of transformations. So one has 
been seeking various feasible methods to solve linear matrix equations without using the Kronecker product. As far 
as the author knows, Theorem 27.9 could be regarded as one of the most general conclusions on solvability
of linear matrix equations up to now.   

Just as what  we did in Chapter 21, the two results (27.6) and 
(26.7) can be used to establish various types of rank equalities for generalized inverses of matrices. We next 
list some of them. 

\medskip
  
\noindent {\bf Theorem  27.10.}\, {\em Let $ A \in  {\cal F}^{m \times k}$ and  $ B \in  {\cal F}^{l \times m}$ be given. Then

{\rm (a)}\, The maximal and the minimal ranks of $ AA^- + B^-B$  with
 respect to $ A^-$  and $ B^-$ are 
$$ \displaylines{ 
\hspace*{1.5cm}
\max_{A^-,\, B^-} r( \, AA^- + B^-B \, ) =
\min \{ \ m,  \ \ r(A) + r(B) \ \}, \hfill (27.26)
\cr
\hspace*{1.5cm}
\min_{A^-,\, B^-} r( \, AA^- + B^-B \, ) = r(A) + r(B) - r(BA). \hfill (27.27)
\cr}
$$ 

{\rm (b)}\, There are $A^-$  and $B^-$ such that $ AA^- + B^-B$ is
nonsingular if and only if $r(A) + r(B) \geq m$.  

{\rm (c)} \, The rank of $ AA^- + B^-B$ is invariant with respect to the
choice of $A^- $  and $B^-$ if and only if $ BA = 0 $ or $r(BA) = r(A) + r(B) - m$. 

{\rm (d)} \, The rank of $ AA^- + B^-B$ is invariant with respect to the
choice of $A^- $  and $B^-$ if and only if the rank of $ AA^- - B^-B$ is invariant with respect to the choice of $A^- $  and $B^-$. 
} 

\medskip

\noindent {\bf Proof.}\, Note that 
$$ 
AA^- + B^-B = AA^{\sim} + B^{\sim}B + AV_1E_A + F_BV_2B. 
$$
This is a matrix expression with two independent variant matrices. Applying 
Theorem 27.2 to it and simplify  we can get Part (a). The detailed is omitted here.  Parts (b) and (c) are 
direct  consequences of Part (a). Contracting Part (c) and 
Theorem 21.16(d) we get Part (d).  \qquad $\Box$

\medskip 

Applying (27.27), (1.11) and (1.12), we can get the following 
$$ 
\displaylines{ 
\hspace*{0.5cm}
\min_{A^-,\, (I_m - A)^-} r[ \, AA^- + (I_m - A)^-(I_m - A) \, ] = r(A) + r(I_m - A) - r(A - A^2) = m, 
\hfill(27.28)
\cr
 \hspace*{0.5cm}
\min_{(I_m + A)^-,\, (I_m - A)^-} r[ \, (I_m + A)(I_m + A)^- + (I_m - A)^-(I_m - A) \, ] 
= r(I_m + A) + r(I_m - A) - r(I_m - A^2) = m, 
\cr
\hfill (27.29)
\cr}
$$
which imply that the matrices $ AA^- + (I_m - A)^-(I_m - A)$ and $ (I_m + A)(I_m + A)^- + 
(I_m - A)^-(I_m - A)$ are nonsingular for any $ A^-, \ (I_m - A)^-$  and  $ (I_m +  A)^-$. 

\medskip

By the similar approach, we can obtain the following.

\medskip
  
\noindent {\bf Theorem  27.11.}\, {\em Let $ A \in  {\cal F}^{m \times n}$ and  $ B \in  {\cal F}^{m \times k}$ be 
given. Then

{\rm (a)}\, The maximal and the minimal ranks of $ AA^- + BB^-$  with
 respect to $ A^-$  and $ B^-$ are 
$$ \displaylines{ 
\hspace*{1.5cm}
\max_{A^-,\, B^-} r( \, AA^- + BB^- \, ) = r[\, A, \ B \,], \hfill  (27.30)
\cr
\hspace*{1.5cm}
\min_{A^-,\, B^-} r( \, AA^- + BB^- \, ) =  \max \{ \, r(A), \ \ \  r(B) \, \}. \hfill (27.31)
\cr}
$$ 

{\rm (b)}\, The maximal and the minimal ranks of $ AA^- - BB^-$  with
 respect to $ A^-$  and $ B^-$ are 
$$ \displaylines{ 
\hspace*{1.5cm}
\max_{A^-,\, B^-} r( \, AA^-  -BB^- \, ) = \min \{ \, r[\, A, \ B \,], \ \  r[\, A, \ B \,] + m - r(A) - r(B) \, \}, 
\hfill (27.32)
\cr
\hspace*{1.5cm}
\min_{A^-,\, B^-} r( \, AA^- - BB^- \, ) =  \max \{ \,  r[\, A, \ B \,] - r(A), \ \ \   r[\, A, \ B \,] - r(B) \, \}. \hfill (27.33)
\cr}
$$ 

{\rm (c)}\, There are $A^-$  and $B^-$ such that $ AA^- = BB^-$  if and only if $R(A) = R(B)$.  
} 

\medskip

The rank equality in (27.30) can be extended to  
$$ 
\max_{A^-_1,\, \cdots, \, A_k^-} r( \, A_1A^-_1 + \cdots + A_kA_k^- \, ) = r[\, A_1, \, \cdots, \,  A_k \,]. 
\eqno (27.34)
$$
Notice that $ A_1A^-_1 + \cdots + A_kA_k^-$ is in fact a  matrix expression with $k$ independent variant matrices. 
Hence we have no rank formula at present for determining the minimal rank of $ A_1A^-_1 + \cdots + A_kA_k^-$. 
Nevertheless, we can guess from (27.31) the following 
$$ 
\min_{A^-_1,\, \cdots, \, A_k^-} r( \, A_1A^-_1 + \cdots + A_kA_k^- \, ) = \max \{\, r(A_1), \ \cdots, \ r(A_k) \, \}. \eqno (27.35)
$$

\medskip
  
\noindent {\bf Theorem  27.12.}\, {\em Let $ A \in  {\cal F}^{n \times m}$ and  $ B \in  {\cal F}^{k \times m}$ be 
given. Then
$$
 \displaylines{ 
\hspace*{1.5cm}
\min_{A^-,\, B^-} r[ \, A^-, \  B^- \, ] 
= \min_{A^-,\, B^-} r[ \, A^-A, \  B^-B \, ] = \max \{ \, r(A), \ \ \ r(B) \, \}. \hfill  (27.36)
\cr}
$$}

In general, we can guess from (27.36) the following  
$$ 
\min_{A^-_1,\, \cdots, \, A_k^-}r[ \, A^-_1, \, \cdots, \, A_k^- \, ]
 = \min_{A^-_1,\, \cdots, \, A_k^-} r[ \, A_1^-A_1, \, \cdots, \, A_k^-A_k \,]
 = \max \{ \, r(A_1), \ \cdots, \ r(A_k) \, \}. \eqno (27.37)
$$ 

\noindent {\bf Theorem  27.13.}\, {\em Let $ A \in  {\cal F}^{m \times n}, \  
B \in  {\cal F}^{m \times k}$ and $ C \in {\cal F}^{l \times n}$ be given. 
Then
$$ \displaylines{ 
\hspace*{1cm}
\min_{B^-,\, C^-} r( \, A -  BB^-A - AC^-C \, )  \hfill (27.38)
\cr
\hspace*{1cm}
= \max \left\{ r \left[ \begin{array}{cc}  A & B \\ C & 0 \end{array}\right] - r(B) - 
r(C),  \ \  r \left[ \begin{array}{cc}  A & B \\ C & 0 \end{array}
 \right] + r(A) - r[\, A, \ B\,] - r \left[ \begin{array}{c} A  \\ C \end{array}\right] \right\}. \hfill
(27.39)
\cr}
$$ 
In particular$,$ there are $B^-$  and $C^-$ such that $ BB^-A + AC^-C = A,$ i.e.$,$ the matrix equation 
$ BX + YC = A$ has a a solution with the form $ X = B^-A$ and  $ Y = AC^-,$ if and only if
$$
r \left[ \begin{array}{cc}  A & B \\ C & 0 \end{array}\right] =  r(B)+  
r(C) =  r \left[ \begin{array}{c} A  \\ C \end{array}\right] +
 r[\, A, \ B\,] - r(A).  \eqno (27.40)
$$
} 

Moreover, one can also find extreme ranks of matrix expressions 
$ A^kA^- + B^-B^k$, $ A^kA^- \pm B^kB^-$, $ A^-A^k \pm B^-B^k$, $ A - BB^- \pm C^-C$, and so on. The reader can try them and 
establish some more general results.      

\markboth{YONGGE  TIAN }
{28. EXTREME RANKS OF $ A - B_1XC_1$ SUBJECT TO  $ B_2XC_2 = A_2$ and $ B_3XC_3 = A_3$} 

\chapter{Extreme ranks of $ A - B_1XC_1$ subject to $ B_2XC_2 = A_2$ and $ B_3XC_3 = A_3$}

\noindent This chapter considers extreme ranks of the matrix 
expression $ A - B_1XC_1$ subject to a pair of
 consistent matrix equations $ B_2XC_2 = A_2$ and $ B_3XC_3 = A_3$ over an 
arbitrary filed ${\cal F}$. A direct
 motivation for this work comes from considering consistency of the
 triple matrix equations $ B_1XC_1 = A_1, \ B_2XC_2 = A_2$ and
$ B_3XC_3 = A_3.$ To do so,  we need to know expression of general solution to the pair of matrix equations 
$ B_2XC_2 = A_1, \ B_3XC_3 = A_3$.   

\medskip

\noindent {\bf Lemma 28.1.}\, {\em Suppose that 
$$\displaylines{
\hspace*{2cm} 
 B_2XC_2 = A_2, \qquad  B_3XC_3 = A_3 \hfill (28.1)
\cr}
$$
is a pair of matrix equations over an arbitrary filed ${\cal F}$. Then

{\rm (a)}\, The general common solution of the pair of homogeneous matrix 
equations $B_2 XC_2 = 0$ and $B_3 XC_3= 0$  can factor as 
$$\displaylines{
\hspace*{2cm} 
            X = X_1 + X_2 + X_3 + X_4,              \hfill (28.2)
\cr}
$$
where $ X_1, \  X_2, \  X_3 $ and $ X_4$ are$,$ respectively$,$ the general
solutions of the following four systems of homogeneous linear matrix equations
$$\displaylines{
\hspace*{2cm} 
\left\{ \begin{array}{l}  B_2 X_1 = 0 \\  B_3 X_1 = 0, \end{array} \right.\ \ \ 
\left\{ \begin{array}{l}   X_2C_2 = 0 \\   X_2C_3 = 0, \end{array} \right. \ \ \  
\left\{ \begin{array}{l}  B_2 X_3 = 0 \\  X_3C_3= 0, \end{array} \right. \ \ \ 
\left\{ \begin{array}{l}  X_4C_2 = 0 \\  B_3X_4 = 0. \end{array} \right. \ \ \   \hfill (28.3)
\cr}
$$
Written in an explicit form$,$  it is 
$$\displaylines{
\hspace*{2cm} 
  X =  F_B V_1 + V_2 E_C + F_{B_2}V_3 E_{C_3}  + F_{B_3}V_4 E_{C_2},
 \hfill (28.4)
\cr}
$$  
where  $ B = \left[ \begin{array}{c} B_2  \\ B_3  \end{array} \right], \ C = [ \, C_2, \ C_3 \, ],$ and $ V_1$---$V_4$ are four arbitrary matrices.  

{\rm (b)}\, Suppose that the pair of matrix equations $ B_2 XC_2 = A_2$ and 
$ B_3 XC_3 = A_3$ have a common solution. Then the general common solution can be written as
$$\displaylines{
\hspace*{2cm} 
  X = X_0 +  F_B V_1 + V_2 E_C + F_{B_2}V_3 E_{C_3}  + F_{B_3}V_4 E_{C_2}, \hfill (28.5)
\cr}
$$  
where $X_0$ is a particular common solution to $ B_2 XC_2 = A_2$ and 
$ B_3 XC_3 = A_3$.

{\rm (c)}\, Suppose that the pair of matrix equations $ B_2 XC_2 = A_2$ and 
$ B_3 XC_3 = A_3$ have a common solution$,$ and the given matrices satisfy   
$$ 
\displaylines{
\hspace*{2cm} 
 R(B_2^T) \subseteq R(B_3^T),  \qquad  R(C_3) \subseteq R(C_2),  \hfill 
\cr}
$$   
or equivalently 
$$ 
\displaylines{
\hspace*{2cm} 
   R(F_{B_3}) \subseteq R(F_{B_2}),  \qquad    R( E_{C_2}^T ) \subseteq R( E_{C_3}^T ).  \hfill
\cr}
$$   
Then the general common solution  $ B_2 XC_2 = A_2$ and $ B_3 XC_3 = A_3$ can be written as
$$
\displaylines{
\hspace*{2cm} 
  X = X_0 +  F_{B_3} V_1 + V_2 E_{C_2} + F_{B_2}V_3 E_{C_3}, \hfill
\cr}
$$  
where $X_0$ is a particular common solution to $ B_2 XC_2 = A_2$ and 
$ B_3 XC_3 = A_3$.}

\medskip
   
\noindent {\bf Proof.}\, According to Lemma 18.1, the general solution of 
$ B_2 XC_2= 0$  can be written as 
$$\displaylines{
\hspace*{2cm} 
             X = F_{B_2}W_1 + W_2 E_{C_2},         \hfill (28.6)
\cr}
$$
where $ W_1, \  W_2$ are arbitrary. Substituting it into $B_3XC_3 = 0$  yields 
$$\displaylines{
\hspace*{2cm} 
         B_3 XC_3 = B_3 F_{B_2} W_1C_3 + B_3 W_2 E_{C_2}C_3 = 0.          \hfill (28.7)
\cr}
$$
Observe that $ R( B_3 F_{B_2} ) \subseteq R( B_3 )$ and $ R[(E_{C_2}C_3)^T]
\subseteq R( C_3^T ).$ We can find by Lemma 26.2(a) that the general solutions for $W_1 $ and $W_2$ of (28.7) 
can be written as 
$$\displaylines{
\hspace*{2cm} 
  W_1 = U E_{C_2} + F_G V_1 + V_3E_{C_3}, \hfill
\cr
\hspace*{2cm} 
 W_2 = -F_{B_2}U + V_2 E_H +F_{B_3} V_4, \hfill
\cr}
$$
where $ H = B_3 F_{B_2}, \ G =  E_{C_2}C_3$, and $U, \  V_1$---$V_4$ are arbitrary.  Substituting both of them 
into (28.6) produces the general common solution of  $ B_2 XC_2 = 0 $ and 
$ B_3 XC_3 = 0 $ as follows
$$\displaylines{
\hspace*{2cm} 
 X = F_{B_2}F_G V_1 + V_2 E_HE_{C_2}  + F_{B_2}V_3 E_{C_3}  + F_{B_3}V_4 E_{C_2}.   \hfill (28.8)
\cr}
$$ 
It is easy to verify that the four terms in (28.8)  are, in turn, the general common solutions of the four 
 pairs of homogeneous equations in (28.3). Thus we have (28.2) and (28.4).  The result in Part (b) is obvious 
from Part (a).   \qquad  $ \Box $  

\medskip

 Putting (28.5) in $A_1 - B_1XC_1,$ we get 
$$
A_1 - B_1XC_1 = A_1 - B_1X_0C_1 - B_1F_B V_1C_1 - B_1V_2 E_CC_1 -  B_1F_{B_2}V_3 E_{C_3}C_1 - 
B_1F_{B_3}V_4 E_{C_2}C_1. \eqno (28.9) 
$$
Thus the maximal and the minimal ranks of the matrix 
expression $ A - B_1XC_1$ subject to $ B_2XC_2 = A_2$ and $ B_3XC_3 = A_3$ can be determined by the  matrix expression
 (28.9).  For convenience of representation, we write (28.9) as 
$$\displaylines{
\hspace*{2cm} 
A_1 - B_1XC_1 = A - G_1V_1H_1 - G_1V_2H_2 - G_3V_3H_3 - G_4V_4H_4, \hfill (28.10)
\cr
\hspace*{0cm}
where  \hfill
\cr
\hspace*{2cm} 
 A = A_1 - B_1X_0C_1, \ \ \
G_1 =  B_1F_B, \ \ \  G_2 =  B_1, \ \ \ G_3 = B_1F_{B_2},   \ \ \ G_4 = B_1F_{B_3}, \hfill (28.11)
\cr
\hspace*{2cm} 
H_1 = C_1 ,\ \ \ H_2 = E_CC_1, \ \ \ H_3 =  E_{C_3}C_1, \ \ \ H_4 =  E_{C_2}C_1.  \hfill (28.12)
\cr}
$$

Observe that (28.10) involves four independent variant matrices $ V_1$---$V_4$. Moreover it is not difficult to derive
 that the above matrices satisfy the following conditions 
$$\displaylines{
\hspace*{2cm} 
R(G_1) \subseteq R(G_i) \subseteq R(G_2),  \ \ {\rm and}   \ \  R(H_2^T) \subseteq R(H_i^T) \subseteq R(H_1^T), \ \ \ \ 
i = \ 3, \ 4, \hfill (28.13)
\cr}
$$ 
Thus (28.10) can be regarded as a special case of the matrix expression in Theorem 27.4. In that case, 
applying the two formulas in (27.10) and (27.11) to (28.10), we get  the main results of the chapter. 

\medskip

\noindent {\bf Theorem 28.2.}\, {\em Suppose that the pair of matrix equations
 $ B_2 XC_2 = A_2$ and $ B_3 XC_3 = A_3$ have a common solution. Then the maximal rank of $ A_1 - B_1XC_1$ subject to 
$ A_2 XB_2 = C_2$ and $ A_3 XB_3 = C_3$  is 
$$
\displaylines{
\hspace*{1cm}
\max_{\begin{array}{c} B_2 XC_2 = A_2 \\ B_3 XC_3 = A_3  \end{array}} r( \, A_1 - B_1XC_1 \,)
 = \min  \left\{  r[ \, A_1, \ B_1 \, ],    \ \ r\left[ \begin{array}{c} A_1 \\ C_1 \end{array} \right], \ \ 
s_1, \ \ s_2, \ \ s_3, \ \ s_4  \right\}, \hfill (28.14)
\cr
\hspace*{0cm}
where \hfill
\cr
\hspace*{1cm} 
s_1 =  r\left[ \begin{array}{cccc} A_1  & 0 & 0 & B_1 \\
 0 & -A_2  & 0 & B_2  \\ 0 & 0 & -A_3 & B_3 \\ C_1   & C_2  & 0 & 0  \\
 C_1 & 0 & C_3 & 0 \end{array} \right]
 - r\left[ \begin{array}{cc} B_2 \\
 B_3 \end{array} \right] - r(C_2) - r(C_3), \hfill
\cr
\hspace*{1cm}
s_2 = r\left[ \begin{array}{ccccc} A_1  & 0 & 0 & B_1 & B_1 \\
 0 & -A_2  & 0 & B_2 & 0  \\ 0 & 0 & -A_3  & 0 & B_3 \\ C_1   & C_2  & C_3
  & 0 & 0 \end{array} \right] -r[\, C_2, \ C_4 \,] -r(B_2) - r(B_3),
  \hfill
\cr
\hspace*{1cm}
s_3 = r\left[ \begin{array}{ccc}  A_1  & 0  & B_1  \\ 0 &  -A_2  & B_2  \\
C_1 & C_2 & 0 \end{array} \right] - r(B_2) -r(C_2),  \ \ \ s_4 = r \left[ \begin{array}{ccc} A_1  & 0  & B_1  \\ 0 &  -A_3  & B_3  \\
C_1 & C_3 & 0 \end{array} \right] - r(B_3) - r(C_3). \hfill
\cr}
$$
}
{\bf Proof.}\, Under (28.13), we first find  by (27.10) that
$$
\displaylines{
\hspace*{0cm} 
\max_{\begin{array}{c} B_2 XC_2 = A_2 \\ B_3 XC_3 = A_3  \end{array}} r( \,
A_1 - B_1XC_1 \,) = \max_{\{V_i\}} r( \, A - G_1V_1H_1 - G_1V_2H_2 - G_3V_3H_3 - G_4V_4H_4 \, )
\hfill
\cr}
$$ 
$$
= \min \left\{  r[ \, A, \  G_2 \, ], \   r \left[ \begin{array}{c}  A  \\
H_1 \end{array} \right], \ 
r \left[ \begin{array}{cc}  A  & G_1 \\ H_3 & 0 \\ H_4 & 0  \end{array}
\right],  \  r \left[ \begin{array}{ccc} A  & G_3 & G_4 \\ H_2  & 0 & 0
\end{array} \right],  \  r \left[ \begin{array}{cc}  A  & G_3 \\ H_4 & 0  \end{array} \right],  \ \
r \left[ \begin{array}{cc}  A & G_4 \\ H_3 & 0 \end{array} \right] \right\}. \eqno (28.15)
$$
Simplifying the ranks of the  block matrices in (28.15) by (1.2)---(1.4), as well as $ B_2X_0C_2 = A_2,$ 
$ B_3X_0C_3 = A_3$,  we have 
$$
\displaylines{
\hspace*{0.5cm} 
r[ \, A, \  G_2 \, ] = r[ \, A_1 - B_1X_0C_1 , \  B_1  \, ] =
r[ \, A_1, \ B_1 \, ],  \ \  \ \  
r \left[ \begin{array}{c}  A  \\ H_1 \end{array} \right] = 
r \left[ \begin{array}{c}  A_1 - B_1X_0C_1  \\ C_1 \end{array} \right] 
 = r \left[ \begin{array}{c}  A_1 \\ C_1 \end{array} \right], \hfill
\cr
\hspace*{1cm}
r \left[ \begin{array}{cc} A  & G_1 \\ H_3 & 0 \\ H_4 & 0  \end{array}
\right] = r \left[ \begin{array}{cc} A_1 - B_1X_0C_1  & B_1F_B \\ E_{C_3}C_1
& 0 \\  E_{C_2}C_1 & 0 \end{array} \right] \hfill
\cr
\hspace*{1cm}
= r \left[ \begin{array}{cccc} A_1 - B_1X_0C_1  & B_1 & 0 & 0
\\ C_1 & 0  & C_3 & 0 \\  C_1 & 0 & 0 &C_2 \\
 0 & B_2 & 0 & 0 \\ 0 & B_3 & 0 & 0 \end{array} \right] -
 r\left[ \begin{array}{cc} B_2 \\ B_3 \end{array} \right] - r(C_2) - r(C_3)
\hfill
\cr
\hspace*{1cm}
= r \left[ \begin{array}{cccc} A_1 & B_1 & 0 & 0
\\ C_1 & 0  & C_3 & 0 \\  C_1 & 0 & 0 & C_2 \\
 0 & B_2 & 0 & -A_2 \\ 0 & B_3 & -A_3 & 0 \end{array} \right] -
 r\left[ \begin{array}{cc} B_2 \\ B_3 \end{array} \right] - r(C_2)
 - r(C_3)
 \hfill
\cr
\hspace*{1cm}
 =  r\left[ \begin{array}{cccc} A_1  & 0 & 0 & B_1 \\
 0 & -A_2  & 0 & B_2  \\ 0 & 0 & -A_3 & B_3 \\ C_1   & C_2  & 0 & 0  \\
 C_1 & 0 & C_3 & 0 \end{array} \right]
 - r\left[ \begin{array}{cc} B_2 \\
 B_3 \end{array} \right] - r(C_2) - r(C_3). \hfill
\cr}
$$
Similarly we can get 
$$
\displaylines{
\hspace*{1cm}
r \left[ \begin{array}{ccc} A  & G_3 & G_4 \\ H_2  & 0 & 0
\end{array} \right] = r\left[ \begin{array}{ccccc} A_1  & 0 & 0 & B_1 &
B_1 \\ 0 & -A_2  & 0  & B_2 & 0  \\ 0 & 0 & -A_3  & 0 & B_3 \\ C_1
& C_2  & C_3  & 0 & 0 \end{array} \right] - r[\, C_2, \ C_4 \,] -r(B_2)
- r(B_3), \hfill
\cr
\hspace*{1cm}
 r \left[ \begin{array}{cc}  A  & G_3 \\ H_4 & 0  \end{array}
 \right] = r\left[ \begin{array}{ccc}  A_1  & 0  & B_1  \\ 0 &  -A_2
 & B_2  \\ C_1 & C_2 & 0 \end{array} \right] - r(B_2) -r(C_2), \hfill
\cr
\hspace*{1cm}
r \left[ \begin{array}{cc}  A & G_4 \\ H_3 & 0 \end{array} \right]
= r \left[ \begin{array}{ccc} A_1  & 0  & B_1  \\ 0 &  -A_3  & B_3  \\
C_1 & C_3 & 0 \end{array} \right] - r(B_3) - r(C_3).  \hfill
\cr}
$$
Putting them in (28.15) yields (28.14). \qquad $ \Box$ 

\medskip

\noindent {\bf Theorem 28.3.}\, {\em Suppose that the pair of matrix
equations $ B_2 XC_2 = A_2$ and $ B_3 XC_3 = A_3$ have a common solution.
Then the minimal rank of $ A_1 - B_1XC_1$ subject to
$ B_2 XC_2 = A_2$ and $ B_3 XC_3 = A_3$  is 
$$
\displaylines{ 
\hspace*{0cm}
\min_{\begin{array}{c}  B_2 XC_2 = A_2 \\ B_3 XC_3 = A_3  \end{array}}
r( \, A_1 - B_1XC_1 \,) =  r\left[ \begin{array}{cccc} A_1  & 0 & 0 & B_1 \\
 0 & -A_2  & 0 & B_2  \\ 0 & 0 & -A_3 & B_3 \\ C_1   & C_2  & 0 & 0  \\ C_1
 & 0 & C_3 & 0 \end{array} \right] +  r\left[ \begin{array}{ccccc} A_1  & 0
 & 0 & B_1 & B_1 \\ 0 & -A_2   & 0 & B_2 & 0  \\ 0 & 0 & -A_3  & 0 & B_3 \\ C_1  
 & C_2  & C_3  & 0 & 0 \end{array} \right] \hfill
\cr
\hspace*{0.5cm} 
- \ r\left[ \begin{array}{cc} A & B_1 \\ C_1 & 0 \\ 0 & B_2 \\  0 &  B_3
\end{array} \right] - r\left[ \begin{array}{cccc} A & B_1 & 0 & 0 \\ C_1   & 0 & C_2 & C_3 \end{array} \right] + 
\left[ \begin{array}{cc} A_1 \\ C_1 \end{array} \right] + r[\, A_1,  \ B_1 \,]  \hfill
\cr
\hspace*{0.5cm}
+ \max \left\{ r\left[ \begin{array}{ccc}  A_1  & 0  & B_1  \\ 0 &  -A_2  & B_2  \\
C_1 & C_2 & 0 \end{array} \right] - r\left[ \begin{array}{cccc}  A_1  & 0  & B_1 & B_1  
\\ 0 &  -A_2  & B_2 & 0   \\
C_1 & C_2 & 0 & 0 \\ 0 & 0 & 0 & B_3 \end{array} \right]
 -  
 r\left[ \begin{array}{cccc}  A_1  & 0  & B_1 & 0   \\ 0 &  -A_2  & B_2  & 0  \\
C_1 & C_2 & 0  & 0  \\ C_1 & 0 & 0 & C_3 \end{array} \right],  \right.
\hfill
\cr
\hspace*{1.5cm}
\left.  r \left[ \begin{array}{ccc} A_1  & 0  & B_1  \\ 0 &  -A_3  & B_3  \\
C_1 & C_3 & 0 \end{array} \right] - r\left[ \begin{array}{cccc}  A_1  & 0  & B_1 & B_1  
\\ 0 &  -A_3  & B_3 & 0   \\
C_1 & C_3 & 0 & 0 \\ 0 & 0 & 0 & B_2 \end{array} \right] - 
 r\left[ \begin{array}{cccc}  A_1  & 0  & B_1 & 0   \\ 0 &  -A_3  & B_3  & 0  \\
C_1 & C_3 & 0  & 0  \\ C_1 & 0 & 0 & C_2 \end{array} \right] \right\}.
\hfill (28.16)
\cr}
$$
}
{\bf Proof.}\, Under (28.13), applying (27.11) to (28.10) yields
$$
\displaylines{
\hspace*{0cm} 
\min_{\begin{array}{c} B_2 XC_2 = A_2 \\ B_3 XC_3 = A_3 \end{array}} r( \,
A_1 - B_1XC_1 \,) \hfill
\cr
\hspace*{0cm} 
= \min_{\{V_i\}} r( \, A - G_1V_1H_1 - G_1V_2H_2 - G_3V_3H_3 - G_4V_4H_4 \, )
\hfill
\cr 
\hspace*{0cm}
= r[ \, A, \  G_2 \, ] +  r \left[ \begin{array}{c}  A  \\
H_1 \end{array} \right] + r \left[ \begin{array}{cc}  A  & G_1 \\ H_3 & 0 \\
H_4 & 0  \end{array}
\right] + r \left[ \begin{array}{ccc} A  & G_3 & G_4 \\ H_2  & 0 & 0
 \end{array} \right] + r[ \, A, \  G_2 \, ] +  r \left[ \begin{array}{c}
  A  \\ H_1 \end{array} \right]  - \ r \left[ \begin{array}{cc}  A  & G_1 \\
  H_1 & 0 \end{array} \right] \hfill
\cr 
\hspace*{0.5cm}
-  \ r \left[ \begin{array}{cc}  A  & G_2 \\
  H_2 & 0  \end{array} \right] + \max \left\{ r\left[ \begin{array}{cc} A & G_3 \\ H_4 & 0
\end{array} \right] - r \left[ \begin{array}{ccc}  A & G_3 & G_4 \\
H_3 & 0 & 0 \end{array} \right] - r \left[ \begin{array}{cc}  A  & G_3 \\
 H_3 & 0 \\ H_4 & 0 \end{array} \right], \right. \hfill
\cr 
\hspace*{4.5cm} 
\left.   \ r \left[ \begin{array}{cc} A & G_4 \\ H_3 & 0
\end{array} \right] - r \left[ \begin{array}{ccc}  A & G_3 & G_4 \\
H_3 & 0 & 0 \end{array} \right] - r \left[ \begin{array}{cc}  A  & G_4 \\
 H_3 & 0 \\ H_4 & 0 \end{array} \right] \right\}.  \hfill  (28.17)
\cr}
$$
Simplifying the ranks of block matrices in (28.17) by (1.2)---(1.4), as well as $ B_2X_0C_2 = A_2, \ B_3X_0C_3 
= A_3$,  we can eventually get the rank formula (28.16). But we omit here the tedious steps. \qquad $ \Box $ 

\medskip

\noindent {\bf Corollary 28.4.}\, {\em Suppose that the three matrix equations   $ B_1 XC_1 = A_1,$ 
 $ B_2 XC_2 = A_2$ and $ B_3 XC_3 = A_3$ are consistent$,$ respectively. Also suppose that any pair 
of the  three matrix equations has a common solution.  Then 
$$
\displaylines{ 
\hspace*{0cm}
\min_{\begin{array}{c}  B_2 XC_2 = A_2 \\ B_3 XC_3 = A_3  \end{array}}
r( \, A_1 - B_1XC_1 \,) =  r\left[ \begin{array}{cccc} A_1  & 0 & 0 & B_1 \\
 0 & -A_2  & 0 & B_2  \\ 0 & 0 & -A_3 & B_3 \\ C_1   & C_2  & 0 & 0  \\ C_1 & 0 & C_3 & 0 
\end{array} \right] +  r\left[ \begin{array}{ccccc} A_1  & 0 & 0 & B_1 & B_1 \\
 0 & -A_2  &  0 & B_2  & 0 \\ 0 & 0 & -A_3 & 0 & B_3 \\ C_1 & C_2  & C_3 & 0
 & 0 \end{array} \right] \hfill
\cr
\hspace*{3.5cm} 
- \ r\left[ \begin{array}{cc} B_1 & B_1 \\ B_2 & 0 \\ 0 & B_3 \end{array} \right] 
- r\left[ \begin{array}{cc} B_1 \\  B_2 \\  B_3 \end{array} \right]
- r\left[ \begin{array}{cccc} C_1 & C_2 & 0  \\ C_1   & 0 & C_3  \end{array} \right] 
- r[\, C_1,  \ C_2, \ C_3 \,].  \hfill (28.18)
\cr}
$$}
{\bf Proof.}\, Under the assumption of the corollary, we know by
Corollary 20.3 that the given matrices in the three equations satisfy the
conditions
$$
R(A_i) \subseteq R(B_i),  \ \ \  R(A_i^T) \subseteq R(C_i^T),  \ \ \ i = 1,  \ 2, \ 3,.
$$
$$
r\left[ \begin{array}{ccc}  A_i  & 0  & B_i  \\ 0 &  -A_j  & B_j  \\
C_i & C_j & 0 \end{array} \right]  = r \left[ \begin{array}{c} B_i  \\ B_j \end{array} \right] +
r[ \, C_i, \  C_j \,].  \ \ \ i = 1,  \ 2, \ 3. 
$$
In that case, the formula (28.16) reduces to (28.18). \qquad $\Box$ 

\medskip

Based on the formula (28.18), one  can easily verify that under the assumption
of Corollary 28.4, the following identity holds
$$
\min_{\begin{array}{c} B_2 XC_2 = A_2 \\ B_3 XC_3 = A_3  \end{array}} r( \, A_1 - B_1XC_1 \,) 
=\min_{\begin{array}{c} B_1 XC_1 = A_1 \\ B_3 XC_3 = A_3  \end{array}} r( \, A_2 - B_2XC_2 \,) 
= \min_{\begin{array}{c} B_1 XC_1 = A_1 \\ B_2 XC_2 = A_2  \end{array}} r( \, A_3 - B_3XC_3 \,).  
$$
\hspace*{0.3cm} One of the most important consequences of (28.18) is concerning the consistency of a triple matrix equations.  

\medskip

\noindent {\bf Corollary 28.5.}\, {\em The triple linear matrix equations  $ B_1 XC_1 = A_1,$ 
 $ B_2 XC_2 = A_2$ and $ B_3 XC_3 = A_3$ have a common solution if and only if  any pair 
 of the three equations has a common solution$,$ meanwhile the given matrices  satisfy the two rank equalities  
$$
\displaylines{ 
\hspace*{1cm}
r\left[ \begin{array}{ccccc} A_1  & 0 & 0 & B_1 & B_1 \\
 0 & -A_2 & 0 &  B_2 & 0 \\ 0 & 0 & -A_3  & 0 & B_3 \\ C_1   & C_2  & C_3  & 0 & 0 
\end{array} \right] =  r\left[ \begin{array}{cc} B_1 & B_1 \\ B_2 & 0 \\ 0 & 
B_3 \end{array} \right] +  r[\, C_1,  \ C_2, \ C_3 \,], \hfill (28.19)
\cr
\hspace*{1cm} 
r\left[ \begin{array}{cccc} A_1  & 0 & 0 & B_1 \\
 0 & -A_2  & 0 & B_2  \\ 0 & 0 & -A_3 & B_3 \\ C_1   & C_2  & 0 & 0  \\ C_1 & 0 & C_3 & 0 
\end{array} \right] = r\left[ \begin{array}{ccc} C_1 & C_2 & 0  \\ C_1   & 0 & C_3  \end{array} \right] + r\left[\begin{array}{c} B_1 \\  B_2 \\  B_3 \end{array} \right]. \hfill
 (28.20)
\cr}
$$}
\hspace*{0.3cm} This result can also be alternatively stated as follows.

\medskip

\noindent {\bf Corollary 28.6.}\, {\em The triple matrix equations  $ B_1 XC_1 = A_1,$ 
 $ B_2 XC_2 = A_2$ and $ B_3 XC_3 = A_3$ have a common solution if and only if  the following eight 
independent simple matrix equations are all solvable
$$
B_1 X_1C_1 = A_1,  \ \ \ B_2 X_2C_2 = A_2,  \ \ \  B_3 X_3C_3 = A_3, 
$$
$$
\left[ \begin{array}{c} B_1  \\ B_2  \end{array} \right] X_4 + Y_4[\, C_1, \ C_2 \, ]
 = \left[ \begin{array}{ccc} A_1 & 0  \\ 0  & -A_2  \end{array} \right],  \ \ \
\left[ \begin{array}{c} B_1  \\ B_3  \end{array} \right] X_5 + Y_5[\, C_1, \ C_3 \, ]
 = \left[ \begin{array}{ccc} A_1 & 0  \\ 0  & -A_3  \end{array} \right],  
$$
$$
\left[ \begin{array}{c} B_2  \\ B_3  \end{array} \right] X_6 + Y_6[\, C_2, \ C_3 \, ]
 = \left[ \begin{array}{ccc} A_2 & 0  \\ 0  & -A_3  \end{array} \right], 
$$
$$
\left[ \begin{array}{cc} B_1  & B_1 \\ B_2 & 0 \\ 0 & B_3  \end{array} \right] X_7 + 
Y_7[\, C_1, \ C_2, \ C_3 \, ]
 = \left[ \begin{array}{ccc} A_1 & 0   & 0 \\ 0  & -A_2  & 0 \\ 0 & 0 & -A_3 \end{array} \right], 
$$
$$ 
\left[ \begin{array}{c} B_1 \\ B_2  \\ B_3  \end{array} \right] X_8 + 
Y_8\left[ \begin{array}{ccc}  C_1  & C_2  & 0 \\ C_1 & 0 & C_3 \end{array} \right] 
= \left[ \begin{array}{ccc} A_1 & 0   & 0 \\ 0  & -A_2  & 0 \\ 0 & 0 & -A_3 \end{array} \right]. 
$$}
\hspace*{0.3cm} Of course, one can also equivalently write the consistency condition for $ B_1 XC_1 = A_1,$ 
 $ B_2 XC_2 = A_2$ and $ B_3 XC_3 = A_3$ in Theorems 28.5 and 28.6 in term of equivalence of matrices,
column or row spaces of matrices, generalized inverses of matrices, and so on.     

As a simple consequence  of Theorems 28.2 and 28.3 we can also get the maximal and  the minimal ranks of common 
solutions to a pair of linear matrix equations. This problem was examined by Mitra \cite{Mi6}. 

\medskip

\noindent {\bf  Corollary 28.7.}\, {\em Suppose that the pair of matrix
equations $ B_2 XC_2 = A_2$ and $ B_3 XC_3 = A_3$ have a common solution$,$
where $ X $ is a $ p \times q $ matrix.  Then the maximal rank of common
solutions to the pair of equations is
$$\displaylines{
\hspace*{2cm} 
\max_{\begin{array}{c}B_2 XC_2 = A_2 \\ B_3 XC_3 = A_3  \end{array}} r(X) 
 = \min \{ \, p, \ \ q, \ \  s_1, \ \ s_2, \ \ s_3, \ \ s_4 \, \}, \hfill (28.21) 
\cr
\hspace*{0cm}
where  \hfill
\cr
\hspace*{2cm} 
s_1 = r(A_2) - r(B_2) - r(C_2) + p + q,  \hfill
\cr
\hspace*{2cm}
 s_2 = r(A_3) - r(B_3) -
r(C_3) + p + q, \hfill
\cr
\hspace*{2cm} 
s_3 = r\left[ \begin{array}{cc} A_2 &  0 \\ 0 &  A_3 \\ C_2 & C_3 \end{array} \right] 
- r[\, C_2,  \  C_3 \,] - r(C_2) - r(C_3) + p + q, \hfill 
\cr
\hspace*{2cm} 
s_4 = r\left[ \begin{array}{ccc} A_2 &  0 & B_2  \\ 0 &  A_3  &  B_3
\end{array} \right] - r\left[ \begin{array}{c}  B_2 \\  B_3 \end{array}
\right] - r(B_2) - r(B_3) + p + q. \hfill
\cr}
$$
The minimal rank of common solutions to the pair of equations is
$$
\displaylines{ 
\hspace*{0cm}
\min_{\begin{array}{c}  B_2 XC_2 = A_2 \\ B_3 XC_3 = A_3  \end{array}} r(X)
= r\left[ \begin{array}{cc} A_2 &  0 \\ 0 &  A_3 \\ C_2 & C_3 \end{array} \right] 
+ r\left[ \begin{array}{ccc} A_2 &  0 & B_2  \\ 0 &  A_3  &  B_3 \end{array} \right] \hfill
\cr
\hspace*{0cm}
 +  \max \left\{ r(A_2) - r\left[ \begin{array}{cc} A_2 & B_2 \\ 0 & B_3 \end{array} \right]   
- r\left[ \begin{array}{cc} A_2 & 0  \\ C_2 & C_3  \end{array} \right],  \ \
r(A_3) - r\left[ \begin{array}{cc} B_2 & 0 \\ B_3 & A_3 \end{array} \right]   
- r\left[ \begin{array}{cc} C_2 & C_3  \\ 0 & A_3  \end{array} \right]
\right\}.  \hfill (28.22)
\cr}
$$}
\noindent{\bf Corollary 28.8.}\, {\em Let  $ A, \, B , \, C \in {\cal F}^{ m \times n}$ be given. Then 
 $ A, \, B $ and $C$ have a common inner inverse if and only if 
$$
\displaylines{ 
\hspace*{2cm}
r( \, A - B \, ) = r\left[ \begin{array}{c} A \\ B \end{array} \right] + r[ \, A, \  B \,] - r(A) - r(B),\hfill 
\cr
\hspace*{2cm}
r( \, A - C \, ) = r\left[ \begin{array}{c} A \\ C \end{array} \right] + r[ \, A, \ C \, ] - r(A) - r(C), \hfill
\cr 
\hspace*{2cm}
r( \, B - C \, ) = r\left[ \begin{array}{c} B \\ C \end{array} \right] + r[ \, B, \ C \,] - r(B) - r(C), \hfill
\cr
\hspace*{2cm}
r[ \, A - B ,\ A - C ] = r\left[ \begin{array}{cc} A & A \\ B & 0 \\ 0 & C \end{array} \right] + 
r[ \, A, \ B, \ C \,] - r(A) - r(B) - r(C),  \hfill
\cr
\hspace*{2cm}
r\left[ \begin{array}{cc} A - B \\ A - C  \end{array} \right] 
= r\left[ \begin{array}{c} A \\ B \\ C \end{array} \right] + 
 r\left[ \begin{array}{ccc} A & B & 0  \\ A & 0 & C \end{array} \right] - r(A) - r(B) - r(C). \hfill
\cr} 
$$
In particular$,$ if 
$$\displaylines{
\hspace*{2cm} 
r\left[ \begin{array}{c} A \\ B \\ C \end{array} \right]  = r[ \, A, \ B, \ C \,] = r(A) + r(B) + r(C),\hfill
\cr} 
$$
then $ A, \ B $ and $C$ have a common inner inverse. }

\medskip

\noindent {\bf Proof.}\, Consider the three matrix equations $ AXA = A, \ BXB = B $ and $ CXC = C$. Then the result in 
the corollary follows directly from Theorems 21.10(a) and 28.5.   \qquad $\Box$

\medskip

When the matrices  $ A, \ B $ and $C$ are all idempotent, they have identity matrix as their common inner inverse. Thus 
the five rank equalities in Corollary 28.7 are all satisfied, the first three occurred in Theorem 3.1, the fourth and 
the fifth are two new rank equalities for idempotent matrices.

Another work related to a  triple matrix equations   $ B_1 XC_1 = A_1,$  $ B_2 XC_2 = A_2$ and $ B_3 XC_3 = A_3$
 is to determine 
$$
\max_{\begin{array}{c}  B_1XC_1 = A_1  \\  B_2 YC_2 = A_2 \\ B_3 YC_3 = A_3  \end{array}} 
r(\, X - Y \, ),    \qquad   \ \ \ 
\min_{\begin{array}{c}  B_1XC_1 = A_1  \\  B_2 YC_2 = A_2 \\ B_3 YC_3 = A_3  \end{array}} 
r(\, X - Y \, ). \eqno (28.23)
$$
Based on Lemma 28.1 and Corollary 27.5, one can routinely find the two ranks in (28.23). From them one can also 
establish  Corollaries 28.5 and 28.6. We leave this work to the reader. 

A more general work than those for a triple matrix equations  is  to consider  common solution to 
a quadruple of matrix equations  
$$
 B_1 XC_1 = A_1,  \ \ \  B_2 XC_2 = A_2,  \ \ \  B_3 XC_3 = A_3, \ \ \   B_4 XC_4 = A_4. \eqno (28.24)
$$
Clearly, the  quadruple  matrix equations have a common solution if and only if 
$$ 
\min_{\begin{array}{c}  B_1XC_1 = A_1 \\  B_2 XC_2 = A_2 \\  B_3 YC_3 = A_3 \\ B_4 YC_4 = A_4 \end{array}} 
r(\, X - Y \, ) = 0. \eqno (28.25)
$$
If the two  pairs  $B_1 XC_1 = A_1, \ B_2 XC_2 = A_2$ and $B_3 YC_3 = A_3, \ B_4YC_4 = A_4$ are consistent 
respectively, then the deference $ X - Y $ in (28.25), according to Lemma 28.1(b), is a linear matrix expression with 
eight independent variant matrices, four of them are one-sided and other four are two-sided. Unforturnately we can not 
find in general the minimal rank of such a matrix expression. However if the quadruple  matrix equations  satisfy some
restrictions, and the expressions for $ X $ and $ Y $ are reduced to some simple forms, then we can find (28.25).    
One such a case is when (28.24) satisfy the conditions     
$$ 
 R( B_1^T ) \subseteq R( B_2^T ),  \qquad  R(C_2) \subseteq R( C_1 ), \qquad   R( B_3^T ) \subseteq R( B_4^T ),  \qquad R(C_4) \subseteq R( C_3 ).
 \eqno (28.26)
$$   
or equivalently 
$$ 
   R( F_{B_2} ) \subseteq R(F_{B_1}),  \qquad    R( E_{C_1}^T ) \subseteq R( E_{C_2}^T ), 
\qquad  R( F_{B_4} ) \subseteq R( F_{B_3} ),  \qquad R( E_{C_3}^T ) \subseteq R( E_{C_4}^T ).
 \eqno  (28.27)
$$   
In that case, the general common solution to $B_1XC_1 = A_1 $ and  $B_2 XC_2 = A_2$, according to Lemma 28.1(c), 
is 
$$
X = X_0 + F_{B_2} V_1 + V_2 E_{C_1} + F_{B_1}V_3E_{C_2},     
$$
where $X_0$ is a particular common solution to the pair $ B_1XC_1 = A_1$ and $ B_2XC_2 = A_2$,
$ V_1$---$V_3$ are arbitrary; the general common solution to $B_3YC_3 = A_3 $ and  $B_4 YC_4 = A_4$ is 
is
$$
Y = Y_0 - F_{B_4} W_1 - W_2 E_{C_3} - F_{B_3}W_3E_{C_4},  
$$
where  $Y_0$ is a particular common solution of  the pair $B_3YC_3 = A_3 $ and $B_4 YC_4 = A_4,$  $W_1$---$W_3$ are arbitrary. Hence
\begin{eqnarray*} 
X - Y & = & X_0 - Y_0 + F_{B_2} V_1  + F_{B_4} W_1 + V_2 E_{C_1} + W_2 E_{C_3} 
+ F_{B_1}V_3E_{C_2} + F_{B_3}W_3E_{C_4}  \\
& = & Z  + [\, F_{B_2}, \ F_{B_4} \,] \left[ \begin{array}{c} V_1  \\ W_1 \end{array} \right] + [\, V_2, \ W_2\,] 
\left[ \begin{array}{c} E_{C_1}  \\ E_{C_3}\end{array} \right] + F_{B_1}V_3E_{C_2} + F_{B_3}W_3E_{C_4},  
\end{eqnarray*}
where $ Z = X_0 - Y_0$. Applying (27.15) to it, one can determine (28.25), we leave the routine work to the reader. 
Furthermore, we have the following useful consequence.  

\medskip

\noindent {\bf Theorem 28.9.}\, {\em  Suppose that the quadruple matrix equations {\rm (28.24)} satisfy
the condition {\rm (28.26)}. Then they have  a common solution if and only if the following fourteen rank equalities 
are all satisfied
$$\displaylines{
\hspace*{2cm}
 r[ \, B_i , \ A_i \,] = r( B_i ) , \qquad   r \left[ \begin{array}{c} C_i  \\ A_i  \end{array} \right]
 = r( C_i ), \ \ \ \ i = 1, \, 2, \, 3, \, 4,   \hfill (28.28)
\cr
\hspace*{2cm}
r\left[ \begin{array}{ccc}  A_1  & 0 & B_1  \\ 0 &  -A_2  & B_2  \\ C_1 & C_2 & 0  \end{array} \right] =
  r \left[ \begin{array}{c} B_1  \\ B_2  \end{array} \right] + r[\, C_1, \  C_2 \,],  \hfill (28.29)
\cr
\hspace*{2cm} 
r\left[ \begin{array}{ccc}  A_3  & 0  & B_3  \\ 0 &  -A_4  & B_4  \\ C_3 & C_4 & 0  \end{array} \right]  = 
r \left[ \begin{array}{c} B_3  \\ B_4  \end{array} \right] + r[\, C_3, \  C_4 \,],  \hfill (28.30)
\cr
\hspace*{2cm} 
r\left[ \begin{array}{ccc}  A_i  & 0  & B_i  \\ 0 &  -A_j  & B_j  \\ C_i & C_j & 0  \end{array} \right]  = 
r \left[ \begin{array}{c} B_i  \\ B_j  \end{array} \right] + r[\, C_i, \  C_j \,], \ \ \ \  i = 1, \, 2, \ \ 
j=  3, \, 4.  \hfill (28.31) 
\cr}
$$}
 
In fact, it is obvious that (28.24) has a common solution if and only if  the two pairs $ B_1 XC_1 = A_1, \, B_2 XC_2 = A_2$ and $B_3 YC_3 = A_3, \, B_4YC_4 = A_4$ are consistent, respectively, and the  equation 
$$
 [\, F_{B_2}, \ F_{B_4} \,] \left[ \begin{array}{c} V_1  \\ W_1 \end{array} \right] + [\, V_2, \ W_2\,] 
\left[ \begin{array}{c} E_{C_1}  \\ E_{C_3}\end{array} \right] + F_{B_1}V_3E_{C_2} + F_{B_3}W_3E_{C_4}
= Y_0 - X_0 \eqno (28.32)
$$
is consistent. According to corollary 20.3, the consistency conditions for the two pairs $ B_1 XC_1 = A_1,$ $ 
B_2 XC_2 = A_2$ and $B_3 YC_3 = A_3, \, B_4YC_4 = A_4$ are the ten rank equalities  in (28.28)---(28.30). 
Next applying the four rank equalities in (27.17) and (27.18) to the equation (28.32) and simplifying, we can eventually 
find that (28.32) is consistent  if and only if the four rank equalities in (28.31) hold.   Based on Theorem 28.9, we are
 now able to establish a consistency condition for the matrix equation 
$$ 
B_1 X_1 C_1 + B_2 X_2C_2 +  B_3 X_3 C_3  = A, \eqno (28.33)
$$
which was  presented in Theorem 27.9.

\medskip

\noindent{\bf The Proof of Theorem 27.9.}\, Write first (28.33) as
$$
         B_1 X_1 C_1 + B_2 X_2C_2 = A  - B_3 X_3 C_3 .   \eqno (28.34) 
$$     
Then by Corollary 27.3 we know that this equation is solvable if and only if there exists an $ X_3$ satisfying 
the following four rank equalities
$$
\displaylines{
\hspace*{1cm}
 r[\, B_1, \ B_2, \ A - B_3 X_3 C_3 \,] = r[\,B_1, \ B_2 \,], \ \ \ \ r \left[ \begin{array}{c} C_1  \\ C_2 \\
 A - B_3 X_3 C_3 \end{array} \right] = r \left[ \begin{array}{c} C_1  \\ C_2 \end{array} \right],
 \hfill (28.35) 
\cr
\hspace*{1cm}
r\left[ \begin{array}{cc}  A - B_3 X_3 C_3  & B_1  \\ C_2 &  0  \end{array} \right] = r(B_1) +
 r(C_2),    \qquad  r\left[ \begin{array}{cc}  A - B_3 X_3 C_3  & B_2  \\ C_1 &  0 \end{array} \right] = 
r(B_2) + r(C_1).  \hfill (28.36) 
\cr}
$$  
Applying (1.2)---(1.5) to the right hand sides of  these four rank equalities, we see that they are equivalent to
the following four matrix equations 
$$ 
\displaylines{
\hspace*{2cm}
E_P B_3 X_3 C_3 = E_PA,   \qquad    E_{B_1}B_3 X_3 C_3F_{C_2}  = E_{B_1}AF_{C_2}, \hfill
\cr
\hspace*{2cm}
E_{B_2}B_3 X_3 C_3F_{C_1}  = E_{B_2}AF_{C_1} ,  \qquad    B_3 X_3 C_3F_Q  = AF_Q,   \hfill
\cr}
$$ 
where $ P = [\, B_1, \  B_2 \,]$ and  $Q = \left[ \begin{array}{c} C_1  \\ C_2 \end{array} \right],$  
which can be simply written as 
$$ 
\displaylines{
\hspace*{2cm}
G_1 X H_1 = L_1 ,   \qquad G_2 X H_2 = L_2 ,   \qquad G_3 X H_3 = L_3 ,  \qquad   G_34X H_4 = L_4,  \hfill (28.37) 
\cr}
$$    
where 
$$ 
\displaylines{
\hspace*{2cm}
G_1 =  E_P B_3, \ \ \ \  G_2 = E_{B_1}B_3, \ \ \ \  G_3 = E_{B_2}B_3, \ \ \ \  G_4 = B_3,  \hfill (28.38)
\cr
\hspace*{2cm} 
H_1 = C_3, \ \ \ \ H_2 = C_3F_{C_2}, \ \ \ \  H_3 = C_3F_{C_1}, \ \ \ \ H_4 = C_3F_Q,  \hfill (28.39) 
\cr
\hspace*{2cm}
L_1 = E_PA, \ \ \ \  L_2 = E_{B_1}CF_{B_2}, \ \ \ \  L_3 = E_{B_2}AF_{C_1}, \ \ \ \ L_4 = AF_Q.  \hfill (28.40)
\cr}
 $$ 
It is not difficult to deduce that the given matrices in (28.37) satisfy the following four range inclusions 
$$ \displaylines{
\hspace*{1cm}
R( G_1^T ) \subseteq  R( G^T_2 ), \qquad  R( H_2)  \subseteq   R( H_1 ), \qquad
 R( G_3^T) \subseteq  R( G^T_4 ), \qquad  R( H_4)  \subseteq   R( H_3 ). \hfill (28.41)
\cr}
$$ 
Thus by Theorem 28.9 we know that the four equations in (28.37) have a common solution if and only if the following
 fourteen rank equalities all hold
$$\displaylines{
\hspace*{2cm}
 r[\, G_i , \ L_i \,] = r( G_i ) , \qquad  r \left[ \begin{array}{c} H_i  \\ L_i  \end{array} \right]
 = r( H_i ), \ \ \ \ i = 1, \, 2, \, 3, \, 4,   \hfill (28.42) 
\cr
\hspace*{2cm}
r\left[ \begin{array}{ccc}  L_1  & 0  & G_1  \\ 0 &  -L_i  & G_i  \\ H_1 & H_i & 0  \end{array} \right] =
  r \left[ \begin{array}{c} G_1  \\ G_i \end{array} \right] + r[\, H_1, \  H_i\,],  \ \ \ \ i =  2, \, 3, \, 4,  
 \hfill (28.43) 
\cr
\hspace*{2cm}
r\left[ \begin{array}{ccc}  L_i  & 0  & G_i  \\ 0 &  -L_4  & G_4  \\ H_i & H_4 & 0  \end{array} \right] =
  r \left[ \begin{array}{c} G_i  \\ G_4 \end{array} \right] + r[\, H_i, \  H_4\,],  \ \ \ \ i =  2, \, 3, 
  \hfill (28.44) 
\cr
\hspace*{2cm}
r\left[ \begin{array}{ccc}  L_2  & 0 & G_2  \\ 0 &  -L_3  & G_3  \\ H_2 & H_3 & 0  \end{array} \right] =
  r \left[ \begin{array}{c} G_2  \\ G_3 \end{array} \right] + r[\, H_2, \  H_3 \,]. \hfill (28.45) 
\cr}
$$ 
Substituting the explicit expressions of $G_i, \  H_i $ and $ L_i (i = 1, \ 2, \ 3, \ 4)$  into
the eight rank equalities in (28.42) and  simplifying  by (1.2)---(1.4),  we can find that they are  
 equivalent to the first eight rank  equalities in Theorem 27.9, respectively. Next substituting (28.38)---(28.40) 
into the five rank equalities in  (28.43) and (28.44) and simplifying by (1.2)---(1.4), we can also find that 
they are equivalent to the first eight rank  equalities in Theorem 27.9, respectively. We omit the routine processes  
here for simplicity. As for (28.45), we have by (1.2)---(1.4) that 
\begin{eqnarray*}
r\left[ \begin{array}{ccc}  L_2  & 0  & G_2  \\ 0 &  -L_3  & G_3  \\ H_2 & H_3 & 0 \end{array} \right] 
& = & r\left[ \begin{array}{ccc}   E_{B_1}AF_{C_2} & 0  & E_{B_1}B_3  \\ 0 &  -E_{B_2}CF_{C_1}  & E_{B_2}B_3 
 \\ C_3F_{C_2} & C_3F_{C_1} & 0 \end{array} \right]  \\
& = & r \left[ \begin{array}{crccc} A  & 0 & B_1 &  0 & B_3  \\ 0  &  -A & 0 & B_2 & B_3  \\  C_2  & 0 & 0 & 0 & 0  \\ 
0  & C_1 & 0 & 0  & 0  \\ C_3  & C_3  & 0 &  0 & 0  \end{array} \right]  -  r \left[ \begin{array}{cc}  B_1  & 0 \\ 0  &  B_2 \end{array} \right] - 
r \left[ \begin{array}{cc}  C_2  & 0 \\ 0  &  C_1 \end{array} \right], 
\end{eqnarray*}
\begin{eqnarray*}
r \left[ \begin{array}{c} G_2  \\ G_3 \end{array} \right] + r[\, H_2, \  H_3 \,] 
& = & r \left[ \begin{array}{c} E_{B_1}B_3 \\ E_{B_2}B_3 \end{array} \right] 
+ r[\, C_3F_{C_2}, \  C_3F_{C_1} \, ] \\ 
& = & r\left[ \begin{array}{ccc}  B_1  & 0  & B_3  \\ 0 &  B_2  & B_3   \end{array} \right] 
+ r \left[ \begin{array}{cc}  C_2  & 0 \\ 0  &  C_1 \\  C_3 &  C_3   \end{array} \right] - 
r \left[ \begin{array}{cc}  B_1  & 0  \\ 0  &  B_2 \end{array} \right] - r \left[ \begin{array}{cc}  C_2  & 0 \\ 0  & 
 C_1 \end{array} \right]. 
\end{eqnarray*}
Thus (28.45) is equivalent to the last rank  equality in Theorem 27.9. Summing up, we find that  
(28.37) has a common solution, or equivalently  (28.33) is consistent, if and only if the nine rank 
 equalities in Theorem 27.9 all hold. \qquad $\Box$

\markboth{YONGGE  TIAN }
{29. EXTREME RANKS OF $ A - BX - XC$ SUBJECT TO  $BX C = D$} 

\chapter{Extreme ranks of $ A - BX - XC$ subject to $BX C = D$} 

\noindent As a simple application of the rank formulas in Chapter 27, we determine in this chapter  extreme ranks of a linear matrix equation $A - BX - XC$ subject to a consistent matrix $BXC = D$. This work is motivated by factoring a matrix $D$
as $ A = BB^- \pm B^-C$, and some related topics. Another motivation is from considering extreme ranks of
 $ A- BX- XC$ subject to $X.$ Quite different to the matrix expressions in the previous chapters, the same variant term $ X$ occurs two places in $ A- BX- XC$. Although it is quite simple in form, we fail to establish a general method for expressing its extreme ranks except some special cases. 
An interesting exception is that when $X$ is restricted by a consistent matrix equation $ BXC = D$, extreme ranks of $ A- BX- XC$ can completely be determined.    
 
\medskip

\noindent {\bf Theorem 29.1.}\, {\em Let  $ A, \, D  \in {\cal F}^{ m \times n}, \
B \in {\cal F}^{ m \times m}$ and $ C \in {\cal F}^{ n \times n}$ be
given$,$  and the matrix equation $ BXC = D $ is consistent. Then
$$
\displaylines{
\hspace*{0.5cm}
\max_{BXC = D}r(\, A - BX - XC \,) = \min \left\{  m +  r[ \, BA - D, \  B^2 \, ] - r( B ), \ \  
 n +  r \left[ \begin{array}{c} AC - D \\  C^2 \end{array} \right] - r(C),  \right.  \hfill
\cr
\hspace*{4.5cm}
\left.  r \left[ \begin{array}{cc} A & B  \\  C & 0 \end{array} \right],  
\ \   m + n +  r( \, BAC - BD - DC \,) - r(B) - r(C) \right\}, \hfill (29.1)
\cr
\hspace*{0cm}
and  \hfill
\cr
\hspace*{0.5cm}
\min_{BXC = D}r(\, A - BX - XC \,) =  r[ \, BA - D, \  B^2 \, ] +  r \left[ \begin{array}{c} AC - D \\ 
 C^2 \end{array} \right] +  \max\{  \, s_1, \ \ s_2 \, \},  \hfill (29.2)
\cr}
$$
where
$$
\displaylines{
\hspace*{1.5cm}
s_1  = r \left[ \begin{array}{cc} A & B  \\  C & 0 \end{array} \right]
 - r \left[ \begin{array}{cc}    C & 0 \\ BA & B^2 \end{array} \right] -
 r \left[ \begin{array}{cc} B & AC  \\  C^2 & 0 \end{array} \right],  \hfill
\cr
\hspace*{1.5cm}
s _2 =  r( \, BAC - BD - DC \,) - r[ \, BAC - DC, \  B^2 \, ] - r \left[ \begin{array}{c} BAC - BD \\  C^2 \end{array} 
\right]. \hfill
\cr}
$$
}
{\bf Proof.}\, Putting the general solution $X = B^-DC^- + F_B V_1 + V_2 E_C$ of $ BXC = D$ 
in $ A - BX - XC$ we first get 
$$
\displaylines{
\hspace*{0.5cm}
A - BX - XC  = A - DC^-  B^-D - F_BV_1C - BV_2 E_C  = p( V_1, \, V_2 ),  \hfill (29.4)
\cr}
$$
Clearly, this is a linear matrix expression involving two independent variant matrices $ V_1$ and $V_2$. In that case, 
we get by (27.6) and (27.7) that 
$$
\displaylines{
\hspace*{0.5cm} 
\max_{V_1, \, V_2} r[ \, p( V_1, \, V_2 )\,] = \min  \left\{ r[ \, A_1, \,  B, \,  F_B, \, ], \ \
r \left[ \begin{array}{c}  A_1  \\ C \\ E_C  \end{array} \right], \ \
r \left[ \begin{array}{cc}  A_1 & B \\ C  & 0 \end{array} \right], \ \ r \left[ \begin{array}{cc}  A_1   
& F_B \\ E_C   & 0 \end{array} \right]  \right\}, \hfill (29.5)
\cr
\hspace*{0cm}
and  \hfill
\cr 
\hspace*{0.5cm} 
\min_{V_1, \, V_2} r[ \, p( V_1, \, V_2 )\,] \hfill
\cr
\hspace*{0.5cm} 
= r \left[ \begin{array}{c}  A_1  \\ C \\ E_C  \end{array} \right] + 
 r[ \, A, \, B, \, F_B \, ] + \max \left\{ r\left[ \begin{array}{cc} A_1 & B \\ C& 0 \end{array} \right] 
- r\left[ \begin{array}{ccc} A_1 & B & F_B  \\ C   & 0 & 0  \end{array} \right] - r\left[ \begin{array}{cc} A_1   
& B \\ C   & 0 \\ E_C & 0  \end{array} \right] \right., \hfill
\cr
\hspace*{5cm} 
 \left. r \left[ \begin{array}{cc}  A_1   
& F_B \\ E_C & 0 \end{array} \right]  - r\left[ \begin{array}{ccc} A_1  & E_B & B  \\ E_C   & 0 & 0  \end{array} 
\right] - r\left[ \begin{array}{cc} A_1 & E_B  \\ E_C   & 0 \\ C & 0  \end{array} \right] \ \right\}, \hfill (29.6)
\cr}
$$
where $ A_1 =  A - DC^- - B^-D$. Simplifying the ranks of the block matrix in them by Lemma 1.1, we have
$$
\displaylines{
 \hspace*{1cm}
r[ \, A_1, \, B, \, F_B \, ] = r \left[ \begin{array}{ccc} A - B^-D & B  &
 I_m  \\ 0 &  0 & B  \end{array} \right] -r(B) =  r[ \, BA - D, \  B^2 \, ] + m - r( B ), \hfill
\cr
\hspace*{1cm}
 r \left[ \begin{array}{c}  A_1  \\  C   \\ E_C \end{array} \right] =
 r \left[ \begin{array}{cc} A -DC^-   & 0 \\ C & 0  \\  I_n  & C
 \end{array} \right] - r(C) =  r \left[ \begin{array}{c} AC - D \\
 C^2 \end{array} \right]  + n - r(C), \hfill
\cr
\hspace*{1cm}
r \left[ \begin{array}{cc} A_1 & B   \\ C & 0 \end{array} \right]
 = r \left[ \begin{array}{cc} A - DC^- -B^-D & B  \\  C &  0 \end{array}
 \right] = r \left[ \begin{array}{cc} A & B   \\ C & 0 \end{array} \right],
  \hfill
\cr }
$$
\begin{eqnarray*}
r \left[ \begin{array}{cc}  A_1 & F_B   \\ E_C &  0 \end{array} \right]
& = & r \left[ \begin{array}{ccc} A &  I_m  & 0  \\ I_n  &
 0  & C  \\  0 & B & 0 \end{array} \right]  - r(B) - r(C)  \\  
& = & r \left[ \begin{array}{ccc} 0 & I_m  & 0  \\ I_n  &  0  & 0  \\
0 & 0 &  BAC - BD - DC \end{array} \right] - r(B) - r(C) \\ 
 &= & m + n +  r( \, BAC - BD - DC \,) - r(B) - r(C),  
\end{eqnarray*}
\begin{eqnarray*}
r \left[ \begin{array}{ccc}  A_1 & B & F_B   \\ E_C &  0  & 0 \end{array} \right]
& = & r \left[ \begin{array}{ccc}  A  &  B  & I_m  \\ C  & 0  & 0  \\  0 & 0 & B \end{array} \right] - r(B)  \\  
& = & r \left[ \begin{array}{ccc} 0 & 0  & I_m  \\ C  &  0  & 0  \\
BA & B^2 & 0 \end{array} \right] - r(B) =  m +  r \left[ \begin{array}{cc}  C & 0 \\ BA  & B^2  \end{array} \right]  - r(B),  
\end{eqnarray*}
\begin{eqnarray*}
r\left[ \begin{array}{ccc}  A_1 & B \\ C &  0 \\ E_C & 0 \end{array} \right]
& = & r \left[ \begin{array}{ccc}  A  &  B  & 0 \\ C  & 0  & 0  \\  I_n & 0 & C \end{array} \right] - r(C)  \\  
& = & r \left[ \begin{array}{ccc} 0 & B  & AC  \\ 0  &  0  & C^2  \\
I_n & 0 & 0 \end{array} \right] - r(C) = n +  r \left[ \begin{array}{cc} B & AC \\ 0  & C^2  \end{array} \right]  - r(C),  
\end{eqnarray*}
\begin{eqnarray*}
r \left[ \begin{array}{ccc}  A_1 & F_B & B   \\ E_C &  0  & 0 \end{array} \right]
& = & r \left[ \begin{array}{cccc}  A - B^-D  & I_m  & B & 0  \\ I_n  & 0  & 0  & C \\  0 & B  & 0 & 0 \end{array} 
\right] - r(B) - r(C)  \\  
& = & r \left[ \begin{array}{cccc} 0  & I_m  & 0  & 0  \\ I_n  & 0  & 0  & 0 \\  0 & 0  & B^2 & 
BAC - DC \end{array} 
\right] - r(B) - r(C) \\
& = &  r[\, B^2, \ BAC- DC \,]  - r(B ) - r(C) + m + n,  
\end{eqnarray*}
\begin{eqnarray*}
r\left[ \begin{array}{cc}  A_1 & F_B \\ E_C &  0 \\ C & 0 \end{array} \right]
& = & r \left[ \begin{array}{ccc}  A DC^-  &  I_m  & 0 \\ I_n  & 0  & C  \\  C & 0 & 0 \\ 0 & B & 0
 \end{array} \right] - r(B) -r(C)  \\  
& = & r \left[ \begin{array}{ccc}  0  &  I_m  & 0 \\ I_n  & 0  & 0  \\  0 & 0 & C^2 \\ 0 & 0  & BAC - BD
\end{array} \right] - r(B) -r(C)  \\  
& = & r \left[ \begin{array}{c} C^2  \\ BAC - BD  \end{array} \right] - r(B) - r(C) + m + n.  
\end{eqnarray*}
Putting them in (29.3) and (29.4) yields (29.1) and (29.2).  \qquad $\Box$

\medskip

\noindent {\bf Corollary 29.2.}\, {\em Let  $ A, \, D  \in {\cal F}^{ m \times n}, \
B \in {\cal F}^{ m \times m}$ and $ C \in {\cal F}^{ n \times n}$ be
given$,$ 
 and the two matrix equations $ BX + YC = A$ and $ BXC = D $ are  consistent$,$ respectively. Then
$$
\displaylines{
\hspace*{0cm}
\max_{BXC = D}r(\, A - BX - XC \,) = \min \left\{  m +  r[ \, BA - D, \  B^2 \, ] - r( B ), \ \  
 n +  r \left[ \begin{array}{c} AC - D \\  C^2 \end{array} \right] - r(C),  \right.  \hfill
\cr
\hspace*{4.5cm}
\left.  r(B) + r(C),  \ \  m + n +  r( \, BAC - BD - DC \,) - r(B) - r(C) \right\}, \hfill (29.7)
\cr
\hspace*{0cm}
and \hfill
\cr
\hspace*{0cm}
\min_{BXC = D}r(\, A - BX - XC \,) =  r[ \, BA - D, \  B^2 \, ] +  r \left[ \begin{array}{c} AC - D \\ 
 C^2 \end{array} \right] \hfill
\cr
\hspace*{0cm}
+  \max \left\{  - r(B^2) - r(C^2),  \ \   r( \, BAC - BD - DC \,) - r[ \, BAC - DC, \  B^2 \, ] - 
r \left[ \begin{array}{c} BAC - BD \\  C^2 \end{array} \right] \right\}. \hfill (29.8)
\cr}
$$
}
{\bf Proof.}\, The consistency of $ BX + YC = A $ implies that 
$ r \left[ \begin{array}{cc} A  & B \\  C & 0  \end{array} \right] = r(B) +  r(C)$. Thus (29.1) and (29.2) reduce
 to (29.7) and (29.8)  \qquad $\Box$

\medskip

\noindent {\bf Corollary 29.3.}\, {\em Let  $ A, \, D  \in {\cal F}^{ m \times n}, \,
B \in {\cal F}^{ m \times m}$ and $ C \in {\cal F}^{ n \times n}$ be
given$,$  and $ BXC = D $ is consistent. If $r(B^2) = r(B)$ and $ r(C^2) = r(C),$ then
$$
\displaylines{
\hspace*{0cm}
\max_{BXC = D}r(\, A - BX - XC \,) = \min \left\{  m +  r[ \, D, \  B \, ] - r( B ), \ \  
 n +  r \left[ \begin{array}{c} D \\  C \end{array} \right] - r(C),  \right.  \hfill
\cr
\hspace*{4.5cm}
\left.  r(B) + r(C),  \ \  m + n +  r( \, BAC - BD - DC \,) - r(B) - r(C) \right\}, \hfill (29.9)
\cr
\hspace*{0cm}
and  \hfill
\cr
\hspace*{0cm}
\min_{BXC = D}r(\, A - BX - XC \,) =  r[ \, D, \  B \, ] +  r \left[ \begin{array}{c} D \\ 
 C \end{array} \right] + r( \, BAC - BD - DC \,) - r(B) - r(C). \hfill (29.10)
\cr}
$$
}
{\bf Proof.}\, Under $r(B^2) = r(B) $ and $ r(C^2) = r(C)$, there are 
$$
r[ \, BA - D, \  B^2 \, ] =  r[ \, D, \  B \,],  \ \ \ \ r \left[ \begin{array}{c} AC - D \\  C^2 \end{array} \right] 
 = r \left[ \begin{array}{c} D \\  C \end{array} \right],  \ \ \  \ r[ \, BAC - DC, \  B^2 \, ] = r(B),  
$$
$$
 r\left[ \begin{array}{c} BAC - BD \\  C^2 \end{array} \right] = r(C),  \ \ \ \  \left[ \begin{array}{cc} BA & B^2  \\  C & 0 \end{array} \right]  = \left[ \begin{array}{cc} B & AC 
 \\ C^2 & 0 \end{array} \right] = r(B) + r(C). 
$$ 
Thus we have (29.9) and (29.10). \qquad $\Box$

\medskip

\noindent {\bf Corollary 29.4.}\, {\em Let $ A, \, D  \in {\cal F}^{ m \times n}, \
B \in {\cal F}^{ m \times m}$ and $ C \in {\cal F}^{ n \times n}$ be
given$,$ and the two matrix equations $ BX + YC = A$ and $ BXC = D $ are  
consistent$,$ respectively. Then the pair of matrix equations   
$$
 BX + XC = A \qquad  BXC = D \eqno (29.11)
$$
have a common solution if and only if the following three conditions hold 
$$
R(\, BA - D \, ) \subseteq R(B^2), \qquad  R[(\, AC - D \, )^T] \subseteq R[(C^2)^T],  \qquad  BD + DC = BAC. 
\eqno (29.12)
$$ }
{\bf Proof.}\, Letting the right hand side of (29.8) be zero and  simplifying yield (29.12).
 \qquad $\Box$

\medskip

If $ BX + XC = A$ and $ BXC = D$ have a common solution, their general common solution can be simply  found by the 
 following two steps: First solve the matrix equation
$$
 F_BV_1C + BV_2 E_C = A - DC^-  - B^-D  \eqno(29.13)
$$ 
for $ V_1 $ and $ V_2$. Then put $V_1$ and $V_2$ in $X = B^-DC^- + F_B V_1 + V_2 E_C$ to yield the general common 
solution to the pair of equations. Based on the results in Lemma 26.2, we find that their general common 
solution can be written as 
$$
X = X_0 +  [\, F_B, \  0 \, ] F_G UE_H \left[ \begin{array}{c} I_n  \\ 0
 \end{array} \right] + [\,  0, \  I_m \, ] F_G UE_H
 \left[ \begin{array}{c} 0 \\ E_C \end{array} \right] + F_BSE_C,
\eqno (29.14) 
$$ 
where $ X_0$ is a particular common solution to (29.11),  $G = [\, F_B, \ -B \, ],
\ H = \left[ \begin{array}{c} C \\ E_C \end{array} \right],$  $U $ and $S$
are arbitrary. From Corollary 29.4 and (29.14) we also see that (29.11) has a unique common solution 
if and only if both $B$ and $ C $ are nonsingular and  $BD + DC = BDC$.
 In that case$,$ the unique common solution is $ X = B^{-1}DC^{-1}.$ 

\medskip

Theorem 29.1 can apply to determine extreme ranks of the matrix expressions $ A - BB^- - B^-B$  and 
$ A - BB^- + B^-B$  with respect to $B^-$. In fact 
$$\displaylines{
\hspace*{2cm}
\max_{B^-}r(\, A - BB^- - B^-B \,) = \max_{BXB = B}r(\, A - BX - XB \, ), \hfill 
\cr
\hspace*{2cm}
\min_{B^-}r(\, A - BB^- - B^-B \,) = \min_{BXB = B}r(\, A - BX - XB \, ), \hfill
\cr
\hspace*{2cm}
\max_{B^-}r(\, A - BB^- + B^-B \,) = \max_{BX(-B) = -B} r[\, A - BX - X(-B) \, ], \hfill
\cr
\hspace*{2cm}
\min_{B^-}r(\, A - BB^- + B^-B \,) = \min_{BX(-B) = -B }r[\, A - BX - X(-B) \, ]. \hfill
\cr}
$$
Applying Theorem 29.1, one can easily get the maximal and the minimal ranks of the matrix expressions, as well as 
necessary and sufficient conditions for the factorization $ A = BB^- +  B^-B$ or $ A = BB^- - B^-B$ to hold.

As we mentioned in beginning the chapter, it is a quite difficult problem to find in general extreme ranks 
of the linear matrix expressions $A - BX - XC$ as well as  $A - B_1XC_1 - B_2XC_2$. However, if the given matrices 
in them satisfy conditions, we can find their extreme ranks. Here we present two special results related to 
the minimal ranks of  $  A - BX + XC$ and  $  A - X + BXC$ when both  $ B$ and $ C$ are idempotent. 

\medskip

\noindent {\bf Theorem 29.5.}\,  {\em Suppose that $ B $ and $ C$ are $ m \times m$ and
 $ n \times n$ idempotent matrices$,$ respectively.  Then
$$\displaylines{
\hspace*{2cm}
\min_Xr(\, A - BX + XC \, ) = \max \{ \, r(BAC), \ \ \   r(\, I_m - B \,)A(\, I_n - C \,) \, \}, \hfill (29.15)
\cr
\hspace*{2cm}
\min_X r (\, A - X  + BXC \, ) = r(BAC). \hfill (29.16)
\cr}
$$
In particular$,$ 

{\rm (a)}\, The matrix equation $BX - XC = A$ is consistent if and only if $ABC = 0$ and 
$r(\, I_m - B \,)A(\, I_n - C \,) = 0.$ 

{\rm (b)}\, The matrix equation $X - BXC = A$ is consistent if and only if $ABC = 0$.}

\medskip

\noindent {\bf Proof.}\, Observe that 
$$\displaylines{
\hspace*{2cm}
B(\, A - X + BXC \,)C = BAC - BXC + B^2XC^2 = BAC. \hfill
\cr}
$$
We first see that $r(\, A - X + BXC \, ) \geq r(BAC)$ holds for all
$X$. On the other hand, let $ X = A$, then $ A  - A  + BAC = BAC$.
The combination of the above two facts yields (29.16).

To prove (29.15), we use the simple result
$$
\displaylines{
\hspace*{2cm}
\min_{X, \,Y} r \left[ \begin{array}{cc}  M &  X  \\ Y & N \end{array}
\right]
 = \max \{ \, r(M), \ \ \ r(N) \, \}.   \hfill (29.17)
\cr}
$$
Since both $ B$ and $ C$ are idempotent, we can factor them as
$$
\displaylines{
\hspace*{2cm}
B = P^{-1} \left[ \begin{array}{cc} I_k &  0  \\ 0 &  0 \end{array}
\right]P, \ \ \ C = Q \left[ \begin{array}{cc} I_l &  0  \\ 0 &  0
\end{array} \right] Q^{-1}, \hfill
\cr}
$$
where $ k = r(B)$ and $ l = r(C)$. In that case,
$$\displaylines{
\hspace*{2cm}
A - BX + XC = P^{-1} \left( \, PAQ - \left[ \begin{array}{cc} I_k &  0  \\ 0 & 0
\end{array} \right]PXQ  + PXQ \left[ \begin{array}{cc} I_l &  0  \\ 0 &  0
\end{array} \right] \,  \right)  Q^{-1}.   \hfill (29.18)
\cr}
$$
Let $ Y = PXQ  = \left[ \begin{array}{cc} Y_1 &  Y_2  \\ Y_3 & Y_4
\end{array} \right]$ and $PAQ
= \left[ \begin{array}{cc} S_1 &  S_2  \\ S_3
&  S_4 \end{array} \right].$ Then from (29.18) we get
$$\displaylines{
\hspace*{1cm}
r(\, A - BX + XC \,) =  r \left( \, \left[ \begin{array}{cc} S_1 &  S_2  \\
S_3 & S_4 \end{array} \right]  - \left[ \begin{array}{cc}
Y_1 & Y_2 \\ 0 & 0 \end{array} \right]  +  \left[ \begin{array}{cc} Y_1 &  0  \\ Y_3
&  0 \end{array} \right] \right) =
r \left[ \begin{array}{cc} S_1 &  S_2  - Y_2 \\ S_3 + Y_3 & S_4 \end{array} \right]. \hfill
\cr}
$$
Applying (29.17) to it we find
$$\displaylines{
\hspace*{1cm}
\min_X r(\, A - BX + XC \,) = \min_{Y_2, \, Y_3}r \left[ \begin{array}{cc}
S_1 &  S_2 - Y_2 \\ S_3 + Y_3 & S_4 \end{array} \right] = \min \{ \ r(S_1), \ \ 
r(S_4) \ \}, \hfill
\cr
\hspace*{0cm}
where \hfill
\cr
\hspace*{1cm}
r(S_1)  = r \left( \, \left[ \begin{array}{cc} I_k & 0  \\
0 & 0 \end{array} \right] PAQ \left[ \begin{array}{cc}
I_l & 0 \\ 0 & 0 \end{array} \right] \right) = r \left( \, P^{-1}\left[
\begin{array}{cc} I_k & 0  \\
0 & 0 \end{array} \right]PAQ \left[ \begin{array}{cc}
I_l & 0 \\ 0 & 0 \end{array} \right]Q^{-1} \right) = r(BAC). \hfill
\cr
\hspace*{0cm}
and \hfill
\cr
\hspace*{1cm}
r(S_4) = r \left( \, \left[ \begin{array}{cc} 0 & 0  \\
0 & I_{m-k} \end{array} \right] PAQ \left[ \begin{array}{cc}
 0 & 0 \\ 0 & I_{n-l} \end{array} \right] \right) \hfill
\cr
\hspace*{1.9cm}
  =  r \left( \, P^{-1}\left[ \begin{array}{cc} 0 & 0  \\
0 & I_{m-k} \end{array} \right] PAQ \left[ \begin{array}{cc}
0 & 0 \\ 0 & I_{n-l}\end{array}  \right] Q^{-1} \right) = r(\, I_m - B \,)A(\, I_n - C \,). \hfill
\cr}
$$
Thus we have (29.15). \qquad $\Box$

\medskip

Finally we present another interesting result related to the minimal rank of a complex matrix with respect to 
its imaginary part. It could also be regarded a  special case of a  matrix  expression  $A - B_1XC_1 - B_2XC_2$. 
We leave it as an exercise  to the reader.

\medskip

\noindent {\bf Theorem 29.6.}\,  {\em Suppose that $ A$ and $X$ are two real matrices of the same size. Then 
$$
\min_Xr(\, A + iX \, )  = \frac{1}{2} \min_X r\left[ \begin{array}{rr} A  & -X  \\ X & A \end{array} \right]
 =  \left\{ \begin{array}{cl} r(A)/2  &  if \ r(A) \ is   \ even  \\
 \left[\,r(A)+ 1 \,\right]/2  & if \ r(A) \ is \ odd \end{array} \right..
$$ 
}

The problem can also reasonably be considered  for a quaternion matrix.  Here we list a conjecture.   

\medskip

\noindent {\bf Conjecture 29.7.}\, {\em  Let $A + iX + jY + kZ $  be a quaternion matrix$,$
 where $ i^2 = j^2 = k^2 = -1$ and $ ijk =-1,$ the matrices $A, \, X_1, \, X_2$ and $X_3$ are real. Then 
$$
\min_{X, \, Y, \, Z} r( \, A + iX + jY + kZ  \, ) = \left\{
\begin{array}{ll} 
r(A)/4  &   if \ r(A) \ \equiv 0 \  ( {\rm mod} = 4 ) \\
\left[ \, r(A) + 3 \, \right]/4 & if \ r(A) \ \equiv 1 \  ( {\rm mod} = 4 ) \\
\left[ \, r(A) + 2 \, \right]/4 & if \ r(A) \ \equiv 2 \  ( {\rm mod} = 4 ) \\
\left[ \, r(A) + 1 \, \right]/4 & if \ r(A) \ \equiv 3 \  ( {\rm mod} = 4 ).
\end{array} \right.
$$}

\markboth{YONGGE  TIAN }
{30. EXTREME RANKS OF  SOME QUADRATIC MATRIX EXPRESSIONS}
\chapter{Extreme ranks of some quadratic matrix expressions}

\noindent Without much effort, the work in previous chapters can be easily 
extended quadratic matrix expressions involving two independent variant 
matrices. In this chapter we first present the maximal and the minimal 
 ranks of a matrix expression 
$$ 
q(X_1, \, X_2) = A - ( \, A_1 - B_1X_1C_1\, )D( \, A_2 - B_2X_2C_2\, ) 
\eqno (30.1)
$$
subject to $X_1$ and $ X_2$, and then present their various consequences.
The fundamental tool used for coping with (30.1) is the following rank 
formula
$$ 
 r(\, A - PNQ \, ) =  r \left[ \begin{array}{cc} A  & PN \\  NQ  & N 
 \end{array} \right] - r(N). \eqno (30.2)
$$
Applying (30.2) to (30.1), we can get
$$
\displaylines{
\hspace*{1cm} 
 r[q(X_1, \, X_2)] \hfill
\cr 
\hspace*{1cm}
= r\left( \left[ \begin{array}{cc} A  & ( \, A_1 - B_1X_1C_1 \, )D \\  D( \, A_2 - B_2X_2C_2\, )  & D \end{array} \right] \right) -r(D) \hfill
\cr
\hspace*{1cm}
 = r\left( \left[ \begin{array}{cc} A  & A_1D \\  DA_2  & D \end{array} 
\right] -   \left[ \begin{array}{c} B_1 \\  0 \end{array} \right] X_1[\, 0, \ C_1D \, ] - 
  \left[ \begin{array}{c} 0 \\  DB_2 \end{array} \right] X_1[\, C_2, \ 0  \, ] \right) -r(D). \hfill (30.3)
\cr}
$$
Evidently the matrix expression in the right hand side of (30.3) is linear with two independent variant matrices. 
Applying the rank formulas (27.6) and (27.7) to (30.3) and simplifying,  we get the following. 

\medskip

\noindent  {\bf Theorem 30.1.}\, {\em Let $ q(X_1, \, X_2)$ be given by
{\rm (30.1)}.  Then
$$
\displaylines{
\hspace*{0cm} 
\max_{X_1, \, X_2} r[ \, q( X_1, \, X_2 )\,] = \min  \left\{  r[ \, A - A_1DA_2 , \ A_1DB_2, \  B_1 \, ], \ \
r \left[ \begin{array}{c}  A - A_1DA_2 \\ C_1DA_2 \\ C_2  \end{array} \right], \right.  \hfill  
\cr
\hspace*{5cm} 
\left.  r \left[ \begin{array}{cc}  A_1DA_2 - A  & B_1 \\ C_2   & 0 \end{array} \right], \ \ 
r \left[ \begin{array}{cc}  A_1DA_2 - A  & A_1DB_2  \\ C_1DA_2   & C_1DB_2  \end{array} \right]  \
 \right\}, \hfill (30.4)
\cr 
and \hfill
\cr
\hspace*{0cm} 
\min_{X_1, \, X_2} r[ \, q( X_1, \, X_2 )\,] = r \left[ \begin{array}{c}  A - A_1DA_2 \\ C_1DA_2 \\ C_2  \end{array} \right] + 
r[ \, A - A_1DA_2 , \ A_1DB_2, \  B_1 \, ] \hfill 
\cr
\hspace*{0cm}
+  \max \left\{ r \left[ \begin{array}{cc}  A_1DA_2 - A  & B_1 \\ C_2   & 0 \end{array} \right]
 - r \left[ \begin{array}{ccc}  A_1DA_2 - A  & B_1 & A_1DB_2 \\ C_2   & 0 & 0  \end{array} \right] -
r \left[ \begin{array}{cc}  A_1DA_2 - A  & B_1 \\ C_2   & 0 \\ C_1DA_2  & 0 \end{array} \right],   \right. \hfill
\cr}
$$
$$
\displaylines{
\hspace*{0cm} 
 \left. r \left[ \begin{array}{cc}  A_1DA_2 - A  & A_1DB_2  \\ C_1DA_2   & C_1DB_2  \end{array} \right] -
r \left[ \begin{array}{ccc}  A_1DA_2 - A  & A_1DB_2  & B_1  \\ C_1DA_2   & C_1DB_2 &  0  \end{array} \right]  
- r \left[ \begin{array}{cc}  A_1DA_2 - A  & A_1DB_2  \\ C_1DA_2   & C_1DB_2 \\ C_2 & 0 \end{array} \right]  \right\}. \ 
 \hfill (30.5) 
\cr}
$$
}
\hspace*{0.3cm} The  formulas (30.4) and (30.5) can further simplify when the given matrices 
in them  satisfy some conditions, for example, the two equations 
$ B_1X_1C_1 = A_1$ and $ B_2X_2C_2 = A_2$ are solvable, respectively; or 
some of them are identity matrices or zero matrices.  

Two nice results are given below. 

\medskip

\noindent  {\bf Corollary  30.2.}\, {\em Let $ q(X, \, Y) = XAY + XB + CY + D,$ where where $ A, \, B, \, C$ and $D$ are $m \times n$, \, $l \times n,$ 
$m \times k,$ and $l \times k$ matrices$,$ respectively. Then

{\rm (a)}\, The maximal and the minimal ranks of $ q(X_1, \, X_2)$ are    
$$
\displaylines{
\hspace*{2cm} 
\max_{X, \, Y} r[ \, q( X, \, Y )\,] = \min  \left\{ \, m, \ \ \ n, \ \ \
r \left[ \begin{array}{cc}  A
& B \\ C   & D \end{array} \right] \, \right\}, \hfill (30.6)
\cr
\hspace*{2cm} 
\min_{X, \, Y} r[ \, q( X, \, Y )\,] = \max  \left\{ \, 0, \ \ \
r \left[ \begin{array}{cc}  A & B \\ C   & D \end{array} \right] - 
r\left[ \begin{array}{c}  A \\ C \end{array} \right] - r[\, A, \ B \,]
 \, \right\}. \hfill (30.7)
\cr}
$$

{\rm (b)}\,  Let $ m = n.$ Then there are $ X $ and $ Y $ such that 
$q(X, \, Y )$ is nonsingular if and only if
$$
r \left[ \begin{array}{cc}  A & B \\ C & D \end{array} \right] \geq m. 
\eqno (30.8)
$$

{\rm (c)}\, There are $ X $ and $ Y $ such that  $XAY + XB + YC + D = 0$  
 if and only if
$$
r \left[ \begin{array}{cc}  A & B \\ C & D \end{array} \right] \leq  
r\left[ \begin{array}{c}  A \\ C \end{array} \right] + r[\, A, \ B \,].
 \eqno (30.9)
$$ }
\hspace*{0.3cm} The matrix expression $ q(X, \, Y) = XAY + XB + CY + D$ occurs in an elementary operation for 
a $ 2 \times 2$ block matrix  
$$
\left[ \begin{array}{cc}  I_m  & 0 \\ X  & I_l \end{array} \right]
\left[ \begin{array}{cc}  A  & B \\ C  & D \end{array} \right]
\left[ \begin{array}{cc}  I_n  & Y \\ 0  & I_k \end{array} \right]
=\left[ \begin{array}{cc}  A & AY + B  \\ C + XA  & XAY + XB + CY + D
 \end{array} \right]. \eqno (30.10)
$$
Clearly the lower right block in (30.10) is the matrix expression $ q(X, \, Y ).$ If we let $ X = - CA^- $ and $ Y = - A^-B$, where $ A^- $ is an inner inverse of $ A$, then $q(X, \, Y ) = D - CA^-B$, 
the well-known Schur complement $ A $ in $ M = \left[ \begin{array}{cc}  A  & B \\ C  & D \end{array} 
\right]$. Corollary 30.2 acctually gives possible ranks of the lower right block of (30.10) 
after the block elementary operation, including necessary and sufficient conditions for the block to 
be nonsingular or null.

\medskip

\noindent  {\bf Corollary  30.3.}\, {\em Let $ q(X_1, \, X_2) = A - BX_1DX_2C,$  where 
$ A, \ B, \ C$ and $D$ are $m \times n,$ \ $m \times k,$ 
$l \times n,$ and $p \times q$ matrices$,$ respectively. Then the maximal and the minimal ranks of
 $ q(X_1, \, X_2)$ are    
$$
\displaylines{
\hspace*{1cm} 
\max_{X_1, \, X_2} r[ \, q( X_1, \, X_2 )\,] = \min  \left\{  r\left[ \begin{array}{c}  A \\ C 
\end{array} \right],  \ \ r[\, A, \ B \,],  \ \ r(A) + r(D)  \right\}, \hfill (30.10)
\cr
\hspace*{1cm} 
\min_{X_1, \, X_2} r[ \, q( X_1, \, X_2 )\,] = \max \left\{  r(A) - r(D), \ \ \  
r\left[ \begin{array}{c}  A \\ C \end{array} \right] + r[\, A, \ B \,]
 - r \left[ \begin{array}{cc}  A & B \\ C   & 0 \end{array} \right]  \right\}. \hfill (30.11)
\cr}
$$}
\hspace*{0.3cm} Theoretically one can express rank of any nonlinear matrix expression  through  rank of a linear matrix expression. A basic 
transformation formula for this is 
$$
\displaylines{
\hspace*{1cm} 
r(\, A - B_1X_1B_2X_2\cdots B_kX_kB_{k+1} \,) \hfill
\cr
\hspace*{1cm} 
= r \left[ \begin{array}{ccccc}  
B_1X_1B_2  &  0  & \cdots &  0 & (-1)^k A  \\
B_2        & B_2X_2B_3 & \cdots & 0   & 0  \\  
           &   \ddots  & \ddots & \vdots  & \vdots \\
           &           &         & B_{k-1}X_{k-1}B_k & 0 \\
           &           &         & B_k  & B_kX_kB_{k+1}
 \end{array}\right] - r(B_2) - \cdots - r(B_k). \hfill (30.12) 
\cr}
$$
\hspace*{0.3cm} The block matrix on the right side of (30.12) is obviously a linear 
matrix expression.

As an important application we next consider extreme ranks of the Schur 
complement $ D - CA^-_rB$ with respect to an reflexive inner inverse 
 $ A^-_r$ of $ A$. A reflexive inner inverse of $ A $ is a solution
of the  pair of matrix equations $ AXA = A $ and $ XAX = X$. The
 general expression of reflexive inner inverse of $ A $ can be 
written as $ A^-_r = A^-AA^-  = (\, A^{\sim} - F_AV_1 \, )A(\, A^{\sim} - V_2E_A \,),$ 
where $ A^{\sim}$ is a particular inner inverse of $ A $, $ V_1 $ and $V_2$ are arbitrary.    
Therefore,  we have  
$$ 
D - CA^-_rB =D - (\, CA^{\sim} - CF_AV_1 \, )A(\, A^{\sim}B - V_2E_AB \,),
$$ 
Applying (30.4) and (30.5) to it we get the following. 
 
\medskip

\noindent  {\bf Theorem 30.4.}\, {\em The maximal and the minimal ranks of the 
Schur complement $ D - CA^-_rB $  with respect to $ A^-_r $ are given by  
$$ 
\displaylines{
\hspace*{0cm} 
\max_{A_r^-} r( \, D - CA^-_rB \,) = \min  \left\{ r(A) + r(D),  \ \   r[ \, A, \ B \, ], \ \
r \left[ \begin{array}{c}  A \\ C  \end{array} \right],  \ \ 
 r \left[ \begin{array}{cc}  A  & B \\ C   & D \end{array} \right] -r(A)  \right\}, \hfill (30.13)
\cr
\hspace*{0cm}
and \hfill
\cr
\hspace*{0cm}
\min_{A_r^-} r( \, D - CA^-_rB \,)  = r \left[ \begin{array}{c} B \\ D  \end{array} \right] + 
r[ \, C, \  D \, ] + r(A) \hfill 
\cr
\hspace*{0cm}
+  \max \left\{  r \left[ \begin{array}{cc} A  & B \\ C & D \end{array} \right]
 - r \left[ \begin{array}{ccc} A  & 0  & B \\ 0   & C & D  \end{array} \right] -
r \left[ \begin{array}{cc} A  & 0 \\ 0   & B \\ C & D \end{array} \right],  \  r(D) -  r \left[ \begin{array}{cc} A  & 0 \\ C  & D  \end{array} \right] -
r \left[ \begin{array}{cc}  A  & B \\ 0   & D  \end{array} \right]  \right\}.  \hfill (30.14) 
\cr}
$$
} 
\noindent {\bf Remark.} Just as what we did in Chapters 21, 22  and 23,  one can further find many consequences from (30.13) and (30.14), such as, the rank invariance of  $ D - CA^-_rB$ with respect to the choice of $ A^-_r$; various special cases of (30.13) and (30.14) and their interesting consequences when $ A, \, B, \, C$ and $ D$ satisfy some conditions;
reverse order laws for reflexive inner inverses of  products of matrices, rank equalities for sums of reflexive 
inner inverses of matrices,  and so on. The reader can easily list 
them and apply them to find some more interesting results.         
 
\medskip 

Through (30.4) and (30.5) we cal also derive various rank equalities for matrix expressions involving products of inner inverses of matrices. We next list several of them without detailed proofs. 

\medskip

\noindent {\bf Theorem 30.5.}\, {\em  Let $ A \in {\cal F}^{ m \times n}, \, 
 B \in {\cal F}^{ m \times k}$ and $ C \in {\cal F}^{ l \times n}$ be given. Then
$$ 
\displaylines{
\hspace*{1cm} 
\max_{B^-,\, C^-} r[ \, A - (\, I_m - BB^-\,)A(\, I_n - C^-C\,) \,] \hfill
\cr 
\hspace*{1cm} 
= \min  \left\{ r(B) + r(C),  \ \  r[ \, A, \ B \, ], \ \ r \left[ \begin{array}{c}  A \\ C  \end{array} \right],  \ \ 
 r \left[ \begin{array}{cc}  A  & B \\ C   & 0 \end{array} \right] + r(A) - r(B) - r(C) \right\}, \hfill (30.15)
\cr
\hspace*{0cm}
and \hfill
\cr
\hspace*{1cm}
\min_{B^-, \, C^-}r[ \, A - (\, I_m - BB^-\,)A(\, I_n - C^-C\,) \,] = 
r(A) + r(B) + r(C) - \left[ \begin{array}{cc}  A  & B \\ C   & 0  \end{array} \right]. \hfill (30.16) 
\cr}
$$}
{\bf Proof.}\, Notice the two general expressions  $BB^-  = BB^{\sim} + F_BV_1B$ and 
 $C^-C = CC^{\sim} + CV_2E_C$, where $ B^{\sim}$  and  $ C^{\sim}$ are two particular inner inverses of 
$ B $ and $C$, respectively, $ V_1 $ and $V_2$ are arbitrary. Put them in 
$ A - (\, I_m - BB^-\,)A(\, I_n - C^-C\,)$  to yield a quadratic matrix expression. In this case 
applying (30.4) and (30.5) to it and then simplifying  we may trivially give (30.15) and (30.16). \qquad $ \Box$   
 
\medskip

\noindent {\bf Theorem 30.6.}\, {\em  Let $ A \in {\cal F}^{ m \times n}, \, 
 B \in {\cal F}^{ p \times m}$ and $ C \in {\cal F}^{ n \times q}$ be given. Then
$$ 
\displaylines{
\hspace*{1cm} 
\max_{B^-,\, C^-} r[ \, (\, I_m - B^-B\,)A(\, I_n - CC^- \,) \,] = \min  \left\{ \, r(A),  \ \ m -r(B), 
\ \ n - r(C) \, \right\}, \hfill (30.17)
\cr
\hspace*{1cm}
\min_{B^-, \, C^-}r[ \,(\, I_m - B^-B\,)A(\, I_n - CC^- \,) \,] = \max \{ \, 0, \ \  \ r(A) - r(BA) - r(AC) \, \}. 
\hfill (30.18)
\cr}
$$
In particular$,$ there are $ B^-$ and $C^-$ such that $(\, I_m - B^-B\,)A(\, I_n - CC^- \,) = 0,$ if and only if 
 $ r(A) \leq r(BA) + r(AC)$. 
}

\medskip

Some special cases of (30.17) and (30.18) are listed below:
$$ 
\displaylines{
\hspace*{1cm} 
\max_{B^-,\, C^-} r[ \, (\, I_m - B^-B\,)(\, I_m - CC^- \,) \,] = \min  \left\{\,  m -r(B), 
\ \  \ m - r(C) \,  \right\}, \hfill
\cr
\hspace*{1cm}
\min_{B^-, \, C^-}r[ \,(\, I_m - B^-B\,)(\, I_m - CC^- \,) \,] = \max \{ \, 0, \ \  \ m - r(B) - r(C) \, \}, \hfill
\cr
\hspace*{1cm} 
\max_{A^-} r[ \, (\, I_m - A^-A \,)(\, I_m - AA^- \,) \,] = m  - r(A), \hfill
\cr
\hspace*{1cm}
\min_{A^-}r[ \,(\, I_m - A^-A \,)(\, I_m - AA^- \,) \,] = \max \{ \, 0, \ \ \ m - 2r(A) \, \}. \hfill
\cr}
$$
{\bf Theorem 30.7.}\, {\em  Let $ A \in {\cal F}^{ m \times m}$ be given. Then
$$  
\max_{A^-} r[ \, A -  (\, I_m - AA^-\,)(\, I_m - A^-A \,) \,] = 
\min_{A^-} r[ \, A -  (\, I_m - AA^-\,)(\, I_m - A^-A \,) \,] =  m + r(A^2) - r(A). \eqno (30.19)
$$}
{\bf Theorem 30.8.}\, {\em  Let $ A \in {\cal F}^{ m \times m}$ be given. Then
$$ 
\displaylines{
\hspace*{0cm} 
\max_{A^-} r[ \, A -  (\, I_m - A^-A \,)(\, I_m - AA^- \,) \,] =   m + r(A^2) - r(A), \hfill (30.20)
\cr
\hspace*{0cm} 
\min_{A^-} r[ \, A -  (\, I_m - A^-A\,)(\, I_m - AA^- \,) \,] = \max \{ \,  2r(A^2) - r(A^3),  \ \ 
m - 2r(A) + 2r(A^2) - r(A^3) \, \}.  \  \hfill (30.21)
\cr}
$$}
{\bf Theorem 30.9.}\, {\em  Let $ A, \ D \in {\cal F}^{ m \times n},  \ 
 B \in {\cal F}^{m \times k} $ and $C \in {\cal F}^{ l \times n}$ be given. Then
$$ 
\displaylines{
\hspace*{0cm} 
\max_{B^-,\, C^-} r( \, A - BB^-DC^-C \,) = \min  \left\{   r[ \, A, \ B \, ],  \ \  
r \left[ \begin{array}{c}  A \\ C  \end{array} 
\right],  \ \ r \left[ \begin{array}{ccc} A  & 0 & B \\ 0  & -D & B \\ C & C & 0
\end{array} \right] -r(B) - r(C) \right\},  \ \hfill (30.22) 
\cr
\hspace*{0cm}
and \hfill
\cr
\hspace*{0cm}
\min_{B^-, \, C^-} r( \, A - BB^-DC^-C \,)  = \min \{ \,  s_1, \ \  s_2 \, \}, \hfill (30.23) 
\cr
\hspace*{0cm}
where \hfill
\cr
\hspace*{0cm} 
s_1 = r[ \, A, \ B \, ]  +  r \left[ \begin{array}{c}  A \\ C  \end{array} \right] - 
r \left[ \begin{array}{cc} A  &  B \\ C & 0 \end{array} \right], \hfill
\cr
\hspace*{0cm}
s_2 = r(B) + r(C) +  r[ \, A, \ B \, ] + r \left[ \begin{array}{c} A \\ C \end{array} \right] 
+ r \left[ \begin{array}{ccc} A  & 0 & B \\ 0  & -D & B \\ C & C & 0 \end{array} \right]
- r \left[ \begin{array}{cccc} A  & 0 & B & 0  \\ 0  & D & 0 & B \\ C & C & 0 & 0 \end{array}\right] -
r \left[ \begin{array}{ccc} A  & 0 & B \\ 0  & D & B \\ C & 0 & 0 \\ 0 & C & 0 \end{array} \right]. \hfill
\cr}
$$}
\hspace*{0.3cm}Some useful consequences of (30.22) and (30.23) are listed below:

\medskip

\noindent {\bf Theorem 30.10.}\, {\em  Let $ A \in {\cal F}^{ m \times n},  \, 
 B \in {\cal F}^{m \times k} $ and $C \in {\cal F}^{ l \times n}$ be given. Then
$$ 
\displaylines{
\hspace*{0cm} 
\max_{B^-,\, C^-} r( \, A - BB^-AC^-C \,) = \min  \left\{ r \left[ \begin{array}{c}  A \\ C  \end{array} 
\right], \ \  r[ \, A, \ B \, ], \ \  \ r \left[ \begin{array}{c}  A \\ C  \end{array} \right] + 
 r[ \, A, \ B \, ] -r(B) - r(C) \right\},  \hfill (30.24) 
\cr
\hspace*{0cm}
and \hfill
\cr
\hspace*{0cm}
\min_{B^-, \, C^-} r( \, A - BB^-AC^-C \,)  = r[ \, A, \ B \, ]  + r \left[ \begin{array}{c}  A \\ C  \end{array} \right] -  r \left[ \begin{array}{cc} A  & B \\ C & 0 \end{array} \right]. \hfill (30.25)
\cr}
$$ }
\hspace*{0.3cm} Contrasting (30.25) with (18.6), we see that 
$$\displaylines{
\hspace*{2cm} 
\min_{B^-, \, C^-} r( \, A - BB^-AC^-C \,) =  \min_X r( \, A - BXC \,), \hfill 
\cr}
$$ 
which can be stated that there is a matrix $ X$ with form $ X = B^-AC^-$ such that $ A - BXC$ reaches to its minimal rank. We can call this $ X $  a minimal rank solution of $ BXC = A$. In that case, $ X = B^-AC^- + F_AV_1 +V_2E_B$ is also 
minimizing $ r(\ A - BXC \, )$, where $ V_1$ and $V_2$ are arbitrary.     

\medskip

\noindent {\bf Theorem 30.11.}\, {\em  Let $ B, \, C  \in {\cal F}^{ m \times m}$ be given. Then
$$ 
\displaylines{
\hspace*{2cm} 
\max_{B^-, \, C^-} r( \, BC - BB^-C^-C \,)  = \min \{ \, r(B),  \ \ r(C),  \ \  m + r( \, C  - CBC \,) - r(C) \, \}, \hfill (30.26)
\cr
\hspace*{2cm} 
\min_{B^-, \, C^-} r( \, BC - BB^-C^-C \,)  = \max \{ \, 0, \ \ r( \, C  - CBC \,) + r(B) - m \, \}, \hfill
\cr 
\hspace*{0cm} 
and \hfill (30.27)
\cr
\hspace*{2cm} 
\max_{B^-} r( \, B^2 - BB^-B^-B \,)  = \min \{ \, r(B),  \ \  m + r( \, B  - B^3 \,) - r(B) \, \}, \hfill 
(30.28)
\cr
\hspace*{2cm} 
\min_{B^- } r( \, B^2 - BB^-B^-B \,)  = \max \{ \, 0, \ \ r( \, B  - B^3 \,) + r(B) - m \, \}. \hfill 
(30.29)
\cr}
$$    
In particular$,$ there are  $ B^-$ and $C^-$  such that $BC$ can factor as 
$BC = (BB^-)(C^-C)$ if and only if $ r( \, C - CBC \,) \leq  m - r(B).$ There is $ B^-$ such that $ B^2$ can factor as 
$B^2 = (BB^-)(B^-B)$ if and only if $ r( \, B - B^3 \,) \leq m - r(B).$}

\medskip

\noindent {\bf Theorem 30.12.}\, {\em  Let $ A \in {\cal F}^{ m \times m}$ be given. Then
$$\displaylines{
\hspace*{2cm} 
 \max_{A^-} r( \, A - AA^-A^-A \,)  = \max \{ \, r(A), \ \ r( \, I_m - A \,) \, \}, \hfill (30.30)
\cr
\hspace*{2cm} 
 \min_{A^-} r( \, A - AA^-A^-A \,)  = \max \{ \, 0, \ \  r( \, I_m - A \,) + 2r(A) - 2m \, \}. \hfill (30.31)
\cr}
$$    
In particular$,$ there is an $ A^-$ such that $ A $ can factor as $A = (AA^-)(A^-A)$ if and only if 
$ r( \, I_m - A \,) \leq  2m - 2r(A).$ Moreover$,$  when $ r(A) \leq m/2,$  there must exist an  $ A^-$ such that $ A $ can factor as $A = (AA^-)(A^-A)$.}
  
\medskip

\noindent {\bf Theorem 30.13.}\, {\em  Let $ B, \, C  \in {\cal F}^{ m \times m}$ be given. Then
$$
\displaylines{
\hspace*{2cm} 
 \max_{B^-, \, C^-} r( \, BB^-C^-C \,)  = \max \{ \, r(B), \ \ r(C) \, \}, \hfill (30.32)   
\cr
\hspace*{2cm} 
\min_{B^-, \, C^-} r( \,B^-C^- \,) = \min_{B^-, \, C^-} r( \,BB^-C^-C \,)  = \max \{ \, 0, \ \  r(B) + r(C) - m  \, \}.
\hfill (30.33)  
\cr}
$$    
In particular$,$ there are $B^-$ and $C^-$  such that $ B^-C^- = 0$ if and only if $r(B) + r(C) \leq m.$ }

\medskip

Moreover one can also find extreme ranks of matrix expressions $ A - (\, I_m - BB^-\,)D(\, I_n - C^-C\,),$  
$ A - (\, I_m - B^-B\,)D(\, I_n - CC^-\,),$ $ A - B^-BDCC^-$, $ A - B^kB^-DC^-C^k$, and so on. Based on them 
more consequences can be derived. We leave them to the reader.  In addition,  we present  another 
interesting result for the reader to prove 
$$
\displaylines{
\hspace*{2cm} 
\min_{A^-, \,B^-, \, C^-} r( \,AA^-BB^-CC^- \,) = {\rm dim} [\, R(A) \cap  R(B) \cap R(C)\, ]. \hfill 
(30.34)
\cr}
$$  

Finally we present a conjecture on the minimal rank of multiple product of inner inverses. 

\medskip

\noindent {\bf Conjecture 30.14.}\, {\em  Let $ A_i \in {\cal F}^{ m_{i+1} \times m_i}, \, i = 1, \, 2,
 \, \cdots, \, k.$  Then  
$$
\min_{A^-_1, \, \cdots \, A^-_k} r( \,A^-_1A_2^- \cdots A^-_k \,) 
= \max \{ \, 0, \ \  r(A_1) + r(A_2)  + \cdots  +r(A_k) - m_2 - m_3 - \cdots - m_k \, \}.   \eqno (30.35)    
$$
}
 
\markboth{YONGGE  TIAN }
{31. COMPLETING BLOCK MATRICES WITH EXTREME RANKS }
\chapter{Completing triangular block matrices with extreme ranks}

\noindent Suppose that $A_n$ and $ X_n$ are  triangular block matrices with the forms 
$$
A_n = \left[ \begin{array}{ccccc}  A_{11} & & & & \cr
A_{21}&A_{22}& & &\cr
\vdots &\vdots &\ddots & &\cr
A_{n1}& A_{n2}&\cdots & A_{nn} \end{array} \right] , \qquad  X_n = \left[ \begin{array}{ccccc} 0 & X_{12} &
\cdots & X_{1n} \cr
& \ddots  & \ddots & \vdots \cr
& & 0 & X_{n-1,n}\cr
& & & 0 \end{array} \right] , \eqno (31.1)
$$ 
where $A_{ij}\ (n \geq i \geq j \geq 1)$ is a given $s_i \times t_j$ matrix,
  $X_{ij} \ (1 \leq i < j \leq n ) $ is a variant $s_i\times t_j $ matrix.
  Further let $S(X_n) $ be the collection of all matrices $ X_n$ in (31.1).
  In this article we consider how to choose $X_n \in  S(X_n)$ such
that 
$$
r(\, A_n+X_n \, ) = \left[ \begin{array}{ccccc} A_{11} & X_{12} & \cdots & X_{1n} \cr
A_{21} & A_{22} & \ddots & \vdots \cr
\vdots & \vdots & \ddots & X_{n-1,n}\cr
A_{n1} & A_{n2} & \cdots &  A_{nn} \end{array} \right] = \max \eqno (31.2) 
$$
and
$$
r(\, A_n+ X_n \, ) = \left[ \begin{array}{cccc} A_{11} & X_{12} & \cdots & X_{1n} \cr
A_{21} & A_{22} & \ddots & \vdots \cr
\vdots & \vdots & \ddots & X_{n-1,n}\cr
A_{n1} & A_{n2} & \cdots &  A_{nn}  \end{array} \right] = \min  \eqno (31.3)
$$
hold, respectively.

These two problems are well known in matrix theory as maximal and minimal
rank completion problems, which have been previously examined by lots of
authors from different aspects (see, e.g., \cite{CJRW, GKL, Jo, JW, Wo1, Wo2}). In this chapter, we wish to give a new
investigation to the two problems by making use of the theory of
generalized inverses of matrices.

\medskip

\noindent { \bf Lemma 31.1.}\, {\it Suppose that 
$$ \displaylines{
\hspace*{2cm}
M(X_{12}) = \left[ \begin{array}{cc} A_{11} & X_{12} \cr
A_{21} & A_{22}  \end{array} \right] \hfill (31.4)
\cr}
$$ 
is a $ 2 \times 2$ block matrix$,$ where $ A_{11}, \,   A_{21}$  and
$  A_{22}$ are three given $s_1 \times t_1, s_2\times t_1 \  and \
s_2 \times t_2 $ matrices$,$ respectively$,$ and $X_{12} $ is a variant
$s_1 \times t_2 $ matrix. Then

{\rm(a)}\,  The maximal rank of $M(X_{12}) $ with respect to $X_{12}$ is
$$\displaylines{
\hspace*{2cm}
\max_{X_{12}} r[M(X_{12} )]  = \min \left\{  r\left[ \begin{array}{c}  A_{11} \cr A_{21}
 \end{array} \right]+ t_2, \quad   r[\, A_{21}, \ A_{22} \, ]+ s_1 \right\}, \hfill (31.5)
\cr}
$$
and the matrix $ X_{12} $ satisfying {\rm (31.5)} can be expressed as 
 $$\displaylines{
\hspace*{2cm}
 X_{12} = \widehat{X}_{12} + A_{11} A_{21}^-A_{22} + A_{11}F_{A_{21}}V +
 WE_{A_{21}}A_{22}, \hfill (31.6)
\cr}
$$
where $V$ and $W$ are two arbitrary matrices, $ \widehat{X}_{12} $ is chosen
such that
\begin{eqnarray*}
r( E_G\widehat{X}_{12}F_H ) & = & \min \{ \ r(E_G), \quad r( F_H ) \ \} \\
& = & \min \left\{ \, s_1 + r(A_{21}) - r\left[ \begin{array}{c}  A_{11} \cr A_{21}  \end{array} \right], \ \ \ 
 t_2+ r( A_{21}) - r[\, A_{21}, \ A_{22} \,]  \, \right\}, 
\end{eqnarray*}
where $ G = A_{11}F_{A_{21}} $ and $ H = E_{A_{21}} A_{22}$.
  
{\rm (b)}\, The minimal rank of $ M(X_{12}) $ with respect to $X_{12}$ is
$$ \displaylines{
\hspace*{2cm}
\min_{ X_{12} } r[M(X_{12})]  = r\left[ \begin{array}{c}  A_{11} \cr A_{21}  \end{array} \right] +
 r[ \, A_{21}, \ A_{22} \, ] - r( A_{21}), \hfill (31.7) 
\cr}
$$
and the matrix $X_{12}$ satisfying {\rm (31.7)} is exactly the general
solution of the following consistent linear matrix equation
$$
\displaylines{
\hspace*{2cm} 
E_G( X_{12} - A_{11}A_{21}^-A_{22}) F_H = 0, \hfill
\cr}
$$
which can be written as 
$$ \displaylines{
\hspace*{2cm}
X_{12} = A_{11}A_{21}^-A_{22} + A_{11}F_{A_{21}} V + WE_{A_{21}}A_{22},
\hfill (31.8)
\cr}
$$
where $V$ and $W$ are two arbitrary matrices. 

{\rm(c)}\, The matrix $X_{12}$ satisfying {\rm (31.7)} is unique if and
only if
$$ \displaylines{
\hspace*{2cm}
 r\left[ \begin{array}{cc} A_{11} \cr A_{21}  \end{array} \right] = r[ \,  A_{21}, \ A_{22} \, ] =
r( A_{21} ). \hfill
\cr}
$$
In that  case$,$  the unique matrix is $X_{12} = A_{11}A^-_{21}A_{22}$.

{\rm (d)}\, The rank of $M(X_{12})$ is invariant with respect to the
choice of $X_{12},$ if and only if
$$\displaylines{
\hspace*{2cm}
 r(A_{11}) = s_1  \quad   and \quad R(A_{11}^T) \cap R(A_{21}^T) =\{ 0 \},
 \hfill (31.9)
\cr
\hspace*{0cm}
and \hfill
\cr
\hspace*{2cm} r(A_{22}) = t_2 \quad  and  \quad  R(A_{21}) \cap R(A_{22})
= \{ 0 \}.  \hfill (31.10)
\cr}
$$ } 
{\bf Proof.}\, Applying (1.6) to $ M(X_{12})$ in (31.4), we first obtain
$$ 
r[ M(X_{12})] = r\left[ \begin{array}{cc}  A_{11} & X_{12} \cr A_{21} & A_{22} \end{array} \right] \\
 = r\left[ \begin{array}{c}  A_{11} \cr A_{21} \end{array} \right] + r[\, A_{21} , \ A_{22} \, ] -
 r( A_{21} ) + r[ \, E_G(X_{12} - A_{11}A^-_{21}A_{22})F_H \,], \eqno (31.11)
$$
where $ G = A_{11}F_{A_{21}} $ and $ H = E_{A_{21}}A_{22} $. Thus the maximal
 and the minimal ranks of $M(X_{12})$ subject to $X_{12}$ are, in fact,
 determined by the term $ E_G(X_{12} - A_{11}A^-_{21}A_{22})F_H$. It is quite
 easy to see that
 $$
\displaylines{
\hspace*{2cm}
\max_{X_{12}} r[ \, E_G(\, X_{12} - A_{11}A^-_{21}A_{22} \,)F_H \,] = \min
\{ \, r(E_G), \ \ r(F_H) \, \}, \hfill (31.12)
\cr}
$$
and the matrix $ X_{12}$ is given by (31.6). Moreover
$$
\displaylines{
\hspace*{2cm}
\min_{X_{12}} r[ \, E_G(\, X_{12} - A_{11}A^-_{21}A_{22} \, )F_H \,] = 0, \hfill (31.13)
\cr}
$$
and the matrix $ X_{12}$ is given by (31.8). Putting (31.12) and (31.13) in
(31.11) produces (31.5) and (31.7). The result in Part (c) is
direct consequence of Part (b). The invariance of the rank of $ M(X_{12}) $ 
subject to $X_{12} $ is equivalent to the fact  that (31.5) and (31.7) are 
equal, that is,
$$
\displaylines{
\hspace*{2cm}
 r\left[ \begin{array}{c} A_{11} \cr A_{21} \end{array} \right] = r(A_{21} ) + s_1 \quad  {\rm or} \quad  r[ \, A_{21} , \ A_{22} \,]
 = r(A_{21} ) + t_2. \hfill
\cr}
$$ 
Finally applying Lemma 1.2 to both of them leads to  (31.9) and (31.10).  \qquad $ \Box$

\medskip

Notice that $M(X_{12}) $ is the simplest case of $ A_n + X_n $ in (31.1)
corresponding to $ n = 2$. Thus Lemma 31.1 presents, in fact,  a complete
solution to the two problems in (31.2) and (31.3) when $ n = 2$. Our work
in the next two sections is actually to extend the results in Lemma 1.2 to 
$ n \times n$ case. \\

\noindent {\large {\bf 31.1. The maximal rank completion of $A_n+X_n$ }} \\

\noindent For convenience of representation, we adopt the following
notations for the block matrices in (31.1),
$$\displaylines{
\hspace*{2cm}
 P_i=[ \, A_{i1}, \, A_{i2}, \, \cdots, \, A_{i,i-1}\, ], \qquad
 i= 2, \ 3, \ \cdots, \  n,  \hfill (31.14)
\cr
\hspace*{2cm}
M_i=[ \,  A_{i1}, \, A_{i2} , \, \cdots , \, A_{ii} \, ] , \qquad
i=1, \, 2, \, \cdots , \,  n, \hfill (31.15)
\cr
\hspace*{2cm}
A_i = \left[ \begin{array}{ccc} A_{11} & & \cr
\vdots & \ddots & \cr
A_{i1} & \cdots & A_{ii} \end{array} \right], \qquad  i=1, \, 2, \,  \cdots ,\,  n,
\hfill (31.16)
\cr
\hspace*{2cm}
Q_{ij} = \left[ \begin{array}{ccc}  A_{i1} & \cdots & A_{ij} \cr
 \vdots & & \vdots \cr
A_{n1}& \cdots & A_{nj} \end{array} \right], \qquad  1 \leq i \leq j \leq n,
\hfill (31.17)
\cr
\hspace*{2cm} 
N_i = \left[ \begin{array}{c}  A_{ii} \cr A_{i+1,i} \cr \vdots \cr A_{ni} \end{array} \right],
\qquad Y_i = \left[ \begin{array}{c}  X_{1i} \cr X_{2i} \cr \vdots \cr X_{i-1, i}  \end{array} \right],
\quad i= 2, \ 3,\ \cdots, \  n, \hfill (31.18)
\cr
\hspace*{2cm}
 X_i =\left[ \begin{array}{cccc}  0 & X_{12} & \cdots & X_{1i} \cr
& 0 & \ddots & \vdots \cr
& &  \ddots & X_{i-1,i} \cr
& & & 0  \end{array} \right], \qquad  i = 2, \, 3,  \, \cdots , \,  n-1. \hfill (31.19)
\cr}
$$
>From (31.17) and (31.18) we see that
$$\displaylines{
\hspace*{2cm}
[\, Q_{i+1,i}, \ N_{i+1}  \, ] = Q_{i+1,i+1}, \qquad  i = 1,\,  2, \,
\cdots , \  n-1. \hfill (31.20)
\cr}
$$
Besides, we use 
$$ \displaylines{
\hspace*{2cm}
S_i = \{ \,  X_{1,i+1}, \, X_{ 2,i+2}, \,  \cdots , \, X_{n-i,n} \,  \},
\qquad i =1, \, 2, \, \cdots , \,  n-1 \hfill
\cr}
$$
to denote the set of the $n-i$ variant block entries in the $i$th upper
block subdiagonal of $ X_n $ in (31.1). 

\medskip

\noindent {\bf Theorem 31.2.}\, {\it Let $ A_n + X_n $ be given by
{\rm (31.2)}. Then the maximal rank of $A_n + X_n $ subject to
$X_n \in S( X_n) $ is
$$ \displaylines{
\hspace*{0cm}
\max_{X_n \in S(X_n)} r(\, A_n + X_n\,) \hfill
\cr
}
$$
$$
=  \min \left\{ \, r(Q_{11}) + ( s - k_1) + (t-l_1), \ r( Q_{22}) + (s - k_2) +
(t - l_2), \ \cdots, \ r( Q_{nn}) + (s - k_n) + (t - l_n) \, \right\},
\eqno (31.21)
$$
where $s$ and $t$ are the row number and column number of $A_n + X_n,$
respectively$,$ $ k_i = \sum_{j=i}^n s_j $ and $l_i = \sum^i_{j=1} t_j $
are the row number and column number of $Q_{ii} (i= 1, \, 2, \, \cdots, \, n)$
in {\rm (31.17)}$,$ respectively.}

\medskip

\noindent {\bf Proof.}\, By induction on $n$. When $n = 2, \, A_n + X_n $ in
(31.2) has the same form as $ M(X_{12})$ in Lemma 31.1, and the result
in (31.5) is exactly the result in (31.21) when $ n = 2$. Hence (32.21) is true
for $n=2$. Now suppose that (31.21) is true for $ A_{n-1} + X_{n-1}.$
 Then we next consider $n$. According to (31.14)---(31.19), $A_n + X_n$ in (31.2)
 can be partitioned as
$$
\displaylines{
\hspace*{2cm}
A_n + X_n = \left[ \begin{array}{cc}  A_{n-1} + X_{n-1} & Y_n \cr Q_{n,n-1} & N_n \end{array} \right].
\hfill
\cr}
$$
In that case, the maximal rank of $A_n + X_n $ subject to $ X_n \in S(X_n) $
can be calculated by the following two steps
$$ \displaylines{
\hspace*{2cm}
\max_{X_n \in S(X_n)} r(A_n +X_n )= \max_{X_{n-1} \in S(X_{n-1})}
\max_{Y_n} r(A_n+ X_n).
\hfill (31.22) 
\cr}
$$
applying (31.5),  we first find that
$$ 
\max_{Y_n} r(\,A_n+X_n\,)  =  \max_{Y_n} r\left[ \begin{array}{cc} A_{n-1} +X_{n-1} &
Y_n \cr Q_{n,n-1} & N_n  \end{array} \right] =  \min \left\{ r\left[ \begin{array}{c} A_{n-1} + X_{n-1}
\cr Q_{n,n-1} \end{array} \right]+ t_n, \ \ r(Q_{nn}) + (s - k_n )  \right\},
\eqno (31.23)
$$
and the matrix $Y_n$ satisfying (31.23) can be written as 
$$\displaylines{
\hspace*{2cm}
Y_n = \widehat{Y}_n +(\,A_{n-1} + X_{n-1}\, ) Q_{n,n-1}^- N_n + G V_n +
W_n H , \hfill
\cr}
$$
where $V_n$ and $W_n$ are two arbitrary matrices, $G=(A_{n-1} + X_{n-1} )
F_{Q_{n,n-1}},$ $ H=E_{Q_{n,n-1}} N_n, $ and $\widehat{Y}_n $ is chosen such
that
$$
\displaylines{
\hspace*{2cm} 
r(E_G \widehat{Y}_n F_H ) = \min \{ \, r(E_G), \quad r(F_H) \, \}.\hfill 
\cr}
$$
The next step for continuing (31.22) is to find the maximal rank of the
block matrix in (31.23) subject to $X_{n-1}.$ Observe that
$$
\displaylines{
\hspace*{2cm} 
\left[ \begin{array}{c} A_{n-1} + X_{n-1} \cr Q_{n, n-1} \end{array} \right] = \left[ \begin{array}{cccc}  A_{11} & X_{12} & \cdots & X_{1,n-1} \cr
A_{21} & A_{22} & \ddots & \vdots \cr
\vdots & \vdots & \ddots &  X_{n-2,n-1} \cr
B_1 & B_2 & \cdots & B_{n-1}  \end{array} \right], \hfill (31.24) 
\cr}
$$
where $ B_i = \left[ \begin{array}{c} A_{n-1,i} \cr A_{ni}  \end{array} \right], \, i = 1
, \, 2, \, \cdots, \, n-1.$ Hence (31.24) is, in fact, a new
$ (n-1) \times (n-1)$ block matrix with the same form as
$A_{n-1} + X_{n-1} $ in (31.2). Thus by hypothesis of induction, we know that
$$ \displaylines{
\hspace*{0.5cm} 
\max_{X_{n-1} \in S(X_{n-1})}  r \left[ \begin{array}{c} A_{n-1} + X_{n-1} \cr Q_{n,n-1} 
\end{array} \right] \hfill
\cr
\hspace*{0.5cm}
= \min \{ \, r(Q_{11}) + (\bar{t} - l_1), \quad  r(Q_{22}) + (s - k_2) + 
( \bar{t} -l_2), \quad \cdots, \quad r(Q_{n-1, n-1}) + (s - k_{n-1}) \,  \},
\hfill 
\cr}
$$
where $ \bar{t} = \sum_{i=1}^{n-1} t_i $ . Substituting it into (31.23)
yields
$$
\displaylines{
\hspace*{0.5cm}
\max_{X_n \in S(X_n)} r(\, A_n + X_n \,)
=  \min \{ \, r(Q_{11}) + ( \bar{t} - l_1) + t_n, \ \ 
r( Q_{22}) + ( s - k_2) + ( \bar{t} - l_2)  + t_n, \hfill
\cr
\hspace*{5cm} \cdots, \ \
r( Q_{n-1,n-1}) + ( s - k_{n-1}) + t_n , \ \ r( Q_{nn}) + (s - k_n) \, \}.
\hfill
\cr}
$$
Note that $ t = \bar{t} + t_n, \ s - k_1 = 0$ and $ t- l_n = 0$, thus the
above result is exactly the formula in (31.21).  \qquad  $ \Box$

\medskip

>From the proof of Theorem 31.2 we can also conclude a group of formulas for 
calculating the column block matrices $ Y_2, \, Y_3, \, \cdots, \, Y_n $ in the
matrix $ X_n $ satisfying (31.21).

\medskip

\noindent{\bf Theorem 31.3.}\, {\it The general expressions of the column
block entries
 $Y_2, \ Y_3, \ \cdots, \  Y_n $ in the matrix $X_n $ satisfying {\rm (31.21)}
  can be written in the inductive formulas
$$\displaylines{
\hspace*{2cm}  
Y_2 = \widehat{Y}_2 + A_{11}Q_{21}^- N_2 + G_2 V_2 + W_2 H_2 , \hfill (31.25)
\cr
\hspace*{2cm} 
 Y_i = \widehat{Y}_i + (\, A_{i-1} + X _{n-1} \,) Q^-_{i,i-1} N_i + G_i V_i
  + W_i H_i, \quad i = 3, \,  \cdots, \, n, \hfill (31.26)
\cr
\hspace*{0cm}
where \hfill
\cr
\hspace*{2cm}  
A_{i-1} + X_{i-1} = \left[ \begin{array}{cc}  A_{i-2} + X_{i-2} & Y_{i-1} \cr  P_{i-1} &
A_{i-1,i-1} \end{array} \right], \qquad i= 3, \  \cdots, \ n,  \hfill
\cr}
$$
$V_2, \, V_3, \, \cdots, \, V_n, \, W_2, \, W_3, \,  \cdots, \, W_n $ are arbitrary 
matrices, $G_2 = A_{11}F_{Q_{21}}, \,  H_2 = E_{Q_{21}} N_2 $ and 
$$\displaylines{
\hspace*{2cm} 
G_i = (\, A_{i-1} + X_{i-1}\, )F_{Q_{i,i-1}},  \quad H_i= E_{Q_{i,i-1}} N_i,
\qquad  i= 3, \, \cdots, \, n,  \hfill
\cr}
$$
meanwhile $ \widehat{Y}_2, \, \widehat{Y}_3, \, \cdots, \, \widehat{Y}_n$ are
chosen such that
$$\displaylines{
\hspace*{2cm}  
r( E_{G_i} \widehat{Y}_i F_{H_i}) = \min \{ \, r( E_{G_i}), \quad r(F_{H_i})
\, \}, \qquad i= 2, \, 3, \, \cdots, \, n.  \hfill
\cr}
$$ }
\hspace*{0.25cm} Substituting (31.25) and (31.26) into the matrix $X_n$ in (31.2) will produce
a general expression for the maximal rank completion of $ A_n + X_n $ in
(31.2).

On the basis of Theorems 31.3 and 31.3, now we are able to consider the
nonsingularity of $ A_n + X_n $ in (31.2) when it is a square block
matrix.

\medskip

\noindent {\bf Corollary 31.4.}\, {\it Suppose that $A_n + X_n$ in {\rm (31.2)} is a square 
block matrix of size $ t \times t $. Then there exists an $ X_n \in S(X_n ) $ such that 
$ A_n + X_n $ in {\rm (31.2)} is nonsingular, if and only if the block 
matrices $ Q_{11}, \,  Q_{22}, \,  \cdots,  \,  Q_{nn} $ in $A_n$ satisfy
$$ 
r( Q_{11}) = l_1, \quad  r( Q_{22}) \geq k_2 + l_2 - t, \quad \cdots,
\quad  r( Q_{n-1,n-1}) \geq k_{n-1} + l_{n-1} - t, \quad
r(Q_{nn}) = k_n,
$$
where $k_i$ and $l_i$ are$,$ respectively$,$ the row number and the column
number of $Q_{ii}(i=1,\, 2, \, \cdots, \, n).$ In that case$,$ the column block
matrices $ Y_2, \, Y_3, \, \cdots, \, Y_n $ in the matrix $ X_n $ such that
$ A_n + X_n $ is nonsingular are also given by {\rm (31.25)} and
{\rm (31.26)}. }

If the matrix $A_n $ in (31.1) satisfies some additional conditions, the
results in Theorems 31.2 and 31.3 can  further simplify. In particular,
when $A_n $ in (31.1) is a diagonal block matrix, we have the following
simple result.

\medskip

\noindent {\bf Corollary 31.5.}\, {\it Suppose that $A_n $ in {\rm (31.1)} is a diagonal 
block matrix$,$ i.e.$,$ $ A_{ij} = 0 (i>j) $ in {\rm (31.1)}. Then 
$$\displaylines{
\hspace*{0cm}
\max_{X_n \in S(X_n)} r(\, A_n + X_n \, ) \hfill
\cr
\hspace*{0cm}
=  \min \{ \, r(A_{11}) + ( s - k_1) + ( t - l_1), \ \  r(A_{22}) + ( s - k_2) + ( t - l_2), \  \cdots, \
r( A_{nn}) + ( s- k_n ) + ( t - l_n) \, \}, \hfill (31.27)
\cr}
$$
where $k_i $ and $ l_i(i = 1, \, 2, \, \cdots, \, n)$ are as in {\rm (31.21)}.
The column block entries $ Y_2, \, Y_3, \, \cdots, \, Y_n $  in the matrix
$X_n $ satisfying {\rm(31.27)} are given by inductive formulas
$$ \displaylines{
\hspace*{2cm} 
Y_2 = \widehat{Y}_2 + A_{11}V_1+ W_1A_{22}, \hfill
\cr
\hspace*{2cm} 
Y_i =\widehat{Y}_i + ( A_{n-1} + X_{n-1}) V_i+ W_iA_{ii}, \quad i=3, \, \cdots,\,  n,  \hfill
\cr
\hspace*{0cm}
where \hfill
\cr
\hspace*{2cm} 
A_{i-1} + X_{i-1} = \left[ \begin{array}{cc}  A_{i-2} + X_{i-2} & Y_{i-1} \cr 0 &
A_{i-1,i-1} \end{array} \right], \quad i= 3, \, \cdots, \, n,  \hfill
\cr}
$$ 
the matrices $ V_2, \, V_3, \, \cdots, \, V_n,\,  W_2, \, W_3, \, \cdots, \, W_n $
are arbitrary, and $ \widehat{Y}_2,$ $ \widehat{Y}_3,$ $ \cdots,$
$\widehat{Y}_n $ in them satisfy
$$\displaylines{
\hspace*{2cm}  
r( E_{A_{11}}\widehat{Y}_2 F_{A_{22}} ) = \min \{ \, r( E_{A_{11}}),
\quad r( F_{A_{22}}) \, \}, \hfill
\cr
\hspace*{2cm} 
r( E_{G_i} \widehat{Y}_i F_{A_{ii}}) = \min \{ \, r( E_{G_i}), \quad
r( F_{A_{ii}}) \,  \}, \quad i = 3, \, \cdots, \, n,  \hfill
\cr}
$$
where $G_i = A_{i-1} + X_{i-1},\, i = 3, \, \cdots, \, n. $}  \\

\noindent {\large {\bf 31.2. The minimal rank completion of $A_n+X_n$ }} \\

\noindent From the results in Lemma 31.1(b), we see that the variant entry $X_{12}$
such that $ M(X_{12})$ has its minimal rank is, in fact, the general
solution of a consistent linear matrix equation constructed by the given
matrices $A_{11},\ A_{21}$ and $A_{22}$ in $A_2$. It is not difficult to
find by repeatedly using Lemma 31.1(b) that the column block entries
$ Y_2, \ Y_3, \ \cdots, \ Y_n $ in the minimal rank completion of
$ A_n + X_n $ are also the general solutions of $n-1$ consistent linear
matrix equations constructed by the given block entries in $A_n$.

\medskip

\noindent {\bf Theorem 31.6.}\, {\it Let $A_n + X_n $ be given by
{\rm (31.3)}. Then 

{\rm (a)\cite{Wo1}}\, The minimal rank of $A_n +X_n $  subject to $X_n \in S(X_n) $ is 
$$\displaylines{
\hspace*{2cm}  
\min_{X_n \in S(X_n)} r( \, A_n + X_n \, ) = \sum^n_{i=1} r(G_{ii}) - \sum^
{n-1}_{i=1} r( Q_{i+1,i}) . \hfill (31.29)
\cr}
$$

{\rm (b)}\, The general expressions of the column block entries
$ Y_2, \ Y_3, \ \cdots, \ Y_n $ in the matrix $X_n$ satisfying {\rm (31.3)}
can be calculated by the inductive formulas
$$ \displaylines{
\hspace*{2cm} 
Y_2 = A_{11}Q^-_{21} N_2 + A_{11}F_{Q_{21}}V_2 + W_2E_{Q_{21}}N_2, \hfill (31.30)
\cr
\hspace*{2cm} 
Y_i = ( \, A_{i-1} + X_{i-1} \,) Q^-_{i,i-1}N_i + G_iV_i + W_iH_i, \quad i=3, \, 
\cdots, \, n, \hfill (31.31)
\cr
\hspace*{0cm}
where \hfill
\cr
\hspace*{2cm} 
 A_{i-1} + X_{i-1} = \left[ \begin{array}{cc}  A_{i-2} + X_{i-2} & Y_{i-1} \cr  P_{i-1} & A_{i-1,i-1} \end{array} \right], \quad i= 3, \ \cdots, \ n, \hfill
\cr}
$$
$V_2, \, V_3, \, \cdots, \, V_n, \, W_2, \, W_3, \, \cdots, \, W_n $ are arbitrary matrices, and
$$\displaylines{
\hspace*{2cm}  
G_i = (  \, A_{i-1} + X_{i-1} \, )F_{Q_{i,i-1}}, \quad H_i = E_{Q_{i,i-1}}N_i, \quad i= 3, \, \cdots, \, n.  \hfill
\cr}
$$} 
{\bf Proof.}\, According to  the structure of $X_n$ in (31.1) we determine the minimal rank completion of $A_n + X_n $ by
 the following $n-1$ steps
$$\displaylines{
\hspace*{2cm} 
\min_{X_n \in S(X_n)} r(  \, A_n + X_n \, ) = \min_{Y_2} \min_{Y_3} \, \cdots \, \min_{Y_n} r( A_n + X_n ).  \hfill (31.32)
\cr}
$$
Applying Lemma 31.1(b) we  first find that  
$$
\displaylines{
\hspace*{0cm}  
\min_{Y_n} r( \, A_n + X_n  \,)  =  \min_{Y_n} r\left[ \begin{array}{cc}  A_{n-1} + X_{n-1} & Y_n \cr Q_{n,n-1} & N_n \end{array} \right] \hfill
\cr
\hspace*{2.2cm}  
=  r(\, Q_{n, n-1}, \ N_n \,)-r( Q_{n,n-1}) + r\left[ \begin{array}{c}  A_{n-1} + X_{n-1} \cr Q_{n,n-1} \end{array} \right] \hfill 
\cr
\hspace*{2.2cm}  
=  r(\, Q_{nn}) - r( Q_{n,n-1}\, ) + r\left[ \begin{array}{c} A_{n-1}+ X_{n-1}
\cr Q_{n,n-1}, \end{array} \right],  \hfill (31.33)
\cr}
$$     
and the column block matrix $ Y_n $ satisfying (31.33) is 
$$\displaylines{
\hspace*{2cm} 
 Y_n = ( \, A_{n-1} + X_{n-1} \, )Q_{n,n-1}^-N_n + G_nV_n  + W_nH_n,
 \hfill  (31.34)
\cr
\hspace*{0cm} 
where \hfill
\cr
\hspace*{2cm} 
G_n = ( \, A_n + X_n  \,)F_{Q_{n,n-1}}N_n,  \ \ \  {\rm and}  \ \ \
H_n = E_{Q_{n,n-1}}N_n,  \hfill
\cr}
$$
$V_n $ and $W_n $ are two arbitrary matrices. Clearly (31.34) is exactly the
result in (31.31) when $ i=n.$ Observe that
$$
\displaylines{
\hspace*{2cm}  
 \left[ \begin{array}{c} A_{n-1} + X_{n-1} \cr Q_{n,n-1} \end{array} \right] = \left[ \begin{array}{cc}  A_{n-2}+ X_{n-2}
& Y_{n-1} \cr Q_{n-1,n-2} & N_{n-1} \end{array} \right]. \hfill 
\cr}
$$
Thus we find by Lemma 31.1(b) that
$$ \displaylines{
\hspace*{2cm} 
\min_{Y_{n-1}} r\left[ \begin{array}{c}  A_{n-1} + X_{n-1} \cr Q_{n,n-1} \end{array} \right]
= r( Q_{n-1,n-1})- r( Q_{n-1,n-2} ) + r\left[ \begin{array}{c}  A_{n-2} + X_{n-2} \cr
Q_{n-1,n-2} \end{array} \right], \quad  \hfill (31.35)
\cr}
$$ 
and $Y_{n-1}$ satisfying (31.35) is 
$$\displaylines{
\hspace*{2cm} 
Y_n = ( \,A_{n-2} + X_{n-2}\, )Q^-_{n-1,n-2}N_{n-1} + G_{n-1}V_{n-1} +
W_{n-1}H_{n-1}, \hfill  (31.36)
\cr
\hspace*{0cm}
where \hfill
\cr
\hspace*{2cm}
G_{n-1} = (\, A_{n-2} + X _{n-2}\, ) F_{ Q_{n-1,n-2}}, \ \ \  {\rm and }  \ \ \
 H_{n-1} = E_{Q_{n-1,n-2}}N_{n-1}, \ V_{n-1},  \hfill
\cr}
$$
and $W_{n-1}$ are two arbitrary matrices. Clearly (31.36) is exactly the
result in (31.31) when $i = n-1.$ By the same method, we can derive a group of
inductive formulas as follows
$$ \displaylines{
\hspace*{2cm} 
\min_{Y_i} r \left[ \begin{array}{c} A_i + X_i \cr  Q_{i+1,i}  \end{array} \right] = r( Q_{ii} ) -
r( Q_{i,i-1}) + r\left[ \begin{array}{c}  A_{i-1} + X_{i-1} \cr  Q_{i,i-1} \end{array} \right],
\qquad  i= n-1, \, \cdots, \, 3,  \hfill
\cr 
\hspace*{0cm} and \hfill
\cr
\hspace*{2cm} 
 \min_{Y_2} r\left[ \begin{array}{c}  A_2+ X_2  \cr  Q_{32} \end{array} \right] = r( Q_{11}) +
 r( Q_{22} ) - r( Q_{21}).
 \hfill
\cr}
$$ 
The general expressions of $Y_2, \, Y_3, \, \cdots, \, Y_n $ satisfying
the above group  of equalities are given in (31.30) and (31.31). Finally
substituting the above $ n-1$  equalities into (31.32)  results in (31.29).
\qquad $\Box$

\medskip

A particular concern for  the problem (31.3) is the uniqueness of $X_n$
satisfying (31.3), which was well examined in \cite{GKL}, \cite{Wo1} and
\cite{Wo2}.

\medskip

\noindent {\bf Corollary 31.7.}\, {\it Let $A_n + X_n $ be given by {\rm (31.2)}. 

{\rm (a)\cite{Wo1}}\, The matrix $X_n$ satisfying {\rm (31.1)} is unique$,$
that is$,$ the matrix $ A_n + X_n $ has a unique minimal rank completion$,$
if and only if
$$ \displaylines{
\hspace*{2cm} 
 r( Q_{ii}) = r( Q_{i+1,i}) = r( Q_{i+1, i+1} ), \qquad i= 1, \, 2,  \, \cdots, \, n-1, \hfill (31.37)
 \cr}
$$
In that case$,$ the minimal rank of $ A_n + X_n $ subject to
$ X_n \in S( X_n)$
is
$$\displaylines{
\hspace*{2cm}  
\min_{ X_n \in S(X_n)} r(\, A_n + X_n\, ) = r( Q_{11}).  \hfill (31.38)
\cr}
$$

{\rm (b) }\, Under the conditions in {\rm (31.37)}$,$ the column block entries 
$Y_2, \ Y_3, \, \cdots, \, Y_n $ in the unique minimal rank completion of
$A_n + X_n $ are given by the group of inductive formulas
$$\displaylines{
\hspace*{2cm} 
Y_2 = A_{11}Q_{21}^- N_2 , \qquad  Y_i = ( \,A_{i-1} + X_{i-1}\, )Q^-_{i,i-1}N_i, \quad \quad i= 3, \, \cdots, \, n,  \hfill
\cr
\hspace*{0cm} 
where \hfill
\cr
\hspace*{2cm} 
A_{i-1} + X_{i-1} = \left[ \begin{array}{cc}  A_{i-2} + X_{i-2}  & Y_{i-1} \cr P_{i-1} & 
A_{i-1,i-1} \end{array} \right], \quad \quad  i= 3, \, \cdots, \, n.  \hfill
\cr}
$$ }
{\bf Proof.}\,  Follows directly from (31.30) and (31.31). \qquad $\Box $   

\medskip

In general cases, the maximal rank and the minimal rank of $A_n + X_n$
subject to $X_n \in S(X_n )$ are different. If they are equal,
it implies that the rank of $A_n + X_n $ is invariant for any choice of
$ X_n \in S(X_n )$. The combination of Theorem 31.2 with Theorem 31.6
 yields  the following result.

\medskip

\noindent {\bf Corollary 31.8.}\, {\it Let $A_n + X_n $ be given by {\rm (31.2)}. 
Then the rank of $A_n + X_n $ is invariant with respect to the choice of 
$ X_n \in S(X_n) $ if and only if 
$$ \displaylines{
\hspace*{0cm} 
\sum_{i=1}^n r( Q_{ii}) - \sum_{ i=1}^{n-1} r( Q_{i+1, i}) \hfill
\cr
\hspace*{0cm}
= \min \{ r( Q_{11}) + (s - k_1 ) + ( t - l_1), \, 
r( Q_{22} ) + ( s - k_2) + ( t - l_2), \, \cdots, \,  r( Q_{nn}) +
( s - k_n ) + ( t - l_n )  \},  \ \hfill  (31.39)
\cr}
$$
where $s$ and $t$ are the row and the column numbers of $A_n + X_n, $
respectively$;$  $ k_i $ and $l_i $ are the row and the column numbers of
$Q_{ii} ( i = 1, \, 2, \, \cdots, \, n ),$ respectively. In that case$,$ any
sum of of $A_n + X_n $ can be regarded as
a maximal or a minimal rank completion of the triangular block matrix
$ A_n $. }

\medskip

The rank equality (31.39) can be written as some explicit
equivalent expressions. When $n=2$, the rank equality in (31.39) is Lemma 31.1(d). 
We next present a equivalent statement for (31.39) when $n=3$. The proof is analogous to
that of Lemma 31.1(d) and is, therefore, omitted.

\medskip

\noindent {\bf Corollary 31.9.}\, {\it Suppose that $A_3 + X_3 $ is a
$ 3 \times 3$ block matrix as follows
$$ 
A_3 + X_3 = \left[ \begin{array}{ccc}  A_{11} & X_{12} & X_{13} \cr
                         A_{21} & A_{22} & X_{23} \cr
                         A_{31} & A_{32} & A_{33} \end{array} \right],  \eqno  (31.40)
$$
Then the rank of $ A_3 + X_3 $ is invariant with respect to  the choice of
$ X_{12}, \ X_{13} $ and $ X_{23} $, if and only if one of the following
three groups of conditions is satisfied:

{\rm (a) }\, $A_{11}$ and $ [ \, A_{21}, \ A_{22} \,]$ have  full row
ranks$,$ respectively$,$ and
$$ 
R(A_{11}^T) \cap R[ \, A_{21}^T, \ A_{31}^T \, ] = \{  0 \}, \ \ and  \ \  
R \left[ \begin{array}{c}  A_{21}^T \cr A_{22}^T \end{array} \right] \cap R \left[ \begin{array}{c}  A_{31}^T \cr A_{32}^T 
\end{array} \right] = \{ 0 \}.
$$

 {\rm (b) }\, $A_{11}$ has full row rank$,$ $A_{33}$ has full column rank$,$
  and
$$ 
R(A_{11}^T) \cap R[ \, A_{21}^T, \ A_{31}^T \, ] = \{ 0 \}, \ \ and  \ \
R[ \,  A_{31}, \  A_{32} \, ] \cap R( A_{33} ) = \{ 0\}.
$$ 
 
{\rm (c) }\, $ A_{33}$ and $ \left[ \begin{array}{c} A_{22} \cr A_{32} \end{array} \right] $ has full
column ranks$,$ respectively$,$ and
$$ 
R\left[ \begin{array}{c} A_{21} \cr A_{31} \end{array} \right] \cap R\left[ \begin{array}{c}  A_{22} \cr A_{32} \end{array} \right]
= \{  0 \},  \ \ and  \ \  R[ \, A_{31}, \ A_{32} \, ]  \cap R( A_{33} )=
\{  0 \}.
$$ }
{\bf Remark 31.10.}\, Besides the partitioning method for $X_n $
shown in (31.18) and (31.19), we can also inductively partition $ X_n $ in (31.2) into
$$ 
X_n = \left[ \begin{array}{cc}  0 & Z_1 \cr 0 & X_{n-1} \end{array} \right], \qquad \cdots,  \qquad
X_3 = \left[ \begin{array}{cc}  0 & Z_{n-2} \cr 0 & X_2 \end{array} \right], \qquad
 X_2 = \left[ \begin{array}{cc}  0 & X_{n-1,n} \cr 0 & 0
 \end{array} \right], 
$$
where $ Z_i = [\, X_{i,i+1}, \, X_{i,i+2}, \, \cdots, \ X_{in} \,], \, i = 1, \,
2, \, \cdots, \, n-1.$ In the same method for deducing the results in
Sections 31.1 and 31.2, we can also find general expressions of the row block
matrices $ Z_1, \ Z_2, \ \cdots, \ Z_{n-1} $  to the two problems in
(31.2) and (31.3). The corresponding results are much analogous to those in
Theorems 31.2 and 31.3, and are, therefore, omitted here.

In addition to the two methods used or mentioned above, another method
available for constructing  maximal and minimal rank completions of
$A_n + X_n $ in (31.2) and (31.3) is inductively calculating the block
entries in each upper block subdiagonal of $ X_n $. We next illustrate this
method by constructing  maximal and minimal rank completions of
$A_3 + X_3 $ in (31.40).

According to Lemma 31.1(a) and (b), the maximal and the minimal ranks of
$A_3 +X_3$ in (31.40) with respect to $X_{13}$ respectively are
$$ \displaylines{
\hspace*{2cm} 
 \max_{X_{13}} r(\, A_3 + X_3\, ) = \min  \left\{ \, r\left[ \begin{array}{cc}  A_{11} &
 X_{12} \cr Q_{21} & N_2 \end{array} \right] + t_3, \quad r\left[ \begin{array}{cc} M_2 & X_{23} \cr
 Q_{32} & A_{33} \end{array} \right]+ s_1 \,  \right\},  \hfill (31.41)
\cr}
$$
and
$$ \displaylines{
\hspace*{2cm} 
\min_{X_{13}}r(\, A_3+ X_3\,) = r\left[ \begin{array}{cc}  A_{11} & X_{12} \cr Q_{21} & N_2 \end{array} \right]
+r\left[ \begin{array}{cc}  M_2 & X_{23} \cr Q_{32} & A_{33} \end{array} \right] -r(Q_{22}). \hfill (31.42)
\cr}
$$
The general expressions of $X_{13} $ satisfying (31.41) and (31.42) can
respectively be derived from (31.6) and (31.8). Observe that the variant
entries $X_{12}$ and $X_{23} $ occur in two independent block matrices in
(31.41) and (31.42).  Applying Lemma 31.1(a) and (b) to the block matrices in
(31.41) and (31.42), we easily obtain
\begin{eqnarray*} 
\max_{X_3 \in S(X_3)} r(A_3 + X_3) & = & \min \left\{ \max_{X_{12}}
 r\left[ \begin{array}{cc}  A_{11} & X_{12} \cr Q_{21} & N_2 \end{array} \right]  + t_3, \quad  \max_{X_{23}} 
r\left[ \begin{array}{cc}  M_2 & X_{23} \cr Q_{32} & A_{33} \end{array} \right] +s_1  \right\} \\
  &=& \min \{ r( Q_{11}) + t_2 + t_3, \quad r( Q_{22} ) + s_1 + t_3,
  \quad r(Q_{33}) + s_1 + s_2 \},
\end{eqnarray*} 
and   
\begin{eqnarray*} 
\min_{X_3 \in S(X_3)} r( A_3+ X_3) &=& \min_{X_{21}} r\left[ \begin{array}{cc}  A_{11} & X_{12} \cr Q_{21} & N_2 \end{array} \right] + \min_{X_{23}} r \left[ \begin{array}{cc}  M_2 & X_{23} \cr Q_{32} & A_{33} \end{array} \right] - r( Q_{22} ) \\
& = & r( Q_{11} ) + r( Q_{22} ) + r( Q_{33} ) - r( Q_{21} ) - r( Q_{32} ).
\end{eqnarray*} 
The matrices $X_{12} $ and $ X_{23} $ satisfying the above two equalities
can respectively be derived from (31.6) and (31.7). Clearly the above
two equalities are exactly (31.27) and (31.29) when $n = 3$ in them.

It is easy to conclude from the above example that the two completion
problems can also be constructed by inductively calculating the block
entries in the  upper block subdiagonals in $X_n$. Speaking precisely, the
first step of this work is to find the $n-1$ block entries in the set
$ S_1 = \{ \, X_{21}, \, X_{23}, \, \cdots, \, X_{n-1,n} \, \}$ of $X_n $
 such that
$$
r\left[ \begin{array}{cc}  M_1 & X_{12} \cr Q_{21} & N_2 \end{array} \right] = \max, \qquad
r \left[ \begin{array}{cc}  M_2 & X_{23} \cr Q_{32} & N_3 \end{array} \right] = \max, \ \ \ \cdots,  \ \ \
r \left[ \begin{array}{cc}  M_{n-1} & X_{n-1,n} \cr Q_{n, n-1} & N_n \end{array} \right] = \max,
$$ 
and 
$$ 
r\left[ \begin{array}{cc} M_1 & X_{12} \cr Q_{21} & N_2 \end{array} \right] = \min,  \qquad
r \left[ \begin{array}{cc} M_2 & X_{23} \cr Q_{32} & N_3 \end{array} \right]= \min,  \ \ \ \cdots,
\ \ \ r\left[ \begin{array}{cc} M_{n-1} & X_{n-1,n} \cr Q_{n,n-1} & N_n  \end{array} \right] = \min
$$
hold, respectively, where $ M_i, \ Q_{i,i-1}$ and $N_i$ are defined in
(31.14)---(31.17). The second step is to substitute the $n-1$ given block
entries in $S_1$ into $ A_n + X_n $ and then to determine the $n-2$ block
entries in  the  set $S_2 = \{ X_{13}, \ X_{24}, \ \cdots, $  $ X_{n-2,n} \}
$ of $ X_n $ such that
$$ 
r \left[ \begin{array}{cc} * & X_{13} \cr * & * \end{array} \right] = \max, \qquad
r \left[ \begin{array}{cc} * & X_{24} \cr * & * \end{array} \right] = \max, \qquad  \cdots,  \qquad
r \left[ \begin{array}{cc} * & X_{n-2,n} \cr * & * \end{array} \right] = \max
$$
and
$$ 
r \left[ \begin{array}{cc} * & X_{13} \cr * & * \end{array} \right] = \min, \qquad
r \left[ \begin{array}{cc} * & X_{24} \cr
* & * \end{array} \right] = \min, \qquad  \cdots, \qquad  r \left[ \begin{array}{cc} * & X_{n-2,n }
\cr * & * \end{array} \right] =   \min
$$
hold, respectively, where $ \left[ \begin{array}{cc} * & X_{ij} \cr * & * \end{array} \right] $ denotes
the submatrix in $ A_n + X_n $ as follows
$$ 
\left[ \begin{array}{cc} * & X_{ij} \cr * & *  \end{array} \right] = \left[ \begin{array}{ccc} A_{i1} & \cdots & X_{ij} 
\cr \vdots & & \vdots \cr A_{n1} & \cdots & A_{nj}  \end{array} \right], \quad  1 \leq i < j \leq n. 
$$
Next substituting the $ n-2 $ given block matrices in $ S_2$ into
$A_n + X_n$ and repeating the analogous calculations produces inductively
the block entries in the 3rd, 4th, ..., $(n-1)$th upper block subdiagonals
in $ X_n $. Finally substituting all the given entries in the $n-1$ upper
block subdiagonals of $X_n$ into $A_n + X_n $ yields maximal and minimal
rank completions of $A_n + X_n $, respectively.

Summing up all the results and discussion given in Sections 31.1 and 31.2,
we see that there are three kinds of general methods available for
constructing  maximal and minimal rank completions of  a triangular block
matrix:
\begin{enumerate}
\item[(i)]  Calculating inductively the unspecified column block entries
$ Y_2, \, Y_3, \cdots.$ $ \, Y_n $ in $ A_n +X_n.$

\item[(ii)] Calculating inductively the unspecified row block entries
$ Z_1, \, Z_2,  \, \cdots,$ $ \, Z_{n-1}$ in $ A_n + X_n $.

\item[(iii)] Calculating inductively the unspecified upper block
subdiagonal entries in $A_n + X_n.$ 
\end{enumerate} 

\bigskip

\noindent {\large {\bf 31.3. Extreme ranks of $ A - BXC $ when
$ X$ is a triangular block matrix}} \\

\noindent Let $ A - BXC$ be a matrix expression over an arbitrary field
 ${\cal F}$, and suppose $ X $ is a variant $ p \times p$ upper triangular
 block matrix. In that case, $ A- BXC$ can also be written as
 $$
 A - BXC = A - [\, B_1, \, B_2, \, \cdots, \, B_p \,]
  \left[ \begin{array}{cccc}
 X_{11} & X_{12} & \cdots & X_{1p} \\
 & X_{22} & \cdots & X_{2p} \\
& & \ddots & \vdots \\
& &  & X_{pp} \end{array} \right]
\left[ \begin{array}{c}
 C_1 \\ C_2  \\ \vdots \\ C_p \end{array} \right]  , \eqno (31.43)
$$
where $ A \in {\cal F}^{m \times n}, \, B \in {\cal F}^{m \times k}, \, C \in {\cal F}^{l \times n}, \, 
B_i \in {\cal F}^{m \times k_i}, \,  C_i \in {\cal F}^{l_i \times n}, \,  X_{ij} \in {\cal F}^{k_i \times l_j}
(1 \leq i \leq j \leq p)$.

In this section, we consider how to determine the maximal and the minimal ranks of the
matrix expression (31.43) with respect to  all variant blocks $X_{ij}$. This work is 
motivated by  some earlier work on triangular block solutions of matrix equations, and triangular block generalized
inverses of matrices (see, e.g., \cite{MJ0}, \cite{MJ1}). 

It is easy to verify that the rank of $ A - BXC $ can be expressed as
the rank of a block matrix as follows
$$
r(\, A - BXC \, ) = r \left[ \begin{array}{ccc}
 0 & I_k & X \\
 C & 0 & I_l \\ -A & B & 0 \end{array} \right] - k - l.
$$
Putting (31.43) in it, we get
$$\displaylines{
\hspace*{1cm}
r(\,A - BXC \,) = r\left[ \begin{array}{ccccccccc}
0 & I_{k_1} & 0 & \cdots& 0 & X_{11} & X_{12} & \cdots & X_{1p} \\
0 & 0 & I_{k_2} & \cdots & 0 & 0 & X_{22} & \cdots & X_{2p} \\
\vdots & \vdots & \vdots & \ddots & \vdots & \vdots & \vdots & \ddots &
\vdots \\
0 & 0 & 0 & \cdots & I_{k_p} & 0 & 0 & \cdots & X_{pp}\\
C_1 & 0 & 0 & \cdots & 0 & I_{l_1} & 0 & \cdots & 0 \\
C_2 & 0 & 0 & \cdots & 0 & 0 & I_{l_2} & \cdots & 0 \\
\vdots & \vdots & \vdots & \ddots & \vdots & \vdots & \vdots& \ddots &
\vdots \\
C_p & 0 & 0 & \cdots & 0 & 0 & 0 & \cdots & I_{l_p} \\
-A & B_1 & B_2 & \cdots & B_p & 0 & 0 & \cdots & 0
\end{array} \right] - k -l \hfill
\cr
\hspace*{3.2cm}
= r\left[ \begin{array}{ccccc} G_{11} & X_{11} & X_{12} & \cdots & X_{1p} \\
G_{21} & G_{22} & X_{22} & \cdots & X_{2p} \\
\vdots & \vdots & \vdots  & \ddots & \vdots \\
G_{p1} & G_{p2} & G_{p3} & \cdots & X_{pp} \\
G_{p+1,1} & G_{p+1,2}  & G_{p+1,3}  & \cdots & G_{p+1,p+1} \end{array} \right] -k -l.
 \hfill(31.44)
\cr}
$$
\hspace*{0.25cm} This equality shows that the  maximal and the minimal ranks of $ A - BXC$ in (31.43) can
completely determined by those of the block matrix in
(31.43). Applying (31.27) and (31.29) to the block matrix in (31.44) and then simplifying, we obtain  the following.

\medskip

\noindent { \bf Theorem 31.11.}\, {\em Let $A - BXC$ be given by {\rm (31.43)}$,$ and let
$$
\displaylines{
\hspace*{2cm}
\widehat{B}_i = [\, B_1, \, \cdots, \, B_i \,], \ \ \ \ \widehat{C}_i = 
\left[ \begin{array}{c} C_i \\ \vdots \\ C_p \end{array} \right],
\ \ \ \  i = 1, \, 2  \, \cdots, \,  p. \hfill
\cr 
\hspace*{0cm}
Then  \hfill
\cr
\hspace*{0cm}
\max_{X_{ij}} r( \, A - BXC \, )  = \min\left\{
r \left[ \begin{array}{c} A \\ \widehat{C}_1 \end{array} \right], \ \
r\left[ \begin{array}{cc} A & \widehat{B}_1 \\ \widehat{C}_2 & 0 \end{array}
\right], \ \  \cdots, \ \ r\left[ \begin{array}{cc} A & \widehat{B}_{p-1}
\\ \widehat{C}_p  & 0 \end{array} \right],  \ \ r[\, A, \ \widehat{B}_p \,]
\right\}, \hfill (31.45)
\cr
\hspace*{0cm} 
and \hfill
\cr
\hspace*{0cm}
\min_{X_{ij}} r( \, A - BXC \, ) \hfill
\cr}
$$
$$
= r \left[ \begin{array}{c} A \\
\widehat{C}_1 \end{array} \right] + r\left[\begin{array}{cc} A &
\widehat{B}_1 \\ \widehat{C}_2  & 0 \end{array} \right] + \cdots +
r\left[ \begin{array}{cc} A & \widehat{B}_{p-1} \\ \widehat{C}_p & 0 
\end{array} \right] + r[\, A, \ \widehat{B}_p \,] - r\left[
\begin{array}{cc} A & \widehat{B}_1 \\ \widehat{C}_1 & 0
\end{array} \right] - \cdots - r\left[ \begin{array}{cc} A & \widehat{B}_p
 \\ \widehat{C}_p & 0 \end{array} \right]. \eqno (31.46)
$$
In particular$,$ the upper triangular block matrix $ X $ satisfying {\rm (31.45)} is unique if and only if
$$
r(G_{ii}) =r(G_{i+1, i}) =r(G_{i+1, i+1}), \ \ \ \ i = 1, \, 2, \, \cdots, \, p, \eqno (31.47)
$$ 
where $ G_{ij}$ is defined in {\rm (31.44)}.
}

\medskip

For simplicity, the steps for presenting (31.45) and (31.46) are omitted. 
Furthermore the upper triangular block matrix $ X $ satisfying (31.45) and (31.46) 
can respectively be determined by the three general methods presented in Sections 31.1 and 2. 
We also omit them here for simplicity.   

Two direct consequences of the formula (31.46) are given below.

\medskip

\noindent { \bf Corollary 31.12.}\, {\em The matrix equation
$$\displaylines{
\hspace*{1.5cm}
[\, B_1, \, B_2, \, \cdots, \, B_p \,]
  \left[ \begin{array}{cccc}
 X_{11} & X_{12} & \cdots & X_{1p} \\
 & X_{22} & \cdots & X_{2p} \\
& & \ddots & \vdots \\
& &  & X_{pp} \end{array} \right]
\left[ \begin{array}{c}
 C_1 \\ C_2  \\ \vdots \\ C_p \end{array} \right] = A \hfill (31.48)
\cr}
$$
is consistent if and only if $ R( A) \subseteq R(B), $ $ R( A^T)
\subseteq R(C^T),$ and 
$$\displaylines{
\hspace*{1.5cm}
r\left[ \begin{array}{cccc} A  & B_1 & \cdots & B_i \\
C_{i+1} & 0  & \cdots & 0 \\
\vdots & \vdots &  & \vdots \\
C_p & 0 & \cdots & 0 \end{array} \right]  = r\left[ \begin{array}{c}
C_{i+1} \\ \vdots \\ C_p \end{array} \right] +
r[\, B_1, \ \cdots, \ B_i \,], \ \ \ \  i = 1, \, 2  \, \cdots, \,  p-1. 
 \hfill (31.49)
\cr}
$$}
{\bf Proof.}\, Let the right-hand side of (31.46) be zero and then
simplify to yield the desired result. \qquad $ \Box$

The general expressions of $ X_{ij} $ satisfying (31.48) can also be determined by 
the general method presented in Sections 31.2. 

\medskip

It is well known that an inner inverse of a matrix $ A $ is a solution
 to the matrix equation $ AXA = A$. Thus applying Corollary 31.12 to the equation
$ AXA = A$, we obtain the following.

\medskip

\noindent { \bf Corollary 31.13.}\,  {\em Let
$$\displaylines{
\hspace*{2cm}
 A = \left[ \begin{array}{cccc}
 A_{11} & A_{12} & \cdots & A_{1p} \\
A_{21} & A_{22} & \cdots & A_{2p} \\
\vdots & \vdots & \ddots & \vdots \\
A_{p1} & A_{p2} & \cdots & A_{pp} \end{array} \right]. \hfill (31.50)
\cr}
$$
Then $ A $ has an inner inverse with the upper triangular block form
$$\displaylines{
\hspace*{2cm}
A^- = \left[ \begin{array}{cccc}
 S_{11} & S_{12} & \cdots & S_{1p} \\
 & S_{22} & \cdots & S_{2p} \\
 &  & \ddots & \vdots \\
 &  &  & S_{pp} \end{array} \right], \hfill(31.51)
\cr}
$$ 
if and only if
$$
\displaylines{
\hspace*{2cm}
r(A) = r(Q_{1k}) + r(Q_{k+1,p}) - r(Q_{k+1,k}), \ \ \ \ k = 1, \, 2, \,
\cdots, \, p-1,  \hfill (31.52)
\cr}
$$
where
$$
\displaylines{
\hspace*{2cm}
Q_{ij} = \left[ \begin{array}{ccc}
 A_{i1} & \cdots & A_{ij} \\
\vdots & \ddots & \vdots \\
A_{p1} & \cdots  & A_{pj}
\end{array} \right] \in {\cal F}^{ s_i \times t_j} ,
\ \ \ \ 1 \leq i ,\ j \leq p. \hfill
\cr}
$$ 
In particular$,$ an upper triangular block matrix
$$\displaylines{
\hspace*{2cm}
A = \left[ \begin{array}{cccc}
 A_{11} & A_{12} & \cdots & A_{1p} \\
 & A_{22} & \cdots & A_{2p} \\
 &  & \ddots & \vdots \\
 &  &  & A_{pp} \end{array} \right] \hfill
\cr}
$$
has an upper triangular block inner inverse of the form {\rm (31.51)}$,$ if and only if  
$$
r(A) = r(\widehat{A}_{ii}) + r(\tilde{A}_{ii}), \ \ \ i = 1, \, 2, \,
\cdots, \, p-1.
$$
where
$$\displaylines{
\hspace*{2cm}
\widehat{A}_{ii} = \left[ \begin{array}{ccc}
  A_{11} & \cdots & A_{1i} \\
  & \ddots & \vdots \\
  &  & A_{ii} \end{array} \right], \ \ \ and  \ \ \
\tilde{A}_{ii} = \left[ \begin{array}{ccc}
 A_{i+1, i+1} & \cdots & A_{i+1,p} \\
  & \ddots & \vdots \\
  &  & A_{pp} \end{array} \right]. \hfill
\cr}
$$
}

\vspace*{2cm}

\noindent  Department of Mathematics and Statistics \\
Queen's University \\
Kingston,  Ontario,  Canada K7L 3N6\\
e-mail address: ytian@mast.queensu.ca

\end{document}